\documentclass{amsart}
\usepackage{amssymb}
\usepackage{amsfonts}
\usepackage{geometry}

\newtheorem{theorem}{Theorem}
\theoremstyle{plain}

\newtheorem{conclusion}[theorem]{Conclusion}

\newtheorem{corollary}[theorem]{Corollary}

\newtheorem{definition}[theorem]{Definition}

\newtheorem{lemma}[theorem]{Lemma}
\newtheorem{notation}[theorem]{Notation}
\newtheorem{problem}[theorem]{Problem}
\newtheorem{proposition}[theorem]{Proposition}
\newtheorem{remark}[theorem]{Remark}

\numberwithin{equation}{section}
\input{tcilatex}
\begin{document}
\title[$Tb$ theorem]{A two weight local $Tb$ theorem for the Hilbert
transform}

\begin{abstract}
We obtain a two weight local $Tb$ theorem for any elliptic and gradient
elliptic fractional singular integral operator $T^{\alpha }$ on the real
line $\mathbb{R}$, and any pair of locally finite positive Borel measures $%
\left( \sigma ,\omega \right) $ on $\mathbb{R}$. The Hilbert transform is
included in the case $\alpha =0$, and is bounded from $L^{2}\left( \sigma
\right) $ to $L^{2}\left( \omega \right) $ if and only if the Muckenhoupt
and energy conditions hold, as well as $b_{Q}$ and $b_{Q}^{\ast }$ testing
conditions over intervals $Q$, where the families $\left\{ b_{Q}\right\} $
and $\left\{ b_{Q}^{\ast }\right\} $ are $p$-weakly accretive for some $p>2$%
. A number of new ideas are needed to accommodate weak goodness, including a
new method for handling the stubborn nearby form, and an additional corona
construction to deal with the stopping form. In a sense, this theorem
improves the $T1$ theorem obtained by the authors and M. Lacey in the two
part paper \cite{LaSaShUr3},\cite{Lac}.
\end{abstract}

\author[E.T. Sawyer]{Eric T. Sawyer}
\address{ Department of Mathematics \& Statistics, McMaster University, 1280
Main Street West, Hamilton, Ontario, Canada L8S 4K1 }
\email{sawyer@mcmaster.ca}
\thanks{Research supported in part by NSERC}
\author[C.-Y. Shen]{Chun-Yen Shen}
\address{Department of Mathematics,\\
National Taiwan University, 10617, Taipei}
\email{cyshen@math.ntu.edu.tw}
\thanks{C.-Y. Shen supported in part by MOST, through grant 104-2628-M-002
-015 -MY4}
\author[I. Uriarte-Tuero]{Ignacio Uriarte-Tuero}
\address{ Department of Mathematics \\
Michigan State University \\
East Lansing MI }
\email{ignacio@math.msu.edu}
\thanks{ I. Uriarte-Tuero has been partially supported by grant
MTM2015-65792-P (MINECO, Spain).}
\date{September 25, 2017}
\maketitle
\tableofcontents

\section{Introduction}

The original $T1$ theorem of David and Journ\'{e} \cite{DaJo}, which
characterized boundedness of a singular integral operator by testing over
indicators $\mathbf{1}_{Q\text{ }}$ of cubes $Q$, was quickly extended to a $%
Tb$ theorem by David, Journ\'{e} and Semmes \cite{DaJoSe}, in which the
indicators $\mathbf{1}_{Q\text{ }}$ were replaced by testing functions $b%
\mathbf{1}_{Q\text{ }}$ for an accretive function $b$, i.e. $0<c\leq \func{Re%
}b\leq \left\vert b\right\vert \leq C<\infty $. Here the accretive function $%
b$ could be chosen to adapt well to the operator at hand, resulting in
almost immediate verification of the $b$-testing conditions, despite
difficulty in verifying the $1$-testing conditions. One motivating example
of this phenomenon is the boundedness of the Cauchy integral on Lipschitz
curves\footnote{%
The problem reduces to boundedness on $L^{2}\left( \mathbb{R}\right) $ of
the singular integral operator $C_{A}$ with kernel $K_{A}\left( x,y\right)
\equiv \frac{1}{x-y+i\left( A\left( x\right) -A\left( y\right) \right) }$,
where the curve has graph $\left\{ x+iA\left( x\right) :x\in \mathbb{R}%
\right\} $. Now $b\left( x\right) \equiv 1+iA^{\prime }\left( x\right) $ is
accretive and the $b$-testing condition $\int_{I}\left\vert C_{A}\left( 
\mathbf{1}_{I}b\right) \left( x\right) \right\vert ^{2}dx\leq \mathfrak{T}%
_{H}^{b}\left\vert I\right\vert $ follows from $\left\vert C_{A}\left( 
\mathbf{1}_{I}b\right) \left( x\right) \right\vert ^{2}\approx \ln \frac{%
x-\alpha }{\beta -x}$, for $x\in I=\left[ \alpha ,\beta \right] $. In the
case of a $C^{1,\delta }$ curve, the kernel $K_{A}$ is $C^{1,\delta }$ and
any $Tb$ theorem applies with $T=C_{A}$ and $\sigma =\omega =dx$, to show
that $C_{A}$ is bounded on $L^{2}\left( \mathbb{R}\right) $.}. See e.g. \cite%
[pages 310-316]{Ste}.

Subsequently, M. Christ \cite{Chr} obtained a far more robust \emph{local} $%
Tb$ theorem in the setting of homogeneous spaces, in which the testing
functions could be further specialized to $b_{Q}\mathbf{1}_{Q\text{ }}$,
where now the accretive functions $b_{Q}$ can be chosen to \emph{differ} for 
\emph{each} cube $Q$. Applications of the local $Tb$ theorem included
boundedness of layer potentials, see e.g. \cite{AAAHK} and references there;
and the Kato problem, see \cite{HoMc}, \cite{HoLaMc} and \cite{AuHoLaMcTc}:
and many authors, including G. David \cite{Dav1}; Nazarov, Treil and Volberg 
\cite{NTV3}, \cite{NTV2}; Auscher, Hofmann, Muscalu, Tao and Thiele \cite%
{AuHoMuTaTh}, Hyt\"{o}nen and Martikainen \cite{HyMa}, and more recently
Lacey and Martikainen \cite{LaMa}, set about proving extensions of the local 
$Tb$ theorem, for example to include a single upper doubling weight together
with weaker upper bounds on the function $b$. But these extensions were
modelled on the `nondoubling' methods that arose in connection with upper
doubling measures in the analytic capacity problem, see Mattila, Melnikov
and Verdera \cite{MaMeVe}, G. David \cite{Dav1}, \cite{Dav2}, X. Tolsa \cite%
{Tol}, and alsoVolberg \cite{Vol}, and were thus constrained to a single
weight - a setting in which both the Muckenhoupt and energy conditions
follow from the upper doubling condition.

In this paper, we consider only the case of dimension $n=1$, and we adapt
methods from the theory of two weight $T1$ theorems, which arose from \cite%
{NTV4}, \cite{Vol}, \cite{LaSaShUr3}, \cite{Lac}, \cite{SaShUr7} and \cite%
{SaShUr9}, and were used in \cite{HyMa} as well, to prove a two weight local 
$Tb$ theorem. These methods involve the `testing' perspective toward
characterizing two weight norm inequalities for an operator $T$. As
suggested by work originating in \cite{DaJo} and \cite{Saw3}, it is
plausible to conjecture that a given operator $T$ is bounded from one
weighted space to another if and only if both it and its dual are bounded
when tested over a suitable family of functions related geometrically to $T$%
, e.g. testing over indicators of intervals for fractional integrals $T$ as
in \cite{Saw3}.

\textbf{Muckenhoupt conditions}: However, for even the simplest singular
integral, the Hilbert transform, testing over indicators of intervals no
longer suffices\footnote{%
consider e.g. $d\omega \left( x\right) =\delta _{0}\left( x\right) $ and $%
d\sigma \left( x\right) =\left\vert x\right\vert dx$.}, and an additional
`side condition' on the weight pair is required - namely the Muckenhoupt $%
\mathfrak{A}_{2}$ condition, a simpler form of which was shown by Hunt,
Muckenhoupt and Wheeden \cite{HuMuWh} to characterize the one weight
inequality for the Hilbert transform. This side condition is a size
condition on the weight pair that is typically shown to be necessary by
testing over so-called tails of indicators of intervals, and indeed is known
to be necessary for boundedness of a broad class of fractional singular
integrals that are `strongly elliptic'. Using this side condition of
Muckenhoupt, the solution of the NTV conjecture, due to the authors and M.
Lacey in the two part paper \cite{LaSaShUr3}-\cite{Lac}, shows that the
Hilbert transform $H$ is bounded between weighted $L^{2}$ spaces if and only
if the Muckenhoupt condition and the two testing conditions over indicators
of intervals all hold. However, the testing conditions for singular
integrals, unlike those for positive operators such as fractional integrals,
are extremely unstable and in principle difficult to check \cite{LaSaUr2}.
On the other hand, given a weight pair, it may be possible to produce a
family of testing functions adapted to intervals on which the boundedness of
the operator is evident. In such a case, one would like to conclude that
finding an appropriately \emph{nondegenerate} family of such testing
functions, for which the corresponding testing conditions hold, is enough to
guarantee boundedness of the operator - bringing us back to a local $Tb$
theorem. In any event, one would in general like to understand the weakest
testing conditions that are sufficient for two weight boundedness of a given
operator.

\textbf{Energy conditions}: Our $Tb$ theorem lies in this direction, but the
method of proof requires in addition a second `side condition', namely the
so-called energy condition, introduced in \cite{LaSaUr2}. The energy
condition is necessary for the boundedness of the Hilbert transform, and
actually follows there from testing over indicators of intervals and,
through the Muckenhoupt condition, testing over tails of indicators of
intervals as well. More generally, it is known that the energy condition is
necessary for boundedness of \emph{gradient elliptic} fractional singular
integrals on the real line \cite{SaShUr11}, but fails to be necessary for
certain elliptic singular integrals on the line.

\textbf{Failure of sufficiency of Muckenhoupt and Energy conditions}:
However, the weight pair $\left( \omega ,\ddot{\sigma}\right) $ constructed
in \cite{LaSaUr2} satisfies the Muckenhoupt and energy conditions, yet fails
to satisfy the norm inequality for the Hilbert transform\footnote{%
The interested reader can easily verify this, or see previous versions of
the current paper where details are included.}. This shows that, even
assuming the necessary conditions of Muckenhoupt and energy, we still need
some sort of testing conditions, and our $Tb$ theorem essentially leaves the
choice of testing conditions at our disposal - subject only to nondegeneracy
and size conditions.

For example, in the case of the Hilbert transform, Theorem \ref{dim one}
below roughly says this. As we are dealing with the case of general locally
finite positive Borel measures, all intervals appearing in this paper should
be assumed to closed on the left and open on the right, except when
otherwise noted.

\begin{theorem}[$Tb$ for Hilbert transform]
Let $H$ denote the Hilbert transform on the real line $\mathbb{R}$, let $%
\sigma $ and $\omega $ be locally finite positive Borel measures on $\mathbb{%
R}$. Then $H_{\sigma }$, where $H_{\sigma }f\equiv H\left( f\sigma \right) $%
, is bounded from \thinspace $L^{2}\left( \sigma \right) $ to $L^{2}\left(
\omega \right) $ \emph{if and only if} the Muckenhoupt and energy side
conditions hold, as well as the $\mathbf{b}$-testing and $\mathbf{b}^{\ast }$%
-testing conditions%
\begin{equation*}
\int_{I}\left\vert T_{\sigma }^{\alpha }b_{I}\right\vert ^{2}d\omega \leq
\left( \mathfrak{T}_{T^{\alpha }}^{\mathbf{b}}\right) ^{2}\left\vert
I\right\vert _{\sigma }\text{ and }\int_{J}\left\vert T_{\omega }^{\alpha
,\ast }b_{J}^{\ast }\right\vert ^{2}d\sigma \leq \left( \mathfrak{T}%
_{T^{\alpha }}^{\mathbf{b}^{\ast },\ast }\right) ^{2}\left\vert J\right\vert
_{\omega }\ ,
\end{equation*}%
taken over two families of test functions $\left\{ b_{I}\right\} _{I\in 
\mathcal{P}}$ and $\left\{ b_{J}^{\ast }\right\} _{J\in \mathcal{P}}$, where 
$b_{I}$ and $b_{J}^{\ast }$ are only required to be nondegenerate in an
average sense, and to be just slightly better than $L^{2}$ functions
themselves, namely $L^{p}$ for some $p>2$.
\end{theorem}

The families of test functions $\left\{ b_{I}\right\} _{I\in \mathcal{P}}$
and $\left\{ b_{J}^{\ast }\right\} _{J\in \mathcal{P}}$ in the $Tb$ theorem
above are nondegenerate and slightly better than $L^{2}$ functions, but
otherwise remain at the disposal of the reader. It is this flexibility in
choosing families of test functions that distinguishes this characterization
as compared to the corresponding $T1$ theorem\footnote{%
The energy conditions in (\ref{strong b* energy}) and (\ref{strong b energy}%
) below are relatively simple, stable and checkable conditions on a weight
pair, that are in addition, an almost immediate consequence of the
Muckenhoupt side condition and the testing conditions for $H$ over
indicators of intervals \cite{LaSaUr2}.}. The $Tb$ theorem here generalizes
many of the one-weight $Tb$ theorems in one dimension, since in the upper
doubling case, the Muckenhoupt $\mathfrak{A}_{2}$ condition and the energy
condition easily follow from the upper doubling condition. Recall that in
the one-weight case with doubling and upper doubling measures $\mu $, there
has been a long and sustained effort to relax the integrability conditions
of the testing functions: see e.g. S. Hofmann \cite{Hof} and Alfonseca,
Auscher, Axelsson, Hofmann and Kim \cite{AAAHK}. Subsequently, Hyt\"{o}nen-
Martikainen \cite{HyMa} assumed $Tb$ in $L^{s}\left( \mu \right) $ for some $%
s>2$, and the one weight theorem with testing functions $b$ in $L^{2}\left(
\mu \right) $ was attained by Lacey-Martikainen \cite{LaMa}, but their
argument strongly uses methods not immediately available in the two weight
setting.

Finally, we point out that the proof of our $Tb$ theorem is mostly
self-contained, but at the expense of considerable length. This is not just
for the convenience of the reader, but mainly because we must repeat much of
the proof strategy from \cite{NTV4}, \cite{LaSaShUr3}, \cite{Lac}, \cite%
{SaShUr7}, \cite{SaShUr9} and \cite{SaShUr10}, as the new ideas used here
force redevelopment of many of the previous arguments in these papers. We
now turn to a brief discussion of these new ideas for those readers already
acquainted with the theory of $T1$ and local $Tb$ theorems. See the brief
schematic (\ref{schematic}) below for a picture summary of the
decompositions involved.

\subsection{New ideas}

For those readers already familiar with the theory of local $Tb$ theorems,
we describe here some of the new techniques introduced in this paper to
handle the two weight situation. There are many difficulties to be overcome
in proving a local $Tb$ theorem, even in the one weight setting, as compared
to the corresponding $T1$ theorem, and we indicate four of them now.

\begin{enumerate}
\item \textbf{First difficulty:} In order to control the dual martingale
differences for `breaking' children, i.e. when the testing function
corresponding to a child is \textbf{not} the restriction of the testing
function of the parent, we need to construct coronas in which the
restrictions don't change, and for which the `breaking' intervals satisfy a
Carleson condition. This makes the so-called `nearby' inner products $%
\left\langle T_{\sigma }^{\alpha }b_{I},b_{J}^{\ast }\right\rangle _{\omega
} $, i.e. those in which the intervals $I$ and $J$ are close in both
position and scale, extremely difficult to estimate due to the fact that the
testing conditions are lost in the corona, except at the tops of coronas,
and are replaced with just a \textbf{weak} testing condition. Ironically,
these nearby inner products are the easiest to estimate in the proof of a $%
T1 $ theorem since the testing conditions remain in force in the coronas
there. In the one weight setting in \cite{NTV3}, \cite{HyMa} and \cite{LaMa}%
, special considerations, such as boundedness of Poisson integrals, are
taken into account in handling nearby inner products with random surgery,
and are unavailable to us here.\newline
\textbf{Resolution:} We develop a new recursive method for controlling the
nearby form by the \emph{energy conditions} and testing at the tops of the
coronas, in which we resurrect the original testing functions discarded
during the corona construction. This is presented in Section \ref{Sec nearby}%
.

\item \textbf{Second difficulty:} Both dual martingale and martingale
differences fail to satisfy two-sided frame-like and Riesz-like inequalities
in the setting of a $Tb$ theorem when $p=2$, complicating the treatment of
paraproducts.\newline
\textbf{Resolution:} We assume $p>2$ in the upper $L^{p}$ control of testing
functions, and then reduce this case to that of \emph{bounded} testing
functions using an absorption and recursion argument. For families of
bounded testing functions, we prove two-sided \emph{Riesz-like} inequalities
for dual martingale differences that are more robust than frame inequalities
(but only one-sided \emph{Riesz-like} inequalities for martingale
differences), and that enable many of the $T1$ two weight techniques to
carry over here in the $Tb$ setting. In particular these are key to
controlling paraproducts here.

\item \textbf{Third difficulty:} Only a weaker form of goodness due to Hyt%
\"{o}nen and Martikainen \cite{HyMa} is available for use in two weight $Tb$
theorems, complicating Lacey's treatment of the stopping form.\newline
\textbf{Resolution:} We adapt the two weight $T1$ arguments to accommodate 
\emph{weak goodness} in two ways, the first highly nontrivial and second
more straightforward:\newline
(\textbf{1}) in bounding the stopping form by Lacey's size functional on
admissible collections and bottom/up corona construction in \cite{Lac}, as
adapted in \cite{SaShUr9} and \cite{SaShUr10}, but using an additional
top/down `indented' corona construction, along with an enlargement of the
skeleton $\limfunc{skel}\left( I\right) $ of an interval $I$ to the body $%
\limfunc{body}\left( I\right) $ of $I$, in order to deal with the lack of
goodness in telescoping intervals - see Section \ref{Sec stop},\newline
(\textbf{2}) in controlling functional energy as in \cite{SaShUr9} and \cite%
{SaShUr10}, but with a different decomposition of the stopping intervals
into `Whitney' intervals, and modified pseudoprojections to accommodate two
independent families of grids - see Appendix B.

\item \textbf{Fourth difficulty:} The dual martingale differences are not in
general projections when some of the children `break', and the Monotonicity
Lemma fails to hold in any of the traditional forms in the setting of $T1$
theorems.\newline
\textbf{Resolution:} We introduce an additional square function bound on the
right hand side involving an infimum of averages, $\inf_{z\in \mathbb{R}%
}\sum_{J^{\prime }\in \mathfrak{C}_{\limfunc{broken}}\left( J\right)
}\left\vert J^{\prime }\right\vert _{\omega }\left( E_{J^{\prime }}^{\omega
}\left\vert x-z\right\vert \right) ^{2}$, summed over broken children. We
also use the fact that the corresponding `unbroken' dual martingale
differences form projections, but then we also need to modify the testing
function at the top of a corona, and also refine the triple corona
construction, so that dual martingale differences have controlled averages
on children throughout the corona.
\end{enumerate}

\subsection{Standard fractional singular integrals}

Let $0\leq \alpha <1$. We define a standard $\alpha $-fractional CZ kernel $%
K^{\alpha }(x,y)$ to be a real-valued function defined on $\mathbb{R}\times 
\mathbb{R}$ satisfying the following fractional size and smoothness
conditions of order $1+\delta $ for some $\delta >0$: For $x\neq y$,%
\begin{eqnarray}
\left\vert K^{\alpha }\left( x,y\right) \right\vert &\leq &C_{CZ}\left\vert
x-y\right\vert ^{\alpha -1}\text{ and }\left\vert \nabla K^{\alpha }\left(
x,y\right) \right\vert \leq C_{CZ}\left\vert x-y\right\vert ^{\alpha -2},
\label{sizeandsmoothness'} \\
\left\vert \nabla K^{\alpha }\left( x,y\right) -\nabla K^{\alpha }\left(
x^{\prime },y\right) \right\vert &\leq &C_{CZ}\left( \frac{\left\vert
x-x^{\prime }\right\vert }{\left\vert x-y\right\vert }\right) ^{\delta
}\left\vert x-y\right\vert ^{\alpha -2},\ \ \ \ \ \frac{\left\vert
x-x^{\prime }\right\vert }{\left\vert x-y\right\vert }\leq \frac{1}{2}, 
\notag
\end{eqnarray}%
and the last inequality also holds for the adjoint kernel in which $x$ and $%
y $ are interchanged. We note that a more general definition of kernel has
only order of smoothness $\delta >0$, rather than $1+\delta $, but the use
of the Monotonicity and Energy Lemmas in arguments below involves first
order Taylor approximations to the kernel functions $K^{\alpha }\left( \cdot
,y\right) $.

\subsubsection{Defining the norm inequality}

We now turn to a precise definition of the weighted norm inequality%
\begin{equation}
\left\Vert T_{\sigma }^{\alpha }f\right\Vert _{L^{2}\left( \omega \right)
}\leq \mathfrak{N}_{T_{\sigma }^{\alpha }}\left\Vert f\right\Vert
_{L^{2}\left( \sigma \right) },\ \ \ \ \ f\in L^{2}\left( \sigma \right) .
\label{two weight'}
\end{equation}%
For this we follow the lead of \cite{NTV3} and introduce a family $\left\{
\eta _{\delta ,R}^{\alpha }\right\} _{0<\delta <R<\infty }$ of nonnegative
functions on $\left[ 0,\infty \right) $ so that the truncated kernels $%
K_{\delta ,R}^{\alpha }\left( x,y\right) =\eta _{\delta ,R}^{\alpha }\left(
\left\vert x-y\right\vert \right) K^{\alpha }\left( x,y\right) $ are bounded
with compact support for fixed $x$ or $y$. Then the truncated operators 
\begin{equation}
T_{\sigma ,\delta ,R}^{\alpha }f\left( x\right) \equiv \int_{\mathbb{R}%
}K_{\delta ,R}^{\alpha }\left( x,y\right) f\left( y\right) d\sigma \left(
y\right) ,\ \ \ \ \ x\in \mathbb{R},  \label{def truncation}
\end{equation}%
are pointwise well-defined, and we will refer to the pair $\left( K^{\alpha
},\left\{ \eta _{\delta ,R}^{\alpha }\right\} _{0<\delta <R<\infty }\right) $
as an $\alpha $-fractional singular integral operator, which we typically
denote by $T^{\alpha }$, suppressing the dependence on the truncations.

\begin{definition}
\label{truncated op}We say that an $\alpha $-fractional singular integral
operator $T^{\alpha }=\left( K^{\alpha },\left\{ \eta _{\delta ,R}^{\alpha
}\right\} _{0<\delta <R<\infty }\right) $ satisfies the norm inequality (\ref%
{two weight'}) provided%
\begin{equation*}
\left\Vert T_{\sigma ,\delta ,R}^{\alpha }f\right\Vert _{L^{2}\left( \omega
\right) }\leq \mathfrak{N}_{T_{\sigma }^{\alpha }}\left\Vert f\right\Vert
_{L^{2}\left( \sigma \right) },\ \ \ \ \ f\in L^{2}\left( \sigma \right)
,0<\delta <R<\infty .
\end{equation*}
\end{definition}

It turns out that, in the presence of the Muckenhoupt conditions (\ref{def
A2}) below, the norm inequality (\ref{two weight'}) is essentially
independent of the choice of truncations used, and this is explained in some
detail in \cite{NTV3}, \cite{LaSaShUr3} and \cite{SaShUr10}. Thus, as in 
\cite{SaShUr10}, we are free to use the tangent line truncations described
there throughout the proofs of our results.

\subsection{Weakly accretive $\mathbf{b}$-testing and $\mathbf{b}^{\ast }$%
-testing conditions}

Denote by $\mathcal{P}$ the collection of intervals in $\mathbb{R}$. Note
that we include an $L^{p}$ upper bound in our definition of `$p$-weakly
accretive family' of functions.

\begin{definition}
Let $p\geq 2$ and let $\mu $ be a locally finite positive Borel measure on $%
\mathbb{R}$. We say that a family $\mathbf{b}=\left\{ b_{Q}\right\} _{Q\in 
\mathcal{P}}$ of functions indexed by $\mathcal{P}$ is a $p$\emph{-weakly }$%
\mu $\emph{-accretive} family of functions on $\mathbb{R}$ if%
\begin{eqnarray}
&&\limfunc{support}b_{Q}\subset Q\ ,\ \ \ \ \ Q\in \mathcal{P},
\label{local accretive} \\
1 &\leq &\left\vert \frac{1}{\left\vert Q\right\vert _{\mu }}%
\int_{Q}b_{Q}d\mu \right\vert \leq \left( \frac{1}{\left\vert Q\right\vert
_{\mu }}\int_{Q}\left\vert b_{Q}\right\vert ^{p}d\mu \right) ^{\frac{1}{p}%
}\leq C_{\mathbf{b}}<\infty ,\ \ \ \ \ Q\in \mathcal{P}\ .  \notag
\end{eqnarray}
\end{definition}

Suppose $\sigma $ and $\omega $ are locally finite positive Borel measures
on $\mathbb{R}$. The $\mathbf{b}$-testing conditions for $T^{\alpha }$ and $%
\mathbf{b}^{\ast }$-testing conditions for the dual $T^{\alpha ,\ast }$ are
given by%
\begin{eqnarray}
\int_{Q}\left\vert T_{\sigma }^{\alpha }b_{Q}\right\vert ^{2}d\omega &\leq
&\left( \mathfrak{T}_{T^{\alpha }}^{\mathbf{b}}\right) ^{2}\left\vert
Q\right\vert _{\sigma }\ ,\ \ \ \ \ \text{for all intervals }Q,
\label{b testing cond} \\
\int_{Q}\left\vert T_{\omega }^{\alpha ,\ast }b_{Q}^{\ast }\right\vert
^{2}d\sigma &\leq &\left( \mathfrak{T}_{T^{\alpha ,\ast }}^{\mathbf{b}^{\ast
}}\right) ^{2}\left\vert Q\right\vert _{\omega }\ ,\ \ \ \ \ \text{for all
intervals }Q,  \notag
\end{eqnarray}%
where these inequalities are interpreted as holding uniformly over
truncations of $T_{\sigma }^{\alpha }$ and $T^{\alpha ,\ast }$.

\subsection{Poisson integrals and the Muckenhoupt condition $\mathfrak{A}%
_{2}^{\protect\alpha }$}

Let $\mu $ be a locally finite positive Borel measure on $\mathbb{R}$, and
suppose $Q$ is an interval in $\mathbb{R}$. Recall that $\left\vert
Q\right\vert =\ell \left( Q\right) $ for an interval $Q$. The two $\alpha $%
-fractional Poisson integrals of $\mu $ on an interval $Q$ are given by the
following expressions:%
\begin{eqnarray}
\mathrm{P}^{\alpha }\left( Q,\mu \right) &\equiv &\int_{\mathbb{R}}\frac{%
\left\vert Q\right\vert }{\left( \left\vert Q\right\vert +\left\vert
x-x_{Q}\right\vert \right) ^{2-\alpha }}d\mu \left( x\right) ,
\label{def Poisson} \\
\mathcal{P}^{\alpha }\left( Q,\mu \right) &\equiv &\int_{\mathbb{R}}\left( 
\frac{\left\vert Q\right\vert }{\left( \left\vert Q\right\vert +\left\vert
x-x_{Q}\right\vert \right) ^{2}}\right) ^{1-\alpha }d\mu \left( x\right) , 
\notag
\end{eqnarray}%
where $\left\vert x-x_{Q}\right\vert $ denotes distance between $x$ and $%
x_{Q}$ and $\left\vert Q\right\vert $ denotes the Lebesgue measure of the
interval $Q$. We refer to $\mathrm{P}^{\alpha }$ as the \emph{standard}
Poisson integral and to $\mathcal{P}^{\alpha }$ as the \emph{reproducing}
Poisson integral. Note that these two kernels satisfy 
\begin{equation*}
0\leq \mathrm{P}^{\alpha }\left( Q,\mu \right) \leq C\mathcal{P}^{\alpha
}\left( Q,\mu \right) ,\ \ \ \ \text{for all intervals }Q\text{ and positive
measures }\mu .
\end{equation*}

We now define the \emph{one-tailed Muckenhoupt constant with holes }$%
\mathcal{A}_{2}^{\alpha }$ using the reproducing Poisson kernel $\mathcal{P}%
^{\alpha }$. On the other hand, the standard Poisson integral $\mathrm{P}%
^{\alpha }$ arises naturally throughout the proof of the $Tb$ theorem in
estimating oscillation of the fractional singular integral $T^{\alpha }$,
and in the definition of the energy conditions below.

\begin{definition}
Suppose $\sigma $ and $\omega $ are locally finite positive Borel measures
on $\mathbb{R}$. The one-tailed Muckenhoupt constants $\mathcal{A}%
_{2}^{\alpha }$ and $\mathcal{A}_{2}^{\alpha ,\ast }$ with holes for the
weight pair $\left( \sigma ,\omega \right) $ are given by%
\begin{eqnarray}
\mathcal{A}_{2}^{\alpha } &\equiv &\sup_{Q\in \mathcal{P}}\mathcal{P}%
^{\alpha }\left( Q,\mathbf{1}_{Q^{c}}\sigma \right) \frac{\left\vert
Q\right\vert _{\omega }}{\left\vert Q\right\vert ^{1-\alpha }}<\infty ,
\label{def call A2} \\
\mathcal{A}_{2}^{\alpha ,\ast } &\equiv &\sup_{Q\in \mathcal{P}}\mathcal{P}%
^{\alpha }\left( Q,\mathbf{1}_{Q^{c}}\omega \right) \frac{\left\vert
Q\right\vert _{\sigma }}{\left\vert Q\right\vert ^{1-\alpha }}<\infty . 
\notag
\end{eqnarray}
\end{definition}

Note that these definitions are the conditions with `holes' introduced by Hyt%
\"{o}nen \cite{Hyt} - the supports of the measures $\mathbf{1}_{Q^{c}}\sigma 
$ and $\mathbf{1}_{Q}\omega $ in the definition of $\mathcal{A}_{2}^{\alpha
} $ are disjoint, and so any common point masses of $\sigma $ and $\omega $
do not appear simultaneously in the factors of any of the products $\mathcal{%
P}^{\alpha }\left( Q,\mathbf{1}_{Q^{c}}\sigma \right) \frac{\left\vert
Q\right\vert _{\omega }}{\left\vert Q\right\vert ^{1-\alpha }}$. We will
also use the smaller `offset' Muckenhoupt condition%
\begin{equation*}
A_{2}^{\alpha }\equiv \sup_{\substack{ Q,Q^{\prime }\in \mathcal{P}  \\ Q%
\text{ and }Q^{\prime }\text{ are adjacent, }\ell \left( Q\right) =\ell
\left( Q^{\prime }\right) }}\frac{\left\vert Q\right\vert _{\omega }}{%
\left\vert Q\right\vert ^{1-\alpha }}\frac{\left\vert Q^{\prime }\right\vert
_{\sigma }}{\left\vert Q^{\prime }\right\vert ^{1-\alpha }}<\infty ,
\end{equation*}%
but the classical Muckenhoupt condition $A_{2}^{\alpha ,\limfunc{class}%
}\equiv \sup_{Q\in \mathcal{P}}\frac{\left\vert Q\right\vert _{\omega
}\left\vert Q\right\vert _{\sigma }}{\left\vert Q\right\vert ^{2-2\alpha }}%
<\infty $ will find no use in the two weight setting with common point
masses permitted.

\begin{remark}
Initially, these definitions of Muckenhoupt type were given in the following
`one weight' case, $d\omega \left( x\right) =w\left( x\right) dx$ and $%
d\sigma \left( x\right) =\frac{1}{w\left( x\right) }dx$, where $\mathcal{A}%
_{2}^{\alpha }\left( \lambda w,\left( \lambda w\right) ^{-1}\right) =%
\mathcal{A}_{2}^{\alpha }\left( w,w^{-1}\right) $ is homogeneous of degree $%
0 $. Of course the two weight version is homogeneous of degree $2$ in the
weight pair, $\mathcal{A}_{2}^{\alpha }\left( \lambda \sigma ,\lambda \omega
\right) =\lambda ^{2}\mathcal{A}_{2}^{\alpha }\left( \sigma ,\omega \right) $%
, while all of the other conditions we consider in connection with two
weight norm inequalities, including the operator norm $\mathfrak{N}%
_{T^{\alpha }}\left( \sigma ,\omega \right) $ itself, are homogeneous of
degree $1$ in the weight pair. This awkwardness regarding the homogeneity of
Muckenhoupt conditions could be rectified by simply taking the square root
of $\mathcal{A}_{2}^{\alpha }$ and renaming it, but the current definition
is so entrenched in the literature, in particular in connection with the $%
A_{2}$ conjecture, that we will leave it as is.
\end{remark}

\subsubsection{Punctured $A_{2}^{\protect\alpha }$ conditions}

The \emph{classical} $A_{2}^{\alpha }$ characteristic $\sup_{Q\in \mathcal{P}%
}\frac{\left\vert Q\right\vert _{\omega }}{\left\vert Q\right\vert
^{1-\alpha }}\frac{\left\vert Q\right\vert _{\sigma }}{\left\vert
Q\right\vert ^{1-\alpha }}$ fails to be finite when the measures $\sigma $
and $\omega $ have a common point mass - simply let $Q$ in the $\sup $ above
shrink to a common mass point. But there is a substitute that is quite
similar in character that is motivated by the fact that for large intervals $%
Q$, the $\sup $ above is problematic only if just \emph{one} of the measures
is \emph{mostly} a point mass when restricted to $Q$.

Given an at most countable set $\mathfrak{P}=\left\{ p_{k}\right\}
_{k=1}^{\infty }$ in $\mathbb{R}$, an interval $Q\in \mathcal{P}$, and a
positive locally finite Borel measure $\mu $, define 
\begin{equation}
\mu \left( Q,\mathfrak{P}\right) \equiv \left\vert Q\right\vert _{\mu }-\sup
\left\{ \mu \left( p_{k}\right) :p_{k}\in Q\cap \mathfrak{P}\right\} ,
\label{puncture}
\end{equation}%
where the supremum is actually achieved since $\sum_{p_{k}\in Q\cap 
\mathfrak{P}}\mu \left( p_{k}\right) <\infty $ as $\mu $ is locally finite.
The quantity $\mu \left( Q,\mathfrak{P}\right) $ is simply the $\widetilde{%
\mu }$ measure of $Q$ where $\widetilde{\mu }$ is the measure $\mu $ with
its largest point mass from $\mathfrak{P}$ in $Q$ removed. Given a locally
finite positive measure pair $\left( \sigma ,\omega \right) $, let 
\begin{equation}
\mathfrak{P}_{\left( \sigma ,\omega \right) }=\left\{ p_{k}\right\}
_{k=1}^{\infty }  \label{def common point mass}
\end{equation}%
be the at most countable set of common point masses of $\sigma $ and $\omega 
$. Then the weighted norm inequality (\ref{two weight'}) typically implies
finiteness of the following \emph{punctured} Muckenhoupt conditions:%
\begin{eqnarray}
A_{2}^{\alpha ,\limfunc{punct}}\left( \sigma ,\omega \right) &\equiv
&\sup_{Q\in \mathcal{P}}\frac{\omega \left( Q,\mathfrak{P}_{\left( \sigma
,\omega \right) }\right) }{\left\vert Q\right\vert ^{1-\alpha }}\frac{%
\left\vert Q\right\vert _{\sigma }}{\left\vert Q\right\vert ^{1-\alpha }},
\label{def punct} \\
A_{2}^{\alpha ,\ast ,\limfunc{punct}}\left( \sigma ,\omega \right) &\equiv
&\sup_{Q\in \mathcal{P}}\frac{\left\vert Q\right\vert _{\omega }}{\left\vert
Q\right\vert ^{1-\alpha }}\frac{\sigma \left( Q,\mathfrak{P}_{\left( \sigma
,\omega \right) }\right) }{\left\vert Q\right\vert ^{1-\alpha }}.  \notag
\end{eqnarray}

All of the above Muckenhuopt conditions $\mathcal{A}_{2}^{\alpha }$, $%
\mathcal{A}_{2}^{\alpha ,\ast }$, $A_{2}^{\alpha ,\limfunc{punct}}$ and $%
A_{2}^{\alpha ,\ast ,\limfunc{punct}}$ are necessary for boundedness of an
elliptic $\alpha $-fractional singular integral $T_{\sigma }^{\alpha }$ on
the line from\thinspace $L^{2}\left( \sigma \right) $ to $L^{2}\left( \omega
\right) $ (see \cite{SaShUr10}). It is convenient to define%
\begin{equation}
\mathfrak{A}_{2}^{\alpha }\equiv \mathcal{A}_{2}^{\alpha }+\mathcal{A}%
_{2}^{\alpha ,\ast }+A_{2}^{\alpha ,\limfunc{punct}}+A_{2}^{\alpha ,\ast ,%
\limfunc{punct}}\ .  \label{def A2}
\end{equation}

\subsection{Energy Conditions}

Here is the definition of the strong energy conditions, which we sometimes
refer to simply as the energy conditions. Let $m_{I}^{\mu }\equiv \frac{1}{%
\left\vert I\right\vert }\int xd\mu \left( x\right) $ be the average of $x$
over $I$ with respect to the measure $\mu $, which we often abbreviate to $%
m_{I}$ when the measure $\mu $ is understood.

\begin{definition}
\label{def strong quasienergy}Let $0\leq \alpha <1$. Suppose $\sigma $ and $%
\omega $ are locally finite positive Borel measures on $\mathbb{R}$. Then
the \emph{strong} energy constant $\mathcal{E}_{2}^{\alpha }$ is defined by 
\begin{equation}
\left( \mathcal{E}_{2}^{\alpha }\right) ^{2}\equiv \sup_{I=\dot{\cup}I_{r}}%
\frac{1}{\left\vert I\right\vert _{\sigma }}\sum_{r=1}^{\infty }\left( \frac{%
\mathrm{P}^{\alpha }\left( I_{r},\mathbf{1}_{I}\sigma \right) }{\left\vert
I_{r}\right\vert }\right) ^{2}\left\Vert x-m_{I_{r}}^{\omega }\right\Vert
_{L^{2}\left( \mathbf{1}_{I_{r}}\omega \right) }^{2}\ ,
\label{strong b* energy}
\end{equation}%
where the supremum is taken over arbitrary decompositions of an interval $I$
using a pairwise disjoint union of subintervals $I_{r}$. Similarly, we
define the dual \emph{strong} energy constant $\mathcal{E}_{\alpha }^{%
\limfunc{strong},\mathbf{b},\ast }$ by switching the roles of $\sigma $ and $%
\omega $:%
\begin{equation}
\left( \mathcal{E}_{2}^{\alpha ,\ast }\right) ^{2}\equiv \sup_{I=\dot{\cup}%
I_{r}}\frac{1}{\left\vert I\right\vert _{\omega }}\sum_{r=1}^{\infty }\left( 
\frac{\mathrm{P}^{\alpha }\left( I_{r},\mathbf{1}_{I}\omega \right) }{%
\left\vert I_{r}\right\vert }\right) ^{2}\left\Vert x-m_{I_{r}}^{\sigma
}\right\Vert _{L^{2}\left( \mathbf{1}_{I_{r}}\sigma \right) }^{2}\ .
\label{strong b energy}
\end{equation}
\end{definition}

These energy conditions are necessary for boundedness of elliptic and
gradient elliptic operators, including the Hilbert transform (but not for
certain elliptic singular operators that fail to be gradient elliptic) - see 
\cite{SaShUr11}, and see also (\ref{energy condition is necd}) below for
control of the energy constants $\mathcal{E}_{2}^{\alpha }$ and $\mathcal{E}%
_{2}^{\alpha ,\ast }$\ by the $\mathbf{1}$-testing and Muckenhoupt constants 
$\mathfrak{T}_{T^{\alpha }}^{\mathbf{1}}$, $\mathfrak{T}_{T^{\alpha }}^{%
\mathbf{1},\ast }$ and $\sqrt{\mathfrak{A}_{2}^{\alpha }}$. It is convenient
to define%
\begin{equation}
\mathfrak{E}_{2}^{\alpha }\equiv \mathcal{E}_{2}^{\alpha }+\mathcal{E}%
_{2}^{\alpha ,\ast },  \label{def frak energy}
\end{equation}%
as well as%
\begin{equation}
\mathcal{NTV}_{\alpha }\equiv \mathfrak{T}_{T^{\alpha }}^{\mathbf{b}}+%
\mathfrak{T}_{T^{\alpha }}^{\mathbf{b}^{\ast },\ast }+\sqrt{\mathfrak{A}%
_{2}^{\alpha }}+\mathfrak{E}_{2}^{\alpha }\ .  \label{def NTV}
\end{equation}

\section{The local $Tb$ theorem and proof preliminaries}

We derive a local $Tb$ theorem based in part on our \emph{proof} of the $T1$
theorem in \cite{SaShUr7}, \cite{SaShUr6}, \cite{SaShUr9} and \cite{SaShUr10}%
, in turn based on prior work in \cite{NTV4}, \cite{LaSaShUr3} and \cite{Lac}%
, and in part on the \emph{proof} of one weight $Tb$ theorems in Nazarov,
Treil and Volberg \cite{NTV3} and Hyt\"{o}nen and Martikainen \cite{HyMa}.
Recall from \cite{SaShUr11} that an $\alpha $-fractional singular integral $%
T^{\alpha }$ with kernel $K^{\alpha }$ is said to be \emph{elliptic} if $%
\left\vert K^{\alpha }\left( x,y\right) \right\vert \geq c\left\vert
x-y\right\vert ^{\alpha -1}$, and \emph{gradient elliptic} if the kernel $%
K^{\alpha }\left( x,y\right) $ satisfies%
\begin{equation}
\frac{d}{dx}K^{\alpha }\left( x,y\right) ,-\frac{d}{dy}K^{\alpha }\left(
x,y\right) \geq c\left\vert x-y\right\vert ^{\alpha -2}.
\label{def grad elliptic}
\end{equation}%
The Hilbert transform kernel $K\left( x,y\right) =\frac{1}{y-x}$ satisfies (%
\ref{def grad elliptic}) with $\alpha =0$. In dimension $n=1$ the
Muckenhoupt conditions are necessary for norm boundedness of elliptic
operators by results in \cite{LaSaUr2}, \cite{Hyt2} and \cite{SaShUr9}, and
the energy conditions are necessary for norm boundedness of gradient
elliptic operators by results in \cite{SaShUr11}. Moreover, in dimension $%
n=1 $, Hyt\"{o}nen \cite[Corollary 3.10]{Hyt2} proves that full testing is
controlled by testing and the Muckenhoupt conditions for the Hilbert
transform, and this is easily extended to $0\leq \alpha <1$ - see (\ref{full
proved}) below. Here is our two weight local $Tb$ theorem.

\begin{theorem}
\label{dim one}Suppose that $\sigma $ and $\omega $ are locally finite
positive Borel measures on the real line $\mathbb{R}$.

\begin{enumerate}
\item Assume that $T^{\alpha }$ is a standard $\alpha $-fractional elliptic
and gradient elliptic singular integral operator on $\mathbb{R}$, and set $%
T_{\sigma }^{\alpha }f=T^{\alpha }\left( f\sigma \right) $ for any smooth
truncation of $T_{\sigma }^{\alpha }$, so that $T_{\sigma }^{\alpha }$ is 
\emph{apriori} bounded from $L^{2}\left( \sigma \right) $ to $L^{2}\left(
\omega \right) $.

\item Let $p>2$ and let $\mathbf{b}=\left\{ b_{Q}\right\} _{Q\in \mathcal{P}%
} $ be a $p$-weakly $\sigma $-accretive family of functions on $\mathbb{R}$,
and let $\mathbf{b}^{\ast }=\left\{ b_{Q}^{\ast }\right\} _{Q\in \mathcal{P}%
} $ be a $p$-weakly $\omega $-accretive family of functions on $\mathbb{R}$.

\item Then for $0\leq \alpha <1$, the operator $T_{\sigma }^{\alpha }$ is
bounded from $L^{2}\left( \sigma \right) $ to $L^{2}\left( \omega \right) $
with operator norm $\mathfrak{N}_{T_{\sigma }^{\alpha }}$, i.e. 
\begin{equation*}
\left\Vert T_{\sigma }^{\alpha }f\right\Vert _{L^{2}\left( \omega \right)
}\leq \mathfrak{N}_{T_{\sigma }^{\alpha }}\left\Vert f\right\Vert
_{L^{2}\left( \sigma \right) },\ \ \ \ \ f\in L^{2}\left( \sigma \right) ,
\end{equation*}%
\textbf{uniformly} in smooth truncations of $T^{\alpha }$ \emph{if and only
if} the Muckenhoupt and energy conditions hold, i.e. $\mathcal{A}%
_{2}^{\alpha },\mathcal{A}_{2}^{\alpha ,\ast },A_{2}^{\alpha ,\limfunc{punct}%
},A_{2}^{\alpha ,\ast ,\limfunc{punct}},\mathcal{E}_{2}^{\alpha },\mathcal{E}%
_{2}^{\alpha ,\ast }<\infty $, and the $\mathbf{b}$-testing conditions for $%
T^{\alpha }$ and the $\mathbf{b}^{\ast }$-testing conditions for the dual $%
T^{\alpha ,\ast }$ both hold. Moreover, we have the equivalence,%
\begin{equation*}
\mathfrak{N}_{T^{\alpha }}\approx \mathcal{NTV}_{\alpha }=\mathfrak{T}_{%
\mathbf{R}^{\alpha }}^{\mathbf{b}}+\mathfrak{T}_{\mathbf{R}^{\alpha }}^{%
\mathbf{b}^{\ast }}+\sqrt{\mathfrak{A}_{2}^{\alpha }}+\mathfrak{E}%
_{2}^{\alpha }\ .
\end{equation*}
\end{enumerate}
\end{theorem}

\begin{remark}
\label{special}In the special case that $\sigma =\omega =\mu $, the
classical Muckenhoupt $A_{2}^{\alpha }$ condition is%
\begin{equation*}
\sup_{Q\in \mathcal{P}}\frac{\left\vert Q\right\vert _{\mu }}{\left\vert
Q\right\vert ^{1-\alpha }}\frac{\left\vert Q\right\vert _{\mu }}{\left\vert
Q\right\vert ^{1-\alpha }}<\infty ,
\end{equation*}%
which is the upper doubling measure\ condition with exponent $1-\alpha $,
i.e. 
\begin{equation*}
\left\vert Q\right\vert _{\mu }\leq C\ell \left( Q\right) ^{1-\alpha },\ \ \
\ \ \text{for all intervals }Q,\text{ }
\end{equation*}%
which of course prohibits point masses in $\mu $. Both Poisson integrals are
then bounded, 
\begin{eqnarray*}
\mathrm{P}^{\alpha }\left( Q,\mu \right) &\lesssim &\sum_{k=0}^{\infty }%
\frac{\left\vert Q\right\vert }{\left( 2^{k}\left\vert Q\right\vert \right)
^{2-\alpha }}\left\vert 2^{k}Q\right\vert _{\mu }\lesssim \sum_{k=0}^{\infty
}\frac{\left\vert Q\right\vert }{\left( 2^{k}\left\vert Q\right\vert \right)
^{2-\alpha }}\left( 2^{k}\left\vert Q\right\vert \right) ^{1-\alpha }=2, \\
\mathcal{P}^{\alpha }\left( Q,\mu \right) &\lesssim &\sum_{k=0}^{\infty
}\left( \frac{\left\vert Q\right\vert }{\left( 2^{k}\left\vert Q\right\vert
\right) ^{2}}\right) ^{1-\alpha }\left\vert 2^{k}Q\right\vert _{\mu
}\lesssim \sum_{k=0}^{\infty }\left( \frac{\left\vert Q\right\vert }{\left(
2^{k}\left\vert Q\right\vert \right) ^{2}}\right) ^{1-\alpha }\left(
2^{k}\left\vert Q\right\vert \right) ^{1-\alpha }=C_{\alpha }<\infty ,
\end{eqnarray*}%
and it follows easily that the equal weight pair $\left( \mu ,\mu \right) $
satisfies not only the Muckenhoupt $\mathfrak{A}_{2}^{\alpha }$ condition,
but also the strong energy condition $\mathfrak{E}_{2}^{\alpha }$:%
\begin{eqnarray*}
&&\sum_{r=1}^{\infty }\left( \frac{\mathrm{P}^{\alpha }\left( I_{r},\mathbf{1%
}_{I}\sigma \right) }{\left\vert I_{r}\right\vert }\right) ^{2}\left\Vert x-%
\mathbb{E}_{I_{r}}^{\omega }x\right\Vert _{L^{2}\left( \omega \right)
}^{2}\leq C\sum_{r=1}^{\infty }\left\Vert \frac{x-\mathbb{E}_{I_{r}}^{\omega
}x}{\left\vert I_{r}\right\vert }\right\Vert _{L^{2}\left( \omega \right)
}^{2} \\
&&\ \ \ \ \ \ \ \ \ \ \ \ \ \ \ \ \ \ \ \ \ \ \ \ \ \ \ \ \ \leq
C\sum_{r=1}^{\infty }\left\vert I_{r}\right\vert _{\omega }\leq C\left\vert
I\right\vert _{\omega }=C\left\vert I\right\vert _{\sigma }\ ,
\end{eqnarray*}%
since $\omega =\sigma $. Thus Theorem \ref{dim one}, when restricted to a
single weight $\sigma =\omega $, recovers a weaker version of the one weight
theorem of Lacey and Martikainen \cite[Theorem 1.1]{LaMa} for dimension $n=1$
- weaker due to our assumption that $p>2$. On the other hand, the
possibility of a two weight theorem for a $2$-weakly $\mu $-accretive family
is highly problematic, as one of the key proof strategies used in \cite{LaMa}
in the one weight case is a reduction to testing over $f$ and $g$ with
controlled $L^{\infty }$ norm via interpolation, a strategy that appears to
be unavailable in the two weight setting.
\end{remark}

\begin{problem}
Does Theorem \ref{dim one} remain true in the case $p=2$, i.e. when $\mathbf{%
b}=\left\{ b_{Q}\right\} _{Q\in \mathcal{P}}$ is a $2$-weakly $\sigma $%
-accretive family of functions, and $\mathbf{b}^{\ast }=\left\{ b_{Q}^{\ast
}\right\} _{Q\in \mathcal{P}}$ is a $2$-weakly $\omega $-accretive family of
functions?
\end{problem}

\begin{problem}
Are the energy conditions in Theorem \ref{dim one} already implied by the
Muckenhoupt, $\mathbf{b}$-testing and dual $\mathbf{b}^{\ast }$-testing
conditions for a pair of $p$-weakly accretive families when $p>2$?
\end{problem}

In order to prove Theorem \ref{dim one}, we first establish some improved
properties for a $p$-weakly $\mu $-accretive family, and establish some
improved energy conditions related to the families of testing functions $%
\mathbf{b}$ and $\mathbf{b}^{\ast }$. We turn to these matters in the next
three subsections.

\subsection{Reduction to the pointwise lower bound property}

Here we show that we may assume without loss of generality that the $p$%
-weakly accretive families of testing functions $b_{Q}$ and $b_{Q}^{\ast }$
for $Q\in \mathcal{P}$ satisfy the \emph{pointwise lower bound property},
written $PLBP$:%
\begin{equation}
\left\vert b_{Q}\left( x\right) \right\vert \geq c_{1}>0\ \ \ \ \ \text{for}%
\ Q\in \mathcal{P}\text{ and }\sigma \text{-a.e. }x\in \mathbb{R},
\label{plb}
\end{equation}%
for some positive constant $c_{1}$. Of course if $b_{Q}=\mathbf{1}_{Q}b$ for
some globally defined $b$, then the $PLBP$ is immediate from Lebesgue's
dyadic differentiation theorem. We make the following definition of a $p$%
\emph{-strongly} $\mu $-accretive family in $\mathbb{R}$.

\begin{definition}
\label{def strongly accretive}Let $\mu $ be a positive Borel measure on $%
\mathbb{R}$. We say that a family $\mathbf{b}=\left\{ b_{Q}\right\} _{Q\in 
\mathcal{P}}$ of functions indexed by $\mathcal{P}$ is a $p$\emph{-strongly }%
$\mu $\emph{-accretive} family of functions on $\mathbb{R}$ if the $b_{Q}$
are real-valued and there are positive constants $C_{\mathbf{b}}$, and $%
c_{1} $ such that 
\begin{eqnarray*}
&&\limfunc{support}b_{Q}\subset Q\ ,\ \ \ \ \ Q\in \mathcal{P}, \\
&&0<1\leq \left\vert \frac{1}{\left\vert Q\right\vert _{\mu }}%
\int_{Q}b_{Q}d\mu \right\vert \leq \left( \frac{1}{\left\vert Q\right\vert
_{\mu }}\int_{Q}\left\vert b_{Q}\right\vert ^{p}d\mu \right) ^{\frac{1}{p}%
}\leq C_{\mathbf{b}}<\infty ,\ \ \ \ \ Q\in \mathcal{P}\ , \\
&&\left\vert b_{Q}\left( x\right) \right\vert \geq c_{1}>0\ \ \ \ \ \text{%
for }\sigma \text{-a.e. }x\in \mathbb{R}.
\end{eqnarray*}
\end{definition}

To obtain that the families $\mathbf{b}=\left\{ b_{Q}\right\} _{Q\in 
\mathcal{P}^{n}}$ and $\mathbf{b}^{\ast }=\left\{ b_{Q}^{\ast }\right\}
_{Q\in \mathcal{P}^{n}}$ can be assumed to satisfy the $PLBP$ requires some
effort. But first, let us make a simple observation (essentially in \cite%
{HyMa}) under the additional assumption that the \emph{breaking intervals} $%
Q $, those for which there is a dyadic child $Q^{\prime }$ of $Q$ with $%
b_{Q^{\prime }}\neq \mathbf{1}_{Q^{\prime }}b_{Q}$, satisfy a $\mu $%
-Carleson condition. If $G_{k}\equiv \cup \mathcal{G}_{k}$ where%
\begin{equation*}
\mathcal{G}_{k}\equiv \left\{ A\in \mathcal{A}:A\text{ is a }k^{th}\text{
generation breaking interval}\right\} ,
\end{equation*}%
then $\left\vert \dbigcap\limits_{k=1}^{\infty }G_{k}\right\vert _{\sigma
}=0 $ since $\left\vert G_{k}\right\vert _{\sigma }\lesssim 2^{-\delta k}$
for some $\delta >0$ by the Carleson condition on breaking intervals. Thus
for $\sigma $-almost every $x$, the sequence of test functions $\left\{
b_{Q}\right\} _{Q:\ x\in Q}$, when arranged in order of decreasing length of 
$Q$, has the property that all sufficiently small $Q$ with $x\in Q$ belong
to the same corona $\mathcal{C}_{A}$ with $x\in A$, and hence $b_{Q}=\mathbf{%
1}_{Q}b_{A}$ for sufficiently small intervals $Q$ containing $x$. Suppose
that $A\in \mathcal{G}_{k}$. Then by Lebesgue's dyadic differentiation
theorem, we have 
\begin{equation*}
\left\vert b_{A}\left( x\right) \right\vert =\left\vert \lim_{Q\searrow x}%
\frac{1}{\left\vert Q\right\vert _{\mu }}\int_{Q}b_{A}d\mu \right\vert
=\lim_{Q\searrow x}\left\vert \frac{1}{\left\vert Q\right\vert _{\mu }}%
\int_{Q}b_{Q}d\mu \right\vert \geq c>0,
\end{equation*}%
for $\sigma $-a.e. $x\in A\setminus \left( \cup \mathcal{G}_{k+1}\right) $.
But this misses showing that $\left\vert b_{A}\left( x\right) \right\vert
\geq c>0$ on $A\cap \left( \cup \mathcal{G}_{k+1}\right) $, and for this we
must work harder.

\begin{proposition}
\label{lower bound}Let $p\geq 2$. Suppose that $\mathbf{b}=\left\{
b_{Q}\right\} _{Q\in \mathcal{P}}$ is a $p$-weakly $\sigma $-accretive
family of complex-valued functions on $\mathbb{R}$ that satisfies the $%
\mathbf{b}$-testing condition (\ref{b testing cond}) for a fractional
singular integral operator $T^{\alpha }$. Then there is a $p$-weakly $\sigma 
$-accretive family $\widetilde{\mathbf{b}}=\left\{ \widetilde{b}_{Q}\right\}
_{Q\in \mathcal{P}}$ that satisfies the $PLBP$. Moreover, the full $%
\widetilde{\mathbf{b}}$-testing condition (\ref{b testing cond}) for $%
T^{\alpha }$ holds and we have the estimate, 
\begin{eqnarray*}
&&\left\{ \widetilde{b}_{Q}\right\} _{Q\in \mathcal{P}^{n}}\text{ is }p\text{%
-strongly }\sigma \text{-accretive}, \\
\mathfrak{FT}_{T^{\alpha }}^{\widetilde{\mathbf{b}}}\left( \sigma ,\omega
\right) &\leq &C_{\alpha }\left( \mathfrak{T}_{T^{\alpha }}^{\mathbf{b}%
}\left( \sigma ,\omega \right) +\mathfrak{T}_{T^{\alpha }}^{\mathbf{b}^{\ast
},\ast }\left( \sigma ,\omega \right) +\sqrt{\mathcal{A}_{2}^{\alpha }\left(
\sigma ,\omega \right) }+\sqrt{\mathcal{A}_{2}^{\alpha ,\ast }\left( \sigma
,\omega \right) }\right) .
\end{eqnarray*}
\end{proposition}

\begin{proof}
For every interval $Q\in \mathcal{P}$ define%
\begin{equation*}
E\left( Q\right) \equiv \left\{ x\in Q:\left\vert b_{Q}\left( x\right)
\right\vert <\frac{1}{4}\right\} .
\end{equation*}%
Momentarily fix an interval $Q\in \mathcal{P}$ and $\delta >0$, and let $%
\left\{ I_{j}\left( Q\right) \right\} _{j=1}^{\infty }$ be a collection of
pairwise disjoint intervals such that%
\begin{eqnarray*}
&&E\left( Q\right) =\left\{ x\in Q:\left\vert b_{Q}\left( x\right)
\right\vert <\frac{1}{4}\right\} \subset \dbigcup\limits_{j=1}^{\infty
}I_{j}\left( Q\right) \ ; \\
&&\left\vert \dbigcup\limits_{j=1}^{\infty }I_{j}\left( Q\right) \setminus
E\left( Q\right) \right\vert _{\sigma }<\delta \left\vert Q\right\vert
_{\sigma }\ .
\end{eqnarray*}%
Note that%
\begin{eqnarray*}
\left\vert Q\right\vert _{\sigma } &\leq &\left\vert \int_{Q}b_{Q}d\sigma
\right\vert =\left\vert \int_{E\left( Q\right) }b_{Q}d\sigma \right\vert
+\left\vert \int_{Q\setminus E\left( Q\right) }b_{Q}d\sigma \right\vert \\
&\leq &\frac{1}{4}\left\vert E\left( Q\right) \right\vert _{\sigma }+\eta
\int_{Q\setminus E\left( Q\right) }\left\vert b_{Q}\right\vert ^{2}d\sigma +%
\frac{1}{\eta }\int_{Q\setminus E\left( Q\right) }d\sigma \\
&\leq &\frac{1}{4}\left\vert E\left( Q\right) \right\vert _{\sigma }+\eta
\int_{Q}\left\vert b_{Q}\right\vert ^{2}d\sigma +\frac{1}{\eta }\left\vert
Q\setminus E\left( Q\right) \right\vert _{\sigma } \\
&\leq &\frac{1}{4}\left\vert E\left( Q\right) \right\vert _{\sigma }+\eta
C\left\vert Q\right\vert _{\sigma }+\frac{1}{\eta }\left\vert Q\setminus
E\left( Q\right) \right\vert _{\sigma }\ .
\end{eqnarray*}%
Thus taking $\eta =\frac{1}{2C}$, and dividing through by $\left\vert
Q\right\vert _{\sigma }$, we get 
\begin{eqnarray*}
1 &\leq &\frac{1}{4}\frac{\left\vert E\left( Q\right) \right\vert _{\sigma }%
}{\left\vert Q\right\vert _{\sigma }}+\frac{1}{2}+2C\frac{\left\vert
Q\setminus E\left( Q\right) \right\vert _{\sigma }}{\left\vert Q\right\vert
_{\sigma }}\leq \frac{3}{4}+2C\frac{\left\vert Q\setminus E\left( Q\right)
\right\vert _{\sigma }}{\left\vert Q\right\vert _{\sigma }}; \\
&\Longrightarrow &\frac{1}{4}\leq 2C\frac{\left\vert Q\setminus E\left(
Q\right) \right\vert _{\sigma }}{\left\vert Q\right\vert _{\sigma }}%
\Longrightarrow \left\vert Q\setminus E\left( Q\right) \right\vert _{\sigma
}\geq \frac{1}{8C}\left\vert Q\right\vert _{\sigma } \\
&\Longrightarrow &\left\vert E\left( Q\right) \right\vert _{\sigma }\leq
\left( 1-\frac{1}{8C}\right) \left\vert Q\right\vert _{\sigma }=\beta
\left\vert Q\right\vert _{\sigma }\ ,
\end{eqnarray*}%
where $\beta =1-\frac{1}{8C}$. Now we note that since $\delta >0$ can be
taken arbitrarily small, we may without loss of generality take $\delta =0$%
\footnote{%
A rigorous limiting argument can be modelled after that given here.}.
Altogether then, we have shown that for every $Q\in \mathcal{P}$, there is a
pairwise disjoint collection of intervals $\left\{ I_{j}^{Q}\right\} _{j}$
such that%
\begin{eqnarray*}
E\left( Q\right) &\equiv &\left\{ x\in Q:\left\vert b_{Q}\left( x\right)
\right\vert <\frac{1}{4}\right\} =\overset{\cdot }{\dbigcup }_{j=1}^{\infty
}I_{j}\left( Q\right) \ , \\
\text{and }\left\vert E\left( Q\right) \right\vert _{\sigma }
&=&\sum_{j=1}^{\infty }\left\vert I_{j}\left( Q\right) \right\vert _{\sigma
}\leq \beta \left\vert Q\right\vert _{\sigma }\ ,\text{ where }0<\beta =1-%
\frac{1}{8C}<1.
\end{eqnarray*}

Now we begin the first step of the construction of a new family $\left\{ 
\widetilde{b}_{Q}\right\} _{Q\in \mathcal{D}}$ that satisfies both the
accretivity conditions and the testing conditions, as well as the pointwise
lower bound condition. We start by defining for $\mathbf{\epsilon }=\left\{
\epsilon _{j}\right\} _{j=1}^{\infty }$,%
\begin{equation*}
\widetilde{b}_{Q}^{\mathbf{\epsilon }}\left( x\right) \equiv b_{Q}\left(
x\right) +\sum_{j=1}^{\infty }\epsilon _{j}b_{I_{j}\left( Q\right) }\left(
x\right) \ ,
\end{equation*}%
where $\epsilon _{j}\in \left\{ -1,1\right\} $ for all $j\geq 1$.

We first note that we can assume that the collection of intervals $\left\{
I_{i}\left( Q\right) \right\} _{i}$ is subordinate to the collection of
children of $Q$, i.e. $I_{i}\left( Q\right) \subset Q^{\prime }$ for some $%
Q^{\prime }\in \mathfrak{C}\left( Q\right) $ depending on $i$. Then we have
for each $Q^{\prime }\in \mathfrak{C}\left( Q\right) $ that%
\begin{eqnarray*}
\int_{Q^{\prime }}\left\vert \widetilde{b}_{Q}^{\mathbf{\epsilon }%
}\right\vert ^{p}d\sigma &=&\int_{Q^{\prime }}\left\vert
b_{Q}+\sum_{j=1}^{\infty }\epsilon _{j}b_{I_{j}}\right\vert ^{p}d\sigma \\
&\leq &2^{p}\left\{ \int_{Q^{\prime }}\left\vert b_{Q}\right\vert
^{2p}d\sigma +\sum_{j:\ I_{j}\subset Q^{\prime }}\int \left\vert
b_{I_{j}}\right\vert ^{p}d\sigma \right\} \\
&\leq &C_{p}\left\{ \left\vert Q^{\prime }\right\vert _{\sigma }+\sum_{j:\
I_{j}\subset Q^{\prime }}\left\vert I_{j}\right\vert _{\sigma }\right\} \leq
C_{p}\left( 1+\beta \right) \left\vert Q^{\prime }\right\vert _{\sigma }\ .
\end{eqnarray*}

Now with $\mathbb{E}$ denoting expectation with respect to the uniform
probability measure on $\Omega _{\infty }\equiv \left\{ -1,1\right\} ^{%
\mathbb{N}}$, we have%
\begin{equation*}
\mathbb{E}\int \widetilde{b}_{Q}^{\mathbf{\epsilon }}d\sigma =\int
b_{Q}d\sigma \geq 1>0\ ,
\end{equation*}%
and%
\begin{eqnarray*}
&&\mathbb{E}\int \left\vert T_{\sigma }^{\alpha }\left( \mathbf{1}_{Q}%
\widetilde{b}_{Q}^{\mathbf{\epsilon }}\right) \right\vert ^{2}d\omega =%
\mathbb{E}\int \left\vert T_{\sigma }^{\alpha }\left(
b_{Q}+\sum_{j=1}^{\infty }\epsilon _{j}b_{I_{j}}\right) \right\vert
^{2}d\omega \\
&=&\mathbb{E}\int \left\{ \left\vert T_{\sigma }^{\alpha }b_{Q}\right\vert
^{2}+2\func{Re}\sum_{j=1}^{\infty }\epsilon _{j}\left( T_{\sigma }^{\alpha
}b_{Q}\right) \overline{\left( T_{\sigma }^{\alpha }b_{I_{j}}\right) }%
+\sum_{j,k=1}^{\infty }\epsilon _{j}\epsilon _{k}\left( T_{\sigma }^{\alpha
}b_{I_{j}}\right) \overline{\left( T_{\sigma }^{\alpha }b_{I_{k}}\right) }%
\right\} d\omega \\
&\leq &\int \left\vert T_{\sigma }^{\alpha }b_{Q}\right\vert ^{2}d\omega
+\sum_{j=1}^{\infty }\int \left\vert T_{\sigma }^{\alpha
}b_{I_{j}}\right\vert ^{2}d\omega \\
&\leq &\mathfrak{FT}_{T^{\alpha }}^{\mathbf{b}}\left\vert Q\right\vert
_{\sigma }+\sum_{j=1}^{\infty }\mathfrak{FT}_{T^{\alpha }}^{\mathbf{b}%
}\left\vert I_{i}\left( Q\right) \right\vert _{\sigma }\leq \mathfrak{F}%
_{T^{\alpha }}^{\mathbf{b}}\left[ 1+\beta \right] \left\vert Q\right\vert
_{\sigma }\ .
\end{eqnarray*}

So altogether, at this point in the first step of the construction, we have
for each pair $\left( Q,E\left( Q\right) \right) $ consisting of an interval 
$Q$ and a subset $E\left( Q\right) $ having measure at most $\beta
\left\vert Q\right\vert _{\sigma }$:%
\begin{eqnarray*}
\left\vert E\left( Q\right) \right\vert _{\sigma } &\leq &\beta \left\vert
Q\right\vert _{\sigma }\ , \\
\mathbb{E}\int \widetilde{b}_{Q}^{\mathbf{\epsilon }}d\sigma &=&\int
b_{Q}d\sigma \geq \left\vert Q\right\vert _{\sigma }>0\ , \\
\int_{Q^{\prime }}\left\vert \widetilde{b}_{Q}^{\mathbf{\epsilon }%
}\right\vert ^{p}d\sigma &\leq &C_{p}\left( 1+\beta \right) \left\vert
Q^{\prime }\right\vert _{\sigma }\ ,\ \ \ \ \ Q^{\prime }\in \mathfrak{C}%
\left( Q\right) \ , \\
\mathbb{E}\int \left\vert T_{\sigma }^{\alpha }\left( \mathbf{1}_{Q}%
\widetilde{b}_{Q}^{\mathbf{\epsilon }}\right) \right\vert ^{2}d\omega &\leq &%
\mathfrak{FT}_{T^{\alpha }}^{\mathbf{b}}\left[ 1+\beta \right] \left\vert
Q\right\vert _{\sigma }\ .
\end{eqnarray*}

Now we choose a positive constant $A$ large enough so that with
probabilities $\frac{1}{2}$ for $\int \widetilde{b}_{Q}^{\mathbf{\epsilon }%
}d\sigma $ and $\frac{3}{4}$ for $\int_{Q}\left\vert T_{\sigma }^{\alpha
}\left( \mathbf{1}_{Q}\widetilde{b}_{Q}^{\mathbf{\epsilon }}\right)
\right\vert ^{2}d\omega $ (note that $\left( 1-\frac{1}{2}\right) +\left( 1-%
\frac{3}{4}\right) =\frac{3}{4}<1$), there exists $\mathbf{\epsilon }\in
\Omega _{\infty }$ so that $\widetilde{b}_{Q}^{1}\equiv \widetilde{b}_{Q}^{%
\mathbf{\epsilon }}$ satisfies%
\begin{eqnarray}
\int \widetilde{b}_{Q}^{1}d\sigma &\geq &\left\vert Q\right\vert _{\sigma
}>0\ ,  \label{exists} \\
\int_{Q^{\prime }}\left\vert \widetilde{b}_{Q}^{1}\right\vert ^{p}d\sigma
&\leq &AC_{p}\left( 1+\beta \right) \left\vert Q^{\prime }\right\vert
_{\sigma }\ ,\ \ \ \ \ Q^{\prime }\in \mathfrak{C}\left( Q\right) \ ,  \notag
\\
\int \left\vert T_{\sigma }^{\alpha }\left( \mathbf{1}_{Q}\widetilde{b}%
_{Q}^{1}\right) \right\vert ^{2}d\omega &\leq &A\mathfrak{FT}_{T^{\alpha }}^{%
\mathbf{b}}\left[ 1+\beta \right] \beta \left\vert Q\right\vert _{\sigma }\ .
\notag
\end{eqnarray}

To see how big $A$ must be taken to achieve (\ref{exists}), we use
Chebyshev's inequality as follows. Take $N$ large, set%
\begin{equation*}
\widetilde{b}_{Q}^{\mathbf{\epsilon },N}\left( x\right) \equiv b_{Q}\left(
x\right) +\sum_{j=1}^{N}\epsilon _{j}b_{I_{j}\left( Q\right) }\left(
x\right) ,
\end{equation*}%
and equip $\Omega _{N}=\left\{ -1,1\right\} ^{N}$ with the uniform
probability measure that assigns mass $\frac{1}{2^{N}}$ to each $\mathbf{%
\epsilon }\in \Omega _{N}$. Then for each child $Q^{\prime }$ of $Q$ we have
as above that%
\begin{equation*}
\frac{1}{\#\Omega _{N}}\sum_{\mathbf{\epsilon }\in \Omega _{N}}\int
\left\vert T_{\sigma }^{\alpha }\left( \mathbf{1}_{Q}\widetilde{b}_{Q}^{%
\mathbf{\epsilon },N}\right) \right\vert ^{2}d\omega =\mathbb{E}\int
\left\vert T_{\sigma }^{\alpha }\left( \mathbf{1}_{Q}\widetilde{b}_{Q}^{%
\mathbf{\epsilon },N}\right) \right\vert ^{2}d\omega \leq C\left( 1+\beta
\right) \left\vert Q\right\vert _{\sigma }\ ,
\end{equation*}%
which by Chebyshev's inequality implies that%
\begin{eqnarray*}
&&AC\left( 1+\beta \right) \left\vert Q\right\vert _{\sigma }\ \#\left\{ 
\mathbf{\epsilon }\in \Omega _{N}:\int \left\vert T_{\sigma }^{\alpha
}\left( \mathbf{1}_{Q}\widetilde{b}_{Q}^{\mathbf{\epsilon },N}\right)
\right\vert ^{2}d\omega >AC\left( 1+\beta \right) \left\vert Q\right\vert
_{\sigma }\right\} \\
&&\ \ \ \ \ \ \ \ \ \ \ \ \ \ \ \leq \sum_{\mathbf{\epsilon }\in \Omega
_{N}}\int \left\vert T_{\sigma }^{\alpha }\left( \mathbf{1}_{Q}\widetilde{b}%
_{Q}^{\mathbf{\epsilon },N}\right) \right\vert ^{2}d\omega \leq \#\Omega
_{N}\ C\left( 1+\beta \right) \left\vert Q\right\vert _{\sigma }\ ,
\end{eqnarray*}%
which gives in turn that%
\begin{equation*}
\frac{1}{\#\Omega _{N}}\#\left\{ \mathbf{\epsilon }\in \Omega _{N}:\int
\left\vert T_{\sigma }^{\alpha }\left( \mathbf{1}_{Q}\widetilde{b}_{Q}^{%
\mathbf{\epsilon },N}\right) \right\vert ^{2}d\omega >AC\left( 1+\beta
\right) \left\vert Q\right\vert _{\sigma }\right\} \leq \frac{1}{A}.
\end{equation*}%
This of course just says that the probability that $\int \left\vert
T_{\sigma }^{\alpha }\left( \mathbf{1}_{Q}\widetilde{b}_{Q}^{\mathbf{%
\epsilon },N}\right) \right\vert ^{2}d\omega $ exceeds $A$ times $C\left(
1+\beta \right) \left\vert Q\right\vert _{\sigma }$ is less than $\frac{1}{A}
$. So in order to achieve that this latter probability is at most $\frac{1}{4%
}$, we take $A\geq 4$. Then we get that%
\begin{equation*}
\int \left\vert T_{\sigma }^{\alpha }\left( \mathbf{1}_{Q}\widetilde{b}_{Q}^{%
\mathbf{\epsilon },N}\right) \right\vert ^{2}d\omega \leq AC\left( 1+\beta
\right) \left\vert Q\right\vert _{\sigma }\ ,\ \ \ \ \ Q^{\prime }\in 
\mathfrak{C}\left( Q\right) \ ,
\end{equation*}%
holds for a set of $\mathbf{\epsilon }$ in $\Omega _{N}$ of probability at
least $1-\frac{1}{4}=\frac{3}{4}$.

So altogether, provided that we take $A=4$, we obtain that all inequalities, 
\begin{eqnarray*}
\int \widetilde{b}_{Q}^{\mathbf{\epsilon },N}d\sigma &\geq &\left\vert
Q\right\vert _{\sigma }>0\ , \\
\int_{Q^{\prime }}\left\vert \widetilde{b}_{Q}^{\mathbf{\epsilon }%
,N}\right\vert ^{p}d\sigma &\leq &C_{p}\left( 1+\beta \right) \left\vert
Q^{\prime }\right\vert _{\sigma }\ ,\ \ \ \ \ Q^{\prime }\in \mathfrak{C}%
\left( Q\right) \ , \\
\int \left\vert T_{\sigma }^{\alpha }\left( \mathbf{1}_{Q}\widetilde{b}_{Q}^{%
\mathbf{\epsilon },N}\right) \right\vert ^{2}d\omega &\leq &A\mathfrak{FT}%
_{T^{\alpha }}^{\mathbf{b}}\left[ 1+\beta \right] \beta \left\vert
Q\right\vert _{\sigma }\ ,
\end{eqnarray*}%
hold simultaneously for a set of $\mathbf{\epsilon }$ in $\Omega _{N}$ of
probability at least%
\begin{equation*}
1-\left\{ \left( 1-\frac{1}{2}\right) +\left( 1-\frac{3}{4}\right) \right\}
=1-\left\{ \frac{1}{2}+\frac{1}{4}\right\} =\frac{1}{4}.
\end{equation*}%
Since these estimates are independent of $N$, we can let $N\rightarrow
\infty $ to obtain that there is at least one $\mathbf{\epsilon }\in \Omega
_{\infty }$ for which (\ref{exists}) holds with $\widetilde{b}_{Q}^{1}\equiv 
\widetilde{b}_{Q}^{\mathbf{\epsilon },N}$.

We also have $\left\vert \widetilde{b}_{Q}^{1}\right\vert \geq \frac{3}{4}$
in $Q$ except on the exceptional set%
\begin{equation*}
G\equiv \overset{\cdot }{\dbigcup }_{j=1}^{\infty }E\left( I_{j}^{Q}\right)
\end{equation*}%
whose $\sigma $-measure satisfies 
\begin{equation*}
\left\vert G\right\vert _{\sigma }=\sum_{j=1}^{\infty }\left\vert E\left(
I_{j}^{Q}\right) \right\vert _{\sigma }\leq \sum_{j=1}^{\infty }\beta
\left\vert I_{j}^{Q}\right\vert _{\sigma }=\beta \left\vert E\left( Q\right)
\right\vert _{\sigma }\leq \beta ^{2}\left\vert Q\right\vert _{\sigma }\ .
\end{equation*}%
Now we consider the set%
\begin{equation*}
F\left( I_{j}\left( Q\right) \right) \equiv \left\{ x\in E\left( I_{j}\left(
Q\right) \right) :\left\vert b_{Q}\left( x\right) +\epsilon
_{j}b_{I_{j}^{Q}}\left( x\right) \right\vert <\frac{1}{4}\right\} .
\end{equation*}%
We may assume that%
\begin{equation*}
F\left( I_{j}\left( Q\right) \right) =\overset{\cdot }{\dbigcup }%
_{k=1}^{\infty }I_{k}\left( I_{j}\left( Q\right) \right)
\end{equation*}%
where $\left\{ I_{k}\left( I_{j}\left( Q\right) \right) \right\}
_{k=1}^{\infty }$ is a pairwise disjoint collection of intervals for each $j$%
. Now we apply the above step to each pair $\left( I_{j}\left( Q\right)
,F\left( I_{j}\left( Q\right) \right) \right) $ consisting of an interval $%
I_{j}^{Q}$ and a subset $F\left( I_{j}^{Q}\right) $ having measure at most $%
\beta \left\vert I_{j}^{Q}\right\vert _{\sigma }$. Then arguing as above for
each such pair, and adding results, we obtain that there exists $\mathbf{%
\epsilon }_{j}$ so that $\widetilde{b}_{Q}^{2}\equiv \sum_{j=1}^{\infty }%
\widetilde{b^{1}}_{Q}^{\mathbf{\epsilon }_{j}}$ satisfies%
\begin{eqnarray*}
\int \widetilde{b}_{Q}^{2}d\sigma &\geq &\left\vert Q\right\vert _{\sigma
}>0\ , \\
\int_{Q^{\prime }}\left\vert \widetilde{b}_{Q}^{2}\right\vert ^{p}d\sigma
&\leq &AC_{p}\left( 1+\beta +\beta ^{2}\right) \left\vert Q^{\prime
}\right\vert _{\sigma }\ ,\ \ \ \ \ Q^{\prime }\in \mathfrak{C}\left(
Q\right) \ , \\
\int \left\vert T_{\sigma }^{\alpha }\left( \mathbf{1}_{Q}\widetilde{b}%
_{Q}^{2}\right) \right\vert ^{2}d\omega &\leq &A\mathfrak{FT}_{T^{\alpha }}^{%
\mathbf{b}}\left[ 1+\beta +\beta ^{2}\right] \beta \left\vert Q\right\vert
_{\sigma }\ ,
\end{eqnarray*}%
as well as $\left\vert \widetilde{b}_{Q}^{2}\right\vert \geq \frac{3}{4}$ in 
$Q$ except on the exceptional set%
\begin{equation*}
G\equiv \overset{\cdot }{\dbigcup }_{j,k=1}^{\infty }E\left( I_{k}\left(
I_{j}\left( Q\right) \right) \right) ,
\end{equation*}%
whose $\sigma $-measure satisfies 
\begin{eqnarray*}
\left\vert G\right\vert _{\sigma } &=&\sum_{j,k=1}^{\infty }\left\vert
E\left( I_{k}\left( I_{j}\left( Q\right) \right) \right) \right\vert
_{\sigma }\leq \sum_{j=1}^{\infty }\beta \left\vert E\left( I_{j}\left(
Q\right) \right) \right\vert _{\sigma } \\
&\leq &\sum_{j=1}^{\infty }\beta ^{2}\left\vert I_{j}^{Q}\right\vert
_{\sigma }=\beta ^{2}\left\vert E\left( Q\right) \right\vert _{\sigma }\leq
\beta ^{3}\left\vert Q\right\vert _{\sigma }\ .
\end{eqnarray*}%
Continuing in this way, we end up with the desired function $\widetilde{b}%
_{Q}\left( x\right) =\lim_{n\rightarrow \infty }\widetilde{b}_{Q}^{n}\left(
x\right) $, since $\beta <1$ implies that the collection of nested intervals 
$\left\{ I_{j_{N}}\left( ...I_{j_{2}}\left( I_{j_{1}}\left( Q\right) \right)
...\right) \right\} $ satisfy a $\sigma $-Carleson condition. We emphasize
that for each interval $Q$, we then have $\left\vert \widetilde{b}_{Q}\left(
x\right) \right\vert \geq \frac{3}{4}$ for $\sigma $-a.e. $x\in Q$, as well
as $\int_{Q^{\prime }}\left\vert \widetilde{b}_{Q}\right\vert ^{2}d\sigma
\leq C^{\prime }\left\vert Q^{\prime }\right\vert _{\sigma }\ $for all $%
Q^{\prime }\in \mathfrak{C}\left( Q\right) $. The full $\widetilde{\mathbf{b}%
}$-testing condition constant that we have obtained satisfies%
\begin{equation*}
\mathfrak{FT}_{T^{\alpha }}^{\widetilde{\mathbf{b}}}\left( \sigma ,\omega
\right) \leq C_{\alpha }\left( \mathfrak{FT}_{T^{\alpha }}^{\mathbf{b}%
}\left( \sigma ,\omega \right) +\mathfrak{FT}_{T^{\alpha }}^{\mathbf{b}%
^{\ast },\ast }\left( \sigma ,\omega \right) +\sqrt{\mathcal{A}_{2}^{\alpha
}\left( \sigma ,\omega \right) }+\sqrt{\mathcal{A}_{2}^{\alpha ,\ast }\left(
\sigma ,\omega \right) }\right) ,
\end{equation*}%
and since the full testing constant is controlled by the testing constant
and the Muckenhoupt constant in (\ref{full proved})\ below, the proof is
complete.
\end{proof}

\subsection{Reduction to real bounded accretive families\label{str acc}}

Recall that a vector of `complex-valued testing functions' $\mathbf{b}\equiv
\left\{ b_{Q}\right\} _{Q\in \mathcal{D}}$ is a $p$-\emph{strongly }$\mu $-%
\emph{accretive} family if 
\begin{equation*}
\limfunc{support}b_{Q}\subset Q\ ,\ \ \ \ \ Q\in \mathcal{P},
\end{equation*}%
\begin{equation}
1\leq \left\vert \frac{1}{\left\vert Q\right\vert _{\mu }}\int_{Q}b_{Q}d\mu
\right\vert \leq \left( \frac{1}{\left\vert Q\right\vert _{\mu }}%
\int_{Q}\left\vert b_{Q}\right\vert ^{p}d\mu \right) ^{\frac{1}{p}}\leq C_{%
\mathbf{b}}\left( p\right) <\infty ,\ \ \ \ \ Q\in \mathcal{P}\ ,
\label{acc}
\end{equation}%
and if $\mathbf{b}$ satisfies the $PLBP$, i.e. 
\begin{equation*}
\left\vert b_{Q}\left( x\right) \right\vert \geq c_{1}>0\ \ \ \ \ \text{for}%
\ Q\in \mathcal{D}\text{ and }\sigma -a.e.x\in \mathbb{R}^{n}.
\end{equation*}%
We begin by noting that if $b_{Q}$ satisfies (\ref{acc}) with $\mu =\sigma $%
, and satisfies a given $\mathbf{b}$-testing condition for a weight pair $%
\left( \sigma ,\omega \right) $, then $\func{Re}b_{Q}$ satisfies $\left( 
\frac{1}{\left\vert Q\right\vert _{\mu }}\int_{Q}\left\vert \func{Re}%
b_{Q}\right\vert ^{p}d\mu \right) ^{\frac{1}{p}}\leq C_{\mathbf{b}}\left(
p\right) $ and the given $\mathbf{b}$-testing condition for $\left( \sigma
,\omega \right) $ with $\func{Re}b_{Q}$ in place of $b_{Q}$.

\begin{conclusion}
We may assume throughout the proof of Theorem \ref{dim one} that our $p$%
-weakly\emph{\ }$\mu $-accretive families $\mathbf{b}\equiv \left\{
b_{Q}\right\} _{Q\in \mathcal{D}}$ and $\mathbf{b}^{\ast }\equiv \left\{
b_{Q}^{\ast }\right\} _{Q\in \mathcal{G}}$ consist of \textbf{real-valued}
functions that in addition satisfy the $PLBP$ and $1\leq \frac{1}{\left\vert
Q\right\vert _{\sigma }}\int_{Q}b_{Q}d\sigma $ and $1\leq \frac{1}{%
\left\vert Q\right\vert _{\sigma }}\int_{Q}b_{Q}^{\ast }d\sigma $.
\end{conclusion}

Next we show that the assumption of testing conditions for a Calder\'{o}%
n-Zygmund operator $T$ and $p$-strongly $\mu $-accretive testing functions $%
\mathbf{b}=\left\{ b_{Q}\right\} _{Q\in \mathcal{P}}$ and $\mathbf{b}^{\ast
}=\left\{ b_{Q}^{\ast }\right\} _{Q\in \mathcal{P}}$ with $p>2$ can always
be replaced with real-valued $\infty $-strongly $\mu $-accretive testing
functions, thus reducing the $Tb$ theorem for the case $p>2$ to the case
when $p=\infty $ and the $PLBP$ (\ref{plb}) holds. We now proceed to develop
a precise statement. We extend (\ref{acc}) to $2\leq p\leq \infty $ by%
\begin{eqnarray}
&&\limfunc{support}b_{Q}\subset Q\ ,\ \ \ \ \ Q\in \mathcal{P},
\label{acc infinity} \\
1 &\leq &\left\vert \frac{1}{\left\vert Q\right\vert _{\mu }}%
\int_{Q}b_{Q}d\mu \right\vert \leq \left\{ 
\begin{array}{cc}
\left( \frac{1}{\left\vert Q\right\vert _{\mu }}\int_{Q}\left\vert
b_{Q}\right\vert ^{p}d\mu \right) ^{\frac{1}{p}}\leq C_{\mathbf{b}}\left(
p\right) <\infty & \text{for }2\leq p<\infty \\ 
\left\Vert b_{Q}\right\Vert _{L^{\infty }\left( \mu \right) }\leq C_{\mathbf{%
b}}\left( \infty \right) <\infty & \text{for }p=\infty%
\end{array}%
\right. ,\ \ \ \ \ Q\in \mathcal{P}\ .  \notag
\end{eqnarray}

\begin{proposition}
\label{conditional}Let $0\leq \alpha <1$, and let $\sigma $ and $\omega $ be
locally finite positive Borel measures on the real line $\mathbb{R}$, and
let $T^{\alpha }$ be a standard $\alpha $-fractional elliptic and gradient
elliptic singular integral operator on $\mathbb{R}$. Set $T_{\sigma
}^{\alpha }f=T^{\alpha }\left( f\sigma \right) $ for any smooth truncation
of $T_{\sigma }^{\alpha }$, so that $T_{\sigma }^{\alpha }$ is \emph{apriori}
bounded from $L^{2}\left( \sigma \right) $ to $L^{2}\left( \omega \right) $.
Finally, define the sequence of positive extended real numbers%
\begin{equation*}
\left\{ p_{n}\right\} _{n=0}^{\infty }=\left\{ \frac{2}{1-\left( \frac{2}{3}%
\right) ^{n}}\right\} _{n=0}^{\infty }=\left\{ \infty ,6,\frac{18}{5},\frac{%
162}{65},...\right\} .
\end{equation*}%
Suppose that the following statement is true:\medskip

\begin{description}
\item[$\left( \mathcal{S}_{\infty }\right) $] If $\mathbf{b}=\left\{
b_{Q}\right\} _{Q\in \mathcal{P}}$ is an $\infty $-strongly $\sigma $%
-accretive family of functions on $\mathbb{R}$, and if $\mathbf{b}^{\ast
}=\left\{ b_{Q}^{\ast }\right\} _{Q\in \mathcal{P}}$ is an $\infty $%
-strongly $\omega $-accretive family of functions on $\mathbb{R}$, then the
operator norm $\mathfrak{N}_{T_{\sigma }^{\alpha }}$ of $T_{\sigma }^{\alpha
}$ from $L^{2}\left( \sigma \right) $ to $L^{2}\left( \omega \right) $, i.e.
the best constant in%
\begin{equation*}
\left\Vert T_{\sigma }^{\alpha }f\right\Vert _{L^{2}\left( \omega \right)
}\leq \mathfrak{N}_{T_{\sigma }^{\alpha }}\left\Vert f\right\Vert
_{L^{2}\left( \sigma \right) },\ \ \ \ \ f\in L^{2}\left( \sigma \right) ,
\end{equation*}%
\textbf{uniformly} in smooth truncations of $T^{\alpha }$, satisfies%
\begin{equation*}
\mathfrak{T}_{T^{\alpha }}^{\mathbf{b}}+\mathfrak{T}_{T^{\alpha }}^{\mathbf{b%
}^{\ast }}+\sqrt{\mathfrak{A}_{2}^{\alpha }}+\mathfrak{E}_{2}^{\alpha
}\lesssim \mathfrak{N}_{T^{\alpha }}\lesssim \left( C_{\mathbf{b}}\left(
\infty \right) +C_{\mathbf{b}^{\ast }}\left( \infty \right) \right) \left( 
\mathfrak{T}_{T^{\alpha }}^{\mathbf{b}}+\mathfrak{T}_{T^{\alpha }}^{\mathbf{b%
}^{\ast }}+\sqrt{\mathfrak{A}_{2}^{\alpha }}+\mathfrak{E}_{2}^{\alpha
}\right) \ ,
\end{equation*}%
where $C_{\mathbf{b}}\left( \infty \right) ,C_{\mathbf{b}^{\ast }}\left(
\infty \right) $ are the accretivity constants in (\ref{acc infinity}), and
the constants implied by $\lesssim $ depend on $\alpha $ and the Calder\'{o}%
n-Zygmund constant $C_{CZ}$ in (\ref{sizeandsmoothness'}).\\[0.15cm]
Then for each $n\geq 0$, the following statements hold:\medskip

\item[$\left( \mathcal{S}_{n}\right) $] Let $p\in \left( p_{n+1},p_{n}\right]
$. If $\mathbf{b}=\left\{ b_{Q}\right\} _{Q\in \mathcal{P}}$ is a $p$%
-strongly $\sigma $-accretive family of functions on $\mathbb{R}$, and if $%
\mathbf{b}^{\ast }=\left\{ b_{Q}^{\ast }\right\} _{Q\in \mathcal{P}}$ is a $%
p $-strongly $\omega $-accretive family of functions on $\mathbb{R}$, then
the operator norm $\mathfrak{N}_{T_{\sigma }^{\alpha }}$ of $T_{\sigma
}^{\alpha }$ from $L^{2}\left( \sigma \right) $ to $L^{2}\left( \omega
\right) $, i.e. the best constant in%
\begin{equation*}
\left\Vert T_{\sigma }^{\alpha }f\right\Vert _{L^{2}\left( \omega \right)
}\leq \mathfrak{N}_{T_{\sigma }^{\alpha }}\left\Vert f\right\Vert
_{L^{2}\left( \sigma \right) },\ \ \ \ \ f\in L^{2}\left( \sigma \right) ,
\end{equation*}%
\textbf{uniformly} in smooth truncations of $T^{\alpha }$, satisfies%
\begin{equation*}
\mathfrak{T}_{T^{\alpha }}^{\mathbf{b}}+\mathfrak{T}_{T^{\alpha }}^{\mathbf{b%
}^{\ast }}+\sqrt{\mathfrak{A}_{2}^{\alpha }}+\mathfrak{E}_{2}^{\alpha
}\lesssim \mathfrak{N}_{T^{\alpha }}\lesssim \left( C_{\mathbf{b}}\left(
p\right) +C_{\mathbf{b}^{\ast }}\left( p\right) \right) ^{3^{n+1}}\left( 
\mathfrak{T}_{T^{\alpha }}^{\mathbf{b}}+\mathfrak{T}_{T^{\alpha }}^{\mathbf{b%
}^{\ast }}+\sqrt{\mathfrak{A}_{2}^{\alpha }}+\mathfrak{E}_{2}^{\alpha
}\right) \ ,
\end{equation*}%
where $C_{\mathbf{b}}\left( p\right) ,C_{\mathbf{b}^{\ast }}\left( p\right) $
are the accretivity constants in (\ref{acc}), and the constants implied by $%
\lesssim $ depend on $p$, $\alpha $, and the Calder\'{o}n-Zygmund constant $%
C_{CZ}$ in (\ref{sizeandsmoothness'}).
\end{description}
\end{proposition}

\begin{proof}[Proof of Proposition \protect\ref{conditional}]
We first prove $\left( \mathcal{S}_{0}\right) $. So fix $p\in \left(
p_{1},p_{0}\right) =\left( 6,\infty \right) $, and let $\mathbf{b}=\left\{
b_{Q}\right\} _{Q\in \mathcal{P}}$ be a $p$-weakly $\sigma $-accretive
family of functions on $\mathbb{R}$, and let $\mathbf{b}^{\ast }=\left\{
b_{Q}^{\ast }\right\} _{Q\in \mathcal{P}}$ be a $p$-weakly $\omega $%
-accretive family of functions on $\mathbb{R}$. Let $0<\varepsilon <1$ (to
be chosen differently at various points in the argument below) and define%
\begin{equation}
\lambda =\lambda \left( \varepsilon \right) =\left( \frac{p}{p-2}C_{\mathbf{b%
}}\left( p\right) ^{p}\frac{1}{\varepsilon }\right) ^{\frac{1}{p-2}},
\label{lambda choice}
\end{equation}%
and a new collection of test functions,%
\begin{equation}
\widehat{b}_{Q}\equiv 2b_{Q}\left( \mathbf{1}_{\left\{ \left\vert
b_{Q}\right\vert \leq \lambda \right\} }+\frac{\lambda }{\left\vert
b_{Q}\right\vert }\mathbf{1}_{\left\{ \left\vert b_{Q}\right\vert >\lambda
\right\} }\right) ,\ \ \ \ \ Q\in \mathcal{P},  \label{new}
\end{equation}%
which continue to satisfy the $PLBP$ (\ref{plb}). We compute%
\begin{eqnarray*}
\int_{\left\{ \left\vert b_{Q}\right\vert >\lambda \right\} }\left\vert
b_{Q}\right\vert ^{2}d\sigma &=&\int_{\left\{ \left\vert b_{Q}\right\vert
>\lambda \right\} }\left[ \int_{0}^{\left\vert b_{Q}\right\vert }2tdt\right]
d\sigma \\
&=&\int \int_{\left\{ \left( x,t\right) \in \mathbb{R}\times \left( 0,\infty
\right) :\max \left\{ t,\lambda \right\} <\left\vert b_{Q}\left( x\right)
\right\vert \right\} }2tdtd\sigma \left( x\right) \\
&=&\int_{\left[ 0,\lambda \right] }\int_{\left\{ x\in \mathbb{R}:\lambda
<\left\vert b_{Q}\left( x\right) \right\vert \right\} }d\sigma \left(
x\right) 2tdt+\int_{\left( \lambda ,\infty \right) }\int_{\left\{ x\in 
\mathbb{R}:t<\left\vert b_{Q}\left( x\right) \right\vert \right\} }d\sigma
\left( x\right) 2tdt \\
&=&\lambda ^{2}\left\vert \left\{ \left\vert b_{Q}\right\vert >\lambda
\right\} \right\vert _{\sigma }+\int_{\lambda }^{\infty }\left\vert \left\{
\left\vert b_{Q}\right\vert >t\right\} \right\vert _{\sigma }2tdt,
\end{eqnarray*}%
and hence%
\begin{eqnarray}
\int_{\left\{ \left\vert b_{Q}\right\vert >\lambda \right\} }\left\vert
b_{Q}\right\vert ^{2}d\sigma &\leq &\lambda ^{2}\frac{1}{\lambda ^{p}}\left(
\int \left\vert b_{Q}\right\vert ^{p}d\sigma \right) +\int_{\lambda
}^{\infty }\frac{1}{t^{p}}\left( \int \left\vert b_{Q}\right\vert
^{p}d\sigma \right) 2tdt  \label{hence} \\
&=&\left\{ \lambda ^{2-p}+\int_{\lambda }^{\infty }2t^{1-p}dt\right\} C_{%
\mathbf{b}}\left( p\right) ^{p}\left\vert Q\right\vert _{\sigma }  \notag \\
&=&\frac{p}{p-2}\lambda ^{2-p}C_{\mathbf{b}}\left( p\right) ^{p}\left\vert
Q\right\vert _{\sigma }=\varepsilon \left\vert Q\right\vert _{\sigma }\ , 
\notag
\end{eqnarray}%
by (\ref{lambda choice}). Thus we have the lower bound,%
\begin{eqnarray}
\left\vert \frac{1}{\left\vert Q\right\vert _{\sigma }}\int_{Q}\widehat{b}%
_{Q}d\sigma \right\vert &=&2\left\vert \frac{1}{\left\vert Q\right\vert
_{\sigma }}\int_{Q}b_{Q}d\sigma -\frac{1}{\left\vert Q\right\vert _{\sigma }}%
\int_{Q}b_{Q}\left( \frac{\lambda }{\left\vert b_{Q}\right\vert }-1\right) 
\mathbf{1}_{\left\{ \left\vert b_{Q}\right\vert >\lambda \right\} }d\sigma
\right\vert  \label{low} \\
&\geq &2\left\vert \frac{1}{\left\vert Q\right\vert _{\sigma }}%
\int_{Q}b_{Q}d\mu \right\vert -2\left( \frac{1}{\left\vert Q\right\vert
_{\sigma }}\int_{Q}\left\vert b_{Q}\right\vert ^{2}\mathbf{1}_{\left\{
\left\vert b_{Q}\right\vert >\lambda \right\} }d\sigma \right) ^{\frac{1}{2}}
\notag \\
&\geq &2-2\left( \frac{1}{\left\vert Q\right\vert _{\sigma }}\varepsilon
\left\vert Q\right\vert _{\sigma }\right) ^{\frac{1}{2}}=2-2\sqrt{%
\varepsilon }\geq 1>0,\ \ \ \ \ Q\in \mathcal{P}\ ,  \notag
\end{eqnarray}%
if we choose $0<\varepsilon \leq \frac{1}{4}$. For an upper bound we have 
\begin{equation*}
\left\Vert \widehat{b}_{Q}\right\Vert _{L^{\infty }\left( \sigma \right)
}\leq 2\lambda =2\lambda \left( \varepsilon \right) =2\left( \frac{p}{p-2}C_{%
\mathbf{b}}\left( p\right) ^{p}\frac{1}{\varepsilon }\right) ^{\frac{1}{p-2}%
},
\end{equation*}%
which altogether shows that%
\begin{eqnarray*}
C_{\widehat{\mathbf{b}}}\left( \infty \right) &\leq &2\left( \frac{p}{p-2}C_{%
\mathbf{b}}\left( p\right) ^{p}\frac{1}{\varepsilon }\right) ^{\frac{1}{p-2}%
}=2\left( \frac{p}{p-2}\right) ^{\frac{1}{p-2}}C_{\mathbf{b}}\left( p\right)
^{\frac{p}{p-2}}\varepsilon ^{-\frac{1}{p-2}}, \\
&&\ \ \ \ \ \ \ \ \ \ \ \ \ \ \ \ \ \ \ \ \text{for }0<\varepsilon \leq 
\frac{1}{4}.
\end{eqnarray*}%
Similarly we have%
\begin{eqnarray*}
C_{\widehat{\mathbf{b}}^{\ast }}\left( \infty \right) &\leq &2\left( \frac{p%
}{p-2}C_{\mathbf{b}^{\ast }}\left( p\right) ^{p}\frac{1}{\varepsilon ^{\ast }%
}\right) ^{\frac{1}{p-2}}=2\left( \frac{p}{p-2}\right) ^{\frac{1}{p-2}}C_{%
\mathbf{b}^{\ast }}\left( p\right) ^{\frac{p}{p-2}}\left( \varepsilon ^{\ast
}\right) ^{-\frac{1}{p-2}} \\
&&\ \ \ \ \ \ \ \ \ \ \ \ \ \ \ \ \ \ \ \ \text{for }0<\varepsilon ^{\ast
}\leq \frac{1}{4}.
\end{eqnarray*}

Moreover, we also have%
\begin{eqnarray*}
\sqrt{\int_{Q}\left\vert T_{\sigma }^{\alpha }\widehat{b}_{Q}\right\vert
^{2}d\omega } &\leq &2\sqrt{\int_{Q}\left\vert T_{\sigma }^{\alpha
}b_{Q}\right\vert ^{2}d\omega }+2\sqrt{\int_{Q}\left\vert T_{\sigma
}^{\alpha }\mathbf{1}_{\left\{ \left\vert b_{Q}\right\vert >\lambda \right\}
}\left( \frac{\lambda }{\left\vert b_{Q}\right\vert }-1\right)
b_{Q}\right\vert ^{2}d\omega } \\
&\leq &2\mathfrak{T}_{T^{\alpha }}^{\mathbf{b}}\sqrt{\left\vert Q\right\vert
_{\sigma }}+2\mathfrak{N}_{T^{\alpha }}\sqrt{\int_{\left\{ \left\vert
b_{Q}\right\vert >\lambda \right\} }\left\vert b_{Q}\right\vert ^{2}d\sigma }
\\
&=&2\left\{ \mathfrak{T}_{T^{\alpha }}^{\mathbf{b}}+\sqrt{\varepsilon }%
\mathfrak{N}_{T^{\alpha }}\right\} \sqrt{\left\vert Q\right\vert _{\sigma }}%
\ ,\ \ \ \ \ \text{for all intervals }Q,
\end{eqnarray*}%
which shows that%
\begin{equation}
\mathfrak{T}_{T^{\alpha }}^{\widehat{\mathbf{b}}}\leq 2\mathfrak{T}%
_{T^{\alpha }}^{\mathbf{b}}+2\sqrt{\varepsilon }\mathfrak{N}_{T^{\alpha }}\ .
\label{test}
\end{equation}

Now we take $\varepsilon =\varepsilon ^{\ast }$ and apply the fact that $%
\left( \mathcal{S}_{\infty }\right) $ holds to obtain%
\begin{eqnarray*}
\mathfrak{N}_{T^{\alpha }} &\lesssim &\left( C_{\widehat{\mathbf{b}}}\left(
\infty \right) +C_{\widehat{\mathbf{b}}^{\ast }}\left( \infty \right)
\right) \left\{ \mathfrak{T}_{T^{\alpha }}^{\widehat{\mathbf{b}}}+\mathfrak{T%
}_{T^{\alpha ,\ast }}^{\widehat{\mathbf{b}}^{\ast }}+\sqrt{\mathfrak{A}%
_{2}^{\alpha }}+\mathfrak{E}_{2}^{\alpha }\right\} \\
&\lesssim &\left( C_{\mathbf{b}}\left( p\right) +C_{\mathbf{b}^{\ast
}}\left( p\right) \right) ^{\frac{p}{p-2}}\varepsilon ^{-\frac{1}{p-2}%
}\left\{ \left[ \mathfrak{T}_{T^{\alpha }}^{\mathbf{b}}+\sqrt{\varepsilon }%
\mathfrak{N}_{T^{\alpha }}\right] +\left[ \mathfrak{T}_{T^{\alpha ,\ast }}^{%
\mathbf{b}^{\ast }}+\sqrt{\varepsilon }\mathfrak{N}_{T^{\alpha }}\right] +%
\sqrt{\mathfrak{A}_{2}^{\alpha }}+\mathfrak{E}_{2}^{\alpha }\right\} \\
&\lesssim &\left( C_{\mathbf{b}}\left( p\right) +C_{\mathbf{b}^{\ast
}}\left( p\right) \right) ^{\frac{p}{p-2}}\varepsilon ^{-\frac{1}{p-2}%
}\left\{ \mathfrak{T}_{T^{\alpha }}^{\mathbf{b}}+\mathfrak{T}_{T^{\alpha
,\ast }}^{\mathbf{b}^{\ast }}+\sqrt{\mathfrak{A}_{2}^{\alpha }}+\mathfrak{E}%
_{2}^{\alpha }\right\} \\
&&\ \ \ \ \ \ \ \ \ \ \ \ \ \ \ \ \ \ \ \ +\left( C_{\mathbf{b}}\left(
p\right) +C_{\mathbf{b}^{\ast }}\left( p\right) \right) ^{\frac{p}{p-2}%
}\varepsilon ^{\frac{1}{2}-\frac{1}{p-2}}\mathfrak{N}_{T^{\alpha }}\ .
\end{eqnarray*}%
Now we choose 
\begin{equation*}
\varepsilon =\frac{1}{\Gamma }\left( C_{\mathbf{b}}\left( p\right) +C_{%
\mathbf{b}^{\ast }}\left( p\right) \right) ^{-\frac{\frac{p}{p-2}}{\frac{1}{2%
}-\frac{1}{p-2}}}
\end{equation*}%
with $\Gamma $ large enough, depending only on the implied constant $C_{%
\limfunc{implied}}$ (where $\lesssim $ is written $\leq C_{\limfunc{implied}%
} $), so that the final term on the right satisfies%
\begin{eqnarray*}
&&C_{\limfunc{implied}}\left( C_{\mathbf{b}}\left( p\right) +C_{\mathbf{b}%
^{\ast }}\left( p\right) \right) ^{\frac{p}{p-2}}\ \varepsilon ^{\frac{1}{2}-%
\frac{1}{p-2}}\ \mathfrak{N}_{T^{\alpha }} \\
&=&C_{\limfunc{implied}}\left( C_{\mathbf{b}}\left( p\right) +C_{\mathbf{b}%
^{\ast }}\left( p\right) \right) ^{\frac{p}{p-2}}\ \left( \frac{1}{\Gamma }%
\right) ^{\frac{1}{2}-\frac{1}{p-2}}\left( C_{\mathbf{b}}\left( p\right) +C_{%
\mathbf{b}^{\ast }}\left( p\right) \right) ^{-\frac{p}{p-2}}\ \mathfrak{N}%
_{T^{\alpha }} \\
&\leq &C_{\limfunc{implied}}\ \left( \frac{1}{\Gamma }\right) ^{\frac{1}{4}%
}\ \mathfrak{N}_{T^{\alpha }}=\frac{1}{2}\mathfrak{N}_{T^{\alpha }}\ ,
\end{eqnarray*}%
i.e, we choose $\Gamma =\left( 2C_{\limfunc{implied}}\right) ^{4}$, and
where we have used $\frac{1}{2}-\frac{1}{p-2}\geq \frac{1}{4}$ for $p>6$.
This term can then be absorbed into the left hand side to obtain%
\begin{eqnarray*}
\mathfrak{N}_{T^{\alpha }} &\lesssim &\left( C_{\mathbf{b}}\left( p\right)
+C_{\mathbf{b}^{\ast }}\left( p\right) \right) ^{\frac{p}{p-2}}\left( \left(
C_{\mathbf{b}}\left( p\right) +C_{\mathbf{b}^{\ast }}\left( p\right) \right)
^{-\frac{\frac{p}{p-2}}{\frac{1}{2}-\frac{1}{p-2}}}\right) ^{-\frac{1}{p-2}}
\\
&&\ \ \ \ \ \ \ \ \ \ \ \ \ \ \ \ \ \ \ \ \times \left\{ \mathfrak{T}%
_{T^{\alpha }}^{\mathbf{b}}+\mathfrak{T}_{T^{\alpha ,\ast }}^{\mathbf{b}%
^{\ast }}+\sqrt{\mathfrak{A}_{2}^{\alpha }}+\mathfrak{E}_{2}^{\alpha
}\right\} \\
&\lesssim &\left( C_{\mathbf{b}}\left( p\right) +C_{\mathbf{b}^{\ast
}}\left( p\right) \right) ^{\frac{p}{p-2}\left\{ 1+\frac{\frac{1}{p-2}}{%
\frac{1}{2}-\frac{1}{p-2}}\right\} }\left\{ \mathfrak{T}_{T^{\alpha }}^{%
\mathbf{b}}+\mathfrak{T}_{T^{\alpha ,\ast }}^{\mathbf{b}^{\ast }}+\sqrt{%
\mathfrak{A}_{2}^{\alpha }}+\mathfrak{E}_{2}^{\alpha }\right\}
\end{eqnarray*}%
Since 
\begin{equation*}
\frac{p}{p-2}\left\{ 1+\frac{\frac{1}{p-2}}{\frac{1}{2}-\frac{1}{p-2}}%
\right\} =\left( 1+\frac{2}{p-2}\right) \left( 1+\frac{2}{p-4}\right) \leq 3%
\text{ for }p>6,
\end{equation*}%
we get%
\begin{equation*}
\mathfrak{N}_{T^{\alpha }}\lesssim \left( C_{\mathbf{b}}\left( p\right) +C_{%
\mathbf{b}^{\ast }}\left( p\right) \right) ^{3}\left\{ \mathfrak{T}%
_{T^{\alpha }}^{\mathbf{b}}+\mathfrak{T}_{T^{\alpha ,\ast }}^{\mathbf{b}%
^{\ast }}+\sqrt{\mathfrak{A}_{2}^{\alpha }}+\mathfrak{E}_{2}^{\alpha
}\right\} ,
\end{equation*}%
which completes the proof of $\left( \mathcal{S}_{0}\right) $.

Suppose now, in order to derive a contradiction, that $\left( \mathcal{S}%
_{n}\right) $ fails for some $n\geq 1$. Then there exists an integer $n\geq
1 $ and an exponent $r\in \left( p_{n+1},p_{n}\right] $ such that%
\begin{equation*}
\left( \mathcal{S}_{r}\right) \text{ fails and }\left( \mathcal{S}%
_{q}\right) \text{ holds\ for }q>p_{n}.
\end{equation*}%
We now show that $\left( \mathcal{S}_{p}\right) $ holds for all $p\in \left(
p_{n+1},p_{n}\right] $, contradicting the fact that $\left( \mathcal{S}%
_{r}\right) $ fails. So fix $n\geq 1$, $p\in \left( p_{n+1},p_{n}\right] $,
and suppose that $\mathbf{b}=\left\{ b_{Q}\right\} _{Q\in \mathcal{P}}$ is a 
$p$-weakly $\sigma $-accretive family of functions on $\mathbb{R}$, and that 
$\mathbf{b}^{\ast }=\left\{ b_{Q}^{\ast }\right\} _{Q\in \mathcal{P}}$ is a $%
p$-weakly $\omega $-accretive family of functions on $\mathbb{R}$. Note that
the sequence $\left\{ p_{n}\right\} _{n=0}^{\infty }=\left\{ \frac{2}{%
1-\left( \frac{2}{3}\right) ^{n}}\right\} _{n=0}^{\infty }$ satisfies the
recursion relation%
\begin{equation*}
p_{n+1}=\frac{6}{1+\frac{4}{p_{n}}},\text{ equivalently }p_{n}=\frac{4}{%
\frac{6}{p_{n+1}}-1},\text{\ \ \ \ \ }n\geq 0.
\end{equation*}%
Choose $q\in \left( p_{n},p_{n-1}\right] $ so that%
\begin{equation}
p>\frac{6}{1+\frac{4}{q}},\text{ i.e. }q<\frac{4}{\frac{6}{p}-1},
\label{restrict}
\end{equation}%
which can be done since $p>p_{n+1}=\frac{2}{1-\left( \frac{2}{3}\right)
^{n+1}}$ is equivalent to $p_{n}=\frac{2}{1-\left( \frac{2}{3}\right) ^{n}}<%
\frac{4}{\frac{6}{p}-1}$, which leaves room to choose $q$ satisfying $%
p_{n}<q<\frac{4}{\frac{6}{p}-1}$.

Now let $0<\varepsilon <1$ (to be fixed later), define $\lambda =\lambda
\left( \varepsilon \right) $ as in (\ref{lambda choice}), and define $%
\widehat{b}_{Q}$ as in (\ref{new}). Recall from (\ref{hence}) and (\ref{low}%
) that we then have%
\begin{equation*}
\int_{\left\{ \left\vert b_{Q}\right\vert >\lambda \right\} }\left\vert
b_{Q}\right\vert ^{2}d\sigma \leq \varepsilon \left\vert Q\right\vert
_{\sigma }\ ,
\end{equation*}%
and%
\begin{equation*}
\left\vert \frac{1}{\left\vert Q\right\vert _{\sigma }}\int_{Q}\widehat{b}%
_{Q}d\sigma \right\vert \geq 1,\ \ \ \ \ Q\in \mathcal{P}\ ,
\end{equation*}%
if we choose $0<\varepsilon \leq \frac{1}{4}$. We of course have the
previous upper bound 
\begin{equation*}
\left\Vert \widehat{b}_{Q}\right\Vert _{L^{\infty }\left( \sigma \right)
}\leq 2\lambda =2\lambda \left( \varepsilon \right) =2\left( \frac{p}{p-2}C_{%
\mathbf{b}}\left( p\right) ^{p}\frac{1}{\varepsilon }\right) ^{\frac{1}{p-2}%
},
\end{equation*}%
and while this turned out to be sufficient in the case $n=0$, we must do
better than $O\left( \frac{1}{\varepsilon }\right) ^{\frac{1}{p-2}}$ in the
case $n\geq 1$. In fact we compute the $L^{q}$ norm instead, recalling that $%
q>p$:%
\begin{eqnarray*}
&&\left( \frac{1}{\left\vert Q\right\vert _{\mu }}\int_{Q}\left\vert 
\widehat{b}_{Q}\right\vert ^{q}d\mu \right) ^{\frac{1}{q}}=2\left( \frac{1}{%
\left\vert Q\right\vert _{\mu }}\int_{Q}\left\vert b_{Q}\left( \mathbf{1}%
_{\left\{ \left\vert b_{Q}\right\vert \leq \lambda \right\} }+\frac{\lambda 
}{\left\vert b_{Q}\right\vert }\mathbf{1}_{\left\{ \left\vert
b_{Q}\right\vert >\lambda \right\} }\right) \right\vert ^{q}d\mu \right) ^{%
\frac{1}{q}} \\
&=&2\left( \frac{1}{\left\vert Q\right\vert _{\mu }}\int_{\left\{ \left\vert
b_{Q}\right\vert \leq \lambda \right\} }\left[ \int_{0}^{\left\vert
b_{Q}\right\vert }qt^{q-1}dt\right] d\sigma +\frac{\lambda ^{q}\left\vert
\left\{ \left\vert b_{Q}\right\vert >\lambda \right\} \right\vert _{\mu }}{%
\left\vert Q\right\vert _{\mu }}\right) ^{\frac{1}{q}} \\
&\leq &2\left( \frac{1}{\left\vert Q\right\vert _{\mu }}\int_{0}^{\lambda }%
\left[ \int_{\left\{ t<\left\vert b_{Q}\right\vert \leq \lambda \right\}
}d\sigma \right] qt^{q-1}dt+C_{\mathbf{b}}\left( p\right) ^{p}\lambda
^{q-p}\right) ^{\frac{1}{q}} \\
&\leq &2\left( \frac{1}{\left\vert Q\right\vert _{\mu }}\int_{0}^{\lambda }%
\left[ \frac{1}{t^{p}}\int \left\vert b_{Q}\right\vert ^{p}d\sigma \right]
qt^{q-1}dt+C_{\mathbf{b}}\left( p\right) ^{p}\lambda ^{q-p}\right) ^{\frac{1%
}{q}} \\
&\leq &2C_{\mathbf{b}}\left( p\right) ^{\frac{p}{q}}\left( \int_{0}^{\lambda
}qt^{q-p-1}dt+\lambda ^{q-p}\right) ^{\frac{1}{q}}=2\left( C_{\mathbf{b}%
}\left( p\right) ^{p}\frac{2q-p}{q-p}\lambda ^{q-p}\right) ^{\frac{1}{q}},
\end{eqnarray*}%
which shows that $C_{\widehat{\mathbf{b}}}\left( q\right) $ satisfies the
estimate%
\begin{eqnarray*}
C_{\widehat{\mathbf{b}}}\left( q\right) &\leq &2C_{\mathbf{b}}\left(
p\right) ^{\frac{p}{q}}\left( \frac{2q-p}{q-p}\right) ^{\frac{1}{q}}\left[
\left( \frac{p}{p-2}C_{\mathbf{b}}\left( p\right) ^{p}\frac{1}{\varepsilon }%
\right) ^{\frac{1}{p-2}}\right] ^{1-\frac{p}{q}} \\
&\lesssim &C_{\mathbf{b}}\left( p\right) ^{\frac{p}{q}\left( 1+\frac{q-p}{p-2%
}\right) }\varepsilon ^{-\frac{1-\frac{p}{q}}{p-2}}\lesssim C_{\mathbf{b}%
}\left( p\right) ^{\frac{3}{2}}\varepsilon ^{-\frac{1-\frac{p}{q}}{p-2}},
\end{eqnarray*}%
a significant improvement over the bound $O\left( \varepsilon ^{-\frac{1}{p-2%
}}\right) $. Here we have used that if $p>\frac{6}{1+\frac{4}{q}}=\frac{6q}{%
q+4}$, then%
\begin{eqnarray*}
\frac{p}{q}\left( 1+\frac{q-p}{p-2}\right) &=&\frac{p}{p-2}\frac{q-2}{q}<%
\frac{\frac{6q}{q+4}}{\frac{6q}{q+4}-2}\frac{q-2}{q} \\
&=&\frac{6q}{6q-\left( 2q+8\right) }\frac{q-2}{q}=\frac{6\left( q-2\right) }{%
4q-8}=\frac{3}{2}.
\end{eqnarray*}

Moreover, from (\ref{test}) we also have%
\begin{equation*}
\mathfrak{T}_{T^{\alpha }}^{\widehat{\mathbf{b}}}\leq 2\mathfrak{T}%
_{T^{\alpha }}^{\mathbf{b}}+2\sqrt{\varepsilon }\mathfrak{N}_{T^{\alpha }}\ .
\end{equation*}%
We can do the same for the dual testing functions $\mathbf{b}^{\ast
}=\left\{ b_{Q}^{\ast }\right\} _{Q\in \mathcal{P}}$, and then altogether we
have both%
\begin{eqnarray*}
1 &\leq &\left\vert \frac{1}{\left\vert Q\right\vert _{\sigma }}\int_{Q}%
\widehat{b}_{Q}d\sigma \right\vert \leq \left\Vert \widehat{b}%
_{Q}\right\Vert _{L^{q}\left( \sigma \right) }\leq C_{\mathbf{b}}\left(
p\right) ^{\frac{3}{2}}\varepsilon ^{-\frac{1-\frac{p}{q}}{p-2}},\ \ \ \ \
Q\in \mathcal{P}\ , \\
&&\ \ \ \ \ \ \ \ \ \ \mathfrak{T}_{T^{\alpha }}^{\widehat{\mathbf{b}}}\leq 2%
\mathfrak{T}_{T^{\alpha }}^{\mathbf{b}}+2\sqrt{\varepsilon }\mathfrak{N}%
_{T^{\alpha }}\ ,
\end{eqnarray*}%
as well as%
\begin{eqnarray*}
1 &\leq &\left\vert \frac{1}{\left\vert Q\right\vert _{\omega }}\int_{Q}%
\widehat{b^{\ast }}_{Q}d\omega \right\vert \leq \left\Vert \widehat{b^{\ast }%
}_{Q}\right\Vert _{L^{q}\left( \omega \right) }\leq C_{\mathbf{b}^{\ast
}}\left( p\right) ^{\frac{3}{2}}\varepsilon ^{-\frac{1-\frac{p}{q}}{p-2}},\
\ \ \ \ Q\in \mathcal{P}\ , \\
&&\ \ \ \ \ \ \ \ \ \ \mathfrak{T}_{T^{\alpha }}^{\widehat{\mathbf{b}^{\ast }%
}}\leq 2\mathfrak{T}_{T^{\alpha }}^{\mathbf{b}^{\ast }}+2\sqrt{\varepsilon
^{\ast }}\mathfrak{N}_{T^{\alpha }}\ ,
\end{eqnarray*}%
provided%
\begin{equation}
0<\varepsilon =\varepsilon ^{\ast }\leq \frac{1}{4}\text{ }.  \label{eps}
\end{equation}

We now use these estimates, together with the fact that $\left( \mathcal{S}%
_{n-1}\right) $ holds, to obtain%
\begin{eqnarray*}
\mathfrak{N}_{T^{\alpha }} &\lesssim &\left( C_{\widehat{\mathbf{b}}}\left(
q\right) +C_{\widehat{\mathbf{b}}^{\ast }}\left( q\right) \right)
^{3^{n}}\left\{ \mathfrak{T}_{T^{\alpha }}^{\widehat{\mathbf{b}}}+\mathfrak{T%
}_{T^{\alpha ,\ast }}^{\widehat{\mathbf{b}}^{\ast }}+\sqrt{\mathfrak{A}%
_{2}^{\alpha }}+\mathfrak{E}_{2}^{\alpha }\right\} \\
&\lesssim &\left( C_{\mathbf{b}}\left( p\right) +C_{\mathbf{b}^{\ast
}}\left( p\right) \right) ^{\frac{3}{2}3^{n}}\varepsilon ^{-\frac{1-\frac{p}{%
q}}{p-2}} \\
&&\ \ \ \ \ \ \ \ \ \ \ \ \ \ \ \times \left\{ \left[ \mathfrak{T}%
_{T^{\alpha }}^{\mathbf{b}}+\sqrt{\varepsilon }\mathfrak{N}_{T^{\alpha }}%
\right] +\left[ \mathfrak{T}_{T^{\alpha ,\ast }}^{\mathbf{b}^{\ast }}+\sqrt{%
\varepsilon }\mathfrak{N}_{T^{\alpha }}\right] +\sqrt{\mathfrak{A}%
_{2}^{\alpha }}+\mathfrak{E}_{2}^{\alpha }\right\} \\
&\lesssim &\left( C_{\mathbf{b}}\left( p\right) +C_{\mathbf{b}^{\ast
}}\left( p\right) \right) ^{\frac{3}{2}3^{n}}\varepsilon ^{-\frac{1-\frac{p}{%
q}}{p-2}}\left\{ \mathfrak{T}_{T^{\alpha }}^{\mathbf{b}}+\mathfrak{T}%
_{T^{\alpha ,\ast }}^{\mathbf{b}^{\ast }}+\sqrt{\mathfrak{A}_{2}^{\alpha }}+%
\mathfrak{E}_{2}^{\alpha }\right\} \\
&&\ \ \ \ \ \ \ \ \ \ \ \ \ \ \ \ \ \ \ \ +\left( C_{\mathbf{b}}\left(
p\right) +C_{\mathbf{b}^{\ast }}\left( p\right) \right) ^{\frac{3}{2}3^{n}}%
\sqrt{\varepsilon }\varepsilon ^{-\frac{1-\frac{p}{q}}{p-2}}\mathfrak{N}%
_{T^{\alpha }}\ .
\end{eqnarray*}%
We can absorb the last term on the right hand side above into the left hand
side for $\varepsilon >0$ sufficiently small, since (\ref{restrict}) gives $%
\frac{\frac{6}{p}-1}{2}<\frac{2}{q}$, and hence 
\begin{equation}
\frac{1}{2}-\frac{1-\frac{p}{q}}{p-2}=\frac{p\left( 1+\frac{2}{q}\right) -4}{%
2p-4}>\frac{p\left( 1+\frac{\frac{6}{p}-1}{2}\right) -4}{2p-4}=\frac{1}{4}.
\label{quarter}
\end{equation}%
In fact, we choose 
\begin{equation*}
\varepsilon =\frac{1}{\Gamma }\left( C_{\mathbf{b}}\left( p\right) +C_{%
\mathbf{b}^{\ast }}\left( p\right) \right) ^{\left[ \frac{3}{2}3^{n}\right] %
\left[ \frac{1}{\frac{1-\frac{p}{q}}{p-2}-\frac{1}{2}}\right] }
\end{equation*}%
with $\Gamma $ sufficiently large, depending only on the implied constant,
to get%
\begin{eqnarray*}
\mathfrak{N}_{T^{\alpha }} &\lesssim &\left( C_{\mathbf{b}}\left( p\right)
+C_{\mathbf{b}^{\ast }}\left( p\right) \right) ^{\left[ \frac{3}{2}3^{n}%
\right] }\left( \left( \left( C_{\mathbf{b}}\left( p\right) +C_{\mathbf{b}%
^{\ast }}\left( p\right) \right) ^{\left[ \frac{3}{2}3^{n}\right] \left[ 
\frac{1}{\frac{1-\frac{p}{q}}{p-2}-\frac{1}{2}}\right] }\right) ^{-\frac{1-%
\frac{p}{q}}{p-2}}\right) \\
&&\ \ \ \ \ \ \ \ \ \ \ \ \ \ \ \ \ \ \ \ \times \left\{ \mathfrak{T}%
_{T^{\alpha }}^{\mathbf{b}}+\mathfrak{T}_{T^{\alpha ,\ast }}^{\mathbf{b}%
^{\ast }}+\sqrt{\mathfrak{A}_{2}^{\alpha }}+\mathfrak{E}_{2}^{\alpha
}\right\} \\
&\lesssim &\left( C_{\mathbf{b}}\left( p\right) +C_{\mathbf{b}^{\ast
}}\left( p\right) \right) ^{\left[ \frac{3}{2}3^{n}\right] \left[ 1+1\right]
}\left\{ \mathfrak{T}_{T^{\alpha }}^{\mathbf{b}}+\mathfrak{T}_{T^{\alpha
,\ast }}^{\mathbf{b}^{\ast }}+\sqrt{\mathfrak{A}_{2}^{\alpha }}+\mathfrak{E}%
_{2}^{\alpha }\right\} .
\end{eqnarray*}%
Here we have used that (\ref{quarter}) applied twice implies%
\begin{equation*}
\frac{\frac{1-\frac{p}{q}}{p-2}}{\frac{1}{2}-\frac{1-\frac{p}{q}}{p-2}}<4%
\frac{1-\frac{p}{q}}{p-2}\leq 1.
\end{equation*}%
So we finally have%
\begin{equation*}
\mathfrak{N}_{T^{\alpha }}\lesssim \left( C_{\mathbf{b}}\left( p\right) +C_{%
\mathbf{b}^{\ast }}\left( p\right) \right) ^{3^{n+1}}\left\{ \mathfrak{T}%
_{T^{\alpha }}^{\mathbf{b}}+\mathfrak{T}_{T^{\alpha ,\ast }}^{\mathbf{b}%
^{\ast }}+\sqrt{\mathfrak{A}_{2}^{\alpha }}+\mathfrak{E}_{2}^{\alpha
}\right\} ,
\end{equation*}%
which completes the proof of Proposition \ref{conditional}.
\end{proof}

\begin{remark}
Propositions \ref{lower bound}\ and \ref{conditional} extend to higher
dimensions with analogous proofs.
\end{remark}

\begin{conclusion}
\label{bounded PLBP}We may assume for the proof of Theorem \ref{dim one}
given below that $p=\infty $ and that the testing functions are real-valued,
satisfy the $PLBP$ and satisfy%
\begin{eqnarray}
&&\limfunc{support}b_{Q}\subset Q\ ,\ \ \ \ \ Q\in \mathcal{P},
\label{acc infinity'} \\
1 &\leq &\frac{1}{\left\vert Q\right\vert _{\mu }}\int_{Q}b_{Q}d\mu \leq
\left\Vert b_{Q}\right\Vert _{L^{\infty }\left( \mu \right) }\leq C_{\mathbf{%
b}}\left( \infty \right) <\infty ,\ \ \ \ \ Q\in \mathcal{P}\ .  \notag
\end{eqnarray}
\end{conclusion}

\subsection{Reverse H\"{o}lder control of children}

Here we begin to further reduce the proof of Theorem \ref{dim one} to the
case of bounded real testing functions $\mathbf{b}=\left\{ b_{Q}\right\}
_{Q\in \mathcal{P}}$ having reverse H\"{o}lder control 
\begin{equation}
\left\vert \frac{1}{\left\vert Q^{\prime }\right\vert _{\sigma }}%
\int_{Q^{\prime }}b_{Q}d\sigma \right\vert \geq c\left\Vert \mathbf{1}%
_{Q^{\prime }}b_{Q}\right\Vert _{L^{\infty }\left( \sigma \right) }>0,
\label{rev Hol con}
\end{equation}%
for all children $Q^{\prime }\in \mathfrak{C}\left( Q\right) $ with $%
\left\vert Q^{\prime }\right\vert _{\sigma }>0$ and $Q\in \mathcal{P}$.

\subsubsection{Control of averages over children}

Here we address the case of a single interval $Q$.

\begin{lemma}
\label{further red}Suppose that $\sigma $ and $\omega $ are locally finite
positive Borel measures on the real line $\mathbb{R}$. Assume that $%
T^{\alpha }$ is a standard $\alpha $-fractional elliptic and gradient
elliptic singular integral operator on $\mathbb{R}$, and set $T_{\sigma
}^{\alpha }f=T^{\alpha }\left( f\sigma \right) $ for any smooth truncation
of $T_{\sigma }^{\alpha }$, so that $T_{\sigma }^{\alpha }$ is \emph{apriori}
bounded from $L^{2}\left( \sigma \right) $ to $L^{2}\left( \omega \right) $.
Let $Q\in \mathcal{P}$ and let $\mathfrak{N}_{T^{\alpha }}\left( Q\right) $
be the best constant in the local inequality%
\begin{equation*}
\sqrt{\int_{Q^{\prime }}\left\vert T_{\sigma }^{\alpha }\left( \mathbf{1}%
_{Q}f\right) \right\vert ^{2}d\omega }\leq \mathfrak{N}_{T^{\alpha }}\left(
Q\right) \sqrt{\int_{Q}\left\vert f\right\vert ^{2}d\sigma }\ ,\ \ \ \ \
f\in L^{2}\left( \mathbf{1}_{Q}\sigma \right) .
\end{equation*}%
Suppose that $b_{Q}$ is a real-valued function supported in $Q$ such that 
\begin{eqnarray*}
&&1\leq \frac{1}{\left\vert Q\right\vert _{\sigma }}\int_{Q}b_{Q}d\sigma
\leq \left\Vert \mathbf{1}_{Q}b_{Q}\right\Vert _{L^{\infty }\left( \sigma
\right) }\leq C_{b_{Q}}\ , \\
&&\sqrt{\int_{Q}\left\vert T_{\sigma }^{\alpha }b_{Q}\right\vert ^{2}d\omega 
}\leq \mathfrak{T}_{T^{\alpha }}^{b_{Q}}\left( Q\right) \sqrt{\left\vert
Q\right\vert _{\sigma }}\ .
\end{eqnarray*}%
Then for every $0<\delta <\frac{1}{4C_{\mathbf{b}}^{3}}$, there exists a
real-valued function $\widetilde{b}_{Q}$ supported in $Q$ such that 
\begin{eqnarray*}
&&1\leq \frac{1}{\left\vert Q\right\vert _{\sigma }}\int_{Q}\widetilde{b}%
_{Q}d\sigma \leq \left\Vert \mathbf{1}_{Q}\widetilde{b}_{Q}\right\Vert
_{L^{\infty }\left( \sigma \right) }\leq 2\left( 1+\sqrt{C_{b_{Q}}}\right)
C_{b_{Q}}\ , \\
&&\sqrt{\int_{Q}\left\vert T_{\sigma }^{\alpha }\widetilde{b}_{Q}\right\vert
^{2}d\omega }\leq \left[ 2\mathfrak{T}_{T^{\alpha }}^{b_{Q}}\left( Q\right)
+4C_{b_{Q}}^{\frac{3}{2}}\delta ^{\frac{1}{4}}\mathfrak{N}_{T^{\alpha
}}\left( Q\right) \right] \sqrt{\left\vert Q\right\vert _{\sigma }}\ , \\
&&0<\left\Vert \mathbf{1}_{Q_{i}}\widetilde{b}_{Q}\right\Vert _{L^{\infty
}\left( \sigma \right) }\leq \frac{16C_{b_{Q}}}{\delta }\left\vert \frac{1}{%
\left\vert Q_{i}\right\vert _{\sigma }}\int_{Q_{i}}\widetilde{b}_{Q}d\sigma
\right\vert \ ,\ \ \ \ \ Q_{i}\in \mathfrak{C}\left( Q\right) .
\end{eqnarray*}%
Moreover, if $\left\vert b_{Q}\right\vert \geq c_{1}>0$, then we may take $%
\left\vert \widetilde{b}_{Q}\right\vert \geq c_{1}$ as well.
\end{lemma}

\begin{proof}
Let $0<\delta <1$ and fix $Q\in \mathcal{P}$. By assumption we have%
\begin{equation*}
1\leq \frac{1}{\left\vert Q\right\vert _{\sigma }}\int_{Q}b_{Q}d\sigma \leq
\left\Vert \mathbf{1}_{Q}b_{Q}\right\Vert _{L^{\infty }\left( \sigma \right)
}\leq C_{b_{Q}}.
\end{equation*}%
Let $Q_{\limfunc{left}}$ and $Q_{\limfunc{right}}$ be the children of $Q$.
We now define $\widetilde{b}_{Q}$. First we note that the inequality%
\begin{equation}
\left\vert \frac{1}{\left\vert Q^{\prime }\right\vert _{\sigma }}%
\int_{Q^{\prime }}b_{Q}d\sigma \right\vert <\frac{\delta }{C_{b_{Q}}}%
\left\Vert \mathbf{1}_{Q^{\prime }}b_{Q}\right\Vert _{L^{\infty }\left(
\sigma \right) }  \label{Q' big}
\end{equation}%
cannot hold for $Q^{\prime }$ equal to both $Q_{\limfunc{left}}$ and $Q_{%
\limfunc{right}}$, since otherwise we obtain the contradiction 
\begin{eqnarray*}
\left\vert \int_{Q}b_{Q}d\sigma \right\vert &\leq &\left\vert \int_{Q_{%
\limfunc{left}}}b_{Q}d\sigma \right\vert +\left\vert \int_{Q_{\limfunc{right}%
}}b_{Q}d\sigma \right\vert \\
&<&\frac{\delta }{C_{\mathbf{b}}}\left( \left\vert Q_{\limfunc{left}%
}\right\vert _{\sigma }\left\Vert \mathbf{1}_{Q_{\limfunc{left}%
}}b_{Q}\right\Vert _{L^{\infty }\left( \sigma \right) }+\left\vert Q_{%
\limfunc{right}}\right\vert _{\sigma }\left\Vert \mathbf{1}_{Q_{\limfunc{%
right}}}b_{Q}\right\Vert _{L^{\infty }\left( \sigma \right) }\right) \\
&\leq &\frac{\delta }{C_{b_{Q}}}\left\vert Q\right\vert _{\sigma }\left\Vert 
\mathbf{1}_{Q}b_{Q}\right\Vert _{L^{\infty }\left( \sigma \right) }\leq
\delta \left\vert \int_{Q}b_{Q}d\sigma \right\vert <\left\vert
\int_{Q}b_{Q}d\sigma \right\vert .
\end{eqnarray*}%
If (\ref{Q' big}) holds for neither $Q_{\limfunc{left}}$ nor $Q_{\limfunc{%
right}}$, then we simply define $\widetilde{b}_{Q}=b_{Q}$. If (\ref{Q' big})
holds for just one of the children, say $Q_{\limfunc{left}}$, then we define 
$\widetilde{b}_{Q}$ differently according to how large the $L^{1}\left(
\sigma \right) $-average $\frac{1}{\left\vert Q_{\limfunc{left}}\right\vert
_{\sigma }}\int_{Q_{\limfunc{left}}}\left\vert b_{Q}\right\vert d\sigma $ is.

\textbf{Case (0)} $\frac{1}{\left\vert Q_{\limfunc{left}}\right\vert
_{\sigma }}\int_{Q_{\limfunc{left}}}\left\vert b_{Q}\right\vert d\sigma =0$:
In this case we define%
\begin{equation*}
\widetilde{b}_{Q}\equiv \delta \mathbf{1}_{Q_{\limfunc{left}}}+b_{Q}\mathbf{1%
}_{Q_{\limfunc{right}}}\ ,
\end{equation*}%
and the reader can easily verify that the conclusions of Lemma \ref{further
red} hold.

\textbf{Case (1)} $0<\frac{1}{\left\vert Q_{\limfunc{left}}\right\vert
_{\sigma }}\int_{Q_{\limfunc{left}}}\left\vert b_{Q}\right\vert d\sigma \leq 
\sqrt{C_{b_{Q}}\delta }$: In this case we define 
\begin{equation*}
\widetilde{b}_{Q}\equiv \left( \frac{1}{\left\vert Q_{\limfunc{left}%
}\right\vert _{\sigma }}\int_{Q_{\limfunc{left}}}\left\vert b_{Q}\right\vert
d\sigma \right) \mathbf{1}_{Q_{\limfunc{left}}}+b_{Q}\mathbf{1}_{Q_{\limfunc{%
right}}}\ .
\end{equation*}

With this definition we then have%
\begin{eqnarray*}
1 &\leq &\frac{1}{\left\vert Q\right\vert _{\sigma }}\int_{Q}b_{Q}d\sigma
\leq \frac{1}{\left\vert Q\right\vert _{\sigma }}\left( \int_{Q_{\limfunc{%
left}}}\left\vert b_{Q}\right\vert d\sigma +\int_{Q_{\limfunc{right}%
}}b_{Q}d\sigma \right) \\
&=&\frac{1}{\left\vert Q\right\vert _{\sigma }}\int_{Q}\widetilde{b}%
_{Q}d\sigma \leq \left\Vert \widetilde{b}_{Q}\right\Vert _{L^{\infty }\left(
\sigma \right) }\leq \left\Vert b_{Q}\right\Vert _{L^{\infty }\left( \sigma
\right) }\leq C_{b_{Q}}\ ,
\end{eqnarray*}%
and both%
\begin{eqnarray*}
\frac{\left\Vert \mathbf{1}_{Q_{\limfunc{left}}}\widetilde{b}_{Q}\right\Vert
_{L^{\infty }\left( \sigma \right) }}{\left\vert \frac{1}{\left\vert Q_{%
\limfunc{left}}\right\vert _{\sigma }}\int_{Q_{\limfunc{left}}}\widetilde{b}%
_{Q}d\sigma \right\vert } &=&\frac{\frac{1}{\left\vert Q_{\limfunc{left}%
}\right\vert _{\sigma }}\int_{Q_{\limfunc{left}}}\left\vert b_{Q}\right\vert
d\sigma }{\left\vert \frac{1}{\left\vert Q_{\limfunc{left}}\right\vert
_{\sigma }}\int_{Q_{\limfunc{left}}}\left\vert b_{Q}\right\vert d\sigma
\right\vert }=1<\frac{1}{\delta }C_{b_{Q}}\text{\ }, \\
\frac{\left\Vert \mathbf{1}_{Q_{\limfunc{right}}}\widetilde{b}%
_{Q}\right\Vert _{L^{\infty }\left( \sigma \right) }}{\left\vert \frac{1}{%
\left\vert Q_{\limfunc{right}}\right\vert _{\sigma }}\int_{Q_{\limfunc{right}%
}}\widetilde{b}_{Q}d\sigma \right\vert } &=&\frac{\left\Vert \mathbf{1}_{Q_{%
\limfunc{right}}}b_{Q}\right\Vert _{L^{\infty }\left( \sigma \right) }}{%
\left\vert \frac{1}{\left\vert Q_{\limfunc{right}}\right\vert _{\sigma }}%
\int_{Q_{\limfunc{right}}}b_{Q}d\sigma \right\vert }<\frac{1}{\delta }%
C_{b_{Q}}\ ,
\end{eqnarray*}%
where the second line follows since (\ref{Q' big}) fails for $Q^{\prime }=Q_{%
\limfunc{right}}$.

Finally we check the testing condition in this case. We have from
Minkowski's inequality, 
\begin{eqnarray*}
\sqrt{\frac{1}{\left\vert Q\right\vert _{\sigma }}\int_{Q}\left\vert
T_{\sigma }^{\alpha }\widetilde{b}_{Q}\right\vert ^{2}d\omega } &\leq &\sqrt{%
\frac{1}{\left\vert Q\right\vert _{\sigma }}\int_{Q}\left\vert T_{\sigma
}^{\alpha }b_{Q}\right\vert ^{2}d\omega }+\sqrt{\frac{1}{\left\vert
Q\right\vert _{\sigma }}\int_{Q}\left\vert T_{\sigma }^{\alpha }\left( 
\widetilde{b}_{Q}-b_{Q}\right) \right\vert ^{2}d\omega } \\
&\leq &\mathfrak{T}_{T^{\alpha }}^{b_{Q}}\left( Q\right) +\mathfrak{N}%
_{T^{\alpha }}\left( Q\right) \sqrt{\frac{1}{\left\vert Q\right\vert
_{\sigma }}\int_{Q}\left\vert \widetilde{b}_{Q}-b_{Q}\right\vert ^{2}d\sigma 
}.
\end{eqnarray*}%
In the case (\ref{Q' big}) holds for neither $Q_{\limfunc{left}}$ nor $Q_{%
\limfunc{right}}$, then $\widetilde{b}_{Q}-b_{Q}=0$ and so%
\begin{equation*}
\sqrt{\frac{1}{\left\vert Q\right\vert _{\sigma }}\int_{Q}\left\vert
T_{\sigma }^{\alpha }\widetilde{b}_{Q}\right\vert ^{2}d\omega }\leq 
\mathfrak{T}_{T^{\alpha }}^{b_{Q}}\left( Q\right) \ .
\end{equation*}%
In the case (\ref{Q' big}) holds for just one child, say $Q_{\limfunc{left}}$%
, then%
\begin{eqnarray*}
\sqrt{\frac{1}{\left\vert Q\right\vert _{\sigma }}\int_{Q}\left\vert 
\widetilde{b}_{Q}-b_{Q}\right\vert ^{2}d\omega } &=&\sqrt{\frac{1}{%
\left\vert Q\right\vert _{\sigma }}\int_{Q_{\limfunc{left}}}\left\vert
\left( \frac{1}{\left\vert Q_{\limfunc{left}}\right\vert _{\sigma }}\int_{Q_{%
\limfunc{left}}}\left\vert b_{Q}\right\vert d\sigma \right)
-b_{Q}\right\vert ^{2}d\sigma } \\
&\leq &\sqrt{\frac{1}{\left\vert Q\right\vert _{\sigma }}\int_{Q_{\limfunc{%
left}}}\left\vert \frac{1}{\left\vert Q_{\limfunc{left}}\right\vert _{\sigma
}}\int_{Q_{\limfunc{left}}}\left\vert b_{Q}\right\vert d\sigma \right\vert
^{2}d\sigma }+\sqrt{\frac{1}{\left\vert Q\right\vert _{\sigma }}\int_{Q_{%
\limfunc{left}}}\left\vert b_{Q}\right\vert ^{2}d\sigma } \\
&\leq &\sqrt{\frac{1}{\left\vert Q\right\vert _{\sigma }}\int_{Q_{\limfunc{%
left}}}C_{\mathbf{b}}\delta d\sigma }+\sqrt{C_{\mathbf{b}}\frac{1}{%
\left\vert Q\right\vert _{\sigma }}\int_{Q_{\limfunc{left}}}\left\vert
b_{Q}\right\vert d\sigma } \\
&\leq &\sqrt{C_{b_{Q}}\delta }+\sqrt{C_{b_{Q}}\sqrt{C_{b_{Q}}\delta }}\leq
2C_{b_{Q}}^{\frac{3}{4}}\delta ^{\frac{1}{4}}.
\end{eqnarray*}

\textbf{Case (2)} $\frac{1}{\left\vert Q_{\limfunc{left}}\right\vert
_{\sigma }}\int_{Q_{\limfunc{left}}}\left\vert b_{Q}\right\vert d\sigma >%
\sqrt{C_{b_{Q}}\delta }$: Let 
\begin{eqnarray*}
\mathbf{1}_{Q_{\limfunc{left}}}\left( x\right) b_{Q}\left( x\right)
&=&p\left( x\right) -n\left( x\right) , \\
\mathbf{1}_{Q_{\limfunc{left}}}\left( x\right) \left\vert b_{Q}\left(
x\right) \right\vert &=&p\left( x\right) +n\left( x\right) ,
\end{eqnarray*}%
where $p\left( x\right) $ and $n\left( x\right) $ are the positive and
negative parts respectively of $b_{Q}$ on $Q_{\limfunc{left}}$. Then define $%
\widetilde{b}_{Q}$ by%
\begin{equation*}
\widetilde{b}_{Q}\equiv \left\{ 
\begin{array}{ccc}
\left( \frac{1}{\left\vert Q_{\limfunc{left}}\right\vert _{\sigma }}\int_{Q_{%
\limfunc{left}}}\left[ p-n\left( 1+\sqrt{C_{\mathbf{b}}\delta }\right) %
\right] d\sigma \right) \mathbf{1}_{Q_{\limfunc{left}}}+b_{Q}\mathbf{1}_{Q_{%
\limfunc{right}}} & \text{ if } & \int_{Q_{\limfunc{left}}}pd\sigma
<\int_{Q_{\limfunc{left}}}nd\sigma \\ 
\left( \frac{1}{\left\vert Q_{\limfunc{left}}\right\vert _{\sigma }}\int_{Q_{%
\limfunc{left}}}\left[ \left( 1+\sqrt{C_{\mathbf{b}}\delta }\right) p-n%
\right] d\sigma \right) \mathbf{1}_{Q_{\limfunc{left}}}+b_{Q}\mathbf{1}_{Q_{%
\limfunc{right}}} & \text{ if } & \int_{Q_{\limfunc{left}}}pd\sigma \geq
\int_{Q_{\limfunc{left}}}nd\sigma%
\end{array}%
\right. .
\end{equation*}

\textbf{Subcase (2a)} $\int_{Q_{\limfunc{left}}}pd\sigma <\int_{Q_{\limfunc{%
left}}}nd\sigma $: In this case we have%
\begin{eqnarray*}
1 &\leq &\frac{1}{\left\vert Q\right\vert _{\sigma }}\int_{Q}b_{Q}d\sigma =%
\frac{1}{\left\vert Q\right\vert _{\sigma }}\left( \int_{Q_{\limfunc{left}%
}}\left( p-n\right) d\sigma +\int_{Q_{\limfunc{right}}}b_{Q}d\sigma \right)
\\
&&-\frac{\sqrt{C_{b_{Q}}\delta }}{\left\vert Q\right\vert _{\sigma }}%
\int_{Q_{\limfunc{left}}}nd\sigma +\frac{\sqrt{C_{b_{Q}}\delta }}{\left\vert
Q\right\vert _{\sigma }}\int_{Q_{\limfunc{left}}}nd\sigma \\
&=&\frac{1}{\left\vert Q\right\vert _{\sigma }}\int_{Q}\widetilde{b}%
_{Q}d\sigma +\frac{\sqrt{C_{b_{Q}}\delta }}{\left\vert Q\right\vert _{\sigma
}}\int_{Q_{\limfunc{left}}}nd\sigma \leq \frac{1}{\left\vert Q\right\vert
_{\sigma }}\int_{Q}\widetilde{b}_{Q}d\sigma +\frac{\sqrt{C_{b_{Q}}\delta }}{%
\left\vert Q\right\vert _{\sigma }}C_{b_{Q}}\left\vert Q_{\limfunc{left}%
}\right\vert _{\sigma } \\
&\leq &\frac{1}{\left\vert Q\right\vert _{\sigma }}\int_{Q}\widetilde{b}%
_{Q}d\sigma +\sqrt{C_{b_{Q}}\delta }C_{b_{Q}}\leq \frac{1}{\left\vert
Q\right\vert _{\sigma }}\int_{Q}\widetilde{b}_{Q}d\sigma +\frac{1}{2},
\end{eqnarray*}%
if we choose $0<\delta <\frac{1}{4C_{\mathbf{b}}^{3}}$. This gives the lower
bound $\frac{1}{\left\vert Q\right\vert _{\sigma }}\int_{Q}\widetilde{b}%
_{Q}d\sigma \geq \frac{1}{2}$, and for an upper bound we have%
\begin{equation*}
\left\Vert \widetilde{b}_{Q}\right\Vert _{L^{\infty }\left( \sigma \right)
}\leq \left( 1+\sqrt{C_{b_{Q}}}\right) C_{b_{Q}}\ .
\end{equation*}%
Since we are taking $\delta <\frac{1}{4C_{\mathbf{b}}^{3}}$, we have $1+%
\sqrt{C_{b_{Q}}\delta }\leq 2$, and so we also have both 
\begin{eqnarray*}
\frac{\left\Vert \mathbf{1}_{Q_{\limfunc{left}}}\widetilde{b}_{Q}\right\Vert
_{L^{\infty }\left( \sigma \right) }}{\left\vert \frac{1}{\left\vert Q_{%
\limfunc{left}}\right\vert _{\sigma }}\int_{Q_{\limfunc{left}}}\widetilde{b}%
_{Q}d\sigma \right\vert } &\leq &\frac{\left( 1+\sqrt{C_{b_{Q}}\delta }%
\right) C_{b_{Q}}}{\left\vert \frac{1}{\left\vert Q_{\limfunc{left}%
}\right\vert _{\sigma }}\int_{Q_{\limfunc{left}}}\left[ p-n\left( 1+\sqrt{%
C_{b_{Q}}\delta }\right) \right] d\sigma \right\vert } \\
&\leq &\frac{\left( 1+\sqrt{C_{b_{Q}}\delta }\right) C_{b_{Q}}}{\left\vert 
\sqrt{C_{b_{Q}}\delta }\frac{1}{\left\vert Q_{\limfunc{left}}\right\vert
_{\sigma }}\int_{Q_{\limfunc{left}}}nd\sigma \right\vert }\leq \frac{%
4C_{b_{Q}}}{\sqrt{C_{b_{Q}}\delta }\frac{1}{\left\vert Q_{\limfunc{left}%
}\right\vert _{\sigma }}\int_{Q_{\limfunc{left}}}\left\vert b_{Q}\right\vert
d\sigma } \\
&\leq &\frac{4C_{b_{Q}}}{C_{b_{Q}}\delta }=\frac{4}{\delta }\text{\ }, \\
\frac{\left\Vert \mathbf{1}_{Q_{\limfunc{right}}}\widetilde{b}%
_{Q}\right\Vert _{L^{\infty }\left( \sigma \right) }}{\left\vert \frac{1}{%
\left\vert Q_{\limfunc{right}}\right\vert _{\sigma }}\int_{Q_{\limfunc{right}%
}}\widetilde{b}_{Q}d\sigma \right\vert } &=&\frac{\left\Vert \mathbf{1}_{Q_{%
\limfunc{right}}}b_{Q}\right\Vert _{L^{\infty }\left( \sigma \right) }}{%
\left\vert \frac{1}{\left\vert Q_{\limfunc{right}}\right\vert _{\sigma }}%
\int_{Q_{\limfunc{right}}}b_{Q}d\sigma \right\vert }<\frac{1}{\delta }%
C_{b_{Q}}\ ,
\end{eqnarray*}%
where the second line follows since (\ref{Q' big}) fails for $Q^{\prime }=Q_{%
\limfunc{right}}$.

Finally we check the testing condition in this case. We have from
Minkowski's inequality, 
\begin{eqnarray*}
\sqrt{\frac{1}{\left\vert Q\right\vert _{\sigma }}\int_{Q}\left\vert
T_{\sigma }^{\alpha }\widetilde{b}_{Q}\right\vert ^{2}d\omega } &\leq &\sqrt{%
\frac{1}{\left\vert Q\right\vert _{\sigma }}\int_{Q}\left\vert T_{\sigma
}^{\alpha }b_{Q}\right\vert ^{2}d\omega }+\sqrt{\frac{1}{\left\vert
Q\right\vert _{\sigma }}\int_{Q}\left\vert T_{\sigma }^{\alpha }\left( 
\widetilde{b}_{Q}-b_{Q}\right) \right\vert ^{2}d\omega } \\
&\leq &\mathfrak{T}_{T^{\alpha }}^{b_{Q}}\left( Q\right) +\mathfrak{N}%
_{T^{\alpha }}\left( Q\right) \sqrt{\frac{1}{\left\vert Q\right\vert
_{\sigma }}\int_{Q}\left\vert \widetilde{b}_{Q}-b_{Q}\right\vert ^{2}d\sigma 
}.
\end{eqnarray*}%
Now recall we are assuming $\int_{Q_{\limfunc{left}}}nd\sigma >\int_{Q_{%
\limfunc{left}}}pd\sigma $, so that%
\begin{eqnarray*}
\sqrt{\frac{1}{\left\vert Q\right\vert _{\sigma }}\int_{Q}\left\vert 
\widetilde{b}_{Q}-b_{Q}\right\vert ^{2}d\sigma } &=&\sqrt{\frac{1}{%
\left\vert Q\right\vert _{\sigma }}\int_{Q_{\limfunc{left}}}\left\vert \sqrt{%
C_{\mathbf{b}}\delta }n\right\vert ^{2}d\sigma }\leq \sqrt{C_{b_{Q}}\delta }%
\sqrt{\frac{1}{\left\vert Q\right\vert _{\sigma }}\int_{Q_{\limfunc{left}%
}}\left\vert n\right\vert ^{2}d\sigma } \\
&\lesssim &\sqrt{C_{b_{Q}}\delta }C_{b_{Q}}=C_{b_{Q}}^{\frac{3}{2}}\sqrt{%
\delta }.
\end{eqnarray*}%
Thus in any case we have%
\begin{equation*}
\sqrt{\frac{1}{\left\vert Q\right\vert _{\sigma }}\int_{Q}\left\vert
T_{\sigma }^{\alpha }\widetilde{b}_{Q}\right\vert ^{2}d\omega }\leq 
\mathfrak{T}_{T^{\alpha }}^{b_{Q}}\left( Q\right) +2C_{b_{Q}}^{\frac{3}{2}%
}\delta ^{\frac{1}{2}}\mathfrak{N}_{T^{\alpha }}\left( Q\right) \ .
\end{equation*}%
\textbf{Subcase (2b)} $\int_{Q_{\limfunc{left}}}pd\sigma \geq \int_{Q_{%
\limfunc{left}}}nd\sigma $: The same estimates arise in this case, except
that we get the better lower bound $\frac{1}{\left\vert Q\right\vert
_{\sigma }}\int_{Q}\widetilde{b}_{Q}d\sigma \geq 1$.

Collecting all of our estimates for $\widetilde{b}_{Q}$ in the various cases
above we have%
\begin{eqnarray*}
&&\frac{1}{2}\leq \frac{1}{\left\vert Q\right\vert _{\sigma }}\int_{Q}%
\widetilde{b}_{Q}d\sigma \leq \left\Vert \mathbf{1}_{Q}\widetilde{b}%
_{Q}\right\Vert _{L^{\infty }\left( \sigma \right) }\leq \left( 1+\sqrt{%
C_{b_{Q}}}\right) C_{b_{Q}}\ ,\ \ \ \ \ Q\in \mathcal{P}\ , \\
&&\frac{\left\Vert \mathbf{1}_{Q^{\prime }}\widetilde{b}_{Q}\right\Vert
_{L^{\infty }\left( \sigma \right) }}{\left\vert \frac{1}{\left\vert
Q^{\prime }\right\vert _{\sigma }}\int_{Q^{\prime }}\widetilde{b}_{Q}d\sigma
\right\vert },\frac{\left\Vert \mathbf{1}_{Q^{\prime }}\widetilde{b}%
_{Q}^{\ast }\right\Vert _{L^{\infty }\left( \omega \right) }}{\left\vert 
\frac{1}{\left\vert Q^{\prime }\right\vert _{\omega }}\int_{Q^{\prime }}%
\widetilde{b}_{Q}^{\ast }d\omega \right\vert }<\frac{4}{\delta }C_{b_{Q}}\
,\ \ \ \ \ Q\in \mathcal{P}\text{ and }Q^{\prime }\in \mathfrak{C}\left(
Q\right) \ , \\
&&\sqrt{\frac{1}{\left\vert Q\right\vert _{\sigma }}\int_{Q}\left\vert
T_{\sigma }^{\alpha }\widetilde{b}_{Q}\right\vert ^{2}d\omega }\leq 
\mathfrak{T}_{T^{\alpha }}^{b_{Q}}\left( Q\right) +2C_{b_{Q}}^{\frac{3}{2}%
}\delta ^{\frac{1}{4}}\mathfrak{N}_{T^{\alpha }}\left( Q\right) \ \ \ \ \
Q\in \mathcal{P}\ .
\end{eqnarray*}%
In order to obtain the inequalities for $\widetilde{b}_{Q}$ in the
conclusion of Lemma \ref{further red}, we simply multiply the above function 
$\widetilde{b}_{Q}$ by a factor of $2$.

Finally, if $\left\vert b_{Q}\right\vert \geq c_{1}>0$, we need only
consider \textbf{Case (2)} above, in which case we have $\left\vert
b_{Q}\right\vert \geq \left\vert \widetilde{b}_{Q}\right\vert $. This
completes the proof of Lemma \ref{further red}.
\end{proof}

\subsubsection{Control of averages in coronas}

Let $\mathcal{D}_{Q}$ be the grid of dyadic subintervals of $Q$. In the
construction of the triple corona below, we will need to repeat the
construction in the previous subsubsection for a subdecomposition $\left\{
Q_{i}\right\} _{i=1}^{\infty }$ of dyadic subintervals $Q_{i}\in \mathcal{D}%
_{Q}$ of an interval $Q$. Define the corona corresponding to the
subdecomposition $\left\{ Q_{i}\right\} _{i=1}^{\infty }$ by 
\begin{equation*}
\mathcal{C}_{Q}\equiv \mathcal{D}_{Q}\setminus \bigcup_{i=1}^{\infty }%
\mathcal{D}_{Q_{i}}\ .
\end{equation*}

\begin{lemma}
\label{prelim control of corona}Suppose that $\sigma $ and $\omega $ are
locally finite positive Borel measures on the real line $\mathbb{R}$. Assume
that $T^{\alpha }$ is a standard $\alpha $-fractional elliptic and gradient
elliptic singular integral operator on $\mathbb{R}$, and set $T_{\sigma
}^{\alpha }f=T^{\alpha }\left( f\sigma \right) $ for any smooth truncation
of $T_{\sigma }^{\alpha }$, so that $T_{\sigma }^{\alpha }$ is \emph{apriori}
bounded from $L^{2}\left( \sigma \right) $ to $L^{2}\left( \omega \right) $.
Let $Q\in \mathcal{P}$ and let $\mathfrak{N}_{T^{\alpha }}\left( Q\right) $
be the best constant in the local inequality%
\begin{equation*}
\sqrt{\int_{Q^{\prime }}\left\vert T_{\sigma }^{\alpha }\left( \mathbf{1}%
_{Q}f\right) \right\vert ^{2}d\omega }\leq \mathfrak{N}_{T^{\alpha }}\left(
Q\right) \sqrt{\int_{Q}\left\vert f\right\vert ^{2}d\sigma }\ ,\ \ \ \ \
f\in L^{2}\left( \mathbf{1}_{Q}\sigma \right) .
\end{equation*}%
Let $\left\{ Q_{i}\right\} _{i=1}^{\infty }\subset \mathcal{D}_{Q}$\ be a
collection of pairwise disjoint dyadic subintervals of $Q$. Suppose that $%
b_{Q}$ is a real-valued function supported in $Q$ such that 
\begin{eqnarray*}
&&1\leq \frac{1}{\left\vert Q^{\prime }\right\vert _{\sigma }}%
\int_{Q^{\prime }}b_{Q}d\sigma \leq \left\Vert \mathbf{1}_{Q^{\prime
}}b_{Q}\right\Vert _{L^{\infty }\left( \sigma \right) }\leq C_{\mathbf{b}}\
,\ \ \ \ \ Q^{\prime }\in \mathcal{C}_{Q}\ , \\
&&\sqrt{\int_{Q}\left\vert T_{\sigma }^{\alpha }b_{Q}\right\vert ^{2}d\omega 
}\leq \mathfrak{T}_{T^{\alpha }}^{b_{Q}}\left( Q\right) \sqrt{\left\vert
Q\right\vert _{\sigma }}\ .
\end{eqnarray*}%
Then for every $0<\delta <\frac{1}{4C_{\mathbf{b}}^{3}}$, there exists a
real-valued function $\widetilde{b}_{Q}$ supported in $Q$ such that%
\begin{eqnarray*}
&&1\leq \frac{1}{\left\vert Q^{\prime }\right\vert _{\sigma }}%
\int_{Q^{\prime }}\widetilde{b}_{Q}d\sigma \leq \left\Vert \mathbf{1}%
_{Q^{\prime }}\widetilde{b}_{Q}\right\Vert _{L^{\infty }\left( \sigma
\right) }\leq 2\left( 1+\sqrt{C_{\mathbf{b}}}\right) C_{\mathbf{b}}\ ,\ \ \
\ \ Q^{\prime }\in \mathcal{C}_{Q}\ , \\
&&\sqrt{\int_{Q}\left\vert T_{\sigma }^{\alpha }\widetilde{b}_{Q}\right\vert
^{2}d\omega }\leq \left[ 2\mathfrak{T}_{T^{\alpha }}^{b_{Q}}\left( Q\right)
+4C_{\mathbf{b}}^{\frac{3}{2}}\delta ^{\frac{1}{4}}\mathfrak{N}_{T^{\alpha
}}\left( Q\right) \right] \sqrt{\left\vert Q\right\vert _{\sigma }}\ , \\
&&0<\left\Vert \mathbf{1}_{Q_{i}}\widetilde{b}_{Q}\right\Vert _{L^{\infty
}\left( \sigma \right) }\leq \frac{16C_{\mathbf{b}}}{\delta }\left\vert 
\frac{1}{\left\vert Q_{i}\right\vert _{\sigma }}\int_{Q_{i}}\widetilde{b}%
_{Q}d\sigma \right\vert \ ,\ \ \ \ \ 1\leq i<\infty .
\end{eqnarray*}%
Moreover, if $\left\vert b_{Q}\right\vert \geq c_{1}>0$, then we may take $%
\left\vert \widetilde{b}_{Q}\right\vert \geq c_{1}$ as well.
\end{lemma}

The additional gain in the lemma is in the final line that controls the
degeneracy of $\widetilde{b}_{Q}$ at the `bottom' of the corona $\mathcal{C}%
_{Q}$ by establishing a reverse H\"{o}lder control. Note that if we combine
this control with the accretivity control in the corona $\mathcal{C}_{Q}$,
namely 
\begin{equation*}
\left\Vert \mathbf{1}_{Q^{\prime }}\widetilde{b}_{Q}\right\Vert _{L^{\infty
}\left( \sigma \right) }\leq 2\left( 1+\sqrt{C_{\mathbf{b}}}\right) C_{%
\mathbf{b}}\leq 2\left( 1+\sqrt{C_{\mathbf{b}}}\right) C_{\mathbf{b}}\frac{1%
}{\left\vert Q^{\prime }\right\vert _{\sigma }}\int_{Q^{\prime }}\widetilde{b%
}_{Q}d\sigma ,
\end{equation*}%
we obtain reverse H\"{o}lder control throughout the entire collection $%
\mathcal{C}_{Q}\cup \left\{ Q_{i}\right\} _{i=1}^{\infty }$: 
\begin{equation*}
\left\Vert \mathbf{1}_{I}\widetilde{b}_{Q^{\prime }}\right\Vert _{L^{\infty
}\left( \sigma \right) }\leq C_{\delta ,\mathbf{b}}\left\vert \frac{1}{%
\left\vert I\right\vert _{\sigma }}\int_{I}\widetilde{b}_{Q^{\prime
}}d\sigma \right\vert ,\ \ \ \ \ I\in \mathfrak{C}\left( Q^{\prime }\right)
,Q^{\prime }\in \mathcal{C}_{Q}\text{ }.
\end{equation*}%
This has the crucial consequence that the martingale and dual martingale
differences $\bigtriangleup _{Q^{\prime }}^{\sigma ,\mathbf{b}}$ and $%
\square _{Q^{\prime }}^{\sigma ,\mathbf{b}}$ associated with these functions
as defined in (\ref{def diff}) of Appendix A, satisfy%
\begin{equation}
\left\vert \bigtriangleup _{Q^{\prime }}^{\sigma ,\mathbf{b}}h\right\vert
,\left\vert \square _{Q^{\prime }}^{\sigma ,\mathbf{b}}h\right\vert \leq
C_{\delta ,\mathbf{b}}\sum_{I\in \mathfrak{C}\left( Q^{\prime }\right)
}\left( \frac{1}{\left\vert I\right\vert _{\sigma }}\int_{I}\left\vert
h\right\vert d\sigma +\frac{1}{\left\vert Q^{\prime }\right\vert _{\sigma }}%
\int_{Q^{\prime }}\left\vert h\right\vert d\sigma \right) \mathbf{1}_{I}\ .
\label{cruc conseq}
\end{equation}%
See Appendix A for more detail on this. However, the defect in this lemma is
that we lose the weak testing condition for $\widetilde{b}_{Q}$ in the
corona even if we had assumed it at the outset for $b_{Q}$.

\begin{proof}
The proof of Lemma \ref{prelim control of corona} is similar to that of the
special case given by Lemma \ref{further red}. Indeed, we define%
\begin{eqnarray*}
\widetilde{b}_{Q} &\equiv &\sum_{i\in G_{0}}\delta \mathbf{1}%
_{Q_{i}}+\sum_{i\in G_{+}}\left( \frac{1}{\left\vert Q_{i}\right\vert
_{\sigma }}\int_{Q_{i}}\left\vert b_{Q}\right\vert d\sigma \right) \mathbf{1}%
_{Q_{i}} \\
&&+\sum_{i\in B_{-}}\left( \frac{1}{\left\vert Q_{i}\right\vert _{\sigma }}%
\int_{Q_{i}}\left[ p_{i}-n_{i}\left( 1+\sqrt{C_{\mathbf{b}}\delta }\right) %
\right] d\sigma \right) \mathbf{1}_{Q_{i}} \\
&&+\sum_{i\in B_{+}}\left( \frac{1}{\left\vert Q_{i}\right\vert _{\sigma }}%
\int_{Q_{i}}\left[ \left( 1+\sqrt{C_{\mathbf{b}}\delta }\right) p_{i}-n_{i}%
\right] d\sigma \right) \mathbf{1}_{Q_{i}} \\
&&+b_{Q}\mathbf{1}_{Q\setminus \cup _{i=1}^{\infty }Q_{i}}\ ,
\end{eqnarray*}%
where%
\begin{eqnarray*}
G_{0} &\equiv &\left\{ i:\frac{1}{\left\vert Q_{i}\right\vert _{\sigma }}%
\int_{Q_{i}}\left\vert b_{Q}\right\vert d\sigma =0\right\} \\
G_{+} &\equiv &\left\{ i:0<\frac{1}{\left\vert Q_{i}\right\vert _{\sigma }}%
\int_{Q_{i}}\left\vert b_{Q}\right\vert d\sigma \leq \sqrt{C_{\mathbf{b}%
}\delta }\right\} , \\
B_{-} &\equiv &\left\{ i:\frac{1}{\left\vert Q_{i}\right\vert _{\sigma }}%
\int_{Q_{i}}\left\vert b_{Q}\right\vert d\sigma >\sqrt{C_{\mathbf{b}}\delta }%
\text{ and }\int_{Q_{i}}nd\sigma >\int_{Q_{i}}pd\sigma \right\} , \\
B_{+} &\equiv &\left\{ i:\frac{1}{\left\vert Q_{i}\right\vert _{\sigma }}%
\int_{Q_{i}}\left\vert b_{Q}\right\vert d\sigma >\sqrt{C_{\mathbf{b}}\delta }%
\text{ and }\int_{Q_{i}}pd\sigma \geq \int_{Q_{i}}nd\sigma \right\} .
\end{eqnarray*}%
First we note that%
\begin{eqnarray*}
1 &\leq &\frac{1}{\left\vert Q^{\prime }\right\vert _{\sigma }}%
\int_{Q^{\prime }}b_{Q}d\sigma \leq \frac{1}{\left\vert Q^{\prime
}\right\vert _{\sigma }}\int_{Q^{\prime }}\widetilde{b}_{Q}d\sigma +\frac{1}{%
\left\vert Q^{\prime }\right\vert _{\sigma }}\int_{Q^{\prime }}\left( b_{Q}-%
\widetilde{b}_{Q}\right) d\sigma \\
&\leq &\frac{1}{\left\vert Q^{\prime }\right\vert _{\sigma }}\int_{Q^{\prime
}}\widetilde{b}_{Q}d\sigma +\sum_{i:\ Q_{i}\subset Q^{\prime }}\frac{1}{%
\left\vert Q^{\prime }\right\vert _{\sigma }}\left\{ 
\begin{array}{ccc}
\int_{Q_{i}}\left\vert b_{Q}\right\vert d\sigma -\delta & \text{ in } & 
\text{\textbf{Case (0)}} \\ 
2\int_{Q_{i}}\left\vert b_{Q}\right\vert d\sigma & \text{ in } & \text{%
\textbf{Case (1)}} \\ 
\sqrt{C_{\mathbf{b}}\delta }\int_{Q_{i}}\left\vert b_{Q}\right\vert d\sigma
& \text{ in } & \text{\textbf{Case (2)}}%
\end{array}%
\right. \\
&\leq &\frac{1}{\left\vert Q^{\prime }\right\vert _{\sigma }}\int_{Q^{\prime
}}\widetilde{b}_{Q}d\sigma +2\sqrt{C_{\mathbf{b}}\delta }C_{\mathbf{b}%
}\sum_{i:\ Q_{i}\subset Q^{\prime }}\frac{\left\vert Q_{i}\right\vert
_{\sigma }}{\left\vert Q^{\prime }\right\vert _{\sigma }}\leq \frac{1}{%
\left\vert Q^{\prime }\right\vert _{\sigma }}\int_{Q^{\prime }}\widetilde{b}%
_{Q}d\sigma +\frac{1}{2},
\end{eqnarray*}%
if $0<\delta <\frac{1}{4C_{\mathbf{b}}^{3}}$.

Then we estimate the testing condition for the interval $Q$ by%
\begin{eqnarray*}
\sqrt{\frac{1}{\left\vert Q\right\vert _{\sigma }}\int_{Q}\left\vert
T_{\sigma }^{\alpha }\widetilde{b}_{Q}\right\vert ^{2}d\omega } &\leq &\sqrt{%
\frac{1}{\left\vert Q\right\vert _{\sigma }}\int_{Q}\left\vert T_{\sigma
}^{\alpha }b_{Q}\right\vert ^{2}d\omega }+\sqrt{\frac{1}{\left\vert
Q\right\vert _{\sigma }}\int_{Q}\left\vert T_{\sigma }^{\alpha }\left( 
\widetilde{b}_{Q}-b_{Q}\right) \right\vert ^{2}d\omega } \\
&\leq &\mathfrak{T}_{T^{\alpha }}^{b_{Q}}\left( Q\right) +\mathfrak{N}%
_{T^{\alpha }}\left( Q\right) \sqrt{\frac{1}{\left\vert Q\right\vert
_{\sigma }}\int_{Q}\left\vert \widetilde{b}_{Q}-b_{Q}\right\vert ^{2}d\sigma 
},
\end{eqnarray*}%
and note that the arguments above show that 
\begin{eqnarray*}
\widetilde{b}_{Q}-b_{Q} &=&\sum_{i\in G}\left( \frac{1}{\left\vert
Q_{i}\right\vert _{\sigma }}\int_{Q_{i}}\left\vert b_{Q}\right\vert d\sigma
\right) \mathbf{1}_{Q_{i}}-\sum_{i\in G}b_{Q}\mathbf{1}_{Q_{i}} \\
&&-\sum_{i\in B_{-}}\left( \frac{\sqrt{C_{\mathbf{b}}\delta }}{\left\vert
Q_{i}\right\vert _{\sigma }}\int_{Q_{i}}n_{i}d\sigma \right) \mathbf{1}%
_{Q_{i}}+\sum_{i\in B_{+}}\left( \frac{\sqrt{C_{\mathbf{b}}\delta }}{%
\left\vert Q_{i}\right\vert _{\sigma }}\int_{Q_{i}}p_{i}d\sigma \right) 
\mathbf{1}_{Q_{i}},
\end{eqnarray*}%
satisfies an inequality of the form%
\begin{equation*}
\int_{Q}\left\vert \widetilde{b}_{Q}-b_{Q}\right\vert ^{2}d\sigma \leq
C\left( C_{\mathbf{b}}\right) \delta ^{\frac{1}{4}}\sum_{i=1}^{\infty
}\left\vert Q_{i}\right\vert _{\sigma }\leq C\left( C_{\mathbf{b}}\right)
\delta ^{\frac{1}{4}}\left\vert Q\right\vert _{\sigma }\ .
\end{equation*}
\end{proof}

\begin{remark}
The estimate $\int_{Q}\left\vert \widetilde{b}_{Q}-b_{Q}\right\vert
^{2}d\sigma \leq C\left( C_{\mathbf{b}}\right) \delta ^{\frac{1}{4}%
}\sum_{i=1}^{\infty }\left\vert Q_{i}\right\vert _{\sigma }$ in the last
line of the above proof is of course too large in general to be dominated by
a fixed multiple of $\left\vert Q^{\prime }\right\vert _{\sigma }$ for $%
Q^{\prime }\in \mathcal{C}_{Q}$, and this is the reason we have no control
of weak testing for $\widetilde{b}_{Q}$ in the rest of the corona even if we
assume weak testing for $b_{Q}$ in the corona $\mathcal{C}_{Q}$. This defect
is addressed in the next subsection below.
\end{remark}

\subsection{Three corona decompositions}

We will use multiple corona constructions, namely a Calder\'{o}n-Zygmund
decomposition, an accretive decomposition, a weak testing decomposition, and
an energy decomposition, in order to reduce matters to the stopping form,
which is treated in Section \ref{Sec stop} by adapting the bottom/up
stopping time in the argument of M. Lacey in \cite{Lac}, and using an
additional `indented' top/down corona construction, in order to accommodate
weak goodness. We will then iterate these corona decompositions into a
single corona decomposition, which we refer to as the \emph{triple corona}.
More precisely, we iterate the first generation of common stopping times
with an infusion of the reverse H\"{o}lder condition on children, followed
by another iteration of the first generation of weak testing stopping times.
Recall that we must show the bilinear inequality%
\begin{equation*}
\left\vert \int \left( T_{\sigma }^{\alpha }f\right) gd\omega \right\vert
\leq \mathfrak{N}_{T^{\alpha }}\left\Vert f\right\Vert _{L^{2}\left( \sigma
\right) }\left\Vert g\right\Vert _{L^{2}\left( \omega \right) },\ \ \ \ \
f\in L^{2}\left( \sigma \right) \text{ and }g\in L^{2}\left( \omega \right) .
\end{equation*}

\subsubsection{The Calder\'{o}n-Zygmund corona decomposition}

We first introduce the Calder\'{o}n-Zygmund stopping times $\mathcal{F}$ for
a function $\phi \in L^{2}\left( \mu \right) $ relative to an interval $%
S_{0} $ and a positive constant $C_{0}\geq 4$. Let $\mathcal{F}=\left\{
F\right\} _{F\in \mathcal{F}}$ be the collection of Calder\'{o}n-Zygmund
stopping intervals for $\phi $ defined so that $F\subset S_{0}$, $S_{0}\in 
\mathcal{F} $, and for all $F\in \mathcal{F}$ with $F\subsetneqq S_{0}$ we
have%
\begin{eqnarray*}
\frac{1}{\left\vert F\right\vert _{\mu }}\int_{F}\left\vert \phi \right\vert
d\mu &>&C_{0}\frac{1}{\left\vert \pi _{\mathcal{F}}F\right\vert _{\mu }}%
\int_{\pi _{\mathcal{F}}F}\left\vert \phi \right\vert d\mu ; \\
\frac{1}{\left\vert F^{\prime }\right\vert _{\mu }}\int_{F^{\prime
}}\left\vert \phi \right\vert d\mu &\leq &C_{0}\frac{1}{\left\vert \pi _{%
\mathcal{F}}F\right\vert _{\mu }}\int_{\pi _{\mathcal{F}}F}\left\vert \phi
\right\vert d\mu \ \ \ \ \text{ for }F\subsetneqq F^{\prime }\subset \pi _{%
\mathcal{F}}F.
\end{eqnarray*}%
To achieve this construction we use the following definition.

\begin{definition}
\label{CZ stopping times}Let $C_{0}\geq 4$. Given a dyadic grid $\mathcal{D}$
an interval $S_{0}\in \mathcal{D}$, define $\mathcal{S}\left( S_{0}\right) $
to be the \emph{maximal} $\mathcal{D}$-subintervals $I\subset S_{0}$ such
that%
\begin{equation*}
\frac{1}{\left\vert I\right\vert _{\mu }}\int_{I}\left\vert \phi \right\vert
d\mu >C_{0}\frac{1}{\left\vert S_{0}\right\vert _{\mu }}\int_{S_{0}}\left%
\vert \phi \right\vert d\mu \ ,
\end{equation*}%
and then define the $CZ$ stopping intervals of $S_{0}$ to be the collection 
\begin{equation*}
\mathcal{S}=\left\{ S_{0}\right\} \cup \dbigcup\limits_{n=0}^{\infty }%
\mathcal{S}_{n}
\end{equation*}%
where $\mathcal{S}_{0}=\mathcal{S}\left( S_{0}\right) $ and $\mathcal{S}%
_{n+1}=\dbigcup\limits_{S\in \mathcal{S}_{n}}\mathcal{S}\left( S\right) $
for $n\geq 0$.
\end{definition}

Let $\mathcal{D}=\dbigcup\limits_{F\in \mathcal{F}}\mathcal{C}_{F}$ be the
associated corona decomposition of the dyadic grid $\mathcal{D}$ where%
\begin{equation*}
\mathcal{C}_{F}\equiv \left\{ F^{\prime }\in \mathcal{D}:F\supset F^{\prime
}\supsetneqq H\text{ for some }H\in \mathfrak{C}_{\mathcal{F}}\left(
F\right) \right\} .
\end{equation*}%
We now recall some of the definitions just used above. See \cite{SaShUr7}
and/or \cite{SaShUr6} for more detail. For an interval $I\in \mathcal{D}$
let $\pi _{\mathcal{D}}I$ be the $\mathcal{D}$-parent of $I$ in the grid $%
\mathcal{D}$, and let $\pi _{\mathcal{F}}I$ be the smallest member of $%
\mathcal{F}$ that \emph{strictly} contains $I$. For $F,F^{\prime }\in 
\mathcal{F}$, we say that $F^{\prime }$ is an $\mathcal{F}$-child of $F$ if $%
\pi _{\mathcal{F}}\left( F^{\prime }\right) =F$ (it could be that $F=\pi _{%
\mathcal{D}}F^{\prime }$), and we denote by $\mathfrak{C}_{\mathcal{F}%
}\left( F\right) $ the set of $\mathcal{F}$-children of $F$. We call $\pi _{%
\mathcal{F}}\left( F^{\prime }\right) $ the $\mathcal{F}$-parent of $%
F^{\prime }\in \mathcal{F}$.

The stopping intervals $\mathcal{F}$ above satisfy a Carleson condition:%
\begin{equation}
\sum_{F\in \mathcal{F}:\ F\subset \Omega }\left\vert F\right\vert _{\mu
}\leq C\left\vert \Omega \right\vert _{\mu }\ ,\ \ \ \ \ \text{for all open
sets }\Omega .  \label{CZ Car}
\end{equation}%
Indeed, 
\begin{equation*}
\sum_{F^{\prime }\in \mathfrak{C}_{\mathcal{F}}\left( F\right) }\left\vert
F^{\prime }\right\vert _{\mu }\leq \sum_{F^{\prime }\in \mathfrak{C}_{%
\mathcal{F}}\left( F\right) }\frac{\int_{F^{\prime }}\left\vert \phi
\right\vert d\mu }{C_{0}\frac{1}{\left\vert F\right\vert _{\mu }}%
\int_{F}\left\vert \phi \right\vert d\mu }\leq \frac{1}{C_{0}},
\end{equation*}%
and standard arguments now complete the proof of the Carleson condition.

We emphasize that accretive functions $b$ play no role in the Calder\'{o}%
n-Zygmund corona decomposition.

\subsubsection{The $\mathbf{b}$-accretive / weak testing corona decomposition%
}

Recall that we are assuming $p=\infty $, and that our testing functions $%
\mathbf{b}$ and $\mathbf{b}^{\ast }$ are real-valued, in the proof of
Theorem \ref{dim one}. We use a corona construction modelled after that of
Hyt\"{o}nen and Martikainen \cite{HyMa}, that delivers a \emph{weak corona
testing condition} that coincides with the testing condition itself \textbf{%
only} at the tops of the coronas. This corona decomposition is developed to
optimize the choice of a new family of testing functions $\left\{ \widehat{b}%
_{Q}\right\} _{Q\in \mathcal{D}}$ taken from the vector $\mathbf{b}\equiv
\left\{ b_{Q}\right\} _{Q\in \mathcal{D}}$ so that we have

\begin{enumerate}
\item the telescoping property at our disposal in each accretive corona,

\item a weak corona testing condition remains in force for the new testing
functions $\widehat{b}_{Q}$ that coincides with the usual testing condition
at the tops of the coronas,

\item the tops of the coronas, i.e. the stopping intervals, enjoy a Carleson
condition.
\end{enumerate}

We will sometimes refer to the old family as the \emph{original} family, and
denote it by $\left\{ b_{Q}^{\limfunc{orig}}\right\} _{Q\in \mathcal{D}}$.
The original family will reappear later in helping to estimate the nearby
form in Section \ref{Sec nearby}.

Let $\sigma $ and $\omega $ be locally finite Borel measures on $\mathbb{R}$%
. We assume that the vector of `testing functions' $\mathbf{b}\equiv \left\{
b_{Q}\right\} _{Q\in \mathcal{D}}$ is an $\infty $-strongly $\sigma $%
-accretive real-valued family, i.e. 
\begin{equation*}
\limfunc{support}b_{Q}\subset Q\ ,\ \ \ \ \ Q\in \mathcal{D},
\end{equation*}%
and%
\begin{equation*}
1\leq \frac{1}{\left\vert Q\right\vert _{\mu }}\int_{Q}b_{Q}d\sigma \leq
\left\Vert b_{Q}\right\Vert _{L^{\infty }\left( \sigma \right) }\leq C_{%
\mathbf{b}}<\infty ,\ \ \ \ \ Q\in \mathcal{D}\ ,
\end{equation*}%
and also that $\mathbf{b}^{\ast }\equiv \left\{ b_{Q}\right\} _{Q\in 
\mathcal{D}}$ is an $\infty $-strongly $\omega $-accretive real-valued
family, and we assume in addition the testing conditions%
\begin{eqnarray*}
\int_{Q}\left\vert T_{\sigma }^{\alpha }\left( \mathbf{1}_{Q}b_{Q}\right)
\right\vert ^{2}d\omega &\leq &\mathfrak{T}_{T^{\alpha }}^{\mathbf{b}%
}\left\vert Q\right\vert _{\sigma }\ ,\ \ \ \ \ \text{for all intervals }Q,
\\
\int_{Q}\left\vert T_{\omega }^{\alpha ,\ast }\left( \mathbf{1}%
_{Q}b_{Q}^{\ast }\right) \right\vert ^{2}d\sigma &\leq &\mathfrak{T}%
_{T^{\alpha }}^{\mathbf{b}^{\ast },\ast }\left\vert Q\right\vert _{\omega }\
,\ \ \ \ \ \text{for all intervals }Q.
\end{eqnarray*}

\begin{definition}
\label{accretive stopping times gen}Given a dyadic grid $\mathcal{D}$ an
interval $S_{0}\in \mathcal{D}$, define $\mathcal{S}\left( S_{0}\right) $ to
be the \emph{maximal} $\mathcal{D}$-subintervals $I\subset S_{0}$ such that%
\begin{eqnarray*}
\text{\textbf{either} }\left\vert \frac{1}{\left\vert I\right\vert _{\sigma }%
}\int_{I}b_{S_{0}}d\sigma \right\vert &<&\gamma \ , \\
\ \text{\textbf{or} }\int_{I}\left\vert T_{\sigma }^{\alpha }\left(
b_{S_{0}}\right) \right\vert ^{2}d\omega &>&\Gamma \left( \mathfrak{T}%
_{T^{\alpha }}^{\mathbf{b}}\right) ^{2}\left\vert I\right\vert _{\sigma }%
\text{ },
\end{eqnarray*}%
where the positive constants $\gamma ,\Gamma $ satisfy $0<\gamma <1<\Gamma
<\infty $. Then define the $\mathbf{b}$-accretive stopping intervals of $%
S_{0}$ to be the collection 
\begin{equation*}
\mathcal{S}=\left\{ S_{0}\right\} \cup \dbigcup\limits_{n=0}^{\infty }%
\mathcal{S}_{n}
\end{equation*}%
where $\mathcal{S}_{0}=\mathcal{S}\left( S_{0}\right) $ and $\mathcal{S}%
_{n+1}=\dbigcup\limits_{S\in \mathcal{S}_{n}}\mathcal{S}\left( S\right) $
for $n\geq 0$.
\end{definition}

For $\gamma <1$ chosen small enough and $\Gamma >1$ chosen large enough, the 
$\mathbf{b}$-accretive stopping intervals satisfy a $\sigma $-Carleson
condition relative to the measure $\sigma $, and the corresponding stopping
functions $b_{S_{0}}$ satisfy \emph{weak} testing inequalities in the
corona. The following lemma is essentially in \cite{HyMa}, but we include a
proof for completeness.

\begin{lemma}[\protect\cite{HyMa}]
\label{Car and Test gen}For $\gamma <1$ small enough and $\Gamma >1$ large
enough we have the following:

\begin{enumerate}
\item For every open set $\Omega $ we have we have the inequality,%
\begin{equation}
\sum_{S\in \mathcal{S}:\ S\subset \Omega }\left\vert S\right\vert _{\sigma
}\leq C\left( \mathfrak{T}_{T^{\alpha }}^{\mathbf{b}}\right) ^{2}\left\vert
\Omega \right\vert _{\sigma }\ .  \label{Car gen}
\end{equation}

\item For every interval $S\in \mathcal{C}_{S_{0}}$ we have the weak corona
testing inequality,%
\begin{equation}
\int_{S}\left\vert T_{\sigma }^{\alpha }b_{S_{0}}\right\vert ^{2}d\omega
\leq C\left( \mathfrak{T}_{T^{\alpha }}^{\mathbf{b}}\right) ^{2}\left\vert
S\right\vert _{\sigma }\ .  \label{Test gen}
\end{equation}
\end{enumerate}
\end{lemma}

\begin{proof}
We first address the Carleson condition (\ref{Car gen}). A standard argument
reduces matters to the case where $\Omega $ is an interval $Q\in \mathcal{S}$
with $\left\vert Q\right\vert _{\sigma }>0$. It suffices to consider each of
the two stopping criteria separately. We first address the stopping
condition $\left\vert \frac{1}{\left\vert I\right\vert _{\sigma }}%
\int_{I}b_{S_{0}}d\sigma \right\vert <\gamma $. Throughout this proof we
will denote the union of these children $\mathcal{S}\left( Q\right) $ of $Q$
by $E\left( Q\right) \equiv \dbigcup\limits_{S\in \mathcal{S}\left( Q\right)
}S$. Then we have%
\begin{equation*}
\left\vert \int_{E\left( Q\right) }b_{Q}d\sigma \right\vert \leq \sum_{S\in 
\mathcal{S}\left( Q\right) }\left\vert \int_{S}b_{Q}d\sigma \right\vert
<\gamma \sum_{S\in \mathcal{S}\left( Q\right) }\left\vert S\right\vert
_{\sigma }\leq \gamma \left\vert Q\right\vert _{\sigma }\ ,
\end{equation*}%
which together with our hypotheses on $b_{Q}$ gives%
\begin{eqnarray*}
\left\vert Q\right\vert _{\sigma } &<&\left\vert \int_{Q}b_{Q}d\sigma
\right\vert =\left\vert \int_{E\left( Q\right) }b_{Q}d\sigma \right\vert
+\left\vert \int_{Q\setminus E\left( Q\right) }b_{Q}d\sigma \right\vert \\
&\leq &\gamma \left\vert Q\right\vert _{\sigma }+\sqrt{\int_{Q\setminus
E\left( Q\right) }\left\vert b_{Q}\right\vert ^{2}d\sigma }\sqrt{\left\vert
Q\setminus E\left( Q\right) \right\vert _{\sigma }} \\
&\leq &\gamma \left\vert Q\right\vert _{\sigma }+\Gamma C_{\mathbf{b}}\sqrt{%
\left\vert Q\right\vert _{\sigma }}\sqrt{\left\vert Q\setminus E\left(
Q\right) \right\vert _{\sigma }}.
\end{eqnarray*}%
Rearranging the inequality yields successively%
\begin{eqnarray*}
\left( 1-\gamma \right) \left\vert Q\right\vert _{\sigma } &\leq &\Gamma C_{%
\mathbf{b}}\sqrt{\left\vert Q\right\vert _{\sigma }}\sqrt{\left\vert
Q\setminus E\left( Q\right) \right\vert _{\sigma }}; \\
\left( 1-\gamma \right) ^{2}\left\vert Q\right\vert _{\sigma }^{2} &\leq
&\Gamma ^{2}C_{\mathbf{b}}^{2}\left\vert Q\right\vert _{\sigma }\left\vert
Q\setminus E\left( Q\right) \right\vert _{\sigma }; \\
\frac{\left( 1-\gamma \right) ^{2}}{\Gamma ^{2}C_{\mathbf{b}}^{2}}\left\vert
Q\right\vert _{\sigma } &\leq &\left\vert Q\setminus E\left( Q\right)
\right\vert _{\sigma }\ ,
\end{eqnarray*}%
which in turn gives%
\begin{eqnarray*}
\sum_{S\in \mathcal{S}\left( Q\right) }\left\vert S\right\vert _{\sigma }
&=&\left\vert Q\right\vert _{\sigma }-\left\vert Q\setminus E\left( Q\right)
\right\vert _{\sigma } \\
&\leq &\left\vert Q\right\vert _{\sigma }-\frac{\left( 1-\gamma \right) ^{2}%
}{\Gamma ^{2}C_{\mathbf{b}}^{2}}\left\vert Q\right\vert _{\sigma }=\left( 1-%
\frac{\left( 1-\gamma \right) ^{2}}{\Gamma ^{2}C_{\mathbf{b}}^{2}}\right)
\left\vert Q\right\vert _{\sigma }\equiv \beta \left\vert Q\right\vert
_{\sigma }\ ,
\end{eqnarray*}%
where $0<\beta <1$ since $1\leq C_{\mathbf{b}}$. If we now iterate this
inequality, we obtain for each $k\geq 1$,%
\begin{eqnarray*}
\sum_{\substack{ S\in \mathcal{S}:\ S\subset Q  \\ \pi _{\mathcal{S}%
}^{\left( k\right) }\left( S\right) =Q}}\left\vert S\right\vert _{\sigma }
&=&\sum _{\substack{ S\in \mathcal{S}:\ S\subset Q  \\ \pi _{\mathcal{S}%
}^{\left( k-1\right) }\left( S\right) =Q}}\sum_{S^{\prime }\in \mathcal{S}%
\left( S\right) }\left\vert S^{\prime }\right\vert _{\sigma }\leq \sum 
_{\substack{ S\in \mathcal{S}:\ S\subset Q  \\ \pi _{\mathcal{S}}^{\left(
k-1\right) }\left( S\right) =Q}}\beta \left\vert S\right\vert _{\sigma } \\
&&\vdots \\
&\leq &\sum_{\substack{ S\in \mathcal{S}:\ S\subset Q  \\ \pi _{\mathcal{S}%
}^{\left( 1\right) }\left( S\right) =Q}}\left( 1-\gamma ^{2}\right)
^{k-1}\left\vert S\right\vert _{\sigma }\leq \beta ^{k}\left\vert
Q\right\vert _{\sigma }\ .
\end{eqnarray*}%
Finally then%
\begin{equation*}
\sum_{S\in \mathcal{S}:\ S\subset Q}\left\vert S\right\vert _{\sigma }\leq
\left\vert Q\right\vert _{\sigma }+\sum_{k=1}^{\infty }\sum_{\substack{ S\in 
\mathcal{S}:\ S\subset Q  \\ \pi _{\mathcal{S}}^{\left( k\right) }\left(
S\right) =Q}}\left\vert S\right\vert _{\sigma }\leq \sum_{k=0}^{\infty
}\beta ^{k}\left\vert Q\right\vert _{\sigma }=\frac{1}{1-\beta }\left\vert
Q\right\vert _{\sigma }=\frac{\Gamma ^{2}C_{\mathbf{b}}^{2}}{\left( 1-\gamma
\right) ^{2}}\left\vert Q\right\vert _{\sigma }\ .
\end{equation*}

Now we turn to the second stopping criterion 
\begin{equation*}
\int_{I}\left\vert T_{\sigma }^{\alpha }\left( b_{S_{0}}\right) \right\vert
^{2}d\omega >\Gamma \left( \mathfrak{T}_{T^{\alpha }}^{\mathbf{b}}\right)
^{2}\left\vert I\right\vert _{\sigma }.
\end{equation*}%
We have%
\begin{eqnarray*}
\sum_{S\in \mathfrak{C}_{\mathcal{S}}\left( S_{0}\right) }\left\vert
S\right\vert _{\sigma } &\leq &\frac{1}{\Gamma \left( \mathfrak{T}%
_{T^{\alpha }}^{\mathbf{b}}\right) ^{2}}\sum_{S\in \mathfrak{C}_{\mathcal{S}%
}\left( S_{0}\right) }\int_{S}\left\vert T_{\sigma }^{\alpha }\left(
b_{S_{0}}\right) \right\vert ^{2}d\omega \\
&\leq &\frac{1}{\Gamma \left( \mathfrak{T}_{T^{\alpha }}^{\mathbf{b}}\right)
^{2}}\sum_{S_{0}}\int_{S}\left\vert T_{\sigma }^{\alpha }\left(
b_{S_{0}}\right) \right\vert ^{2}d\omega \leq \frac{1}{\Gamma }\left\vert
S_{0}\right\vert _{\sigma }.
\end{eqnarray*}%
Iterating this inequality gives%
\begin{equation*}
\sum_{\substack{ S\in \mathcal{S}  \\ S\subset S_{0}}}\left\vert
S\right\vert _{\sigma }\leq \sum_{k=0}^{\infty }\frac{1}{\Gamma ^{k}}%
\left\vert S_{0}\right\vert _{\sigma }=\frac{\Gamma }{\Gamma -1}\left\vert
S_{0}\right\vert _{\sigma },
\end{equation*}%
and then%
\begin{equation*}
\sum_{\substack{ S\in \mathcal{S}  \\ S\subset \Omega }}\left\vert
S\right\vert _{\sigma }=\sum_{\substack{ \text{maximal }S_{0}\in \mathcal{S} 
\\ S_{0}\subset \Omega }}\sum_{\substack{ S\in \mathcal{S}  \\ S\subset
S_{0} }}\left\vert S\right\vert _{\sigma }\leq \frac{\Gamma }{\Gamma -1}\sum 
_{\substack{ \text{maximal }S_{0}\in \mathcal{S}  \\ S_{0}\subset \Omega }}%
\left\vert S_{0}\right\vert _{\sigma }=\frac{\Gamma }{\Gamma -1}\left\vert
\Omega \right\vert _{\sigma }\ .
\end{equation*}

Finally, for $I\in \mathcal{C}_{S_{0}}$ we have the weak testing inequality%
\begin{equation*}
\int_{I}\left\vert T_{\sigma }^{\alpha }\left( b_{S_{0}}\right) \right\vert
^{2}d\omega \leq \Gamma \left( \mathfrak{T}_{T^{\alpha }}^{\mathbf{b}%
}\right) ^{2}\left\vert I\right\vert _{\sigma }\ ,
\end{equation*}

and this completes the proof of Lemma \ref{Car and Test gen}.
\end{proof}

\subsubsection{The energy corona decompositions}

Given a weight pair $\left( \sigma ,\omega \right) $, we now construct an
energy corona decomposition for $\sigma $, and an energy corona
decomposition for $\omega $, that uniformize estimates (c.f. \cite{NTV3}, 
\cite{LaSaUr2}, \cite{SaShUr6} and \cite{SaShUr7}). In order to define these
constructions, we recall that the energy condition constant $\mathcal{E}%
_{2}^{\alpha }$ in Definition \ref{def strong quasienergy} is given by%
\begin{equation*}
\left( \mathcal{E}_{2}^{\alpha }\right) ^{2}\equiv \sup_{\substack{ Q\in 
\mathcal{P}  \\ Q\supset \dot{\cup}J_{r}}}\frac{1}{\left\vert Q\right\vert
_{\sigma }}\sum_{r=1}^{\infty }\left( \frac{\mathrm{P}^{\alpha }\left( J_{r},%
\mathbf{1}_{Q}\sigma \right) }{\left\vert J_{r}\right\vert }\right)
^{2}\left\Vert x-m_{J_{r}}\right\Vert _{L^{2}\left( \mathbf{1}_{J_{r}}\omega
\right) }^{2}\ ,
\end{equation*}%
where $\dot{\cup}J_{r}$ is an arbitary subdecomposition of $Q$ into
intervals $J_{r}\in \mathcal{P}$. In the next definition we restrict the
intervals $Q$ to a dyadic grid $\mathcal{D}$, but keep the subintervals $%
J_{r}$ unrestricted.

\begin{definition}
\label{def energy corona 3}Given a dyadic grid $\mathcal{D}$ and an interval 
$S_{0}\in \mathcal{D}$, define $\mathcal{S}\left( S_{0}\right) $ to be the 
\emph{maximal} $\mathcal{D}$-subintervals $I\subset S_{0}$ such that%
\begin{equation}
\sup_{I\supset \dot{\cup}J_{r}}\sum_{r=1}^{\infty }\left( \frac{\mathrm{P}%
^{\alpha }\left( J_{r},\mathbf{1}_{S_{0}}\sigma \right) }{\left\vert
J_{r}\right\vert }\right) ^{2}\left\Vert x-m_{J_{r}}\right\Vert
_{L^{2}\left( \mathbf{1}_{J_{r}}\omega \right) }^{2}\geq C_{\limfunc{energy}}%
\left[ \left( \mathfrak{E}_{2}^{\alpha }\right) ^{2}+\mathfrak{A}%
_{2}^{\alpha }\right] \ \left\vert I\right\vert _{\sigma },
\label{def stop 3}
\end{equation}%
where the intervals $J_{r}\in \mathcal{P}$ are pairwise disjoint in $I$, $%
\mathcal{E}_{2}^{\alpha }$ is the energy condition constant, and $C_{%
\limfunc{energy}}$ is a sufficiently large positive constant depending only
on $\alpha $. Then define the $\sigma $-energy stopping intervals of $S_{0}$
to be the collection 
\begin{equation*}
\mathcal{S}=\left\{ S_{0}\right\} \cup \dbigcup\limits_{n=0}^{\infty }%
\mathcal{S}_{n}
\end{equation*}%
where $\mathcal{S}_{0}=\mathcal{S}\left( S_{0}\right) $ and $\mathcal{S}%
_{n+1}=\dbigcup\limits_{S\in \mathcal{S}_{n}}\mathcal{S}\left( S\right) $
for $n\geq 0$.
\end{definition}

We now claim that from the energy condition $\mathcal{E}_{2}^{\alpha }%
\mathcal{<\infty }$, we obtain the $\sigma $-Carleson estimate,%
\begin{equation}
\sum_{S\in \mathcal{S}:\ S\subset I}\left\vert S\right\vert _{\sigma }\leq
2\left\vert I\right\vert _{\sigma },\ \ \ \ \ I\in \mathcal{D}^{\sigma }.
\label{sigma Carleson 3}
\end{equation}%
Indeed, for any $S_{1}\in \mathcal{S}$ we have%
\begin{eqnarray*}
\sum_{S\in \mathfrak{C}_{\mathcal{S}}\left( S_{1}\right) }\left\vert
S\right\vert _{\sigma } &\leq &\frac{1}{C_{\limfunc{energy}}\left( \mathcal{E%
}_{2}^{\alpha }\right) ^{2}}\sum_{S\in \mathfrak{C}_{\mathcal{S}}\left(
S_{1}\right) }\sup_{S\supset \dot{\cup}J_{r}}\sum_{r=1}^{\infty }\left( 
\frac{\mathrm{P}^{\alpha }\left( J_{r},\mathbf{1}_{S_{1}}\sigma \right) }{%
\left\vert J_{r}\right\vert }\right) ^{2}\left\Vert x-m_{J_{r}}\right\Vert
_{L^{2}\left( \mathbf{1}_{J_{r}}\omega \right) }^{2} \\
&\leq &\frac{1}{C_{\limfunc{energy}}\left( \mathcal{E}_{2}^{\alpha }\right)
^{2}}\left( \mathcal{E}_{2}^{\alpha }\right) ^{2}\left\vert S_{1}\right\vert
_{\sigma }=\frac{1}{C_{\limfunc{energy}}}\left\vert S_{1}\right\vert
_{\sigma }\ ,
\end{eqnarray*}%
upon noting that the union of the subdecompositions $\dot{\cup}J_{r}\subset
S $ over $S\in \mathfrak{C}_{\mathcal{S}}\left( S_{1}\right) $ is a
subdecomposition of $S_{1}$, and the proof of the Carleson estimate is now
finished by iteration in the standard way.

Finally, we record the reason for introducing energy stopping times. If 
\begin{equation}
\mathbf{X}_{\alpha }\left( \mathcal{C}_{S}\right) ^{2}\equiv \sup_{I\in 
\mathcal{C}_{S}}\frac{1}{\left\vert I\right\vert _{\sigma }}\sup_{I\supset 
\dot{\cup}J_{r}}\sum_{r=1}^{\infty }\left( \frac{\mathrm{P}^{\alpha }\left(
J_{r},\mathbf{1}_{S}\sigma \right) }{\left\vert J_{r}\right\vert }\right)
^{2}\left\Vert x-m_{J_{r}}\right\Vert _{L^{2}\left( \mathbf{1}_{J_{r}}\omega
\right) }^{2}  \label{def stopping energy 3}
\end{equation}%
is (the square of) the $\alpha $\emph{-stopping energy} of the weight pair $%
\left( \sigma ,\omega \right) $ with respect to the corona $\mathcal{C}_{S}$%
, then we have the \emph{stopping energy bounds}%
\begin{equation}
\mathbf{X}_{\alpha }\left( \mathcal{C}_{S}\right) \leq \sqrt{C_{\limfunc{%
energy}}}\sqrt{\left( \mathfrak{E}_{2}^{\alpha }\right) ^{2}+\mathfrak{A}%
_{2}^{\alpha }},\ \ \ \ \ S\in \mathcal{S},  \label{def stopping bounds 3}
\end{equation}%
where $\mathfrak{A}_{2}^{\alpha }$ and the energy constant $\mathfrak{E}%
_{2}^{\alpha }$ are controlled by the assumptions in Theorem \ref{dim one}.

\subsection{Iterated coronas and general stopping data}

We will use a construction that permits \emph{iteration} of the above three
corona decompositions by combining Definitions \ref{CZ stopping times}, \ref%
{accretive stopping times gen} and \ref{def energy corona 3} into a single
stopping condition. However, there is one remaining difficulty with the
triple corona constructed in this way, namely if a stopping interval $I\in 
\mathcal{A}$ is a child of an interval $Q$ in the corona $\mathcal{C}_{A}$,
then the modulus of the average $\left\vert \frac{1}{\left\vert I\right\vert
_{\sigma }}\int_{I}b_{Q}d\sigma \right\vert $ of $b_{Q}$ on $I$ may be far
smaller than the sup norm of $\left\vert b_{Q}\right\vert $ on the child $I$%
, indeed it may be that $\frac{1}{\left\vert I\right\vert _{\sigma }}%
\int_{I}b_{Q}d\sigma =0$. This of course destroys any reasonable estimation
of the martingale and dual martingale differences $\bigtriangleup
_{Q}^{\sigma ,\mathbf{b}}f$ and $\square _{Q}^{\sigma ,\mathbf{b}}f$ used in
the proof of Theorem \ref{dim one}, and so we will use Lemma \ref{prelim
control of corona} on the function $b_{A}$ to obtain a new function $%
\widetilde{b}_{A}$ for which this problem is circumvented at the `bottom' of
the corona, i.e. for those $A^{\prime }\in C_{\mathcal{A}}\left( A\right) $.
We then refer to the stopping times $A^{\prime }\in C_{\mathcal{A}}\left(
A\right) $ as `shadow' stopping times since we have lost control of the weak
testing condition relative to the new function $\widetilde{b}_{A}$. Thus we
must redo the weak testing stopping times for the new function $\widetilde{b}%
_{A}$, but also stopping if we hit one of the shadow stopping times. Here
are the details.

\begin{definition}
\label{def shadow}Let $C_{0}\geq 4$, $0<\gamma <1$ and $1<\Gamma <\infty $.
Suppose that $\mathbf{b}=\left\{ b_{Q}\right\} _{Q\in \mathcal{P}}$ is an $%
\infty $-strongly $\sigma $-accretive family on $\mathbb{R}$. Given a dyadic
grid $\mathcal{D}$ and an interval $Q\in \mathcal{D}$, define the collection
of `shadow' stopping times $\mathcal{S}_{\limfunc{shadow}}\left( Q\right) $
to be the \emph{maximal} $\mathcal{D}$-subintervals $I\subset Q$ such that
either%
\begin{equation*}
\frac{1}{\left\vert I\right\vert _{\sigma }}\int_{I}\left\vert f\right\vert
d\sigma >C_{0}\frac{1}{\left\vert Q\right\vert _{\sigma }}\int_{Q}\left\vert
f\right\vert d\sigma \ ,
\end{equation*}%
or%
\begin{equation*}
\left\vert \frac{1}{\left\vert I\right\vert _{\mu }}\int_{I}b_{Q}d\sigma
\right\vert <\gamma \ \text{or }\int_{I}\left\vert T_{\sigma }^{\alpha
}\left( b_{Q}\right) \right\vert ^{2}d\omega >\Gamma \left( \mathfrak{T}%
_{T^{\alpha }}^{\mathbf{b}}\right) ^{2}\left\vert I\right\vert _{\sigma }%
\text{ },
\end{equation*}%
or%
\begin{equation*}
\sup_{I\supset \dot{\cup}J_{r}}\sum_{r=1}^{\infty }\left( \frac{\mathrm{P}%
^{\alpha }\left( J_{r},\sigma \right) }{\left\vert J_{r}\right\vert }\right)
^{2}\left\Vert x-m_{J_{r}}\right\Vert _{L^{2}\left( \mathbf{1}_{J_{r}}\omega
\right) }^{2}\geq C_{\limfunc{energy}}\left[ \left( \mathcal{E}_{2}^{\alpha ,%
\mathbf{b},\mathbf{b}^{\ast }}\right) ^{2}+\mathfrak{A}_{2}^{\alpha }\right]
\ \left\vert I\right\vert _{\sigma }\ .
\end{equation*}
\end{definition}

Now we apply Lemma \ref{prelim control of corona}\ to the function $b_{Q}$
with the subdecomposition $\mathcal{S}_{\limfunc{shadow}}\left( Q\right)
\equiv \left\{ Q_{i}\right\} _{i=1}^{\infty }$ to obtain a new function $%
\widetilde{b}_{Q}$ satisfying the properties%
\begin{eqnarray}
&&\limfunc{support}\widetilde{b}_{Q}\subset Q\ ,  \label{props} \\
&&1\leq \frac{1}{\left\vert Q^{\prime }\right\vert _{\sigma }}%
\int_{Q^{\prime }}\widetilde{b}_{Q}d\sigma \leq \left\Vert \mathbf{1}%
_{Q^{\prime }}\widetilde{b}_{Q}\right\Vert _{L^{\infty }\left( \sigma
\right) }\leq 2\left( 1+\sqrt{C_{\mathbf{b}}}\right) C_{\mathbf{b}}\ ,\ \ \
\ \ Q^{\prime }\in \mathcal{C}_{Q}\ ,  \notag \\
&&\sqrt{\int_{Q}\left\vert T_{\sigma }^{\alpha }b_{Q}\right\vert ^{2}d\omega 
}\leq \left[ 2\mathfrak{T}_{T^{\alpha }}^{\mathbf{b}}\left( Q\right) +4C_{%
\mathbf{b}}^{\frac{3}{2}}\delta ^{\frac{1}{4}}\mathfrak{N}_{T^{\alpha
}}\left( Q\right) \right] \sqrt{\left\vert Q\right\vert _{\sigma }}\ , 
\notag \\
&&\left\Vert \mathbf{1}_{Q_{i}}\widetilde{b}_{Q}\right\Vert _{L^{\infty
}\left( \sigma \right) }\leq \frac{16C_{\mathbf{b}}}{\delta }\left\vert 
\frac{1}{\left\vert Q_{i}\right\vert _{\sigma }}\int_{Q_{i}}\widetilde{b}%
_{Q}d\sigma \right\vert \ ,\ \ \ \ \ 1\leq i<\infty .  \notag
\end{eqnarray}%
Note that each of the functions $\widetilde{b}_{Q^{\prime }}\equiv \mathbf{1}%
_{Q^{\prime }}\widetilde{b}_{Q}$, for $Q^{\prime }\in \mathcal{C}_{Q}$, now
satisfies the crucial reverse H\"{o}lder property%
\begin{equation*}
\left\Vert \mathbf{1}_{I}\widetilde{b}_{Q^{\prime }}\right\Vert _{L^{\infty
}\left( \sigma \right) }\leq C_{\delta ,\mathbf{b}}\left\vert \frac{1}{%
\left\vert I\right\vert _{\sigma }}\int_{I}\widetilde{b}_{Q^{\prime
}}d\sigma \right\vert \ ,\ \ \ \ \ \text{for all }I\in \mathfrak{C}\left(
Q^{\prime }\right) ,\ Q^{\prime }\in \mathcal{C}_{Q}.
\end{equation*}%
Indeed, if $I$ equals one of the $Q_{i}$ then the reverse H\"{o}lder
condition in the last line of (\ref{props}) applies, while if $I\in \mathcal{%
C}_{Q}$ then the accretivity in the second line of (\ref{props}) applies.

Since we have lost the weak testing condition in the corona for this new
function $\widetilde{b}_{Q}$, the next step is to run again the weak testing
construction of stopping times, but this time starting with the new function 
$\widetilde{b}_{Q}$, and also stopping if we hit one of the `shadow'
stopping times $Q_{i}$. Here is the new stopping criterion.

\begin{definition}
\label{def iterated}Let $C_{0}\geq 4$ and $1<\Gamma <\infty $. Let $\mathcal{%
S}_{\limfunc{shadow}}\left( Q\right) \equiv \left\{ Q_{i}\right\}
_{i=1}^{\infty }$ be as in Definition \ref{def shadow}. Define $\mathcal{S}_{%
\limfunc{iterated}}\left( Q\right) $ to be the \emph{maximal} $\mathcal{D}$%
-subintervals $I\subset Q$ such that either%
\begin{equation*}
\int_{I}\left\vert T_{\sigma }^{\alpha }\left( \widetilde{b}_{Q}\right)
\right\vert ^{2}d\omega >\Gamma \left( \mathfrak{T}_{T^{\alpha }}^{%
\widetilde{\mathbf{b}}}\right) ^{2}\left\vert I\right\vert _{\sigma }\text{ }%
,
\end{equation*}%
or%
\begin{equation*}
I=Q_{i}\text{ for some }1\leq i<\infty .
\end{equation*}
\end{definition}

Thus for each interval $Q$ we have now constructed \emph{iterated stopping
children} $\mathcal{S}_{\limfunc{iterated}}\left( Q\right) $ by first
constructing shadow stopping times $\mathcal{S}_{\limfunc{shadow}}\left(
Q\right) $ using one step of the triple corona construction, then modifying
the testing function to have reverse H\"{o}lder controlled children, and
finally running again the weak testing stopping time construction to get $%
\mathcal{S}_{\limfunc{iterated}}\left( Q\right) $. These iterated stopping
times $\mathcal{S}_{\limfunc{iterated}}\left( Q\right) $ have control of CZ
averages of $f$ and energy control of $\sigma $ and $\omega $, simply
because these controls were achieved in the shadow construction, and were
unaffected by either the application of Lemma \ref{prelim control of corona}
or the rerunning of the weak testing stopping criterion for $\widetilde{b}%
_{Q}$. And of course we now have weak testing within the corona determined
by $Q$ and $\mathcal{S}_{\limfunc{iterated}}\left( Q\right) $, and we also
have the crucial reverse H\"{o}lder condition on all the children of
intervals in the corona. With all of this in hand, here then is the
definition of the construction of iterated coronas.

\begin{definition}
\label{iterated stopping times}Let $C_{0}\geq 4$, $0<\gamma <1$ and $%
1<\Gamma <\infty $. Suppose that $\mathbf{b}=\left\{ b_{Q}\right\} _{Q\in 
\mathcal{P}}$ is an $\infty $-strongly $\sigma $-accretive family on $%
\mathbb{R}$. Given a dyadic grid $\mathcal{D}$ and an interval $S_{0}$ in $%
\mathcal{D}$, define the iterated stopping intervals of $S_{0}$ to be the
collection 
\begin{equation*}
\mathcal{S}=\left\{ S_{0}\right\} \cup \dbigcup\limits_{n=0}^{\infty }%
\mathcal{S}_{n}
\end{equation*}%
where $\mathcal{S}_{0}=\mathcal{S}_{\limfunc{iterated}}\left( S_{0}\right) $
and $\mathcal{S}_{n+1}=\dbigcup\limits_{S\in \mathcal{S}_{n}}\mathcal{S}_{%
\limfunc{iterated}}\left( S\right) $ for $n\geq 0$, and where $\mathcal{S}_{%
\limfunc{iterated}}\left( Q\right) $ is defined in Definition \ref{def
iterated}.
\end{definition}

It is useful to append to the notion of stopping times $\mathcal{S}$ in the
above $\sigma $-iterated corona decomposition a positive constant $A_{0}$
and an additional structure $\alpha _{\mathcal{S}}$ called stopping bounds
for a function $f$. We will refer to the resulting\ triple $\left( A_{0},%
\mathcal{F},\alpha _{\mathcal{F}}\right) $ as constituting stopping data for 
$f$. If $\mathcal{F}$ is a grid, we define $F^{\prime }\prec F$ if $%
F^{\prime }\subsetneqq F$ and $F^{\prime },F\in \mathcal{F}$. Recall that $%
\pi _{\mathcal{F}}F^{\prime }$ is the smallest $F\in \mathcal{F}$ such that $%
F^{\prime }\prec F$.

\begin{definition}
\label{general stopping data}Suppose we are given a positive constant $%
A_{0}\geq 4$, a subset $\mathcal{F}$ of the dyadic grid $\mathcal{D}$
(called the stopping times), and a corresponding sequence $\alpha _{\mathcal{%
F}}\equiv \left\{ \alpha _{\mathcal{F}}\left( F\right) \right\} _{F\in 
\mathcal{F}}$ of nonnegative numbers $\alpha _{\mathcal{F}}\left( F\right)
\geq 0$ (called the stopping bounds). Let $\left( \mathcal{F},\prec ,\pi _{%
\mathcal{F}}\right) $ be the tree structure on $\mathcal{F}$ inherited from $%
\mathcal{D}$, and for each $F\in \mathcal{F}$ denote by $\mathcal{C}%
_{F}=\left\{ I\in \mathcal{D}:\pi _{\mathcal{F}}I=F\right\} $ the corona
associated with $F$: 
\begin{equation*}
\mathcal{C}_{F}=\left\{ I\in \mathcal{D}:I\subset F\text{ and }I\not\subset
F^{\prime }\text{ for any }F^{\prime }\prec F\right\} .
\end{equation*}%
We say the triple $\left( A_{0},\mathcal{F},\alpha _{\mathcal{F}}\right) $
constitutes \emph{stopping data} for a function $f\in L_{loc}^{1}\left(
\sigma \right) $ if

\begin{enumerate}
\item $\mathbb{E}_{I}^{\sigma }\left\vert f\right\vert \leq \alpha _{%
\mathcal{F}}\left( F\right) $ for all $I\in \mathcal{C}_{F}$ and $F\in 
\mathcal{F}$,

\item $\sum_{F^{\prime }\preceq F}\left\vert F^{\prime }\right\vert _{\sigma
}\leq A_{0}\left\vert F\right\vert _{\sigma }$ for all $F\in \mathcal{F}$,

\item $\sum_{F\in \mathcal{F}}\alpha _{\mathcal{F}}\left( F\right)
^{2}\left\vert F\right\vert _{\sigma }\mathbf{\leq }A_{0}^{2}\left\Vert
f\right\Vert _{L^{2}\left( \sigma \right) }^{2}$,

\item $\alpha _{\mathcal{F}}\left( F\right) \leq \alpha _{\mathcal{F}}\left(
F^{\prime }\right) $ whenever $F^{\prime },F\in \mathcal{F}$ with $F^{\prime
}\subset F$.
\end{enumerate}
\end{definition}

Property (1) says that $\alpha _{\mathcal{F}}\left( F\right) $ bounds the
averages of $f$ in the corona $\mathcal{C}_{F}$, and property (2) says that
the intervals at the tops of the coronas satisfy a Carleson condition
relative to the weight $\sigma $. Note that a standard `maximal interval'
argument extends the Carleson condition in property (2) to the inequality%
\begin{equation}
\sum_{F^{\prime }\in \mathcal{F}:\ F^{\prime }\subset A}\left\vert F^{\prime
}\right\vert _{\sigma }\leq A_{0}\left\vert A\right\vert _{\sigma }\text{
for all open sets }A\subset \mathbb{R}.  \label{Car ext}
\end{equation}%
Property (3) is the quasi-orthogonality condition that says the sequence of
functions $\left\{ \alpha _{\mathcal{F}}\left( F\right) \mathbf{1}%
_{F}\right\} _{F\in \mathcal{F}}$ is in the vector-valued space $L^{2}\left(
\ell ^{2};\sigma \right) $ with control, and is often referred to as a
Carleson embedding theorem, and property (4) says that the control on
stopping data is nondecreasing on the stopping tree $\mathcal{F}$. We
emphasize that we are \emph{not} assuming in this definition the stronger
property that there is $C>1$ such that $\alpha _{\mathcal{F}}\left(
F^{\prime }\right) >C\alpha _{\mathcal{F}}\left( F\right) $ whenever $%
F^{\prime },F\in \mathcal{F}$ with $F^{\prime }\subsetneqq F$. Instead, the
properties (2) and (3) substitute for this lack. Of course the stronger
property \emph{does} hold for the familiar \emph{Calder\'{o}n-Zygmund}
stopping data determined by the following requirements for $C>1$,%
\begin{equation*}
\mathbb{E}_{F^{\prime }}^{\sigma }\left\vert f\right\vert >C\mathbb{E}%
_{F}^{\sigma }\left\vert f\right\vert \text{ whenever }F^{\prime },F\in 
\mathcal{F}\text{ with }F^{\prime }\subsetneqq F,\ \ \ \ \ \mathbb{E}%
_{I}^{\sigma }\left\vert f\right\vert \leq C\mathbb{E}_{F}^{\sigma
}\left\vert f\right\vert \text{ for }I\in \mathcal{C}_{F},
\end{equation*}%
which are themselves sufficiently strong to automatically force properties
(2) and (3) with $\alpha _{\mathcal{F}}\left( F\right) =\mathbb{E}%
_{F}^{\sigma }\left\vert f\right\vert $.

We have the following useful consequence of (2) and (3) that says the
sequence $\left\{ \alpha _{\mathcal{F}}\left( F\right) \mathbf{1}%
_{F}\right\} _{F\in \mathcal{F}}$ has a \emph{quasi-orthogonal} property
relative to $f$ with a constant $A_{0}^{\prime }$ depending only on $A_{0}$
(see e.g. \cite{SaShUr7}):%
\begin{equation}
\left\Vert \sum_{F\in \mathcal{F}}\alpha _{\mathcal{F}}\left( F\right) 
\mathbf{1}_{F}\right\Vert _{L^{2}\left( \sigma \right) }^{2}\leq
A_{0}^{\prime }\left\Vert f\right\Vert _{L^{2}\left( \sigma \right) }^{2}.
\label{q orth}
\end{equation}

\begin{proposition}
\label{data}Let $f\in L^{2}\left( \sigma \right) $, let $\mathcal{F}$ be the
iterated corona $\mathcal{S}\left( S_{0}\right) $ in Definition \ref%
{iterated stopping times}, and define stopping data $\alpha _{\mathcal{F}}$
by $\alpha _{F}=\frac{1}{\left\vert F\right\vert _{\sigma }}%
\int_{F}\left\vert f\right\vert d\sigma $. Then there is $A_{0}\geq 4$,
depending only on the constant $C_{0}$ in Definition \ref{CZ stopping times}%
, such that the triple $\left( A_{0},\mathcal{F},\alpha _{\mathcal{F}%
}\right) $ constitutes \emph{stopping data} for the function $f$.
\end{proposition}

\begin{proof}
This is an easy exercise using (\ref{CZ Car}), (\ref{Car gen}) and (\ref%
{sigma Carleson 3}), and is left for the reader.
\end{proof}

\subsection{Grid parameterizations}

It is important to use \emph{two} independent random grids, one for each
function $f$ and $g$ simultaneously, as this is necessary in order to apply
probabilistic methods to the dual martingale averages $\square _{I}^{\mu ,%
\mathbf{b}}$ that depend, not only on $I$, but also on the underlying grid
in which $I$ lives. The proof methods for functional energy from \cite%
{SaShUr7} and \cite{SaShUr6} relied heavily on the use of a single grid, and
this must now be modified to accommodate two independent grids.

Now we recall the construction from our paper \cite{SaShUr10}. We
momentarily fix a large positive integer $M\in \mathbb{N}$, and consider the
tiling of $\mathbb{R}$ by the family of intervals $\mathbb{D}_{M}\equiv
\left\{ I_{\alpha }^{M}\right\} _{\alpha \in \mathbb{Z}}$ having side length 
$2^{-M}$ and given by $I_{\alpha }^{M}\equiv I_{0}^{M}+2^{-M}\alpha $ where $%
I_{0}^{M}=\left[ 0,2^{-M}\right) $. A \emph{dyadic grid} $\mathcal{D}$ built
on $\mathbb{D}_{M}$ is\ defined to be a family of intervals $\mathcal{D}$
satisfying:

\begin{enumerate}
\item Each $I\in \mathcal{D}$ has side length $2^{-\ell }$ for some $\ell
\in \mathbb{Z}$ with $\ell \leq M$, and $I$ is a union of $2^{M-\ell }$
intervals from the tiling $\mathbb{D}_{M}$,

\item For $\ell \leq M$, the collection $\mathcal{D}_{\ell }$ of intervals
in $\mathcal{D}$ having side length $2^{-\ell }$ forms a pairwise disjoint
decomposition of the space $\mathbb{R}$,

\item Given $I\in \mathcal{D}_{i}$ and $J\in \mathcal{D}_{j}$ with $j\leq
i\leq M$, it is the case that either $I\cap J=\emptyset $ or $I\subset J$.
\end{enumerate}

We now momentarily fix a \emph{negative} integer $N\in -\mathbb{N}$, and
restrict the above grids to intervals of side length at most $2^{-N}$:%
\begin{equation*}
\mathcal{D}^{N}\equiv \left\{ I\in \mathcal{D}:\text{side length of }I\text{
is at most }2^{-N}\right\} \text{.}
\end{equation*}%
We refer to such grids $\mathcal{D}^{N}$ as a (truncated) dyadic grid $%
\mathcal{D}$ built on $\mathbb{D}_{M}$ of size $2^{-N}$. There are now two
traditional means of constructing probability measures on collections of
such dyadic grids, namely parameterization by choice of parent, and
parameterization by translation.

\textbf{Construction \#1}: For any 
\begin{equation*}
\beta =\{\beta _{i}\}_{i\in _{M}^{N}}\in \omega _{M}^{N}\equiv \left\{
0,1\right\} ^{\mathbb{Z}_{M}^{N}},
\end{equation*}%
where $\mathbb{Z}_{M}^{N}\equiv \left\{ \ell \in \mathbb{Z}:N\leq \ell \leq
M\right\} $, define the dyadic grid $\mathcal{D}_{\beta }$ built on $\mathbb{%
D}_{m}$ of size $2^{-N}$ by 
\begin{equation}
\mathcal{D}_{\beta }=\left\{ 2^{-\ell }\left( [0,1)+k+\sum_{i:\ \ell <i\leq
m}2^{-i+\ell }\beta _{i}\right) \right\} _{N\leq \ell \leq m,\,k\in {\mathbb{%
Z}}}\ .  \label{def dyadic grid}
\end{equation}%
Place the uniform probability measure $\rho _{M}^{N}$ on the finite index
space $\omega _{M}^{N}=\left\{ 0,1\right\} ^{\mathbb{Z}_{M}^{N}}$, namely
that which charges each $\beta \in \omega _{M}^{N}$ equally.

\textbf{Construction \#2}: Momentarily fix a (truncated) dyadic grid $%
\mathcal{D}$ built on $\mathbb{D}_{M}$ of size $2^{-N}$. For any 
\begin{equation*}
\gamma \in \Gamma _{M}^{N}\equiv \left\{ \gamma \in 2^{-M}\mathbb{Z}%
_{+}:\left\vert \gamma \right\vert <2^{-N}\right\} ,
\end{equation*}%
define the dyadic grid $\mathcal{D}^{\gamma }$ built on $\mathbb{D}_{M}$ of
size $2^{-N}$ by%
\begin{equation*}
\mathcal{D}^{\gamma }\equiv \mathcal{D}+\gamma .
\end{equation*}%
Place the uniform probability measure $\nu _{M}^{N}$ on the finite index set 
$\Gamma _{M}^{N}$, namely that which charges each multiindex $\gamma $ in $%
\Gamma _{M}^{N}$ equally.

The two probability spaces $\left( \left\{ \mathcal{D}_{\beta }\right\}
_{\beta \in \Omega _{M}^{N}},\mu _{M}^{N}\right) $ and $\left( \left\{ 
\mathcal{D}^{\gamma }\right\} _{\gamma \in \Gamma _{M}^{N}},\nu
_{M}^{N}\right) $ are isomorphic since both collections $\left\{ \mathcal{D}%
_{\beta }\right\} _{\beta \in \Omega _{M}^{N}}$ and $\left\{ \mathcal{D}%
^{\gamma }\right\} _{\gamma \in \Gamma _{M}^{N}}$ describe the set $%
\boldsymbol{A}_{M}^{N}$ of \textbf{all} (truncated) dyadic grids $\mathcal{D}%
^{\gamma }$ built on $\mathbb{D}_{m}$ of size $2^{-N}$, and since both
measures $\mu _{M}^{N}$ and $\nu _{M}^{N}$ are the uniform measure on this
space. The first construction may be thought of as being \emph{parameterized
by scales} - each component $\beta _{i}$ in $\beta =\{\beta _{i}\}_{i\in
_{M}^{N}}\in \omega _{M}^{N}$ amounting to a choice of the two possible
tilings at level $i$ that respect the choice of tiling at the level below -
and since any grid in $\boldsymbol{A}_{M}^{N}$ is determined by a choice of
scales , we see that $\left\{ \mathcal{D}_{\beta }\right\} _{\beta \in
\Omega _{M}^{N}}=\boldsymbol{A}_{M}^{N}$. The second construction may be
thought of as being \emph{parameterized by translation} - each $\gamma \in
\Gamma _{M}^{N}$ amounting to a choice of translation of the grid $\mathcal{D%
}$ fixed in construction \#2\ - and since any grid in $\boldsymbol{A}%
_{M}^{N} $ is determined by any of the intervals at the top level, i.e. with
side length $2^{-N}$, we see that $\left\{ \mathcal{D}^{\gamma }\right\}
_{\gamma \in \Gamma _{M}^{N}}=\boldsymbol{A}_{M}^{N}$ as well, since every
interval at the top level in $\boldsymbol{A}_{M}^{N}$ has the form $Q+\gamma 
$ for some $\gamma \in \Gamma _{M}^{N}$ and $Q\in \mathcal{D}$ at the top
level in $\boldsymbol{A}_{M}^{N}$ (i.e. every interval at the top level in $%
\boldsymbol{A}_{M}^{N}$ is a union of small intervals in $\mathbb{D}_{m}$,
and so must be a translate of some $Q\in \mathcal{D}$ by an amount $2^{-M}$
times an element of $\mathbb{Z}_{+}$). Note also that $\#\Omega
_{M}^{N}=\#\Gamma _{M}^{N}=2^{M-N}$. We will use $\boldsymbol{E}_{\Omega
_{M}^{N}}$ to denote expectation with respect to this common probability
measure on $\boldsymbol{A}_{M}^{N}$.

\begin{notation}
\label{suppress M and N}For purposes of notation and clarity, we now
suppress all reference to $M$ and $N$ in our families of grids, and in the
notations $\Omega $ and $\Gamma $ for the parameter sets, and we use $%
\boldsymbol{P}_{\Omega }$ and $\boldsymbol{E}_{\Omega }$ to denote
probability and expectation with respect to families of grids, and instead
proceed as if all grids considered are unrestricted. The careful reader can
supply the modifications necessary to handle the assumptions made above on
the grids $\mathcal{D}$ and the functions $f$ and $g$ regarding $M$ and $N$.
\end{notation}

\subsection{The Monotonicity Lemma}

As in virtually all proofs of a two weight $T1$ theorem (see e.g. \cite%
{LaSaShUr3}, \cite{Lac}, \cite{SaShUr7} and/or \cite{SaShUr6}), the key to
starting an estimate for any of the forms we consider below, is the
Monotonicity Lemma and the Energy Lemma, to which we now turn. In dimension $%
n=1$ (\cite{LaSaShUr3}, \cite{Lac}) the Haar functions have opposite sign on
their children, and this was exploited in a simple but powerful monotonicity
argument. In higher dimensions, this simple argument no longer holds and
that Monotonicity Lemma is replaced with the Lacey-Wick formulation of the
Monotonicity Lemma (see \cite{LaWi}, and also \cite{SaShUr6}) involving the
smaller Poisson operator. As the martingale differences with test functions $%
b_{Q\,}$ here are no longer of one sign on children, we will adapt the
Lacey-Wick formulation of the Monotonicity Lemma to the operator $T^{\alpha
} $ and the dual martingale differences $\left\{ \square _{J}^{\omega ,%
\mathbf{b}^{\ast }}\right\} _{J\in \mathcal{G}}$, bearing in mind that the
operators $\square _{J}^{\omega ,\mathbf{b}^{\ast }}$ are no longer
projections, which results in only a one-sided estimate with additional
terms on the right hand side. It is here that we need the crucial property
that $\limfunc{Range}\square _{J}^{\omega ,\mathbf{b}^{\ast }}$ is
orthogonal to constants, $\int \left( \square _{J}^{\omega ,\mathbf{b}^{\ast
}}\Psi \right) d\sigma =\int \left( \triangle _{J}^{\sigma ,\mathbf{b}^{\ast
}}1\right) \Psi d\omega =\int \left( 0\right) \Psi d\omega =0$ (see Appendix
A). See Definition \ref{controlled accretive} in Appendix A for the
terminology `$p$-weakly $\mu $-controlled accretive family' along with more
detail on martingale and dual martingale expansions.

Recall from Appendix A that%
\begin{eqnarray*}
\mathbb{E}_{Q}^{\mu ,\mathbf{b}}f\left( x\right) &\equiv &\mathbf{1}%
_{Q}\left( x\right) \frac{1}{\int_{Q}b_{Q}d\mu }\int_{Q}fb_{Q}d\mu ,\ \ \ \
\ Q\in \mathcal{P}\ , \\
\mathbb{F}_{Q}^{\mu ,\mathbf{b}}f\left( x\right) &\equiv &\mathbf{1}%
_{Q}\left( x\right) b_{Q}\left( x\right) \frac{1}{\int_{Q}b_{Q}d\mu }%
\int_{Q}fd\mu ,\ \ \ \ \ Q\in \mathcal{P}\ ,
\end{eqnarray*}%
and 
\begin{eqnarray*}
\bigtriangleup _{Q}^{\mu ,\mathbf{b}}f\left( x\right) &\equiv &\left(
\sum_{Q^{\prime }\in \mathfrak{C}\left( Q\right) }\mathbb{E}_{Q^{\prime
}}^{\mu ,\mathbf{b}}f\left( x\right) \right) -\mathbb{E}_{Q}^{\mu ,\mathbf{b}%
}f\left( x\right) =\sum_{Q^{\prime }\in \mathfrak{C}\left( Q\right) }\mathbf{%
1}_{Q^{\prime }}\left( x\right) \left( \mathbb{E}_{Q^{\prime }}^{\mu ,%
\mathbf{b}}f\left( x\right) -\mathbb{E}_{Q}^{\mu ,\mathbf{b}}f\left(
x\right) \right) , \\
\square _{Q}^{\mu ,\mathbf{b}}f\left( x\right) &\equiv &\left(
\sum_{Q^{\prime }\in \mathfrak{C}\left( Q\right) }\mathbb{F}_{Q^{\prime
}}^{\mu ,\mathbf{b}}f\left( x\right) \right) -\mathbb{F}_{Q}^{\mu ,\mathbf{b}%
}f\left( x\right) =\sum_{Q^{\prime }\in \mathfrak{C}\left( Q\right) }\mathbf{%
1}_{Q^{\prime }}\left( x\right) \left( \mathbb{F}_{Q^{\prime }}^{\mu ,%
\mathbf{b}}f\left( x\right) -\mathbb{F}_{Q}^{\mu ,\mathbf{b}}f\left(
x\right) \right) ,
\end{eqnarray*}%
and from (\ref{Carleson avg op}),%
\begin{equation*}
\nabla _{Q}^{\mu }h=\sum_{Q^{\prime }\in \mathfrak{C}_{\limfunc{broken}%
}\left( Q\right) }\left( E_{Q^{\prime }}^{\mu }\left\vert h\right\vert
\right) ^{2}\mathbf{1}_{Q^{\prime }}\ .
\end{equation*}

We will also need the smaller Poisson integral used in the Lacey-Wick
formulation of the Monontonicity Lemma,%
\begin{equation*}
\mathrm{P}_{1+\delta }^{\alpha }\left( J,\mu \right) \equiv \int \frac{%
\left\vert J\right\vert ^{\frac{1+\delta }{n}}}{\left( \left\vert
J\right\vert +\left\vert y-c_{J}\right\vert \right) ^{n+1+\delta -\alpha }}%
d\mu \left( y\right) ,
\end{equation*}%
which is discussed in more detail below.

\begin{lemma}[Monotonicity Lemma]
\label{mono}Suppose that$\ I$ and $J$ are intervals in $\mathbb{R}$ such
that $J\subset \gamma J\subset I$ for some $\gamma >1$, and that $\mu $ is a
signed measure on $\mathbb{R}$ supported outside $I$. Let $0<\delta <1$ and
let $\Psi \in L^{2}\left( \omega \right) $. Finally suppose that $T^{\alpha
} $ is a standard fractional singular integral on $\mathbb{R}$ as in \cite%
{SaShUr6}, \cite{SaShUr7} and \cite{SaShUr9} with $0\leq \alpha <1$, and
suppose that $\mathbf{b}^{\ast }$ is an $\infty $-weakly $\mu $-controlled
accretive family on $\mathbb{R}$. Then we have the estimate%
\begin{equation}
\left\vert \left\langle T^{\alpha }\mu ,\square _{J}^{\omega ,\mathbf{b}%
^{\ast }}\Psi \right\rangle _{\omega }\right\vert \lesssim C_{\mathbf{b}%
^{\ast }}C_{CZ}\ \Phi ^{\alpha }\left( J,\left\vert \mu \right\vert \right)
\ \left\Vert \square _{J}^{\omega ,\mathbf{b}^{\ast }}\Psi \right\Vert
_{L^{2}\left( \omega \right) }^{\bigstar },  \label{estimate}
\end{equation}%
where%
\begin{eqnarray*}
\Phi ^{\alpha }\left( J,\left\vert \mu \right\vert \right) &\equiv &\frac{%
\mathrm{P}^{\alpha }\left( J,\left\vert \mu \right\vert \right) }{\left\vert
J\right\vert }\left\Vert \bigtriangleup _{J}^{\omega ,\mathbf{b}^{\ast
}}x\right\Vert _{L^{2}\left( \omega \right) }^{\spadesuit }+\frac{\mathrm{P}%
_{1+\delta }^{\alpha }\left( J,\left\vert \mu \right\vert \right) }{%
\left\vert J\right\vert }\left\Vert x-m_{J}\right\Vert _{L^{2}\left( \mathbf{%
1}_{J}\omega \right) }, \\
\left\Vert \bigtriangleup _{J}^{\omega ,\mathbf{b}^{\ast }}x\right\Vert
_{L^{2}\left( \omega \right) }^{\spadesuit 2} &\equiv &\left\Vert
\bigtriangleup _{J}^{\omega ,\mathbf{b}^{\ast }}x\right\Vert _{L^{2}\left(
\omega \right) }^{2}+\inf_{z\in \mathbb{R}}\sum_{J^{\prime }\in \mathfrak{C}%
_{\limfunc{broken}}\left( J\right) }\left\vert J^{\prime }\right\vert
_{\omega }\left( E_{J^{\prime }}^{\omega }\left\vert x-z\right\vert \right)
^{2}, \\
\left\Vert \square _{J}^{\omega ,\mathbf{b}^{\ast }}\Psi \right\Vert
_{L^{2}\left( \mu \right) }^{\bigstar 2} &\equiv &\left\Vert \square
_{J}^{\omega ,\mathbf{b}^{\ast }}\Psi \right\Vert _{L^{2}\left( \mu \right)
}^{2}+\sum_{J^{\prime }\in \mathfrak{C}_{\limfunc{broken}}\left( J\right)
}\left\vert J^{\prime }\right\vert _{\omega }\left[ E_{J^{\prime }}^{\omega
}\left\vert \Psi \right\vert \right] ^{2}.
\end{eqnarray*}%
All of the implied constants above depend only on $\gamma >1$, $0<\delta <1$
and $0<\alpha <1$.
\end{lemma}

Using $\bigtriangledown _{J}^{\omega }h=\sum_{J^{\prime }\in \mathfrak{C}_{%
\limfunc{broken}}\left( JQ\right) }\left( E_{J^{\prime }}^{\omega
}\left\vert h\right\vert \right) ^{2}\mathbf{1}_{J^{\prime }}$ defined in (%
\ref{Carleson avg op}) in Appendix A, we can rewrite the expressions $%
\left\Vert \bigtriangleup _{J}^{\omega ,\mathbf{b}^{\ast }}x\right\Vert
_{L^{2}\left( \omega \right) }^{\spadesuit 2}$ and $\left\Vert \square
_{J}^{\omega ,\mathbf{b}^{\ast }}\Psi \right\Vert _{L^{2}\left( \mu \right)
}^{\bigstar 2}$ as%
\begin{eqnarray*}
\left\Vert \bigtriangleup _{J}^{\omega ,\mathbf{b}^{\ast }}x\right\Vert
_{L^{2}\left( \omega \right) }^{\spadesuit 2} &\equiv &\left\Vert
\bigtriangleup _{J}^{\omega ,\mathbf{b}^{\ast }}x\right\Vert _{L^{2}\left(
\omega \right) }^{2}+\inf_{z\in \mathbb{R}}\left\Vert \bigtriangledown
_{J}^{\omega }\left( x-z\right) \right\Vert _{L^{2}\left( \omega \right)
}^{2}, \\
\left\Vert \square _{J}^{\omega ,\mathbf{b}^{\ast }}\Psi \right\Vert
_{L^{2}\left( \mu \right) }^{\bigstar 2} &\equiv &\left\Vert \square
_{J}^{\omega ,\mathbf{b}^{\ast }}\Psi \right\Vert _{L^{2}\left( \mu \right)
}^{2}+\left\Vert \bigtriangledown _{J}^{\omega }\Psi \right\Vert
_{L^{2}\left( \omega \right) }^{2}.
\end{eqnarray*}

\begin{proof}
We also use formulas (\ref{def pi box}), (\ref{square of delta}) and the
estimate (\ref{F est}) from Appendix A:%
\begin{eqnarray*}
\square _{Q}^{\mu ,\pi ,\mathbf{b}}f &=&\left[ \sum_{Q^{\prime }\in 
\mathfrak{C}\left( Q\right) }\mathbb{F}_{Q^{\prime }}^{\mu ,\pi ,\mathbf{b}}f%
\right] -\mathbb{F}_{Q}^{\mu ,\mathbf{b}}f=\sum_{Q^{\prime }\in \mathfrak{C}%
\left( Q\right) }\mathbb{F}_{Q^{\prime }}^{\mu ,b_{Q}}f-\mathbb{F}_{Q}^{\mu
,b_{Q}}f, \\
\mathbb{F}_{Q}^{\mu ,\pi ,\mathbf{b}}f &=&\mathbf{1}_{Q}\frac{b_{\pi Q}}{%
\int_{Q}b_{\pi Q}d\mu }\int_{Q}fd\mu , \\
\square _{Q}^{\mu ,\mathbf{b}} &=&\square _{Q}^{\mu ,\pi ,\mathbf{b}}\square
_{Q}^{\mu ,\pi ,\mathbf{b}}+\square _{Q,\limfunc{broken}}^{\mu ,\mathbf{b}}%
\text{ and }\square _{Q,\limfunc{broken}}^{\mu ,\mathbf{b}}f=\sum_{Q^{\prime
}\in \mathfrak{C}_{\limfunc{broken}}\left( Q\right) }\mathbb{F}_{Q^{\prime
}}^{\mu ,b_{Q^{\prime }}}f-\mathbb{F}_{Q^{\prime }}^{\mu ,b_{Q}}f, \\
\left\vert \square _{Q,\limfunc{broken}}^{\mu ,\mathbf{b}}f\right\vert
&\lesssim &\left\vert \bigtriangledown _{Q}^{\mu }f\right\vert ,
\end{eqnarray*}%
with similar equalities and inequalities for $\bigtriangleup $ and $\mathbb{E%
}$. Here $\mathfrak{C}_{\limfunc{broken}}\left( Q\right) $ denotes the set
of broken children, i.e. those $Q^{\prime }\in \mathfrak{C}\left( Q\right) $
for which $b_{Q^{\prime }}\neq \mathbf{1}_{Q^{\prime }}b_{Q}$, and more
generally and typically, $\mathfrak{C}_{\limfunc{broken}}\left( Q\right) =%
\mathfrak{C}\left( Q\right) \cap \mathcal{A}$ where $\mathcal{A}$ is a
collection of stopping intervals that includes the broken children and
satisfies a $\sigma $-Carleson condition. Using $\square _{J}^{\omega ,%
\mathbf{b}^{\ast }}=\square _{J}^{\omega ,\pi ,\mathbf{b}^{\ast }}\square
_{J}^{\omega ,\pi ,\mathbf{b}^{\ast }}+\square _{J,\limfunc{broken}}^{\omega
,\mathbf{b}^{\ast }}$, we write%
\begin{eqnarray*}
\left\vert \left\langle T^{\alpha }\mu ,\square _{J}^{\omega ,\mathbf{b}%
^{\ast }}\Psi \right\rangle _{\omega }\right\vert &=&\left\vert \left\langle
T^{\alpha }\mu ,\left( \square _{J}^{\omega ,\pi ,\mathbf{b}^{\ast }}\square
_{J}^{\omega ,\pi ,\mathbf{b}^{\ast }}+\square _{J,\limfunc{broken}}^{\omega
,\mathbf{b}^{\ast }}\right) \Psi \right\rangle _{\omega }\right\vert \\
&\leq &\left\vert \left\langle T^{\alpha }\mu ,\square _{J}^{\omega ,\pi ,%
\mathbf{b}^{\ast }}\square _{J}^{\omega ,\pi ,\mathbf{b}^{\ast }}\Psi
\right\rangle _{\omega }\right\vert +\left\vert \left\langle T^{\alpha }\mu
,\square _{J,\limfunc{broken}}^{\omega ,\mathbf{b}^{\ast }}\Psi
\right\rangle _{\omega }\right\vert \equiv I+II.
\end{eqnarray*}%
Since $\left\langle 1,\square _{J}^{\omega ,\pi ,\mathbf{b}^{\ast
}}h\right\rangle _{\omega }=0$, we use $m_{J}=\frac{1}{\left\vert
J\right\vert _{\omega }}\int_{J}xd\omega \left( x\right) $ to obtain%
\begin{equation*}
T^{\alpha }\mu \left( x\right) -T^{\alpha }\mu \left( m_{J}\right) =\int 
\left[ \left( K^{\alpha }\right) \left( x,y\right) -\left( K^{\alpha
}\right) \left( m_{J},y\right) \right] d\mu \left( y\right) =\int \left[
\left( K_{y}^{\alpha }\right) ^{\prime }\left( \theta \left( x,m_{J}\right)
\right) \left( x-m_{J}\right) \right] d\mu \left( y\right)
\end{equation*}%
for some $\theta \left( x,m_{J}\right) \in J$ to obtain%
\begin{eqnarray*}
I &=&\left\vert \int \left[ T^{\alpha }\mu \left( x\right) -T^{\alpha }\mu
\left( m_{J}\right) \right] \ \square _{J}^{\omega ,\pi ,\mathbf{b}^{\ast
}}\square _{J}^{\omega ,\pi ,\mathbf{b}^{\ast }}\Psi \left( x\right) d\omega
\left( x\right) \right\vert \\
&=&\left\vert \int \left\{ \int \left[ \left( K_{y}^{\alpha }\right)
^{\prime }\left( \theta \left( x,m_{J}\right) \right) \right] d\mu \left(
y\right) \right\} \ \left( x-m_{J}\right) \ \square _{J}^{\omega ,\pi ,%
\mathbf{b}^{\ast }}\square _{J}^{\omega ,\pi ,\mathbf{b}^{\ast }}\Psi \left(
x\right) d\omega \left( x\right) \right\vert \\
&\leq &\left\vert \int \left\{ \int \left[ \left( K_{y}^{\alpha }\right)
^{\prime }\left( m_{J}\right) \right] d\mu \left( y\right) \right\} \ \left(
x-m_{J}\right) \ \square _{J}^{\omega ,\pi ,\mathbf{b}^{\ast }}\square
_{J}^{\omega ,\pi ,\mathbf{b}^{\ast }}\Psi \left( x\right) d\omega \left(
x\right) \right\vert \\
&&+\left\vert \int \left\{ \int \left[ \left( K_{y}^{\alpha }\right)
^{\prime }\left( \theta \left( x,m_{J}\right) \right) -\left( K_{y}^{\alpha
}\right) ^{\prime }\left( m_{J}\right) \right] d\mu \left( y\right) \right\}
\ \left( x-m_{J}\right) \ \square _{J}^{\omega ,\pi ,\mathbf{b}^{\ast
}}\square _{J}^{\omega ,\pi ,\mathbf{b}^{\ast }}\Psi \left( x\right) d\omega
\left( x\right) \right\vert \\
&\equiv &I_{1}+I_{2}.
\end{eqnarray*}%
Now we estimate%
\begin{eqnarray*}
I_{1} &=&\left\vert \int \left[ \left( K_{y}^{\alpha }\right) ^{\prime
}\left( m_{J}\right) \right] d\mu \left( y\right) \right\vert \ \left\vert
\int \left( x-m_{J}\right) \ \square _{J}^{\omega ,\pi ,\mathbf{b}^{\ast
}}\square _{J}^{\omega ,\pi ,\mathbf{b}^{\ast }}\Psi \left( x\right) d\omega
\left( x\right) \right\vert \\
&=&\left\vert \int \left[ \left( K_{y}^{\alpha }\right) ^{\prime }\left(
m_{J}\right) \right] d\mu \left( y\right) \right\vert \ \left\vert \int
\left( \bigtriangleup _{J}^{\omega ,\pi ,\mathbf{b}^{\ast }}x\right) \
\left( \square _{J}^{\omega ,\pi ,\mathbf{b}^{\ast }}\Psi \left( x\right)
\right) \ d\omega \left( x\right) \right\vert \\
&\lesssim &C_{CZ}\frac{\mathrm{P}^{\alpha }\left( J,\left\vert \mu
\right\vert \right) }{\left\vert J\right\vert }\ \left\Vert \bigtriangleup
_{J}^{\omega ,\pi ,\mathbf{b}^{\ast }}x\right\Vert _{L^{2}\left( \omega
\right) }\left\Vert \square _{J}^{\omega ,\pi ,\mathbf{b}^{\ast }}\Psi
\right\Vert _{L^{2}\left( \omega \right) }\ ,
\end{eqnarray*}%
and%
\begin{eqnarray*}
I_{2} &\lesssim &C_{CZ}\frac{\mathrm{P}_{1+\delta }^{\alpha }\left(
J,\left\vert \mu \right\vert \right) }{\left\vert J\right\vert }\int
\left\vert x-m_{J}\right\vert \left\vert \square _{J}^{\omega ,\pi ,\mathbf{b%
}^{\ast }}\square _{J}^{\omega ,\pi ,\mathbf{b}^{\ast }}\Psi \left( x\right)
\right\vert d\omega \left( x\right) \\
&\lesssim &C_{CZ}\frac{\mathrm{P}_{1+\delta }^{\alpha }\left( J,\left\vert
\mu \right\vert \right) }{\left\vert J\right\vert }\sqrt{\int_{J}\left\vert
x-m_{J}\right\vert ^{2}d\omega \left( x\right) }\left\Vert \square
_{J}^{\omega ,\pi ,\mathbf{b}^{\ast }}\square _{J}^{\omega ,\pi ,\mathbf{b}%
^{\ast }}\Psi \right\Vert _{L^{2}\left( \omega \right) } \\
&\lesssim &C_{CZ}\frac{\mathrm{P}_{1+\delta }^{\alpha }\left( J,\left\vert
\mu \right\vert \right) }{\left\vert J\right\vert }\left\Vert
x-m_{J}\right\Vert _{L^{2}\left( \mathbf{1}_{J}\omega \right) }\left\Vert
\square _{J}^{\omega ,\pi ,\mathbf{b}^{\ast }}\Psi \right\Vert _{L^{2}\left(
\omega \right) }\ .
\end{eqnarray*}%
For term $II$ we fix $z\in \overline{J}$ for the moment. Then since $%
\left\langle 1,\square _{J,\limfunc{broken}}^{\omega ,\mathbf{b}^{\ast
}}h\right\rangle _{\omega }=\left\langle 1,\square _{J}^{\omega ,\mathbf{b}%
^{\ast }}h-\square _{J}^{\omega ,\pi ,\mathbf{b}^{\ast }}h\right\rangle
_{\omega }=0$, we have%
\begin{equation*}
II=\left\vert \left\langle T^{\alpha }\mu ,\square _{J,\limfunc{broken}%
}^{\omega ,\mathbf{b}^{\ast }}\Psi \right\rangle _{\omega }\right\vert
=\left\vert \int \left\{ \int \left[ \left( K_{y}^{\alpha }\right) ^{\prime
}\left( \theta \left( x,z\right) \right) \right] d\mu \left( y\right)
\right\} \left( x-z\right) \ \square _{J,\limfunc{broken}}^{\omega ,\mathbf{b%
}^{\ast }}\Psi \left( x\right) d\omega \left( x\right) \right\vert .
\end{equation*}

Using reverse H\"{o}lder control of children (\ref{rev Hol con}), we obtain
the estimate (\ref{F est}) from Appendix A,%
\begin{equation*}
\left\vert \square _{J,\limfunc{broken}}^{\omega ,\mathbf{b}^{\ast }}\Psi
\right\vert =\left\vert \sum_{J^{\prime }\in \mathfrak{C}_{\limfunc{broken}%
}\left( JQ\right) }\left( \mathbb{F}_{J^{\prime }}^{\omega ,b_{J^{\prime }}}-%
\mathbb{F}_{J^{\prime }}^{\omega ,b_{J}}\right) \Psi \right\vert \lesssim
\sum_{J^{\prime }\in \mathfrak{C}_{\limfunc{broken}}\left( J\right) }\mathbf{%
1}_{J^{\prime }}E_{J^{\prime }}^{\omega }\left\vert \Psi \right\vert ,
\end{equation*}%
and so%
\begin{equation*}
II\lesssim C_{CZ}\frac{\mathrm{P}^{\alpha }\left( J,\left\vert \mu
\right\vert \right) }{\left\vert J\right\vert }\sqrt{\sum_{J^{\prime }\in 
\mathfrak{C}_{\limfunc{broken}}\left( J\right) }\left\vert J^{\prime
}\right\vert _{\omega }\left( E_{J^{\prime }}^{\omega }\left\vert
x-z\right\vert \right) ^{2}}\sqrt{\sum_{J^{\prime }\in \mathfrak{C}_{%
\limfunc{broken}}\left( J\right) }\left\vert J^{\prime }\right\vert _{\omega
}\left[ E_{J^{\prime }}^{\omega }\left\vert \Psi \right\vert \right] ^{2}}.
\end{equation*}

Combining the estimates for terms $I$ and $II$, we obtain%
\begin{eqnarray*}
&&\left\vert \left\langle T^{\alpha }\mu ,\square _{J}^{\omega ,\mathbf{b}%
^{\ast }}\Psi \right\rangle _{\omega }\right\vert \\
&\lesssim &C_{CZ}\frac{\mathrm{P}^{\alpha }\left( J,\left\vert \mu
\right\vert \right) }{\left\vert J\right\vert }\ \left\Vert \bigtriangleup
_{J}^{\omega ,\pi ,\mathbf{b}^{\ast }}x\right\Vert _{L^{2}\left( \omega
\right) }\left\Vert \square _{J}^{\omega ,\pi ,\mathbf{b}^{\ast }}\Psi
\right\Vert _{L^{2}\left( \omega \right) } \\
&&+C_{CZ}\frac{\mathrm{P}_{1+\delta }^{\alpha }\left( J,\left\vert \mu
\right\vert \right) }{\left\vert J\right\vert }\left\Vert x-m_{J}\right\Vert
_{L^{2}\left( \mathbf{1}_{J}\omega \right) }\left\Vert \square _{J}^{\omega
,\pi ,\mathbf{b}^{\ast }}\Psi \right\Vert _{L^{2}\left( \omega \right) } \\
&&+C_{CZ}\frac{\mathrm{P}^{\alpha }\left( J,\left\vert \mu \right\vert
\right) }{\left\vert J\right\vert }\inf_{z\in \overline{J}}\sqrt{%
\sum_{J^{\prime }\in \mathfrak{C}_{\limfunc{broken}}\left( J\right)
}\left\vert J^{\prime }\right\vert _{\omega }\left( E_{J^{\prime }}^{\omega
}\left\vert x-z_{2}\right\vert \right) ^{2}}\sqrt{\sum_{J^{\prime }\in 
\mathfrak{C}_{\limfunc{broken}}\left( J\right) }\left\vert J^{\prime
}\right\vert _{\omega }\left[ E_{J^{\prime }}^{\omega }\left\vert \Psi
\right\vert +E_{J}^{\omega }\left\vert \Psi \right\vert \right] ^{2}},
\end{eqnarray*}%
and then noting that the infimum over $z\in \mathbb{R}$ is achieved for $%
z\in \overline{J}$, and using the triangle inequality on $\square
_{J}^{\omega ,\pi ,\mathbf{b}^{\ast }}=\square _{J}^{\omega ,\mathbf{b}%
^{\ast }}-\square _{J,\limfunc{broken}}^{\omega ,\mathbf{b}^{\ast }}$ we get
(\ref{estimate}).
\end{proof}

The right hand side of (\ref{estimate}) in the Monotonicity Lemma will be
typically estimated in what follows using the frame inequalities (see
Appendix A) for any interval $K$,%
\begin{eqnarray*}
\sum_{J\subset K}\left\Vert \square _{J}^{\omega ,\mathbf{b}^{\ast }}\Psi
\right\Vert _{L^{2}\left( \omega \right) }^{\bigstar 2} &\lesssim
&\left\Vert \Psi \right\Vert _{L^{2}\left( \omega \right) }^{2}\ , \\
\sum_{J\subset K}\left\Vert \bigtriangleup _{J}^{\omega ,\mathbf{b}^{\ast
}}x\right\Vert _{L^{2}\left( \omega \right) }^{\spadesuit 2} &\lesssim
&\int_{K}\left\vert x-m_{K}\right\vert ^{2}d\omega \left( x\right) \ ,
\end{eqnarray*}%
together with these inequalities for the square function expressions.

\begin{lemma}
For any interval $K$ we have%
\begin{eqnarray}
\sum_{J\subset K}\sum_{J^{\prime }\in \mathfrak{C}_{\limfunc{broken}}\left(
J\right) }\left\vert J^{\prime }\right\vert _{\omega }\left[ E_{J^{\prime
}}^{\omega }\left\vert \Psi \right\vert \left( x\right) \right] ^{2}
&\lesssim &\int_{K}\left\vert \Psi \left( x\right) \right\vert ^{2}d\omega
\left( x\right) ,  \label{with both} \\
\text{and }\sum_{J\subset K}\inf_{z\in \mathbb{R}}\sum_{J^{\prime }\in 
\mathfrak{C}_{\limfunc{broken}}\left( J\right) }\left\vert J^{\prime
}\right\vert _{\omega }\left( E_{J^{\prime }}^{\omega }\left\vert
x-z\right\vert \right) ^{2} &\lesssim &\int_{K}\left\vert x-m_{K}\right\vert
^{2}d\omega \left( x\right) .  \notag
\end{eqnarray}
\end{lemma}

\begin{proof}
The first inequality in (\ref{with both}) is just the Carleson embedding
theorem since the intervals $\left\{ J^{\prime }\in \mathfrak{C}_{\limfunc{%
broken}}\left( J\right) :J\subset K\right\} $ satisfy an $\omega $-Carleson
condition, and the second inequality in (\ref{with both}) follows by
choosing $z=m_{K}$ to obtain%
\begin{equation*}
\inf_{z\in \mathbb{R}}\sum_{J^{\prime }\in \mathfrak{C}_{\limfunc{broken}%
}\left( J\right) }\left\vert J^{\prime }\right\vert _{\omega }\left(
E_{J^{\prime }}^{\omega }\left\vert x-z\right\vert \right) ^{2}\leq
\sum_{J^{\prime }\in \mathfrak{C}_{\limfunc{broken}}\left( J\right)
}\left\vert J^{\prime }\right\vert _{\omega }\left( E_{J^{\prime }}^{\omega
}\left\vert x-m_{K}\right\vert \right) ^{2},
\end{equation*}%
and then applying the Carleson embedding theorem again:%
\begin{equation*}
\sum_{J\subset K}\sum_{J^{\prime }\in \mathfrak{C}_{\limfunc{broken}}\left(
J\right) }\left\vert J^{\prime }\right\vert _{\omega }\left( E_{J^{\prime
}}^{\omega }\left\vert x-m_{K}\right\vert \right) ^{2}\lesssim
\int_{K}\left\vert x-m_{K}\right\vert ^{2}d\omega \left( x\right) .
\end{equation*}
\end{proof}

\subsubsection{The smaller Poisson integral}

The expressions $\inf_{z\in \mathbb{R}}\frac{\mathrm{P}_{1+\delta }^{\alpha
}\left( J,\left\vert \mu \right\vert \right) }{\left\vert J\right\vert }%
\left\Vert x-z\right\Vert _{L^{2}\left( \mathbf{1}_{J}\omega \right)
}\left\Vert \square _{J}^{\omega ,\mathbf{b}^{\ast }}\Psi \right\Vert
_{L^{2}\left( \omega \right) }^{\bigstar }$ are typically easier to sum due
to the small Poisson operator $\mathrm{P}_{1+\delta }^{\alpha }\left(
J,\left\vert \mu \right\vert \right) $. To illlustrate, we show here one way
in which we can exploit the additional decay in the Poisson integral $%
\mathrm{P}_{1+\delta }^{\alpha }$. Suppose that $J$ is good in $I$ with $%
\ell \left( J\right) =2^{-s}\ell \left( I\right) $ (see Definition \ref{good
arb} below for `goodness'). We then compute%
\begin{eqnarray*}
\frac{\mathrm{P}_{1+\delta }^{\alpha }\left( J,\mathbf{1}_{A\setminus
I}\sigma \right) }{\left\vert J\right\vert ^{\frac{1}{n}}} &\approx
&\int_{A\setminus I}\frac{\left\vert J\right\vert ^{\frac{\delta }{n}}}{%
\left\vert y-c_{J}\right\vert ^{n+1+\delta -\alpha }}d\sigma \left( y\right)
\\
&\leq &\int_{A\setminus I}\left( \frac{\left\vert J\right\vert ^{\frac{1}{n}}%
}{\limfunc{qdist}\left( c_{J},I^{c}\right) }\right) ^{\delta }\frac{1}{%
\left\vert y-c_{J}\right\vert ^{n+1-\alpha }}d\sigma \left( y\right) \\
&\lesssim &\left( \frac{\left\vert J\right\vert ^{\frac{1}{n}}}{\limfunc{%
qdist}\left( c_{J},I^{c}\right) }\right) ^{\delta }\frac{\mathrm{P}^{\alpha
}\left( J,\mathbf{1}_{A\setminus I}\sigma \right) }{\left\vert J\right\vert
^{\frac{1}{n}}},
\end{eqnarray*}%
and use the goodness inequality,%
\begin{equation*}
d\left( c_{J},I^{c}\right) \geq \frac{1}{2}\ell \left( I\right)
^{1-\varepsilon }\ell \left( J\right) ^{\varepsilon }\geq \frac{1}{2}%
2^{s\left( 1-\varepsilon \right) }\ell \left( J\right) ,
\end{equation*}%
to conclude that%
\begin{equation}
\left( \frac{\mathrm{P}_{1+\delta }^{\alpha }\left( J,\mathbf{1}_{A\setminus
I}\sigma \right) }{\left\vert J\right\vert ^{\frac{1}{n}}}\right) \lesssim
2^{-s\delta \left( 1-\varepsilon \right) }\frac{\mathrm{P}^{\alpha }\left( J,%
\mathbf{1}_{A\setminus I}\sigma \right) }{\left\vert J\right\vert ^{\frac{1}{%
n}}}.  \label{Poisson decay}
\end{equation}%
Now we can estimate%
\begin{eqnarray*}
&&\sum_{J\subset K:\ J\text{ is }\limfunc{good}\text{ in }K}\inf_{z\in 
\mathbb{R}}\frac{\mathrm{P}_{1+\delta }^{\alpha }\left( J,\mathbf{1}%
_{K^{c}}\left\vert \mu \right\vert \right) }{\left\vert J\right\vert }%
\left\Vert x-z\right\Vert _{L^{2}\left( \mathbf{1}_{J}\omega \right)
}\left\Vert \square _{J}^{\omega ,\mathbf{b}^{\ast }}\Psi \right\Vert
_{L^{2}\left( \omega \right) }^{\bigstar } \\
&\leq &\sqrt{\sum_{J\subset K:\ J\text{ is }\limfunc{good}\text{ in }%
K}\left( \frac{\mathrm{P}_{1+\delta }^{\alpha }\left( J,\mathbf{1}%
_{K^{c}}\left\vert \mu \right\vert \right) }{\left\vert J\right\vert }%
\right) ^{2}\inf_{z\in \mathbb{R}}\left\Vert x-z\right\Vert _{L^{2}\left( 
\mathbf{1}_{J}\omega \right) }^{2}}\sqrt{\sum_{J\subset K:\ J\text{ is }%
\limfunc{good}\text{ in }K}\left\Vert \square _{J}^{\omega ,\mathbf{b}^{\ast
}}\Psi \right\Vert _{L^{2}\left( \omega \right) }^{\bigstar 2}},
\end{eqnarray*}%
where%
\begin{eqnarray*}
&&\sum_{J\subset K:\ J\text{ is }\limfunc{good}\text{ in }K}\left( \frac{%
\mathrm{P}_{1+\delta }^{\alpha }\left( J,\mathbf{1}_{K^{c}}\left\vert \mu
\right\vert \right) }{\left\vert J\right\vert }\right) ^{2}\inf_{z\in 
\mathbb{R}}\left\Vert x-z\right\Vert _{L^{2}\left( \mathbf{1}_{J}\omega
\right) }^{2} \\
&=&\sum_{s=0}^{\infty }\sum_{\substack{ J\subset K:\ J\text{ is }\limfunc{%
good}\text{ in }K  \\ \ell \left( J\right) =2^{-s}\ell \left( I\right) }}%
\left( \frac{\mathrm{P}_{1+\delta }^{\alpha }\left( J,\mathbf{1}%
_{K^{c}}\left\vert \mu \right\vert \right) }{\left\vert J\right\vert }%
\right) ^{2}\inf_{z\in \mathbb{R}}\left\Vert x-z\right\Vert _{L^{2}\left( 
\mathbf{1}_{J}\omega \right) }^{2} \\
&\leq &\sum_{s=0}^{\infty }\sum_{\substack{ J\subset K:\ J\text{ is }%
\limfunc{good}\text{ in }K  \\ \ell \left( J\right) =2^{-s}\ell \left(
I\right) }}\left( 2^{-s\delta \left( 1-\varepsilon \right) }\frac{\mathrm{P}%
^{\alpha }\left( J,\mathbf{1}_{K^{c}}\sigma \right) }{\left\vert
J\right\vert ^{\frac{1}{n}}}\right) ^{2}\inf_{z\in \mathbb{R}}\left\Vert
x-z\right\Vert _{L^{2}\left( \mathbf{1}_{J}\omega \right) }^{2} \\
&\leq &\left( \frac{\mathrm{P}^{\alpha }\left( K,\mathbf{1}_{K^{c}}\sigma
\right) }{\left\vert K\right\vert ^{\frac{1}{n}}}\right)
^{2}\sum_{s=0}^{\infty }\sum_{\substack{ J\subset K:\ J\text{ is }\limfunc{%
good}\text{ in }K  \\ \ell \left( J\right) =2^{-s}\ell \left( I\right) }}%
2^{-2s\delta \left( 1-\varepsilon \right) }\inf_{z\in \mathbb{R}}\left\Vert
x-z\right\Vert _{L^{2}\left( \mathbf{1}_{K}\omega \right) }^{2} \\
&\lesssim &\left( \frac{\mathrm{P}^{\alpha }\left( K,\mathbf{1}%
_{K^{c}}\sigma \right) }{\left\vert K\right\vert ^{\frac{1}{n}}}\right)
^{2}\inf_{z\in \mathbb{R}}\left\Vert x-z\right\Vert _{L^{2}\left( \mathbf{1}%
_{K}\omega \right) }^{2}\ ,
\end{eqnarray*}%
and where we have used (\ref{Poisson inequality}), which gives in particular%
\begin{equation*}
\mathrm{P}^{\alpha }(J,\mu \mathbf{1}_{I^{c}})\lesssim \left( \frac{\ell
\left( J\right) }{\ell \left( I\right) }\right) ^{1-\varepsilon \left(
2-\alpha \right) }\mathrm{P}^{\alpha }(I,\mu \mathbf{1}_{I^{c}}).
\end{equation*}%
for $J\subset I$ and $d\left( J,\partial I\right) >\tfrac{1}{2}\ell \left(
J\right) ^{\varepsilon }\ell \left( I\right) ^{1-\varepsilon }$. We will use
such arguments repeatedly in the sequel.

Armed with the Monotonicity Lemma and the lower frame inequality%
\begin{equation*}
\sum_{I\in \mathcal{D}}\left\Vert \square _{I}^{\omega ,\mathbf{b}^{\ast
}}g\right\Vert _{L^{2}\left( \mu \right) }^{\bigstar 2}\lesssim \left\Vert
g\right\Vert _{L^{2}\left( \omega \right) }^{2}\ ,
\end{equation*}%
we can obtain a $\mathbf{b}^{\ast }$-analogue of the Energy Lemma as in \cite%
{SaShUr7} and/or \cite{SaShUr6}.

\subsubsection{The Energy Lemma}

Suppose now we are given a subset $\mathcal{H}$ of the dyadic grid $\mathcal{%
G}$.

\begin{notation}
\label{nonstandard norm}Due to the failure of both martingale and dual
martingale pseudoprojections $\mathsf{Q}_{\mathcal{H}}^{\omega ,\mathbf{b}%
^{\ast }}x$ and $\mathsf{P}_{\mathcal{H}}^{\omega ,\mathbf{b}^{\ast }}g$, as
in Definition \ref{Psi op} in Appendix A, to satisfy inequalities of the
form $\left\Vert \mathsf{P}_{\mathcal{H}}^{\omega ,\mathbf{b}^{\ast
}}g\right\Vert _{L^{2}\left( \omega \right) }\lesssim \left\Vert
g\right\Vert _{L^{2}\left( \omega \right) }$ when the children `break', it
is convenient to define the `square function norms' $\left\Vert \mathsf{Q}_{%
\mathcal{H}}^{\omega ,\mathbf{b}^{\ast }}x\right\Vert _{L^{2}\left( \omega
\right) }^{\spadesuit }$ and $\left\Vert \mathsf{P}_{\mathcal{H}}^{\omega ,%
\mathbf{b}^{\ast }}g\right\Vert _{L^{2}\left( \omega \right) }^{\bigstar }$
of the pseudoprojections 
\begin{equation*}
\mathsf{Q}_{\mathcal{H}}^{\omega ,\mathbf{b}^{\ast }}x=\sum_{J\in \mathcal{H}%
}\bigtriangleup _{J}^{\omega ,\mathbf{b}^{\ast }}x\text{ and }\mathsf{P}_{%
\mathcal{H}}^{\omega ,\mathbf{b}^{\ast }}g=\sum_{J\in \mathcal{H}}\square
_{J}^{\omega ,\mathbf{b}^{\ast }}g\ ,
\end{equation*}%
by%
\begin{eqnarray*}
\left\Vert \mathsf{Q}_{\mathcal{H}}^{\omega ,\mathbf{b}^{\ast }}x\right\Vert
_{L^{2}\left( \omega \right) }^{\spadesuit 2} &\equiv &\sum_{J\in \mathcal{H}%
}\left\Vert \bigtriangleup _{J}^{\omega ,\mathbf{b}^{\ast }}x\right\Vert
_{L^{2}\left( \omega \right) }^{\spadesuit 2}=\sum_{J\in \mathcal{H}%
}\left\Vert \bigtriangleup _{J}^{\omega ,\mathbf{b}^{\ast }}x\right\Vert
_{L^{2}\left( \omega \right) }^{2}+\sum_{J\in \mathcal{H}}\inf_{z\in \mathbb{%
R}}\sum_{J^{\prime }\in \mathfrak{C}_{\limfunc{broken}}\left( J\right)
}\left\vert J^{\prime }\right\vert _{\omega }\left( E_{J^{\prime }}^{\omega
}\left\vert x-z\right\vert \right) ^{2}, \\
\left\Vert \mathsf{P}_{\mathcal{H}}^{\omega ,\mathbf{b}^{\ast }}g\right\Vert
_{L^{2}\left( \omega \right) }^{\bigstar 2} &\equiv &\sum_{J\in \mathcal{H}%
}\left\Vert \square _{J}^{\omega ,\mathbf{b}^{\ast }}g\right\Vert
_{L^{2}\left( \omega \right) }^{\bigstar 2}=\sum_{J\in \mathcal{H}%
}\left\Vert \square _{J}^{\omega ,\mathbf{b}^{\ast }}g\right\Vert
_{L^{2}\left( \omega \right) }^{2}+\sum_{J\in \mathcal{H}}\sum_{J^{\prime
}\in \mathfrak{C}_{\limfunc{broken}}\left( J\right) }\left\vert J^{\prime
}\right\vert _{\omega }\left[ E_{J^{\prime }}^{\omega }\left\vert
g\right\vert +E_{J}^{\omega }\left\vert g\right\vert \right] ^{2},
\end{eqnarray*}%
for any subset $\mathcal{H}\subset \mathcal{G}$. The average $E_{J}^{\omega
}\left\vert x-z\right\vert $ above is taken with respect to the variable $x$%
, i.e. $E_{J}^{\omega }\left\vert x-z\right\vert =\frac{1}{\left\vert
J\right\vert _{\omega }}\int \left\vert x-z\right\vert d\omega \left(
x\right) $, and it is important that the infimum $\inf_{z\in \mathbb{R}}$ is
taken \emph{inside} the sum $\sum_{J\in \mathcal{H}}$.
\end{notation}

Note that we are defining here square function expressions related to
pseudoprojections, which depend not only on the functions $\mathsf{Q}_{%
\mathcal{H}}^{\omega ,\mathbf{b}^{\ast }}x$ and $\mathsf{P}_{\mathcal{H}%
}^{\omega ,\mathbf{b}^{\ast }}g$, but also on the particular representations 
$\sum_{J\in \mathcal{H}}\bigtriangleup _{J}^{\omega ,\mathbf{b}^{\ast }}x$
and $\sum_{J\in \mathcal{H}}\square _{J}^{\omega ,\mathbf{b}^{\ast }}g$.
This slight abuse of notation should not cause confusion, and it provides a
useful way of bookkeeping the sums of squares of norms of martingale and
dual martingale differences $\left\Vert \bigtriangleup _{J}^{\omega ,\mathbf{%
b}^{\ast }}x\right\Vert _{L^{2}\left( \omega \right) }^{2}$ and $\left\Vert
\square _{J}^{\omega ,\mathbf{b}^{\ast }}g\right\Vert _{L^{2}\left( \omega
\right) }^{2}$, along with the norms of the associated Carleson square
function expressions 
\begin{eqnarray*}
\sum_{J\in \mathcal{H}}\inf_{z\in \mathbb{R}}\left\Vert \nabla _{J}^{\omega
}\left( x-z\right) \right\Vert _{L^{2}\left( \omega \right) }^{2}
&=&\sum_{J\in \mathcal{H}}\inf_{z\in \mathbb{R}}\sum_{J^{\prime }\in 
\mathfrak{C}_{\limfunc{broken}}\left( J\right) }\left\vert J^{\prime
}\right\vert _{\omega }\left( E_{J^{\prime }}^{\omega }\left\vert
x-z\right\vert \right) ^{2} \\
\sum_{J\in \mathcal{H}}\left\Vert \nabla _{J}^{\omega }\Psi \right\Vert
_{L^{2}\left( \omega \right) }^{2} &=&\sum_{J\in \mathcal{H}}\sum_{J^{\prime
}\in \mathfrak{C}_{\limfunc{broken}}\left( J\right) }\left\vert J^{\prime
}\right\vert _{\omega }\left[ E_{J^{\prime }}^{\omega }\left\vert \Psi
\right\vert \right] ^{2}.
\end{eqnarray*}%
Note also that the upper weak Riesz inequalities in Appendix A below yield
the inequalities%
\begin{eqnarray*}
\left\Vert \mathsf{Q}_{\mathcal{H}}^{\omega ,\mathbf{b}^{\ast }}x\right\Vert
_{L^{2}\left( \omega \right) }^{2} &\lesssim &\sum_{J\in \mathcal{H}%
}\left\Vert \bigtriangleup _{J}^{\omega ,\mathbf{b}^{\ast }}x\right\Vert
_{L^{2}\left( \omega \right) }^{2}\leq \left\Vert \mathsf{Q}_{\mathcal{H}%
}^{\omega ,\mathbf{b}^{\ast }}x\right\Vert _{L^{2}\left( \omega \right)
}^{\spadesuit 2}\ , \\
\left\Vert \mathsf{P}_{\mathcal{H}}^{\omega ,\mathbf{b}^{\ast }}g\right\Vert
_{L^{2}\left( \omega \right) }^{2} &\lesssim &\sum_{J\in \mathcal{H}%
}\left\Vert \square _{J}^{\omega ,\mathbf{b}^{\ast }}g\right\Vert
_{L^{2}\left( \omega \right) }^{2}\leq \left\Vert \mathsf{P}_{\mathcal{H}%
}^{\omega ,\mathbf{b}^{\ast }}g\right\Vert _{L^{2}\left( \omega \right)
}^{\bigstar 2}\ .
\end{eqnarray*}%
We will exclusively use $\left\Vert \mathsf{Q}_{\mathcal{H}}^{\omega ,%
\mathbf{b}^{\ast }}x\right\Vert _{L^{2}\left( \omega \right) }^{\spadesuit
2} $ in connection with energy terms, and use $\left\Vert \mathsf{P}_{%
\mathcal{H}}^{\sigma ,\mathbf{b}^{\ast }}f\right\Vert _{L^{2}\left( \sigma
\right) }^{\bigstar 2}$ and $\left\Vert \mathsf{P}_{\mathcal{H}}^{\omega ,%
\mathbf{b}^{\ast }}g\right\Vert _{L^{2}\left( \omega \right) }^{\bigstar 2}$%
\ in connection with functions $f\in L^{2}\left( \sigma \right) $ and $g\in
L^{2}\left( \omega \right) $. Finally, note that $\mathsf{Q}_{\mathcal{H}%
}^{\omega ,\mathbf{b}^{\ast }}x=\mathsf{Q}_{\mathcal{H}}^{\omega ,\mathbf{b}%
^{\ast }}\left( x-m\right) $ for any constant $m$. We also define the `large 
$\mathcal{G}$-pseudoprojections'%
\begin{equation*}
\mathsf{Q}_{L}^{\omega ,\mathbf{b}^{\ast }}\equiv \sum_{J^{\prime }\in 
\mathcal{G}:\ J^{\prime }\subset L}\bigtriangleup _{J^{\prime }}^{\omega ,%
\mathbf{b}^{\ast }},\ \ \ \ \ \text{for any interval }L\text{.}
\end{equation*}

Recall that%
\begin{equation*}
\Phi ^{\alpha }\left( J,\nu \right) \equiv \frac{\mathrm{P}^{\alpha }\left(
J,\nu \right) }{\left\vert J\right\vert }\left\Vert \bigtriangleup
_{J}^{\omega ,\mathbf{b}^{\ast }}x\right\Vert _{L^{2}\left( \omega \right)
}^{\spadesuit }+\frac{\mathrm{P}_{1+\delta }^{\alpha }\left( J,\nu \right) }{%
\left\vert J\right\vert }\left\Vert x-m_{J}\right\Vert _{L^{2}\left( \mathbf{%
1}_{J}\omega \right) }\ .
\end{equation*}

\begin{lemma}[\textbf{Energy Lemma}]
\label{ener}Let $J\ $be an interval in $\mathcal{G}$. Let $\Psi _{J}$ be an $%
L^{2}\left( \omega \right) $ function supported in $J$ with vanishing $%
\omega $-mean, and let $\mathcal{H}\subset \mathcal{G}$ be such that $%
J^{\prime }\subset J$ for every $J^{\prime }\in \mathcal{H}$. Let $\nu $ be
a positive measure supported in $\mathbb{R}\setminus \gamma J$ with $\gamma
>1$, and for each $J^{\prime }\in \mathcal{H}$, let $d\nu _{J^{\prime
}}=\varphi _{J^{\prime }}d\nu $ with $\left\vert \varphi _{J^{\prime
}}\right\vert \leq 1$. Suppose that $\mathbf{b}^{\ast }$ is an $\infty $%
-weakly $\mu $-controlled accretive family on $\mathbb{R}$. Let $T^{\alpha }$
be a standard $\alpha $-fractional singular integral operator with $0\leq
\alpha <1$. Then we have%
\begin{eqnarray*}
&&\left\vert \sum_{J^{\prime }\in \mathcal{H}}\left\langle T^{\alpha }\left(
\nu _{J^{\prime }}\right) ,\square _{J^{\prime }}^{\omega ,\mathbf{b}^{\ast
}}\Psi _{J}\right\rangle _{\omega }\right\vert \lesssim C_{\gamma
}\sum_{J^{\prime }\in \mathcal{H}}\Phi ^{\alpha }\left( J^{\prime },\nu
\right) \left\Vert \square _{J^{\prime }}^{\omega ,\mathbf{b}^{\ast }}\Psi
_{J}\right\Vert _{L^{2}\left( \mu \right) }^{\bigstar } \\
&\lesssim &C_{\gamma }\sqrt{\sum_{J^{\prime }\in \mathcal{H}}\Phi ^{\alpha
}\left( J^{\prime },\nu \right) ^{2}}\sqrt{\sum_{J^{\prime }\in \mathcal{H}%
}\left\Vert \square _{J^{\prime }}^{\omega ,\mathbf{b}^{\ast }}\Psi
_{J}\right\Vert _{L^{2}\left( \mu \right) }^{\bigstar 2}} \\
&\leq &C_{\gamma }\left( \frac{\mathrm{P}^{\alpha }\left( J,\nu \right) }{%
\left\vert J\right\vert }\left\Vert \mathsf{Q}_{\mathcal{H}}^{\omega ,%
\mathbf{b}^{\ast }}x\right\Vert _{L^{2}\left( \omega \right) }^{\spadesuit }+%
\frac{\mathrm{P}_{1+\delta }^{\alpha }\left( J,\nu \right) }{\left\vert
J\right\vert }\left\Vert x-m_{J}\right\Vert _{L^{2}\left( \mathbf{1}%
_{J}\omega \right) }\right) \left\Vert \mathsf{P}_{\mathcal{H}}^{\omega ,%
\mathbf{b}^{\ast }}\Psi _{J}\right\Vert _{L^{2}\left( \mu \right)
}^{\bigstar },
\end{eqnarray*}%
and in particular the `energy' estimate%
\begin{equation*}
\left\vert \left\langle T^{\alpha }\varphi \nu ,\Psi _{J}\right\rangle
_{\omega }\right\vert \lesssim C_{\gamma }\left( \frac{\mathrm{P}^{\alpha
}\left( J,\nu \right) }{\left\vert J\right\vert }\left\Vert \mathsf{Q}%
_{J}^{\omega ,\mathbf{b}^{\ast }}x\right\Vert _{L^{2}\left( \omega \right)
}^{\spadesuit }+\frac{\mathrm{P}_{1+\delta }^{\alpha }\left( J,\nu \right) }{%
\left\vert J\right\vert }\left\Vert x-m_{J}\right\Vert _{L^{2}\left( \mathbf{%
1}_{J}\omega \right) }\right) \left\Vert \sum_{J^{\prime }\subset J}\square
_{J^{\prime }}^{\omega ,\mathbf{b}^{\ast }}\Psi _{J}\right\Vert
_{L^{2}\left( \mu \right) }^{\bigstar }\ ,
\end{equation*}%
where $\left\Vert \sum_{J^{\prime }\subset J}\square _{J^{\prime }}^{\omega ,%
\mathbf{b}^{\ast }}\Psi _{J}\right\Vert _{L^{2}\left( \mu \right)
}^{\bigstar }\lesssim \left\Vert \Psi _{J}\right\Vert _{L^{2}\left( \mu
\right) }$, and the `pivotal' bound%
\begin{equation*}
\left\vert \left\langle T^{\alpha }\left( \varphi \nu \right) ,\Psi
_{J}\right\rangle _{\omega }\right\vert \lesssim C_{\gamma }\mathrm{P}%
^{\alpha }\left( J,\left\vert \nu \right\vert \right) \sqrt{\left\vert
J\right\vert _{\omega }}\left\Vert \Psi _{J}\right\Vert _{L^{2}\left( \omega
\right) }\ ,
\end{equation*}%
for any function $\varphi $ with $\left\vert \varphi \right\vert \leq 1$.
\end{lemma}

\begin{proof}
Beginning with the first display in the conclusion of the Energy Lemma, we
need only prove the first line since the next two follow from Poisson
inequalities and the definitions in Notation \ref{nonstandard norm}. Using
the Monotonicity Lemma \ref{mono}, followed by $\left\vert \nu _{J^{\prime
}}\right\vert \leq \nu $, the Poisson equivalence 
\begin{equation}
\frac{\mathrm{P}^{\alpha }\left( J^{\prime },\nu \right) }{\left\vert
J^{\prime }\right\vert }\approx \frac{\mathrm{P}^{\alpha }\left( J,\nu
\right) }{\left\vert J\right\vert },\ \ \ \ \ J^{\prime }\subset J\subset
\gamma J,\ \ \ \limfunc{supp}\nu \cap \gamma J=\emptyset ,
\label{Poisson equiv}
\end{equation}%
and the weak frame inequalities for dual martingale differences in Appendix
A, we have%
\begin{eqnarray*}
&&\left\vert \sum_{J^{\prime }\in \mathcal{H}}\left\langle T^{\alpha }\left(
\nu _{J^{\prime }}\right) ,\square _{J^{\prime }}^{\omega ,\mathbf{b}^{\ast
}}\Psi _{J}\right\rangle _{\omega }\right\vert \lesssim \sum_{J^{\prime }\in 
\mathcal{H}}\Phi ^{\alpha }\left( J^{\prime },\left\vert \mu \right\vert
\right) \left\Vert \square _{J^{\prime }}^{\omega ,\mathbf{b}^{\ast }}\Psi
_{J}\right\Vert _{L^{2}\left( \mu \right) }^{\bigstar } \\
&\leq &\sum_{J^{\prime }\in \mathcal{H}}\left( \frac{\mathrm{P}^{\alpha
}\left( J^{\prime },\nu \right) }{\left\vert J^{\prime }\right\vert }%
\left\Vert \bigtriangleup _{J^{\prime }}^{\omega ,\mathbf{b}^{\ast
}}x\right\Vert _{L^{2}\left( \omega \right) ,z}^{\spadesuit }+\frac{\mathrm{P%
}_{1+\delta }^{\alpha }\left( J^{\prime },\nu \right) }{\left\vert J^{\prime
}\right\vert }\left\Vert x-m_{J^{\prime }}\right\Vert _{L^{2}\left( \mathbf{1%
}_{J^{\prime }}\omega \right) }\right) \left\Vert \square _{J^{\prime
}}^{\omega ,\mathbf{b}^{\ast }}\Psi _{J}\right\Vert _{L^{2}\left( \omega
\right) }^{\bigstar } \\
&\lesssim &\left( \sum_{J^{\prime }\in \mathcal{H}}\left( \frac{\mathrm{P}%
^{\alpha }\left( J^{\prime },\nu \right) }{\left\vert J^{\prime }\right\vert 
}\right) ^{2}\left\Vert \bigtriangleup _{J^{\prime }}^{\omega ,\mathbf{b}%
^{\ast }}x\right\Vert _{L^{2}\left( \omega \right) ,z}^{\spadesuit 2}\right)
^{\frac{1}{2}}\left( \sum_{J^{\prime }\in \mathcal{H}}\left\Vert \square
_{J^{\prime }}^{\omega ,\mathbf{b}^{\ast }}\Psi _{J}\right\Vert
_{L^{2}\left( \omega \right) }^{\bigstar 2}\right) ^{\frac{1}{2}} \\
&&+\left( \sum_{J^{\prime }\in \mathcal{H}}\left( \frac{\mathrm{P}_{1+\delta
}^{\alpha }\left( J^{\prime },\left\vert \mu \right\vert \right) }{%
\left\vert J^{\prime }\right\vert }\right) ^{2}\left\Vert x-m_{J^{\prime
}}\right\Vert _{L^{2}\left( \mathbf{1}_{J^{\prime }}\omega \right)
}^{2}\right) ^{\frac{1}{2}}\left( \sum_{J^{\prime }\in \mathcal{H}%
}\left\Vert \square _{J^{\prime }}^{\omega ,\mathbf{b}^{\ast }}\Psi
_{J}\right\Vert _{L^{2}\left( \omega \right) }^{\bigstar 2}\right) ^{\frac{1%
}{2}} \\
&\lesssim &\left( \frac{\mathrm{P}^{\alpha }\left( J,\nu \right) }{%
\left\vert J\right\vert }\right) \left\Vert \mathsf{Q}_{\mathcal{H}}^{\omega
,\mathbf{b}^{\ast }}x\right\Vert _{L^{2}\left( \omega \right) }^{\spadesuit
}\left\Vert \Psi _{J}\right\Vert _{L^{2}\left( \omega \right) }+\frac{1}{%
\gamma ^{\delta ^{\prime }}}\left( \frac{\mathrm{P}_{1+\delta ^{\prime
}}^{\alpha }\left( J,\nu \right) }{\left\vert J\right\vert }\right)
\left\Vert x-m_{J}\right\Vert _{L^{2}\left( \mathbf{1}_{J}\omega \right)
}\left\Vert \Psi _{J}\right\Vert _{L^{2}\left( \omega \right) }\ .
\end{eqnarray*}%
The last inequality follows from the following calculation using Haar
projections $\bigtriangleup _{K}^{\omega }$:

\begin{eqnarray}
&&\sum_{J^{\prime }\in \mathcal{H}}\left( \frac{\mathrm{P}_{1+\delta
}^{\alpha }\left( J^{\prime },\nu \right) }{\left\vert J^{\prime
}\right\vert }\right) ^{2}\left\Vert x-m_{J^{\prime }}\right\Vert
_{L^{2}\left( \mathbf{1}_{J^{\prime }}\omega \right) }^{2}
\label{Haar trick} \\
&=&\sum_{J^{\prime }\in \mathcal{H}}\left( \frac{\mathrm{P}_{1+\delta
}^{\alpha }\left( J^{\prime },\nu \right) }{\left\vert J^{\prime
}\right\vert }\right) ^{2}\sum_{J^{\prime \prime }\subset J^{\prime
}}\left\Vert \bigtriangleup _{J^{\prime \prime }}^{\omega }x\right\Vert
_{L^{2}\left( \omega \right) }^{2}=\sum_{J^{\prime \prime }\subset J}\left\{
\sum_{J^{\prime }:\ J^{\prime \prime }\subset J^{\prime }\subset J}\left( 
\frac{\mathrm{P}_{1+\delta }^{\alpha }\left( J^{\prime },\nu \right) }{%
\left\vert J^{\prime }\right\vert }\right) ^{2}\right\} \left\Vert
\bigtriangleup _{J^{\prime \prime }}^{\omega }x\right\Vert _{L^{2}\left(
\omega \right) }^{2}  \notag \\
&\lesssim &\frac{1}{\gamma ^{2\delta ^{\prime }}}\sum_{J^{\prime \prime
}\subset J}\left( \frac{\mathrm{P}_{1+\delta ^{\prime }}^{\alpha }\left(
J^{\prime \prime },\nu \right) }{\left\vert J^{\prime \prime }\right\vert }%
\right) ^{2}\left\Vert \bigtriangleup _{J^{\prime \prime }}^{\omega
}x\right\Vert _{L^{2}\left( \omega \right) }^{2}\leq \frac{1}{\gamma
^{2\delta ^{\prime }}}\left( \frac{\mathrm{P}_{1+\delta ^{\prime }}^{\alpha
}\left( J,\nu \right) }{\left\vert J\right\vert }\right) ^{2}\sum_{J^{\prime
\prime }\subset J}\left\Vert \bigtriangleup _{J^{\prime \prime }}^{\omega
}x\right\Vert _{L^{2}\left( \omega \right) }^{2}\ ,  \notag
\end{eqnarray}%
which in turn follows from (recalling $\delta =2\delta ^{\prime }$ and using 
$\left\vert J^{\prime }\right\vert +\left\vert y-c_{J^{\prime }}\right\vert
\approx \left\vert J\right\vert +\left\vert y-c_{J}\right\vert $ and $\frac{%
\left\vert J\right\vert }{\left\vert J\right\vert +\left\vert
y-c_{J}\right\vert }\leq \frac{1}{\gamma }$ for $y\in \mathbb{R}\setminus
\gamma J$)%
\begin{eqnarray*}
&&\sum_{J^{\prime }:\ J^{\prime \prime }\subset J^{\prime }\subset J}\left( 
\frac{\mathrm{P}_{1+\delta }^{\alpha }\left( J^{\prime },\nu \right) }{%
\left\vert J^{\prime }\right\vert }\right) ^{2}=\sum_{J^{\prime }:\
J^{\prime \prime }\subset J^{\prime }\subset J}\left\vert J^{\prime
}\right\vert ^{2\delta }\left( \int_{\mathbb{R}\setminus \gamma J}\frac{1}{%
\left( \left\vert J^{\prime }\right\vert +\left\vert y-c_{J^{\prime
}}\right\vert \right) ^{2+\delta -\alpha }}d\nu \left( y\right) \right) ^{2}
\\
&\lesssim &\sum_{J^{\prime }:\ J^{\prime \prime }\subset J^{\prime }\subset
J}\frac{1}{\gamma ^{2\delta ^{\prime }}}\frac{\left\vert J^{\prime
}\right\vert ^{2\delta }}{\left\vert J\right\vert ^{2\delta }}\left( \int_{%
\mathbb{R}\setminus \gamma J}\frac{\left\vert J\right\vert ^{\delta ^{\prime
}}}{\left( \left\vert J\right\vert +\left\vert y-c_{J}\right\vert \right)
^{2+\delta ^{\prime }-\alpha }}d\nu \left( y\right) \right) ^{2} \\
&=&\frac{1}{\gamma ^{2\delta ^{\prime }}}\left( \sum_{J^{\prime }:\
J^{\prime \prime }\subset J^{\prime }\subset J}\frac{\left\vert J^{\prime
}\right\vert ^{2\delta }}{\left\vert J\right\vert ^{2\delta }}\right) \left( 
\frac{\mathrm{P}_{1+\delta ^{\prime }}^{\alpha }\left( J,\nu \right) }{%
\left\vert J\right\vert }\right) ^{2}\lesssim \frac{1}{\gamma ^{2\delta
^{\prime }}}\left( \frac{\mathrm{P}_{1+\delta ^{\prime }}^{\alpha }\left(
J,\nu \right) }{\left\vert J\right\vert }\right) ^{2}.
\end{eqnarray*}%
Finally we obtain the `energy' estimate from the equality%
\begin{equation*}
\Psi _{J}=\sum_{J^{\prime }\subset J}\square _{J^{\prime }}^{\omega ,\mathbf{%
b}^{\ast }}\Psi _{J}\ ,\ \ \ \ \ (\text{since }\Psi _{J}\text{ has vanishing 
}\omega \text{-mean)},\text{ }
\end{equation*}%
and we obtain the `pivotal' bound from the inequality%
\begin{equation*}
\sum_{J^{\prime \prime }\subset J}\left\Vert \bigtriangleup _{J^{\prime
\prime }}^{\omega ,\mathbf{b}^{\ast }}x\right\Vert _{L^{2}\left( \omega
\right) }^{\spadesuit 2}\lesssim \left\Vert \left( x-m_{J}\right)
\right\Vert _{L^{2}\left( \mathbf{1}_{J}\omega \right) }^{2}\leq \left\vert
J\right\vert ^{2}\left\vert J\right\vert _{\omega }\ .
\end{equation*}
\end{proof}

\subsection{Organization of the proof}

We adapt the proof of the main theorem in \cite{SaShUr7}, \cite{SaShUr9} and 
\cite{SaShUr10}, but beginning instead with the decomposition of Hyt\"{o}nen
and Martikainen \cite{HyMa}, to obtain the norm inequality%
\begin{equation*}
\mathfrak{N}_{T^{\alpha }}\lesssim \mathcal{NTV}_{\alpha }=\mathfrak{T}%
_{T^{\alpha }}^{\mathbf{b}}+\mathfrak{T}_{T^{\alpha }}^{\mathbf{b}^{\ast }}+%
\sqrt{\mathfrak{A}_{2}^{\alpha }}+\mathfrak{E}_{2}^{\alpha },
\end{equation*}%
under the \emph{apriori} assumption $\mathfrak{N}_{T^{\alpha }}<\infty $,
which is achieved by considering one of the truncations $T_{\sigma ,\delta
,R}^{\alpha }$ defined in (\ref{def truncation}) above. This will be carried
out in the next five sections of this paper. In the next section we consider
the various form splittings and reduce matters to the \emph{disjoint} form,
the \emph{nearby} form and the \emph{main below} form. Then these latter
three forms are taken up in the subsequent three sections, using material
from the appendices. Finally, the stopping form is treated in the section
following these three.

A major source of difficulty will arise in the infusion of goodness for the
intervals $J$ into the main below form where the sum is taken over all pairs 
$\left( I,J\right) $ such that $\ell \left( J\right) \leq \ell \left(
I\right) $. We will infuse goodness in a weak way pioneered by Hyt\"{o}nen
and Martikainen in a one weight setting. This weak form of goodness is then
exploited in all subsequent constructions by typically replacing $J$ with $%
J^{\maltese }$ in defining relations, where $J^{\maltese }$ is the smallest
interval $K$ for which $J$ is good in $K$ and beyond (see the next section
for terminology, in particular Definition \ref{def sharp cross}).

Another source of difficulty arises in the treatment of the nearby form in
the setting of two weights. The one weight proofs in \cite{HyMa} and \cite%
{LaMa} relied strongly on a property peculiar to the one weight setting -
namely the fact already pointed out in Remark \ref{special}\ above, that
both of the Poisson integrals are bounded, namely $\mathrm{P}^{\alpha
}\left( Q,\mu \right) \lesssim 1$ and $\mathcal{P}^{\alpha }\left( Q,\mu
\right) \lesssim 1$. We will circumvent this difficulty by combining a
recursive energy argument with the full testing conditions assumed for the $%
\infty $-weakly accretive family of \emph{original} testing functions $%
b_{Q}^{\limfunc{orig}}$, before these conditions were suppressed by corona
constructions that delivered only weak testing conditions for the new family
of testing functions $b_{Q}$.

In Section \ref{Sec stop} we bound the stopping form using the arguments
from \cite{SaShUr7}, \cite{SaShUr9} and \cite{SaShUr10}, which were in turn
based on the bottom/up stopping time and recursion of M. Lacey in \cite{Lac}%
. Here we introduce an additional top/down `indented' corona construction to
handle the lack of goodness in size testing intervals, and we use an
absorption in place of recursion. Finally, the treatment of various
`straddling lemmas' is complicated by weak goodness, and we use a stronger
form of weak goodness defined with the three point `skeleton' of an interval
replaced by an infinite `body', coupled with two geometric `Key Facts' to
establish these lemmas.

Of particular importance will be two independent results proved in
Appendices A and B that follow from known work with some new twists. In
Appendix A we establish convergence and frame inequalities for martingale
and dual martingale differences, and derive certain weak Riesz inequalities
associated with $\infty $-weakly $\mu $-accretive families of testing
functions, which will find application in treating the paraproduct form
below. The boundedness of testing functions, and the reverse H\"{o}lder
condition on their children, is important here.

In Appendix B we show that the functional energy for an arbitrary pair of
grids is controlled by the Muckenhoupt and energy side conditions. The
somewhat lengthy proof of this latter assertion is similar to the
corresponding proof in the $T1$ setting - see e.g. \cite{SaShUr9} - but
requires a different decomposition of the stopping intervals into `Whitney
intervals' in order to accommodate the weaker notion of goodness used here,
as well as the usual decomposition into maximal deeply embedded intervals
that is used to control expressions involving the `small' Poisson integral.

Finally, we include in Appendix C an up-to-date list of errata for our most
often referred to paper \cite{SaShUr7}.

\section{Form splittings}

\begin{notation}
Fix grids $\mathcal{D}$ and $\mathcal{G}$. We will use $\mathcal{D}$ to
denote the grid associated with $f\in L^{2}\left( \sigma \right) $, and we
will use $\mathcal{G}$ to denote the grid associated with $g\in L^{2}\left(
\omega \right) $.
\end{notation}

We have defined corona decompositions of $f$ and $g$ in the $\sigma $%
-iterated triple corona construction above, but in order to start these
corona decompositions for $f$ and $g$ respectively within the dyadic grids $%
\mathcal{D}$ and $\mathcal{G}$, we need to first restrict $f$ and $g$ to be
supported in a large common interval $Q_{\infty }$. Then we cover $Q_{\infty
}$ with two pairwise disjoint intervals $I_{\infty }\in \mathcal{D}$ with $%
\ell \left( I_{\infty }\right) =\ell \left( Q_{\infty }\right) $, and
similarly cover $Q_{\infty }$ with two pairwise disjoint intervals $%
J_{\infty }\in \mathcal{G}$ with $\ell \left( J_{\infty }\right) =\ell
\left( Q_{\infty }\right) $. We can now use the broken martingale
decompositions from Appendix A, together with full $\mathbf{b}$-testing (see
(\ref{full b testing}) and (\ref{full proved}) below), to reduce matters to
consideration of the four forms%
\begin{equation*}
\sum_{I\in \mathcal{D}:\ I\subset I_{\infty }}\sum_{J\in \mathcal{G}:\
J\subset J_{\infty }}\int \left( T_{\sigma }^{\alpha }\square _{I}^{\sigma ,%
\mathbf{b}}f\right) \square _{J}^{\omega ,\mathbf{b}^{\ast }}gd\omega ,
\end{equation*}%
with $I_{\infty }$ and $J_{\infty }$ as above, and where we can then use the
intervals $I_{\infty }$ and $J_{\infty }$ as the starting intervals in our
corona constructions below. Indeed, the identities in Lemma \ref{conv prop}
from Appendix A below, give%
\begin{eqnarray*}
f &=&\sum_{I\in \mathcal{D}_{N}}\mathbb{F}_{I}^{\sigma ,\mathbf{b}%
}f+\sum_{I\in \mathcal{D}:\ I\subset I_{\infty },\ \ell \left( I\right) \geq
N+1}\square _{I}^{\sigma ,\mathbf{b}}f+\mathbb{F}_{I_{\infty }}^{\sigma ,%
\mathbf{b}}f, \\
g &=&\sum_{J\in \mathcal{G}_{N}}\mathbb{F}_{J}^{\omega ,\mathbf{b}^{\ast
}}g+\sum_{J\in \mathcal{G}:\ J\subset J_{\infty },\ \ell \left( J\right)
\geq N+1}\square _{J}^{\omega ,\mathbf{b}^{\ast }}g+\mathbb{F}_{J_{\infty
}}^{\omega ,\mathbf{b}^{\ast }}g,
\end{eqnarray*}%
which can then be used to write the bilinear form $\int \left( T_{\sigma
}f\right) gd\omega $ as a sum of the forms%
\begin{eqnarray}
&&\ \ \ \ \ \ \ \ \ \ \ \ \ \ \ \int \left( T_{\sigma }f\right) gd\omega
=\sum_{\substack{ \text{four pairs}  \\ \left( I_{\infty },J_{\infty
}\right) }}\left\{ \sum_{I\in \mathcal{D}:\ I\subset I_{\infty }}\sum_{J\in 
\mathcal{G}:\ J\subset J_{\infty }}\int \left( T_{\sigma }^{\alpha }\square
_{I}^{\sigma ,\mathbf{b}}f\right) \square _{J}^{\omega ,\mathbf{b}^{\ast
}}gd\omega \right.  \label{sum of forms} \\
&&\left. +\sum_{I\in \mathcal{D}:\ I\subset I_{\infty }}\int \left(
T_{\sigma }^{\alpha }\square _{I}^{\sigma ,\mathbf{b}}f\right) \mathbb{F}%
_{J_{\infty }}^{\omega ,\mathbf{b}^{\ast }}gd\omega +\sum_{J\in \mathcal{G}%
:\ J\subset J_{\infty }}\int \left( T_{\sigma }^{\alpha }\mathbb{F}%
_{I_{\infty }}^{\sigma ,\mathbf{b}}f\right) \square _{J}^{\omega ,\mathbf{b}%
^{\ast }}gd\omega +\int \left( T_{\sigma }^{\alpha }\mathbb{F}_{I_{\infty
}}^{\sigma ,\mathbf{b}}f\right) \mathbb{F}_{J_{\infty }}^{\omega ,\mathbf{b}%
^{\ast }}gd\omega \right\} ,  \notag
\end{eqnarray}%
taken over the four pairs of intervals $\left( I_{\infty },J_{\infty
}\right) $ above, plus the limit of the sum of terms involving $\sum_{I\in 
\mathcal{D}_{N}}\mathbb{F}_{I}^{\sigma ,\mathbf{b}}f$ and $\sum_{J\in 
\mathcal{G}_{N}}\mathbb{F}_{J}^{\omega ,\mathbf{b}^{\ast }}g$. This latter
limit is easily shown to vanish due to the strong convergence of the dual
martingale differences $\square _{I}^{\sigma ,\mathbf{b}}f$ and $\square
_{J}^{\omega ,\mathbf{b}^{\ast }}g$ in $L^{2}\left( \sigma \right) $ and $%
L^{2}\left( \omega \right) $ respectively. More precisely, we have%
\begin{eqnarray*}
&&\left\vert \int \left( T_{\sigma }^{\alpha }\sum_{I\in \mathcal{D}_{N}}%
\mathbb{F}_{I}^{\sigma ,\mathbf{b}}f\right) \sum_{J\in \mathcal{G}_{N}}%
\mathbb{F}_{J}^{\omega ,\mathbf{b}^{\ast }}g\ d\omega \right\vert \lesssim 
\mathfrak{N}_{T^{\alpha }}\left\Vert \sum_{I\in \mathcal{D}_{N}}\mathbb{F}%
_{I}^{\sigma ,\mathbf{b}}f\right\Vert _{L^{2}\left( \sigma \right)
}\left\Vert \sum_{J\in \mathcal{G}_{N}}\mathbb{F}_{J}^{\omega ,\mathbf{b}%
^{\ast }}g\right\Vert _{L^{2}\left( \omega \right) } \\
&=&\mathfrak{N}_{T^{\alpha }}\left\Vert \sum_{I\in \mathcal{D}:\ \ell \left(
I\right) \geq 2^{N}}\square _{I}^{\sigma ,\mathbf{b}}f\right\Vert
_{L^{2}\left( \sigma \right) }\left\Vert \sum_{J\in \mathcal{G}:\ \ell
\left( J\right) \geq 2^{N}}\square _{J}^{\omega ,\mathbf{b}^{\ast
}}g\right\Vert _{L^{2}\left( \omega \right) } \\
&\lesssim &\mathfrak{N}_{T^{\alpha }}\left( \sum_{I\in \mathcal{D}:\ \ell
\left( I\right) \geq 2^{N}}\left\Vert \square _{I}^{\sigma ,\mathbf{b}%
}f\right\Vert _{L^{2}\left( \sigma \right) }^{2}\right) ^{\frac{1}{2}}\left(
\sum_{J\in \mathcal{G}:\ \ell \left( J\right) \geq 2^{N}}\left\Vert \square
_{J}^{\omega ,\mathbf{b}^{\ast }}g\right\Vert _{L^{2}\left( \omega \right)
}^{2}\right) ^{\frac{1}{2}},
\end{eqnarray*}%
which tends to $0$ as $N\rightarrow \infty $ since%
\begin{equation*}
\left( \sum_{I\in \mathcal{D}}\left\Vert \square _{I}^{\sigma ,\mathbf{b}%
}f\right\Vert _{L^{2}\left( \sigma \right) }^{2}\right) ^{\frac{1}{2}}\left(
\sum_{J\in \mathcal{G}}\left\Vert \square _{J}^{\omega ,\mathbf{b}^{\ast
}}g\right\Vert _{L^{2}\left( \omega \right) }^{2}\right) ^{\frac{1}{2}%
}\lesssim \left\Vert f\right\Vert _{L^{2}\left( \sigma \right) }\left\Vert
g\right\Vert _{L^{2}\left( \omega \right) }\ .
\end{equation*}

\begin{remark}
In particular, 
\begin{equation*}
\lim_{N\rightarrow \infty }\sup_{I\in \mathcal{D}_{N}}\left\Vert \mathbb{F}%
_{I}^{\sigma ,\mathbf{b}}f\right\Vert _{L^{2}\left( \sigma \right)
}=0=\lim_{N\rightarrow \infty }\sup_{J\in \mathcal{G}_{N}}\left\Vert \mathbb{%
F}_{J}^{\omega ,\mathbf{b}^{\ast }}g\right\Vert _{L^{2}\left( \sigma \right)
}
\end{equation*}%
and so we can use the pointwise telescoping identities%
\begin{equation*}
\mathbb{F}_{I}^{\sigma ,\mathbf{b}}f\left( x\right) =\sum_{I^{\prime }\in 
\mathcal{D}:\ I\subset I^{\prime }}\square _{I^{\prime }}^{\sigma ,\mathbf{b}%
}f\left( x\right) \text{ and }\mathbb{F}_{J}^{\omega ,\mathbf{b}^{\ast
}}g\left( x\right) =\sum_{J\in \mathcal{G}:\ J\subset J^{\prime }}\square
_{J^{\prime }}^{\omega ,\mathbf{b}^{\ast }}g\left( x\right) .
\end{equation*}
\end{remark}

The second, third and fourth sums in (\ref{sum of forms}) can be controlled
by the full testing conditions, e.g.%
\begin{eqnarray}
&&\left\vert \sum_{I\in \mathcal{D}:\ I\subset I_{\infty }}\int \left(
T_{\sigma }^{\alpha }\square _{I}^{\sigma ,\mathbf{b}}f\right) \mathbb{F}%
_{J_{\infty }}^{\omega ,\mathbf{b}^{\ast }}gd\omega \right\vert =\left\vert
\int \left( \sum_{I\in \mathcal{D}:\ I\subset I_{\infty }}\square
_{I}^{\sigma ,\mathbf{b}}f\right) T_{\omega }^{\alpha ,\ast }\left( \mathbb{F%
}_{J_{\infty }}^{\omega ,\mathbf{b}^{\ast }}g\right) d\sigma \right\vert
\label{top control} \\
&\leq &\left\Vert \sum_{I\in \mathcal{D}:\ I\subset I_{\infty }}\square
_{I}^{\sigma ,\mathbf{b}}f\right\Vert _{L^{2}\left( \sigma \right)
}\left\Vert T_{\omega }^{\alpha ,\ast }\left( \mathbb{F}_{J_{\infty
}}^{\omega ,\mathbf{b}^{\ast }}g\right) \right\Vert _{L^{2}\left( \sigma
\right) }\lesssim \left\Vert f\right\Vert _{L^{2}\left( \sigma \right)
}\left( \mathfrak{FT}_{T_{\omega }^{\alpha ,\ast }}+\sqrt{\mathfrak{A}%
_{2}^{\alpha }}\right) \left\Vert g\right\Vert _{L^{2}\left( \omega \right) }
\notag
\end{eqnarray}%
since $\mathbb{F}_{J_{\infty }}^{\omega ,\mathbf{b}^{\ast }}g=b_{J_{\infty
}}^{\ast }\frac{E_{J_{\infty }}^{\omega }g}{E_{J_{\infty }}^{\omega
}b_{J_{\infty }}^{\ast }}$ is $b_{J_{\infty }}^{\ast }$ times an `accretive'
average of $g$ on $J_{\infty }$, and similarly for the third and fourth sum.
Finally, the full testing conditions $\mathfrak{FT}_{T_{\sigma }^{\alpha }}^{%
\mathbf{b}}$ and $\mathfrak{FT}_{T_{\omega }^{\alpha ,\ast }}^{\mathbf{b}%
^{\ast }}$ are controlled by the usual testing conditions $\mathfrak{T}%
_{T_{\sigma }^{\alpha }}^{\mathbf{b}}$ and $\mathfrak{T}_{T_{\omega
}^{\alpha ,\ast }}^{\mathbf{b}^{\ast }}$ together with the Muckenhoupt
condition $\mathfrak{A}_{2}^{\alpha }$, by (\ref{full proved}) below.

\begin{description}
\item[Important] In the $\sigma $-iterated triple corona construction we
redefined the family $\mathbf{b}=\left\{ b_{Q}\right\} _{Q\in \mathcal{D}}$
so that the new functions $b_{Q}^{\limfunc{new}}$ are given in terms of the
original functions $b_{Q}^{\limfunc{orig}}$ by $b_{Q}^{\limfunc{new}}=%
\mathbf{1}_{Q}b_{A}^{\limfunc{orig}}$ for $Q\in \mathcal{C}_{A}^{\sigma }$,
and of course we then dropped the superscript `$\limfunc{new}$'. We continue
to refer to the triple stopping intervals $A$ as `breaking' intervals even
if $b_{A}$ happens to equal $\mathbf{1}_{A}b_{\pi A}$. The results of
Appendix A apply with this more inclusive definition of `breaking'
intervals, and the associated definition of `broken' children, since only
the Carleson condition\ on stopping intervals is relevant here.
\end{description}

Altogether this and Proposition \ref{data} give us the \emph{triple corona
decomposition} of $f=\sum_{A\in \mathcal{A}}\mathsf{P}_{\mathcal{C}%
_{A}}^{\sigma ,\mathbf{b}}f$, where the pseudoprojection $\mathsf{P}_{%
\mathcal{C}_{A}}^{\sigma ,\mathbf{b}}$ is defined in Appendix A:%
\begin{equation*}
\mathsf{P}_{\mathcal{C}_{A}}^{\mu ,\mathbf{b}}f=\sum_{I\in \mathcal{C}%
_{A}}\square _{I}^{\mu ,\mathbf{b}}f\ .
\end{equation*}%
We now record the main facts proved above, and in Appendix A below, for the
triple corona.

\begin{lemma}
Let $f\in L^{2}\left( \sigma \right) $. We have%
\begin{eqnarray*}
f &=&\sum_{A\in \mathcal{A}}\mathsf{P}_{\mathcal{C}_{A}}^{\sigma }f,\ \ \ \
\ \text{both in the sense of norm convergence in }L^{2}\left( \sigma \right)
\\
&&\ \ \ \ \ \ \ \ \ \ \ \ \ \ \ \ \ \ \ \text{ and pointwise }\sigma \text{%
-a.e}.
\end{eqnarray*}%
The corona tops $\mathcal{A}$ and stopping bounds $\left\{ \alpha _{\mathcal{%
A}}\left( A\right) \right\} _{A\in \mathcal{A}}$ satisfy properties (1),
(2), (3) and (4) in Definition \ref{general stopping data}, hence constitute
stopping data for $f$. Moreover, $\mathbf{b}=\left\{ b_{I}\right\} _{I\in 
\mathcal{D}}$ is an $\infty $-strongly $\sigma $-controlled accretive family
on $\mathcal{D}$ with corona tops $\mathcal{A\subset D}$, where $b_{I}=%
\mathbf{1}_{I}b_{A}$ has reverse H\"{o}lder control on children for all $%
I\in \mathcal{C}_{A}$, and the weak corona forward testing condition holds
uniformly in coronas, i.e.%
\begin{equation*}
\frac{1}{\left\vert I\right\vert _{\sigma }}\int_{I}\left\vert T_{\sigma
}^{\alpha }b_{A}\right\vert ^{2}d\sigma \leq C,\ \ \ \ \ I\in \mathcal{C}%
_{A}^{\sigma }\ .
\end{equation*}%
Similar statements hold for $g\in L^{2}\left( \omega \right) $.
\end{lemma}

Now we turn to the various splittings of forms, beginning with the two
weight analogue of the decomposition of Hyt\"{o}nen and Martikainen \cite%
{HyMa}. Fix the stopping data $\mathcal{A}$ and $\left\{ \alpha _{\mathcal{A}%
}\left( A\right) \right\} _{A\in \mathcal{A}}$ and dual martingale
differences $\square _{I}^{\sigma ,\mathbf{b}}$ constructed above with the
triple iterated coronas, as well as the corresponding data for $g$. Here is
a brief schematic diagram of the splittings and decompositions we will
describe below, with associated bounds given in $\fbox{}$. We split the form 
$\left\langle T_{\sigma }^{\alpha }f,g\right\rangle _{\omega }$ into the sum
of two essentially symmetric forms by interval size,%
\begin{equation}
\int \left( T_{\sigma }f\right) gd\omega =\left\{ \sum_{\substack{ I\in 
\mathcal{D}:\ J\in \mathcal{G}  \\ \ell \left( J\right) \leq \ell \left(
I\right) }}+\sum_{\substack{ I\in \mathcal{D}:\ J\in \mathcal{G}  \\ \ell
\left( J\right) >\ell \left( I\right) }}\right\} \int \left( T_{\sigma
}\square _{I}^{\sigma ,\mathbf{b}}f\right) \square _{J}^{\omega ,\mathbf{b}%
^{\ast }}gd\omega ,  \label{ess symm}
\end{equation}%
and focus on the first sum,%
\begin{equation*}
\Theta \left( f,g\right) =\sum_{I\in \mathcal{D}\text{ and }J\in \mathcal{G}%
:\ \ell \left( J\right) \leq \ell \left( I\right) }\left\langle T_{\sigma
}^{\alpha }\square _{I}^{\sigma ,\mathbf{b}}f,\square _{J}^{\omega ,\mathbf{b%
}^{\ast }}\right\rangle _{\omega },
\end{equation*}%
since the second sum is handled dually, but is easier due to the missing
diagonal.%
\begin{equation}
\fbox{$%
\begin{array}{ccccccc}
\Theta \left( f,g\right) &  &  &  &  &  &  \\ 
\downarrow &  &  &  &  &  &  \\ 
\Theta _{2}^{\limfunc{good}}\left( =\mathsf{B}_{\Subset _{\mathbf{r}%
}}\right) \left( f,g\right) & + & \Theta _{1}\left( =\mathsf{B}_{\cap
}\right) \left( f,g\right) & + & \Theta _{3}\left( =\mathsf{B}_{\diagup
}\right) \left( f,g\right) & + & \Theta _{2}^{\limfunc{bad}}\left( f,g\right)
\\ 
\downarrow &  & \fbox{$\mathcal{NTV}_{\alpha }$} &  & \fbox{$\mathcal{NTV}%
_{\alpha }+\sqrt{\theta }\mathfrak{N}_{T^{\alpha }}$} &  & \fbox{$2^{-%
\mathbf{r}\varepsilon }\mathfrak{N}_{T^{\alpha }}$} \\ 
\downarrow &  &  &  &  &  &  \\ 
\mathsf{T}_{\limfunc{diagonal}}\left( f,g\right) & + & \mathsf{T}_{\limfunc{%
far}\limfunc{below}}\left( f,g\right) & + & \mathsf{T}_{\limfunc{far}%
\limfunc{above}}\left( f,g\right) & + & \mathsf{T}_{\limfunc{disjoint}%
}\left( f,g\right) \\ 
\downarrow &  & \fbox{$\mathcal{NTV}_{\alpha }$} &  & \fbox{$\emptyset $} & 
& \fbox{$\emptyset $} \\ 
\mathsf{B}_{\Subset _{\mathbf{r}}}^{A}\left( f,g\right) &  &  &  &  &  &  \\ 
\downarrow &  &  &  &  &  &  \\ 
\mathsf{B}_{\limfunc{stop}}^{A}\left( f,g\right) & + & \mathsf{B}_{\limfunc{%
paraproduct}}^{A}\left( f,g\right) & + & \mathsf{B}_{\limfunc{neighbour}%
}^{A}\left( f,g\right) & + & \mathsf{B}_{\limfunc{broken}}^{A}\left(
f,g\right) \\ 
\fbox{$\mathcal{E}_{2}^{\alpha }+\sqrt{A_{2}^{\alpha }}+\sqrt{A_{2}^{\alpha ,%
\limfunc{punct}}}$} &  & \fbox{$\mathfrak{T}_{T^{\alpha }}^{\mathbf{b}}$} & 
& \fbox{$\sqrt{A_{2}^{\alpha }}$} &  & \fbox{$\mathfrak{T}_{T^{\alpha }}^{%
\mathbf{b}}$}%
\end{array}%
$}  \label{schematic}
\end{equation}

For the reader's convenience we now collect the various martingale and
probability estimates that will be used in the proof that follows. First we
summarize the martingale identities and estimates from Appendix A that we
will use in our proof, noting in particular that \emph{lower weak Riesz}
inequalities are \textbf{not} used in the proof of our $Tb$ theorem. Suppose 
$\mu $ is a positive locally finite Borel measure, and that $\mathbf{b}$ is
an $\infty $-strongly $\mu $-controlled accretive family. Then the following
martingale identities and estimates hold:

\begin{description}
\item[Martingale identities] Both of the following identities hold pointwise 
$\mu $-almost everywhere, as well as in the sense of strong convergence in $%
L^{2}\left( \mu \right) $:%
\begin{eqnarray*}
f &=&\sum_{I\in \mathcal{D}_{N}}\mathbb{F}_{I}^{\mu ,\mathbf{b}}f+\sum_{I\in 
\mathcal{D}:\ \ell \left( I\right) \geq N+1}\square _{I}^{\mu ,\mathbf{b}}f\
, \\
f &=&\sum_{I\in \mathcal{D}_{N}}\mathbb{E}_{I}^{\mu ,\mathbf{b}}f+\sum_{I\in 
\mathcal{D}:\ \ell \left( I\right) \geq N+1}\bigtriangleup _{I}^{\mu ,%
\mathbf{b}}f\ .
\end{eqnarray*}

\item[Frame estimates] Both of the following frame estimates hold:%
\begin{eqnarray}
\left\Vert f\right\Vert _{L^{2}\left( \mu \right) }^{2} &\approx &\sum_{Q\in 
\mathcal{D}}\left\{ \left\Vert \square _{Q}^{\mu ,\mathbf{b}}f\right\Vert
_{L^{2}\left( \mu \right) }^{2}+\left\Vert \bigtriangledown _{Q}^{\mu
}f\right\Vert _{L^{2}\left( \mu \right) }^{2}\right\}  \label{FRAME} \\
&\approx &\sum_{Q\in \mathcal{D}}\left\{ \left\Vert \bigtriangleup _{Q}^{\mu
,\mathbf{b}}f\right\Vert _{L^{2}\left( \mu \right) }^{2}+\left\Vert
\bigtriangledown _{Q}^{\mu }f\right\Vert _{L^{2}\left( \mu \right)
}^{2}\right\} \ .  \notag
\end{eqnarray}

\item[Weak upper Riesz estimates] Define the pseudoprojections, 
\begin{eqnarray*}
\Psi _{\mathcal{B}}^{\mu ,\mathbf{b}}f &\equiv &\sum_{I\in \mathcal{B}%
}\square _{I}^{\mu ,\mathbf{b}}f, \\
\left( \Psi _{\mathcal{B}}^{\mu ,\mathbf{b}}\right) ^{\ast }f &\equiv
&\sum_{I\in \mathcal{B}}\left( \square _{I}^{\mu ,\mathbf{b}}\right) ^{\ast
}f=\sum_{I\in \mathcal{B}}\bigtriangleup _{I}^{\mu ,\mathbf{b}}f.
\end{eqnarray*}%
We have the `upper Riesz' inequalities for pseudoprojections $\Psi _{%
\mathcal{B}}^{\mu ,\mathbf{b}}$ and $\left( \Psi _{\mathcal{B}}^{\mu ,%
\mathbf{b}}\right) ^{\ast }$:%
\begin{eqnarray}
\left\Vert \Psi _{\mathcal{B}}^{\mu ,\mathbf{b}}f\right\Vert _{L^{2}\left(
\mu \right) }^{2} &\leq &C\sum_{I\in \mathcal{B}}\left\Vert \square
_{I}^{\mu ,\mathbf{b}}f\right\Vert _{L^{2}\left( \mu \right)
}^{2}+\sum_{I\in \mathcal{B}}\left\Vert \nabla _{I}^{\mu }f\right\Vert
_{L^{2}\left( \mu \right) }^{2},  \label{UPPER RIESZ} \\
\left\Vert \left( \Psi _{\mathcal{B}}^{\mu ,\mathbf{b}}\right) ^{\ast
}f\right\Vert _{L^{2}\left( \mu \right) }^{2} &\leq &C\sum_{I\in \mathcal{B}%
}\left\Vert \bigtriangleup _{I}^{\mu ,\mathbf{b}}f\right\Vert _{L^{2}\left(
\mu \right) }^{2}+\sum_{I\in \mathcal{B}}\left\Vert \nabla _{I}^{\mu
}f\right\Vert _{L^{2}\left( \mu \right) }^{2},  \notag
\end{eqnarray}%
for all $f\in L^{2}\left( \mu \right) $ and all subsets $\mathcal{B}$ of the
grid $\mathcal{D}$, and where the positive constant $C$ is independent of
the subset $\mathcal{B}$. Here $\nabla _{I}^{\mu }$ is the Carleson
averaging operator defined in (\ref{Carleson avg op}) in Appendix A.
\end{description}

Now we turn to the probability estimates for martingale differences and
halos that we will use. Recall that given $0<\lambda <\frac{1}{2}$, the $%
\lambda $-halo of $J$ is defined to be 
\begin{equation*}
\partial _{\lambda }J\equiv \left( 1+\lambda \right) J\setminus \left(
1-\lambda \right) J.
\end{equation*}%
Suppose $\mu $ is a positive locally finite Borel measure, and that $\mathbf{%
b}$ is an $\infty $-weakly $\mu $-controlled accretive family. Then the
following probability estimates hold. See Definition \ref{def Gbad} below
for the notation $\mathcal{G}_{k-\limfunc{bad}}^{\mathcal{D}}$.

\begin{description}
\item[Bad cube probability estimates] Suppose that $\mathcal{D}$ and $%
\mathcal{G}$ are independent random dyadic grids. With $\Psi _{\mathcal{G}%
_{k-\limfunc{bad}}^{\mathcal{D}}}^{\mu ,\mathbf{b}^{\ast }}g\equiv
\sum_{J\in \mathcal{G}_{k-\limfunc{bad}}^{\mathcal{D}}}\square _{J}^{\mu ,%
\mathbf{b}^{\ast }}g$ equal to the pseudoprojection of $g$ onto $k$-bad $%
\mathcal{G}$-intervals, we have%
\begin{equation*}
\boldsymbol{E}_{\Omega }^{\mathcal{D}}\left( \left\Vert \Psi _{\mathcal{G}%
_{k-\limfunc{bad}}^{\mathcal{D}}}^{\mu ,\mathbf{b}^{\ast }}g\right\Vert
_{L^{2}\left( \mu \right) }^{2}\right) \lesssim \boldsymbol{E}_{\Omega }^{%
\mathcal{D}}\left( \sum_{J\in \mathcal{G}_{k-\limfunc{bad}}^{\mathcal{D}}}%
\left[ \left\Vert \square _{J,\mathcal{G}}^{\mu ,\mathbf{b}^{\ast
}}g\right\Vert _{L^{2}\left( \mu \right) }^{2}+\left\Vert \nabla _{J,%
\mathcal{G}}^{\mu }g\right\Vert _{L^{2}\left( \mu \right) }^{2}\right]
\right) \leq Ce^{-k\varepsilon }\left\Vert g\right\Vert _{L^{2}\left( \mu
\right) }^{2}\ ,
\end{equation*}%
where the first inequality is the `weak upper half Riesz' inequality for the
pseudoprojection $\Psi _{\mathcal{G}_{k-\limfunc{bad}}^{\mathcal{D}}}^{\mu ,%
\mathbf{b}^{\ast }}$, and the second inequality is proved using the frame
inequality in (\ref{main bad prob}) below.

\item[Halo probability estimates] Suppose that $\mathcal{D}$ and $\mathcal{G}
$ are independent random grids. Using the \emph{parameterization by
translations}\ of grids and taking the average over certain translates $\tau
+\mathcal{D}$ of the grid $\mathcal{D}$ we have%
\begin{eqnarray}
\boldsymbol{E}_{\Omega }^{\mathcal{D}}\sum_{I^{\prime }\in \mathcal{D}:\
\ell \left( I^{\prime }\right) \approx \ell \left( J^{\prime }\right)
}\int_{J^{\prime }\cap \partial _{\delta }I^{\prime }}d\omega &\lesssim &%
\mathbb{\delta }\int_{J^{\prime }}d\omega ,\ \ \ \ \ J^{\prime }\in 
\mathfrak{C}\left( J\right) ,J\in \mathcal{G},  \label{hand'} \\
\boldsymbol{E}_{\Omega }^{\mathcal{G}}\sum_{J^{\prime }\in \mathcal{G}:\
\ell \left( J^{\prime }\right) \approx \ell \left( I^{\prime }\right)
}\int_{I^{\prime }\cap \partial _{\delta }J^{\prime }}d\sigma &\lesssim &%
\mathbb{\delta }\int_{I^{\prime }}d\sigma ,\ \ \ \ \ I^{\prime }\in 
\mathfrak{C}\left( I\right) ,I\in \mathcal{D},  \notag
\end{eqnarray}%
and where the expectations $\boldsymbol{E}_{\Omega }^{\mathcal{D}}$ and $%
\boldsymbol{E}_{\Omega }^{\mathcal{G}}$ are taken over grids $\mathcal{D}$
and $\mathcal{G}$ respectively. Indeed, it is geometrically evident that for
any fixed pair of side lengths $\ell _{1}\approx \ell _{2}$, the average of
the measure $\left\vert J^{\prime }\cap \partial _{\delta }I^{\prime
}\right\vert _{\omega }$ of the set $J^{\prime }\cap \partial _{\delta
}I^{\prime }$, as an interval $I^{\prime }\in \mathcal{D}$ with side length $%
\ell \left( I^{\prime }\right) =\ell _{1}$ is translated across an interval $%
J^{\prime }\in \mathcal{G}$ of side length $\ell \left( J^{\prime }\right)
=\ell _{2}$, is at most $C\left\vert J^{\prime }\right\vert _{\omega }$.
Using this observation it is now easy to see that (\ref{hand'}) holds.
\end{description}

\subsection{The Hyt\"{o}nen-Martikainen decomposition and a weak variant of
NTV goodness\label{Subsec HM}}

Let $\mathbf{b}$ (respectively $\mathbf{b}^{\ast }$) be $\infty $-weakly $%
\sigma $-controlled (respectively $\omega $-controlled) accretive families.
At the beginning of this section, we reduced the estimation of the bilinear
form $\int_{\mathbb{R}}\left( T_{\sigma }f\right) gd\omega $ to that of the
sum%
\begin{equation*}
\sum_{I\in \mathcal{D}}\sum_{J\in \mathcal{G}}\int \left( T_{\sigma
}^{\alpha }\square _{I}^{\sigma ,\mathbf{b}}f\right) \square _{J}^{\omega ,%
\mathbf{b}^{\ast }}gd\omega ,
\end{equation*}%
and then we decomposed this sum by interval side length,%
\begin{eqnarray*}
\sum_{I\in \mathcal{D}}\sum_{J\in \mathcal{G}}\int \left( T_{\sigma
}^{\alpha }\square _{I}^{\sigma ,\mathbf{b}}f\right) \square _{J}^{\omega ,%
\mathbf{b}^{\ast }}gd\omega &=&\left\{ \sum_{\substack{ I\in \mathcal{D}:\
J\in \mathcal{G}  \\ \ell \left( J\right) \leq \ell \left( I\right) }}+\sum 
_{\substack{ I\in \mathcal{D}:\ J\in \mathcal{G}  \\ \ell \left( J\right)
>\ell \left( I\right) }}\right\} \int \left( T_{\sigma }^{\alpha }\square
_{I}^{\sigma ,\mathbf{b}}f\right) \square _{J}^{\omega ,\mathbf{b}^{\ast
}}gd\omega \\
&\equiv &\Theta \left( f,g\right) +\Theta ^{\ast }\left( f,g\right) ,
\end{eqnarray*}%
and noted that by symmetry, it suffices to estimate the first form $\Theta
\left( f,g\right) $. Before introducing goodness into the sum, we follow 
\cite{HyMa} and split the form $\Theta \left( f,g\right) $ into 3 pieces:%
\begin{eqnarray*}
&&\Theta \left( f,g\right) \equiv \sum_{\substack{ I\in \mathcal{D}:\ J\in 
\mathcal{G}  \\ \ell \left( J\right) \leq \ell \left( I\right) }}\int \left(
T_{\sigma }^{\alpha }\square _{I}^{\sigma ,\mathbf{b}}f\right) \square
_{J}^{\omega ,\mathbf{b}^{\ast }}gd\omega \\
&=&\sum_{I\in \mathcal{D}}\left\{ \sum_{\substack{ J\in \mathcal{G}:\ \ell
\left( J\right) \leq \ell \left( I\right)  \\ d\left( J,I\right) >2\ell
\left( J\right) ^{\varepsilon }\ell \left( I\right) ^{1-\varepsilon }}}+\sum 
_{\substack{ J\in \mathcal{G}:\ \ell \left( J\right) \leq 2^{-\mathbf{r}%
}\ell \left( I\right)  \\ d\left( J,I\right) \leq 2\ell \left( J\right)
^{\varepsilon }\ell \left( I\right) ^{1-\varepsilon }}}+\sum_{\substack{ %
J\in \mathcal{G}:\ 2^{-\mathbf{r}}\ell \left( I\right) <\ell \left( J\right)
\leq \ell \left( I\right)  \\ d\left( J,I\right) \leq 2\ell \left( J\right)
^{\varepsilon }\ell \left( I\right) ^{1-\varepsilon }}}\right\} \int \left(
T_{\sigma }^{\alpha }\square _{I}^{\sigma ,\mathbf{b}}f\right) \square
_{J}^{\omega ,\mathbf{b}^{\ast }}gd\omega \\
&\equiv &\Theta _{1}\left( f,g\right) +\Theta _{2}\left( f,g\right) +\Theta
_{3}\left( f,g\right) \ ,
\end{eqnarray*}%
where $\varepsilon >0$ will be chosen to satisfy $0<\varepsilon <\frac{1}{2}%
\leq \frac{1}{2-\alpha }$ later, and the goodness parameter $\mathbf{r}$ is
then determined in (\ref{choice of r}) below. Now the disjoint form $\Theta
_{1}\left( f,g\right) $ can be handled by `long-range' and `short-range'
arguments which we give in the next section below, and the nearby form $%
\Theta _{3}\left( f,g\right) $ will be handled in the subsequent section
using probabilistic surgery methods and a new deterministic surgery
involving energy conditions and the `original' testing functions discarded
in the corona construction. The remaining form $\Theta _{2}\left( f,g\right) 
$ will be treated further in this section after introducing weak goodness.

\subsubsection{Good intervals with `$\limfunc{body}$'}

We begin with the weaker extension of goodness introduced in \cite{HyMa},
except that we will make it a bit stronger by replacing the skeleton `$%
\limfunc{skel}K$' of an interval $K$, as used in \cite{HyMa}, by a larger
collection of points `$\limfunc{body}K$', which we call the dyadic body of $%
K $. This modification will prove useful in establishing the Straddling
Lemma in the treatment of the stopping form in Section \ref{Sec stop} below.
Let $\mathcal{P}$ denote the collection of all intervals in $\mathbb{R}$.
The content of the next four definitions is inspired by, or sometimes
identical with, that already appearing in the work of Nazarov, Treil and
Volberg in \cite{NTV1} and \cite{NTV3}.

\begin{definition}
\label{skel spray body}Let $K\in \mathcal{P}$.

\begin{enumerate}
\item Define the \emph{skeleton} `$\limfunc{skel}K$' of $K$ to consist of
its two endpoints and its midpoint.

\item For a point $x$ in $\mathbb{R}$, define the \emph{dyadic spray} $%
\mathbb{S}_{x}^{\limfunc{dy}}\equiv \left\{ x\right\} \cup \left\{ x\pm 
\frac{1}{2^{j}}\right\} _{k\in \mathbb{Z}}$ of $x$ to consist of $x$ and all
points $y$ in $\mathbb{R}$ that have distance $\frac{1}{2^{j}}$ from $x$ for
some $j\in \mathbb{Z}$.

\item Then define the \emph{dyadic body} `$\limfunc{body}K$' of an interval $%
K\in \mathcal{P}$ to be the intersection of $\overline{K}$ with the union of
the dyadic sprays of its two endpoints, i.e. if $K=\left[ a,b\right) $, then%
\begin{equation*}
\limfunc{body}K=\overline{K}\cap \left( \mathbb{S}_{a}^{\limfunc{dy}}\cup 
\mathbb{S}_{b}^{\limfunc{dy}}\right) .
\end{equation*}
\end{enumerate}
\end{definition}

Thus the body of the unit interval $\left[ 0,1\right) $ consists of the
points 
\begin{equation*}
\left\{ 0\right\} \dot{\cup}\left\{ \frac{1}{2^{j}}\right\} _{j=1}^{\infty }%
\dot{\cup}\left\{ 1-\frac{1}{2^{j}}\right\} _{j=2}^{\infty }\dot{\cup}%
\left\{ 1\right\} \ ,
\end{equation*}%
which have the endpoints of $\left[ 0,1\right) $ as cluster points.

\begin{definition}
\label{good arb}Let $0<\varepsilon <1$ (to be chosen later). For intervals $%
J,K\in \mathcal{P}$ with $\ell \left( J\right) \leq \ell \left( K\right) $,
we define $J$ to be $\varepsilon -\limfunc{good}$ \emph{with respect to} an
interval $K$ if 
\begin{equation}
d\left( J,\limfunc{body}K\right) >2\left\vert J\right\vert ^{\varepsilon
}\left\vert K\right\vert ^{1-\varepsilon },  \label{eps far}
\end{equation}%
and we say $J$ is $\varepsilon -\limfunc{bad}$ \emph{with respect to} $K$ if
(\ref{eps far}) fails. We also say that $J$ is $\varepsilon -\limfunc{good}$ 
\emph{inside} an interval $K$ if $J$ is $\varepsilon -\limfunc{good}$ with
respect to $K$ and $J\subset K$.
\end{definition}

A key consequence of an interval $J$ being $\varepsilon -\limfunc{good}$
inside an interval $S$, is that $J$ must then be contained in some \emph{%
dyadic subinterval} $K$ of $S$ with $3K\subset S$:%
\begin{equation}
\text{If }J\text{ is }\varepsilon -\limfunc{good}\text{inside }S\text{, then 
}J\subset K\text{ for some }K\in \mathcal{W}\left( S\right) ,
\label{key contain}
\end{equation}%
where $\mathcal{W}\left( S\right) $ is the collection of maximal dyadic
subintervals of $S$ whose triples are contained in $S$. Indeed, the
endpoints of the intervals in $\mathcal{W}\left( S\right) $ are precisely
the $\limfunc{body}$ of $S$. Note that this property can fail if we use the
smaller set $\limfunc{skel}S$ in place of $\limfunc{body}S$ in Definition %
\ref{good arb}, since then an $\varepsilon -\limfunc{good}$ interval $J$
could intersect one of the sprays. Of course we will also need to know that $%
\limfunc{body}S$ is not so much larger than $\limfunc{skel}S$ that the
crucial probability estimate for good intervals fails - namely we need to
know that given $k\gg 1$, an interval $J\subset S$ of side length $\ell
\left( J\right) =2^{-k}\ell \left( S\right) $ is $\varepsilon -\limfunc{good}
$ inside $S$ with `large probability'. This will be made precise below using
random dyadic grids.

\begin{definition}
\label{good two grids}Let $\mathcal{D}$ and $\mathcal{G}$ be dyadic grids.
Define $\mathcal{G}_{\left( k,\varepsilon \right) -\limfunc{good}}^{\mathcal{%
D}}$ to consist of those $J\in \mathcal{G}$ such that $J$ is $\varepsilon -%
\limfunc{good}$ \textbf{inside} every interval $K\in \mathcal{D}$ with $%
K\cap J\neq \emptyset $ that lies at least $k$ levels `above' $J$, i.e. $%
\ell \left( K\right) \geq 2^{k}\ell \left( J\right) $ (note that the use of 
\textbf{inside} forces such $K$ with $K\cap J\neq \emptyset $ to actually
contain $J$). We also define $J$ to be `$\varepsilon -\limfunc{good}$ inside
an interval $K$ and beyond' if $J\in \mathcal{G}_{\left( k,\varepsilon
\right) -\limfunc{good}}^{\mathcal{D}}$ where $k=\log _{2}\frac{\ell \left(
K\right) }{\ell \left( J\right) }$ and where $K\cap J\neq \emptyset $,
equivalently in this situation $K\supset J$. As the goodness parameter $%
\varepsilon $ will eventually be fixed throughout the proof, we sometimes
suppress it, and simply say `$J$ is $\limfunc{good}$ inside an interval $K$
and beyond' instead of `$J$ is $\varepsilon -\limfunc{good}$ inside an
interval $K$ and beyond'. When $\varepsilon >0$ is understood, we will often
write $\mathcal{G}_{k-\limfunc{good}}^{\mathcal{D}}=\mathcal{G}_{\left(
k,\varepsilon \right) -\limfunc{good}}^{\mathcal{D}}$.
\end{definition}

\begin{remark}
Note that%
\begin{equation*}
\mathcal{G}_{\left( k,\varepsilon \right) -\limfunc{good}}^{\mathcal{D}%
}\equiv \left\{ J\in \mathcal{G}:J\text{ is }\varepsilon -\limfunc{good}%
\text{ with respect to every }K\in \mathcal{D}\text{ with }\ell \left(
K\right) \geq 2^{k}\ell \left( J\right) .\right\}
\end{equation*}%
Indeed, if $J$ is $\varepsilon -\limfunc{bad}$ with respect to some $K\in 
\mathcal{D}$ with $K\cap J=\emptyset $, then $J$ is also $\varepsilon -%
\limfunc{bad}$ with respect to one of the two neighbours (of the same side
length) of $K$ in $\mathcal{D}$.
\end{remark}

\subsubsection{Grid probability}

As pointed out on page 14 of \cite{HyMa} by Hyt\"{o}nen and Martikainen,
there are subtle difficulties associated in using dual martingale
decompositions of functions which depend on the entire dyadic grid, rather
than on just the local interval in the grid. We will proceed at first in the
spirit of \cite{HyMa}, and the goodness that we will infuse below into the
main `below' form $\mathsf{B}_{\Subset _{\mathbf{r}}}\left( f,g\right) $
will be the Hyt\"{o}nen-Martikainen `weak' version of NTV goodness, but
using the body `$\limfunc{body}I$' of an interval rather than\ its skeleton `%
$\limfunc{skel}I$': every pair $\left( I,J\right) \in \mathcal{D}\times 
\mathcal{G}$ that arises in the form $\mathsf{B}_{\Subset _{\mathbf{r}%
}}\left( f,g\right) $ will satisfy $J\in \mathcal{G}_{\left( k,\varepsilon
\right) -\limfunc{good}}^{\mathcal{D}}$ where $\ell \left( I\right)
=2^{k}\ell \left( J\right) $.

Now we return to the martingale differences $\square _{I}^{\sigma ,\mathbf{b}%
}$ and $\square _{J}^{\omega ,\mathbf{b}^{\ast }}$ with controlled families $%
\mathbf{b}$ and $\mathbf{b}^{\ast }$ in the real line $\mathbb{R}$. When we
want to emphasize that the grid in use is $\mathcal{D}$ or $\mathcal{G}$, we
will denote the martingale difference by $\square _{I,\mathcal{D}}^{\sigma ,%
\mathbf{b}}$, and similarly for $\square _{J,\mathcal{G}}^{\omega ,\mathbf{b}%
^{\ast }}$. Recall Definition \ref{good arb} for the meaning of when an
interval $J$ is $\varepsilon $-$\limfunc{bad}$ with respect to another
interval $K$.

\begin{definition}
\label{bad in grid}We say that $J\in \mathcal{P}$ is $k$-$\limfunc{bad}$ in
a grid $\mathcal{D}$ if there is an interval $K\in \mathcal{D}$ with $\ell
\left( K\right) =2^{k}\ell \left( J\right) $ such that $J$ is $\varepsilon $-%
$\limfunc{bad}$ with respect to $K$ (context should eliminate any ambiguity
between the different use of $k$-$\limfunc{bad}$ when $k\in \mathbb{N}$ and $%
\varepsilon $-$\limfunc{bad}$ when $0<\varepsilon <\frac{1}{2}$).
\end{definition}

A key observation here (see \cite{NTV1}, \cite{NTV2}, \cite{NTV3} or \cite%
{NTV4} for the case when goodness is defined using the skeleton instead of
the body) is that for any $J\in \mathcal{G}$ where $\mathcal{D}$ and $%
\mathcal{G}$ are independent random grids,%
\begin{equation}
\boldsymbol{P}_{\Omega }^{\mathcal{D}}\left( \mathcal{D}:J\text{ is }k\text{-%
}\limfunc{bad}\text{ in }\mathcal{D}\right) \equiv \int_{\Omega }\mathbf{1}%
_{\left\{ \mathcal{D}:\ J\text{ is }k\text{-}\limfunc{bad}\text{ in }%
\mathcal{D}\right\} }d\mu _{\Omega }\left( \mathcal{D}\right) \leq
C\varepsilon k2^{-\varepsilon k}.  \label{key prob}
\end{equation}%
Indeed, it suffices to consider the case when $J\in \mathcal{G}$ with $%
J\subset \left[ 0,1\right) $ and $\ell \left( J\right) =2^{-k}$. So fix such
an interval $J$. For each $m\in \mathbb{Z}_{2^{k}}\equiv \left\{ \ell \in 
\mathbb{Z}:0\leq \ell \leq 2^{k}-1\right\} $, consider the collection $%
\mathfrak{D}_{m}$ of all grids $\mathcal{D}$ that contain the interval $%
I_{m}\equiv \left[ 0,1\right) +\frac{m}{2^{k}}=\left[ \frac{m}{2^{k}},1+%
\frac{m}{2^{k}}\right) $. Then for every $m$, it is the case that

\begin{enumerate}
\item \textbf{either} $J$ is $\varepsilon $-$\limfunc{bad}$ in $\mathcal{D}$
for all $\mathcal{D}\in \mathfrak{D}_{m}$,

\item \textbf{or} $J$ is $\varepsilon $-$\limfunc{good}$ in $\mathcal{D}$
for all $\mathcal{D}\in \mathfrak{D}_{m}$.
\end{enumerate}

We will say that the \emph{collection} $\mathfrak{D}_{m}$ is $k$-$\limfunc{%
bad}$ if the first case holds. We have the same dichotomy for $\mathfrak{D}%
_{m}+s$ if we replace $\left[ 0,1\right) $ with the translate $\left[
0,1\right) +s=\left[ s,1+s\right) $ where $0\leq s<2^{-k}$. We now claim
that for any fixed $0\leq s<2^{-k}$, the number of $k$-$\limfunc{bad}$
collections $\mathfrak{D}_{m}+s$ is at most $C\varepsilon k2^{\left(
1-\varepsilon \right) k}$, hence the proportion of $k$-$\limfunc{bad}$
collections is $\frac{C\varepsilon k2^{\left( 1-\varepsilon \right) k}}{2^{k}%
}=C\varepsilon k2^{-\varepsilon k}$, from which we obtain the estimate (\ref%
{key prob}) as follows. Every grid $\mathcal{D}\in \Omega $ is contained in
exactly one of the collections $\left\{ \mathfrak{D}_{m}+s\right\} _{m\in 
\mathbb{Z}_{2^{k}}\text{ and }s\in \left[ 0,2^{-k}\right) }$, and so%
\begin{eqnarray*}
\boldsymbol{P}_{\Omega }^{\mathcal{D}}\left( \mathcal{D}:J\text{ is }k\text{-%
}\limfunc{bad}\text{ in }\mathcal{D}\right) &=&\frac{1}{2^{-k}}\int_{\left[
0,2^{-k}\right) }\left\{ \frac{\#\left\{ m\in \mathbb{Z}_{2^{k}}:\ \mathfrak{%
D}_{m}+s\text{ is }k-\limfunc{bad}\right\} }{\#\mathbb{Z}_{2^{k}}}\right\} ds
\\
&\leq &\frac{1}{2^{-k}}\int_{\left[ 0,2^{-k}\right) }\left\{ \frac{%
C\varepsilon k2^{\left( 1-\varepsilon \right) k}}{2^{k}}\right\}
ds=C\varepsilon k2^{-\varepsilon k}.
\end{eqnarray*}

To see our claim, it suffices to consider the case $s=0$, to keep the
interval $I_{0}=\left[ 0,1\right) $ fixed, and consider instead the
translates $J_{m}\equiv J+\frac{m}{2^{k}}$ of the interval $J$ for $0\leq
m\leq 2^{k}-1$. Moreover we can assume without loss of generality that $J$
intersects the point $\frac{1}{2^{k}}$ so that all of the intervals $J_{m}$
in the collection $\left\{ J_{m}\right\} _{m=0}^{2^{k}-1}$ lie in $I_{0}$
except for the last one $J_{2^{k}-1}=J+1-\frac{1}{2^{k}}$, which intersects
the point $1$. In this situation our claim becomes%
\begin{equation}
\#\left\{ m\in \mathbb{Z}_{2^{k}}:J_{m}\text{ is }\varepsilon -\limfunc{bad}%
\text{ in }\left[ 0,1\right) \right\} \leq C\varepsilon k2^{\left(
1-\varepsilon \right) k}.  \label{claim becomes}
\end{equation}

To prove (\ref{claim becomes}), we begin by defining%
\begin{equation*}
d\equiv \varepsilon k-1\text{ and }L\equiv 2\ell \left( J\right)
^{\varepsilon }\ell \left( I_{0}\right) ^{1-\varepsilon }=2^{1-\varepsilon
k}=2^{-d},
\end{equation*}%
where we may assume $k>\frac{1}{\varepsilon }$ so that $d>0$. Then if $J_{m}$
is $\varepsilon -\limfunc{bad}$ in $\left[ 0,1\right) $, at least one of the
following two inequalities must hold:%
\begin{equation*}
\limfunc{dist}\left( J_{m},\mathbb{S}_{0}^{\limfunc{dy}}\right) \leq L,\ 
\limfunc{dist}\left( J_{m},\mathbb{S}_{1}^{\limfunc{dy}}\right) \leq L,
\end{equation*}%
where we recall that $\mathbb{S}_{a}^{\limfunc{dy}}$ is the dyadic spray of $%
a$. Now if $\limfunc{dist}\left( J_{m},\mathbb{S}_{0}^{\limfunc{dy}}\right)
\leq L$, then%
\begin{eqnarray*}
\text{either }\limfunc{dist}\left( J_{m},\left\{ 0\right\} \cup \left\{ 
\frac{1}{2^{j}}\right\} _{j>d}\right) &\leq &L, \\
\text{or }\limfunc{dist}\left( J_{m},\frac{1}{2^{j}}\right) &\leq &L,\ \ \ \
\ \text{for some }0\leq j\leq d.
\end{eqnarray*}%
However, if 
\begin{equation}
\limfunc{dist}\left( J_{m},\left\{ 0\right\} \cup \left\{ \frac{1}{2^{j}}%
\right\} _{j>d}\right) \leq L,  \label{first case}
\end{equation}%
then we must have 
\begin{equation*}
m\ell \left( J\right) \leq \frac{1}{2^{d}}+L=2L,
\end{equation*}%
and if $\limfunc{dist}\left( J_{m},\frac{1}{2^{j}}\right) \leq L$ for some $%
0\leq j\leq d$, then we must have%
\begin{equation*}
\left\vert \frac{m}{2^{k}}-\frac{1}{2^{j}}\right\vert \leq 2L.
\end{equation*}%
So altogether the number of indices $m\in \mathbb{Z}_{2^{k}}$ for which $%
\limfunc{dist}\left( J_{m},\mathbb{S}_{0}^{\limfunc{dy}}\right) \leq L$
holds is at most%
\begin{eqnarray*}
2\frac{L}{\ell \left( J\right) }+1+\left( d+1\right) \left(
2^{k+2}L+1\right) &=&\left( 2d+3\right) 2^{k+1}L+d+2 \\
&=&\left( 2\varepsilon k+1\right) \cdot 2^{k+1}\cdot 2^{1-\varepsilon
k}+\varepsilon k+1\leq 20\varepsilon k2^{\left( 1-\varepsilon \right) k}.
\end{eqnarray*}%
Similarly the number of indices $m\in \mathbb{Z}_{2^{k}}$ for which $%
\limfunc{dist}\left( J_{m},\mathbb{S}_{1}^{\limfunc{dy}}\right) \leq L$
holds is at most $20\varepsilon k2^{\left( 1-\varepsilon \right) k}$. Thus
we conclude that (\ref{claim becomes}) holds with $C=40$.

Then we obtain from (\ref{key prob}), using the lower frame inequality, the
expectation estimate%
\begin{eqnarray*}
&&\int_{\Omega }\sum_{J\in \mathcal{G}_{k-\limfunc{bad}}^{\mathcal{D}}}\left[
\left\Vert \square _{J,\mathcal{G}}^{\omega ,\mathbf{b}^{\ast }}g\right\Vert
_{L^{2}\left( \omega \right) }^{2}+\left\Vert \nabla _{J,\mathcal{G}%
}^{\omega }g\right\Vert _{L^{2}\left( \omega \right) }^{2}\right] d\mu
_{\Omega }\left( \mathcal{D}\right) \\
&=&\sum_{J\in \mathcal{G}}\left[ \left\Vert \square _{J,\mathcal{G}}^{\omega
,\mathbf{b}^{\ast }}g\right\Vert _{L^{2}\left( \omega \right)
}^{2}+\left\Vert \nabla _{J,\mathcal{G}}^{\omega }g\right\Vert _{L^{2}\left(
\omega \right) }^{2}\right] \int_{\Omega }\mathbf{1}_{\left\{ \mathcal{D}:\ J%
\text{ is }k\text{-}\limfunc{bad}\text{ in }\mathcal{D}\right\} }d\mu
_{\Omega }\left( \mathcal{D}\right) \\
&\leq &Ck\varepsilon 2^{-k\varepsilon }\sum_{J\in \mathcal{G}}\left[
\left\Vert \square _{J,\mathcal{G}}^{\omega ,\mathbf{b}^{\ast }}g\right\Vert
_{L^{2}\left( \omega \right) }^{2}+\left\Vert \nabla _{J,\mathcal{G}%
}^{\omega }g\right\Vert _{L^{2}\left( \omega \right) }^{2}\right] \leq
Ck\varepsilon 2^{-k\varepsilon }\left\Vert g\right\Vert _{L^{2}\left( \omega
\right) }^{2}\ ,
\end{eqnarray*}%
where $\nabla _{J,\mathcal{G}}^{\omega }$ denotes the `broken' Carleson
averaging operator in (\ref{Carleson avg op}) that depends on the broken
children in the grid $\mathcal{G}$. Altogether then it follows easily that%
\begin{equation}
\boldsymbol{E}_{\Omega }^{\mathcal{D}}\left( \sum_{J\in \bigcup_{\ell
=k}^{\infty }\mathcal{G}_{\ell -\limfunc{bad}}^{\mathcal{D}}}\left[
\left\Vert \square _{J,\mathcal{G}}^{\omega ,\mathbf{b}^{\ast }}g\right\Vert
_{L^{2}\left( \omega \right) }^{2}+\left\Vert \nabla _{J,\mathcal{G}%
}^{\omega }g\right\Vert _{L^{2}\left( \omega \right) }^{2}\right] \right)
\leq Ck\varepsilon 2^{-k\varepsilon }\left\Vert g\right\Vert _{L^{2}\left(
\omega \right) }^{2}\ ,  \label{main bad prob}
\end{equation}%
for some large positive constant $C$.

From such inequalities summed for $k\geq \mathbf{r}$, it can be concluded as
in \cite{NTV3} that there is an absolute choice of $\mathbf{r}$ depending on 
$0<\varepsilon <\frac{1}{2}$ so that the following holds. Let $%
T\;:\;L^{2}(\sigma )\rightarrow L^{2}(\omega )$ be a bounded linear
operator. We then have the following traditional inequality for two random
grids in the case that $\mathbf{b}$ is an $\infty $-strongly $\mu $%
-controlled accretive family: 
\begin{equation}
\left\Vert T\right\Vert _{L^{2}(\sigma )\rightarrow L^{2}(\omega )}\leq
2\sup_{\left\Vert f\right\Vert _{L^{2}(\sigma )}=1}\sup_{\left\Vert
g\right\Vert _{L^{2}(\omega )}=1}\boldsymbol{E}_{\Omega }\boldsymbol{E}%
_{\Omega ^{\prime }}\left\vert \left\langle \sum_{I,J\in \mathcal{D}_{%
\mathbf{r}-\limfunc{good}}^{\mathcal{G}}}T\left( \square _{I,\mathcal{D}%
}^{\sigma ,\mathbf{b}}f\right) f,\square _{J,\mathcal{D}}^{\omega ,\mathbf{b}%
^{\ast }}g\right\rangle _{\omega }\right\vert \,.  \label{e.Tgood'}
\end{equation}

However, this traditional method of introducing goodness is flawed here in
the general setting of dual martingale differences, since these differences
are no longer orthogonal projections, and as emphasized in \cite{HyMa}, we
cannot simply add back in bad intervals whenever we want telescoping
identities to hold - but these are needed in order to control the right hand
side of (\ref{e.Tgood'}). In fact, in the analysis of the form $\Theta
\left( f,g\right) $ above, it is necessary to have goodness for the
intervals $J$ and telescoping for the intervals $I$. On the other hand, in
the analysis of the form $\Theta ^{\ast }\left( f,g\right) $ above, it is
necessary to have just the opposite - namely goodness for the intervals $I$
and telescoping for the intervals $J$.

Thus, because in this unfortunate set of circumstances we can no longer `add
back in' bad cubes to achieve telescoping, we are prevented from introducing
goodness in the \emph{full} sum (\ref{ess symm})\ over all $I$ and $J$,
prior to splitting according to side lengths of $I$ and $J$. Thus the
infusion of goodness must come \emph{after} the splitting by side length,
but one must work much harder to introduce goodness directly into the form $%
\Theta \left( f,g\right) $ \emph{after} we have restricted the sum to
intervals $J$ that have smaller side length than $I$. This is accomplished
in the next subsubsection using the \emph{weaker form of NTV goodness}
introduced by Hyt\"{o}nen and Martikainen in \cite{HyMa} (that permits
certain additional pairs $\left( I,J\right) $ in the good forms where $\ell
\left( J\right) \leq 2^{-\mathbf{r}}\ell \left( I\right) $ and yet $J$ is $%
\limfunc{bad}$ in the traditional sense), and that will prevail later in the
treatment of the far below forms $\mathsf{T}_{\limfunc{far}\limfunc{below}%
}^{1}\left( f,g\right) $, and of the local forms $\mathsf{B}_{\Subset _{%
\mathbf{r}}}^{A}\left( f,g\right) $ (see Subsection \ref{Sub wrapup}) where
the need for using the `body' of an interval will become apparent in dealing
with the stopping form, and also in the treatment of the functional energy
in Appendix B.

\subsubsection{Weak goodness}

Let $\mathcal{D}$ and $\mathcal{G}$ be dyadic grids. It remains to estimate
the form $\Theta _{2}\left( f,g\right) $ which, following \cite{HyMa}, we
will split into a `bad' part and a `good' part. For this we introduce our
main definition associated with the above modification of the weak goodness
of Hyt\"{o}nen and Martikainen, namely the definition of the interval $%
R^{\maltese }$ in a grid $\mathcal{D}$, given an arbitrary interval $R\in 
\mathcal{P}$.

\begin{definition}
\label{def sharp cross}Let $\mathcal{D}$ be a dyadic grid. Given $R\in 
\mathcal{P}$, let $R^{\maltese }$ be the smallest (if any such exist) $%
\mathcal{D}$-dyadic superinterval $Q$ of $R$ such that $R$ is good inside 
\textbf{all} $\mathcal{D}$-dyadic superintervals $K$ of $Q$. Of course $%
R^{\maltese }$ will not exist if there is no $\mathcal{D}$-dyadic interval $%
Q $ containing $R$ in which $R$ is good. For intervals $R,Q\in \mathcal{P}$
let $\kappa \left( Q,R\right) =\log _{2}\frac{\ell \left( Q\right) }{\ell
\left( R\right) }$. For $R\in \mathcal{P}$ for which $R^{\maltese }$ exists,
let $\kappa \left( R\right) \equiv \kappa \left( R^{\maltese },R\right) $.
\end{definition}

Note that we typically suppress the dependence of $R^{\maltese }$ on the
grid $\mathcal{D}$, since the grid is usually understood from context. If $%
R^{\maltese }$ exists, we thus have that $R$ is good inside all $\mathcal{D}$%
-dyadic superintervals $K$ of $R$ with $\ell \left( K\right) \geq \ell
\left( R^{\maltese }\right) $. Note in particular the monotonicity property
for $J^{\prime },J\in \mathcal{P}$:%
\begin{equation*}
J^{\prime }\subset J\Longrightarrow \left( J^{\prime }\right) ^{\maltese
}\subset J^{\maltese }.
\end{equation*}%
Here now is the decomposition:%
\begin{eqnarray*}
\Theta _{2}\left( f,g\right) &=&\sum_{I\in \mathcal{D}}\sum_{\substack{ J\in 
\mathcal{G}:\ J^{\maltese }\not\subsetneqq I\text{, }\ell \left( J\right)
\leq 2^{-\mathbf{r}}\ell \left( I\right)  \\ d\left( J,I\right) \leq 2\ell
\left( J\right) ^{\varepsilon }\ell \left( I\right) ^{1-\varepsilon }}}\int
\left( T_{\sigma }^{\alpha }\square _{I}^{\sigma ,\mathbf{b}}f\right)
\square _{J}^{\omega ,\mathbf{b}^{\ast }}gd\omega \\
&&+\sum_{I\in \mathcal{D}}\sum_{\substack{ J\in \mathcal{G}:\ J^{\maltese
}\subsetneqq I\text{, }\ell \left( J\right) \leq 2^{-\mathbf{r}}\ell \left(
I\right)  \\ d\left( J,I\right) \leq 2\ell \left( J\right) ^{\varepsilon
}\ell \left( I\right) ^{1-\varepsilon }}}\int \left( T_{\sigma }^{\alpha
}\square _{I}^{\sigma ,\mathbf{b}}f\right) \square _{J}^{\omega ,\mathbf{b}%
^{\ast }}gd\omega \\
&\equiv &\Theta _{2}^{\limfunc{bad}}\left( f,g\right) +\Theta _{2}^{\limfunc{%
good}}\left( f,g\right) \ ,
\end{eqnarray*}%
and where if $J^{\maltese }$ fails to exist, we assume by convention that $%
J^{\maltese }\not\subsetneqq I$, i.e. $J^{\maltese }$ is \emph{not} strictly
contained in $I$, so that the pair $\left( I,J\right) $ is then included in
the bad form $\Theta _{2}^{\limfunc{bad}}\left( f,g\right) $. We will in
fact estimate a larger quantity corresponding to the bad form, namely%
\begin{equation}
\Theta _{2}^{\limfunc{bad}\natural }\left( f,g\right) \equiv \sum_{I\in 
\mathcal{D}}\sum_{\substack{ J\in \mathcal{G}:\ J^{\maltese }\not\subsetneqq
I\text{, }\ell \left( J\right) \leq 2^{-\mathbf{r}}\ell \left( I\right)  \\ %
d\left( J,I\right) \leq 2\ell \left( J\right) ^{\varepsilon }\ell \left(
I\right) ^{1-\varepsilon }}}\left\vert \int \left( T_{\sigma }^{\alpha
}\square _{I}^{\sigma ,\mathbf{b}}f\right) \square _{J}^{\omega ,\mathbf{b}%
^{\ast }}gd\omega \right\vert  \label{Theta_2^bad sharp}
\end{equation}%
with absolute value signs \emph{inside} the sum.

\begin{remark}
We now make some general comments on where we now stand and where we are
going.

\begin{enumerate}
\item In the first sum $\Theta _{2}^{\limfunc{bad}}\left( f,g\right) $
above, we are roughly keeping the pairs of intervals $\left( I,J\right) $
such that $J$ is $\limfunc{bad}$ with respect to some `nearby' interval
having side length larger than that of $I$.

\item We have defined energy and dual energy conditions that are independent
of the testing families (because the definition of $\mathsf{E}\left(
J,\omega \right) =\mathbb{E}_{J}^{\omega ,x}\mathbb{E}_{J}^{\omega
,x^{\prime }}\left( \left\vert \frac{x-x^{\prime }}{\ell \left( J\right) }%
\right\vert ^{2}\right) $ does not involve pseudoprojections $\square _{J,%
\mathcal{D}}^{\omega ,\mathbf{b}^{\ast }}$), but the functional energy
condition defined below \emph{does} involve the dual martingale
pseudoprojections $\square _{J,\mathcal{D}}^{\omega ,\mathbf{b}^{\ast }}$.

\item Using the notion of weak goodness above, we will be able to eliminate
all pairs of intervals with $J$ bad in $I$, which then permits control of
the short range form in Section \ref{Sec disj form} and the neighbour form
in Section \ref{Sec Main below} provided $0<\varepsilon <\frac{1}{2-\alpha }$%
. Defining shifted coronas in terms of $J^{\maltese }$ will then allow
existing arguments to prove the Intertwining Proposition and obtain control
of the functional energy in Appendix B, as well as permitting control of the
stopping form in Section \ref{Sec stop}, but all of this with some new
twists, for example the introduction of a top/down `indented corona' in the
analysis of the stopping form.

\item The nearby form $\Theta _{3}\left( f,g\right) $ is handled in Section %
\ref{Sec nearby} using the energy condition assumption along with the
original testing functions $b_{Q}^{\limfunc{orig}}$ discarded during the
construction of the testing/accretive corona.
\end{enumerate}
\end{remark}

These remarks will become clear in this and the following sections. Recall
that we earlier defined in Definition \ref{good two grids}, the set $%
\mathcal{G}_{k-\limfunc{good}}^{\mathcal{D}}=\mathcal{G}_{\left(
k,\varepsilon \right) -\limfunc{good}}^{\mathcal{D}}$ to consist of those $%
J\in \mathcal{G}$ such that $J$ is $\varepsilon -\limfunc{good}$ inside
every interval $K\in \mathcal{D}$ with $K\cap J\neq \emptyset $ that lies at
least $k$ levels `above' $J$, i.e. $\ell \left( K\right) \geq 2^{k}\ell
\left( J\right) $. We now define an analogous notion of $\mathcal{G}_{k-%
\limfunc{bad}}^{\mathcal{D}}$.

\begin{definition}
\label{def Gbad}Let $\varepsilon >0$. Define the set $\mathcal{G}_{k-%
\limfunc{bad}}^{\mathcal{D}}=\mathcal{G}_{\left( k,\varepsilon \right) -%
\limfunc{bad}}^{\mathcal{D}}$ to consist of all $J\in \mathcal{G}$ such that
there is a $\mathcal{D}$-interval $K$ with sidelength $\ell \left( K\right)
=2^{k}\ell \left( J\right) $ for which $J$ is $\varepsilon -\limfunc{bad}$
with respect to $K$.
\end{definition}

Note that for grids $\mathcal{D}$ and $\mathcal{G}$, the complement of $%
\mathcal{G}_{k-\limfunc{good}}^{\mathcal{D}}$ is the union of $\mathcal{G}%
_{\ell -\limfunc{bad}}^{\mathcal{D}}$ for $\ell \geq k$, i.e.%
\begin{equation*}
\mathcal{G\setminus G}_{k-\limfunc{good}}^{\mathcal{D}}=\bigcup_{\ell \geq k}%
\mathcal{G}_{\ell -\limfunc{bad}}^{\mathcal{D}}\ .
\end{equation*}%
Now assume $\varepsilon >0$. We then have the following important property,
namely for all intervals $R$, and all $k\geq \mathbf{r}$ (where the goodness
parameter $\mathbf{r}$ will be fixed given $\varepsilon >0$ in (\ref{choice
of r}) below): 
\begin{equation}
\#\left\{ Q:\kappa \left( Q,R\right) =k\text{ and }d\left( R,Q\right) \leq
2\ell \left( R\right) ^{\varepsilon }\ell \left( Q\right) ^{1-\varepsilon
}\right\} \lesssim 1.  \label{imp}
\end{equation}%
As in \cite{HyMa}, set%
\begin{equation*}
\mathcal{G}_{\limfunc{bad},n}^{\mathcal{D}}\equiv \left\{ J\in \mathcal{G}:J%
\text{ is }\varepsilon -\limfunc{bad}\ \text{with respect to some }K\in 
\mathcal{D}\text{ with }\ell \left( K\right) \geq n\right\} .
\end{equation*}%
We will now use the set equality%
\begin{eqnarray}
&&\left\{ J\in \mathcal{G}:\ J^{\maltese }\not\subset I,\text{ }\ell \left(
J\right) \leq 2^{-\mathbf{r}}\ell \left( I\right) ,\ d\left( J,I\right) \leq
2\ell \left( J\right) ^{\varepsilon }\ell \left( I\right) ^{1-\varepsilon
}\right\}  \label{set equ} \\
&=&\left\{ R\in \mathcal{G}_{\limfunc{bad},\ell \left( Q\right) }^{\mathcal{D%
}}:\ \mathbf{r}\leq \kappa \left( Q,R\right) <\kappa \left( R\right) ,\
d\left( R,Q\right) \leq 2\ell \left( R\right) ^{\varepsilon }\ell \left(
Q\right) ^{1-\varepsilon }\right\} ,  \notag
\end{eqnarray}%
which the careful reader can prove by painstakingly verifying both
containments.

Assuming only that $\mathbf{b}$ is $2$-weakly $\mu $-controlled accretive
(recall we are assuming the stronger condition that $\mathbf{b}$ is $\infty $%
-strongly $\mu $-controlled accretive in our proof here), and following the
proof in \cite{HyMa}, we use (\ref{set equ}) to show that for any fixed
grids $\mathcal{D}$ and $\mathcal{G}$, and any bounded linear operator $%
T_{\sigma }^{\alpha }$ we have the following inequality for the form $\Theta
_{2}^{\limfunc{bad}\natural ,\limfunc{strict}}\left( f,g\right) $, defined
to be $\Theta _{2}^{\limfunc{bad}\natural }\left( f,g\right) $ as in (\ref%
{Theta_2^bad sharp}) with the pairs $\left( I,J\right) $ removed when $%
J^{\maltese }=I$. We use $\varepsilon _{Q,R}=\pm 1$ to obtain 
\begin{eqnarray*}
&&\Theta _{2}^{\limfunc{bad}\natural ,\limfunc{strict}}\left( f,g\right)
=\sum_{Q\in \mathcal{D}}\sum_{\substack{ R\in \mathcal{G}_{\limfunc{bad}%
,\ell \left( Q\right) }^{\mathcal{D}}:\ \mathbf{r}\leq \kappa \left(
Q,R\right) <\kappa \left( R\right)  \\ d\left( R,Q\right) \leq 2\ell \left(
R\right) ^{\varepsilon }\ell \left( Q\right) ^{1-\varepsilon }}}\left\vert
\left\langle T_{\sigma }^{\alpha }\left( \square _{Q,\mathcal{D}}^{\sigma ,%
\mathbf{b}}f\right) ,\square _{R,\mathcal{G}}^{\omega ,\mathbf{b}^{\ast
}}g\right\rangle \right\vert \\
&=&\sum_{Q\in \mathcal{D}}\sum_{\substack{ R\in \mathcal{G}_{\limfunc{bad}%
,\ell \left( Q\right) }^{\mathcal{D}}:\ \mathbf{r}\leq \kappa \left(
Q,R\right) <\kappa \left( R\right)  \\ d\left( R,Q\right) \leq 2\ell \left(
R\right) ^{\varepsilon }\ell \left( Q\right) ^{1-\varepsilon }}}\varepsilon
_{Q,R}\left\langle T_{\sigma }^{\alpha }\left( \square _{Q,\mathcal{D}%
}^{\sigma ,\mathbf{b}}f\right) ,\square _{R,\mathcal{G}}^{\omega ,\mathbf{b}%
^{\ast }}g\right\rangle \\
&\leq &\sum_{Q\in \mathcal{D}}\left\vert \left\langle T_{\sigma }^{\alpha
}\left( \square _{Q,\mathcal{D}}^{\sigma ,\mathbf{b}}f\right) ,\sum 
_{\substack{ R\in \mathcal{G}_{\limfunc{bad},\ell \left( Q\right) }^{%
\mathcal{D}}:\ \mathbf{r}\leq \kappa \left( Q,R\right) <\kappa \left(
R\right)  \\ d\left( R,Q\right) \leq 2\ell \left( R\right) ^{\varepsilon
}\ell \left( Q\right) ^{1-\varepsilon }}}\varepsilon _{Q,R}\square _{R,%
\mathcal{G}}^{\omega ,\mathbf{b}^{\ast }}g\right\rangle \right\vert \\
&\leq &\mathfrak{N}_{T^{\alpha }}\sum_{Q\in \mathcal{D}}\left\Vert \square
_{Q,\mathcal{D}}^{\sigma ,\mathbf{b}}f\right\Vert _{L^{2}\left( \sigma
\right) }\left\Vert \sum_{\substack{ R\in \mathcal{G}_{\limfunc{bad},\ell
\left( Q\right) }^{\mathcal{D}}:\ \mathbf{r}\leq \kappa \left( Q,R\right)
<\kappa \left( R\right)  \\ d\left( R,Q\right) \leq 2\ell \left( R\right)
^{\varepsilon }\ell \left( Q\right) ^{1-\varepsilon }}}\varepsilon
_{Q,R}\square _{R,\mathcal{G}}^{\omega ,\mathbf{b}^{\ast }}g\right\Vert
_{L^{2}\left( \omega \right) } \\
&\leq &\mathfrak{N}_{T^{\alpha }}\sum_{Q\in \mathcal{D}}\left\Vert \square
_{Q,\mathcal{D}}^{\sigma ,\mathbf{b}}f\right\Vert _{L^{2}\left( \sigma
\right) }\sum_{k=\mathbf{r}}^{\infty }\left\Vert \sum_{\substack{ R\in 
\mathcal{G}_{\limfunc{bad},\ell \left( Q\right) }^{\mathcal{D}}:k=\kappa
\left( Q,R\right) <\kappa \left( R\right)  \\ d\left( R,Q\right) \leq 2\ell
\left( R\right) ^{\varepsilon }\ell \left( Q\right) ^{1-\varepsilon }}}%
\varepsilon _{Q,R}\square _{R,\mathcal{G}}^{\omega ,\mathbf{b}^{\ast
}}g\right\Vert _{L^{2}\left( \omega \right) }\ ,
\end{eqnarray*}%
by Minkowski's inequality, and we continue with%
\begin{eqnarray*}
&\leq &2\mathfrak{N}_{T^{\alpha }}\sum_{k=\mathbf{r}}^{\infty }\left(
\sum_{Q\in \mathcal{D}}\left\Vert \square _{Q,\mathcal{D}}^{\sigma ,\mathbf{b%
}}f\right\Vert _{L^{2}\left( \sigma \right) }^{2}\right) ^{\frac{1}{2}%
}\left( \sum_{Q\in \mathcal{D}}\sum_{\substack{ R\in \mathcal{G}_{\limfunc{%
bad},\ell \left( Q\right) }^{\mathcal{D}}:\ k=\kappa \left( Q,R\right)
<\kappa \left( R\right)  \\ d\left( R,Q\right) \leq 2\ell \left( R\right)
^{\varepsilon }\ell \left( Q\right) ^{1-\varepsilon }}}\left( \left\Vert
\square _{R,\mathcal{G}}^{\omega ,\mathbf{b}^{\ast }}g\right\Vert
_{L^{2}\left( \omega \right) }^{2}+\left\Vert \nabla _{R,\mathcal{G}%
}^{\omega }g\right\Vert _{L^{2}\left( \omega \right) }^{2}\right) \right) ^{%
\frac{1}{2}} \\
&\lesssim &\mathfrak{N}_{T^{\alpha }}\left\Vert f\right\Vert _{L^{2}\left(
\sigma \right) }\sum_{k=\mathbf{r}}^{\infty }\left( \sum_{R\in \mathcal{G}_{%
\limfunc{bad},2^{k}\ell \left( R\right) }^{\mathcal{D}}}\left( \left\Vert
\square _{R,\mathcal{G}}^{\omega ,\mathbf{b}^{\ast }}g\right\Vert
_{L^{2}\left( \omega \right) }^{2}+\left\Vert \nabla _{R,\mathcal{G}%
}^{\omega }g\right\Vert _{L^{2}\left( \omega \right) }^{2}\right) \right) ^{%
\frac{1}{2}},
\end{eqnarray*}%
where $\nabla _{R,\mathcal{G}}^{\omega }$ denotes the `broken' Carleson
averaging operator in (\ref{Carleson avg op}) that depends on the grid $%
\mathcal{G}$, and

\begin{enumerate}
\item the penultimate inequality uses Cauchy-Schwarz in $Q$ and the weak
upper Riesz inequalities (\ref{UPPER RIESZ}) for $\sum_{\substack{ R\in 
\mathcal{G}_{\limfunc{bad},\ell \left( Q\right) }^{\mathcal{D}}:\ k=\kappa
\left( Q,R\right) <\kappa \left( R\right)  \\ d\left( R,Q\right) \leq 2\ell
\left( R\right) ^{\varepsilon }\ell \left( Q\right) ^{1-\varepsilon }}}%
\varepsilon _{Q,R}\square _{R,\mathcal{G}}^{\omega ,\mathbf{b}^{\ast }}$,
once for the sum when $\varepsilon _{Q,R}=1$, and again for the sum when $%
\varepsilon _{Q,R}=-1$. However, we note that since the sum in $R$ is
pigeonholed by $k=\kappa \left( Q,R\right) $, the $R$'s are pairwise
disjoint intervals and the pseudoprojections $\square _{R,\mathcal{G}%
}^{\omega ,\mathbf{b}^{\ast }}g$ are pairwise orthogonal. Thus we could
instead apply Cauchy-Schwarz first in $R$, and then in $Q$ as was done in 
\cite{HyMa}, but we must still apply weak upper Riesz inequalities as above.

\item and the final inequality uses the frame inequality (\ref{FRAME})
together with (\ref{imp}), namely the fact that there are at most $C$
intervals $Q$ such that $\kappa \left( Q,R\right) \geq \mathbf{r}$ is fixed
and $d\left( R,Q\right) \leq 2\ell \left( R\right) ^{\varepsilon }\ell
\left( Q\right) ^{1-\varepsilon }$.
\end{enumerate}

Now it is easy to verify that we have the same inequality for the pairs $%
\left( J^{\maltese },J\right) $ that were removed, and then we take grid
expectations and use the probability estimate (\ref{main bad prob}) to
obtain for $\varepsilon ^{\prime }=\frac{1}{2}\varepsilon $ that 
\begin{eqnarray}
&&\boldsymbol{E}_{\Omega }^{\mathcal{D}}\left( \Theta _{2}^{\limfunc{bad}%
\natural }\left( f,g\right) \right)  \label{HM bad} \\
&\leq &\boldsymbol{E}_{\Omega }^{\mathcal{D}}\mathfrak{N}_{T^{\alpha
}}\left\Vert f\right\Vert _{L^{2}\left( \sigma \right) }\sum_{k=\mathbf{r}%
}^{\infty }\left( \sum_{R\in \mathcal{G}_{\limfunc{bad},2^{k}\ell \left(
R\right) }^{\mathcal{D}}}\left( \left\Vert \square _{R,\mathcal{G}}^{\omega ,%
\mathbf{b}^{\ast }}g\right\Vert _{L^{2}\left( \omega \right)
}^{2}+\left\Vert \nabla _{R,\mathcal{G}}^{\omega }g\right\Vert _{L^{2}\left(
\omega \right) }^{2}\right) \right) ^{\frac{1}{2}}  \notag \\
&\leq &\mathfrak{N}_{T^{\alpha }}\left\Vert f\right\Vert _{L^{2}\left(
\sigma \right) }\sum_{k=\mathbf{r}}^{\infty }\left( \boldsymbol{E}_{\Omega
}^{\mathcal{D}}\sum_{R\in \mathcal{G}_{\limfunc{bad},2^{k}\ell \left(
R\right) }^{\mathcal{D}}}\left( \left\Vert \square _{R,\mathcal{G}}^{\omega ,%
\mathbf{b}^{\ast }}g\right\Vert _{L^{2}\left( \omega \right)
}^{2}+\left\Vert \nabla _{R,\mathcal{G}}^{\omega }g\right\Vert _{L^{2}\left(
\omega \right) }^{2}\right) \right) ^{\frac{1}{2}}  \notag \\
&\lesssim &2^{-\frac{1}{2}\varepsilon ^{\prime }\mathbf{r}}\mathfrak{N}%
_{T^{\alpha }}\left\Vert f\right\Vert _{L^{2}\left( \sigma \right) }\sum_{k=%
\mathbf{r}}^{\infty }\left( C_{1}2^{-\varepsilon k}\left\Vert g\right\Vert
_{L^{2}\left( \omega \right) }^{2}\right) ^{\frac{1}{2}}\leq C_{\limfunc{good%
}}2^{-\frac{1}{2}\varepsilon \mathbf{r}}\mathfrak{N}_{T^{\alpha }}\left\Vert
f\right\Vert _{L^{2}\left( \sigma \right) }\left\Vert g\right\Vert
_{L^{2}\left( \omega \right) }\ .  \notag
\end{eqnarray}%
Clearly we can now fix $\mathbf{r}$ sufficiently large depending on $%
\varepsilon >0$ so that%
\begin{equation}
C_{\limfunc{good}}2^{-\frac{1}{2}\varepsilon \mathbf{r}}<\frac{1}{100},
\label{choice of r}
\end{equation}%
and then the final term above, namely $C_{\limfunc{good}}2^{-\frac{1}{2}%
\varepsilon \mathbf{r}}\mathfrak{N}_{T^{\alpha }}\left\Vert f\right\Vert
_{L^{2}\left( \sigma \right) }\left\Vert g\right\Vert _{L^{2}\left( \omega
\right) }$, can be absorbed at the end of the proof in Subsection \ref{Sub
wrapup}. Note that (\ref{choice of r}) fixes our choice of the parameter $%
\mathbf{r}$ for any given $\varepsilon >0$. Later we will choose $%
0<\varepsilon <\frac{1}{2}\leq \frac{1}{2-\alpha }$. It is this type of weak
goodness that we will exploit in the local forms $\mathsf{B}_{\Subset _{%
\mathbf{r}}}^{A}\left( f,g\right) $ treated below in Section \ref{Sec Main
below}.

We are now left with the following `good' form to control:%
\begin{equation*}
\Theta _{2}^{\limfunc{good}}\left( f,g\right) =\sum_{I\in \mathcal{D}}\sum 
_{\substack{ J^{\maltese }\subsetneqq I:\ \ell \left( J\right) \leq 2^{-%
\mathbf{r}}\ell \left( I\right)  \\ d\left( J,I\right) \leq 2\ell \left(
J\right) ^{\varepsilon }\ell \left( I\right) ^{1-\varepsilon }}}\int \left(
T_{\sigma }^{\alpha }\square _{I}^{\sigma ,\mathbf{b}}f\right) \square
_{J}^{\omega ,\mathbf{b}^{\ast }}gd\omega .
\end{equation*}%
The first thing we observe regarding this form is that the intervals $J$
which arise in the sum for $\Theta _{2}^{\limfunc{good}}\left( f,g\right) $
must lie entirely inside $I$ since $J\subset J^{\maltese }\subsetneqq I$.
Then in the remainder of the paper, we proceed to analyze 
\begin{equation}
\Theta _{2}^{\limfunc{good}}\left( f,g\right) =\sum_{I\in \mathcal{D}%
}\sum_{J^{\maltese }\subsetneqq I:\ \ell \left( J\right) \leq 2^{-\mathbf{r}%
}\ell \left( I\right) }\int \left( T_{\sigma }^{\alpha }\square _{I}^{\sigma
,\mathbf{b}}f\right) \square _{J}^{\omega ,\mathbf{b}^{\ast }}gd\omega ,
\label{def Theta 2 good}
\end{equation}%
in the same way we analyzed the below term $\mathsf{B}_{\Subset _{\mathbf{r}%
}}\left( f,g\right) $ in \cite{SaShUr6}; namely, by implementing the
canonical corona splitting and the decomposition into paraproduct, neighbour
and stopping forms, but now with an additional broken form. We have $\left(
\kappa ,\varepsilon \right) $-goodness available for all the intervals $J\in 
\mathcal{G}$ arising in the form $\Theta _{2}^{\limfunc{good}}\left(
f,g\right) $, and moreover, the intervals $I\in \mathcal{D}$ arising in the
form $\Theta _{2}^{\limfunc{good}}\left( f,g\right) $ for a fixed $J$ are
tree-connected, so that telescoping identities hold for these intervals $I$.
This will prove decisive in the following three sections of the paper.

The forms $\Theta _{1}\left( f,g\right) $ and $\Theta _{3}\left( f,g\right) $
are analogous to the disjoint and nearby forms $\mathsf{B}_{\cap }\left(
f,g\right) $ and $\mathsf{B}_{/}\left( f,g\right) $ in \cite{SaShUr6}
respectively. In the next\ two sections, we control the disjoint form $%
\Theta _{1}\left( f,g\right) $ in essentially the same way that the disjoint
form $\mathsf{B}_{\cap }\left( f,g\right) $ was treated in \cite{SaShUr6}
and in earlier papers of many authors beginning with Nazarov, Treil and
Volberg (see e.g. \cite{Vol}), and we control the nearby form $\Theta
_{3}\left( f,g\right) $ using the probabilistic surgery of Hyt\"{o}nen and
Martikainen building on that of NTV, together with a new deterministic
surgery involving the energy condition and the original testing functions.
But first we recall, in the following subsection, the characterization of
boundedness of one-dimensional forms supported on disjoint intervals \cite%
{Hyt2}.

\subsection{A characterization of bilinear forms supported on disjoint
intervals \label{disjoint sets}}

Matters here in the one-dimensional setting are greatly simplified by a
generalization to fractional integrals of Hyt\"{o}nen's characterization of
the \emph{restricted bilinear inequality,}%
\begin{equation}
\left\vert \int_{\mathbb{R}\setminus I}\left( \int_{I}\frac{f\left( y\right) 
}{\left\vert x-y\right\vert }d\sigma \left( y\right) \right) g\left(
x\right) d\omega \left( x\right) \right\vert \lesssim \mathfrak{D}\left\Vert
f\right\Vert _{L^{2}\left( \sigma \right) }\left\Vert g\right\Vert
_{L^{2}\left( \omega \right) }\ ,  \label{bilin disjoint}
\end{equation}%
for all intervals $I$, in terms of the Muckenhoupt conditions, namely%
\begin{equation*}
\mathfrak{D}\approx \sqrt{\mathcal{A}_{2}}+\sqrt{\mathcal{A}_{2}^{\ast }},
\end{equation*}%
where $\mathfrak{D}$ is the best constant in (\ref{bilin disjoint}). In \cite%
{HyMa} this inequality was proved for complementary half-lines, where it was
pointed out that the passage to an interval and its complement is then
routine.

We claim that Hyt\"{o}nen's characterization extends immediately to
fractional integrals on the line with the same proof. Namely, we have, 
\begin{eqnarray}
\left\vert \int_{\left( -\infty ,a\right) }\left( \int_{\left( a,\infty
\right) }\frac{f\left( y\right) }{\left\vert x-y\right\vert ^{1-\alpha }}%
d\sigma \left( y\right) \right) g\left( x\right) d\omega \left( x\right)
\right\vert &\lesssim &\left( \sqrt{\mathcal{A}_{2}^{\alpha }}+\sqrt{%
\mathcal{A}_{2}^{\alpha ,\ast }}\right) \left\Vert f\right\Vert
_{L^{2}\left( \sigma \right) }\left\Vert g\right\Vert _{L^{2}\left( \omega
\right) }\ ,  \label{disj supp} \\
\left\vert \int_{\mathbb{R}\setminus I}\left( \int_{I}\frac{f\left( y\right) 
}{\left\vert x-y\right\vert ^{1-\alpha }}d\sigma \left( y\right) \right)
g\left( x\right) d\omega \left( x\right) \right\vert &\lesssim &\left( \sqrt{%
\mathcal{A}_{2}^{\alpha }}+\sqrt{\mathcal{A}_{2}^{\alpha ,\ast }}\right)
\left\Vert f\right\Vert _{L^{2}\left( \sigma \right) }\left\Vert
g\right\Vert _{L^{2}\left( \omega \right) }\ ,  \notag
\end{eqnarray}%
and that $\sqrt{\mathcal{A}_{2}^{\alpha }}+\sqrt{\mathcal{A}_{2}^{\alpha
,\ast }}\approx \mathfrak{D}^{\alpha }$ where $\mathfrak{D}^{\alpha }$ is
the best constant in the inequality above (a limiting argument shows that we
may take one of the half-lines to be closed in (\ref{disj supp})). First,
the proof that $\sqrt{\mathcal{A}_{2}^{\alpha }}+\sqrt{\mathcal{A}%
_{2}^{\alpha ,\ast }}\lesssim \mathfrak{D}^{\alpha }$ is the standard proof
of necessity of the one-tailed Muckenhoupt conditions. In the other
direction, we use the general two weight Hardy inequality of Muckenhoupt as
presented in \cite[Theorem 3.3]{Hyt2}, see also \cite{LaSaUr2}: if $\sigma $
and $\omega $ are locally finite positive Borel measures on the interval $%
\left( 0,\infty \right) $, then%
\begin{equation*}
\int_{0}^{\infty }\left( \int_{\left( 0,x\right] }fd\sigma \right)
^{2}d\omega \left( x\right) \leq C\int_{0}^{\infty }f\left( y\right)
^{2}d\sigma \left( y\right) ,
\end{equation*}%
holds for all $f\in L^{2}\left( \sigma \right) $ if and only if%
\begin{equation*}
A\equiv \sup_{t>0}\left( \int_{\left( 0,t\right] }d\sigma \right) \left(
\int_{\left[ t,\infty \right) }d\omega \right) <\infty .
\end{equation*}%
Moreover, if $C$ is the best constant above, then%
\begin{equation*}
A\leq C\leq 4A.
\end{equation*}%
We easily obtain the following characterization of an intermediate
inequality:

\begin{equation*}
\int_{0}^{\infty }\int_{0}^{\infty }\frac{f\left( y\right) g\left( x\right) 
}{\left( x+y\right) ^{1-\alpha }}d\sigma \left( y\right) d\omega \left(
x\right) \leq C\left\Vert f\right\Vert _{L^{2}\left( \sigma \right)
}\left\Vert g\right\Vert _{L^{2}\left( \omega \right) },
\end{equation*}%
if and only if%
\begin{equation*}
A\equiv \sup_{t>0}\sqrt{\left( \int_{\left( 0,t\right] }d\sigma \right)
\left( \int_{\left[ t,\infty \right) }\frac{d\omega \left( x\right) }{%
x^{2-2\alpha }}\right) }+\sup_{t>0}\sqrt{\left( \int_{\left( 0,t\right]
}d\omega \right) \left( \int_{\left[ t,\infty \right) }\frac{d\sigma \left(
y\right) }{y^{2-2\alpha }}\right) }<\infty ,
\end{equation*}%
and moreover the best constant $C$ satisfies $\frac{1}{4}A\leq C\leq 2A$. To
see this we simply use the estimates 
\begin{equation*}
\frac{1}{2}\max \left\{ \frac{1}{x^{1-\alpha }}\mathbf{1}_{\left( 0,x\right]
}\left( y\right) ,\frac{1}{y^{1-\alpha }}\mathbf{1}_{\left( 0,y\right]
}\left( x\right) \right\} \leq \frac{1}{\left( x+y\right) ^{1-\alpha }}\leq 
\frac{1}{x^{1-\alpha }}\mathbf{1}_{\left( 0,x\right] }\left( y\right) +\frac{%
1}{y^{1-\alpha }}\mathbf{1}_{\left( 0,y\right] }\left( x\right) ,
\end{equation*}%
together with Hardy's inequality. From this and duality, we immediately
obtain (\ref{disj supp}).

\subsubsection{Control of triple testing and triple energy}

We also define the \emph{triple} $\mathbf{b}$-testing conditions for $%
T^{\alpha }$ and \emph{triple} $\mathbf{b}^{\ast }$-testing conditions for
the dual $T^{\alpha ,\ast }$ given by%
\begin{eqnarray}
\int_{3Q}\left\vert T_{\sigma }^{\alpha }b_{Q}\right\vert ^{2}d\omega &\leq
&\left( \mathfrak{3T}_{T^{\alpha }}^{\mathbf{b}}\right) ^{2}\left\vert
Q\right\vert _{\sigma }\ ,\ \ \ \ \ \text{for all intervals }Q,
\label{triple b testing cond} \\
\int_{3Q}\left\vert T_{\omega }^{\alpha ,\ast }b_{Q}^{\ast }\right\vert
^{2}d\sigma &\leq &\left( \mathfrak{3T}_{T^{\alpha }}^{\mathbf{b}^{\ast
},\ast }\right) ^{2}\left\vert Q\right\vert _{\omega }\ ,\ \ \ \ \ \text{for
all intervals }Q,  \notag
\end{eqnarray}%
as well as the \emph{full} $\mathbf{b}$-testing conditions for $T^{\alpha }$
and \emph{full} $\mathbf{b}^{\ast }$-testing conditions for the dual $%
T^{\alpha ,\ast }$ given by%
\begin{eqnarray}
\int_{\mathbb{R}}\left\vert T_{\sigma }^{\alpha }b_{Q}\right\vert
^{2}d\omega &\leq &\left( \mathfrak{FT}_{T^{\alpha }}^{\mathbf{b}}\right)
^{2}\left\vert Q\right\vert _{\sigma }\ ,\ \ \ \ \ \text{for all intervals }%
Q,  \label{full b testing} \\
\int_{\mathbb{R}}\left\vert T_{\omega }^{\alpha ,\ast }b_{Q}^{\ast
}\right\vert ^{2}d\sigma &\leq &\left( \mathfrak{FT}_{T^{\alpha }}^{\mathbf{b%
}^{\ast },\ast }\right) ^{2}\left\vert Q\right\vert _{\omega }\ ,\ \ \ \ \ 
\text{for all intervals }Q.  \notag
\end{eqnarray}%
Note that the full testing conditions are implied by the triple testing
conditions and the Muckenhoupt conditions (e.g. use the above
characterization on complementary half-lines),%
\begin{equation*}
\mathfrak{FT}_{T^{\alpha }}^{\mathbf{b}}\lesssim \mathfrak{3T}_{T^{\alpha
}}^{\mathbf{b}}+\sqrt{\mathcal{A}_{2}^{\alpha }}\text{ and }\mathfrak{FT}%
_{T^{\alpha ,\ast }}^{\mathbf{b}^{\ast }}\lesssim \mathfrak{3T}_{T^{\alpha
,\ast }}^{\mathbf{b}^{\ast }}+\sqrt{\mathcal{A}_{2}^{\alpha ,\ast }}.
\end{equation*}%
Since dimension $n=1$, the full testing conditions are controlled by testing
and Muckenhoupt, as we now show. Indeed, if we now set $f=b_{I}$ in the
second line of (\ref{disj supp}), and take the supremum over all $g\in
L^{2}\left( \omega \right) $ with $\left\Vert g\right\Vert _{L^{2}\left(
\omega \right) }=1$, we obtain%
\begin{eqnarray*}
\sqrt{\int_{\mathbb{R}\setminus I}\left( \int_{I}\frac{f\left( y\right) }{%
\left\vert x-y\right\vert ^{1-\alpha }}d\sigma \left( y\right) \right)
^{2}d\omega \left( x\right) } &\lesssim &\left( \sqrt{\mathcal{A}%
_{2}^{\alpha }}+\sqrt{\mathcal{A}_{2}^{\alpha ,\ast }}\right) \left\Vert
b_{I}\right\Vert _{L^{2}\left( \sigma \right) } \\
&\lesssim &\left( \sqrt{\mathcal{A}_{2}^{\alpha }}+\sqrt{\mathcal{A}%
_{2}^{\alpha ,\ast }}\right) \sqrt{\left\vert I\right\vert _{\sigma }},
\end{eqnarray*}%
which gives%
\begin{eqnarray*}
\int_{\mathbb{R}}\left\vert T_{\sigma }^{\alpha }b_{I}\right\vert
^{2}d\omega \left( x\right) &=&\int_{I}\left\vert T_{\sigma }^{\alpha
}b_{I}\right\vert ^{2}d\omega \left( x\right) +\int_{\mathbb{R}\setminus
I}\left\vert T_{\sigma }^{\alpha }b_{I}\right\vert ^{2}d\omega \left(
x\right) \\
&\lesssim &\left( \mathfrak{T}_{T^{\alpha }}^{\mathbf{b}}\right)
^{2}\left\vert I\right\vert _{\sigma }+\int_{\mathbb{R}\setminus I}\left(
\int_{I}\frac{f\left( y\right) }{\left\vert x-y\right\vert ^{1-\alpha }}%
d\sigma \left( y\right) \right) ^{2}d\omega \left( x\right) \\
&\lesssim &\left\{ \left( \mathfrak{T}_{T^{\alpha }}^{\mathbf{b}}\right)
^{2}+\mathcal{A}_{2}^{\alpha }+\mathcal{A}_{2}^{\alpha ,\ast }\right\}
\left\vert I\right\vert _{\sigma }\ .
\end{eqnarray*}%
Thus we have obtained control of full $\mathbf{b}$-testing by just $\mathbf{b%
}$-testing and the Muckenhoupt conditions in dimension $n=1$:%
\begin{equation}
\mathfrak{FT}_{T^{\alpha }}^{\mathbf{b}}\lesssim \mathfrak{T}_{T^{\alpha }}^{%
\mathbf{b}}+\sqrt{\mathcal{A}_{2}^{\alpha }}+\sqrt{\mathcal{A}_{2}^{\alpha
,\ast }}\text{ and }\mathfrak{FT}_{T^{\alpha ,\ast }}^{\mathbf{b}^{\ast
}}\lesssim \mathfrak{T}_{T^{\alpha ,\ast }}^{\mathbf{b}^{\ast }}+\sqrt{%
\mathcal{A}_{2}^{\alpha }}+\sqrt{\mathcal{A}_{2}^{\alpha ,\ast }}\text{ }.
\label{full proved}
\end{equation}

Now we turn to the analogous notion of \emph{triple} energy conditions
defined in analogy with the triple testing condtions. Namely, the sum over
the intervals $I_{r}$ in the energy condition in Definition \ref{def strong
quasienergy} is permitted to extend to the triple $3I$ of the interval $I$,
but with the additional proviso that the distance of $I_{r}$ from the
boundary of $I$ is at least a positive multiple $\delta $ (the exact value
of which is immaterial) of the side length of $I_{r}$.

\begin{definition}
\label{def triple energy}Let $0\leq \alpha <1$ and $0<\delta \leq \frac{1}{2}
$. Suppose $\sigma $ and $\omega $ are locally finite positive Borel
measures on $\mathbb{R}$. Then the \emph{triple} energy constant $\mathcal{E}%
_{2}^{\alpha ,\limfunc{triple}}$ is defined by 
\begin{equation*}
\left( \mathcal{E}_{2}^{\alpha ,\limfunc{triple}}\right) ^{2}\equiv \sup 
_{\substack{ 3I=\dot{\cup}I_{r}  \\ d\left( I_{r},\partial I\right) \geq
\delta \ell \left( I_{r}\right) \text{ when }I_{r}\cap I^{c}=\emptyset }}%
\frac{1}{\left\vert I\right\vert _{\sigma }}\sum_{r=1}^{\infty }\left( \frac{%
\mathrm{P}^{\alpha }\left( I_{r},\mathbf{1}_{I}\sigma \right) }{\left\vert
I_{r}\right\vert }\right) ^{2}\left\Vert x-m_{I_{r}}\right\Vert
_{L^{2}\left( \mathbf{1}_{I_{r}}\omega \right) }^{2}\ ,
\end{equation*}%
where the supremum is taken over arbitrary decompositions of the triple $3I$
of an interval $I$ using a pairwise disjoint union of subintervals $I_{r}$
whose distance to the boundary of $I$ is at least a positive multiple of $%
\ell \left( I_{r}\right) $ when $I_{r}$ is not contained in $I$. Similarly,
we define the \emph{dual triple} energy constant $\mathcal{E}_{2}^{\alpha ,%
\limfunc{triple},\ast }$ by switching the roles of $\sigma $ and $\omega $:%
\begin{equation*}
\left( \mathcal{E}_{2}^{\alpha ,\limfunc{triple},\ast }\right) ^{2}\equiv
\sup_{\substack{ 3I=\dot{\cup}I_{r}  \\ d\left( I_{r},\partial I\right) \geq
\delta \ell \left( I_{r}\right) \text{ when }I_{r}\cap I^{c}=\emptyset }}%
\sum_{r=1}^{\infty }\left( \frac{\mathrm{P}^{\alpha }\left( I_{r},\mathbf{1}%
_{I}\omega \right) }{\left\vert I_{r}\right\vert }\right) ^{2}\left\Vert
x-m_{I_{r}}\right\Vert _{L^{2}\left( \mathbf{1}_{I_{r}}\sigma \right) }^{2}\
.
\end{equation*}
\end{definition}

We now show that in dimension $n=1$, the triple energy conditions are
controlled by the energy and Muckenhoupt conditions, namely%
\begin{equation}
\mathcal{E}_{2}^{\alpha ,\limfunc{triple}}+\mathcal{E}_{2}^{\alpha ,\limfunc{%
triple},\ast }\lesssim \mathcal{E}_{2}^{\alpha }+\mathcal{E}_{2}^{\alpha
,\ast }+\sqrt{\mathfrak{A}_{2}^{\alpha }}  \label{triple energy control}
\end{equation}%
Indeed, assuming for convenience that $\delta =1$, we need only control by $%
C\left\vert I\right\vert _{\sigma }$ the sum,%
\begin{equation*}
\sum_{r=1}^{\infty }\left( \frac{\mathrm{P}^{\alpha }\left( I_{r},\mathbf{1}%
_{I}\sigma \right) }{\left\vert I_{r}\right\vert }\right) ^{2}\left\Vert
x-m_{I_{r}}\right\Vert _{L^{2}\left( \mathbf{1}_{I_{r}}\omega \right) }^{2}\
,
\end{equation*}
over adjacent intervals $J$ and $I$ of equal length where $\left\{
I_{r}\right\} _{r=1}^{\infty }$ is a disjoint decomposition of $J=\dot{\cup}%
I_{r}$ with $d\left( I_{r},\partial I\right) \geq \ell \left( I_{r}\right) $
for all $r\geq 1$. However, using reversal of energy for the standard
gradient elliptic operator $T^{\alpha }$ with convolution kernel $K^{\alpha
}\left( x\right) =\frac{x}{\left\vert x\right\vert ^{2-\alpha }}$ (see e.g. 
\cite{SaShUr10}), we have since $2I_{r}\cap I=\emptyset $, 
\begin{eqnarray*}
&&\sum_{r=1}^{\infty }\left( \frac{\mathrm{P}^{\alpha }\left( I_{r},\mathbf{1%
}_{I}\sigma \right) }{\left\vert I_{r}\right\vert }\right) ^{2}\left\Vert
x-m_{I_{r}}\right\Vert _{L^{2}\left( \mathbf{1}_{I_{r}}\omega \right) }^{2}
\\
&=&\frac{1}{2}\sum_{r=1}^{\infty }\left( \frac{\mathrm{P}^{\alpha }\left(
I_{r},\mathbf{1}_{I}\sigma \right) }{\left\vert I_{r}\right\vert }\right)
^{2}\frac{1}{\left\vert I_{r}\right\vert _{\omega }}\int_{I_{r}}\int_{I_{r}}%
\left\vert x-z\right\vert ^{2}d\omega \left( x\right) d\omega \left( z\right)
\\
&\lesssim &\sum_{r=1}^{\infty }\frac{1}{\left\vert I_{r}\right\vert _{\omega
}}\int_{I_{r}}\int_{I_{r}}\left\vert T_{\sigma }^{\alpha }\mathbf{1}%
_{I}\left( x\right) -T_{\sigma }^{\alpha }\mathbf{1}_{I}\left( z\right)
\right\vert ^{2}d\omega \left( x\right) d\omega \left( z\right) \\
&\lesssim &\sum_{r=1}^{\infty }\frac{1}{\left\vert I_{r}\right\vert _{\omega
}}\int_{I_{r}}\int_{I_{r}}\left\vert T_{\sigma }^{\alpha }\mathbf{1}%
_{I}\left( x\right) \right\vert ^{2}d\omega \left( x\right) d\omega \left(
z\right) +\sum_{r=1}^{\infty }\frac{1}{\left\vert I_{r}\right\vert _{\omega }%
}\int_{I_{r}}\int_{I_{r}}\left\vert T_{\sigma }^{\alpha }\mathbf{1}%
_{I}\left( z\right) \right\vert ^{2}d\omega \left( x\right) d\omega \left(
z\right) \\
&\lesssim &\sum_{r=1}^{\infty }\int_{I_{r}}\left\vert T_{\sigma }^{\alpha }%
\mathbf{1}_{I}\right\vert ^{2}d\omega \lesssim \int_{J}\left\vert T_{\sigma
}^{\alpha }\mathbf{1}_{I}\right\vert ^{2}d\omega \lesssim \mathfrak{A}%
_{2}^{\alpha }\left\vert I\right\vert _{\sigma }\ ,
\end{eqnarray*}%
which completes the proof of (\ref{triple energy control}).

Finally, we note that a modification of this last argument also shows that
the energy condition itself is controlled by the $\mathbf{1}$-testing
condition and the Muckenhoupt condition. Indeed, as shown in \cite{SaShUr11}%
, 
\begin{eqnarray*}
&&\sum_{r=1}^{\infty }\left( \frac{\mathrm{P}^{\alpha }\left( I_{r},\mathbf{1%
}_{I}\sigma \right) }{\left\vert I_{r}\right\vert }\right) ^{2}\left\Vert
x-m_{I_{r}}\right\Vert _{L^{2}\left( \mathbf{1}_{I_{r}}\omega \right) }^{2}
\\
&=&\frac{1}{2}\sum_{r=1}^{\infty }\left( \frac{\mathrm{P}^{\alpha }\left(
I_{r},\left[ \mathbf{1}_{I\setminus I_{r}}+\mathbf{1}_{I_{r}}\right] \sigma
\right) }{\left\vert I_{r}\right\vert }\right) ^{2}\frac{1}{\left\vert
I_{r}\right\vert _{\omega }}\int_{I_{r}}\int_{I_{r}}\left\vert
x-z\right\vert ^{2}d\omega \left( x\right) d\omega \left( z\right) \\
&\lesssim &\sum_{r=1}^{\infty }\left( \frac{\mathrm{P}^{\alpha }\left( I_{r},%
\mathbf{1}_{I\setminus I_{r}}\sigma \right) }{\left\vert I_{r}\right\vert }%
\right) ^{2}\frac{1}{\left\vert I_{r}\right\vert _{\omega }}%
\int_{I_{r}}\int_{I_{r}}\left\vert x-z\right\vert ^{2}d\omega \left(
x\right) d\omega \left( z\right) +\sum_{r=1}^{\infty }A_{2}^{\alpha ,%
\limfunc{energy}}\left\vert I_{r}\right\vert _{\sigma } \\
&\lesssim &\sum_{r=1}^{\infty }\frac{1}{\left\vert I_{r}\right\vert _{\omega
}}\int_{I_{r}}\int_{I_{r}}\left\vert T_{\sigma }^{\alpha }\mathbf{1}%
_{I\setminus I_{r}}\left( x\right) -T_{\sigma }^{\alpha }\mathbf{1}%
_{I\setminus I_{r}}\left( z\right) \right\vert ^{2}d\omega \left( x\right)
d\omega \left( z\right) +\mathfrak{A}_{2}^{\alpha }\left\vert I\right\vert
_{\sigma } \\
&\lesssim &\sum_{r=1}^{\infty }\frac{1}{\left\vert I_{r}\right\vert _{\omega
}}\int_{I_{r}}\int_{I_{r}}\left( \left\vert T_{\sigma }^{\alpha }\mathbf{1}%
_{I\setminus I_{r}}\left( x\right) \right\vert ^{2}+\left\vert T_{\sigma
}^{\alpha }\mathbf{1}_{I\setminus I_{r}}\left( z\right) \right\vert
^{2}\right) d\omega \left( x\right) d\omega \left( z\right) +\mathfrak{A}%
_{2}^{\alpha }\left\vert I\right\vert _{\sigma } \\
&\lesssim &\sum_{r=1}^{\infty }\int_{I_{r}}\left\vert T_{\sigma }^{\alpha }%
\mathbf{1}_{I\setminus I_{r}}\left( x\right) \right\vert ^{2}d\omega \left(
x\right) +\mathfrak{A}_{2}^{\alpha }\left\vert I\right\vert _{\sigma }\ ,
\end{eqnarray*}%
and now we `plug the hole' to continue with%
\begin{eqnarray*}
\sum_{r=1}^{\infty }\int_{I_{r}}\left\vert T_{\sigma }^{\alpha }\mathbf{1}%
_{I\setminus I_{r}}\left( x\right) \right\vert ^{2}d\omega \left( x\right)
&\lesssim &\sum_{r=1}^{\infty }\int_{I_{r}}\left\vert T_{\sigma }^{\alpha }%
\mathbf{1}_{I}\left( x\right) \right\vert ^{2}d\omega \left( x\right)
+\sum_{r=1}^{\infty }\int_{I_{r}}\left\vert T_{\sigma }^{\alpha }\mathbf{1}%
_{I_{r}}\left( x\right) \right\vert ^{2}d\omega \left( x\right) \\
&\lesssim &\int_{I}\left\vert T_{\sigma }^{\alpha }\mathbf{1}_{I}\right\vert
^{2}d\omega +\left( \mathfrak{T}_{T^{\alpha }}\right) ^{2}\sum_{r=1}^{\infty
}\left\vert I_{r}\right\vert _{\sigma }\lesssim \left( \mathfrak{T}%
_{T^{\alpha }}\right) ^{2}\left\vert I\right\vert _{\sigma }\ .
\end{eqnarray*}%
Altogether this gives%
\begin{equation}
\mathcal{E}_{2}^{\alpha }+\mathcal{E}_{2}^{\alpha ,\ast }\lesssim \mathfrak{T%
}_{T^{\alpha }}+\mathfrak{T}_{T^{\alpha }}^{\ast }+\mathfrak{A}_{2}^{\alpha
}\ .  \label{energy condition is necd}
\end{equation}

\section{Disjoint form\label{Sec disj form}}

Here we control the disjoint form $\Theta _{1}\left( f,g\right) $ by further
decomposing it as follows:%
\begin{eqnarray}
\Theta _{1}\left( f,g\right) &=&\sum_{I\in \mathcal{D}}\sum_{\substack{ J\in 
\mathcal{G}:\ \ell \left( J\right) \leq \ell \left( I\right)  \\ d\left(
J,I\right) >2\ell \left( J\right) ^{\varepsilon }\ell \left( I\right)
^{1-\varepsilon }}}\int \left( T_{\sigma }\square _{I}^{\sigma ,\mathbf{b}%
}f\right) \square _{J}^{\omega ,\mathbf{b}^{\ast }}gd\omega
\label{decomp long short} \\
&=&\sum_{I\in \mathcal{D}}\left\{ \sum_{\substack{ J\in \mathcal{G}:\ \ell
\left( J\right) \leq \ell \left( I\right)  \\ d\left( J,I\right) >\max
\left( \ell \left( I\right) ,2\ell \left( J\right) ^{\varepsilon }\ell
\left( I\right) ^{1-\varepsilon }\right) }}+\sum_{\substack{ J\in \mathcal{G}%
:\ \ell \left( J\right) \leq \ell \left( I\right)  \\ \ell \left( I\right)
\geq d\left( J,I\right) >2\ell \left( J\right) ^{\varepsilon }\ell \left(
I\right) ^{1-\varepsilon }}}\right\} \int \left( T_{\sigma }\square
_{I}^{\sigma ,\mathbf{b}}f\right) \square _{J}^{\omega ,\mathbf{b}^{\ast
}}gd\omega  \notag \\
&\equiv &\Theta _{1}^{\limfunc{long}}\left( f,g\right) +\Theta _{1}^{%
\limfunc{short}}\left( f,g\right) ,  \notag
\end{eqnarray}%
where $\Theta _{1}^{\limfunc{long}}\left( f,g\right) $ is a `long range'
form in which $J$ is far from $I$, and where $\Theta _{1}^{\limfunc{short}%
}\left( f,g\right) $ is a short range form. It should be noted that weak
goodness plays no role in treating the disjoint form.

\subsection{Long range form}

\begin{lemma}
\label{delta long}We have%
\begin{equation*}
\sum_{I\in \mathcal{D}}\sum_{\substack{ J\in \mathcal{G}:\ \ell \left(
J\right) \leq \ell \left( I\right)  \\ d\left( J,I\right) >\ell \left(
I\right) }}\left\vert \int \left( T_{\sigma }\square _{I}^{\sigma ,\mathbf{b}%
}f\right) \square _{J}^{\omega ,\mathbf{b}^{\ast }}gd\omega \right\vert
\lesssim \sqrt{A_{2}^{\alpha }}\left\Vert f\right\Vert _{L^{2}\left( \sigma
\right) }\left\Vert g\right\Vert _{L^{2}\left( \omega \right) }
\end{equation*}
\end{lemma}

\begin{proof}
Since $J$ and $I$ are separated by at least $\max \left\{ \ell \left(
J\right) ,\ell \left( I\right) \right\} $, we have the inequality%
\begin{equation*}
\mathrm{P}^{\alpha }\left( J,\left\vert \square _{I}^{\sigma ,\mathbf{b}%
}f\right\vert \sigma \right) \approx \int_{I}\frac{\ell \left( J\right) }{%
\left\vert y-c_{J}\right\vert ^{2-\alpha }}\left\vert \square _{I}^{\sigma ,%
\mathbf{b}}f\left( y\right) \right\vert d\sigma \left( y\right) \lesssim
\left\Vert \square _{I}^{\sigma ,\mathbf{b}}f\right\Vert _{L^{2}\left(
\sigma \right) }\frac{\ell \left( J\right) \sqrt{\left\vert I\right\vert
_{\sigma }}}{d\left( I,J\right) ^{2-\alpha }},
\end{equation*}%
since $\int_{I}\left\vert \square _{I}^{\sigma ,\mathbf{b}}f\left( y\right)
\right\vert d\sigma \left( y\right) \leq \left\Vert \square _{I}^{\sigma ,%
\mathbf{b}}f\right\Vert _{L^{2}\left( \sigma \right) }\sqrt{\left\vert
I\right\vert _{\sigma }}$. Thus if $A\left( f,g\right) $ denotes the left
hand side of the conclusion of Lemma \ref{delta long}, we have%
\begin{eqnarray*}
A\left( f,g\right) &\lesssim &\sum_{I\in \mathcal{D}}\sum_{J\;:\;\ell \left(
J\right) \leq \ell \left( I\right) :\ d\left( I,J\right) \geq \ell \left(
I\right) }\left\Vert \square _{I}^{\sigma ,\mathbf{b}}f\right\Vert
_{L^{2}\left( \sigma \right) }\left\Vert \square _{J}^{\omega ,\mathbf{b}%
^{\ast }}g\right\Vert _{L^{2}\left( \omega \right) } \\
&&\ \ \ \ \ \ \ \ \ \ \ \ \ \ \ \times \frac{\ell \left( J\right) }{d\left(
I,J\right) ^{2-\alpha }}\sqrt{\left\vert I\right\vert _{\sigma }}\sqrt{%
\left\vert J\right\vert _{\omega }} \\
&\equiv &\sum_{\left( I,J\right) \in \mathcal{P}}\left\Vert \square
_{I}^{\sigma ,\mathbf{b}}f\right\Vert _{L^{2}\left( \sigma \right)
}\left\Vert \square _{J}^{\omega ,\mathbf{b}^{\ast }}g\right\Vert
_{L^{2}\left( \omega \right) }A\left( I,J\right) ; \\
\text{with }A\left( I,J\right) &\equiv &\frac{\ell \left( J\right) }{d\left(
I,J\right) ^{2-\alpha }}\sqrt{\left\vert I\right\vert _{\sigma }}\sqrt{%
\left\vert J\right\vert _{\omega }}; \\
\text{ and }\mathcal{P} &\equiv &\left\{ \left( I,J\right) \in \mathcal{D}%
\times \mathcal{G}:\ell \left( J\right) \leq \ell \left( I\right) \text{ and 
}d\left( I,J\right) \geq \ell \left( I\right) \right\} .
\end{eqnarray*}%
Now let $\mathcal{D}_{N}\equiv \left\{ K\in \mathcal{D}:\ell \left( K\right)
=2^{N}\right\} $ for each $N\in \mathbb{Z}$. For $N\in \mathbb{Z}$ and $s\in 
\mathbb{Z}_{+}$, we further decompose $A\left( f,g\right) $ by pigeonholing
the sidelengths of $I$ and $J$ by $2^{N}$ and $2^{N-s}$ respectively: 
\begin{eqnarray*}
A\left( f,g\right) &=&\sum_{s=0}^{\infty }\sum_{N\in \mathbb{Z}%
}A_{N}^{s}\left( f,g\right) ; \\
A_{N}^{s}\left( f,g\right) &\equiv &\sum_{\left( I,J\right) \in \mathcal{P}%
_{N}^{s}}\left\Vert \square _{I}^{\sigma ,\mathbf{b}}f\right\Vert
_{L^{2}\left( \sigma \right) }\left\Vert \square _{J}^{\omega ,\mathbf{b}%
^{\ast }}g\right\Vert _{L^{2}\left( \omega \right) }A\left( I,J\right) \\
\text{where }\mathcal{P}_{N}^{s} &\equiv &\left\{ \left( I,J\right) \in 
\mathcal{D}_{N}\times \mathcal{G}_{N-s}:d\left( I,J\right) \geq \ell \left(
I\right) \right\} .
\end{eqnarray*}

Now let $\mathsf{P}_{M}^{\sigma }=\dsum\limits_{K\in \mathcal{D}_{M}}\square
_{K}^{\sigma ,\mathbf{b}}$ denote the dual martingale pseudoprojection onto $%
\limfunc{Span}\left\{ \square _{K}^{\sigma ,\mathbf{b}}\right\} _{K\in 
\mathcal{D}_{M}}$. Since the intervals $K$ in $\mathcal{D}_{M}$ are pairwise
disjoint, the pseudoprojections $\square _{K}^{\sigma ,\mathbf{b}}$ are
mutually orthogonal, which means that $\left\Vert \mathsf{P}_{M}^{\sigma
}f\right\Vert _{L^{2}\left( \sigma \right) }^{2}=\dsum\limits_{K\in \mathcal{%
D}_{M}}\left\Vert \square _{K}^{\sigma ,\mathbf{b}}f\right\Vert
_{L^{2}\left( \sigma \right) }^{2}$. We claim that%
\begin{equation}
\left\vert A_{N}^{s}\left( f,g\right) \right\vert \leq C2^{-s}\sqrt{%
A_{2}^{\alpha }}\left\Vert \mathsf{P}_{N}^{\sigma }f\right\Vert
_{L^{2}\left( \sigma \right) }^{\bigstar }\left\Vert \mathsf{P}%
_{N-s}^{\omega }g\right\Vert _{L^{2}\left( \omega \right) }^{\bigstar },\ \
\ \ \ \text{for }s\geq 0\text{ and }N\in \mathbb{Z}.  \label{AsN}
\end{equation}%
With this proved, we can then obtain%
\begin{eqnarray*}
A\left( f,g\right) &=&\sum_{s=0}^{\infty }\sum_{N\in \mathbb{Z}%
}A_{N}^{s}\left( f,g\right) \leq C\sqrt{A_{2}^{\alpha }}\sum_{s=0}^{\infty
}2^{-s}\sum_{N\in \mathbb{Z}}\left\Vert \mathsf{P}_{N}^{\sigma }f\right\Vert
_{L^{2}\left( \sigma \right) }^{\bigstar }\left\Vert \mathsf{P}%
_{N-s}^{\omega }g\right\Vert _{L^{2}\left( \omega \right) }^{\bigstar } \\
&\leq &C\sqrt{A_{2}^{\alpha }}\sum_{s=0}^{\infty }2^{-s}\left( \sum_{N\in 
\mathbb{Z}}\left\Vert \mathsf{P}_{N}^{\sigma }f\right\Vert _{L^{2}\left(
\sigma \right) }^{\bigstar 2}\right) ^{\frac{1}{2}}\left( \sum_{N\in \mathbb{%
Z}}\left\Vert \mathsf{P}_{N-s}^{\omega }g\right\Vert _{L^{2}\left( \omega
\right) }^{\bigstar 2}\right) ^{\frac{1}{2}} \\
&\leq &C\sqrt{A_{2}^{\alpha }}\sum_{s=0}^{\infty }2^{-s}\left\Vert
f\right\Vert _{L^{2}\left( \sigma \right) }\left\Vert g\right\Vert
_{L^{2}\left( \omega \right) }=C\sqrt{A_{2}^{\alpha }}\left\Vert
f\right\Vert _{L^{2}\left( \sigma \right) }\left\Vert g\right\Vert
_{L^{2}\left( \omega \right) },
\end{eqnarray*}%
where in the last line we have used the lower frame inequality for $\square
_{I}^{\sigma ,\mathbf{b}}$ in Appendix A.

To prove (\ref{AsN}), we pigeonhole the distance between $I$ and $J$:%
\begin{eqnarray*}
A_{N}^{s}\left( f,g\right) &=&\dsum\limits_{\ell =0}^{\infty }A_{N,\ell
}^{s}\left( f,g\right) ; \\
A_{N,\ell }^{s}\left( f,g\right) &\equiv &\sum_{\left( I,J\right) \in 
\mathcal{P}_{N,\ell }^{s}}\left\Vert \square _{I}^{\sigma ,\mathbf{b}%
}f\right\Vert _{L^{2}\left( \sigma \right) }\left\Vert \square _{J}^{\omega ,%
\mathbf{b}^{\ast }}g\right\Vert _{L^{2}\left( \omega \right) }A\left(
I,J\right) \\
\text{where }\mathcal{P}_{N,\ell }^{s} &\equiv &\left\{ \left( I,J\right)
\in \mathcal{D}_{N}\times \mathcal{G}_{N-s}:d\left( I,J\right) \approx
2^{N+\ell }\right\} .
\end{eqnarray*}%
If we define $\mathcal{H}\left( A_{N,\ell }^{s}\right) $ to be the bilinear
form on $\ell ^{2}\times \ell ^{2}$ with matrix $\left[ A\left( I,J\right) %
\right] _{\left( I,J\right) \in \mathcal{P}_{N,\ell }^{s}}$, then it remains
to show that the norm $\left\Vert \mathcal{H}\left( A_{N,\ell }^{s}\right)
\right\Vert _{\ell ^{2}\rightarrow \ell ^{2}}$ of $\mathcal{H}\left(
A_{N,\ell }^{s}\right) $ on the sequence space $\ell ^{2}$ is bounded by $%
C2^{-s-\ell }\sqrt{A_{2}^{\alpha }}$. In turn, this is equivalent to showing
that the norm $\left\Vert \mathcal{H}\left( B_{N,\ell }^{s}\right)
\right\Vert _{\ell ^{2}\rightarrow \ell ^{2}}$ of the bilinear form $%
\mathcal{H}\left( B_{N,\ell }^{s}\right) \equiv \mathcal{H}\left( A_{N,\ell
}^{s}\right) ^{\limfunc{tr}}\mathcal{H}\left( A_{N,\ell }^{s}\right) $ on
the sequence space $\ell ^{2}$ is bounded by $C^{2}2^{-2s-2\ell }\mathfrak{A}%
_{2}^{\alpha }$. Here $\mathcal{H}\left( B_{N,\ell }^{s}\right) $ is the
quadratic form with matrix kernel $\left[ B_{N,\ell }^{s}\left( J,J^{\prime
}\right) \right] _{J,J^{\prime }\in \mathcal{D}_{N-s}}$ having entries:%
\begin{equation*}
B_{N,\ell }^{s}\left( J,J^{\prime }\right) \equiv \sum_{I\in \mathcal{D}%
_{N}:\ d\left( I,J\right) \approx d\left( I,J^{\prime }\right) \approx
2^{N+\ell }}A\left( I,J\right) A\left( I,J^{\prime }\right) ,\ \ \ \ \ \text{%
for }J,J^{\prime }\in \mathcal{G}_{N-s}.
\end{equation*}

We are reduced to showing the bilinear form inequality,%
\begin{equation*}
\left\Vert \mathcal{H}\left( B_{N,\ell }^{s}\right) \right\Vert _{\ell
^{2}\rightarrow \ell ^{2}}\leq C2^{-2s-2\ell }A_{2}^{\alpha }\ \ \ \ \text{%
for }s\geq 0\text{, }\ell \geq 0\text{ and }N\in \mathbb{Z}.
\end{equation*}%
We begin by computing $B_{N,\ell }^{s}\left( J,J^{\prime }\right) $:%
\begin{eqnarray*}
B_{N,\ell }^{s}\left( J,J^{\prime }\right) &=&\sum_{\substack{ I\in \mathcal{%
D}_{N}  \\ d\left( I,J\right) \approx d\left( I,J^{\prime }\right) \approx
2^{N+\ell }}}\frac{\ell \left( J\right) }{d\left( I,J\right) ^{2-\alpha }}%
\sqrt{\left\vert I\right\vert _{\sigma }}\sqrt{\left\vert J\right\vert
_{\omega }}\frac{\ell \left( J^{\prime }\right) }{d\left( I,J^{\prime
}\right) ^{2-\alpha }}\sqrt{\left\vert I\right\vert _{\sigma }}\sqrt{%
\left\vert J^{\prime }\right\vert _{\omega }} \\
&=&\left\{ \sum_{\substack{ I\in \mathcal{D}_{N}  \\ d\left( I,J\right)
\approx d\left( I,J^{\prime }\right) \approx 2^{N+\ell }}}\left\vert
I\right\vert _{\sigma }\frac{1}{d\left( I,J\right) ^{2-\alpha }d\left(
I,J^{\prime }\right) ^{2-\alpha }}\right\} \ell \left( J\right) \ell \left(
J^{\prime }\right) \sqrt{\left\vert J\right\vert _{\omega }}\sqrt{\left\vert
J^{\prime }\right\vert _{\omega }}.
\end{eqnarray*}%
Now we show that%
\begin{equation}
\left\Vert \mathcal{H}\left( B_{N,\ell }^{s}\right) \right\Vert _{\ell
^{2}\rightarrow \ell ^{2}}=\left\Vert B_{N,\ell }^{s}\right\Vert _{\ell
^{2}\rightarrow \ell ^{2}}\lesssim 2^{-2s-2\ell }A_{2}^{\alpha }\ ,
\label{Schur s}
\end{equation}%
by applying the proof of Schur's lemma. Fix $\ell \geq 0$ and $s\geq 0$.
Choose the Schur function $\beta \left( K\right) =\frac{1}{\sqrt{\left\vert
K\right\vert _{\omega }}}$. Fix $J\in \mathcal{D}_{N-s}$. We now group those 
$I\in \mathcal{D}_{N}$ with $d\left( I,J\right) \approx 2^{N+\ell }$ into
finitely many groups $G_{1},...G_{C}$ for which the union of the $I$ in each
group is contained in an interval of side length roughly $\frac{1}{100}%
2^{N+\ell }$ , and we set $I_{k}^{\ast }\equiv \dbigcup\limits_{I\in G_{k}}I$
for $1\leq k\leq C$. We then have%
\begin{eqnarray*}
&&\sum_{J^{\prime }\in \mathcal{G}_{N-s}}\frac{\beta \left( J\right) }{\beta
\left( J^{\prime }\right) }B_{N,\ell }^{s}\left( J,J^{\prime }\right) \\
&=&\sum_{\substack{ J^{\prime }\in \mathcal{G}_{N-s}  \\ d\left( J^{\prime
},J\right) \leq \frac{1}{100}2^{N+\ell +2}}}\frac{\beta \left( J\right) }{%
\beta \left( J^{\prime }\right) }B_{N,\ell }^{s}\left( J,J^{\prime }\right)
+\sum_{\substack{ J^{\prime }\in \mathcal{G}_{N-s}  \\ d\left( J^{\prime
},J\right) >\frac{1}{100}2^{N+\ell +2}}}\frac{\beta \left( J\right) }{\beta
\left( J^{\prime }\right) }B_{N,\ell }^{s}\left( J,J^{\prime }\right) \\
&=&A+B,
\end{eqnarray*}%
where%
\begin{eqnarray*}
&&A\lesssim \sum_{\substack{ J^{\prime }\in \mathcal{G}_{N-s}  \\ d\left(
J,J^{\prime }\right) \leq \frac{1}{100}2^{N+\ell +2}}}\left\{ \sum 
_{\substack{ I\in \mathcal{D}_{N}  \\ d\left( I,J\right) \approx 2^{N+\ell } 
}}\left\vert I\right\vert _{\sigma }\right\} \ \frac{2^{2\left( N-s\right) }%
}{2^{2\left( \ell +N\right) \left( 2-\alpha \right) }}\left\vert J^{\prime
}\right\vert _{\omega } \\
&=&\sum_{\substack{ J^{\prime }\in \mathcal{G}_{N-s}  \\ d\left( J,J^{\prime
}\right) \leq \frac{1}{100}2^{N+\ell +2}}}\left\{ \sum_{k=1}^{C}\left\vert
I_{k}^{\ast }\right\vert _{\sigma }\right\} \ \frac{2^{2\left( N-s\right) }}{%
2^{2\left( \ell +N\right) \left( 2-\alpha \right) }}\left\vert J^{\prime
}\right\vert _{\omega }=\frac{2^{2\left( N-s\right) }}{2^{2\left( \ell
+N\right) \left( 2-\alpha \right) }}\sum_{k=1}^{C}\sum_{\substack{ J^{\prime
}\in \mathcal{G}_{N-s}  \\ d\left( J,J^{\prime }\right) \leq \frac{1}{100}%
2^{N+\ell +2}}}\left\vert I_{k}^{\ast }\right\vert _{\sigma }\ \left\vert
J^{\prime }\right\vert _{\omega } \\
&\lesssim &2^{-2s-2\ell }\sum_{k=1}^{C}\frac{\left\vert I_{k}^{\ast
}\right\vert _{\sigma }}{2^{\left( \ell +N\right) \left( 1-\alpha \right) }}%
\frac{\left\vert \frac{1}{100}2^{s+\ell +4}J\right\vert _{\omega }}{%
2^{\left( \ell +N\right) \left( 1-\alpha \right) }}\lesssim 2^{-2s-2\ell
}A_{2}^{\alpha },
\end{eqnarray*}%
since the intervals $I_{k}^{\ast }$ and $\frac{1}{100}2^{s+\ell +4}J$ are
well separated.

Define%
\begin{equation*}
E_{k}^{\limfunc{left}}\equiv \dbigcup\limits_{\substack{ J^{\prime }\in 
\mathcal{G}_{N-s}:\ d\left( I_{k}^{\ast },J^{\prime }\right) \approx
2^{N+\ell }  \\ d\left( J,J^{\prime }\right) >\frac{1}{100}2^{N+\ell +2}%
\text{ and }J^{\prime }\text{ is to the left of }I_{k}^{\ast }}}J^{\prime }\ 
\text{\ and }E_{k}^{\limfunc{right}}\equiv \dbigcup\limits_{\substack{ %
J^{\prime }\in \mathcal{G}_{N-s}:\ d\left( I_{k}^{\ast },J^{\prime }\right)
\approx 2^{N+\ell }  \\ d\left( J,J^{\prime }\right) >\frac{1}{100}2^{N+\ell
+2}\text{ and }J^{\prime }\text{ is\ to the right of }I_{k}^{\ast }}}%
J^{\prime },
\end{equation*}%
and let $Q_{k}^{\limfunc{left}}$ (respectively $Q_{k}^{\limfunc{right}}$) be
the smallest interval containing $E_{k}^{\limfunc{left}}$ (respectively $%
E_{k}^{\limfunc{right}}$). Then we have%
\begin{eqnarray*}
B &\lesssim &\sum_{\substack{ J^{\prime }\in \mathcal{G}_{N-s}  \\ d\left(
J,J^{\prime }\right) >\frac{1}{100}2^{N+\ell +2}}}\left\{ \sum_{\substack{ %
I\in \mathcal{D}_{N}  \\ d\left( I,J^{\prime }\right) \approx d\left(
I,J\right) \approx 2^{N+\ell }}}\left\vert I\right\vert _{\sigma }\right\} \ 
\frac{2^{2\left( N-s\right) }}{2^{2\left( \ell +N\right) \left( 2-\alpha
\right) }}\left\vert J^{\prime }\right\vert _{\omega } \\
&\lesssim &\sum_{\substack{ J^{\prime }\in \mathcal{G}_{N-s}  \\ d\left(
J,J^{\prime }\right) >\frac{1}{100}2^{N+\ell +2}}}\left\{ \sum_{k:\ d\left(
I_{k}^{\ast },J^{\prime }\right) \approx 2^{N+\ell }}\left\vert I_{k}^{\ast
}\right\vert _{\sigma }\right\} \ \frac{2^{2\left( N-s\right) }}{2^{2\left(
\ell +N\right) \left( 2-\alpha \right) }}\left\vert J^{\prime }\right\vert
_{\omega } \\
&\lesssim &\frac{2^{2\left( N-s\right) }}{2^{2\left( \ell +N\right) \left(
2-\alpha \right) }}\sum_{k=1}^{C}\left\vert I_{k}^{\ast }\right\vert
_{\sigma }\left\vert E_{k}^{\limfunc{left}}\cup E_{k}^{\limfunc{right}%
}\right\vert _{\omega } \\
&\lesssim &2^{-2s-2\ell }\sum_{k=1}^{C}\frac{\left\vert I_{k}^{\ast
}\right\vert _{\sigma }}{2^{\left( \ell +N\right) \left( 1-\alpha \right) }}%
\frac{\left\vert Q_{k}^{\limfunc{left}}\right\vert _{\omega }+\left\vert
Q_{k}^{\limfunc{right}}\right\vert _{\omega }}{2^{\left( \ell +N\right)
\left( 1-\alpha \right) }}\lesssim 2^{-2s-2\ell }\mathfrak{A}_{2}^{\alpha },
\end{eqnarray*}%
since the interval $I_{k}^{\ast }$ is well separated from each of the
intervals $Q_{k}^{\limfunc{left}}$ and $Q_{k}^{\limfunc{right}}$.

Thus we can now apply Schur's argument with $\sum_{J}\left( a_{J}\right)
^{2}=\sum_{J^{\prime }}\left( b_{J^{\prime }}\right) ^{2}=1$ to obtain%
\begin{eqnarray*}
&&\sum_{J,J^{\prime }\in \mathcal{G}_{N-s}}a_{J}b_{J^{\prime }}B_{N,\ell
}^{s}\left( J,J^{\prime }\right) \\
&=&\sum_{J,J^{\prime }\in \mathcal{G}_{N-s}}a_{J}\beta \left( J\right)
b_{J^{\prime }}\beta \left( J^{\prime }\right) \frac{B_{N,\ell }^{s}\left(
J,J^{\prime }\right) }{\beta \left( J\right) \beta \left( J^{\prime }\right) 
} \\
&\leq &\sum_{J}\left( a_{J}\beta \left( J\right) \right) ^{2}\sum_{J^{\prime
}}\frac{B_{N,\ell }^{s}\left( J,J^{\prime }\right) }{\beta \left( J\right)
\beta \left( J^{\prime }\right) }+\sum_{J^{\prime }}\left( b_{J^{\prime
}}\beta \left( J^{\prime }\right) \right) ^{2}\sum_{J}\frac{B_{N,\ell
}^{s}\left( J,J^{\prime }\right) }{\beta \left( J\right) \beta \left(
J^{\prime }\right) } \\
&=&\sum_{J}\left( a_{J}\right) ^{2}\left\{ \sum_{J^{\prime }}\frac{\beta
\left( J\right) }{\beta \left( J^{\prime }\right) }B_{N,\ell }^{s}\left(
J,J^{\prime }\right) \right\} +\sum_{J^{\prime }}\left( b_{J^{\prime
}}\right) ^{2}\left\{ \sum_{J}\frac{\beta \left( J^{\prime }\right) }{\beta
\left( J\right) }B_{N,\ell }^{s}\left( J,J^{\prime }\right) \right\} \\
&\lesssim &2^{-2s-2\ell }A_{2}^{\alpha }\left( \sum_{J}\left( a_{J}\right)
^{2}+\sum_{J^{\prime }}\left( b_{J^{\prime }}\right) ^{2}\right)
=2^{1-2s-2\ell }A_{2}^{\alpha }.
\end{eqnarray*}%
This completes the proof of (\ref{Schur s}). We can now sum in $\ell $ to
get (\ref{AsN}) and we are done. This completes our proof of the long range
estimate%
\begin{equation*}
\mathcal{A}\left( f,g\right) \lesssim \sqrt{A_{2}^{\alpha }}\left\Vert
f\right\Vert _{L^{2}\left( \sigma \right) }\left\Vert g\right\Vert
_{L^{2}\left( \omega \right) }\ .
\end{equation*}
\end{proof}

\subsection{Short range form}

The form $\Theta _{1}^{\limfunc{short}}\left( f,g\right) $ is handled by the
following lemma.

\begin{lemma}
\label{delta short}We have%
\begin{equation*}
\sum_{I\in \mathcal{D}}\sum_{\substack{ J\in \mathcal{G}:\ \ell \left(
J\right) \leq \ell \left( I\right)  \\ \ell \left( I\right) \geq d\left(
J,I\right) >2\ell \left( J\right) ^{\varepsilon }\ell \left( I\right)
^{1-\varepsilon }}}\left\vert \int \left( T_{\sigma }\square _{I}^{\sigma ,%
\mathbf{b}}f\right) \square _{J}^{\omega ,\mathbf{b}^{\ast }}gd\omega
\right\vert \lesssim \sqrt{A_{2}^{\alpha }}\left\Vert f\right\Vert
_{L^{2}\left( \sigma \right) }\left\Vert g\right\Vert _{L^{2}\left( \omega
\right) }
\end{equation*}
\end{lemma}

\begin{proof}
The pairs $\left( I,J\right) $ that occur in the sum above satisfy $J\subset
4I\setminus I$ and so we consider 
\begin{equation*}
\mathcal{P}\equiv \left\{ \left( I,J\right) \in \mathcal{D}\times \mathcal{G}%
:\ell \left( J\right) \leq \ell \left( I\right) ,\ \ell \left( I\right) \geq
d\left( J,I\right) >2\ell \left( J\right) ^{\varepsilon }\ell \left(
I\right) ^{1-\varepsilon },\text{ }J\subset 4I\setminus I\right\} .
\end{equation*}%
For $\left( I,J\right) \in \mathcal{P}$, the `pivotal' estimate from the
Energy Lemma \ref{ener} gives%
\begin{equation*}
\left\vert \left\langle T_{\sigma }^{\alpha }\left( \square _{I}^{\sigma ,%
\mathbf{b}}f\right) ,\square _{J}^{\omega ,\mathbf{b}^{\ast }}g\right\rangle
_{\omega }\right\vert \lesssim \left\Vert \square _{J}^{\omega ,\mathbf{b}%
^{\ast }}g\right\Vert _{L^{2}\left( \omega \right) }\mathrm{P}^{\alpha
}\left( J,\left\vert \square _{I}^{\sigma ,\mathbf{b}}f\right\vert \sigma
\right) \sqrt{\left\vert J\right\vert _{\omega }}\,.
\end{equation*}%
Now we pigeonhole the lengths of $I$ and $J$ and the distance between them
by defining%
\begin{equation*}
\mathcal{P}_{N,d}^{s}\equiv \left\{ \left( I,J\right) \in \mathcal{P}:\ell
\left( I\right) =2^{N},\ \ell \left( J\right) =2^{N-s},\ 2^{d-1}\leq d\left(
I,J\right) \leq 2^{d},\text{ }J\subset 4I\setminus I\right\} .
\end{equation*}%
Note that the closest an interval $J$ can come to $I$ is determined by: 
\begin{eqnarray*}
&&2^{d}\geq 2\ell \left( I\right) ^{1-\varepsilon }\ell \left( J\right)
^{\varepsilon }=2^{1+N\left( 1-\varepsilon \right) +\left( N-s\right)
\varepsilon }=2^{1+N-\varepsilon s}; \\
&&\text{which implies }N-\varepsilon s+1\leq d\leq N.
\end{eqnarray*}%
Thus we have%
\begin{eqnarray*}
&&\dsum\limits_{\left( I,J\right) \in \mathcal{P}}\left\vert \left\langle
T_{\sigma }^{\alpha }\left( \square _{I}^{\sigma ,\mathbf{b}}f\right)
,\square _{J}^{\omega ,\mathbf{b}^{\ast }}g\right\rangle _{\omega
}\right\vert \lesssim \dsum\limits_{\left( I,J\right) \in \mathcal{P}%
}\left\Vert \square _{J}^{\omega ,\mathbf{b}^{\ast }}g\right\Vert
_{L^{2}\left( \omega \right) }\mathrm{P}^{\alpha }\left( J,\left\vert
\square _{I}^{\sigma ,\mathbf{b}}f\right\vert \sigma \right) \sqrt{%
\left\vert J\right\vert _{\omega }} \\
&&\ \ \ \ \ =\dsum\limits_{s=0}^{\infty }\ \sum_{N\in \mathbb{Z}}\
\sum_{d=N-\varepsilon s+1}^{N}\ \sum_{\left( I,J\right) \in \mathcal{P}%
_{N,d}^{s}}\ \left\Vert \square _{J}^{\omega ,\mathbf{b}^{\ast
}}g\right\Vert _{L^{2}\left( \omega \right) }\mathrm{P}^{\alpha }\left(
J,\left\vert \square _{I}^{\sigma ,\mathbf{b}}f\right\vert \sigma \right) 
\sqrt{\left\vert J\right\vert _{\omega }}.
\end{eqnarray*}%
Now we use%
\begin{eqnarray*}
\mathrm{P}^{\alpha }\left( J,\left\vert \square _{I}^{\sigma ,\mathbf{b}%
}f\right\vert \sigma \right) &=&\int_{I}\frac{\ell \left( J\right) }{\left(
\ell \left( J\right) +\left\vert y-c_{J}\right\vert \right) ^{2-\alpha }}%
\left\vert \square _{I}^{\sigma ,\mathbf{b}}f\left( y\right) \right\vert
d\sigma \left( y\right) \\
&\lesssim &\frac{2^{N-s}}{2^{d\left( 2-\alpha \right) }}\left\Vert \square
_{I}^{\sigma ,\mathbf{b}}f\right\Vert _{L^{2}\left( \sigma \right) }\sqrt{%
\left\vert I\right\vert _{\sigma }}
\end{eqnarray*}%
and apply Cauchy-Schwarz in $J$ and use $J\subset 4I\setminus I$ to get%
\begin{eqnarray*}
&&\dsum\limits_{\left( I,J\right) \in \mathcal{P}}\left\vert \left\langle
T_{\sigma }^{\alpha }\left( \square _{I}^{\sigma ,\mathbf{b}}f\right)
,\square _{J}^{\omega ,\mathbf{b}^{\ast }}g\right\rangle _{\omega
}\right\vert \\
&\lesssim &\dsum\limits_{s=0}^{\infty }\ \sum_{N\in \mathbb{Z}}\
\sum_{d=N-\varepsilon s+1}^{N}\ \sum_{I\in \mathcal{D}_{N}}\frac{%
2^{N-s}2^{N\left( 1-\alpha \right) }}{2^{d\left( 2-\alpha \right) }}%
\left\Vert \square _{I}^{\sigma ,\mathbf{b}}f\right\Vert _{L^{2}\left(
\sigma \right) }\frac{\sqrt{\left\vert I\right\vert _{\sigma }}\sqrt{%
\left\vert 4I\setminus I\right\vert _{\omega }}}{2^{N\left( 1-\alpha \right)
}} \\
&&\ \ \ \ \ \ \ \ \ \ \ \ \ \ \ \ \ \ \ \ \ \ \ \ \ \ \ \ \ \ \times \sqrt{%
\sum_{\substack{ J\in \mathcal{G}_{N-s}  \\ J\subset 4I\setminus I\text{ and 
}d\left( I,J\right) \approx 2^{d}}}\left\Vert \square _{J}^{\omega ,\mathbf{b%
}^{\ast }}g\right\Vert _{L^{2}\left( \omega \right) }^{2}} \\
&\lesssim &\dsum\limits_{s=0}^{\infty }\ \sum_{N\in \mathbb{Z}}\frac{%
2^{N-s}2^{N\left( 1-\alpha \right) }}{2^{\left( N-\varepsilon s\right)
\left( 2-\alpha \right) }}\sqrt{A_{2}^{\alpha }}\sum_{I\in \mathcal{D}%
_{N}}\left\Vert \square _{I}^{\sigma ,\mathbf{b}}f\right\Vert _{L^{2}\left(
\sigma \right) }\sqrt{\sum_{\substack{ J\in \mathcal{G}_{N-s}  \\ J\subset
4I\setminus I}}\left\Vert \square _{J}^{\omega ,\mathbf{b}^{\ast
}}g\right\Vert _{L^{2}\left( \omega \right) }^{2}} \\
&\lesssim &\dsum\limits_{s=0}^{\infty }2^{-s\left[ 1-\varepsilon \left(
2-\alpha \right) \right] }\sqrt{A_{2}^{\alpha }}\left\Vert f\right\Vert
_{L^{2}\left( \sigma \right) }\left\Vert g\right\Vert _{L^{2}\left( \omega
\right) }\lesssim \sqrt{A_{2}^{\alpha }}\left\Vert f\right\Vert
_{L^{2}\left( \sigma \right) }\left\Vert g\right\Vert _{L^{2}\left( \omega
\right) },
\end{eqnarray*}%
where in the third line above we have summed the geometric series in $d$,
and in the last line $\frac{2^{N-s}2^{N\left( 1-\alpha \right) }}{2^{\left(
N-\varepsilon s\right) \left( 2-\alpha \right) }}=2^{-s\left[ 1-\varepsilon
\left( 2-\alpha \right) \right] }$ followed by Cauchy-Schwarz in $I$ and $N$%
, using that we have bounded overlap in the quadruples of $I$ for $I\in 
\mathcal{D}_{N}$, and finally using the lower frame inequality for $\square
_{I}^{\sigma ,\mathbf{b}}$ from Appendix A. More precisely, if we define $%
f_{k}\equiv \Psi _{\mathcal{D}_{k}}^{\sigma ,\mathbf{b}}f=\sum_{I\in 
\mathcal{D}_{k}}\square _{I}^{\sigma ,\mathbf{b}}f$ and $g_{k}\equiv \Psi _{%
\mathcal{G}_{k}}^{\sigma ,\mathbf{b}^{\ast }}g=\sum_{J\in \mathcal{G}%
_{k}}\square _{J}^{\omega ,\mathbf{b}^{\ast }}g$, then we have the
quasi-orthogonality inequality 
\begin{eqnarray*}
\sum_{N\in \mathbb{Z}}\left\Vert f_{N}\right\Vert _{L^{2}\left( \sigma
\right) }\left\Vert g_{N-s}\right\Vert _{L^{2}\left( \omega \right) } &\leq
&\left( \sum_{N\in \mathbb{Z}}\left\Vert f_{N}\right\Vert _{L^{2}\left(
\sigma \right) }^{2}\right) ^{\frac{1}{2}}\left( \sum_{N\in \mathbb{Z}%
}\left\Vert g_{N-s}\right\Vert _{L^{2}\left( \omega \right) }^{2}\right) ^{%
\frac{1}{2}} \\
&\lesssim &\left\Vert f\right\Vert _{L^{2}\left( \sigma \right) }\left\Vert
g\right\Vert _{L^{2}\left( \omega \right) }.
\end{eqnarray*}%
We have assumed that 
\begin{equation}
0<\varepsilon <\frac{1}{2-\alpha }  \label{short requirement}
\end{equation}%
in the calculations above in order that the sum in $s$ converges, and this
completes the proof of Lemma \ref{delta short}.
\end{proof}

\section{Nearby form\label{Sec nearby}}

We dominate the nearby form $\Theta _{3}\left( f,g\right) $ by%
\begin{equation*}
\left\vert \Theta _{3}\left( f,g\right) \right\vert \leq \sum_{I\in \mathcal{%
D}}\sum_{\substack{ J\in \mathcal{G}:\ 2^{-\mathbf{r}}\ell \left( I\right)
<\ell \left( J\right) \leq \ell \left( I\right)  \\ d\left( J,I\right) \leq
2\ell \left( J\right) ^{\varepsilon }\ell \left( I\right) ^{1-\varepsilon }}}%
\left\vert \int \left( T_{\sigma }^{\alpha }\square _{I}^{\sigma ,\mathbf{b}%
}f\right) \square _{J}^{\omega ,\mathbf{b}^{\ast }}gd\omega \right\vert ,
\end{equation*}%
and prove the following lemma that controls the expectation, over two
independent grids, of the nearby form $\Theta _{3}\left( f,g\right) $. It
should be noted that weak goodness plays no role in treating the nearby form.

\begin{lemma}
\label{nearby form}Suppose $T^{\alpha }$ is a standard fractional singular
integral with $0\leq \alpha <1$. Let $\theta \in \left( 0,1\right) $ be
sufficiently small depending only on $0\leq \alpha <1$. Then there is a
constant $C_{\theta }$ such that for $f\in L^{2}\left( \sigma \right) $ and $%
g\in L^{2}\left( \omega \right) $, and dual martingale differences $\square
_{I}^{\sigma ,\mathbf{b}}$ and $\square _{J}^{\omega ,\mathbf{b}^{\ast }}$
with $\infty $-strongly accretive families of test functions $\mathbf{b}$
and $\mathbf{b}^{\ast }$, we have%
\begin{eqnarray}
&&\boldsymbol{E}_{\Omega }^{\mathcal{D}}\boldsymbol{E}_{\Omega }^{\mathcal{G}%
}\sum_{I\in \mathcal{D}}\sum_{\substack{ J\in \mathcal{G}:\ 2^{-\mathbf{r}%
}\ell \left( I\right) <\ell \left( J\right) \leq \ell \left( I\right)  \\ %
d\left( J,I\right) \leq 2\ell \left( J\right) ^{\varepsilon }\ell \left(
I\right) ^{1-\varepsilon }}}\left\vert \left\langle T_{\sigma }^{\alpha
}\left( \square _{I}^{\sigma ,\mathbf{b}}f\right) ,\square _{J}^{\omega ,%
\mathbf{b}^{\ast }}g\right\rangle _{\omega }\right\vert  \label{delta near}
\\
&\lesssim &\left( C_{\theta }\mathcal{NTV}_{\alpha }+\sqrt{\theta }\mathfrak{%
N}_{T^{\alpha }}\right) \left\Vert f\right\Vert _{L^{2}\left( \sigma \right)
}\left\Vert g\right\Vert _{L^{2}\left( \omega \right) }\ .  \notag
\end{eqnarray}
\end{lemma}

\textbf{Proof:} As usual, we continue to write the independent grids for $f$
and $g$ as $\mathcal{D}$ and $\mathcal{G}$ respectively. Write the dual
martingale averages $\square _{I}^{\sigma ,\mathbf{b}}f$ and $\square
_{J}^{\omega ,\mathbf{b}^{\ast }}g$ as linear combinations 
\begin{eqnarray*}
\square _{I}^{\sigma ,\mathbf{b}}f &=&b_{I}\ \sum_{I^{\prime }\in \mathfrak{C%
}_{\limfunc{natural}}\left( I\right) }\mathbf{1}_{I^{\prime }}\ E_{I^{\prime
}}^{\sigma }\left( \widehat{\square }_{I}^{\sigma ,\mathbf{b}}f\right)
+\sum_{I^{\prime }\in \mathfrak{C}_{\limfunc{broken}}\left( I\right)
}b_{I^{\prime }}\ \mathbf{1}_{I^{\prime }}\widehat{\mathbb{F}}_{I^{\prime
}}^{\sigma ,b_{I^{\prime }}}f-b_{I}\ \sum_{I^{\prime }\in \mathfrak{C}_{%
\limfunc{broken}}\left( I\right) }\mathbf{1}_{I^{\prime }}\widehat{\mathbb{F}%
}_{I}^{\sigma ,b_{I}}f, \\
\square _{J}^{\omega ,\mathbf{b}^{\ast }}g &=&b_{J}^{\ast }\ \sum_{J^{\prime
}\in \mathfrak{C}_{\limfunc{natural}}\left( J\right) }\mathbf{1}_{J^{\prime
}}\ E_{J^{\prime }}^{\omega }\left( \widehat{\square }_{J}^{\omega ,\mathbf{b%
}^{\ast }}g\right) +\sum_{J^{\prime }\in \mathfrak{C}_{\limfunc{broken}%
}\left( J\right) }b_{J^{\prime }}^{\ast }\ \mathbf{1}_{J^{\prime }}\widehat{%
\mathbb{F}}_{J^{\prime }}^{\omega ,b_{J^{\prime }}^{\ast }}g-b_{J}^{\ast }\
\sum_{J^{\prime }\in \mathfrak{C}_{\limfunc{broken}}\left( J\right) }\mathbf{%
1}_{J^{\prime }}\widehat{\mathbb{F}}_{J}^{\omega ,b_{J}^{\ast }}g,
\end{eqnarray*}%
of the appropriate function $b$ times the indicators of their children,
denoted $I^{\prime }$ and $J^{\prime }$ respectively. We will regroup the
terms as needed below.

\begin{notation}
On the natural child $I^{\prime }$, the expression $\widehat{\square }%
_{I}^{\sigma ,\mathbf{b}}f=\frac{1}{b_{I}}\square _{I}^{\sigma ,\mathbf{b}}f$
simply denotes the dual martingale average with $b_{I}$ removed, so that we
need not assume $\left\vert b_{I}\right\vert $ is bounded below in order to
make sense of $\frac{1}{b_{I}}\square _{I}^{\sigma ,\mathbf{b}}f$. Similar
comments apply to the expressions $\widehat{\mathbb{F}}_{I^{\prime
}}^{\sigma ,b_{I^{\prime }}}f=\frac{1}{b_{I^{\prime }}}\mathbb{F}_{I^{\prime
}}^{\sigma ,b_{I^{\prime }}}f$ and $\widehat{\mathbb{F}}_{I}^{\sigma
,b_{I}}f=\frac{1}{b_{I}}\mathbb{F}_{I}^{\sigma ,b_{I}}f$. On the other hand,
we are assuming from Proposition \ref{lower bound} that the $PLBP$ (\ref{plb}%
) holds, which shows that $\frac{1}{b_{I^{\prime }}}$ is actually a bounded
function. This latter fact will be used shortly in (\ref{box hat bound})
below.
\end{notation}

Recall that the length of $J$ is at most the length of $I$, i.e. $\ell
\left( J\right) \leq \ell \left( I\right) $. If $J$ and $I$ are \emph{%
separated}, by which we mean here that $J\cap I=\emptyset $, then by (\ref%
{disj supp}) we have the satisfactory estimate%
\begin{equation*}
\left\vert \left\langle T_{\sigma }^{\alpha }\left( \square _{I}^{\sigma ,%
\mathbf{b}}f\right) ,\square _{J}^{\omega ,\mathbf{b}^{\ast }}g\right\rangle
_{\omega }\right\vert \lesssim \sqrt{\mathcal{A}_{2}^{\alpha ,\ast }+%
\mathcal{A}_{2}^{\alpha }}\left\Vert \square _{I}^{\sigma ,\mathbf{b}%
}f\right\Vert _{L^{2}\left( \sigma \right) }\left\Vert \square _{J}^{\omega ,%
\mathbf{b}^{\ast }}g\right\Vert _{L^{2}\left( \omega \right) }.
\end{equation*}%
Suppose now that $J\cap I\neq \emptyset $. Using (\ref{flat broken}) we have%
\begin{eqnarray}
\left\langle T_{\sigma }^{\alpha }\left( \square _{I}^{\sigma ,\mathbf{b}%
}f\right) ,\square _{J}^{\omega ,\mathbf{b}^{\ast }}g\right\rangle _{\omega
} &=&\left\langle T_{\sigma }^{\alpha }\left( \square _{I}^{\sigma ,\flat ,%
\mathbf{b}}f\right) ,\square _{J}^{\omega ,\flat ,\mathbf{b}^{\ast
}}g\right\rangle _{\omega }+\left\langle T_{\sigma }^{\alpha }\left( \square
_{I,\limfunc{broken}}^{\sigma ,\flat ,\mathbf{b}}f\right) ,\square _{J,%
\limfunc{broken}}^{\omega ,\flat ,\mathbf{b}^{\ast }}g\right\rangle _{\omega
}  \label{case d small} \\
&&+\left\langle T_{\sigma }^{\alpha }\left( \square _{I}^{\sigma ,\flat ,%
\mathbf{b}}f\right) ,\square _{J,\limfunc{broken}}^{\omega ,\flat ,\mathbf{b}%
^{\ast }}g\right\rangle _{\omega }+\left\langle T_{\sigma }^{\alpha }\left(
\square _{I,\limfunc{broken}}^{\sigma ,\flat ,\mathbf{b}}f\right) ,\square
_{J}^{\omega ,\flat ,\mathbf{b}^{\ast }}g\right\rangle _{\omega }\ .  \notag
\end{eqnarray}%
The estimation of the latter three inner products, i.e. those in which a
broken operator $\square _{I,\limfunc{broken}}^{\sigma ,\flat ,\mathbf{b}}$
or $\square _{J,\limfunc{broken}}^{\omega ,\flat ,\mathbf{b}^{\ast }}$
arises, is easy. Indeed, recall that 
\begin{eqnarray*}
\square _{I,\limfunc{broken}}^{\sigma ,\flat ,\mathbf{b}}f
&=&\sum_{I^{\prime }\in \mathfrak{C}_{\limfunc{broken}}\left( I\right) }%
\mathbb{F}_{I^{\prime }}^{\sigma ,\mathbf{b}}f=\sum_{I^{\prime }\in 
\mathfrak{C}_{\limfunc{broken}}\left( I\right) }\left( E_{I^{\prime
}}^{\sigma }\widehat{\mathbb{F}}_{I^{\prime }}^{\sigma ,\mathbf{b}}f\right)
b_{I^{\prime }}\ , \\
\square _{J,\limfunc{broken}}^{\omega ,\flat ,\mathbf{b}}g
&=&\sum_{J^{\prime }\in \mathfrak{C}_{\limfunc{broken}}\left( J\right) }%
\mathbb{F}_{J^{\prime }}^{\omega ,\mathbf{b}^{\ast }}g=\sum_{J^{\prime }\in 
\mathfrak{C}_{\limfunc{broken}}\left( J\right) }\left( E_{J^{\prime
}}^{\omega }\widehat{\mathbb{F}}_{J^{\prime }}^{\omega ,\mathbf{b}^{\ast
}}g\right) b_{J^{\prime }}^{\ast }\ ,
\end{eqnarray*}%
so that if at least one broken difference appears in the inner product, as
is the case for the latter three inner products in (\ref{case d small}),
then testing and Cauchy-Schwarz are all that is needed. For example, the
fourth term satisfies 
\begin{eqnarray*}
\left\vert \left\langle T_{\sigma }^{\alpha }\left( \square _{I,\limfunc{%
broken}}^{\sigma ,\flat ,\mathbf{b}}f\right) ,\square _{J}^{\omega ,\flat ,%
\mathbf{b}^{\ast }}g\right\rangle _{\omega }\right\vert &=&\left\vert
\sum_{I^{\prime }\in \mathfrak{C}_{\limfunc{broken}}\left( I\right) }\left(
E_{I^{\prime }}^{\sigma }\widehat{\mathbb{F}}_{I^{\prime }}^{\sigma ,\mathbf{%
b}}f\right) \left\langle T_{\sigma }^{\alpha }b_{I^{\prime }},\square
_{J}^{\omega ,\flat ,\mathbf{b}^{\ast }}g\right\rangle _{\omega }\right\vert
\\
&\lesssim &\sum_{I^{\prime }\in \mathfrak{C}_{\limfunc{broken}}\left(
I\right) }\left\vert E_{I^{\prime }}^{\sigma }\widehat{\mathbb{F}}%
_{I^{\prime }}^{\sigma ,\mathbf{b}}f\right\vert \mathfrak{T}_{T^{\alpha }}^{%
\mathbf{b}}\sqrt{\left\vert I^{\prime }\right\vert _{\sigma }}\left\Vert
\square _{J}^{\omega ,\flat ,\mathbf{b}^{\ast }}g\right\Vert _{L^{2}\left(
\omega \right) } \\
&\lesssim &\mathfrak{T}_{T^{\alpha }}^{\mathbf{b}}\left\Vert \nabla
_{I}^{\sigma }f\right\Vert _{L^{2}\left( \sigma \right) }\sum_{I^{\prime
}\in \mathfrak{C}_{\limfunc{broken}}\left( I\right) }\left( \left\Vert
\square _{J}^{\omega ,\mathbf{b}^{\ast }}g\right\Vert _{L^{2}\left( \omega
\right) }+\left\Vert \square _{J,\limfunc{broken}}^{\omega ,\flat ,\mathbf{b}%
^{\ast }}g\right\Vert _{L^{2}\left( \omega \right) }\right) \\
&\lesssim &\mathcal{NTV}_{\alpha }\left\Vert \square _{I}^{\sigma ,\mathbf{b}%
}f\right\Vert _{L^{2}\left( \sigma \right) }^{\bigstar }\left\Vert \square
_{J}^{\omega ,\mathbf{b}^{\ast }}g\right\Vert _{L^{2}\left( \omega \right)
}^{\bigstar }\ ,
\end{eqnarray*}%
and the third term can be written as $\left\langle \square _{I}^{\sigma
,\flat ,\mathbf{b}}f,T_{\omega }^{\alpha ,\ast }\left( \square _{J,\limfunc{%
broken}}^{\omega ,\flat ,\mathbf{b}^{\ast }}g\right) \right\rangle _{\sigma
} $ and handled similarly.

Thus it remains to consider the first inner product $\left\langle T_{\sigma
}^{\alpha }\left( \square _{I}^{\sigma ,\flat ,\mathbf{b}}f\right) ,\square
_{J}^{\omega ,\flat ,\mathbf{b}^{\ast }}g\right\rangle _{\omega }$ on the
right hand side of (\ref{case d small}), which we call the problematic term,
and write it as%
\begin{eqnarray}
P\left( I,J\right) &\equiv &\left\langle T_{\sigma }^{\alpha }\left( \square
_{I}^{\sigma ,\flat ,\mathbf{b}}f\right) ,\square _{J}^{\omega ,\flat ,%
\mathbf{b}^{\ast }}g\right\rangle _{\omega }  \label{def P(I,J)} \\
&=&\sum_{I^{\prime }\in \mathfrak{C}\left( I\right) \text{ and }J^{\prime
}\in \mathfrak{C}\left( J\right) }\left\langle T_{\sigma }^{\alpha }\left( 
\mathbf{1}_{I^{\prime }}\square _{I}^{\sigma ,\flat ,\mathbf{b}}f\right) ,%
\mathbf{1}_{J^{\prime }}\square _{J}^{\omega ,\flat ,\mathbf{b}^{\ast
}}g\right\rangle _{\omega }  \notag \\
&=&\sum_{I^{\prime }\in \mathfrak{C}\left( I\right) \text{ and }J^{\prime
}\in \mathfrak{C}\left( J\right) }E_{I^{\prime }}^{\sigma }\left( \widehat{%
\square }_{I}^{\sigma ,\flat ,\mathbf{b}}f\right) \ \left\langle T_{\sigma
}^{\alpha }\left( \mathbf{1}_{I^{\prime }}b_{I}\right) ,\mathbf{1}%
_{J^{\prime }}b_{J}^{\ast }\right\rangle _{\omega }\ E_{J^{\prime }}^{\omega
}\left( \widehat{\square }_{J}^{\omega ,\flat ,\mathbf{b}^{\ast }}g\right) .
\notag
\end{eqnarray}

It now remains to show that%
\begin{equation}
\boldsymbol{E}_{\Omega }^{\mathcal{D}}\boldsymbol{E}_{\Omega }^{\mathcal{G}%
}\sum_{I\in \mathcal{D}}\sum_{\substack{ J\in \mathcal{G}:\ 2^{-\mathbf{r}%
}\ell \left( I\right) <\ell \left( J\right) \leq \ell \left( I\right)  \\ %
d\left( J,I\right) \leq 2\ell \left( J\right) ^{\varepsilon }\ell \left(
I\right) ^{1-\varepsilon }}}\left\vert P\left( I,J\right) \right\vert
\lesssim \left( C_{\theta }\mathcal{NTV}_{\alpha }+\sqrt{\theta }\mathfrak{N}%
_{T^{\alpha }}\right) \left\Vert f\right\Vert _{L^{2}\left( \sigma \right)
}\left\Vert g\right\Vert _{L^{2}\left( \omega \right) }.
\label{must show final}
\end{equation}%
We will repeatedly use the inequality $\left\Vert \widehat{\square }%
_{I}^{\sigma ,\flat ,\mathbf{b}}f\right\Vert _{L^{2}\left( \sigma \right)
}\lesssim \left\Vert \square _{I}^{\sigma ,\mathbf{b}}f\right\Vert
_{L^{2}\left( \sigma \right) }^{\bigstar }$ which, upon noting (\ref{flat
broken}), follows from the $PLBP$ (\ref{plb}), 
\begin{eqnarray}
\left\Vert \widehat{\square }_{I}^{\sigma ,\flat ,\mathbf{b}}f\right\Vert
_{L^{2}\left( \sigma \right) } &\lesssim &\left\Vert b_{I}\widehat{\square }%
_{I}^{\sigma ,\flat ,\mathbf{b}}f\right\Vert _{L^{2}\left( \sigma \right)
}=\left\Vert \square _{I}^{\sigma ,\flat ,\mathbf{b}}f\right\Vert
_{L^{2}\left( \sigma \right) }  \label{box hat bound} \\
&\leq &\left\Vert \square _{I}^{\sigma ,\mathbf{b}}f\right\Vert
_{L^{2}\left( \sigma \right) }+\left\Vert \square _{I,\limfunc{broken}%
}^{\sigma ,\flat ,\mathbf{b}}f\right\Vert _{L^{2}\left( \sigma \right)
}\lesssim \left\Vert \square _{I}^{\sigma ,\mathbf{b}}f\right\Vert
_{L^{2}\left( \sigma \right) }^{\bigstar }.  \notag
\end{eqnarray}

Suppose now that $I\in \mathcal{C}_{A}$ for $A\in \mathcal{A}$, and that $%
J\in \mathcal{C}_{B}$ for $B\in \mathcal{B}$. Then the inner product in the
third line of (\ref{def P(I,J)}) becomes%
\begin{equation*}
\left\langle T_{\sigma }^{\alpha }\left( b_{I}\mathbf{1}_{I^{\prime
}}\right) ,b_{J}^{\ast }\mathbf{1}_{J^{\prime }}\right\rangle _{\omega
}=\left\langle T_{\sigma }^{\alpha }\left( b_{A}\mathbf{1}_{I^{\prime
}}\right) ,b_{B}^{\ast }\mathbf{1}_{J^{\prime }}\right\rangle _{\omega }\ ,
\end{equation*}%
and we will write this inner product in either form, depending on context.
We also introduce the following notation:%
\begin{equation*}
P_{\left( I,J\right) }\left( E,F\right) \equiv \left\langle T_{\sigma
}^{\alpha }\left( b_{I}\mathbf{1}_{E}\right) ,b_{J}^{\ast }\mathbf{1}%
_{F}\right\rangle _{\omega },\ \ \ \ \ \text{for any sets }E\text{ and }F,
\end{equation*}%
so that%
\begin{equation*}
P\left( I,J\right) =\sum_{I^{\prime }\in \mathfrak{C}\left( I\right) \text{
and }J^{\prime }\in \mathfrak{C}\left( J\right) }E_{I^{\prime }}^{\sigma
}\left( \widehat{\square }_{I}^{\sigma ,\flat ,\mathbf{b}}f\right) \
P_{\left( I,J\right) }\left( I^{\prime },J^{\prime }\right) \ E_{J^{\prime
}}^{\omega }\left( \widehat{\square }_{J}^{\omega ,\flat ,\mathbf{b}^{\ast
}}g\right) .
\end{equation*}%
The first thing we do is reduce matters to showing inequality (\ref{must
show final}) in the case of \emph{equal intervals}, by which we mean that $%
P_{\left( I,J\right) }\left( I^{\prime },J^{\prime }\right) $ is replaced by 
$P_{\left( I,J\right) }\left( I^{\prime }\cap J^{\prime },I^{\prime }\cap
J^{\prime }\right) $ in the terms $P\left( I,J\right) $ appearing in (\ref%
{must show final}). To see this let $K\equiv I^{\prime }\cap J^{\prime }$,
write%
\begin{eqnarray*}
\left\vert \left\langle T_{\sigma }^{\alpha }\left( b_{I}\mathbf{1}%
_{I^{\prime }}\right) ,b_{J}^{\ast }\mathbf{1}_{J^{\prime }}\right\rangle
_{\omega }\right\vert &=&\left\vert \left\langle T_{\sigma }^{\alpha }\left(
b_{I}\mathbf{1}_{I^{\prime }\setminus J^{\prime }}\right) ,b_{J}^{\ast }%
\mathbf{1}_{J^{\prime }}\right\rangle _{\omega }+\left\langle T_{\sigma
}^{\alpha }\left( b_{I}\mathbf{1}_{I^{\prime }\cap J^{\prime }}\right)
,b_{J}^{\ast }\mathbf{1}_{J^{\prime }\setminus I^{\prime }}\right\rangle
_{\omega }+\left\langle T_{\sigma }^{\alpha }\left( b_{I}\mathbf{1}%
_{I^{\prime }\cap J^{\prime }}\right) ,b_{J}^{\ast }\mathbf{1}_{J^{\prime
}\cap I^{\prime }}\right\rangle _{\omega }\right\vert \\
&\leq &\left\vert \left\langle T_{\sigma }^{\alpha }\left( b_{I}\mathbf{1}%
_{I^{\prime }\setminus J^{\prime }}\right) ,b_{J}^{\ast }\mathbf{1}%
_{J^{\prime }}\right\rangle _{\omega }\right\vert +\left\vert \left\langle
T_{\sigma }^{\alpha }\left( b_{I}\mathbf{1}_{K}\right) ,b_{J}^{\ast }\mathbf{%
1}_{J^{\prime }\setminus I^{\prime }}\right\rangle _{\omega }\right\vert
+\left\vert \left\langle T_{\sigma }^{\alpha }\left( b_{I}\mathbf{1}%
_{K}\right) ,b_{J}^{\ast }\mathbf{1}_{K}\right\rangle \right\vert ,
\end{eqnarray*}%
and use (\ref{disj supp}) to obtain%
\begin{equation*}
\left\vert \left\langle T_{\sigma }^{\alpha }\left( b_{I}\mathbf{1}%
_{I^{\prime }\setminus J^{\prime }}\right) ,b_{J}^{\ast }\mathbf{1}%
_{J^{\prime }}\right\rangle _{\omega }\right\vert +\left\vert \left\langle
T_{\sigma }^{\alpha }\left( b_{I}\mathbf{1}_{K}\right) ,b_{J}^{\ast }\mathbf{%
1}_{J^{\prime }\setminus I^{\prime }}\right\rangle _{\omega }\right\vert
\lesssim \sqrt{\mathfrak{A}_{2}^{\alpha }}\sqrt{\left\vert I^{\prime
}\right\vert _{\sigma }\left\vert J^{\prime }\right\vert _{\omega }}.
\end{equation*}%
It thus remains to consider only the term $P_{\left( I,J\right) }\left(
K,K\right) =\left\langle T_{\sigma }^{\alpha }\left( b_{I}\mathbf{1}%
_{K}\right) ,b_{J}^{\ast }\mathbf{1}_{K}\right\rangle _{\omega }$ where $%
K=I^{\prime }\cap J^{\prime }\neq \emptyset $, and to show that%
\begin{eqnarray}
&&\boldsymbol{E}_{\Omega }^{\mathcal{D}}\boldsymbol{E}_{\Omega }^{\mathcal{G}%
}\sum_{I\in \mathcal{D}}\sum_{\substack{ J\in \mathcal{G}:\ 2^{-\mathbf{r}%
}\ell \left( I\right) <\ell \left( J\right) \leq \ell \left( I\right)  \\ %
d\left( J,I\right) \leq 2\ell \left( J\right) ^{\varepsilon }\ell \left(
I\right) ^{1-\varepsilon }}}\left\vert \sum_{\substack{ I^{\prime }\in 
\mathfrak{C}\left( I\right) \text{ and }J^{\prime }\in \mathfrak{C}\left(
J\right)  \\ K=I^{\prime }\cap J^{\prime }\neq \emptyset }}E_{I^{\prime
}}^{\sigma }\left( \widehat{\square }_{I}^{\sigma ,\flat ,\mathbf{b}%
}f\right) \ P_{\left( I,J\right) }\left( K,K\right) \ E_{J^{\prime
}}^{\omega }\left( \widehat{\square }_{J}^{\omega ,\flat ,\mathbf{b}^{\ast
}}g\right) \right\vert  \label{ineq KK} \\
&\lesssim &\left( C_{\theta }\mathcal{NTV}_{\alpha }+\sqrt{\theta }\mathfrak{%
N}_{T^{\alpha }}\right) \left\Vert f\right\Vert _{L^{2}\left( \sigma \right)
}\left\Vert g\right\Vert _{L^{2}\left( \omega \right) },  \notag
\end{eqnarray}

\subsection{Random surgery}

However, we wish to further reduce matters to the case where $K\in \mathcal{G%
}$ and contained in $I^{\prime }\cap J^{\prime }$, and for this we use
random surgery and (\ref{disj supp}). First, we note that $I^{\prime }$
cannot be strictly contained in $J^{\prime }$ since $\ell \left( J\right)
\leq \ell \left( I\right) $, and in the case that $J^{\prime }\subset
I^{\prime }$, then $K=I^{\prime }\cap J^{\prime }=J^{\prime }\in \mathcal{G}$
and there is nothing more to do in this case. So we may assume that $%
J^{\prime }$ intersects both $I^{\prime }$ and its complement $\left(
I^{\prime }\right) ^{c}$.

Our first step is to reduce matters to showing inequality (\ref{ineq KK}) in
the case where $\ell \left( K\right) \geq \lambda \ell \left( I^{\prime
}\right) $ for a small positive number $\lambda \ll 2^{-\mathbf{r}}$. This
small constant $\lambda $, as well as the constant $\eta _{0}$ introduced in
probability estimates below, will be chosen sufficiently small at the end of
the proof to result in the term $\sqrt{\theta }\mathfrak{N}_{T^{\alpha }}$
appearing on the right hand side of (\ref{ineq KK}). This is accomplished by
writing (recall that $K=I^{\prime }\cap J^{\prime }$ for the moment)%
\begin{eqnarray*}
&&\sum_{I\in \mathcal{D}}\sum_{\substack{ J\in \mathcal{G}:\ 2^{-\mathbf{r}%
}\ell \left( I\right) <\ell \left( J\right) \leq \ell \left( I\right)  \\ %
d\left( J,I\right) \leq 2\ell \left( J\right) ^{\varepsilon }\ell \left(
I\right) ^{1-\varepsilon }}}\left\vert \sum_{\substack{ I^{\prime }\in 
\mathfrak{C}\left( I\right) \text{ and }J^{\prime }\in \mathfrak{C}\left(
J\right)  \\ K\neq \emptyset }}E_{I^{\prime }}^{\sigma }\left( \widehat{%
\square }_{I}^{\sigma ,\flat ,\mathbf{b}}f\right) \ P_{\left( I,J\right)
}\left( K,K\right) \ E_{J^{\prime }}^{\omega }\left( \widehat{\square }%
_{J}^{\omega ,\flat ,\mathbf{b}^{\ast }}g\right) \right\vert \\
&\leq &\sum_{I\in \mathcal{D}}\sum_{\substack{ J\in \mathcal{G}:\ 2^{-%
\mathbf{r}}\ell \left( I\right) <\ell \left( J\right) \leq \ell \left(
I\right)  \\ d\left( J,I\right) \leq 2\ell \left( J\right) ^{\varepsilon
}\ell \left( I\right) ^{1-\varepsilon }}}\left\vert \sum_{\substack{ %
I^{\prime }\in \mathfrak{C}\left( I\right) \text{ and }J^{\prime }\in 
\mathfrak{C}\left( J\right)  \\ \ell \left( K\right) \geq \lambda \ell
\left( I^{\prime }\right) }}E_{I^{\prime }}^{\sigma }\left( \widehat{\square 
}_{I}^{\sigma ,\flat ,\mathbf{b}}f\right) \ P_{\left( I,J\right) }\left(
K,K\right) \ E_{J^{\prime }}^{\omega }\left( \widehat{\square }_{J}^{\omega
,\flat ,\mathbf{b}^{\ast }}g\right) \right\vert \\
&&+\sum_{I\in \mathcal{D}}\sum_{\substack{ J\in \mathcal{G}:\ 2^{-\mathbf{r}%
}\ell \left( I\right) <\ell \left( J\right) \leq \ell \left( I\right)  \\ %
d\left( J,I\right) \leq 2\ell \left( J\right) ^{\varepsilon }\ell \left(
I\right) ^{1-\varepsilon }}}\left\vert \sum_{\substack{ I^{\prime }\in 
\mathfrak{C}\left( I\right) \text{ and }J^{\prime }\in \mathfrak{C}\left(
J\right)  \\ 0<\ell \left( K\right) <\lambda \ell \left( I^{\prime }\right) 
}}E_{I^{\prime }}^{\sigma }\left( \widehat{\square }_{I}^{\sigma ,\flat ,%
\mathbf{b}}f\right) \ P_{\left( I,J\right) }\left( K,K\right) \ E_{J^{\prime
}}^{\omega }\left( \widehat{\square }_{J}^{\omega ,\flat ,\mathbf{b}^{\ast
}}g\right) \right\vert \\
&\equiv &A+B.
\end{eqnarray*}%
Term $B$ is handled using the norm constant $\mathfrak{N}_{T^{\alpha }}$ and
probability, together with the estimate%
\begin{eqnarray*}
&&\sum_{I\in \mathcal{D}}\sum_{\substack{ J\in \mathcal{G}:\ 2^{-\mathbf{r}%
}\ell \left( I\right) <\ell \left( J\right) \leq \ell \left( I\right)  \\ %
d\left( J,I\right) \leq 2\ell \left( J\right) ^{\varepsilon }\ell \left(
I\right) ^{1-\varepsilon }}}\sum_{\substack{ I^{\prime }\in \mathfrak{C}%
\left( I\right) \text{ and }J^{\prime }\in \mathfrak{C}\left( J\right)  \\ %
0<\ell \left( K\right) <\lambda \ell \left( I^{\prime }\right) }}\left\vert 
\sqrt{\left\vert K\right\vert _{\omega }}\ E_{J^{\prime }}^{\omega }\left( 
\widehat{\square }_{J}^{\omega ,\flat ,\mathbf{b}^{\ast }}g\right)
\right\vert ^{2} \\
&=&\sum_{I\in \mathcal{D}}\sum_{\substack{ J\in \mathcal{G}:\ 2^{-\mathbf{r}%
}\ell \left( I\right) <\ell \left( J\right) \leq \ell \left( I\right)  \\ %
d\left( J,I\right) \leq 2\ell \left( J\right) ^{\varepsilon }\ell \left(
I\right) ^{1-\varepsilon }}}\sum_{\substack{ I^{\prime }\in \mathfrak{C}%
\left( I\right) \text{ and }J^{\prime }\in \mathfrak{C}\left( J\right)  \\ %
0<\ell \left( K\right) <\lambda \ell \left( I^{\prime }\right) }}\left\Vert 
\widehat{\square }_{J}^{\omega ,\flat ,\mathbf{b}^{\ast }}g\right\Vert
_{L^{2}\left( \omega \right) }^{2} \\
&\lesssim &\sum_{J\in \mathcal{G}}\left\Vert \widehat{\square }_{J}^{\omega
,\flat ,\mathbf{b}^{\ast }}g\right\Vert _{L^{2}\left( \omega \right)
}^{2}\lesssim \sum_{J\in \mathcal{G}}\left( \left\Vert \square _{J}^{\omega ,%
\mathbf{b}^{\ast }}g\right\Vert _{L^{2}\left( \omega \right)
}^{2}+\left\Vert \widehat{\square }_{J,\limfunc{broken}}^{\omega ,\flat ,%
\mathbf{b}^{\ast }}g\right\Vert _{L^{2}\left( \omega \right) }^{2}\right)
\lesssim \left\Vert g\right\Vert _{L^{2}\left( \omega \right) }^{2}\ ,
\end{eqnarray*}%
to obtain%
\begin{eqnarray*}
\boldsymbol{E}_{\Omega }^{\mathcal{G}}B &\lesssim &\boldsymbol{E}_{\Omega }^{%
\mathcal{G}}\sum_{I\in \mathcal{D}}\sum_{\substack{ J\in \mathcal{G}:\ 2^{-%
\mathbf{r}}\ell \left( I\right) <\ell \left( J\right) \leq \ell \left(
I\right)  \\ d\left( J,I\right) \leq 2\ell \left( J\right) ^{\varepsilon
}\ell \left( I\right) ^{1-\varepsilon }}}\sum_{\substack{ I^{\prime }\in 
\mathfrak{C}\left( I\right) \text{ and }J^{\prime }\in \mathfrak{C}\left(
J\right)  \\ 0<\ell \left( K\right) <\lambda \ell \left( I^{\prime }\right) 
}}\left\vert E_{I^{\prime }}^{\sigma }\left( \widehat{\square }_{I}^{\sigma
,\flat ,\mathbf{b}}f\right) \ \mathfrak{N}_{T^{\alpha }}\sqrt{\left\vert
K\right\vert _{\sigma }\left\vert K\right\vert _{\omega }}\ E_{J^{\prime
}}^{\omega }\left( \widehat{\square }_{J}^{\omega ,\flat ,\mathbf{b}^{\ast
}}g\right) \right\vert \\
&\lesssim &\boldsymbol{E}_{\Omega }^{\mathcal{G}}\mathfrak{N}_{T^{\alpha
}}\left( \sum_{I\in \mathcal{D}}\sum_{I^{\prime }\in \mathfrak{C}\left(
I\right) }\left( \sum_{\substack{ J\in \mathcal{G}:\ 2^{-\mathbf{r}}\ell
\left( I\right) <\ell \left( J\right) \leq \ell \left( I\right)  \\ d\left(
J,I\right) \leq 2\ell \left( J\right) ^{\varepsilon }\ell \left( I\right)
^{1-\varepsilon }}}\sum_{\substack{ J^{\prime }\in \mathfrak{C}\left(
J\right)  \\ 0<\ell \left( I^{\prime }\cap J^{\prime }\right) <\lambda \ell
\left( I^{\prime }\right) }}\left\vert I^{\prime }\cap J^{\prime
}\right\vert _{\sigma }\right) \left\vert E_{I^{\prime }}^{\sigma }\left( 
\widehat{\square }_{I}^{\sigma ,\flat ,\mathbf{b}}f\right) \right\vert
^{2}\right) ^{\frac{1}{2}}\left\Vert g\right\Vert _{L^{2}\left( \omega
\right) } \\
&\leq &\mathfrak{N}_{T^{\alpha }}\left( \sum_{I\in \mathcal{D}%
}\sum_{I^{\prime }\in \mathfrak{C}\left( I\right) }\lambda \left\vert
I^{\prime }\right\vert _{\sigma }\left\vert E_{I^{\prime }}^{\sigma }\left( 
\widehat{\square }_{I}^{\sigma ,\flat ,\mathbf{b}}f\right) \right\vert
^{2}\right) ^{\frac{1}{2}}\left\Vert g\right\Vert _{L^{2}\left( \omega
\right) }\lesssim \mathfrak{N}_{T^{\alpha }}\sqrt{\lambda }\left\Vert
f\right\Vert _{L^{2}\left( \sigma \right) }\left\Vert g\right\Vert
_{L^{2}\left( \omega \right) }
\end{eqnarray*}%
since the probability that, given an interval $I^{\prime }$, a grid $%
\mathcal{G}$ contains an interval $J^{\prime }$ with $2^{-\mathbf{r}}\ell
\left( I^{\prime }\right) <\ell \left( J^{\prime }\right) \leq \ell \left(
I^{\prime }\right) $ and $0<\ell \left( I^{\prime }\cap J^{\prime }\right)
<\lambda \ell \left( I^{\prime }\right) $, is at most $C_{\mathbf{r}}\lambda 
$ for some large constant $C_{\mathbf{r}}$ depending on the goodness
parameter $\mathbf{r}$ (we could also have appealed to the more general halo
estimate (\ref{hand'})). This term of course contributes to the conclusion
of the lemma provided $\lambda >0$ is chosen sufficiently small. Thus we are
left to control term $A$ in which $\ell \left( K\right) \geq \lambda \ell
\left( I^{\prime }\right) $ as required.

Now choose $\eta _{0}\in \left( 0,1\right) \cap \left\{ 2^{-m}\right\}
_{m\in \mathbb{N}}$, say $\eta _{0}=2^{-m}$ for $m$ sufficiently large. We
will assume from now on that $\eta _{0}<\frac{1}{2}\lambda $, where $\lambda 
$ is the small constant above. For any interval $L$ and $0<\eta \leq \frac{1%
}{2}$, define 
\begin{equation}
\partial _{\eta }L\equiv \left( 1+\eta \right) L\setminus \left( 1-\eta
\right) L  \label{def halo}
\end{equation}%
to be the `$\eta $-halo' around the boundary, i.e. endpoints, of $L$. One
should also recall that all intervals are assumed to be closed on the left
and open on the right. For $\eta _{1}\in \left( 0,\eta _{0}\right] $ we write%
\begin{eqnarray*}
I^{\prime }\cap J^{\prime } &=&\left\{ \left( I^{\prime }\setminus \partial
_{\eta _{1}}I^{\prime }\right) \cap J^{\prime }\right\} \cup \left\{ \left[
I^{\prime }\cap J^{\prime }\right] \setminus \left[ \left( I^{\prime
}\setminus \partial _{\eta _{1}}I^{\prime }\right) \cap J^{\prime }\right]
\right\} \\
&\equiv &M\overset{\cdot }{\cup }L,
\end{eqnarray*}%
where $\overset{\cdot }{\cup }$ denotes a disjoint union. At this point we
note that there is precisely one endpoint of $I^{\prime \prime }\equiv
I^{\prime }\setminus \partial _{\eta _{1}}I^{\prime }$ that lies in $%
J^{\prime }$ since $\ell \left( I^{\prime }\cap J^{\prime }\right) \geq
\lambda \ell \left( I^{\prime }\right) $ and $\eta _{1}\leq \eta _{0}<\frac{1%
}{2}\lambda $.

Moreover, and this is the key part of the argument, we can choose $\frac{1}{2%
}\eta _{0}\leq \eta _{1}\leq \eta _{0}$ so that the interval $M=I^{\prime
\prime }\cap J^{\prime }$ is a union of a number $B=B\left( I^{\prime
},J^{\prime }\right) =B\left( I^{\prime },J^{\prime },\eta _{1}\left(
I^{\prime },J^{\prime }\right) \right) $ of intervals $K_{i}\in \mathcal{G}$
each having side length $\ell \left( K_{i}\right) =2^{-m-1}\ell \left(
J^{\prime }\right) \equiv \frac{1}{2}\eta _{0}\ell \left( J^{\prime }\right) 
$ and where $B\leq C\frac{1}{\eta _{0}}$. Indeed, we take $\frac{1}{2}\eta
_{0}\leq \eta _{1}\leq \eta _{0}=2^{-m}$, so that the endpoint of the
interval $I^{\prime \prime }\equiv I^{\prime }\setminus \partial _{\eta
_{1}}I^{\prime }$ that lies in $J^{\prime }$ coincides with an endpoint of
some $K\in \mathcal{G}$ with $\ell \left( K\right) =2^{-m-1}\ell \left(
J^{\prime }\right) $. This can be arranged by varying $\eta _{1}$ between $%
\frac{1}{2}\eta _{0}$ and $\eta _{0}$ until the endpoint in question lies
among the dyadic numbers that form the endpoints of intervals in the grid $%
\mathcal{G}$ with side length $2^{-m-1}\ell \left( J^{\prime }\right) $. The
choice of intervals $\left\{ K_{i}\right\} _{i=1}^{B}$ having common side
length $2^{-m-1}\ell \left( J^{\prime }\right) $ is then uniquely
determined, and it is easy to see that 
\begin{equation*}
B\leq C\frac{1}{\eta _{0}}.
\end{equation*}%
Thus the choice of $\eta _{1}=\eta _{1}\left( I^{\prime },J^{\prime }\right)
>0$ is always at most $\eta _{0}$, and at least $\frac{1}{2}\eta _{0}$, but
changes according to the relative position of $I^{\prime }$ with respect to $%
J^{\prime }$ in order that $M$ is a union of intervals $K_{i}\in \mathcal{G}$
of side length at least $\frac{1}{2}\eta _{0}\ell \left( J^{\prime }\right) $%
, and specified in the manner described above. Define%
\begin{equation}
\mathcal{K}\left( I^{\prime },J^{\prime }\right) \equiv \left\{
K_{i}\right\} _{i=1}^{B\left( I^{\prime },J^{\prime }\right) }=\left\{
K_{i}\right\} _{i=1}^{B}  \label{def K(I',J')}
\end{equation}%
to be this collection of consecutive adjacent intervals $K_{i}$ uniquely
defined here in terms of $I^{\prime }$, $J^{\prime }$ and $\eta _{1}=\eta
_{1}\left( I^{\prime },J^{\prime }\right) $.

Having chosen the parameter $\eta _{1}$ as above, we now momentarily ignore
the decomposition of $M=\dbigcup\limits_{i=1}^{B}K_{i}$ into subintervals,
and return to the representation $K=I^{\prime }\cap J^{\prime }=M\overset{%
\cdot }{\cup }L$ determined by our choice of $\eta _{1}$. We have%
\begin{eqnarray}
\left\langle T_{\sigma }^{\alpha }\left( b_{I}\mathbf{1}_{K}\right)
,b_{J}^{\ast }\mathbf{1}_{K}\right\rangle _{\omega } &=&\left\langle
T_{\sigma }^{\alpha }\left( b_{I}\mathbf{1}_{M}\right) ,b_{J}^{\ast }\mathbf{%
1}_{L}\right\rangle _{\omega }+\left\langle T_{\sigma }^{\alpha }\left( b_{I}%
\mathbf{1}_{L}\right) ,b_{J}^{\ast }\mathbf{1}_{M}\right\rangle _{\omega }
\label{four terms} \\
&&+\left\langle T_{\sigma }^{\alpha }\left( b_{I}\mathbf{1}_{L}\right)
,b_{J}^{\ast }\mathbf{1}_{L}\right\rangle _{\omega }+\left\langle T_{\sigma
}^{\alpha }\left( b_{I}\mathbf{1}_{M}\right) ,b_{J}^{\ast }\mathbf{1}%
_{M}\right\rangle _{\omega }\ .  \notag
\end{eqnarray}%
Now we apply (\ref{disj supp}) to the first two terms in (\ref{four terms})
to obtain that%
\begin{eqnarray*}
&&\left\vert \left\langle T_{\sigma }^{\alpha }\left( b_{I}\mathbf{1}%
_{M}\right) ,b_{J}^{\ast }\mathbf{1}_{L}\right\rangle _{\omega }\right\vert
+\left\vert \left\langle T_{\sigma }^{\alpha }\left( b_{I}\mathbf{1}%
_{L}\right) ,b_{J}^{\ast }\mathbf{1}_{M}\right\rangle _{\omega }\right\vert
\\
&\lesssim &\sqrt{\mathfrak{A}_{2}^{\alpha }}\left\{ \sqrt{\int_{M}\left\vert
b_{I}\right\vert ^{2}d\sigma }\sqrt{\int_{L}\left\vert b_{J}^{\ast
}\right\vert ^{2}d\omega }+\sqrt{\int_{L}\left\vert b_{I}\right\vert
^{2}d\sigma }\sqrt{\int_{M}\left\vert b_{J}^{\ast }\right\vert ^{2}d\omega }%
\right\} \\
&&\ \ \ \ \ \ \ \ \ \ \ \ \ \ \ \ \ \ \ \ \lesssim \sqrt{\mathfrak{A}%
_{2}^{\alpha }}\sqrt{\left\vert I^{\prime }\right\vert _{\sigma }\left\vert
J^{\prime }\right\vert _{\omega }},
\end{eqnarray*}%
which when plugged appropriately into the left hand side of (\ref{must show
final}) is dominated by the right hand side of (\ref{must show final}):%
\begin{eqnarray*}
&&\boldsymbol{E}_{\Omega }^{\mathcal{D}}\boldsymbol{E}_{\Omega }^{\mathcal{G}%
}\sum_{I\in \mathcal{D}}\sum_{\substack{ J\in \mathcal{G}:\ 2^{-\mathbf{r}%
}\ell \left( I\right) <\ell \left( J\right) \leq \ell \left( I\right)  \\ %
d\left( J,I\right) \leq 2\ell \left( J\right) ^{\varepsilon }\ell \left(
I\right) ^{1-\varepsilon }}}\sum_{\substack{ I^{\prime }\in \mathfrak{C}%
\left( I\right)  \\ J^{\prime }\in \mathfrak{C}\left( J\right) }} \\
&&\ \ \ \ \ \ \ \ \ \ \times \left\vert E_{I^{\prime }}^{\sigma }\left( 
\widehat{\square }_{I}^{\sigma ,\flat ,\mathbf{b}}f\right) \ \left(
\left\vert \left\langle T_{\sigma }^{\alpha }\left( b_{I}\mathbf{1}%
_{M}\right) ,b_{J}^{\ast }\mathbf{1}_{L}\right\rangle _{\omega }\right\vert
+\left\vert \left\langle T_{\sigma }^{\alpha }\left( b_{I}\mathbf{1}%
_{L}\right) ,b_{J}^{\ast }\mathbf{1}_{M}\right\rangle _{\omega }\right\vert
\right) \ E_{J^{\prime }}^{\omega }\left( \widehat{\square }_{J}^{\omega
,\flat ,\mathbf{b}^{\ast }}g\right) \right\vert \\
&\lesssim &\boldsymbol{E}_{\Omega }^{\mathcal{D}}\boldsymbol{E}_{\Omega }^{%
\mathcal{G}}\sum_{I\in \mathcal{D}}\sum_{\substack{ J\in \mathcal{G}:\ 2^{-%
\mathbf{r}}\ell \left( I\right) <\ell \left( J\right) \leq \ell \left(
I\right)  \\ d\left( J,I\right) \leq 2\ell \left( J\right) ^{\varepsilon
}\ell \left( I\right) ^{1-\varepsilon }}}\sum_{\substack{ I^{\prime }\in 
\mathfrak{C}\left( I\right)  \\ J^{\prime }\in \mathfrak{C}\left( J\right) }}%
\left\vert E_{I^{\prime }}^{\sigma }\left( \widehat{\square }_{I}^{\sigma
,\flat ,\mathbf{b}}f\right) \right\vert \sqrt{\mathfrak{A}_{2}^{\alpha }}%
\sqrt{\left\vert I^{\prime }\right\vert _{\sigma }\left\vert J^{\prime
}\right\vert _{\omega }}\left\vert E_{J^{\prime }}^{\omega }\left( \widehat{%
\square }_{J}^{\omega ,\flat ,\mathbf{b}^{\ast }}g\right) \right\vert \\
&\lesssim &\sqrt{\mathfrak{A}_{2}^{\alpha }}\left\Vert f\right\Vert
_{L^{2}\left( \sigma \right) }\left\Vert g\right\Vert _{L^{2}\left( \omega
\right) }\ ,
\end{eqnarray*}%
where in the last line we have used $\sum_{\substack{ I^{\prime }\in 
\mathfrak{C}\left( I\right)  \\ J^{\prime }\in \mathfrak{C}\left( J\right) }}%
\left\vert E_{I^{\prime }}^{\sigma }\left( \widehat{\square }_{I}^{\sigma
,\flat ,\mathbf{b}}f\right) \right\vert ^{2}\left\vert I^{\prime
}\right\vert _{\sigma }=\left\Vert \widehat{\square }_{I}^{\sigma ,\flat ,%
\mathbf{b}}f\right\Vert _{L^{2}\left( \sigma \right) }^{2}$ and (\ref{box
hat bound}) and the frame inequalities in Appendix A. Then we apply
Cauchy-Schwarz to the sums in the third term in (\ref{four terms}) using $%
L=\partial _{\eta _{1}}I^{\prime }\cap J^{\prime }$ to get%
\begin{eqnarray*}
&&\boldsymbol{E}_{\Omega }^{\mathcal{D}}\boldsymbol{E}_{\Omega }^{\mathcal{G}%
}\sum_{I\in \mathcal{D}}\sum_{\substack{ J\in \mathcal{G}:\ 2^{-\mathbf{r}%
}\ell \left( I\right) <\ell \left( J\right) \leq \ell \left( I\right)  \\ %
d\left( J,I\right) \leq 2\ell \left( J\right) ^{\varepsilon }\ell \left(
I\right) ^{1-\varepsilon }}}\sum_{\substack{ I^{\prime }\in \mathfrak{C}%
\left( I\right)  \\ J^{\prime }\in \mathfrak{C}\left( J\right) }}\left\vert
E_{I^{\prime }}^{\sigma }\left( \widehat{\square }_{I}^{\sigma ,\flat ,%
\mathbf{b}}f\right) \ \left\langle T_{\sigma }^{\alpha }\left( b_{I}\mathbf{1%
}_{L}\right) ,b_{J}^{\ast }\mathbf{1}_{L}\right\rangle _{\omega }\
E_{J^{\prime }}^{\omega }\left( \widehat{\square }_{J}^{\omega ,\flat ,%
\mathbf{b}^{\ast }}g\right) \right\vert \\
&\lesssim &\boldsymbol{E}_{\Omega }^{\mathcal{G}}\mathfrak{N}_{T^{\alpha
}}\left\Vert f\right\Vert _{L^{2}\left( \sigma \right) }\boldsymbol{E}%
_{\Omega }^{\mathcal{D}}\sqrt{\sum_{I\in \mathcal{D}}\sum_{\substack{ J\in 
\mathcal{G}:\ 2^{-\mathbf{r}}\ell \left( I\right) <\ell \left( J\right) \leq
\ell \left( I\right)  \\ d\left( J,I\right) \leq 2\ell \left( J\right)
^{\varepsilon }\ell \left( I\right) ^{1-\varepsilon }}}\sum_{\substack{ %
I^{\prime }\in \mathfrak{C}\left( I\right)  \\ J^{\prime }\in \mathfrak{C}%
\left( J\right) }}\left( \int_{\partial _{\eta _{1}}I^{\prime }\cap
J^{\prime }}\left\vert b_{J}^{\ast }\right\vert ^{2}d\omega \right)
\left\vert E_{J^{\prime }}^{\omega }\left( \widehat{\square }_{J}^{\omega
,\flat ,\mathbf{b}^{\ast }}g\right) \right\vert ^{2}},
\end{eqnarray*}%
using (\ref{box hat bound}) and the frame inequalities in Appendix A again.
Then using Cauchy-Schwarz on the expectation $\boldsymbol{E}_{\Omega }^{%
\mathcal{D}}$, this is dominated by 
\begin{eqnarray*}
&\lesssim &\boldsymbol{E}_{\Omega }^{\mathcal{G}}\mathfrak{N}_{T^{\alpha
}}\left\Vert f\right\Vert _{L^{2}\left( \sigma \right) }\sqrt{\sum_{J\in 
\mathcal{G}}\sum_{J^{\prime }\in \mathfrak{C}\left( J\right) }\left( 
\boldsymbol{E}_{\Omega }^{\mathcal{D}}\sum_{\substack{ I\in \mathcal{D}:\
2^{-\mathbf{r}}\ell \left( I\right) <\ell \left( J\right) \leq \ell \left(
I\right)  \\ d\left( J,I\right) \leq 2\ell \left( J\right) ^{\varepsilon
}\ell \left( I\right) ^{1-\varepsilon }  \\ I^{\prime }\in \mathfrak{C}%
\left( I\right) }}\left\vert \partial _{\eta _{1}}I^{\prime }\cap J^{\prime
}\right\vert _{\omega }\right) \left\vert E_{J^{\prime }}^{\omega }\left( 
\widehat{\square }_{J}^{\omega ,\flat ,\mathbf{b}^{\ast }}g\right)
\right\vert ^{2}} \\
&\lesssim &\boldsymbol{E}_{\Omega }^{\mathcal{G}}\mathfrak{N}_{T^{\alpha
}}\left\Vert f\right\Vert _{L^{2}\left( \sigma \right) }\sqrt{\sum_{J\in 
\mathcal{G}}\sum_{J^{\prime }\in \mathfrak{C}\left( J\right) }C2^{\mathbf{r}%
}\left( \boldsymbol{E}_{\Omega }^{\mathcal{D}}\sum_{I^{\prime }\in \mathcal{D%
}:\ell \left( J^{\prime }\right) \leq \ell \left( I^{\prime }\right) \leq 2^{%
\mathbf{r}}\ell \left( J^{\prime }\right) }\left\vert \partial _{\eta
_{0}}I^{\prime }\cap J^{\prime }\right\vert _{\omega }\right) \left\vert
E_{J^{\prime }}^{\omega }\left( \widehat{\square }_{J}^{\omega ,\flat ,%
\mathbf{b}^{\ast }}g\right) \right\vert ^{2}} \\
&\lesssim &\sqrt{\eta _{0}}\mathfrak{N}_{T^{\alpha }}\left\Vert f\right\Vert
_{L^{2}\left( \sigma \right) }\left\Vert g\right\Vert _{L^{2}\left( \omega
\right) }\ ,
\end{eqnarray*}%
where in the last line we have used $\eta _{1}\leq \eta _{0}$, and then 
\begin{equation*}
\boldsymbol{E}_{\Omega }^{\mathcal{D}}\sum_{I^{\prime }\in \mathcal{D}:\ell
\left( J^{\prime }\right) \leq \ell \left( I^{\prime }\right) \leq 2^{%
\mathbf{r}}\ell \left( J^{\prime }\right) }\left\vert \partial _{\eta
_{0}}I^{\prime }\cap J^{\prime }\right\vert _{\omega }\lesssim \eta
_{0}\left\vert J^{\prime }\right\vert _{\omega }
\end{equation*}%
from (\ref{hand'}) since $\eta _{0}\leq \frac{1}{2}\lambda \ll 2^{-\mathbf{r}%
}$, and finally quasiorthogonality and (\ref{box hat bound}) yet again, to
obtain 
\begin{eqnarray*}
\sum_{J\in \mathcal{G}}\sum_{J^{\prime }\in \mathfrak{C}\left( J\right)
}\left\vert J^{\prime }\right\vert _{\omega }\left\vert E_{J^{\prime
}}^{\omega }\left( \widehat{\square }_{J}^{\omega ,\flat ,\mathbf{b}^{\ast
}}g\right) \right\vert ^{2} &\lesssim &\sum_{J\in \mathcal{G}}\left\Vert 
\widehat{\square }_{J}^{\omega ,\flat ,\mathbf{b}^{\ast }}g\right\Vert
_{L^{2}\left( J\right) }^{2}\lesssim \sum_{J\in \mathcal{G}}\left\Vert
\square _{J}^{\omega ,\flat ,\mathbf{b}^{\ast }}g\right\Vert _{L^{2}\left(
J\right) }^{2} \\
&\lesssim &\sum_{J\in \mathcal{G}}\left( \left\Vert \square _{J}^{\omega ,%
\mathbf{b}^{\ast }}g\right\Vert _{L^{2}\left( J\right) }^{2}+\left\Vert
\nabla _{J}^{\omega }g\right\Vert _{L^{2}\left( J\right) }^{2}\right)
\lesssim \left\Vert g\right\Vert _{L^{2}\left( \omega \right) }^{2}\ .
\end{eqnarray*}

This leaves us to estimate the fourth term in (\ref{four terms}), i.e. the
inner product $\left\langle T_{\sigma }^{\alpha }\left( b_{I}\mathbf{1}%
_{M}\right) ,b_{J}^{\ast }\mathbf{1}_{M}\right\rangle $. It is at this point
that we will use the decomposition $M=\overset{\cdot }{\dbigcup }_{1\leq
i\leq B}K_{i}$ constructed above. We have%
\begin{equation*}
\left\langle T_{\sigma }^{\alpha }\left( b_{I}\mathbf{1}_{M}\right)
,b_{J}^{\ast }\mathbf{1}_{M}\right\rangle _{\omega
}=\dsum\limits_{i,i^{\prime }=1}^{B}\left\langle T_{\sigma }^{\alpha }\left(
b_{I}\mathbf{1}_{K_{i}}\right) ,b_{J}^{\ast }\mathbf{1}_{K_{i^{\prime
}}}\right\rangle _{\omega }=\dsum\limits_{i=1}^{B}\left\langle T_{\sigma
}^{\alpha }\left( b_{I}\mathbf{1}_{K_{i}}\right) ,b_{J}^{\ast }\mathbf{1}%
_{K_{i}}\right\rangle _{\omega }+\dsum\limits_{i\neq i^{\prime
}}\left\langle T_{\sigma }^{\alpha }\left( b_{I}\mathbf{1}_{K_{i}}\right)
,b_{J}^{\ast }\dsum\limits_{i^{\prime }:\ i^{\prime }\neq i}\mathbf{1}%
_{K_{i^{\prime }}}\right\rangle _{\omega },
\end{equation*}%
and finally, we can use (\ref{disj supp}) once more on the second sum above
to reduce matters, modulo a constant multiple of $\frac{1}{\eta _{0}}$, to
the case of estimating the inner products $\left\langle T_{\sigma }^{\alpha
}\left( b_{I}\mathbf{1}_{K}\right) ,b_{J}^{\ast }\mathbf{1}_{K}\right\rangle 
$ for the intervals $K=K_{i}\in \mathcal{G}$, $1\leq i\leq B\leq C\frac{1}{%
\eta _{0}}$, which are contained in $I^{\prime }\cap J^{\prime }$. Thus it
remains to show%
\begin{eqnarray}
&&  \label{after prob} \\
&&\boldsymbol{E}_{\Omega }^{\mathcal{D}}\boldsymbol{E}_{\Omega }^{\mathcal{G}%
}\sum_{I\in \mathcal{D}}\sum_{\substack{ J\in \mathcal{G}:\ 2^{-\mathbf{r}%
}\ell \left( I\right) <\ell \left( J\right) \leq \ell \left( I\right)  \\ %
d\left( J,I\right) \leq 2\ell \left( J\right) ^{\varepsilon }\ell \left(
I\right) ^{1-\varepsilon }}}\left\vert \sum_{\substack{ I^{\prime }\in 
\mathfrak{C}\left( I\right) \text{ and }J^{\prime }\in \mathfrak{C}\left(
J\right)  \\ K\in \mathcal{K}\left( I^{\prime },J^{\prime }\right) }}%
E_{I^{\prime }}^{\sigma }\left( \widehat{\square }_{I}^{\sigma ,\flat ,%
\mathbf{b}}f\right) \ P_{\left( I,J\right) }\left( K,K\right) \ E_{J^{\prime
}}^{\omega }\left( \widehat{\square }_{J}^{\omega ,\flat ,\mathbf{b}^{\ast
}}g\right) \right\vert  \notag \\
&\lesssim &\left( C_{\theta }\mathcal{NTV}_{\alpha }+\sqrt{\theta }\mathfrak{%
N}_{T^{\alpha }}\right) \left\Vert f\right\Vert _{L^{2}\left( \sigma \right)
}\left\Vert g\right\Vert _{L^{2}\left( \omega \right) },  \notag
\end{eqnarray}%
where we recall that $\mathcal{K}\left( I^{\prime },J^{\prime }\right)
=\left\{ K_{i}\right\} _{i=1}^{B}$, and where $B=B\left( I^{\prime
},J^{\prime }\right) $ depends on the pair $\left( I^{\prime },J^{\prime
}\right) $ but is bounded by $C\frac{1}{\eta _{0}}$ independent of $\left(
I^{\prime },J^{\prime }\right) $ and the choice of $\eta _{1}$, and 
\begin{equation*}
K_{i}\in \mathcal{G},\ K_{i}\subset I^{\prime }\cap J^{\prime },\ \ell
\left( K_{i}\right) =2^{-m-1}\ell \left( J^{\prime }\right) ,\ \ \ \ \ 1\leq
i\leq B\mathfrak{.}
\end{equation*}

There will be just one more use of random probability in dealing with the
nearby form, and that will occur at the end of the finite iteration in
Subsection \ref{Subsection iteration} below.

\subsection{Return of the original testing function}

We now consider the inner product $\left\langle T_{\sigma }^{\alpha }\left(
b_{A}\mathbf{1}_{K}\right) ,b_{B}^{\ast }\mathbf{1}_{K}\right\rangle
_{\omega }$ and estimate the case when%
\begin{equation*}
K\in \mathcal{G},\ K\subset I^{\prime }\cap J^{\prime },\ I^{\prime }\in 
\mathfrak{C}\left( I\right) ,\ J^{\prime }\in \mathfrak{C}\left( J\right) ,\
I\in \mathcal{C}_{A}^{\mathcal{A}},\ J\in \mathcal{C}_{B}^{\mathcal{B}},\
\ell \left( K\right) =2^{-m-1}\ell \left( J^{\prime }\right) \mathfrak{.}
\end{equation*}%
Recall that for $\eta \in \left( 0,\frac{1}{2}\right] $ and any interval $K$
, we defined $\partial _{\eta }K\equiv \left( 1+\eta \right) K\setminus
\left( 1-\eta \right) K$ to be the `$\eta $-halo' around the boundary, i.e.
endpoints, of $K$. In what follows we will now take $\eta =\frac{1}{2}$ and
invoke \emph{deterministic} surgery with $\eta $-halos (which are $\frac{1}{2%
}$-halos), together with the energy condition and one last application of
random surgery, as follows. For subsets $E,F\subset A\cap B$ and intervals $%
K\subset A\cap B$ we define%
\begin{equation}
\left\{ E,F\right\} \equiv \left\langle T_{\sigma }^{\alpha }\left( b_{A}%
\mathbf{1}_{E}\right) ,b_{B}^{\ast }\mathbf{1}_{F}\right\rangle _{\omega }\ ,
\label{def E,F}
\end{equation}%
and%
\begin{equation*}
K_{\limfunc{in}}\equiv K\setminus \partial _{\eta }K\text{ and }K_{\limfunc{%
out}}\equiv K\cap \partial _{\eta }K\ ,
\end{equation*}%
and we write%
\begin{equation}
\left\{ K,K\right\} =\left\{ A,K_{\limfunc{in}}\right\} -\left\{ A\setminus
K,K_{\limfunc{in}}\right\} +\left\{ K_{\limfunc{out}},K_{\limfunc{out}%
}\right\} +\left\{ K_{\limfunc{in}},K_{\limfunc{out}}\right\} .  \label{K,K}
\end{equation}

Note that the first two terms on the right hand side of (\ref{K,K})
decompose the inner product $\left\{ K,K_{\limfunc{in}}\right\} $, which
`includes' the difficult symmetric inner product $\left\{ K_{\limfunc{in}%
},K_{\limfunc{in}}\right\} $. Thus the difficult symmetric inner product is
ultimately controlled by testing on the interval $A$ to handle$\ \left\{
A,K_{\limfunc{in}}\right\} $, and by using a trick that resurrects the
original testing functions $\left\{ b_{J}^{\ast ,\limfunc{orig}}\right\}
_{J\in \mathcal{G}}$, discarded in the corona constructions above, to handle 
$\left\{ A\setminus K,K_{\limfunc{in}}\right\} $. More precisely, these
original testing functions $b_{J}^{\ast ,\limfunc{orig}}$ are the testing
functions obtained after reducing matters to the case of bounded testing
functions with the pointwise lower bound property $PLBP$ as in Conclusion %
\ref{bounded PLBP} above.

The first term on the right side of (\ref{K,K}) satisfies%
\begin{equation}
\left\vert \left\{ A,K_{\limfunc{in}}\right\} \right\vert =\left\vert
\int_{K_{\limfunc{in}}}\left( T_{\sigma }^{\alpha }b_{A}\right) b_{B}^{\ast
}d\omega \right\vert \leq \left\Vert \mathbf{1}_{K_{\limfunc{in}}}T_{\sigma
}^{\alpha }b_{A}\right\Vert _{L^{2}\left( \omega \right) }\left\Vert \mathbf{%
1}_{K_{\limfunc{in}}}b_{B}^{\ast }\right\Vert _{L^{2}\left( \omega \right)
}\leq \left\Vert b_{B}^{\ast }\right\Vert _{\infty }\left\Vert \mathbf{1}%
_{K_{\limfunc{in}}}T_{\sigma }^{\alpha }b_{A}\right\Vert _{L^{2}\left(
\omega \right) }\sqrt{\left\vert K_{\limfunc{in}}\right\vert _{\omega }}\ .
\label{AKin}
\end{equation}%
Before proceeding further it will prove convenient to introduce some
additional notation, namely we will write the energy estimate in the second
display of the Energy Lemma as%
\begin{equation}
\left\vert \left\langle T^{\alpha }\nu ,\Psi _{J}\right\rangle _{\omega
}\right\vert \lesssim C_{\gamma }\ \mathrm{P}_{\delta }^{\alpha }\mathsf{Q}%
^{\omega }\left( J,\upsilon \right) \ \left\Vert \Psi _{J}\right\Vert
_{L^{2}\left( \mu \right) }^{\bigstar },\ \ \ \ \ \text{if }\int \Psi
_{J}d\omega =0\text{ and }\gamma J\cap \func{Supp}\nu =\emptyset ,\gamma >1,
\label{star}
\end{equation}%
where%
\begin{equation}
\mathrm{P}_{\delta }^{\alpha }\mathsf{Q}^{\omega }\left( J,\upsilon \right)
\equiv \frac{\mathrm{P}^{\alpha }\left( J,\nu \right) }{\left\vert
J\right\vert }\left\Vert \mathsf{Q}_{J}^{\omega ,\mathbf{b}^{\ast
}}x\right\Vert _{L^{2}\left( \omega \right) }^{\spadesuit }+\frac{\mathrm{P}%
_{1+\delta }^{\alpha }\left( J,\nu \right) }{\left\vert J\right\vert }%
\left\Vert x-m_{J}\right\Vert _{L^{2}\left( \mathbf{1}_{J}\omega \right) }\ .
\label{def compact}
\end{equation}%
The use of the compact notation $\mathrm{P}_{\delta }^{\alpha }\mathsf{Q}%
^{\omega }\left( J,\upsilon \right) $ to denote the complicated expression
on the right hand side will considerably reduce the size of many subsequent
displays.

Let $K_{\limfunc{in}}^{\limfunc{left}}$ and $K_{\limfunc{in}}^{\limfunc{right%
}}$ denote the left and right children of $K_{\limfunc{in}}$, which until
now have been written as $\left\{ K_{\ell }^{\prime \prime }\right\} _{\ell
=1}^{2}$ in no particular order, and we will continue to use both of these
notations. We now claim that the second term on the right side of (\ref{K,K}%
) satisfies 
\begin{eqnarray}
\left\vert \left\{ A\setminus K,K_{\limfunc{in}}\right\} \right\vert
&\lesssim &\left\{ \mathrm{P}_{\delta }^{\alpha }\mathsf{Q}^{\omega }\left(
K_{\limfunc{in}}^{\limfunc{left}},\mathbf{1}_{A\setminus K}\sigma \right) +%
\mathrm{P}_{\delta }^{\alpha }\mathsf{Q}^{\omega }\left( K_{\limfunc{in}}^{%
\limfunc{right}},\mathbf{1}_{A\setminus K}\sigma \right) \right\} \sqrt{%
\left\vert K_{\limfunc{in}}\right\vert _{\omega }}  \label{second term} \\
&&+\left( \sqrt{\int_{K_{\limfunc{in}}}\left\vert T_{\sigma }^{\alpha
}b_{A}\right\vert ^{2}d\omega }+\left( \mathfrak{T}_{T^{\alpha }}+\sqrt{%
\mathfrak{A}_{2}^{\alpha }}\right) \sqrt{\left\vert K_{\limfunc{in}%
}\right\vert _{\sigma }}\right) \sqrt{\left\vert K_{\limfunc{in}}\right\vert
_{\omega }}  \notag \\
&&+\sum_{\ell =1}^{2}\left\vert \left\langle T_{\sigma }^{\alpha }b_{A}%
\mathbf{1}_{K_{\limfunc{out}}},b_{K_{\ell }^{\prime \prime }}^{\ast ,%
\limfunc{orig}}\right\rangle _{\omega }\right\vert ,  \notag
\end{eqnarray}%
upon using a trick with the \textbf{original} testing functions $b_{K_{\ell
}^{\prime \prime }}^{\ast ,\limfunc{orig}}$ for the two grandchildren $%
\left\{ K_{\ell }^{\prime \prime }\right\} _{\ell =1}^{2}$ of the interval $%
K $ that lie strictly inside $K$, and whose union is $\frac{1}{2}K$. Indeed,
to prove (\ref{second term}), we use the following identity whose proof is
immediate (and whose origin will be made clear in the discussion below):%
\begin{eqnarray}
&&\left\langle T_{\sigma }^{\alpha }\left( b_{A}\mathbf{1}_{A\setminus
K}\right) ,\mathbf{1}_{K_{\limfunc{in}}}b_{B}^{\ast }\right\rangle _{\omega
}-\sum_{\ell =1}^{2}\left( \frac{\frac{1}{\left\vert K_{\ell }^{\prime
\prime }\right\vert _{\omega }}\int_{K_{\ell }^{\prime \prime }}b_{B}^{\ast
}d\omega }{\frac{1}{\left\vert K_{\ell }^{\prime \prime }\right\vert
_{\omega }}\int_{K_{\ell }^{\prime \prime }}b_{K_{\ell }^{\prime \prime
}}^{\ast ,\limfunc{orig}}d\omega }\right) \left\langle T_{\sigma }^{\alpha
}\left( b_{A}\mathbf{1}_{A\setminus K}\right) ,b_{K_{\ell }^{\prime \prime
}}^{\ast ,\limfunc{orig}}\right\rangle _{\omega }  \label{big identity} \\
&=&\left\langle T_{\sigma }^{\alpha }\left( b_{A}\mathbf{1}_{A\setminus
K}\right) ,\mathbf{1}_{K_{\limfunc{in}}}b_{B}^{\ast }\right\rangle _{\omega }
\notag \\
&&-\sum_{\ell =1}^{2}\left( \frac{\frac{1}{\left\vert K_{\ell }^{\prime
\prime }\right\vert _{\omega }}\int_{K_{\ell }^{\prime \prime }}b_{B}^{\ast
}d\omega }{\frac{1}{\left\vert K_{\ell }^{\prime \prime }\right\vert
_{\omega }}\int_{K_{\ell }^{\prime \prime }}b_{K_{\ell }^{\prime \prime
}}^{\ast ,\limfunc{orig}}d\omega }\right) \left\{ \left\langle T_{\sigma
}^{\alpha }b_{A},b_{K_{\ell }^{\prime \prime }}^{\ast ,\limfunc{orig}%
}\right\rangle _{\omega }-\left\langle b_{A}\mathbf{1}_{K_{\limfunc{in}%
}},T_{\omega }^{\alpha ,\ast }b_{K_{\ell }^{\prime \prime }}^{\ast ,\limfunc{%
orig}}\right\rangle _{\sigma }-\left\langle T_{\sigma }^{\alpha }b_{A}%
\mathbf{1}_{K_{\limfunc{out}}},b_{K_{\ell }^{\prime \prime }}^{\ast ,%
\limfunc{orig}}\right\rangle _{\omega }\right\}  \notag \\
&=&\left\langle T_{\sigma }^{\alpha }\left( b_{A}\mathbf{1}_{A\setminus
K}\right) ,\mathbf{1}_{K_{\limfunc{in}}}b_{B}^{\ast }\right\rangle _{\omega
}+\sum_{\ell =1}^{2}\left( \frac{\frac{1}{\left\vert K_{\ell }^{\prime
\prime }\right\vert _{\omega }}\int_{K_{\ell }^{\prime \prime }}b_{B}^{\ast
}d\omega }{\frac{1}{\left\vert K_{\ell }^{\prime \prime }\right\vert
_{\omega }}\int_{K_{\ell }^{\prime \prime }}b_{K_{\ell }^{\prime \prime
}}^{\ast ,\limfunc{orig}}d\omega }\right) \left\langle T_{\sigma }^{\alpha
}b_{A}\mathbf{1}_{K_{\limfunc{out}}},b_{K_{\ell }^{\prime \prime }}^{\ast ,%
\limfunc{orig}}\right\rangle _{\omega }  \notag \\
&&-\sum_{\ell =1}^{2}\left( \frac{\frac{1}{\left\vert K_{\ell }^{\prime
\prime }\right\vert _{\omega }}\int_{K_{\ell }^{\prime \prime }}b_{B}^{\ast
}d\omega }{\frac{1}{\left\vert K_{\ell }^{\prime \prime }\right\vert
_{\omega }}\int_{K_{\ell }^{\prime \prime }}b_{K_{\ell }^{\prime \prime
}}^{\ast ,\limfunc{orig}}d\omega }\right) \left\{ \left\langle T_{\sigma
}^{\alpha }b_{A},b_{K_{\ell }^{\prime \prime }}^{\ast ,\limfunc{orig}%
}\right\rangle _{\omega }-\left\langle b_{A}\mathbf{1}_{K_{\limfunc{in}%
}},T_{\omega }^{\alpha ,\ast }b_{K_{\ell }^{\prime \prime }}^{\ast ,\limfunc{%
orig}}\right\rangle _{\sigma }\right\} .  \notag
\end{eqnarray}%
In fact we have the following estimate, more precise than (\ref{second term}%
).

\begin{lemma}
\label{preiterate}We have 
\begin{eqnarray*}
&&\left\vert \left\{ A\setminus K,K_{\limfunc{in}}\right\} +\sum_{\ell
=1}^{2}\left( \frac{\frac{1}{\left\vert K_{\ell }^{\prime \prime
}\right\vert _{\omega }}\int_{K_{\ell }^{\prime \prime }}b_{B}^{\ast
}d\omega }{\frac{1}{\left\vert K_{\ell }^{\prime \prime }\right\vert
_{\omega }}\int_{K_{\ell }^{\prime \prime }}b_{K_{\ell }^{\prime \prime
}}^{\ast ,\limfunc{orig}}d\omega }\right) \left\{ K_{\limfunc{out}},K_{%
\limfunc{in}}^{\ell }\right\} ^{\limfunc{orig}}\right\vert \\
&\lesssim &\left( \left\{ \mathrm{P}_{\delta }^{\alpha }\mathsf{Q}^{\omega
}\left( K_{\limfunc{in}}^{\limfunc{left}},\mathbf{1}_{A\setminus K}\sigma
\right) +\mathrm{P}_{\delta }^{\alpha }\mathsf{Q}^{\omega }\left( K_{%
\limfunc{in}}^{\limfunc{right}},\mathbf{1}_{A\setminus K}\sigma \right)
\right\} +\sqrt{\int_{K_{\limfunc{in}}}\left\vert T_{\sigma }^{\alpha
}b_{A}\right\vert ^{2}d\omega }+\left( \mathfrak{T}_{T^{\alpha ,\ast }}+%
\sqrt{\mathfrak{A}_{2}^{\alpha }}\right) \sqrt{\left\vert K_{\limfunc{in}%
}\right\vert _{\sigma }}\right) \sqrt{\left\vert K_{\limfunc{in}}\right\vert
_{\omega }},
\end{eqnarray*}%
where $b_{K_{\ell }^{\prime \prime }}^{\ast ,\limfunc{orig}}$ is the
original testing function for the grandchild $K_{\ell }^{\prime \prime }$ of 
$K$ in Conclusion \ref{bounded PLBP} above, and where%
\begin{equation*}
\left\{ K_{\limfunc{out}},K_{\limfunc{in}}^{\ell }\right\} ^{\limfunc{orig}%
}\equiv \left\langle T_{\sigma }^{\alpha }b_{A}\mathbf{1}_{K_{\limfunc{out}%
}},b_{K_{\ell }^{\prime \prime }}^{\ast ,\limfunc{orig}}\right\rangle
_{\omega }\ \text{for }\ell \in \left\{ 1,2\right\} .
\end{equation*}
\end{lemma}

Before starting the proof of the lemma, we motivate the identity (\ref{big
identity}) with the following discussion. A simpler way to start the
analysis for for $\left\{ A\setminus K,K_{\limfunc{in}}\right\}
=\left\langle T_{\sigma }^{\alpha }\left( b_{A}\mathbf{1}_{A\setminus
K}\right) ,\mathbf{1}_{K_{\limfunc{in}}}b_{B}^{\ast }\right\rangle _{\omega
} $ would be to use instead of (\ref{big identity}), the more obvious
decomposition 
\begin{eqnarray}
&&\left\langle T_{\sigma }^{\alpha }\left( b_{A}\mathbf{1}_{A\setminus
K}\right) ,\mathbf{1}_{K_{\limfunc{in}}}b_{B}^{\ast }\right\rangle _{\omega
}=\left\langle T_{\sigma }^{\alpha }\left( b_{A}\mathbf{1}_{A\setminus
K}\right) ,\mathbf{1}_{K_{\limfunc{in}}}\left( b_{B}^{\ast }-\frac{1}{%
\left\vert K_{\limfunc{in}}\right\vert _{\omega }}\int_{K_{\limfunc{in}%
}}b_{B}^{\ast }d\omega \right) \right\rangle _{\omega }  \label{reach''} \\
&&\ \ \ \ \ \ \ \ \ \ \ \ \ \ \ \ \ \ \ \ +\left\langle T_{\sigma }^{\alpha
}\left( b_{A}\mathbf{1}_{A\setminus K}\right) ,\mathbf{1}_{K_{\limfunc{in}%
}}\left( \frac{1}{\left\vert K_{\limfunc{in}}\right\vert _{\omega }}\int_{K_{%
\limfunc{in}}}b_{B}^{\ast }d\omega \right) \right\rangle _{\omega }\ . 
\notag
\end{eqnarray}%
For the first term in (\ref{reach''}), we would like to apply the Energy
Lemma to obtain%
\begin{eqnarray}
&&\left\vert \left\langle T_{\sigma }^{\alpha }\left( b_{A}\mathbf{1}%
_{A\setminus K}\right) ,\mathbf{1}_{K_{\limfunc{in}}}\left( b_{B}^{\ast }-%
\frac{1}{\left\vert K_{\limfunc{in}}\right\vert _{\omega }}\int_{K_{\limfunc{%
in}}}b_{B}^{\ast }d\omega \right) \right\rangle _{\omega }\right\vert
\label{after mon''} \\
&\lesssim &\mathrm{P}_{\delta }^{\alpha }\mathsf{Q}^{\omega }\left( K_{%
\limfunc{in}},\mathbf{1}_{A\setminus K}\sigma \right) \left\Vert \mathbf{1}%
_{K_{\limfunc{in}}}\left( b_{B}^{\ast }-\frac{1}{\left\vert K_{\limfunc{in}%
}\right\vert _{\omega }}\int_{K_{\limfunc{in}}}b_{B}^{\ast }d\omega \right)
\right\Vert _{L^{2}\left( \omega \right) }\ ,  \notag
\end{eqnarray}%
using that the function $h_{K_{\limfunc{in}}}^{\ast }\equiv \mathbf{1}_{K_{%
\limfunc{in}}}\left( b_{B}^{\ast }-\frac{1}{\left\vert K_{\limfunc{in}%
}\right\vert _{\omega }}\int_{K_{\limfunc{in}}}b_{B}^{\ast }d\omega \right) $
has $\omega $-mean value zero and has support $K_{\limfunc{in}}$ that is 
\emph{strictly separated} from the support of $\left\vert b_{A}\right\vert 
\mathbf{1}_{A\setminus K}$. But a problem arises here since $K_{\limfunc{in}%
} $ is not in the dyadic grid $\mathcal{G}$, despite the fact that $K$
itself is. Indeed, the Energy Lemma requires the dual martingale support of $%
h_{K_{\limfunc{in}}}^{\ast }$ to be contained in $K_{\limfunc{in}}$, so that
we can take $\mathcal{H}$ in the Energy Lemma to be pseudoprojection onto $%
K_{\limfunc{in}}$. However, if $K^{\prime }$ is a child of $K$, then $%
\square _{K^{\prime }}^{\omega ,\mathbf{b}^{\ast }}h_{K_{\limfunc{in}%
}}^{\ast }$ could be nonzero, yet $K^{\prime }\not\subset K_{\limfunc{in}}$.
This is easily fixed in two steps as follows.

First, recall $\eta =\frac{1}{2}$ so that 
\begin{equation*}
K_{\limfunc{in}}=\left( 1-\eta \right) K=\frac{1}{2}K=\dbigcup \left\{
K^{\prime \prime }\in \mathfrak{C}^{\left( 2\right) }\left( K\right)
:\partial K^{\prime \prime }\cap \partial K=\emptyset \right\}
\end{equation*}%
is the union of the $2$ grandchildren $K^{\prime \prime }$ of $K$ whose
boundaries are disjoint from the boundary of $K$. Then 
\begin{equation*}
K_{\limfunc{out}}=\dbigcup \left\{ K^{\prime \prime }\in \mathfrak{C}%
^{\left( 2\right) }\left( K\right) :\partial K^{\prime \prime }\cap \partial
K\neq \emptyset \right\}
\end{equation*}%
is the union of the $2$ grandchildren $K^{\prime \prime }$ of $K$ whose
boundaries intersect the boundary of $K$. The only possible dyadic
subintervals $K^{\prime }$ of $K$ for which both $\square _{K^{\prime
}}^{\omega ,\mathbf{b}^{\ast }}h_{K_{\limfunc{in}}}^{\ast }\neq 0$ and $%
K^{\prime }\not\subset K_{\limfunc{in}}$ are the children of $K$. Enumerate
by $\left\{ K_{\ell }^{\prime \prime }\right\} _{\ell =1}^{2}$ the
grandchildren of $K$ whose boundaries are disjoint from the boundary of $K$.

Then second, instead of decomposing $\mathbf{1}_{K_{\limfunc{in}%
}}b_{B}^{\ast }$ as $h_{K_{\limfunc{in}}}^{\ast }$ plus $\mathbf{1}_{K_{%
\limfunc{in}}}\frac{1}{\left\vert K_{\limfunc{in}}\right\vert _{\omega }}%
\int_{K_{\limfunc{in}}}b_{B}^{\ast }d\omega $, we decompose $\mathbf{1}_{K_{%
\limfunc{in}}}b_{B}^{\ast }$ as%
\begin{equation*}
\mathbf{1}_{K_{\limfunc{in}}}b_{B}^{\ast }=\sum_{\ell =1}^{2}\mathbf{1}%
_{K_{\ell }^{\prime \prime }}\left( b_{B}^{\ast }-\frac{1}{\left\vert
K_{\ell }^{\prime \prime }\right\vert _{\omega }}\int_{K_{\ell }^{\prime
\prime }}b_{B}^{\ast }d\omega \right) +\sum_{\ell =1}^{2}\mathbf{1}_{K_{\ell
}^{\prime \prime }}\frac{1}{\left\vert K_{\ell }^{\prime \prime }\right\vert
_{\omega }}\int_{K_{\ell }^{\prime \prime }}b_{B}^{\ast }d\omega ,
\end{equation*}%
and then apply the Energy Lemma to the function%
\begin{equation*}
k_{K_{\limfunc{in}}}^{\ast }\equiv \sum_{\ell =1}^{2}\mathbf{1}_{K_{\ell
}^{\prime \prime }}\left( b_{B}^{\ast }-\frac{1}{\left\vert K_{\ell
}^{\prime \prime }\right\vert _{\omega }}\int_{K_{\ell }^{\prime \prime
}}b_{B}^{\ast }d\omega \right) \equiv k_{K_{\limfunc{in}}}^{\ast ,1}+k_{K_{%
\limfunc{in}}}^{\ast ,2},
\end{equation*}%
which does indeed satisfy $\square _{K^{\prime }}^{\omega ,\mathbf{b}^{\ast
}}k_{K_{\limfunc{in}}}^{\ast }=0$ unless $K^{\prime }$ is a dyadic
subinterval of $K$ that is contained in $K_{\limfunc{in}}$. (Furthermore, we
could even replace grandchildren by $m$-grandchildren in this argument in
order that $\square _{K^{\prime }}^{\omega ,\mathbf{b}^{\ast }}k_{K_{%
\limfunc{in}}}^{\ast }=0$ unless $K^{\prime }$ is a dyadic $m$-grandchild of 
$K$ that is contained in $K_{\limfunc{in}}$, but we will not need this.) If
we now use $k_{K_{\limfunc{in}}}^{\ast }$ instead of $h_{K_{\limfunc{in}%
}}^{\ast }$ in (\ref{reach''}) and (\ref{after mon''}), we obtain%
\begin{eqnarray}
\left\langle T_{\sigma }^{\alpha }\left( b_{A}\mathbf{1}_{A\setminus
K}\right) ,\mathbf{1}_{K_{\limfunc{in}}}b_{B}^{\ast }\right\rangle _{\omega
} &=&\left\langle T_{\sigma }^{\alpha }\left( b_{A}\mathbf{1}_{A\setminus
K}\right) ,k_{K_{\limfunc{in}}}^{\ast }\right\rangle _{\omega }
\label{reach'''} \\
&&+\left\langle T_{\sigma }^{\alpha }\left( b_{A}\mathbf{1}_{A\setminus
K}\right) ,\sum_{\ell =1}^{2}\mathbf{1}_{K_{\ell }^{\prime \prime }}\left( 
\frac{1}{\left\vert K_{\ell }^{\prime \prime }\right\vert _{\omega }}%
\int_{K_{\ell }^{\prime \prime }}b_{B}^{\ast }d\omega \right) \right\rangle
_{\omega }\ ,  \notag
\end{eqnarray}%
and%
\begin{eqnarray}
\left\vert \left\langle T_{\sigma }^{\alpha }\left( b_{A}\mathbf{1}%
_{A\setminus K}\right) ,k_{K_{\limfunc{in}}}^{\ast }\right\rangle _{\omega
}\right\vert &\leq &\left\vert \left\langle T_{\sigma }^{\alpha }\left( b_{A}%
\mathbf{1}_{A\setminus K}\right) ,k_{K_{\limfunc{in}}}^{\ast
,1}\right\rangle _{\omega }\right\vert +\left\vert \left\langle T_{\sigma
}^{\alpha }\left( b_{A}\mathbf{1}_{A\setminus K}\right) ,k_{K_{\limfunc{in}%
}}^{\ast ,2}\right\rangle _{\omega }\right\vert  \label{after mon'''} \\
&\leq &C_{\eta }\left\{ \mathrm{P}_{\delta }^{\alpha }\mathsf{Q}^{\omega
}\left( K_{\limfunc{in}}^{\limfunc{left}},\mathbf{1}_{A\setminus K}\sigma
\right) +\mathrm{P}_{\delta }^{\alpha }\mathsf{Q}^{\omega }\left( K_{%
\limfunc{in}}^{\limfunc{right}},\mathbf{1}_{A\setminus K}\sigma \right)
\right\} \left\Vert k_{K_{\limfunc{in}}}^{\ast }\right\Vert _{L^{2}\left(
\omega \right) }\ ,  \notag
\end{eqnarray}%
where the constant $C_{\eta }$ depends on the constant $C_{\gamma }$ in the
statement of the Monotonicity Lemma with $\gamma =\frac{1}{1-\eta }$ since $%
\frac{1}{1-\eta }K_{\limfunc{in}}\cap \left( A\setminus K\right) =\emptyset $%
, and where we have written $\left\{ K_{\ell }^{\prime \prime }\right\}
_{\ell =1}^{2}=\left\{ K_{\limfunc{in}}^{\limfunc{left}},K_{\limfunc{in}}^{%
\limfunc{right}}\right\} $ with $K_{\limfunc{in}}^{\limfunc{left}}$ and $K_{%
\limfunc{in}}^{\limfunc{right}}$ denoting the left hand child and right hand
child of $K_{\limfunc{in}}$ respectively.

\begin{conclusion}
Thus we see that $\mathsf{P}_{\mathcal{H}}^{\omega ,\mathbf{b}^{\ast }}$ and 
$\mathsf{Q}_{\mathcal{H}}^{\omega ,\mathbf{b}^{\ast }}$ in the Energy Lemma
can be taken to be pseudoprojection onto $K_{\limfunc{in}}$, i.e. $\mathsf{P}%
_{K_{\limfunc{in}}}^{\omega ,\mathbf{b}^{\ast }}=\sum_{J\in \mathcal{G}:\
J\subset K_{\limfunc{in}}}\square _{J}^{\omega ,\mathbf{b}^{\ast }}$ and $%
\mathsf{Q}_{K_{\limfunc{in}}}^{\omega ,\mathbf{b}^{\ast }}=\sum_{J\in 
\mathcal{G}:\ J\subset K_{\limfunc{in}}}\bigtriangleup _{J}^{\omega ,\mathbf{%
b}^{\ast }}$, and we will see below that the intervals $K_{\limfunc{in}}$
that arise in subsequent arguments will be pairwise disjoint. Furthermore,
the energy condition will be used to control these full pseudoprojections $%
\mathsf{P}_{K_{\limfunc{in}}}^{\omega ,\mathbf{b}^{\ast }}$ when taken over
pairwise disjoint decompositions of intervals by subintervals of the form $%
K_{\limfunc{in}}$.
\end{conclusion}

However, the second line of (\ref{reach'''}) remains problematic, and this
is our point of departure for beginning the proof of Lemma \ref{preiterate},
in which we exploit the original testing functions $b_{K_{\ell }^{\prime
\prime }}^{\ast ,\limfunc{orig}}$ in identity (\ref{big identity}) in order
to handle the second line of (\ref{reach'''}).

\begin{proof}[Proof of Lemma \protect\ref{preiterate}]
We begin by rewriting identity (\ref{big identity}) in the form%
\begin{eqnarray*}
\boldsymbol{A} &=&\boldsymbol{B}+\boldsymbol{C}\text{ where} \\
\boldsymbol{A} &\equiv &\left\langle T_{\sigma }^{\alpha }\left( b_{A}%
\mathbf{1}_{A\setminus K}\right) ,\mathbf{1}_{K_{\limfunc{in}}}b_{B}^{\ast
}\right\rangle _{\omega }+\sum_{\ell =1}^{2}\left( \frac{\frac{1}{\left\vert
K_{\ell }^{\prime \prime }\right\vert _{\omega }}\int_{K_{\ell }^{\prime
\prime }}b_{B}^{\ast }d\omega }{\frac{1}{\left\vert K_{\ell }^{\prime \prime
}\right\vert _{\omega }}\int_{K_{\ell }^{\prime \prime }}b_{K_{\ell
}^{\prime \prime }}^{\ast ,\limfunc{orig}}d\omega }\right) \left\langle
T_{\sigma }^{\alpha }b_{A}\mathbf{1}_{K_{\limfunc{out}}},b_{K_{\ell
}^{\prime \prime }}^{\ast ,\limfunc{orig}}\right\rangle _{\omega }\ , \\
\boldsymbol{B} &\equiv &\left\langle T_{\sigma }^{\alpha }\left( b_{A}%
\mathbf{1}_{A\setminus K}\right) ,\mathbf{1}_{K_{\limfunc{in}}}b_{B}^{\ast
}\right\rangle _{\omega }-\sum_{\ell =1}^{2}\left( \frac{\frac{1}{\left\vert
K_{\ell }^{\prime \prime }\right\vert _{\omega }}\int_{K_{\ell }^{\prime
\prime }}b_{B}^{\ast }d\omega }{\frac{1}{\left\vert K_{\ell }^{\prime \prime
}\right\vert _{\omega }}\int_{K_{\ell }^{\prime \prime }}b_{K_{\ell
}^{\prime \prime }}^{\ast ,\limfunc{orig}}d\omega }\right) \left\langle
T_{\sigma }^{\alpha }\left( b_{A}\mathbf{1}_{A\setminus K}\right)
,b_{K_{\ell }^{\prime \prime }}^{\ast ,\limfunc{orig}}\right\rangle _{\omega
}\ , \\
\boldsymbol{C} &\equiv &\sum_{\ell =1}^{2}\left( \frac{\frac{1}{\left\vert
K_{\ell }^{\prime \prime }\right\vert _{\omega }}\int_{K_{\ell }^{\prime
\prime }}b_{B}^{\ast }d\omega }{\frac{1}{\left\vert K_{\ell }^{\prime \prime
}\right\vert _{\omega }}\int_{K_{\ell }^{\prime \prime }}b_{K_{\ell
}^{\prime \prime }}^{\ast ,\limfunc{orig}}d\omega }\right) \left\{
\left\langle T_{\sigma }^{\alpha }b_{A},b_{K_{\ell }^{\prime \prime }}^{\ast
,\limfunc{orig}}\right\rangle _{\omega }-\left\langle b_{A}\mathbf{1}_{K_{%
\limfunc{in}}},T_{\omega }^{\alpha ,\ast }b_{K_{\ell }^{\prime \prime
}}^{\ast ,\limfunc{orig}}\right\rangle _{\sigma }\right\} \ .
\end{eqnarray*}%
From the discussion above, we recall the identity (\ref{reach'''}) and the
estimate (\ref{after mon'''}). We also have the analogous identity and
estimate with $b_{K_{\ell }^{\prime \prime }}^{\ast ,\limfunc{orig}}$ in
place of $\mathbf{1}_{K_{\limfunc{in}}}b_{B}^{\ast }$:%
\begin{eqnarray}
\left\langle T_{\sigma }^{\alpha }\left( b_{A}\mathbf{1}_{A\setminus
K}\right) ,b_{K_{\ell }^{\prime \prime }}^{\ast ,\limfunc{orig}%
}\right\rangle _{\omega } &=&\left\langle T_{\sigma }^{\alpha }\left( b_{A}%
\mathbf{1}_{A\setminus K}\right) ,\mathbf{1}_{K_{\ell }^{\prime \prime
}}\left( b_{K_{\ell }^{\prime \prime }}^{\ast ,\limfunc{orig}}-\frac{1}{%
\left\vert K_{\ell }^{\prime \prime }\right\vert _{\omega }}\int_{K_{\ell
}^{\prime \prime }}b_{K_{\ell }^{\prime \prime }}^{\ast ,\limfunc{orig}%
}d\omega \right) \right\rangle _{\omega }  \label{reach''''} \\
&&+\left\langle T_{\sigma }^{\alpha }\left( b_{A}\mathbf{1}_{A\setminus
K}\right) ,\mathbf{1}_{K_{\ell }^{\prime \prime }}\left( \frac{1}{\left\vert
K_{\ell }^{\prime \prime }\right\vert _{\omega }}\int_{K_{\ell }^{\prime
\prime }}b_{K_{\ell }^{\prime \prime }}^{\ast ,\limfunc{orig}}d\omega
\right) \right\rangle _{\omega }\ ,  \notag
\end{eqnarray}%
and%
\begin{eqnarray}
&&\left\vert \left\langle T_{\sigma }^{\alpha }\left( b_{A}\mathbf{1}%
_{A\setminus K}\right) ,\mathbf{1}_{K_{\ell }^{\prime \prime }}\left(
b_{K_{\ell }^{\prime \prime }}^{\ast ,\limfunc{orig}}-\frac{1}{\left\vert
K_{\ell }^{\prime \prime }\right\vert _{\omega }}\int_{K_{\ell }^{\prime
\prime }}b_{K_{\ell }^{\prime \prime }}^{\ast ,\limfunc{orig}}d\omega
\right) \right\rangle _{\omega }\right\vert  \label{after mon''''} \\
&\lesssim &\mathrm{P}_{\delta }^{\alpha }\mathsf{Q}^{\omega }\left( K_{\ell
}^{\prime \prime },\mathbf{1}_{A\setminus K}\sigma \right) \left\Vert 
\mathbf{1}_{K_{\ell }^{\prime \prime }}\left( b_{K_{\ell }^{\prime \prime
}}^{\ast ,\limfunc{orig}}-\frac{1}{\left\vert K_{\ell }^{\prime \prime
}\right\vert _{\omega }}\int_{K_{\ell }^{\prime \prime }}b_{K_{\ell
}^{\prime \prime }}^{\ast ,\limfunc{orig}}d\omega \right) \right\Vert
_{L^{2}\left( \omega \right) }\ ,  \notag
\end{eqnarray}%
for $1\leq \ell \leq 2$, where the implied constants depend on $L^{\infty }$
norms of testing functions and the constant in the Energy Lemma. We will
typically suppress dependence on $\eta =\frac{1}{2}$ from now on since there
are no other values of $\eta $ used below. Thus we have, using that $K_{\ell
}^{\prime \prime }$ is close to $K_{\limfunc{in}}$ in both scale and
position, 
\begin{eqnarray}
&&\left\vert \left\langle T_{\sigma }^{\alpha }\left( b_{A}\mathbf{1}%
_{A\setminus K}\right) ,\mathbf{1}_{K_{\limfunc{in}}}b_{B}^{\ast
}\right\rangle _{\omega }-\sum_{\ell =1}^{2}\left( \frac{1}{\left\vert
K_{\ell }^{\prime \prime }\right\vert _{\omega }}\int_{K_{\ell }^{\prime
\prime }}b_{B}^{\ast }d\omega \right) \left\langle T_{\sigma }^{\alpha
}\left( b_{A}\mathbf{1}_{A\setminus K}\right) ,\mathbf{1}_{K_{\ell }^{\prime
\prime }}\right\rangle _{\omega }\right\vert  \label{scale and position} \\
&\lesssim &\left\{ \mathrm{P}_{\delta }^{\alpha }\mathsf{Q}^{\omega }\left(
K_{\limfunc{in}}^{\limfunc{left}},\mathbf{1}_{A\setminus K}\sigma \right) +%
\mathrm{P}_{\delta }^{\alpha }\mathsf{Q}^{\omega }\left( K_{\limfunc{in}}^{%
\limfunc{right}},\mathbf{1}_{A\setminus K}\sigma \right) \right\} \left\Vert 
\mathbf{1}_{K_{\limfunc{in}}}\left( b_{B}^{\ast }-\frac{1}{\left\vert K_{%
\limfunc{in}}\right\vert _{\omega }}\int_{K_{\limfunc{in}}}b_{B}^{\ast
}d\omega \right) \right\Vert _{L^{2}\left( \omega \right) }  \notag \\
&\lesssim &\left\{ \mathrm{P}_{\delta }^{\alpha }\mathsf{Q}^{\omega }\left(
K_{\limfunc{in}}^{\limfunc{left}},\mathbf{1}_{A\setminus K}\sigma \right) +%
\mathrm{P}_{\delta }^{\alpha }\mathsf{Q}^{\omega }\left( K_{\limfunc{in}}^{%
\limfunc{right}},\mathbf{1}_{A\setminus K}\sigma \right) \right\} \sqrt{%
\left\vert K_{\limfunc{in}}\right\vert _{\omega }}\ .  \notag
\end{eqnarray}%
Then we obtain, upon applying (\ref{star}) with the function $\Psi
_{J}^{\ell }$ equal to 
\begin{equation*}
\left( \frac{\frac{1}{\left\vert K_{\ell }^{\prime \prime }\right\vert
_{\omega }}\int_{K_{\ell }^{\prime \prime }}b_{B}^{\ast }d\omega }{\frac{1}{%
\left\vert K_{\ell }^{\prime \prime }\right\vert _{\omega }}\int_{K_{\ell
}^{\prime \prime }}b_{K_{\ell }^{\prime \prime }}^{\ast ,\limfunc{orig}%
}d\omega }\right) b_{K_{\ell }^{\prime \prime }}^{\ast ,\limfunc{orig}%
}-\left( \frac{1}{\left\vert K_{\ell }^{\prime \prime }\right\vert _{\omega }%
}\int_{K_{\ell }^{\prime \prime }}b_{B}^{\ast }d\omega \right) \mathbf{1}%
_{K_{\ell }^{\prime \prime }}
\end{equation*}%
for $\ell =1,2$, and also using (\ref{scale and position}), that%
\begin{eqnarray*}
&&\left\vert \boldsymbol{B}\right\vert =\left\vert \left\langle T_{\sigma
}^{\alpha }\left( b_{A}\mathbf{1}_{A\setminus K}\right) ,\mathbf{1}_{K_{%
\limfunc{in}}}b_{B}^{\ast }\right\rangle _{\omega }-\sum_{\ell =1}^{2}\left( 
\frac{\frac{1}{\left\vert K_{\ell }^{\prime \prime }\right\vert _{\omega }}%
\int_{K_{\ell }^{\prime \prime }}b_{B}^{\ast }d\omega }{\frac{1}{\left\vert
K_{\ell }^{\prime \prime }\right\vert _{\omega }}\int_{K_{\ell }^{\prime
\prime }}b_{K_{\ell }^{\prime \prime }}^{\ast ,\limfunc{orig}}d\omega }%
\right) \left\langle T_{\sigma }^{\alpha }\left( b_{A}\mathbf{1}_{A\setminus
K}\right) ,b_{K_{\ell }^{\prime \prime }}^{\ast ,\limfunc{orig}%
}\right\rangle _{\omega }\right\vert \\
&=&\left\vert \left\langle T_{\sigma }^{\alpha }\left( b_{A}\mathbf{1}%
_{A\setminus K}\right) ,\mathbf{1}_{K_{\limfunc{in}}}b_{B}^{\ast
}\right\rangle _{\omega }-\sum_{\ell =1}^{2}\left( \frac{1}{\left\vert
K_{\ell }^{\prime \prime }\right\vert _{\omega }}\int_{K_{\ell }^{\prime
\prime }}b_{B}^{\ast }d\omega \right) \left\langle T_{\sigma }^{\alpha
}\left( b_{A}\mathbf{1}_{A\setminus K}\right) ,\mathbf{1}_{K_{\ell }^{\prime
\prime }}\right\rangle _{\omega }\right\vert \\
&&+O\left\{ \sum_{\ell =1}^{2}\left( \frac{\mathrm{P}^{\alpha }\left(
K_{\ell }^{\prime \prime },\mathbf{1}_{A\setminus K}\sigma \right) }{%
\left\vert K_{\ell }^{\prime \prime }\right\vert }\left\Vert \mathsf{Q}%
_{K_{\ell }^{\prime \prime }}^{\omega ,\mathbf{b}^{\ast }}x\right\Vert
_{L^{2}\left( \omega \right) }^{\spadesuit }+\frac{\mathrm{P}_{1+\delta
}^{\alpha }\left( K_{\ell }^{\prime \prime },\mathbf{1}_{A\setminus K}\sigma
\right) }{\left\vert K_{\ell }^{\prime \prime }\right\vert }\left\Vert
x-m_{K_{\ell }^{\prime \prime }}\right\Vert _{L^{2}\left( \mathbf{1}%
_{K_{\ell }^{\prime \prime }}\omega \right) }\right) \sqrt{\left\vert K_{%
\limfunc{in}}\right\vert _{\omega }}\right\} \\
&\lesssim &\left\{ \mathrm{P}_{\delta }^{\alpha }\mathsf{Q}^{\omega }\left(
K_{\limfunc{in}}^{\limfunc{left}},\mathbf{1}_{A\setminus K}\sigma \right) +%
\mathrm{P}_{\delta }^{\alpha }\mathsf{Q}^{\omega }\left( K_{\limfunc{in}}^{%
\limfunc{right}},\mathbf{1}_{A\setminus K}\sigma \right) \right\} \sqrt{%
\left\vert K_{\limfunc{in}}\right\vert _{\omega }}\ ,
\end{eqnarray*}%
since by the triangle inequality 
\begin{eqnarray*}
&&\left\Vert \Psi _{J}^{\ell }\right\Vert _{L^{2}\left( \omega \right)
}=\left\Vert \left( \frac{\frac{1}{\left\vert K_{\ell }^{\prime \prime
}\right\vert _{\omega }}\int_{K_{\ell }^{\prime \prime }}b_{B}^{\ast
}d\omega }{\frac{1}{\left\vert K_{\ell }^{\prime \prime }\right\vert
_{\omega }}\int_{K_{\ell }^{\prime \prime }}b_{K_{\ell }^{\prime \prime
}}^{\ast ,\limfunc{orig}}d\omega }\right) b_{K_{\ell }^{\prime \prime
}}^{\ast ,\limfunc{orig}}-\left( \frac{1}{\left\vert K_{\ell }^{\prime
\prime }\right\vert _{\omega }}\int_{K_{\ell }^{\prime \prime }}b_{B}^{\ast
}d\omega \right) \mathbf{1}_{K_{\ell }^{\prime \prime }}\right\Vert
_{L^{2}\left( \omega \right) } \\
&\lesssim &\left\vert \frac{\int_{K_{\ell }^{\prime \prime }}b_{B}^{\ast
}d\omega }{\int_{K_{\ell }^{\prime \prime }}b_{K_{\ell }^{\prime \prime
}}^{\ast ,\limfunc{orig}}d\omega }\right\vert \sqrt{\left\vert K_{\ell
}^{\prime \prime }\right\vert _{\omega }}+\sqrt{\left\vert K_{\ell }^{\prime
\prime }\right\vert _{\omega }}\lesssim \sqrt{\left\vert K_{\ell }^{\prime
\prime }\right\vert _{\omega }}\leq \sqrt{\left\vert K_{\limfunc{in}%
}\right\vert _{\omega }},\ \ \ \ \ 1\leq \ell \leq 2.
\end{eqnarray*}

Finally, turning our attention to term $\boldsymbol{C}$, the reason for
using the identity (\ref{big identity}) now becomes clear - namely the terms 
$\left\langle T_{\sigma }^{\alpha }\left( b_{A}\mathbf{1}_{A\setminus
K}\right) ,b_{K_{\ell }^{\prime \prime }}^{\ast ,\limfunc{orig}%
}\right\rangle _{\omega }$, for which multiples are subtracted on the left
side above, involves the original testing functions $b_{K_{\ell }^{\prime
\prime }}^{\ast ,\limfunc{orig}}$ for which we have the\ full testing
condition. Thus we have 
\begin{eqnarray*}
&&\left\vert \sum_{\ell =1}^{2}\left( \frac{\frac{1}{\left\vert K_{\ell
}^{\prime \prime }\right\vert _{\omega }}\int_{K_{\ell }^{\prime \prime
}}b_{B}^{\ast }d\omega }{\frac{1}{\left\vert K_{\ell }^{\prime \prime
}\right\vert _{\omega }}\int_{K_{\ell }^{\prime \prime }}b_{K_{\ell
}^{\prime \prime }}^{\ast ,\limfunc{orig}}d\omega }\right) \left\langle
T_{\sigma }^{\alpha }b_{A},b_{K_{\ell }^{\prime \prime }}^{\ast ,\limfunc{%
orig}}\right\rangle _{\omega }\right\vert \\
&\lesssim &\sum_{\ell =1}^{2}\sqrt{\int_{K_{\ell }^{\prime \prime
}}\left\vert T_{\sigma }^{\alpha }b_{A}\right\vert ^{2}d\omega }\sqrt{%
\int_{K_{\ell }^{\prime \prime }}\left\vert b_{K_{\ell }^{\prime \prime
}}^{\ast ,\limfunc{orig}}\right\vert ^{2}d\omega }\lesssim \sqrt{\int_{K_{%
\limfunc{in}}}\left\vert T_{\sigma }^{\alpha }b_{A}\right\vert ^{2}d\omega }%
\sqrt{\left\vert K_{\limfunc{in}}\right\vert _{\omega }},
\end{eqnarray*}%
and similarly%
\begin{eqnarray*}
&&\left\vert \sum_{\ell =1}^{2}\left( \frac{\frac{1}{\left\vert K_{\ell
}^{\prime \prime }\right\vert _{\omega }}\int_{K_{\ell }^{\prime \prime
}}b_{B}^{\ast }d\omega }{\frac{1}{\left\vert K_{\ell }^{\prime \prime
}\right\vert _{\omega }}\int_{K_{\ell }^{\prime \prime }}b_{K_{\ell
}^{\prime \prime }}^{\ast ,\limfunc{orig}}d\omega }\right) \left\langle b_{A}%
\mathbf{1}_{K_{\limfunc{in}}},T_{\omega }^{\alpha ,\ast }b_{K_{\ell
}^{\prime \prime }}^{\ast ,\limfunc{orig}}\right\rangle _{\sigma
}\right\vert \lesssim \sqrt{\int_{K_{\limfunc{in}}}\left\vert
b_{A}\right\vert ^{2}d\sigma }\sqrt{\int_{K_{\limfunc{in}}}\left\vert
T_{\omega }^{\alpha ,\ast }\sum_{\ell =1}^{2}b_{K_{\ell }^{\prime \prime
}}^{\ast ,\limfunc{orig}}\right\vert ^{2}d\sigma } \\
&\lesssim &\left( \mathfrak{T}_{T^{\alpha ,\ast }}+\sqrt{\mathfrak{A}%
_{2}^{\alpha }}\right) \sqrt{\int_{K_{\limfunc{in}}}\left\vert
b_{A}\right\vert ^{2}d\sigma }\sum_{\ell =1}^{2}\sqrt{\int \left\vert
b_{K_{\ell }^{\prime \prime }}^{\ast ,\limfunc{orig}}\right\vert ^{2}d\omega 
}\lesssim \left( \mathfrak{T}_{T^{\alpha ,\ast }}+\sqrt{\mathfrak{A}%
_{2}^{\alpha }}\right) \sqrt{\int_{K_{\limfunc{in}}}\left\vert
b_{A}\right\vert ^{2}d\sigma }\sqrt{\left\vert K_{\limfunc{in}}\right\vert
_{\omega }},
\end{eqnarray*}%
which together prove%
\begin{eqnarray*}
\left\vert \boldsymbol{C}\right\vert &=&\left\vert \sum_{\ell =1}^{2}\left( 
\frac{\frac{1}{\left\vert K_{\ell }^{\prime \prime }\right\vert _{\omega }}%
\int_{K_{\ell }^{\prime \prime }}b_{B}^{\ast }d\omega }{\frac{1}{\left\vert
K_{\ell }^{\prime \prime }\right\vert _{\omega }}\int_{K_{\ell }^{\prime
\prime }}b_{K_{\ell }^{\prime \prime }}^{\ast ,\limfunc{orig}}d\omega }%
\right) \left\{ \left\langle T_{\sigma }^{\alpha }b_{A},b_{K_{\ell }^{\prime
\prime }}^{\ast ,\limfunc{orig}}\right\rangle _{\omega }-\left\langle b_{A}%
\mathbf{1}_{K_{\limfunc{in}}},T_{\omega }^{\alpha ,\ast }b_{K_{\ell
}^{\prime \prime }}^{\ast ,\limfunc{orig}}\right\rangle _{\sigma }\right\}
\right\vert \\
&\lesssim &\left( \sqrt{\int_{K_{\limfunc{in}}}\left\vert T_{\sigma
}^{\alpha }b_{A}\right\vert ^{2}d\omega }+\left( \mathfrak{T}_{T^{\alpha
,\ast }}+\sqrt{\mathfrak{A}_{2}^{\alpha }}\right) \sqrt{\left\vert K_{%
\limfunc{in}}\right\vert _{\sigma }}\right) \sqrt{\left\vert K_{\limfunc{in}%
}\right\vert _{\omega }}\ .
\end{eqnarray*}%
This completes the proof of Lemma \ref{preiterate}, and hence also that of (%
\ref{second term}) since $\left\vert \frac{\frac{1}{\left\vert K_{\ell
}^{\prime \prime }\right\vert _{\omega }}\int_{K_{\ell }^{\prime \prime
}}b_{B}^{\ast }d\omega }{\frac{1}{\left\vert K_{\ell }^{\prime \prime
}\right\vert _{\omega }}\int_{K_{\ell }^{\prime \prime }}b_{K_{\ell
}^{\prime \prime }}^{\ast ,\limfunc{orig}}d\omega }\right\vert \lesssim 1$.
\end{proof}

The term $\left\{ K_{\limfunc{out}},K_{\limfunc{in}}^{\ell }\right\} ^{%
\limfunc{orig}}$ will be handled below by relatively crude estimates. If we
write%
\begin{eqnarray*}
&&\left\{ K,K\right\} =\left\{ A,K_{\limfunc{in}}\right\} -\left\{
A\setminus K,K_{\limfunc{in}}\right\} +\left\{ K_{\limfunc{out}},K_{\limfunc{%
out}}\right\} +\left\{ K_{\limfunc{in}},K_{\limfunc{out}}\right\} \\
&=&\left\{ A,K_{\limfunc{in}}\right\} -\left( \left\{ A\setminus K,K_{%
\limfunc{in}}\right\} +\sum_{\ell =1}^{2}\left( \frac{\frac{1}{\left\vert
K_{\ell }^{\prime \prime }\right\vert _{\omega }}\int_{K_{\ell }^{\prime
\prime }}b_{B}^{\ast }d\omega }{\frac{1}{\left\vert K_{\ell }^{\prime \prime
}\right\vert _{\omega }}\int_{K_{\ell }^{\prime \prime }}b_{K_{\ell
}^{\prime \prime }}^{\ast ,\limfunc{orig}}d\omega }\right) \left\{ K_{%
\limfunc{out}},K_{\limfunc{in}}^{\ell }\right\} ^{\limfunc{orig}}\right) \\
&&+\sum_{\ell =1}^{2}\left( \frac{\frac{1}{\left\vert K_{\ell }^{\prime
\prime }\right\vert _{\omega }}\int_{K_{\ell }^{\prime \prime }}b_{B}^{\ast
}d\omega }{\frac{1}{\left\vert K_{\ell }^{\prime \prime }\right\vert
_{\omega }}\int_{K_{\ell }^{\prime \prime }}b_{K_{\ell }^{\prime \prime
}}^{\ast ,\limfunc{orig}}d\omega }\right) \left\{ K_{\limfunc{out}},K_{%
\limfunc{in}}^{\ell }\right\} ^{\limfunc{orig}}+\left\{ K_{\limfunc{out}},K_{%
\limfunc{out}}\right\} +\left\{ K_{\limfunc{in}},K_{\limfunc{out}}\right\} ,
\end{eqnarray*}%
then, using (\ref{AKin}) and Lemma \ref{preiterate}, we see that we have
reduced control of $\left\{ K,K\right\} $ to control of%
\begin{equation*}
\sum_{\ell =1}^{2}\left( \frac{\frac{1}{\left\vert K_{\ell }^{\prime \prime
}\right\vert _{\omega }}\int_{K_{\ell }^{\prime \prime }}b_{B}^{\ast
}d\omega }{\frac{1}{\left\vert K_{\ell }^{\prime \prime }\right\vert
_{\omega }}\int_{K_{\ell }^{\prime \prime }}b_{K_{\ell }^{\prime \prime
}}^{\ast ,\limfunc{orig}}d\omega }\right) \left\{ K_{\limfunc{out}},K_{%
\limfunc{in}}^{\ell }\right\} ^{\limfunc{orig}}+\left\{ K_{\limfunc{out}},K_{%
\limfunc{out}}\right\} +\left\{ K_{\limfunc{in}},K_{\limfunc{out}}\right\} .
\end{equation*}

Altogether then, using the above estimates, we have proved%
\begin{eqnarray}
&&  \label{KK bound} \\
&&\left\vert \left\{ K,K\right\} -\left\{ K_{\limfunc{out}},K_{\limfunc{out}%
}\right\} -\left\{ K_{\limfunc{in}},K_{\limfunc{out}}\right\} -\sum_{\ell
=1}^{2}\left( \frac{\frac{1}{\left\vert K_{\ell }^{\prime \prime
}\right\vert _{\omega }}\int_{K_{\ell }^{\prime \prime }}b_{B}^{\ast
}d\omega }{\frac{1}{\left\vert K_{\ell }^{\prime \prime }\right\vert
_{\omega }}\int_{K_{\ell }^{\prime \prime }}b_{K_{\ell }^{\prime \prime
}}^{\ast ,\limfunc{orig}}d\omega }\right) \left\{ K_{\limfunc{out}},K_{%
\limfunc{in}}^{\ell }\right\} ^{\limfunc{orig}}\right\vert  \notag \\
&\lesssim &\left\Vert \mathbf{1}_{K_{\limfunc{in}}}T_{\sigma }^{\alpha
}b_{A}\right\Vert _{L^{2}\left( \omega \right) }\sqrt{\left\vert K_{\limfunc{%
in}}\right\vert _{\omega }}+C\left\{ \mathrm{P}_{\delta }^{\alpha }\mathsf{Q}%
^{\omega }\left( K_{\limfunc{in}}^{\limfunc{left}},\mathbf{1}_{A\setminus
K}\sigma \right) +\mathrm{P}_{\delta }^{\alpha }\mathsf{Q}^{\omega }\left(
K_{\limfunc{in}}^{\limfunc{right}},\mathbf{1}_{A\setminus K}\sigma \right)
\right\} \sqrt{\left\vert K_{\limfunc{in}}\right\vert _{\omega }}  \notag \\
&\lesssim &C\left( \mathfrak{T}_{T^{\alpha ,\ast }}+\mathfrak{A}_{2}^{\alpha
}\right) \sqrt{\left\vert K_{\limfunc{in}}\right\vert _{\sigma }}\sqrt{%
\left\vert K_{\limfunc{in}}\right\vert _{\omega }}.  \notag
\end{eqnarray}%
We emphasize that this bound did not use any special information regarding $%
K $ being in the coronas $\mathcal{C}_{A},\mathcal{C}_{B}$ or not, and thus
holds for any interval $K$ contained in $A\cap B$. For clarity of notation
we define%
\begin{eqnarray}
\Phi ^{A,B}\left( K_{\limfunc{in}}\right) &\equiv &\left\Vert \mathbf{1}_{K_{%
\limfunc{in}}}T_{\sigma }^{\alpha }b_{A}\right\Vert _{L^{2}\left( \omega
\right) }\sqrt{\left\vert K_{\limfunc{in}}\right\vert _{\omega }}
\label{def PHI} \\
&&+C\left\{ \mathrm{P}_{\delta }^{\alpha }\mathsf{Q}^{\omega }\left( K_{%
\limfunc{in}}^{\limfunc{left}},\mathbf{1}_{A\setminus K}\sigma \right) +%
\mathrm{P}_{\delta }^{\alpha }\mathsf{Q}^{\omega }\left( K_{\limfunc{in}}^{%
\limfunc{right}},\mathbf{1}_{A\setminus K}\sigma \right) \right\} \sqrt{%
\left\vert K_{\limfunc{in}}\right\vert _{\omega }}  \notag \\
&&+C\left( \mathfrak{T}_{T^{\alpha ,\ast }}+\mathfrak{A}_{2}^{\alpha
}\right) \sqrt{\left\vert K_{\limfunc{in}}\right\vert _{\sigma }}\sqrt{%
\left\vert K_{\limfunc{in}}\right\vert _{\omega }},  \notag
\end{eqnarray}%
where $\Phi ^{A,B}$ should not be confused with the notation $\Phi ^{\alpha
} $ introduced for the Monotonicity Lemma, and we also define the constants%
\begin{equation*}
A_{K_{\limfunc{in}}^{\ell }}\equiv \frac{\frac{1}{\left\vert K_{\ell
}^{\prime \prime }\right\vert _{\omega }}\int_{K_{\ell }^{\prime \prime
}}b_{B}^{\ast }d\omega }{\frac{1}{\left\vert K_{\ell }^{\prime \prime
}\right\vert _{\omega }}\int_{K_{\ell }^{\prime \prime }}b_{K_{\ell
}^{\prime \prime }}^{\ast ,\limfunc{orig}}d\omega },\ \ \ \ \ \ell \in
\left\{ 1,2\right\} ,
\end{equation*}%
so that we can rewrite (\ref{KK bound}) as%
\begin{equation}
\left\vert \left\{ K,K\right\} -\left\{ K_{\limfunc{out}},K_{\limfunc{out}%
}\right\} -\left\{ K_{\limfunc{in}},K_{\limfunc{out}}\right\} -\sum_{\ell
=1}^{2}A_{K_{\limfunc{in}}^{\ell }}\left\{ K_{\limfunc{out}},K_{\limfunc{in}%
}^{\ell }\right\} ^{\limfunc{orig}}\right\vert \lesssim \Phi ^{A,B}\left( K_{%
\limfunc{in}}\right) ,\ \ \ \ \ K\in \mathcal{G},K\subset A\cap B\ .
\label{KK bound rewrite}
\end{equation}%
We can further simplify notation by defining 
\begin{equation}
\left\{ K_{\limfunc{out}},K_{\limfunc{in}}\right\} ^{\limfunc{orig}}\equiv
\sum_{\ell =1}^{2}A_{K_{\limfunc{in}}^{\ell }}\left\{ K_{\limfunc{out}},K_{%
\limfunc{in}}^{\ell }\right\} ^{\limfunc{orig}},  \label{def orig}
\end{equation}%
which we will often use below. At this point, as we will see below, the only
problematic inner product subtracted from $\left\{ K,K\right\} $ on the left
hand side of (\ref{KK bound rewrite}) is $\left\{ K_{\limfunc{out}},K_{%
\limfunc{out}}\right\} $, and we will handle this by iterating (\ref{KK
bound rewrite}) finitely many times and then appealing to a final
probability argument starting in (\ref{future prob}) below.

\subsection{A finite iteration and final random surgery\label{Subsection
iteration}}

For $K$ an interval, we write $K_{\limfunc{out}}=K_{\limfunc{left}}\cup K_{%
\limfunc{right}}$ where $K_{\limfunc{left}}$ and $K_{\limfunc{right}}$ are
the two small subintervals on the left and right hand sides of $K$
respectively, and then we have 
\begin{equation*}
\left\{ K_{\limfunc{out}},K_{\limfunc{out}}\right\} =\left\{ K_{\limfunc{left%
}},K_{\limfunc{left}}\right\} +\left\{ K_{\limfunc{right}},K_{\limfunc{right}%
}\right\} +\left\{ K_{\limfunc{left}},K_{\limfunc{right}}\right\} +\left\{
K_{\limfunc{right}},K_{\limfunc{left}}\right\} ,
\end{equation*}%
so that (\ref{KK bound rewrite}) can be written using (\ref{def orig}) as%
\begin{eqnarray}
\left\{ K,K\right\} &=&\left\{ K_{\limfunc{out}},K_{\limfunc{out}}\right\}
+\left\{ K_{\limfunc{in}},K_{\limfunc{out}}\right\} +\left\{ K_{\limfunc{out}%
},K_{\limfunc{in}}\right\} ^{\limfunc{orig}}+O\left[ \Phi ^{A,B}\left( K_{%
\limfunc{in}}\right) \right]  \label{KK bound again} \\
&=&\left\{ K_{\limfunc{left}},K_{\limfunc{left}}\right\} +\left\{ K_{%
\limfunc{right}},K_{\limfunc{right}}\right\} +\left\{ K_{\limfunc{in}},K_{%
\limfunc{out}}\right\} +\left\{ K_{\limfunc{out}},K_{\limfunc{in}}\right\} ^{%
\limfunc{orig}}  \notag \\
&&+\left\{ K_{\limfunc{left}},K_{\limfunc{right}}\right\} +\left\{ K_{%
\limfunc{right}},K_{\limfunc{left}}\right\} +O\left[ \Phi ^{A,B}\left( K_{%
\limfunc{in}}\right) \right] .  \notag
\end{eqnarray}

At this point we observe that the terms $\left\{ K_{\limfunc{in}},K_{%
\limfunc{out}}\right\} $, $\left\{ K_{\limfunc{out}},K_{\limfunc{in}%
}\right\} ^{\limfunc{orig}}$, $\left\{ K_{\limfunc{left}},K_{\limfunc{right}%
}\right\} $, and $\left\{ K_{\limfunc{right}},K_{\limfunc{left}}\right\} $
can all be handled somewhat crudely by separation (\ref{disj supp}):%
\begin{eqnarray*}
&&\left\vert \left\{ K_{\limfunc{in}},K_{\limfunc{out}}\right\} \right\vert
=\left\vert \int_{K_{\limfunc{out}}}\left[ T_{\sigma }^{\alpha }\left( 
\mathbf{1}_{K_{\limfunc{in}}}b_{A}\right) \right] b_{B}^{\ast }d\omega
\right\vert \\
&\lesssim &\sqrt{\mathcal{A}_{2}^{\alpha }}\left( \int_{K_{\limfunc{in}%
}}\left\vert b_{A}\right\vert ^{2}d\sigma \right) ^{\frac{1}{2}}\left(
\int_{K_{\limfunc{out}}}\left\vert b_{B}^{\ast }\right\vert ^{2}d\omega
\right) ^{\frac{1}{2}}\lesssim \sqrt{\mathcal{A}_{2}^{\alpha }}\sqrt{%
\left\vert K_{\limfunc{in}}\right\vert _{\sigma }}\sqrt{\left\vert K_{%
\limfunc{out}}\right\vert _{\omega }}\ ,
\end{eqnarray*}%
and similarly 
\begin{eqnarray*}
&&\left\vert \left\{ K_{\limfunc{out}},K_{\limfunc{in}}\right\} ^{\limfunc{%
orig}}\right\vert \lesssim \sqrt{\mathcal{A}_{2}^{\alpha }}\sqrt{\left\vert
K_{\limfunc{out}}\right\vert _{\sigma }}\sqrt{\left\vert K_{\limfunc{in}%
}\right\vert _{\omega }}\ , \\
&&\left\vert \left\{ K_{\limfunc{left}},K_{\limfunc{right}}\right\}
\right\vert +\left\vert \left\{ K_{\limfunc{right}},K_{\limfunc{left}%
}\right\} \right\vert \lesssim \sqrt{\mathcal{A}_{2}^{\alpha }}\sqrt{%
\left\vert K_{\limfunc{out}}\right\vert _{\sigma }}\sqrt{\left\vert K_{%
\limfunc{out}}\right\vert _{\omega }}\ .
\end{eqnarray*}%
Thus we have%
\begin{equation*}
\left\{ K,K\right\} =\left\{ K_{\limfunc{left}},K_{\limfunc{left}}\right\}
+\left\{ K_{\limfunc{right}},K_{\limfunc{right}}\right\} +O\left[ \Phi
^{A,B}\left( K_{\limfunc{in}}\right) \right] +\mathcal{A}_{2}^{\alpha }\sqrt{%
\left\vert K\right\vert _{\sigma }}\sqrt{\left\vert K\right\vert _{\omega }}.
\end{equation*}%
Upon application of a single iteration we obtain%
\begin{eqnarray*}
\left\{ K,K\right\} &=&\left\{ K_{\limfunc{left}\limfunc{left}},K_{\limfunc{%
left}\limfunc{left}}\right\} +\left\{ K_{\limfunc{left}\limfunc{right}},K_{%
\limfunc{left}\limfunc{right}}\right\} +\left\{ K_{\limfunc{right}\limfunc{%
left}},K_{\limfunc{right}\limfunc{left}}\right\} +\left\{ K_{\limfunc{right}%
\limfunc{right}},K_{\limfunc{right}\limfunc{right}}\right\} \\
&&+O\left[ \Phi ^{A,B}\left( K_{\limfunc{in}}\right) +\Phi ^{A,B}\left( K_{%
\limfunc{left}\limfunc{in}}\right) +\Phi ^{A,B}\left( K_{\limfunc{right}%
\limfunc{in}}\right) \right] \\
&&+\mathcal{A}_{2}^{\alpha }\left( \sqrt{\left\vert K\right\vert _{\sigma }}%
\sqrt{\left\vert K\right\vert _{\omega }}+\sqrt{\left\vert K_{\limfunc{left}%
}\right\vert _{\sigma }}\sqrt{\left\vert K_{\limfunc{left}}\right\vert
_{\omega }}+\sqrt{\left\vert K_{\limfunc{right}}\right\vert _{\sigma }}\sqrt{%
\left\vert K_{\limfunc{right}}\right\vert _{\omega }}\right) ,
\end{eqnarray*}%
and then iterating finitely many more times gives for $n\in \mathbb{N}$,%
\begin{eqnarray}
\left\{ K,K\right\} &=&\sum_{M\in \mathcal{M}_{n}}\left\{ M,M\right\}
+O\left( \sum_{M\in \mathcal{M}_{n}^{\ast }}\left[ \Phi ^{A,B}\left( M_{%
\limfunc{in}}\right) \right] \right) +\sqrt{\mathcal{A}_{2}^{\alpha }}%
\sum_{M\in \mathcal{M}_{n}^{\ast }}\sqrt{\left\vert M\right\vert _{\sigma }}%
\sqrt{\left\vert M\right\vert _{\omega }}  \label{K iterated} \\
&\equiv &A\left( K\right) +B\left( K\right) +C\left( K\right) =A_{\left(
I^{\prime },J^{\prime }\right) }\left( K\right) +B_{\left( I^{\prime
},J^{\prime }\right) }\left( K\right) +C_{\left( I^{\prime },J^{\prime
}\right) }\left( K\right) ,  \notag
\end{eqnarray}%
where the collections of intervals $\mathcal{M}_{n}=\mathcal{M}_{n}\left(
K\right) $ and $\mathcal{M}_{n}^{\ast }=\mathcal{M}_{n}^{\ast }\left(
K\right) $ are defined recursively by 
\begin{eqnarray*}
\mathcal{M}_{0} &\equiv &\left\{ K\right\} , \\
\mathcal{M}_{k+1} &\equiv &\dbigcup \left\{ M_{\limfunc{left}},M_{\limfunc{%
right}}:M\in \mathcal{M}_{k}\right\} ,\ \ \ \ \ k\geq 0, \\
\mathcal{M}_{n}^{\ast } &\equiv &\dbigcup_{k=0}^{n}\mathcal{M}_{k}\ .
\end{eqnarray*}%
We will include the subscript $\left( I^{\prime },J^{\prime }\right) $ in
the notation when we want to indicate the pair $\left( I^{\prime },J^{\prime
}\right) $ for which $K\in \mathcal{K}\left( I^{\prime },J^{\prime }\right) $
as defined in (\ref{def K(I',J')}) above. Now the term $C\left( K\right) $
can be estimated by the crude estimate 
\begin{equation}
C\left( K\right) =\sqrt{\mathcal{A}_{2}^{\alpha }}\sum_{M\in \mathcal{M}%
_{n}^{\ast }}\sqrt{\left\vert M\right\vert _{\sigma }}\sqrt{\left\vert
M\right\vert _{\omega }}\leq C_{n}\sqrt{\mathcal{A}_{2}^{\alpha }}\sqrt{%
\left\vert K\right\vert _{\sigma }}\sqrt{\left\vert K\right\vert _{\omega }}%
\ ,  \label{C est}
\end{equation}%
where $n$ is chosen below depending on $\eta _{0}$. For the first term $%
A\left( K\right) $, we will apply the norm inequality and use probability,
namely%
\begin{eqnarray*}
\left\vert A\left( K\right) \right\vert &\leq &\sqrt{C_{\mathbf{b}}C_{%
\mathbf{b}^{\ast }}}\mathfrak{N}_{T^{\alpha }}\sum_{M\in \mathcal{M}_{n}}%
\sqrt{\left\vert M\right\vert _{\sigma }}\sqrt{\left\vert M\right\vert
_{\omega }} \\
&\leq &\sqrt{C_{\mathbf{b}}C_{\mathbf{b}^{\ast }}}\mathfrak{N}_{T^{\alpha }}%
\sqrt{\sum_{M\in \mathcal{M}_{n}}\left\vert M\right\vert _{\sigma }}\sqrt{%
\sum_{M\in \mathcal{M}_{n}}\left\vert M\right\vert _{\omega }}\leq \sqrt{C_{%
\mathbf{b}}C_{\mathbf{b}^{\ast }}}\mathfrak{N}_{T^{\alpha }}\sqrt{\sum_{M\in 
\mathcal{M}_{n}}\left\vert M\right\vert _{\sigma }}\sqrt{\left\vert
K\right\vert _{\omega }},
\end{eqnarray*}%
where $\sqrt{C_{\mathbf{b}}C_{\mathbf{b}^{\ast }}}$ is an upper bound for
the testing functions involved, followed by 
\begin{equation*}
\boldsymbol{E}_{\Omega }^{\mathcal{G}}\left( \sum_{M\in \mathcal{M}%
_{n}}\left\vert M\right\vert _{\sigma }\right) \leq \varepsilon \left\vert
I^{\prime }\right\vert _{\sigma }\ ,
\end{equation*}%
for a sufficiently small $\varepsilon >0$, where \emph{roughly speaking}, we
use the fact that the intervals $M\in \mathcal{M}_{n}$ depend on the grid $%
\mathcal{G}$ and form a relatively small proportion of $I^{\prime }$, which
captures only a small amount of the total mass $\left\vert I^{\prime
}\right\vert _{\sigma }$ as the grid is translated relative to the grid $%
\mathcal{D}$ that contains $I^{\prime }$. To be a bit more precise, recall
that the intervals $K$ are uniquely constructed as consecutive adjacent
intervals of equal length in the grid $\mathcal{G}$ that start out at the
endpoint of $J^{\prime }$ that lies inside $I^{\prime }$, and then progress
toward the constructed endpoint of $I^{\prime \prime }=I^{\prime }\setminus
\partial _{\eta _{1}}I^{\prime }$ lying in the interior of $J^{\prime }$.
Thus translates of the grid $\mathcal{G}$ result in translating the
intervals $K$ across the interval $I^{\prime }$ a distance comparable to at
least the length of $K$, and where $I^{\prime }$ is fixed in the grid $%
\mathcal{D}$. As a consequence the intervals $M\in \mathcal{M}_{n}\left(
K\right) $ are also translated across the fixed interval $I^{\prime }$ a
distance comparable to at least $\ell \left( K\right) $, and the standard
halo estimate (\ref{hand'}) applies since the sum of the lengths of the
intervals $M\in \mathcal{M}_{n}$ is a small proportion of the length of $%
J^{\prime }$, whose length is at most that of $I^{\prime }$.

Here are the specific details. Recall that the intervals $K$ are taken from
the set of consecutive intervals $\left\{ K_{i}\right\} _{i=1}^{B}$ that lie
in $I^{\prime }\cap J^{\prime }$, that the intervals $M\in \mathcal{M}%
_{n}\left( K_{i}\right) $ have length $\frac{1}{4^{n}}\ell \left(
K_{i}\right) $, and that there are $2^{n}$ such intervals in $\mathcal{M}%
_{n}\left( K_{i}\right) $ for each $i$. Thus we have 
\begin{eqnarray}
\sum_{M\in \mathcal{M}_{n}\left( K\right) }\left\vert M\right\vert
&=&\sum_{M\in \mathcal{M}_{n}\left( K\right) }\frac{1}{4^{n}}\left\vert
K\right\vert =2^{n}\frac{1}{4^{n}}\left\vert K\right\vert =\frac{1}{2^{n}}%
\left\vert K\right\vert  \label{future prob} \\
&\Longrightarrow &\boldsymbol{E}_{\Omega }^{\mathcal{G}}\left(
\sum_{i=1}^{B}\sum_{M\in \mathcal{M}_{n}\left( K_{i}\right) }\left\vert
M\right\vert _{\sigma }\right) \leq C\frac{1}{\eta _{0}}\frac{1}{2^{n}}%
\left\vert I^{\prime }\right\vert _{\sigma }\leq \eta _{0}\left\vert
I^{\prime }\right\vert _{\sigma }\ ,  \notag
\end{eqnarray}%
where we have used that the variable $B$ is at most $\frac{1}{\eta _{0}}$,
and where the final inequality holds if $n$ is chosen large enough that $%
\frac{1}{2^{n}}\leq \eta _{0}^{2}$. Then we have by Cauchy-Schwarz applied
first to $\sum_{i=1}^{B}\sum_{M\in \mathcal{M}_{n}\left( K_{i}\right) }$ and
then to $\boldsymbol{E}_{\Omega }^{\mathcal{G}}$, 
\begin{eqnarray}
&&\boldsymbol{E}_{\Omega }^{\mathcal{G}}\left( \sum_{i=1}^{B}\left\vert
A\left( K_{i}\right) \right\vert \right) \leq \boldsymbol{E}_{\Omega }^{%
\mathcal{G}}\sqrt{C_{\mathbf{b}}C_{\mathbf{b}^{\ast }}}\mathfrak{N}%
_{T^{\alpha }}\sqrt{\sum_{i=1}^{B}\sum_{M\in \mathcal{M}_{n}\left(
K_{i}\right) }\left\vert M\right\vert _{\sigma }}\sqrt{\left\vert J^{\prime
}\right\vert _{\omega }}  \label{A est} \\
&\leq &\sqrt{C_{\mathbf{b}}C_{\mathbf{b}^{\ast }}}\mathfrak{N}_{T^{\alpha }}%
\sqrt{\boldsymbol{E}_{\Omega }^{\mathcal{G}}\sum_{i=1}^{B}\sum_{M\in 
\mathcal{M}_{n}\left( K_{i}\right) }\left\vert M\right\vert _{\sigma }}\sqrt{%
\left\vert J^{\prime }\right\vert _{\omega }}  \notag \\
&\leq &\sqrt{C_{\mathbf{b}}C_{\mathbf{b}^{\ast }}}\mathfrak{N}_{T^{\alpha }}%
\sqrt{\eta _{0}\left\vert I^{\prime }\right\vert _{\sigma }}\sqrt{\left\vert
J^{\prime }\right\vert _{\omega }}=\sqrt{C_{\mathbf{b}}C_{\mathbf{b}^{\ast }}%
}\sqrt{\eta _{0}}\mathfrak{N}_{T^{\alpha }}\sqrt{\left\vert I^{\prime
}\right\vert _{\sigma }}\sqrt{\left\vert J^{\prime }\right\vert _{\omega }},
\notag
\end{eqnarray}%
as required.

Now we turn to summing up the remaining terms $B\left( K\right) =C\sum_{M\in 
\mathcal{M}_{n}^{\ast }}\Phi ^{A,B}\left( M_{\limfunc{in}}\right) $ above.
We begin by claiming that in the case when the interval $I^{\prime }$ is a 
\emph{natural} child of $I$, i.e. $I^{\prime }\in \mathfrak{C}_{\limfunc{%
natural}}\left( I\right) $ so that $I^{\prime }\in \mathcal{C}_{A}^{\mathcal{%
A}}$, we have 
\begin{eqnarray}
&&\left( \sum_{M\in \mathcal{M}_{n}^{\ast }\left( K\right) }\left\Vert 
\mathbf{1}_{M_{\limfunc{in}}}T_{\sigma }^{\alpha }b_{A}\right\Vert
_{L^{2}\left( \omega \right) }^{2}\right) ^{\frac{1}{2}}\left( \sum_{M\in 
\mathcal{M}_{n}^{\ast }\left( K\right) }\left\vert M_{\limfunc{in}%
}\right\vert _{\omega }\right) ^{\frac{1}{2}}  \label{begin claim} \\
&&+\left( \sum_{M\in \mathcal{M}_{n}^{\ast }\left( K\right) }\left\{ \mathrm{%
P}_{\delta }^{\alpha }\mathsf{Q}^{\omega }\left( M_{\limfunc{in}}^{\limfunc{%
left}},\mathbf{1}_{A}\sigma \right) ^{2}+\mathrm{P}_{\delta }^{\alpha }%
\mathsf{Q}^{\omega }\left( M_{\limfunc{in}}^{\limfunc{right}},\mathbf{1}%
_{A}\sigma \right) ^{2}\right\} \right) ^{\frac{1}{2}}\left( \sum_{M\in 
\mathcal{M}_{n}^{\ast }\left( K\right) }\left\vert M_{\limfunc{in}%
}\right\vert _{\omega }\right) ^{\frac{1}{2}}  \notag \\
&\lesssim &\left( \mathfrak{T}_{T^{\alpha }}^{\mathbf{b}}+\mathcal{E}%
_{2}^{\alpha }+\sqrt{\mathcal{A}_{2}^{\alpha ,\ast }+\mathcal{A}_{2}^{\alpha
}}\right) \sqrt{\left\vert I^{\prime }\right\vert _{\sigma }\left\vert
J^{\prime }\right\vert _{\omega }}\lesssim \mathcal{NTV}_{\alpha }\sqrt{%
\left\vert I^{\prime }\right\vert _{\sigma }\left\vert J^{\prime
}\right\vert _{\omega }},  \notag
\end{eqnarray}%
where the last line is a consequence of the crucial fact that the intervals $%
\left\{ M_{\limfunc{in}}\right\} _{M\in \mathcal{M}_{n}^{\ast }\left(
K\right) }$ form a pairwise disjoint subdecomposition of $K\subset I^{\prime
}\cap J^{\prime }$ (for any $n\geq 1$). Indeed, we then have the following
inequalities

\begin{enumerate}
\item the first sum on the left hand side satisfies%
\begin{equation*}
\sum_{M\in \mathcal{M}_{n}^{\ast }\left( K\right) }\left\Vert \mathbf{1}_{M_{%
\limfunc{in}}}T_{\sigma }^{\alpha }b_{A}\right\Vert _{L^{2}\left( \omega
\right) }^{2}=\sum_{M\in \mathcal{M}_{n}^{\ast }\left( K\right) }\int_{M_{%
\limfunc{in}}}\left\vert T_{\sigma }^{\alpha }b_{A}\right\vert ^{2}d\omega
\leq \int_{I^{\prime }}\left\vert T_{\sigma }^{\alpha }b_{A}\right\vert
^{2}d\omega \lesssim \left( \mathfrak{T}_{T^{\alpha }}^{\mathbf{b}}\right)
^{2}\left\vert I^{\prime }\right\vert _{\sigma }\ ,
\end{equation*}%
by the weak testing condition for $I^{\prime }$ in the corona $\mathcal{C}%
_{A}$,

\item and%
\begin{equation*}
\sum_{M\in \mathcal{M}_{n}^{\ast }\left( K\right) }\left\vert M_{\limfunc{in}%
}\right\vert _{\omega }\leq \left\vert K\right\vert _{\omega }\leq
\left\vert J^{\prime }\right\vert _{\omega }\ ,
\end{equation*}

\item and, using the definition of $\mathrm{P}_{\delta }^{\alpha }\mathsf{Q}%
^{\omega }\left( J,\upsilon \right) $ in (\ref{def compact}), 
\begin{eqnarray*}
&&\sum_{M\in \mathcal{M}_{n}^{\ast }\left( K\right) }\left\{ \mathrm{P}%
_{\delta }^{\alpha }\mathsf{Q}^{\omega }\left( M_{\limfunc{in}}^{\limfunc{%
left}},\mathbf{1}_{A}\sigma \right) ^{2}+\mathrm{P}_{\delta }^{\alpha }%
\mathsf{Q}^{\omega }\left( M_{\limfunc{in}}^{\limfunc{right}},\mathbf{1}%
_{A}\sigma \right) ^{2}\right\} \\
&\lesssim &\sum_{M\in \mathcal{M}_{n}^{\ast }\left( K\right) }\left\{ \left( 
\frac{\mathrm{P}^{\alpha }\left( M_{\limfunc{in}}^{\limfunc{left}},%
\boldsymbol{1}_{A}\sigma \right) }{\left\vert M_{\limfunc{in}}^{\limfunc{left%
}}\right\vert }\right) ^{2}\left\Vert x-m_{M_{\limfunc{in}}^{\limfunc{left}%
}}\right\Vert _{L^{2}\left( \mathbf{1}_{M_{\limfunc{in}}^{\limfunc{left}%
}}\omega \right) }^{2}+\left( \frac{\mathrm{P}^{\alpha }\left( M_{\limfunc{in%
}}^{\limfunc{right}},\boldsymbol{1}_{A}\sigma \right) }{\left\vert M_{%
\limfunc{in}}^{\limfunc{right}}\right\vert }\right) ^{2}\left\Vert x-m_{M_{%
\limfunc{in}}^{\limfunc{right}}}\right\Vert _{L^{2}\left( \mathbf{1}_{M_{%
\limfunc{in}}^{\limfunc{right}}}\omega \right) }^{2}\right\} \\
&\lesssim &\mathcal{NTV}_{\alpha }\left\vert I^{\prime }\right\vert _{\sigma
}\ ,
\end{eqnarray*}%
upon using the stopping energy condition for $I^{\prime }$ in the corona $%
\mathcal{C}_{A}$, i.e. the failure of (\ref{def stop 3}), in the corona $%
\mathcal{C}_{A}$ with the subdecomposition $I^{\prime }\supset \overset{%
\cdot }{\dbigcup }_{M\in \mathcal{M}_{n}^{\ast }\left( K\right) }\left( M_{%
\limfunc{in}}^{\limfunc{left}}\overset{\cdot }{\dbigcup }M_{\limfunc{in}}^{%
\limfunc{right}}\right) $.
\end{enumerate}

This completes the proof of (\ref{begin claim}). Our next claim is the
inequality%
\begin{eqnarray}
&&\left( \sum_{M\in \mathcal{M}_{n}^{\ast }\left( K\right) }\int_{M_{%
\limfunc{in}}}\left\vert T_{\sigma }^{\alpha }b_{A}\right\vert ^{2}d\omega
+\left( \mathfrak{T}_{T^{\alpha ,\ast }}+\mathfrak{A}_{2}^{\alpha }\right)
^{2}\left\vert M_{\limfunc{in}}\right\vert _{\sigma }\right) ^{\frac{1}{2}%
}\left( \sum_{M\in \mathcal{M}_{n}^{\ast }\left( K\right) }\left\vert M_{%
\limfunc{in}}\right\vert _{\omega }\right) ^{\frac{1}{2}}  \label{next claim}
\\
&\lesssim &\left( \int_{I^{\prime }}\left\vert T_{\sigma }^{\alpha
}b_{A}\right\vert ^{2}d\omega +\left( \mathfrak{T}_{T^{\alpha ,\ast }}+%
\mathfrak{A}_{2}^{\alpha }\right) ^{2}\left\vert I^{\prime }\right\vert
_{\sigma }\right) ^{\frac{1}{2}}\left( \left\vert J^{\prime }\right\vert
_{\omega }\right) ^{\frac{1}{2}}\lesssim \mathcal{NTV}_{\alpha }\sqrt{%
\left\vert I^{\prime }\right\vert _{\sigma }\left\vert J^{\prime
}\right\vert _{\omega }},  \notag
\end{eqnarray}%
and combining this with (\ref{begin claim}), together with the definition of 
$\Phi ^{A,B}$ in (\ref{def PHI}), gives%
\begin{equation*}
\sum_{M\in \mathcal{M}_{n}^{\ast }\left( K\right) }\Phi ^{A,B}\left( M_{%
\limfunc{in}}\right) \lesssim L{\small eft}H{\small and}S{\small ide}\left( %
\ref{begin claim}\right) +L{\small eft}H{\small and}S{\small ide}\left( \ref%
{next claim}\right) \lesssim \mathcal{NTV}_{\alpha }\sqrt{\left\vert
I^{\prime }\right\vert _{\sigma }\left\vert J^{\prime }\right\vert _{\omega }%
}.
\end{equation*}

In order to deal with this sum in the case when the child $I^{\prime }$ is
broken, we must take the estimate one step further and sum over those broken
intervals $I^{\prime }$ whose parents belong to the corona $\mathcal{C}_{A}$%
, i.e. $\left\{ I^{\prime }\in \mathcal{D}:I^{\prime }\in \mathfrak{C}_{%
\limfunc{broken}}\left( I\right) \text{ for some }I\in \mathcal{C}%
_{A}\right\} $. Of course this collection is precisely the set of $\mathcal{A%
}$ -children of $A$, i.e.%
\begin{equation}
\left\{ I^{\prime }\in \mathcal{D}:I^{\prime }\in \mathfrak{C}_{\limfunc{%
broken}}\left( I\right) \text{ for some }I\in \mathcal{C}_{A}\right\} =%
\mathfrak{C}_{\mathcal{A}}\left( A\right) .  \label{precisely}
\end{equation}

To help motivate this, we first recall that we denote the term $B\left(
K\right) $ by $B_{\left( I^{\prime },J^{\prime }\right) }\left( K\right) $
when we wish to indicate the pair $\left( I^{\prime },J^{\prime }\right) $
to which $K$ is associated, i.e. $K\in \mathcal{K}\left( I^{\prime
},J^{\prime }\right) $ as in (\ref{def K(I',J')}). Of course the intervals $%
M\in \mathcal{M}\left( K\right) $ also depend on the pair of intervals $%
I^{\prime }$ and $J^{\prime }$, but we will suppress notation to this
effect, and the reader should keep this in mind. In particular then, if we
now sum over \emph{natural} children $I^{\prime }$ of $I$ $\in \mathcal{C}%
_{A}$ and the associated children $J^{\prime }$ of $J\in \mathcal{C}_{A}^{%
\mathcal{G},\limfunc{nearby}}\left( I\right) $, where%
\begin{equation*}
\mathcal{C}_{A}^{\mathcal{G},\limfunc{nearby}}\left( I\right) \equiv \left\{
J\in \mathcal{G}:2^{-\mathbf{r}}\ell \left( I\right) <\ell \left( J\right)
\leq \ell \left( I\right) \text{ and }d\left( J,I\right) \leq 2\ell \left(
J\right) ^{\varepsilon }\ell \left( I\right) ^{1-\varepsilon }\right\} ,\ \
\ \ \ \text{for }I\in \mathcal{C}_{A}\ ,
\end{equation*}%
we obtain the following corona estimate, using the collection $\mathcal{K}%
\left( I^{\prime },J^{\prime }\right) $ that is defined in (\ref{def
K(I',J')}) above with $B\leq C\frac{1}{\eta _{0}}$, 
\begin{eqnarray}
&&\sum_{I\in \mathcal{C}_{A}\text{ and }J\in \mathcal{C}_{A}^{\mathcal{G},%
\limfunc{nearby}}\left( I\right) }\sum_{\substack{ I^{\prime }\in \mathfrak{C%
}_{\limfunc{natural}}\left( I\right) \text{ and }J^{\prime }\in \mathfrak{C}%
\left( J\right)  \\ K\in \mathcal{K}\left( I^{\prime },J^{\prime }\right) }}%
\left\vert E_{I^{\prime }}^{\sigma }\left( \widehat{\square }_{I}^{\sigma
,\flat ,\mathbf{b}}f\right) \right\vert \ \left\vert B_{\left( I^{\prime
},J^{\prime }\right) }\left( K\right) \right\vert \ \left\vert E_{J^{\prime
}}^{\omega }\left( \widehat{\square }_{J}^{\omega ,\flat ,\mathbf{b}^{\ast
}}g\right) \right\vert  \label{prelim corona natural} \\
&\lesssim &\frac{1}{\eta _{0}}\mathcal{NTV}_{\alpha }\sum_{I\in \mathcal{C}%
_{A}\text{ and }J\in \mathcal{C}_{A}^{\mathcal{G},\limfunc{nearby}}\left(
I\right) }\sum_{I^{\prime }\in \mathfrak{C}_{\limfunc{natural}}\left(
I\right) \text{ and }J^{\prime }\in \mathfrak{C}\left( J\right) }\left\vert
E_{I^{\prime }}^{\sigma }\left( \widehat{\square }_{I}^{\sigma ,\flat ,%
\mathbf{b}}f\right) \right\vert \ \sqrt{\left\vert I^{\prime }\right\vert
_{\sigma }\left\vert J^{\prime }\right\vert _{\omega }}\ \left\vert
E_{J^{\prime }}^{\omega }\left( \widehat{\square }_{J}^{\omega ,\flat ,%
\mathbf{b}^{\ast }}g\right) \right\vert  \notag \\
&\lesssim &\frac{1}{\eta _{0}}\mathcal{NTV}_{\alpha }\left( \sum_{I\in 
\mathcal{C}_{A}}\sum_{I^{\prime }\in \mathfrak{C}_{\limfunc{natural}}\left(
I\right) }\left\vert I^{\prime }\right\vert _{\sigma }\left\vert
E_{I^{\prime }}^{\sigma }\left( \widehat{\square }_{I}^{\sigma ,\flat ,%
\mathbf{b}}f\right) \right\vert ^{2}\right) ^{\frac{1}{2}}\left( \sum_{I\in 
\mathcal{C}_{A}}\sum_{J\in \mathcal{C}_{A}^{\mathcal{G},\limfunc{nearby}%
}\left( I\right) }\sum_{J^{\prime }\in \mathfrak{C}\left( J\right)
}\left\vert J^{\prime }\right\vert _{\omega }\ \left\vert E_{J^{\prime
}}^{\omega }\left( \widehat{\square }_{J}^{\omega ,\flat ,\mathbf{b}^{\ast
}}g\right) \right\vert ^{2}\right) ^{\frac{1}{2}}  \notag \\
&\lesssim &\frac{1}{\eta _{0}}\mathcal{NTV}_{\alpha }\left\Vert \mathsf{P}_{%
\mathcal{C}_{A}}^{\sigma }f\right\Vert _{L^{2}\left( \sigma \right)
}^{\bigstar }\left\Vert \mathsf{P}_{\mathcal{C}_{A}^{\mathcal{G},\limfunc{%
nearby}}}^{\omega }g\right\Vert _{L^{2}\left( \sigma \right) }^{\bigstar }\ ,
\notag
\end{eqnarray}%
where $\mathcal{C}_{A}^{\mathcal{G},\limfunc{nearby}}=\dbigcup\limits_{I\in 
\mathcal{C}_{A}}\mathcal{C}_{A}^{\mathcal{G},\limfunc{nearby}}\left(
I\right) $, and the final line uses (\ref{box hat bound}) to obtain 
\begin{equation*}
\sum_{I\in \mathcal{C}_{A}}\sum_{I^{\prime }\in \mathfrak{C}_{\limfunc{%
natural}}\left( I\right) }\left\vert I^{\prime }\right\vert _{\sigma
}\left\vert E_{I^{\prime }}^{\sigma }\left( \widehat{\square }_{I}^{\sigma
,\flat ,\mathbf{b}}f\right) \right\vert ^{2}=\sum_{I\in \mathcal{C}%
_{A}}\left\Vert \widehat{\square }_{I}^{\sigma ,\flat ,\mathbf{b}%
}f\right\Vert _{L^{2}\left( \sigma \right) }^{2}\lesssim \sum_{I\in \mathcal{%
C}_{A}}\left\Vert \square _{I}^{\sigma ,\mathbf{b}}f\right\Vert
_{L^{2}\left( \sigma \right) }^{2}\leq \left\Vert \mathsf{P}_{\mathcal{C}%
_{A}}^{\sigma }f\right\Vert _{L^{2}\left( \sigma \right) }^{\bigstar 2},
\end{equation*}%
and similarly for the sum in $J$ and $J^{\prime }$, once we note that given $%
J\in \mathcal{C}_{A}^{\mathcal{G},\limfunc{nearby}}$, there are only
boundedly many $I\in \mathcal{C}_{A}$ for which $J\in \mathcal{C}_{A}^{%
\mathcal{G},\limfunc{nearby}}\left( I\right) $.

To obtain the same corona estimate when summing over broken $I^{\prime }$,
we will exploit the fact that the intervals $A^{\prime }\in \mathfrak{C}_{%
\mathcal{A}}\left( A\right) $ are pairwise disjoint. But first we note that
when $I^{\prime }$ is a broken child, neither weak testing nor stopping
energy is available. But if we sum over such broken $I^{\prime }$, and use (%
\ref{precisely}) to see that the broken children are pairwise disjoint, we
obtain the following estimate where for convenience we write%
\begin{equation*}
\mathrm{P}_{\delta }^{\alpha }\mathsf{Q}^{\omega }\left( M_{\limfunc{in}}^{%
\limfunc{left}/\limfunc{right}},\mathbf{1}_{A}\sigma \right) ^{2}\equiv 
\mathrm{P}_{\delta }^{\alpha }\mathsf{Q}^{\omega }\left( M_{\limfunc{in}}^{%
\limfunc{left}},\mathbf{1}_{A}\sigma \right) ^{2}+\mathrm{P}_{\delta
}^{\alpha }\mathsf{Q}^{\omega }\left( M_{\limfunc{in}}^{\limfunc{right}},%
\mathbf{1}_{A}\sigma \right) ^{2},
\end{equation*}%
and we use the notation $\mathcal{M}_{n}^{\ast }\left( I^{\prime },J^{\prime
}\right) \equiv \dbigcup\limits_{K\in \mathcal{K}\left( I^{\prime
},J^{\prime }\right) }\mathcal{M}_{n}^{\ast }\left( K\right) $: 
\begin{eqnarray*}
&&\sum_{I\in \mathcal{C}_{A}\text{ and }J\in \mathcal{C}_{A}^{\mathcal{G},%
\limfunc{nearby}}\left( I\right) }\sum_{\substack{ I^{\prime }\in \mathfrak{C%
}_{\limfunc{broken}}\left( I\right) \text{ and }J^{\prime }\in \mathfrak{C}%
\left( J\right)  \\ K\in \mathcal{K}\left( I^{\prime },J^{\prime }\right) }}%
\left\vert E_{I^{\prime }}^{\sigma }\left( \widehat{\square }_{I}^{\sigma
,\flat ,\mathbf{b}}f\right) \right\vert \ \left\vert B_{\left( I^{\prime
},J^{\prime }\right) }\left( K\right) \right\vert \ \left\vert E_{J^{\prime
}}^{\omega }\left( \widehat{\square }_{J}^{\omega ,\flat ,\mathbf{b}^{\ast
}}g\right) \right\vert \\
&\lesssim &\mathcal{NTV}_{\alpha }\sum_{I\in \mathcal{C}_{A}\text{ and }J\in 
\mathcal{C}_{A}^{\mathcal{G},\limfunc{nearby}}\left( I\right)
}\sum_{I^{\prime }\in \mathfrak{C}_{\limfunc{broken}}\left( I\right) \text{
and }J^{\prime }\in \mathfrak{C}\left( J\right) }\left\vert E_{I^{\prime
}}^{\sigma }\left( \widehat{\square }_{I}^{\sigma ,\flat ,\mathbf{b}%
}f\right) \right\vert \\
&&\times \sqrt{\sum_{M\in \mathcal{M}_{n}^{\ast }\left( I^{\prime
},J^{\prime }\right) }\left\Vert \mathbf{1}_{M_{\limfunc{in}}}T_{\sigma
}^{\alpha }b_{A}\right\Vert _{L^{2}\left( \omega \right) }^{2}+\sum_{M\in 
\mathcal{M}_{n}^{\ast }\left( I^{\prime },J^{\prime }\right) }\mathrm{P}%
_{\delta }^{\alpha }\mathsf{Q}^{\omega }\left( M_{\limfunc{in}}^{\limfunc{%
left}/\limfunc{right}},\mathbf{1}_{A}\sigma \right) ^{2}+\sum_{M\in \mathcal{%
M}_{n}^{\ast }\left( I^{\prime },J^{\prime }\right) }\left\vert M_{\limfunc{%
in}}\right\vert _{\sigma }} \\
&&\ \ \ \ \ \ \ \ \ \ \ \ \ \ \ \ \ \ \ \ \ \ \ \ \ \ \ \ \ \ \ \ \ \ \ \ \
\ \ \ \ \ \ \ \ \ \ \ \ \ \ \ \ \ \ \ \ \ \ \ \ \ \ \ \ \times \sqrt{%
\left\vert J^{\prime }\right\vert _{\omega }}\ \left\vert E_{J^{\prime
}}^{\omega }\left( \widehat{\square }_{J}^{\omega ,\flat ,\mathbf{b}^{\ast
}}g\right) \right\vert ,
\end{eqnarray*}%
which gives%
\begin{eqnarray}
&&  \label{prelim corona broken} \\
&&\sum_{I\in \mathcal{C}_{A}\text{ and }J\in \mathcal{C}_{A}^{\mathcal{G},%
\limfunc{nearby}}\left( I\right) }\sum_{\substack{ I^{\prime }\in \mathfrak{C%
}_{\limfunc{broken}}\left( I\right) \text{ and }J^{\prime }\in \mathfrak{C}%
\left( J\right)  \\ K\in \mathcal{K}\left( I^{\prime },J^{\prime }\right) }}%
\left\vert E_{I^{\prime }}^{\sigma }\left( \widehat{\square }_{I}^{\sigma
,\flat ,\mathbf{b}}f\right) \right\vert \ \left\vert B_{\left( I^{\prime
},J^{\prime }\right) }\left( K\right) \right\vert \ \left\vert E_{J^{\prime
}}^{\omega }\left( \widehat{\square }_{J}^{\omega ,\flat ,\mathbf{b}^{\ast
}}g\right) \right\vert  \notag \\
&\lesssim &\mathcal{NTV}_{\alpha }\left( \sum_{\substack{ I\in \mathcal{C}%
_{A}  \\ I^{\prime }\in \mathfrak{C}_{\limfunc{broken}}\left( I\right) }}%
\sum _{\substack{ J\in \mathcal{C}_{A}^{\mathcal{G},\limfunc{nearby}}\left(
I\right)  \\ J^{\prime }\in \mathfrak{C}\left( J\right) }}\sum_{M\in 
\mathcal{M}_{n}^{\ast }\left( I^{\prime },J^{\prime }\right) }\left\{
\left\Vert \mathbf{1}_{M_{\limfunc{in}}}T_{\sigma }^{\alpha
}b_{A}\right\Vert _{L^{2}\left( \omega \right) }^{2}+\mathrm{P}_{\delta
}^{\alpha }\mathsf{Q}^{\omega }\left( M_{\limfunc{in}}^{\limfunc{left}/%
\limfunc{right}},\mathbf{1}_{A}\sigma \right) ^{2}+\left\vert M_{\limfunc{in}%
}\right\vert _{\sigma }\right\} \right) ^{\frac{1}{2}}  \notag \\
&&\ \ \ \ \ \ \ \ \ \ \ \ \ \ \ \ \ \ \ \ \ \ \ \ \ \ \ \ \ \ \ \ \ \ \
\times \left( \frac{1}{\left\vert A\right\vert _{\sigma }}\int_{A}\left\vert
f\right\vert d\sigma \right) \left( \sum_{\substack{ J\in \mathcal{C}_{A}^{%
\mathcal{G},\limfunc{nearby}}  \\ J^{\prime }\in \mathfrak{C}\left( J\right) 
}}\sum_{\substack{ I\in \mathcal{C}_{A}:\ J\in \mathcal{C}_{A}^{\mathcal{G},%
\limfunc{nearby}}\left( I\right)  \\ I^{\prime }\in \mathfrak{C}_{\limfunc{%
broken}}\left( I\right) }}\left\vert J^{\prime }\right\vert _{\omega
}\left\vert E_{J^{\prime }}^{\omega }\left( \widehat{\square }_{J}^{\omega
,\flat ,\mathbf{b}^{\ast }}g\right) \right\vert ^{2}\right) ^{\frac{1}{2}} 
\notag \\
&\lesssim &\mathcal{NTV}_{\alpha }\sqrt{\left\vert A\right\vert _{\sigma
}\left( \frac{1}{\left\vert A\right\vert _{\sigma }}\int_{A}\left\vert
f\right\vert d\sigma \right) ^{2}}\left\Vert \mathsf{P}_{\mathcal{C}_{A}^{%
\mathcal{G},\limfunc{nearby}}}^{\omega }g\right\Vert _{L^{2}\left( \sigma
\right) }^{\bigstar }\ ,  \notag
\end{eqnarray}%
because%
\begin{equation*}
\left\vert E_{I^{\prime }}^{\sigma }\left( \widehat{\square }_{I}^{\sigma
,\flat ,\mathbf{b}}f\right) \right\vert =\left\vert \frac{1}{%
\int_{I}b_{I}d\sigma }\int_{I}fd\sigma \right\vert \lesssim \frac{1}{%
\left\vert I\right\vert _{\sigma }}\int_{I}\left\vert f\right\vert d\sigma
\lesssim \frac{1}{\left\vert A\right\vert _{\sigma }}\int_{A}\left\vert
f\right\vert d\sigma
\end{equation*}%
if $I^{\prime }\in \mathfrak{C}_{\limfunc{broken}}\left( I\right) $ and $%
I\in \mathcal{C}_{A}$, and because 
\begin{eqnarray}
&&  \label{last inequ} \\
&&\sum_{\substack{ I\in \mathcal{C}_{A}  \\ I^{\prime }\in \mathfrak{C}_{%
\limfunc{broken}}\left( I\right) }}\sum_{\substack{ J\in \mathcal{C}_{A}^{%
\mathcal{G},\limfunc{nearby}}\left( I\right)  \\ J^{\prime }\in \mathfrak{C}%
\left( J\right) }}\sum_{M\in \mathcal{M}_{n}^{\ast }\left( I^{\prime
},J^{\prime }\right) }\left\{ \left\Vert \mathbf{1}_{M_{\limfunc{in}%
}}T_{\sigma }^{\alpha }b_{A}\right\Vert _{L^{2}\left( \omega \right) }^{2}+%
\mathrm{P}_{\delta }^{\alpha }\mathsf{Q}^{\omega }\left( M_{\limfunc{in}}^{%
\limfunc{left}},\mathbf{1}_{A}\sigma \right) ^{2}+\mathrm{P}_{\delta
}^{\alpha }\mathsf{Q}^{\omega }\left( M_{\limfunc{in}}^{\limfunc{right}},%
\mathbf{1}_{A}\sigma \right) ^{2}+\left\vert M_{\limfunc{in}}\right\vert
_{\sigma }\right\}  \notag \\
&\leq &\int_{A}\left\vert T_{\sigma }^{\alpha }b_{A}\right\vert ^{2}d\omega +%
\mathfrak{E}_{2}^{\alpha }\left\vert A\right\vert _{\sigma }+\left\vert
A\right\vert _{\sigma }\leq \left( \mathfrak{T}_{T^{\alpha }}^{\mathbf{b}}+%
\mathfrak{E}_{2}^{\alpha }+1\right) ^{2}\left\vert A\right\vert _{\sigma }\ .
\notag
\end{eqnarray}

Indeed, in this last inequality (\ref{last inequ}), we have used first the
testing condition, which applies since the collection of $\mathcal{G}$%
-dyadic intervals 
\begin{equation*}
\mathcal{R}\equiv \dbigcup\limits_{I\in \mathcal{C}_{A}}\dbigcup\limits_{I^{%
\prime }\in \mathfrak{C}_{\limfunc{broken}}}\dbigcup\limits_{J\in \mathcal{C}%
_{A}^{\mathcal{G},\limfunc{nearby}}\left( I\right)
}\dbigcup\limits_{J^{\prime }\in \mathfrak{C}\left( J\right)
}\dbigcup\limits_{M\in \mathcal{M}_{n}^{\ast }\left( I^{\prime },J^{\prime
}\right) }\left\{ M_{\limfunc{in}}^{\limfunc{left}},M_{\limfunc{in}}^{%
\limfunc{right}}\right\}
\end{equation*}%
has bounded overlap counting repetitions of intervals. Indeed, given an
interval $L\in \mathcal{G}$, there are only a bounded, say $B$, number of
pairs $\left( I^{\prime },J^{\prime }\right) $, comparable in both scale and
position, for which $\mathcal{M}_{n}^{\ast }\left( I^{\prime },J^{\prime
}\right) $ contains an interval $M$ with $M_{\limfunc{in}}^{\limfunc{left}}$
or $M_{\limfunc{in}}^{\limfunc{right}}$ equal to $L$. Thus any tower of such
intervals $M_{\limfunc{in}}^{\limfunc{left}}$ or $M_{\limfunc{in}}^{\limfunc{%
right}}$, that contains a fixed point $x\in \mathbb{R}$, has at most $B$
intervals counting repetitions.

Next we used the energy condition in (\ref{last inequ}), which applies since
if $\mathcal{R}$, considered now without repetitions, has bounded overlap $B$%
, then $\mathcal{R}$ can be decomposed as $B$ pairwise disjoint families $%
\left\{ \mathcal{R}_{i}\right\} _{i=1}^{B}$. Indeed, since all of the
intervals lie in the dyadic grid $\mathcal{G}$ and are contained in a fixed
interval $A$, the family $\mathcal{R}_{1}$ of maximal intervals in $\mathcal{%
R}$ are pairwise disjoint, and after removing them, the remaining collection
of intervals $\mathcal{R}\setminus \mathcal{R}_{1}$ has bounded overlap $B-1$%
. Let $\mathcal{R}_{2}$ be the family of maximal dyadic intervals in $%
\mathcal{R}\setminus \mathcal{R}_{1}$ and continue until all the intervals
are exhausted after removing $\mathcal{R}_{R}$.

The inequality (\ref{prelim corona broken}) is a suitable estimate since%
\begin{equation*}
\sum_{A\in \mathcal{A}}\sqrt{\left\vert A\right\vert _{\sigma }\left( \frac{1%
}{\left\vert A\right\vert _{\sigma }}\int_{A}\left\vert f\right\vert d\sigma
\right) ^{2}}\left\Vert \mathsf{P}_{\mathcal{C}_{A}^{\mathcal{G},\limfunc{%
nearby}}}^{\omega }g\right\Vert _{L^{2}\left( \sigma \right) }^{\bigstar
}\lesssim \left\Vert f\right\Vert _{L^{2}\left( \sigma \right) }\left\Vert
g\right\Vert _{L^{2}\left( \sigma \right) }
\end{equation*}%
by quasiorthogonality and the frame inequalities in Appendix A, (\ref{Car
embed}) and (\ref{low frame}), together with the bounded overlap of the
`nearby' coronas $\left\{ \mathcal{C}_{A}^{\mathcal{G},\limfunc{nearby}%
}\right\} _{A\in \mathcal{A}}$.

Recall that after an initial application of random surgery, we reduced the
proof of Lemma \ref{nearby form} to establishing inequality (\ref{after prob}%
), in which $P_{\left( I,J\right) }\left( K,K\right) =\left\{ K,K\right\} $
in the notation used in (\ref{K iterated}). Now putting all of the above
estimates (\ref{C est}), (\ref{A est}), (\ref{prelim corona natural}) and (%
\ref{prelim corona broken}) together with (\ref{K iterated}) establishes
probabilistic control of the sum of all the inner products $\left\{
K,K\right\} $ taken over appropriate intervals $K$, yielding (\ref{after
prob}) as required if we choose $\lambda $ and $\eta _{0}$ sufficiently
small:%
\begin{eqnarray*}
&&\boldsymbol{E}_{\Omega }^{\mathcal{D}}\boldsymbol{E}_{\Omega }^{\mathcal{G}%
}\sum_{I\in \mathcal{D}}\sum_{\substack{ J\in \mathcal{G}:\ 2^{-\mathbf{r}%
}\ell \left( I\right) <\ell \left( J\right) \leq \ell \left( I\right)  \\ %
d\left( J,I\right) \leq 2\ell \left( J\right) ^{\varepsilon }\ell \left(
I\right) ^{1-\varepsilon }}}\left\vert \sum_{\substack{ I^{\prime }\in 
\mathfrak{C}\left( I\right) \text{ and }J^{\prime }\in \mathfrak{C}\left(
J\right)  \\ K\in \mathcal{K}\left( I^{\prime },J^{\prime }\right) }}%
E_{I^{\prime }}^{\sigma }\left( \widehat{\square }_{I}^{\sigma ,\flat ,%
\mathbf{b}}f\right) \ P_{\left( I,J\right) }\left( K,K\right) \ E_{J^{\prime
}}^{\omega }\left( \widehat{\square }_{J}^{\omega ,\flat ,\mathbf{b}^{\ast
}}g\right) \right\vert \\
&\leq &\boldsymbol{E}_{\Omega }^{\mathcal{D}}\boldsymbol{E}_{\Omega }^{%
\mathcal{G}}\sum_{A\in \mathcal{A}}\sum_{I\in \mathcal{C}_{A}\text{ and }%
J\in \mathcal{C}_{A}^{\mathcal{G},\limfunc{nearby}}\left( I\right) }\sum 
_{\substack{ I^{\prime }\in \mathfrak{C}\left( I\right) \text{ and }%
J^{\prime }\in \mathfrak{C}\left( J\right)  \\ K\in \mathcal{K}\left(
I^{\prime },J^{\prime }\right) }}\left\vert E_{I^{\prime }}^{\sigma }\left( 
\widehat{\square }_{I}^{\sigma ,\flat ,\mathbf{b}}f\right) \right\vert \
\left\vert A_{\left( I^{\prime },J^{\prime }\right) }\left( K\right)
\right\vert \ E_{J^{\prime }}^{\omega }\left( \widehat{\square }_{J}^{\omega
,\flat ,\mathbf{b}^{\ast }}g\right) \\
&&+\boldsymbol{E}_{\Omega }^{\mathcal{D}}\boldsymbol{E}_{\Omega }^{\mathcal{G%
}}\sum_{A\in \mathcal{A}}\sum_{I\in \mathcal{C}_{A}\text{ and }J\in \mathcal{%
C}_{A}^{\mathcal{G},\limfunc{nearby}}\left( I\right) }\sum_{\substack{ %
I^{\prime }\in \mathfrak{C}\left( I\right) \text{ and }J^{\prime }\in 
\mathfrak{C}\left( J\right)  \\ K\in \mathcal{K}\left( I^{\prime },J^{\prime
}\right) }}\left\vert E_{I^{\prime }}^{\sigma }\left( \widehat{\square }%
_{I}^{\sigma ,\flat ,\mathbf{b}}f\right) \right\vert \ \left( \left\vert
B_{\left( I^{\prime },J^{\prime }\right) }\left( K\right) \right\vert
+\left\vert C_{\left( I^{\prime },J^{\prime }\right) }\left( K\right)
\right\vert \right) \ \left\vert E_{J^{\prime }}^{\omega }\left( \widehat{%
\square }_{J}^{\omega ,\flat ,\mathbf{b}^{\ast }}g\right) \right\vert \\
&\lesssim &\left( C_{\theta }\mathcal{NTV}_{\alpha }+\sqrt{\theta }\mathfrak{%
N}_{T^{\alpha }}\right) \left\Vert f\right\Vert _{L^{2}\left( \sigma \right)
}\left\Vert g\right\Vert _{L^{2}\left( \omega \right) }\ .
\end{eqnarray*}%
This completes the proof of Lemma \ref{nearby form}.

\section{Main below form\label{Sec Main below}}

Now we turn to controlling the main below form (\ref{def Theta 2 good})%
\footnote{%
While it remains the case with $Tb$ arguments that this form is the most
challenging, the nearby form also poses great difficulties with $Tb$
arguments, especially in contrast with $T1$ arguments, where the nearby form
was handled almost trivially.},%
\begin{equation*}
\Theta _{2}^{\limfunc{good}}\left( f,g\right) =\sum_{I\in \mathcal{D}%
}\sum_{J^{\maltese }\subsetneqq I:\ \ell \left( J\right) \leq 2^{-\mathbf{r}%
}\ell \left( I\right) }\int \left( T_{\sigma }^{\alpha }\square _{I}^{\sigma
,\mathbf{b}}f\right) \square _{J}^{\omega ,\mathbf{b}^{\ast }}gd\omega .
\end{equation*}%
where we recall that $\mathbf{r}$ is the goodness parameter fixed in (\ref%
{choice of r}) given $0<\varepsilon <\frac{1}{2}$. The corresponding main
below form in \cite{SaShUr6} was denoted $\mathsf{B}_{\Subset _{\mathbf{r}%
}}\left( f,g\right) =\mathsf{B}_{\Subset _{\mathbf{r},\varepsilon }}\left(
f,g\right) $. However, there are significant differences between the forms $%
\Theta _{2}^{\limfunc{good}}\left( f,g\right) $ and $\mathsf{B}_{\Subset _{%
\mathbf{r}}}\left( f,g\right) $. In \cite{SaShUr6}, the Haar martingale
averages $\bigtriangleup _{I}^{\sigma }$ and $\bigtriangleup _{J}^{\omega }$
in $\mathsf{B}_{\Subset _{\mathbf{r}}}\left( f,g\right) $ are orthogonal
projections, and the intervals $I$ and $J$ that appear in $\mathsf{B}%
_{\Subset _{\mathbf{r}}}\left( f,g\right) $ are all good in the traditional
sense. Here, on the other hand, the dual martingale averages $\square
_{I}^{\sigma ,\mathbf{b}}$ and $\square _{J}^{\omega ,\mathbf{b}^{\ast }}$
in $\Theta _{2}^{\limfunc{good}}\left( f,g\right) $ are no longer orthogonal
projections, and while the intervals $J$ paired with $I$ remain $\varepsilon 
$-good inside $I$ and beyond, the lack of orthogonal projections is
compensated by the fact that the collections of intervals $I$ associated
with any fixed $J$ in $\Theta _{2}^{\limfunc{good}}\left( f,g\right) $ are
tree-connected. Nevertheless, in order to efficiently import the methods for
controlling $\mathsf{B}_{\Subset _{\mathbf{r}}}\left( f,g\right) $ from \cite%
{SaShUr6}, we will relabel the main below form $\Theta _{2}^{\limfunc{good}%
}\left( f,g\right) $ as $\mathsf{B}_{\Subset _{\mathbf{r}}}\left( f,g\right)
=\mathsf{B}_{\Subset _{\mathbf{r},\varepsilon }}\left( f,g\right) $ from now
on, keeping in mind the aforementioned differences.

To control $\Theta _{2}^{\limfunc{good}}\left( f,g\right) =\mathsf{B}%
_{\Subset _{\mathbf{r}}}\left( f,g\right) $ we first perform the \emph{%
canonical corona splitting} of $\mathsf{B}_{\Subset _{\mathbf{r}}}\left(
f,g\right) $ into a diagonal form $\mathsf{T}_{\limfunc{diagonal}}\left(
f,g\right) $ and a far below form $\mathsf{T}_{\limfunc{far}\limfunc{below}%
}\left( f,g\right) $ as in \cite{SaShUr6}. This \emph{canonical splitting}
of the form $\mathsf{B}_{\Subset _{\mathbf{r}}}\left( f,g\right) $ involves
the corona pseudoprojections $\mathsf{P}_{\mathcal{C}_{A}^{\mathcal{D}%
}}^{\sigma ,\mathbf{b}}$ acting on $f$ and the \emph{shifted} corona
pseudoprojections $\mathsf{P}_{\mathcal{C}_{B}^{\mathcal{G},\limfunc{shift}%
}}^{\omega ,\mathbf{b}^{\ast }}$ acting on $g$, where $B$ is a stopping
interval in $\mathcal{A}$. The stopping intervals $\mathcal{B}$ constructed
relative to $g\in L^{2}\left( \omega \right) $ play no role in the analysis
here, except to guarantee that the frame and weak Riesz inequalities hold
for $g$ and the pseudprojections $\left\{ \square _{J}^{\omega ,\mathbf{b}%
^{\ast }}g\right\} _{J\in \mathcal{G}}$ and Carleson averaging operators $%
\left\{ \nabla _{J,\mathcal{G}}^{\omega }g\right\} _{J\in \mathcal{G}}$.
Here the shifted corona $\mathcal{C}_{B}^{\mathcal{G},\limfunc{shift}}$ is
defined for $B\in \mathcal{A}$ - and \textbf{not} for $B\in \mathcal{B}$ -
to include those intervals $J\in \mathcal{G}$ such that $J^{\maltese }\in 
\mathcal{C}_{B}^{\mathcal{D}}$.

\begin{definition}
\label{shifted corona}For $B\in \mathcal{A}$ we define the shifted $\mathcal{%
G}$-corona by%
\begin{equation*}
\mathcal{C}_{B}^{\mathcal{G},\limfunc{shift}}=\left\{ J\in \mathcal{G}%
:J^{\maltese }\in \mathcal{C}_{B}^{\mathcal{D}}\right\} .
\end{equation*}
\end{definition}

The Carleson averaging operator $\nabla _{J,\mathcal{G}}^{\omega }$ is taken
over the `broken' children of $J$ which depend on the grid $\mathcal{G}$
(see (\ref{Carleson avg op}) in Appendix A). We will use repeatedly the fact
that the shifted coronas $\mathcal{C}_{B}^{\mathcal{G},\limfunc{shift}}$ are
pairwise disjoint in $B$:%
\begin{equation}
\sum_{B\in \mathcal{A}}\mathbf{1}_{\mathcal{C}_{B}^{\mathcal{G},\limfunc{%
shift}}}\left( J\right) \leq \mathbf{1},\ \ \ \ \ J\in \mathcal{G}.
\label{tau overlap}
\end{equation}%
It is convenient at this point to introduce the following shorthand notation:%
\begin{equation}
\left\langle T_{\sigma }^{\alpha }\left( \mathsf{P}_{\mathcal{C}_{A}^{%
\mathcal{D}}}^{\sigma ,\mathbf{b}}f\right) ,\mathsf{P}_{\mathcal{C}_{B}^{%
\mathcal{G},\limfunc{shift}}}^{\omega ,\mathbf{b}^{\ast }}g\right\rangle
_{\omega }^{\Subset _{\mathbf{r},\varepsilon }}\equiv \sum_{\substack{ I\in 
\mathcal{C}_{A}^{\mathcal{D}}\text{ and }J\in \mathcal{C}_{B}^{\mathcal{G},%
\limfunc{shift}}:\ J^{\maltese }\subsetneqq I  \\ \ell \left( J\right) \leq
2^{-\mathbf{r}}\ell \left( I\right) }}\left\langle T_{\sigma }^{\alpha
}\left( \square _{I}^{\sigma ,\mathbf{b}}f\right) ,\square _{J}^{\omega ,%
\mathbf{b}^{\ast }}g\right\rangle _{\omega }\ .  \label{def shorthand}
\end{equation}

\begin{description}
\item[Caution] One musn't assume, from the notation on the left hand side
above, that the function $T_{\sigma }^{\alpha }\left( \mathsf{P}_{\mathcal{C}%
_{A}^{\mathcal{D}}}^{\sigma ,\mathbf{b}}f\right) $ is simply integrated
against the function $\mathsf{P}_{\mathcal{C}_{B}^{\mathcal{G},\limfunc{shift%
}}}^{\omega ,\mathbf{b}^{\ast }}g$. Indeed, the sum on the right hand side
is taken over pairs $\left( I,J\right) $ such that $J^{\maltese }\in 
\mathcal{C}_{B}^{\mathcal{D}}$ and $J^{\maltese }\subsetneqq I$ and $\ell
\left( J\right) \leq 2^{-\mathbf{r}}\ell \left( I\right) $.
\end{description}

Here is the relevant portion of the brief schematic diagram (\ref{schematic}%
) of the decompositions, with bounds in $\fbox{}$, used in the next
subsections:%
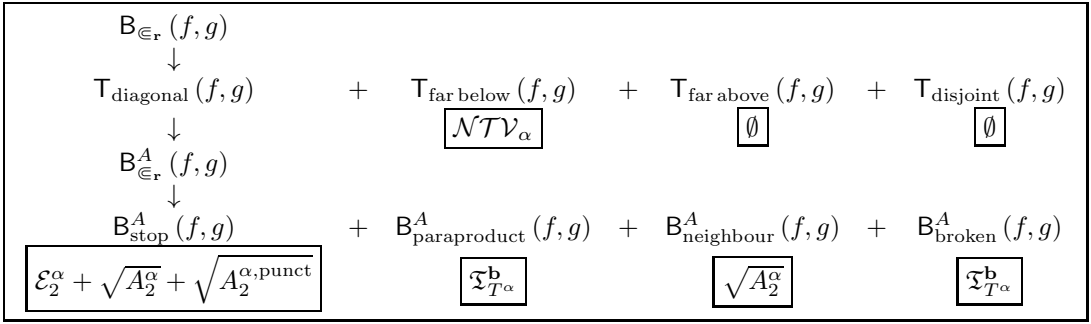
\begin{equation*}
\fbox{$%
\begin{array}{ccccccc}
\mathsf{B}_{\Subset _{\mathbf{r}}}\left( f,g\right) &  &  &  &  &  &  \\ 
\downarrow &  &  &  &  &  &  \\ 
\mathsf{T}_{\limfunc{diagonal}}\left( f,g\right) & + & \mathsf{T}_{\limfunc{%
far}\limfunc{below}}\left( f,g\right) & + & \mathsf{T}_{\limfunc{far}%
\limfunc{above}}\left( f,g\right) & + & \mathsf{T}_{\limfunc{disjoint}%
}\left( f,g\right) \\ 
\downarrow &  & \fbox{$\mathcal{NTV}_{\alpha }$} &  & \fbox{$\emptyset $} & 
& \fbox{$\emptyset $} \\ 
\mathsf{B}_{\Subset _{\mathbf{r}}}^{A}\left( f,g\right) &  &  &  &  &  &  \\ 
\downarrow &  &  &  &  &  &  \\ 
\mathsf{B}_{\limfunc{stop}}^{A}\left( f,g\right) & + & \mathsf{B}_{\limfunc{%
paraproduct}}^{A}\left( f,g\right) & + & \mathsf{B}_{\limfunc{neighbour}%
}^{A}\left( f,g\right) & + & \mathsf{B}_{\limfunc{broken}}^{A}\left(
f,g\right) \\ 
\fbox{$\mathcal{E}_{2}^{\alpha }+\sqrt{A_{2}^{\alpha }}+\sqrt{A_{2}^{\alpha ,%
\limfunc{punct}}}$} &  & \fbox{$\mathfrak{T}_{T^{\alpha }}^{\mathbf{b}}$} & 
& \fbox{$\sqrt{A_{2}^{\alpha }}$} &  & \fbox{$\mathfrak{T}_{T^{\alpha }}^{%
\mathbf{b}}$}%
\end{array}%
$}
\end{equation*}

\subsection{The canonical splitting and local below forms}

We begin with an informal description of decompositions and estimates. The
canonical splitting is determined by the coronas $\mathcal{C}_{A}^{\mathcal{D%
}}$ for $A\in \mathcal{A}$ - note that the stopping times $\mathcal{B}$ play
no explicit role in the canonical splitting of the below form, other than to
guarantee weak Riesz inequalities for $\square _{J}^{\omega ,\mathbf{b}%
^{\ast }}$ and $\nabla _{J,\mathcal{G}}^{\omega }$:%
\begin{eqnarray}
&&\mathsf{B}_{\Subset _{\mathbf{r},\varepsilon }}\left( f,g\right)
\label{parallel corona decomp'} \\
&=&\sum_{A,B\in \mathcal{A}}\left\langle T_{\sigma }^{\alpha }\left( \mathsf{%
P}_{\mathcal{C}_{A}}^{\sigma ,\mathbf{b}}f\right) ,\mathsf{P}_{\mathcal{C}%
_{B}^{\mathcal{G},\limfunc{shift}}}^{\omega ,\mathbf{b}^{\ast
}}g\right\rangle _{\omega }^{\Subset _{\mathbf{r},\varepsilon }}  \notag \\
&=&\sum_{A\in \mathcal{A}}\left\langle T_{\sigma }^{\alpha }\left( \mathsf{P}%
_{\mathcal{C}_{A}}^{\sigma ,\mathbf{b}}f\right) ,\mathsf{P}_{\mathcal{C}%
_{A}^{\mathcal{G},\limfunc{shift}}}^{\omega ,\mathbf{b}^{\ast
}}g\right\rangle _{\omega }^{\Subset _{\mathbf{r},\varepsilon }}+\sum 
_{\substack{ A,B\in \mathcal{A}  \\ B\subsetneqq A}}\left\langle T_{\sigma
}^{\alpha }\left( \mathsf{P}_{\mathcal{C}_{A}}^{\sigma ,\mathbf{b}}f\right) ,%
\mathsf{P}_{\mathcal{C}_{B}^{\mathcal{G},\limfunc{shift}}}^{\omega ,\mathbf{b%
}^{\ast }}g\right\rangle _{\omega }^{\Subset _{\mathbf{r},\varepsilon }} 
\notag \\
&&+\sum_{\substack{ A,B\in \mathcal{A}  \\ B\supsetneqq A}}\left\langle
T_{\sigma }^{\alpha }\left( \mathsf{P}_{\mathcal{C}_{A}}^{\sigma ,\mathbf{b}%
}f\right) ,\mathsf{P}_{\mathcal{C}_{B}^{\mathcal{G},\limfunc{shift}%
}}^{\omega ,\mathbf{b}^{\ast }}g\right\rangle _{\omega }^{\Subset _{\mathbf{r%
},\varepsilon }}+\sum_{\substack{ A,B\in \mathcal{A}  \\ A\cap B=\emptyset }}%
\left\langle T_{\sigma }^{\alpha }\left( \mathsf{P}_{\mathcal{C}%
_{A}}^{\sigma ,\mathbf{b}}f\right) ,\mathsf{P}_{\mathcal{C}_{B}^{\mathcal{G},%
\limfunc{shift}}}^{\omega ,\mathbf{b}^{\ast }}g\right\rangle _{\omega
}^{\Subset _{\mathbf{r},\varepsilon }}  \notag \\
&\equiv &\mathsf{T}_{\limfunc{diagonal}}\left( f,g\right) +\mathsf{T}_{%
\limfunc{far}\limfunc{below}}\left( f,g\right) +\mathsf{T}_{\limfunc{far}%
\limfunc{above}}\left( f,g\right) +\mathsf{T}_{\limfunc{disjoint}}\left(
f,g\right) .  \notag
\end{eqnarray}%
Now the final two terms $\mathsf{T}_{\limfunc{far}\limfunc{above}}\left(
f,g\right) $ and $\mathsf{T}_{\limfunc{disjoint}}\left( f,g\right) $ each
vanish since there are no pairs $\left( I,J\right) \in \mathcal{C}_{A}^{%
\mathcal{D}}\times \mathcal{C}_{B}^{\mathcal{G},\limfunc{shift}}$ with both (%
\textbf{i}) $J^{\maltese }\subsetneqq I$ and (\textbf{ii}) either $%
B\supsetneqq A$ or $B\cap A=\emptyset $. The far below form $\mathsf{T}_{%
\limfunc{far}\limfunc{below}}\left( f,g\right) $ requires functional energy,
which we discuss in a moment.

Next we follow this splitting by a further decomposition of the diagonal
form into local below forms $\mathsf{B}_{\Subset _{\mathbf{r}}}^{A}\left(
f,g\right) $ given by the individual corona pieces 
\begin{equation}
\mathsf{B}_{\Subset _{\mathbf{r},\varepsilon }}^{A}\left( f,g\right) \equiv
\left\langle T_{\sigma }^{\alpha }\left( \mathsf{P}_{\mathcal{C}%
_{A}}^{\sigma ,\mathbf{b}}f\right) ,\mathsf{P}_{\mathcal{C}_{A}^{\mathcal{G},%
\limfunc{shift}}}^{\omega ,\mathbf{b}^{\ast }}g\right\rangle _{\omega
}^{\Subset _{\mathbf{r},\varepsilon }}\ ,  \label{def local}
\end{equation}%
and prove the following estimate where $\mathcal{NTV}_{\alpha }$ is defined
in (\ref{def NTV}):%
\begin{equation}
\left\vert \mathsf{B}_{\Subset _{\mathbf{r},\varepsilon }}^{A}\left(
f,g\right) \right\vert \lesssim \mathcal{NTV}_{\alpha }\ \left( \alpha _{%
\mathcal{A}}\left( A\right) \sqrt{\left\vert A\right\vert _{\sigma }}%
+\left\Vert \mathsf{P}_{\mathcal{C}_{A}}^{\sigma ,\mathbf{b}}f\right\Vert
_{L^{2}\left( \sigma \right) }^{\bigstar }\right) \ \left\Vert \mathsf{P}_{%
\mathcal{C}_{A}^{\mathcal{G},\limfunc{shift}}}^{\omega ,\mathbf{b}^{\ast
}}g\right\Vert _{L^{2}\left( \omega \right) }^{\bigstar }\ .
\label{local est}
\end{equation}%
This reduces matters to the local forms since we then have from
Cauchy-Schwarz that%
\begin{eqnarray*}
&&\sum_{A\in \mathcal{A}}\left\vert \mathsf{B}_{\Subset _{\mathbf{r}%
,\varepsilon }}^{A}\left( f,g\right) \right\vert \lesssim \mathcal{NTV}%
_{\alpha }\ \left( \sum_{A\in \mathcal{A}}\alpha _{\mathcal{A}}\left(
A\right) ^{2}\left\vert A\right\vert _{\sigma }+\left\Vert \mathsf{P}_{%
\mathcal{C}_{A}^{\mathcal{D}}}^{\sigma ,\mathbf{b}}f\right\Vert
_{L^{2}\left( \sigma \right) }^{\bigstar 2}\right) ^{\frac{1}{2}}\left(
\sum_{A\in \mathcal{A}}\left\Vert \mathsf{P}_{\mathcal{C}_{A}^{\mathcal{G},%
\limfunc{shift}}}^{\omega ,\mathbf{b}^{\ast }}g\right\Vert _{L^{2}\left(
\omega \right) }^{\bigstar 2}\right) ^{\frac{1}{2}} \\
&&\ \ \ \ \ \ \ \ \ \ \ \ \ \ \ \ \ \ \ \ \ \ \ \ \ \ \ \ \ \ \lesssim 
\mathcal{NTV}_{\alpha }\ \left\Vert f\right\Vert _{L^{2}\left( \sigma
\right) }\left\Vert g\right\Vert _{L^{2}\left( \omega \right) }\ ,
\end{eqnarray*}%
by the lower frame inequalities $\sum_{A\in \mathcal{A}}\left\Vert \mathsf{P}%
_{\mathcal{C}_{A}}^{\sigma ,\mathbf{b}}f\right\Vert _{L^{2}\left( \sigma
\right) }^{\bigstar 2}\lesssim \left\Vert f\right\Vert _{L^{2}\left( \sigma
\right) }^{2}$ and $\sum_{A\in \mathcal{A}}\left\Vert \mathsf{P}_{\mathcal{C}%
_{A}^{\mathcal{G},\limfunc{shift}}}^{\omega ,\mathbf{b}^{\ast }}g\right\Vert
_{L^{2}\left( \omega \right) }^{\bigstar 2}\lesssim \left\Vert g\right\Vert
_{L^{2}\left( \omega \right) }^{2}$, using quasi-orthogonality $\sum_{A\in 
\mathcal{A}}\alpha _{\mathcal{A}}\left( f\right) ^{2}\left\vert A\right\vert
_{\sigma }\lesssim \left\Vert f\right\Vert _{L^{2}\left( \sigma \right)
}^{2} $ in the stopping intervals $\mathcal{A}$, and the pairwise
disjointedness of the shifted coronas $\mathcal{C}_{A}^{\mathcal{G},\limfunc{%
shift}}$: 
\begin{equation*}
\sum_{A\in \mathcal{A}}\mathbf{1}_{\mathcal{C}_{A}^{\mathcal{G},\limfunc{%
shift}}}\leq \mathbf{1}_{\mathcal{D}}.
\end{equation*}%
From now on we will often write $\mathcal{C}_{A}$ in place of $\mathcal{C}%
_{A}^{\mathcal{D}}$ when no confusion is possible.

Finally, the local forms $\mathsf{B}_{\Subset _{\mathbf{r},\varepsilon
}}^{A}\left( f,g\right) $ are decomposed into stopping forms $\mathsf{B}_{%
\limfunc{stop}}^{A}\left( f,g\right) $, paraproduct forms $\mathsf{B}_{%
\limfunc{paraproduct}}^{A}\left( f,g\right) $, neighbour forms $\mathsf{B}_{%
\limfunc{neighbour}}^{A}\left( f,g\right) $ forms and broken forms $\mathsf{B%
}_{\limfunc{broken}}^{A}\left( f,g\right) $. The paraproduct and neighbour
terms are handled as in \cite{SaShUr6}, which in turn follows the treatment
originating in \cite{NTV3}, and the broken form is handled with Carleson
measure methods, leaving only the stopping form $\mathsf{B}_{\limfunc{stop}%
}^{A}\left( f,g\right) $ to be bounded, which we treat in the next section
below, Section \ref{Sec stop}, by refining the bottom/up stopping time in
the argument of M. Lacey in \cite{Lac} with an additional top/down
`indented' corona construction to handle weak goodness.

However, in order to complete the required bounds of the above forms into
which the below form $\mathsf{B}_{\Subset _{\mathbf{r}}}\left( f,g\right) $
was decomposed, we need functional energy for the far below form. Recall
that the vector-valued function $\mathbf{b}$ in the accretive coronas
`breaks' only at a collection of intervals satisfying a Carleson condition.

\begin{definition}
\label{def Whitney}Define the \emph{Whitney} subintervals $\mathcal{W}\left(
F\right) $ of an interval $F\in \mathcal{D}$ to consist of the \emph{maximal}
dyadic $\mathcal{D}$-subintervals of a $\mathcal{D}$-interval $F$ that have
their triples contained in $F$.
\end{definition}

See (\ref{def pseudo rest}) in Appendix B below for more detail on this and
the remaining terms in (\ref{e.funcEnergy n}) below.

\begin{definition}
\label{functional energy n}Let $\mathfrak{F}_{\alpha }=\mathfrak{F}_{\alpha
}\left( \mathcal{D},\mathcal{G}\right) =\mathfrak{F}_{\alpha }^{\mathbf{b}%
^{\ast }}\left( \mathcal{D},\mathcal{G}\right) $ be the smallest constant in
the `\textbf{f}unctional energy' inequality below, holding for all $h\in
L^{2}\left( \sigma \right) $ and all $\sigma $-Carleson collections $%
\mathcal{F}\subset \mathcal{D}$ with Carleson norm $C_{\mathcal{F}}$ bounded
by a fixed constant $C$: 
\begin{equation}
\sum_{F\in \mathcal{F}}\sum_{M\in \mathcal{W}\left( F\right) }\left( \frac{%
\mathrm{P}^{\alpha }\left( M,h\sigma \right) }{\left\vert M\right\vert }%
\right) ^{2}\left\Vert \mathsf{Q}_{\mathcal{C}_{F}^{\mathcal{G},\limfunc{%
shift}};M}^{\omega ,\mathbf{b}^{\ast }}x\right\Vert _{L^{2}\left( \omega
\right) }^{\spadesuit 2}\leq \left( \mathfrak{F}_{\alpha }^{\mathbf{b}^{\ast
}}\left( \mathcal{D},\mathcal{G}\right) \right) ^{2}\left\Vert h\right\Vert
_{L^{2}\left( \sigma \right) }^{2}\,,  \label{e.funcEnergy n}
\end{equation}
\end{definition}

The main ingredient used in reducing control of the below form $\mathsf{B}%
_{\Subset _{\mathbf{r}}}\left( f,g\right) $ to control of the functional
energy $\mathfrak{F}_{\alpha }=\mathfrak{F}_{\alpha }^{\mathbf{b}^{\ast
}}\left( \mathcal{D},\mathcal{G}\right) $ constant and the stopping form $%
\mathsf{B}_{\limfunc{stop}}^{A}\left( f,g\right) $, is the Intertwining
Proposition from \cite{SaShUr7} and/or \cite{SaShUr6}. The control of the
functional energy condition by the energy and Muckenhoupt conditions must
also be adapted in light of the $\infty $-strongly accretive function $%
\mathbf{b}$ that only `breaks' at a collection of intervals satisfying a
Carleson condition, but this poses no real difficulties. The fact that the
usual Haar bases are orthonormal is here replaced by the weaker condition
that the corresponding broken dual martingale `bases' are merely weak frames
satisfying certain\ weak lower and weak upper Riesz inequalities, but again
this poses no real difference in the arguments. Finally, the fact that
goodness for $J$ has been replaced with weak goodness, namely $J^{\maltese
}\subsetneqq I$ whenever the pair $\left( I,J\right) $ occurs in a sum,
forces the use of a Whitney decomposition $\mathcal{W}$ of intervals instead
of the deeply embedded decomposition $\mathcal{M}_{\left( \mathbf{\rho }%
,\varepsilon \right) -\limfunc{deep}}$ used in \cite{SaShUr7}.

We then use the paraproduct / neighbour / broken / stopping splitting
mentioned above to reduce boundedness of $\mathsf{B}_{\Subset _{\mathbf{r}%
,\varepsilon }}^{A}\left( f,g\right) $ to boundedness of the associated
stopping form 
\begin{equation}
\mathsf{B}_{\limfunc{stop}}^{A}\left( f,g\right) \equiv \sum_{\substack{ %
I\in \mathcal{C}_{A}\text{ and }J\in \mathcal{C}_{A}^{\mathcal{G},\limfunc{%
shift}}  \\ J^{\maltese }\subset I\text{, }J^{\maltese }\neq I\text{ and }%
\ell \left( J\right) \leq 2^{-\mathbf{r}}\ell \left( I\right) }}%
E_{I_{J}}^{\sigma }\left( \widehat{\square }_{I}^{\sigma ,\flat ,\mathbf{b}%
}f\right) \left\langle T_{\sigma }^{\alpha }b_{A}\mathbf{1}_{A\setminus
I_{J}},\square _{J}^{\omega ,\mathbf{b}^{\ast }}g\right\rangle _{\omega }\ ,
\label{bounded stopping form}
\end{equation}%
where $f$ is supported in the interval $A$ and its expectations $\mathbb{E}%
_{I}^{\sigma }\left\vert f\right\vert $ are bounded by $\alpha _{\mathcal{A}%
}\left( A\right) $ for $I\in \mathcal{C}_{A}$, the dual martingale support
of $f$ is contained in the corona $\mathcal{C}_{A}^{\sigma }$, and the dual
martingale support of $g$\ is contained in $\mathcal{C}_{A}^{\mathcal{G},%
\limfunc{shift}}$, and where $I_{J}$ is the $\mathcal{D}$-child of $I$ that
contains $J$.

\subsection{Diagonal and far below forms}

Now we turn to estimating the \emph{diagonal term} $\mathsf{T}_{\limfunc{%
diagonal}}\left( f,g\right) $ and the \emph{far below} term $\mathsf{T}_{%
\limfunc{far}\limfunc{below}}\left( f,g\right) $, where in \cite{SaShUr7}
and/or \cite{SaShUr6}, the far below terms were bounded using the
Intertwining Proposition and the control of functional energy condition by
the energy and Muckenhoupt conditions, but of course under the restriction
there that the intervals $J$ were good. Here we write%
\begin{eqnarray}
&&\mathsf{T}_{\limfunc{far}\limfunc{below}}\left( f,g\right) =\sum 
_{\substack{ A,B\in \mathcal{A}  \\ B\subsetneqq A}}\sum_{\substack{ I\in 
\mathcal{C}_{A}\text{ and }J\in \mathcal{C}_{B}^{\mathcal{G},\limfunc{shift}%
}  \\ J^{\maltese }\subsetneqq I\text{ and }\ell \left( J\right) \leq 2^{-%
\mathbf{r}}\ell \left( I\right) }}\left\langle T_{\sigma }^{\alpha }\left(
\square _{I}^{\sigma ,\mathbf{b}}f\right) ,\square _{J}^{\omega ,\mathbf{b}%
^{\ast }}g\right\rangle _{\omega }  \label{write} \\
&=&\sum_{B\in \mathcal{A}}\sum_{I\in \mathcal{D}:\ B\subsetneqq
I}\left\langle T_{\sigma }^{\alpha }\left( \square _{I}^{\sigma ,\mathbf{b}%
}f\right) ,\sum_{J\in \mathcal{C}_{B}^{\mathcal{G},\limfunc{shift}}}\square
_{J}^{\omega ,\mathbf{b}^{\ast }}g\right\rangle _{\omega }-\sum_{B\in 
\mathcal{A}}\sum_{I\in \mathcal{D}:\ B\subsetneqq I}\left\langle T_{\sigma
}^{\alpha }\left( \square _{I}^{\sigma ,\mathbf{b}}f\right) ,\sum_{\substack{
J\in \mathcal{C}_{B}^{\mathcal{G},\limfunc{shift}}  \\ \ell \left( J\right)
>2^{-\mathbf{r}}\ell \left( I\right) }}\square _{J}^{\omega ,\mathbf{b}%
^{\ast }}g\right\rangle _{\omega }  \notag \\
&=&\sum_{B\in \mathcal{A}}\sum_{I\in \mathcal{D}:\ B\subsetneqq
I}\left\langle T_{\sigma }^{\alpha }\left( \square _{I}^{\sigma ,\mathbf{b}%
}f\right) ,g_{B}\right\rangle _{\omega }-\sum_{B\in \mathcal{A}}\sum_{I\in 
\mathcal{D}:\ B\subsetneqq I}\left\langle T_{\sigma }^{\alpha }\left(
\square _{I}^{\sigma ,\mathbf{b}}f\right) ,\sum_{\substack{ J\in \mathcal{C}%
_{B}^{\mathcal{G},\limfunc{shift}}  \\ \ell \left( J\right) >2^{-\mathbf{r}%
}\ell \left( I\right) }}\square _{J}^{\omega ,\mathbf{b}^{\ast
}}g\right\rangle _{\omega };  \notag \\
&&\text{where }g_{B}=\sum_{J\in \mathcal{C}_{B}^{\mathcal{G},\limfunc{shift}%
}}\square _{J}^{\omega ,\mathbf{b}^{\ast }}g=\mathsf{P}_{\mathcal{C}_{F}^{%
\mathcal{G},\limfunc{shift}}}^{\omega ,\mathbf{b}^{\ast }}g\ ,  \notag
\end{eqnarray}%
since if $I\in \mathcal{C}_{A}$ and $J\in \mathcal{C}_{B}^{\mathcal{G},%
\limfunc{shift}}$, with $J^{\maltese }\subsetneqq I$ and $B\subsetneqq A$,
then we must have $B\subsetneqq I$. First, we note that expectation of the
second sum on the right hand side of (\ref{write}) is controlled by (\ref%
{delta near}) in Lemma \ref{nearby form}, i.e.%
\begin{eqnarray*}
&&\boldsymbol{E}_{\Omega }^{\mathcal{D}}\boldsymbol{E}_{\Omega }^{\mathcal{G}%
}\left\vert \sum_{B\in \mathcal{A}}\sum_{I\in \mathcal{D}:\ B\subsetneqq
I}\left\langle T_{\sigma }^{\alpha }\left( \square _{I}^{\sigma ,\mathbf{b}%
}f\right) ,\sum_{\substack{ J\in \mathcal{C}_{B}^{\mathcal{G},\limfunc{shift}%
}  \\ \ell \left( J\right) >2^{-\mathbf{r}}\ell \left( I\right) }}\square
_{J}^{\omega ,\mathbf{b}^{\ast }}g\right\rangle _{\omega }\right\vert \\
&\lesssim &\boldsymbol{E}_{\Omega }^{\mathcal{D}}\boldsymbol{E}_{\Omega }^{%
\mathcal{G}}\sum_{I\in \mathcal{D}}\sum_{\substack{ J\in \mathcal{G}:\ 2^{-%
\mathbf{r}}\ell \left( I\right) <\ell \left( J\right) \leq \ell \left(
I\right)  \\ d\left( J,I\right) \leq 2\ell \left( J\right) ^{\varepsilon
}\ell \left( I\right) ^{1-\varepsilon }}}\left\vert \left\langle T_{\sigma
}^{\alpha }\left( \square _{I}^{\sigma ,\mathbf{b}}f\right) ,\square
_{J}^{\omega ,\mathbf{b}^{\ast }}g\right\rangle _{\omega }\right\vert \\
&\lesssim &\left( C_{\theta }\mathcal{NTV}_{\alpha }+\sqrt{\theta }\mathfrak{%
N}_{T^{\alpha }}\right) \left\Vert f\right\Vert _{L^{2}\left( \sigma \right)
}\left\Vert g\right\Vert _{L^{2}\left( \omega \right) }\ .
\end{eqnarray*}%
Second, we note that the Intertwining Proposition \ref{strongly adapted},
which controls sums of the form%
\begin{equation*}
\sum_{F\in \mathcal{F}}\ \sum_{I:\ I\supsetneqq F}\ \left\langle T_{\sigma
}^{\alpha }\square _{I}^{\sigma ,\mathbf{b}}f,\mathsf{P}_{\mathcal{C}_{F}^{%
\mathcal{G},\limfunc{shift}}}^{\omega ,\mathbf{b}^{\ast }}g\right\rangle
_{\omega },
\end{equation*}%
can be applied to the first sum on the right hand side of (\ref{write}) to
show that it is bounded by $\left( \mathcal{NTV}_{\alpha }+\mathfrak{F}%
_{\alpha }^{\mathbf{b}^{\ast }}\left( \mathcal{D},\mathcal{G}\right) \right)
\left\Vert f\right\Vert _{L^{2}\left( \sigma \right) }\left\Vert
g\right\Vert _{L^{2}\left( \omega \right) }$, where the goodness parameter $%
\varepsilon >0$ is chosen sufficiently small. Then Proposition \ref{func
ener control}\ can be applied to show that $\mathfrak{F}_{\alpha }^{\mathbf{b%
}^{\ast }}\left( \mathcal{D},\mathcal{G}\right) \lesssim \mathfrak{A}%
_{2}^{\alpha }+\mathfrak{E}_{2}^{\alpha }$, which completes the proof that%
\begin{equation}
\left\vert \mathsf{T}_{\limfunc{far}\limfunc{below}}\left( f,g\right)
\right\vert \lesssim \mathcal{NTV}_{\alpha }\ \left\Vert f\right\Vert
_{L^{2}\left( \sigma \right) }\left\Vert g\right\Vert _{L^{2}\left( \omega
\right) }\ .  \label{far below bound}
\end{equation}

\subsection{Intertwining Proposition}

First we adapt the relevant definitions from \cite{SaShUr6}.

\begin{definition}
\label{sigma carleson n}A collection $\mathcal{F}$ of dyadic intervals is $%
\sigma $\emph{-Carleson} if%
\begin{equation*}
\sum_{F\in \mathcal{F}:\ F\subset S}\left\vert F\right\vert _{\sigma }\leq
C_{\mathcal{F}}\left\vert S\right\vert _{\sigma },\ \ \ \ \ S\in \mathcal{F}.
\end{equation*}%
The constant $C_{\mathcal{F}}$ is referred to as the Carleson norm of $%
\mathcal{F}$.
\end{definition}

\begin{definition}
\label{def shift}Let $\mathcal{F}$ be a collection of dyadic intervals in a
grid $\mathcal{D}$. Then for $F\in \mathcal{F}$, we define the shifted
corona $\mathcal{C}_{F}^{\mathcal{G},\limfunc{shift}}$ in analogy with
Definition \ref{shifted corona} by%
\begin{equation*}
\mathcal{C}_{F}^{\mathcal{G},\limfunc{shift}}\equiv \left\{ J\in \mathcal{G}%
:J^{\maltese }\in \mathcal{C}_{F}\right\} ,
\end{equation*}%
where $J^{\maltese }$ is defined in Definition \ref{def sharp cross}.
\end{definition}

Note that the collections $\mathcal{C}_{F}^{\mathcal{G},\limfunc{shift}}$
are pairwise disjoint in $F$. Let $\mathfrak{C}_{\mathcal{F}}\left( F\right) 
$ denote the set of $\mathcal{F}$-children of $F$. Given any collection $%
\mathcal{H}\subset \mathcal{G}$ of intervals, a family $\mathbf{b}^{\ast }$
of dual testing functions, and an arbitrary interval $K\in \mathcal{P}$, we
define the corresponding dual pseudoprojection $\mathsf{P}_{\mathcal{H}%
}^{\omega ,\mathbf{b}^{\ast }}$ and its localization $\mathsf{P}_{\mathcal{H}%
;K}^{\omega ,\mathbf{b}^{\ast }}$ to $K$ by%
\begin{equation}
\mathsf{Q}_{\mathcal{H}}^{\omega ,\mathbf{b}^{\ast }}=\sum_{H\in \mathcal{H}%
}\bigtriangleup _{H}^{\omega ,\mathbf{b}^{\ast }}\text{ and }\mathsf{Q}_{%
\mathcal{H};K}^{\omega ,\mathbf{b}^{\ast }}=\sum_{H\in \mathcal{H}:\
H\subset K}\bigtriangleup _{H}^{\omega ,\mathbf{b}^{\ast }}\ .
\label{def localization}
\end{equation}%
Recall from Definition \ref{functional energy n} that $\mathfrak{F}_{\alpha
}=\mathfrak{F}_{\alpha }\left( \mathcal{D},\mathcal{G}\right) =\mathfrak{F}%
_{\alpha }^{\mathbf{b}^{\ast }}\left( \mathcal{D},\mathcal{G}\right) $ is
the best constant in (\ref{e.funcEnergy n}), i.e. 
\begin{equation*}
\sum_{F\in \mathcal{F}}\sum_{M\in \mathcal{W}\left( F\right) }\left( \frac{%
\mathrm{P}^{\alpha }\left( M,h\sigma \right) }{\left\vert M\right\vert }%
\right) ^{2}\left\Vert \mathsf{Q}_{\mathcal{C}_{F}^{\mathcal{G},\limfunc{%
shift}};M}^{\omega ,\mathbf{b}^{\ast }}x\right\Vert _{L^{2}\left( \omega
\right) }^{\spadesuit 2}\leq \mathfrak{F}_{\alpha }\lVert h\rVert
_{L^{2}\left( \sigma \right) }\,.
\end{equation*}

\begin{remark}
\label{explaining funct ener}If in (\ref{e.funcEnergy n}), we take $h=%
\mathbf{1}_{I}$ and $\mathcal{F}$ to be the trivial Carleson collection $%
\left\{ I_{r}\right\} _{r=1}^{\infty }$ where the intervals $I_{r}$ are
pairwise disjoint in $I$, then we essentially obtain the Whitney energy
condition in Definition \ref{energy condition}, but with $\mathsf{Q}_{%
\mathcal{C}_{F}^{\mathcal{G},\limfunc{shift}};M}^{\omega ,\mathbf{b}^{\ast
}} $ in place of $\mathsf{Q}_{M}^{\limfunc{weak}\limfunc{good},\omega }$.
However, the pseudoprojection $\mathsf{Q}_{M}^{\limfunc{weak}\limfunc{good}%
,\omega }$ is larger than $\mathsf{Q}_{\mathcal{C}_{F}^{\mathcal{G},\limfunc{%
shift}};J}^{\omega ,\mathbf{b}^{\ast }}$, and so we just miss obtaining the
Whitney energy condition as a consequence of the functional energy
condition. Nevertheless, this near miss with $h=\mathbf{1}_{I}$ explains the
terminology `functional' energy.
\end{remark}

We will need an `indicator' version of the estimate proved above for the
disjoint form%
\begin{equation*}
\Theta _{1}\left( f,g\right) =\sum_{I\in \mathcal{D}}\sum_{\substack{ J\in 
\mathcal{G}:\ \ell \left( J\right) \leq \ell \left( I\right)  \\ d\left(
J,I\right) >2\ell \left( J\right) ^{\varepsilon }\ell \left( I\right)
^{1-\varepsilon }}}\int \left( T_{\sigma }\square _{I}^{\sigma ,\mathbf{b}%
}f\right) \square _{J}^{\omega ,\mathbf{b}^{\ast }}gd\omega .
\end{equation*}

\begin{lemma}
\label{standard indicator}Fix dyadic grids $\mathcal{D}$ and $\mathcal{G}$.
Suppose $T^{\alpha }$ is a standard fractional singular integral with $0\leq
\alpha <1$, that $f\in L^{2}\left( \sigma \right) $ and $g\in L^{2}\left(
\omega \right) $, that $\mathcal{F}\subset \mathcal{D}$ and $\mathcal{H}%
\subset \mathcal{G}$ are $\sigma $-Carleson and $\omega $-Carleson
collections respectively, i.e.,%
\begin{equation*}
\sum_{F^{\prime }\in \mathcal{F}:\ F^{\prime }\subset F}\left\vert F^{\prime
}\right\vert _{\sigma }\lesssim \left\vert F\right\vert _{\sigma },\ \ \ \ \
F\in \mathcal{F},\text{ and }\sum_{H^{\prime }\in \mathcal{H}:\ H^{\prime
}\subset H}\left\vert H^{\prime }\right\vert _{\omega }\lesssim \left\vert
H\right\vert _{\omega },\ \ \ \ \ H\in \mathcal{H},
\end{equation*}%
and that there are numerical sequences $\left\{ \alpha _{\mathcal{F}}\left(
F\right) \right\} _{F\in \mathcal{F}}$ and $\left\{ \beta _{\mathcal{H}%
}\left( H\right) \right\} _{H\in \mathcal{H}}$ such that%
\begin{equation}
\sum_{F\in \mathcal{F}}\alpha _{\mathcal{F}}\left( F\right) ^{2}\left\vert
F\right\vert _{\sigma }\leq \left\Vert f\right\Vert _{L^{2}\left( \sigma
\right) }^{2}\text{ and }\sum_{H\in \mathcal{H}}\beta _{\mathcal{H}}\left(
H\right) ^{2}\left\vert H\right\vert _{\sigma }\leq \left\Vert g\right\Vert
_{L^{2}\left( \sigma \right) }^{2}\ .  \label{qo}
\end{equation}%
Then%
\begin{eqnarray}
&&\sum_{F\in \mathcal{F}}\sum_{\substack{ J\in \mathcal{G}:\ \ell \left(
J\right) \leq \ell \left( F\right)  \\ d\left( J,F\right) >2\ell \left(
J\right) ^{\varepsilon }\ell \left( F\right) ^{1-\varepsilon }}}\left\vert
\int \left( T_{\sigma }^{\alpha }\mathbf{1}_{F}\alpha _{\mathcal{F}}\left(
F\right) \right) \square _{J}^{\omega ,\mathbf{b}^{\ast }}gd\omega
\right\vert  \label{indicator far} \\
&&+\sum_{G\in \mathcal{G}}\sum_{\substack{ I\in \mathcal{D}:\ \ell \left(
I\right) \leq \ell \left( G\right)  \\ d\left( I,G\right) >2\ell \left(
I\right) ^{\varepsilon }\ell \left( G\right) ^{1-\varepsilon }}}\left\vert
\int \left( T_{\sigma }^{\alpha }\square _{I}^{\sigma ,\mathbf{b}}f\right) 
\mathbf{1}_{G}\beta _{\mathcal{G}}\left( G\right) d\omega \right\vert  \notag
\\
&\lesssim &\sqrt{A_{2}^{\alpha }}\left\Vert f\right\Vert _{L^{2}\left(
\sigma \right) }\left\Vert g\right\Vert _{L^{2}\left( \omega \right) }. 
\notag
\end{eqnarray}
\end{lemma}

The proof of this lemma is similar to those of Lemmas \ref{delta long} and %
\ref{delta short} in Section \ref{Sec disj form}\ above, using the square
function inequalities for $\square _{I}^{\sigma ,\mathbf{b}}$, $\nabla _{I,%
\mathcal{F}}^{\sigma }$ and $\square _{J}^{\omega ,\mathbf{b}^{\ast }}$, $%
\nabla _{J,\mathcal{G}}^{\omega }$ in Appendix A, as well as the
quasiorthogonal inequalities assumed in (\ref{qo}), which substitute for the
square function inequalities when dealing with indicators $\mathbf{1}%
_{F}\alpha _{\mathcal{F}}\left( F\right) $ instead instead of dual
martingale differences $\square _{I}^{\sigma ,\mathbf{b}}f$. We note that
there is no explicit restriction of the type $\ell \left( J\right) \leq 2^{-%
\mathbf{\rho }}\ell \left( I\right) $ in any of Lemmas \ref{delta long}, \ref%
{delta short}, or \ref{standard indicator}.

There is one more simple lemma that we will use in the proof of the
Intertwining Proposition, namely that for small $\varepsilon >0$, an
interval $J$ is $\varepsilon $-$\limfunc{good}$ inside an interval $I$ only
if $J$ is many scales smaller in size than $I$. Recall from Definition \ref%
{good arb} that if an interval $J$ is $\varepsilon $-$\limfunc{good}$ inside
an interval $I$, then 
\begin{equation*}
J\subset I\text{ and }d\left( J,\limfunc{skel}I\right) \geq d\left( J,%
\limfunc{body}I\right) >2\ell \left( J\right) ^{\varepsilon }\ell \left(
I\right) ^{1-\varepsilon }.
\end{equation*}

\begin{lemma}
\label{good scale}If $J$ is $\varepsilon $-$\limfunc{good}$ inside $I$, then 
$\ell \left( J\right) <2^{-\frac{3}{\varepsilon }}\ell \left( I\right) $.
\end{lemma}

\begin{proof}
We have%
\begin{equation*}
\frac{1}{4}\ell \left( I\right) \geq d\left( J,\limfunc{skel}I\right) >2\ell
\left( J\right) ^{\varepsilon }\ell \left( I\right) ^{1-\varepsilon
}=2\left( \frac{\ell \left( J\right) }{\ell \left( I\right) }\right)
^{\varepsilon }\ell \left( I\right) ,
\end{equation*}%
which gives $\frac{1}{8}>\left( \frac{\ell \left( J\right) }{\ell \left(
I\right) }\right) ^{\varepsilon }$, i.e. $\frac{\ell \left( J\right) }{\ell
\left( I\right) }<\left( \frac{1}{8}\right) ^{\frac{1}{\varepsilon }}=2^{-%
\frac{3}{\varepsilon }}$.
\end{proof}

\begin{proposition}[The Intertwining Proposition]
\label{strongly adapted}Let $\mathcal{D}$ and $\mathcal{G}$ be grids, and
suppose that $\mathbf{b}$ and $\mathbf{b}^{\ast }$ are $\infty $-strongly $%
\sigma $-accretive families of intervals in $\mathcal{D}$ and $\mathcal{G}$
respectively. Suppose that $\mathcal{F}\subset \mathcal{D}$ is $\sigma $%
-Carleson and that the $\mathcal{F}$-coronas 
\begin{equation*}
\mathcal{C}_{F}\equiv \left\{ I\in \mathcal{D}:I\subset F\text{ but }%
I\not\subset F^{\prime }\text{ for }F^{\prime }\in \mathfrak{C}_{\mathcal{F}%
}\left( F\right) \right\}
\end{equation*}%
satisfy%
\begin{equation*}
E_{I}^{\sigma }\left\vert f\right\vert \lesssim E_{F}^{\sigma }\left\vert
f\right\vert \text{ and }b_{I}=\mathbf{1}_{I}b_{F},\ \ \ \ \ \text{for all }%
I\in \mathcal{C}_{F}\mathfrak{\ ,\ }F\in \mathcal{F}.
\end{equation*}%
Then with the shifted corona in Definition \ref{def shift}, i.e. $\mathcal{C}%
_{F}^{\mathcal{G},\limfunc{shift}}=\left\{ J\in \mathcal{G}:J^{\maltese }\in 
\mathcal{C}_{F}\right\} $ with $J^{\maltese }$ as in Definition \ref{def
sharp cross} that depends on $\varepsilon >0$, we have%
\begin{equation*}
\left\vert \sum_{F\in \mathcal{F}}\ \sum_{I:\ I\supsetneqq F}\ \left\langle
T_{\sigma }^{\alpha }\square _{I}^{\sigma ,\mathbf{b}}f,\mathsf{P}_{\mathcal{%
C}_{F}^{\mathcal{G},\limfunc{shift}}}^{\omega ,\mathbf{b}^{\ast
}}g\right\rangle _{\omega }\right\vert \lesssim \left( \mathfrak{F}_{\alpha
}^{\mathbf{b}^{\ast }}\left( \mathcal{D},\mathcal{G}\right) +\mathcal{NTV}%
_{\alpha }\right) \ \left\Vert f\right\Vert _{L^{2}\left( \sigma \right)
}\left\Vert g\right\Vert _{L^{2}\left( \omega \right) },
\end{equation*}%
where the implied constant depends on the $\sigma $-Carleson norm $C_{%
\mathcal{F}}$ of the family $\mathcal{F}$.
\end{proposition}

\begin{proof}[Proof of Proposition \protect\ref{strongly adapted}]
We write the sum on the left hand side of the display above as%
\begin{eqnarray*}
&&\sum_{F\in \mathcal{F}}\ \sum_{I:\ I_{\infty }\supset I\supsetneqq F}\
\left\langle T_{\sigma }^{\alpha }\square _{I}^{\sigma ,\mathbf{b}}f,\mathsf{%
P}_{\mathcal{C}_{F}^{\mathcal{G},\limfunc{shift}}}^{\omega }g\right\rangle
_{\omega }=\sum_{F\in \mathcal{F}}\ \left\langle T_{\sigma }^{\alpha }\left(
\sum_{I:\ I_{\infty }\supset I\supsetneqq F}\square _{I}^{\sigma ,\mathbf{b}%
}f\right) ,\mathsf{P}_{\mathcal{C}_{F}^{\mathcal{G},\limfunc{shift}%
}}^{\omega }g\right\rangle _{\omega }=\sum_{F\in \mathcal{F}}\ \left\langle
T_{\sigma }^{\alpha }\left( f_{F}^{\ast }\right) ,g_{F}\right\rangle
_{\omega }; \\
&&\text{where }f_{F}^{\ast }\equiv \sum_{I:\ I_{\infty }\supset I\supsetneqq
F}\square _{I}^{\sigma ,\mathbf{b}}f\text{ and }g_{F}\equiv \mathsf{P}_{%
\mathcal{C}_{F}^{\mathcal{G},\limfunc{shift}}}^{\omega ,\mathbf{b}^{\ast }}g,
\end{eqnarray*}%
where $I_{\infty }$ is the starting interval for corona constructions in $%
\mathcal{D}$ as in (\ref{top control}) above. Note that $g_{F}$ is supported
in $F$. By the telescoping identity for $\square _{I}^{\sigma ,\mathbf{b}}$,
the function $f_{F}^{\ast }$ satisfies%
\begin{equation*}
\mathbf{1}_{F}f_{F}^{\ast }=\sum_{I:\ I_{\infty }\supset I\supsetneqq
F}\square _{I}^{\sigma ,\mathbf{b}}f=\mathbb{F}_{F}^{\sigma ,\mathbf{b}}f-%
\mathbf{1}_{F}\mathbb{F}_{I_{\infty }}^{\sigma ,\mathbf{b}}f=b_{F}\frac{%
E_{F}^{\sigma }f}{E_{F}^{\sigma }b_{F}}-\mathbf{1}_{F}b_{I_{\infty }}\frac{%
E_{I_{\infty }}^{\sigma }f}{E_{I_{\infty }}^{\sigma }b_{I_{\infty }}}\ .
\end{equation*}%
However, we cannot apply the testing condition to the function $\mathbf{1}%
_{F}b_{I_{\infty }}$, and since $E_{I_{\infty }}^{\sigma }f$ does not vanish
in general, we will instead add and subtract the term $\mathbb{F}_{I_{\infty
}}^{\sigma ,\mathbf{b}}f$ to get 
\begin{eqnarray*}
\sum_{F\in \mathcal{F}}\ \left\langle T_{\sigma }^{\alpha }\left(
f_{F}^{\ast }\right) ,g_{F}\right\rangle _{\omega } &=&\sum_{F\in \mathcal{F}%
}\ \left\langle T_{\sigma }^{\alpha }\left( \sum_{I:\ I_{\infty }\supset
I\supsetneqq F}\square _{I}^{\sigma ,\mathbf{b}}f\right) ,\mathsf{P}_{%
\mathcal{C}_{F}^{\mathcal{G},\limfunc{shift}}}^{\omega }g\right\rangle
_{\omega } \\
&=&\sum_{F\in \mathcal{F}}\ \left\langle T_{\sigma }^{\alpha }\left( \mathbb{%
F}_{I_{\infty }}^{\sigma ,\mathbf{b}}f+\sum_{I:\ I_{\infty }\supset
I\supsetneqq F}\square _{I}^{\sigma ,\mathbf{b}}f\right) ,\mathsf{P}_{%
\mathcal{C}_{F}^{\mathcal{G},\limfunc{shift}}}^{\omega }g\right\rangle
_{\omega }-\sum_{F\in \mathcal{F}}\ \left\langle T_{\sigma }^{\alpha }\left( 
\mathbb{F}_{I_{\infty }}^{\sigma ,\mathbf{b}}f\right) ,\mathsf{P}_{\mathcal{C%
}_{F}^{\mathcal{G},\limfunc{shift}}}^{\omega }g\right\rangle _{\omega }\ ,
\end{eqnarray*}%
where the second sum on the right hand side of the identity satisfies%
\begin{eqnarray*}
\left\vert \sum_{F\in \mathcal{F}}\ \left\langle T_{\sigma }^{\alpha }\left( 
\mathbb{F}_{I_{\infty }}^{\sigma ,\mathbf{b}}f\right) ,\mathsf{P}_{\mathcal{C%
}_{F}^{\mathcal{G},\limfunc{shift}}}^{\omega }g\right\rangle _{\omega
}\right\vert &=&\left\vert \left\langle T_{\sigma }^{\alpha }\left( \mathbb{F%
}_{I_{\infty }}^{\sigma ,\mathbf{b}}f\right) ,\sum_{F\in \mathcal{F}}\mathsf{%
P}_{\mathcal{C}_{F}^{\mathcal{G},\limfunc{shift}}}^{\omega }g\right\rangle
_{\omega }\right\vert \\
&\leq &\left\Vert T_{\sigma }^{\alpha }\left( \mathbb{F}_{I_{\infty
}}^{\sigma ,\mathbf{b}}f\right) \right\Vert _{L^{2}\left( \omega \right)
}\left\Vert \sum_{F\in \mathcal{F}}\mathsf{P}_{\mathcal{C}_{F;\mathbf{r}}^{%
\mathcal{G},\limfunc{shift}}}^{\omega }g\right\Vert _{L^{2}\left( \omega
\right) } \\
&\lesssim &\mathfrak{FT}_{T_{\alpha }}^{\mathbf{b}}\left\vert E_{I_{\infty
}}^{\sigma }f\right\vert \left\Vert g\right\Vert _{L^{2}\left( \omega
\right) }\lesssim \left( \mathfrak{T}_{T_{\alpha }}^{\mathbf{b}}+\mathfrak{A}%
_{2}^{\alpha }\right) \left\Vert f\right\Vert _{L^{2}\left( \sigma \right)
}\left\Vert g\right\Vert _{L^{2}\left( \omega \right) }
\end{eqnarray*}%
by (\ref{full proved}) above, and the Riesz inequalities in Appendix A. The
advantage now is that with%
\begin{equation*}
f_{F}\equiv \mathbb{F}_{I_{\infty }}^{\sigma ,\mathbf{b}}f+f_{F}^{\ast }=%
\mathbb{F}_{I_{\infty }}^{\sigma ,\mathbf{b}}f+\sum_{I:\ I_{\infty }\supset
I\supsetneqq F}\square _{I}^{\sigma ,\mathbf{b}}f
\end{equation*}%
then in the first term on the right hand side of the identity, the
telescoping identity gives%
\begin{equation*}
\mathbf{1}_{F}f_{F}=\mathbf{1}_{F}\left( \mathbb{F}_{I_{\infty }}^{\sigma ,%
\mathbf{b}}f+\sum_{I:\ I_{\infty }\supset I\supsetneqq F}\square
_{I}^{\sigma ,\mathbf{b}}f\right) =\mathbb{F}_{F}^{\sigma ,\mathbf{b}}f=b_{F}%
\frac{E_{F}^{\sigma }f}{E_{F}^{\sigma }b_{F}},
\end{equation*}%
which shows that $f_{F}$ is a controlled constant times $b_{F}$ on $F$.

The intervals $I$ occurring in this sum are linearly and consecutively
ordered by inclusion, along with the intervals $F^{\prime }\in \mathcal{F}$
that contain $F$. More precisely we can write%
\begin{equation*}
F\equiv F_{0}\subsetneqq F_{1}\subsetneqq F_{2}\subsetneqq ...\subsetneqq
F_{n}\subsetneqq F_{n+1}\subsetneqq ...F_{N}=I_{\infty }
\end{equation*}%
where $F_{m}=\pi _{\mathcal{F}}^{m}F$ for all $m\geq 1$. We can also write%
\begin{equation*}
F=F_{0}\equiv I_{0}\subsetneqq I_{1}\subsetneqq I_{2}\subsetneqq
...\subsetneqq I_{k}\subsetneqq I_{k+1}\subsetneqq ...\subsetneqq
I_{K}=F_{N}=I_{\infty }
\end{equation*}%
where $I_{k}=\pi _{\mathcal{D}}^{k}F$ for all $k\geq 1$. There is a (unique)
subsequence $\left\{ k_{m}\right\} _{m=1}^{N}$ such that%
\begin{equation*}
F_{m}=I_{k_{m}},\ \ \ \ \ 1\leq m\leq N.
\end{equation*}

Then we have%
\begin{eqnarray*}
f_{F}\left( x\right) &\equiv &\mathbb{F}_{I_{\infty }}^{\sigma ,\mathbf{b}%
}f\left( x\right) +\sum_{\ell =1}^{K}\square _{I_{\ell }}^{\sigma ,\mathbf{b}%
}f\left( x\right) , \\
g_{F} &\equiv &\sum_{J\in \mathcal{C}_{F}^{\mathcal{G},\limfunc{shift}%
}}\square _{J}^{\omega ,\mathbf{b}^{\ast }}g.
\end{eqnarray*}%
Assume now that $k_{m}\leq k<k_{m+1}$. We denote the sibling of $I$ by $%
\theta \left( I\right) $, i.e. $\left\{ \theta \left( I\right) \right\} =%
\mathfrak{C}_{\mathcal{D}}\left( \pi _{\mathcal{D}}I\right) \setminus
\left\{ I\right\} $. There are two cases to consider here:%
\begin{equation*}
\theta \left( I_{k}\right) \notin \mathcal{F}\text{ and }\theta \left(
I_{k}\right) \in \mathcal{F}.
\end{equation*}%
We first note that in either case, using a telescoping sum, we compute that
for 
\begin{equation*}
x\in \theta \left( I_{k}\right) =I_{k+1}\setminus I_{k}\subset
F_{m+1}\setminus F_{m},
\end{equation*}%
we have the formula 
\begin{eqnarray*}
f_{F}\left( x\right) &=&\mathbb{F}_{I_{\infty }}^{\sigma ,\mathbf{b}}f\left(
x\right) +\sum_{\ell \geq k+1}\square _{I_{\ell }}^{\sigma ,\mathbf{b}%
}f\left( x\right) \\
&=&\mathbb{F}_{\theta \left( I_{k}\right) }^{\sigma ,\mathbf{b}}f\left(
x\right) -\mathbb{F}_{I_{k+1}}^{\sigma ,\mathbf{b}}f\left( x\right)
+\sum_{\ell =k+1}^{K-1}\left( \mathbb{F}_{I_{\ell }}^{\sigma ,\mathbf{b}%
}f\left( x\right) -\mathbb{F}_{I_{\ell +1}}^{\sigma ,\mathbf{b}}f\left(
x\right) \right) +\mathbb{F}_{I_{\infty }}^{\sigma ,\mathbf{b}}f\left(
x\right) =\mathbb{F}_{\theta \left( I_{k}\right) }^{\sigma ,\mathbf{b}%
}f\left( x\right) \ .
\end{eqnarray*}%
Now fix $x\in \theta \left( I_{k}\right) $. If $\theta \left( I_{k}\right)
\notin \mathcal{F}$, then $\theta \left( I_{k}\right) \in \mathcal{C}%
_{F_{m+1}}^{\sigma }$, and we have 
\begin{equation}
\left\vert f_{F}\left( x\right) \right\vert =\left\vert \mathbb{F}_{\theta
\left( I_{k}\right) }^{\sigma ,\mathbf{b}}f\left( x\right) \right\vert
\lesssim \left\vert b_{\theta \left( I_{k}\right) }\left( x\right)
\right\vert \ \frac{E_{\theta \left( I_{k}\right) }^{\sigma }\left\vert
f\right\vert }{\left\vert E_{\theta \left( I_{k}\right) }^{\sigma }b_{\theta
\left( I_{k}\right) }\right\vert }\lesssim E_{F_{m+1}}^{\sigma }\left\vert
f\right\vert \ ,  \label{bound for f_F}
\end{equation}%
since the testing functions $b_{\theta \left( I_{k}\right) }$ are bounded
and accretive, and $E_{\theta \left( I_{k}\right) }^{\sigma }\left\vert
f\right\vert \lesssim E_{F_{m+1}}^{\sigma }\left\vert f\right\vert $ by
hypothesis. On the other hand, if $\theta \left( I_{k}\right) \in \mathcal{F}
$, then $I_{k+1}\in \mathcal{C}_{F_{m+1}}^{\sigma }$ and we have%
\begin{equation*}
\left\vert f_{F}\left( x\right) \right\vert =\left\vert \mathbb{F}_{\theta
\left( I_{k}\right) }^{\sigma ,\mathbf{b}}f\left( x\right) \right\vert
\lesssim E_{\theta \left( I_{k}\right) }^{\sigma }\left\vert f\right\vert \ .
\end{equation*}%
Note that $F^{c}=\overset{\cdot }{\dbigcup }_{k\geq 0}\theta \left(
I_{k}\right) $. Now we write%
\begin{eqnarray*}
f_{F} &=&\varphi _{F}+\psi _{F}, \\
\varphi _{F} &\equiv &\sum_{k:\ \theta \left( I_{k}\right) \in \mathcal{F}}%
\mathbb{F}_{\theta \left( I_{k}\right) }^{\sigma ,\mathbf{b}}f\text{ and }%
\psi _{F}=f_{F}-\varphi _{F}\ ; \\
\sum_{F\in \mathcal{F}}\ \left\langle T_{\sigma }^{\alpha
}f_{F},g_{F}\right\rangle _{\omega } &=&\sum_{F\in \mathcal{F}}\
\left\langle T_{\sigma }^{\alpha }\varphi _{F},g_{F}\right\rangle _{\omega
}+\sum_{F\in \mathcal{F}}\ \left\langle T_{\sigma }^{\alpha }\psi
_{F},g_{F}\right\rangle _{\omega }\ ,
\end{eqnarray*}%
and note that $\varphi _{F}=0$ on $F$, and $\psi _{F}=b_{F}\frac{%
E_{F}^{\sigma }f}{E_{F}^{\sigma }b_{F}}$ on $F$. We can apply the first line
in (\ref{indicator far}) using $\theta \left( I_{k}\right) \in \mathcal{F}$
to the first sum above since $J\in \mathcal{C}_{F}^{\mathcal{G},\limfunc{%
shift}}$ implies $J\subset J^{\maltese }\subset F\subset \theta \left(
I_{k}\right) ^{c}$, which implies that $d\left( J,\theta \left( I_{k}\right)
\right) >2\ell \left( J\right) ^{\varepsilon }\ell \left( \theta \left(
I_{k}\right) \right) ^{1-\varepsilon }$. Thus we obtain after substituting $%
F^{\prime }$ for $\theta \left( I_{k}\right) $ below, 
\begin{eqnarray*}
\left\vert \sum_{F\in \mathcal{F}}\ \left\langle T_{\sigma }^{\alpha
}\varphi _{F},g_{F}\right\rangle _{\omega }\right\vert &=&\left\vert
\sum_{F\in \mathcal{F}}\sum_{J\in \mathcal{C}_{F}^{\mathcal{G},\limfunc{shift%
}}}\left\langle T_{\sigma }^{\alpha }\left( \sum_{k:\ \theta \left(
I_{k}\right) \in \mathcal{F}}\mathbb{F}_{\theta \left( I_{k}\right)
}^{\sigma ,\mathbf{b}}f\right) ,\square _{J}^{\omega ,\mathbf{b}^{\ast
}}g\right\rangle _{\omega }\right\vert \\
&\leq &\sum_{F\in \mathcal{F}}\sum_{J\in \mathcal{C}_{F}^{\mathcal{G},%
\limfunc{shift}}}\sum_{k:\ \theta \left( I_{k}\right) \in \mathcal{F}%
}\left\vert \left\langle T_{\sigma }^{\alpha }\left( \mathbb{F}_{\theta
\left( I_{k}\right) }^{\sigma ,\mathbf{b}}f\right) ,\square _{J}^{\omega ,%
\mathbf{b}^{\ast }}g\right\rangle _{\omega }\right\vert \\
&\leq &\sum_{F^{\prime }\in \mathcal{F}}\sum_{\substack{ J\in \mathcal{G}:\
\ell \left( J\right) \leq \ell \left( F^{\prime }\right)  \\ d\left(
J,F^{\prime }\right) >2\ell \left( J\right) ^{\varepsilon }\ell \left(
F^{\prime }\right) ^{1-\varepsilon }}}\left\vert \left\langle T_{\sigma
}^{\alpha }\left( \mathbb{F}_{F^{\prime }}^{\sigma ,\mathbf{b}}f\right)
,\square _{J}^{\omega ,\mathbf{b}^{\ast }}g\right\rangle _{\omega
}\right\vert \\
&\lesssim &\sqrt{A_{2}^{\alpha }}\left\Vert f\right\Vert _{L^{2}\left(
\sigma \right) }\left\Vert g\right\Vert _{L^{2}\left( \omega \right) }\ .
\end{eqnarray*}

Turning to the second sum, we note that for $k_{m}\leq k<k_{m+1}$ and $x\in
\theta \left( I_{k}\right) $ with $\theta \left( I_{k}\right) \notin 
\mathcal{F}$, we have%
\begin{equation*}
\left\vert \psi _{F}\left( x\right) \right\vert \lesssim \left\vert
b_{\theta \left( I_{k}\right) }\right\vert \ E_{\theta \left( I_{k}\right)
}^{\sigma }\left\vert f\right\vert \ \mathbf{1}_{F_{m+1}\setminus
F_{m}}\left( x\right) \lesssim \alpha _{\mathcal{F}}\left( F_{m+1}\right) \ 
\mathbf{1}_{F_{m+1}\setminus F_{m}}\left( x\right) \ ,
\end{equation*}%
and hence the following inequality for $x\notin F$, 
\begin{equation}
\left\vert \psi _{F}\left( x\right) \right\vert \lesssim \sum_{F^{\prime
}\in \mathcal{F}:\ F\subset F^{\prime }}\alpha _{\mathcal{F}}\left( \pi _{%
\mathcal{F}}F^{\prime }\right) \ \mathbf{1}_{\pi _{\mathcal{F}}F^{\prime
}\setminus F^{\prime }}\left( x\right) =\Phi \left( x\right) \ \mathbf{1}%
_{F^{c}}\left( x\right) \ ,  \label{Psi_F bound}
\end{equation}%
where%
\begin{equation*}
\Phi \equiv \sum_{F^{\prime \prime }\in \mathcal{F}}\alpha _{\mathcal{F}%
}\left( F^{\prime \prime }\right) \ \mathbf{1}_{F^{\prime \prime }\setminus
\cup \mathfrak{C}_{\mathcal{F}}\left( F^{\prime \prime }\right) }=\sum_{F\in 
\mathcal{F}}\alpha _{\mathcal{F}}\left( F\right) \mathbf{1}_{F\setminus \cup 
\mathfrak{C}_{\mathcal{F}}\left( F\right) }\ .
\end{equation*}

Now we write%
\begin{equation*}
\sum_{F\in \mathcal{F}}\ \left\langle T_{\sigma }^{\alpha }\psi
_{F},g_{F}\right\rangle _{\omega }=\sum_{F\in \mathcal{F}}\ \left\langle
T_{\sigma }^{\alpha }\left( \mathbf{1}_{F}\psi _{F}\right)
,g_{F}\right\rangle _{\omega }+\sum_{F\in \mathcal{F}}\ \left\langle
T_{\sigma }^{\alpha }\left( \mathbf{1}_{F^{c}}\psi _{F}\right)
,g_{F}\right\rangle _{\omega }\equiv I+II,
\end{equation*}%
where $I$ and $II$ are defined at the end of the display. Then by interval
testing, 
\begin{equation*}
\left\vert \left\langle T_{\sigma }^{\alpha }\left( b_{F}\mathbf{1}%
_{F}\right) ,g_{F}\right\rangle _{\omega }\right\vert =\left\vert
\left\langle \mathbf{1}_{F}T_{\sigma }^{\alpha }\left( b_{F}\mathbf{1}%
_{F}\right) ,g_{F}\right\rangle _{\omega }\right\vert \lesssim \mathfrak{T}%
_{T^{\alpha }}\sqrt{\left\vert F\right\vert _{\sigma }}\left\Vert
g_{F}\right\Vert _{L^{2}\left( \omega \right) }^{\bigstar }\ ,
\end{equation*}%
and so quasi-orthogonality, together with the fact that on $F$, $\psi
_{F}=b_{F}\frac{E_{F}^{\sigma }f}{E_{F}^{\sigma }b_{F}}$ is a constant $c=%
\frac{E_{F}^{\sigma }f}{E_{F}^{\sigma }b_{F}}$ times $b_{F}$, where $%
\left\vert c\right\vert $ is bounded by $\alpha _{\mathcal{F}}\left(
F\right) $, give

\begin{eqnarray*}
&&\left\vert I\right\vert =\left\vert \sum_{F\in \mathcal{F}}\ \left\langle
T_{\sigma }^{\alpha }\left( \mathbf{1}_{F}cb_{F}\right) ,g_{F}\right\rangle
_{\omega }\right\vert \lesssim \sum_{F\in \mathcal{F}}\ \alpha _{\mathcal{F}%
}\left( F\right) \ \left\vert \left\langle T_{\sigma }^{\alpha
}b_{F},g_{F}\right\rangle _{\omega }\right\vert \\
&\lesssim &\sum_{F\in \mathcal{F}}\ \alpha _{\mathcal{F}}\left( F\right) 
\mathcal{NTV}_{\alpha }\sqrt{\left\vert F\right\vert _{\sigma }}\left\Vert
g_{F}\right\Vert _{L^{2}\left( \omega \right) }^{\bigstar }\lesssim \mathcal{%
NTV}_{\alpha }\left\Vert f\right\Vert _{L^{2}\left( \sigma \right) }\left[
\sum_{F\in \mathcal{F}}\left\Vert g_{F}\right\Vert _{L^{2}\left( \omega
\right) }^{\bigstar 2}\right] ^{\frac{1}{2}}.
\end{eqnarray*}

Now $\mathbf{1}_{F^{c}}\psi _{F}$ is supported outside $F$, and each $J$ in
the dual martingale support $\mathcal{C}_{F}^{\mathcal{G},\limfunc{shift}}$
of $g_{F}=\mathsf{P}_{\mathcal{C}_{F}^{\mathcal{G},\limfunc{shift}}}^{\omega
}g$ is in particular $\limfunc{good}$ in the interval $F$, and as a
consequence, each such interval $J$ as above is contained in some interval $%
M $ for $M\in \mathcal{W}\left( F\right) $. This containment will be used in
the analysis of the term $II_{G}$ below.

In addition, each $J$ in the dual martingale support $\mathcal{C}_{F}^{%
\mathcal{G},\limfunc{shift}}$ of $g_{F}=\mathsf{P}_{\mathcal{C}_{F}^{%
\mathcal{G},\limfunc{shift}}}^{\omega }g$ is $\left( \left[ \frac{3}{%
\varepsilon }\right] ,\varepsilon \right) $-deeply embedded in $F$, i.e. $%
J\Subset _{\left[ \frac{3}{\varepsilon }\right] ,\varepsilon }F$, by Lemma %
\ref{good scale} and the definition of $\mathcal{C}_{F}^{\mathcal{G},%
\limfunc{shift}}$ in Definition \ref{def shift}. As a consequence, each such
interval $J$ as above is contained in some interval $M$ for $M\in \mathcal{M}%
_{\left( \left[ \frac{3}{\varepsilon }\right] ,\varepsilon \right) -\limfunc{%
deep},\mathcal{D}}\left( F\right) $. This containment will be used in the
analysis of the term $II_{B}$ below.

\begin{notation}
Define $\mathbf{\rho }\equiv \left[ \frac{3}{\varepsilon }\right] $, so that
for every $J\in \mathcal{C}_{F}^{\mathcal{G},\limfunc{shift}}$, there is $%
M\in \mathcal{M}_{\left( \mathbf{\rho },\varepsilon \right) -\limfunc{deep},%
\mathcal{G}}\left( F\right) $ such that $J\subset M$.
\end{notation}

The collections $\mathcal{W}\left( F\right) $ and $\mathcal{M}_{\left( 
\mathbf{\rho },\varepsilon \right) -\limfunc{deep},\mathcal{G}}\left(
F\right) $ used here, and in the display below, are defined in (\ref{def
M_r-deep}) in Appendix B. Finally, since the intervals $M\in \mathcal{W}%
\left( F\right) $, as well as the intervals $M\in \mathcal{M}_{\left( \left[ 
\frac{3}{\varepsilon }\right] ,\varepsilon \right) -\limfunc{deep},\mathcal{G%
}}\left( F\right) $, satisfy $3M\subset F$, we can apply (\ref{estimate}) in
the Monotonicity Lemma \ref{mono} using (\ref{Psi_F bound}) with $\mu =%
\mathbf{1}_{F^{c}}\psi _{F}$ and $J^{\prime }$ in place of $J$ there, to
obtain%
\begin{eqnarray*}
\left\vert II\right\vert &=&\left\vert \sum_{F\in \mathcal{F}}\left\langle
T_{\sigma }^{\alpha }\left( \mathbf{1}_{F^{c}}\psi _{F}\right)
,g_{F}\right\rangle _{\omega }\right\vert =\left\vert \sum_{F\in \mathcal{F}%
}\sum_{J^{\prime }\in \mathcal{C}_{F}^{\mathcal{G},\limfunc{shift}%
}}\left\langle T_{\sigma }^{\alpha }\left( \mathbf{1}_{F^{c}}\psi
_{F}\right) ,\square _{J^{\prime }}^{\omega ,\mathbf{b}^{\ast
}}g\right\rangle _{\omega }\right\vert \\
&\lesssim &\sum_{F\in \mathcal{F}}\sum_{J^{\prime }\in \mathcal{C}_{F}^{%
\mathcal{G},\limfunc{shift}}}\frac{\mathrm{P}^{\alpha }\left( J^{\prime },%
\mathbf{1}_{F^{c}}\Phi \sigma \right) }{\left\vert J^{\prime }\right\vert }%
\left\Vert \bigtriangleup _{J^{\prime }}^{\omega ,\mathbf{b}^{\ast
}}x\right\Vert _{L^{2}\left( \omega \right) }^{\spadesuit }\left\Vert
\square _{J^{\prime }}^{\omega ,\mathbf{b}^{\ast }}g\right\Vert
_{L^{2}\left( \omega \right) }^{\bigstar } \\
&&+\sum_{F\in \mathcal{F}}\sum_{J^{\prime }\in \mathcal{C}_{F}^{\mathcal{G},%
\limfunc{shift}}}\frac{\mathrm{P}_{1+\delta }^{\alpha }\left( J^{\prime },%
\mathbf{1}_{F^{c}}\Phi \sigma \right) }{\left\vert J^{\prime }\right\vert }%
\left\Vert x-m_{J^{\prime }}\right\Vert _{L^{2}\left( \omega \right)
}\left\Vert \square _{J^{\prime }}^{\omega ,\mathbf{b}^{\ast }}g\right\Vert
_{L^{2}\left( \omega \right) }^{\bigstar } \\
&\lesssim &\sum_{F\in \mathcal{F}}\sum_{M\in \mathcal{W}\left( F\right) }%
\frac{\mathrm{P}^{\alpha }\left( M,\mathbf{1}_{F^{c}}\Phi \sigma \right) }{%
\left\vert M\right\vert }\left\Vert \mathsf{Q}_{\mathcal{C}_{F;M}^{\mathcal{G%
},\limfunc{shift}}}^{\omega ,\mathbf{b}^{\ast }}x\right\Vert _{L^{2}\left(
\omega \right) }^{\spadesuit }\left\Vert g_{F;M}\right\Vert _{L^{2}\left(
\omega \right) }^{\bigstar } \\
&&+\sum_{F\in \mathcal{F}}\sum_{J\in \mathcal{M}_{\left( \mathbf{\rho }%
,\varepsilon \right) -\limfunc{deep},\mathcal{G}}\left( F\right)
}\sum_{J^{\prime }\in \mathcal{C}_{F;J}^{\mathcal{G},\limfunc{shift}}}\frac{%
\mathrm{P}_{1+\delta }^{\alpha }\left( J^{\prime },\mathbf{1}_{F^{c}}\Phi
\sigma \right) }{\left\vert J^{\prime }\right\vert }\left\Vert
x-m_{J^{\prime }}\right\Vert _{L^{2}\left( \mathbf{1}_{J^{\prime }}\omega
\right) }\left\Vert \square _{J^{\prime }}^{\omega ,\mathbf{b}^{\ast
}}g\right\Vert _{L^{2}\left( \omega \right) }^{\bigstar } \\
&\equiv &II_{G}+II_{B}\ .
\end{eqnarray*}%
where $g_{F;M}$ denotes the pseudoprojection $g_{F;M}=\sum_{J^{\prime }\in 
\mathcal{C}_{F}^{\mathcal{G},\limfunc{shift}}:\ J^{\prime }\subset M}\square
_{J^{\prime }}^{\omega ,\mathbf{b}^{\ast }}g$.

\smallskip

\textbf{Note}: We could also bound $II_{G}$ by using the decomposition $%
\mathcal{M}_{\left( \mathbf{\rho },\varepsilon \right) -\limfunc{deep},%
\mathcal{G}}\left( F\right) $ of $F$ into certain maximal $\mathcal{G}$%
-intervals, but the `smaller' choice $\mathcal{W}\left( F\right) $ of $%
\mathcal{D}$-intervals is needed for $II_{G}$ in order to bound it by the
corresponding functional energy constant $\mathfrak{F}_{\alpha }^{\mathbf{b}%
^{\ast }}$, which can then be controlled by the energy and Muckenhoupt
constants in Appendix B.

\smallskip

Then from Cauchy-Schwarz, the functional energy condition, and 
\begin{equation*}
\left\Vert \Phi \right\Vert _{L^{2}\left( \sigma \right) }^{2}\leq
\sum_{F\in \mathcal{F}}\alpha _{\mathcal{F}}\left( F\right) ^{2}\left\vert
F\right\vert _{\sigma }\lesssim \left\Vert f\right\Vert _{L^{2}\left( \sigma
\right) }^{2}\ ,
\end{equation*}%
we obtain%
\begin{eqnarray*}
\left\vert II_{G}\right\vert &\leq &\left( \sum_{F\in \mathcal{F}}\sum_{M\in 
\mathcal{W}\left( F\right) }\left( \frac{\mathrm{P}^{\alpha }\left( M,%
\mathbf{1}_{F^{c}}\Phi \sigma \right) }{\left\vert M\right\vert }\right)
^{2}\left\Vert \mathsf{Q}_{\mathcal{C}_{F;M}^{\mathcal{G},\limfunc{shift}%
}}^{\omega ,\mathbf{b}^{\ast }}x\right\Vert _{L^{2}\left( \omega \right)
}^{\spadesuit 2}\right) ^{\frac{1}{2}} \\
&&\times \left( \sum_{F\in \mathcal{F}}\sum_{M\in \mathcal{W}\left( F\right)
}\left\Vert g_{F;M}\right\Vert _{L^{2}\left( \omega \right) }^{\bigstar
2}\right) ^{\frac{1}{2}} \\
&\lesssim &\mathfrak{F}_{\alpha }^{\mathbf{b}^{\ast }}\left\Vert \Phi
\right\Vert _{L^{2}\left( \sigma \right) }\left[ \sum_{F\in \mathcal{F}%
}\left\Vert g_{F}\right\Vert _{L^{2}\left( \omega \right) }^{\bigstar 2}%
\right] ^{\frac{1}{2}}\lesssim \mathfrak{F}_{\alpha }^{\mathbf{b}^{\ast
}}\left\Vert f\right\Vert _{L^{2}\left( \sigma \right) }\left\Vert
g\right\Vert _{L^{2}\left( \omega \right) },
\end{eqnarray*}%
by the pairwise disjointedness of the coronas $\mathcal{C}_{F;M}^{\mathcal{G}%
,\limfunc{shift}}$ jointly in $F$ and $M$, which in turn follows from the
pairwise disjointedness (\ref{tau overlap}) of the shifted coronas $\mathcal{%
C}_{F}^{\mathcal{G},\limfunc{shift}}$ in $F$, together with the pairwise
disjointedness of the cubes $M$. Thus we obtain the pairwise disjointedness
of both of the pseudoprojections $\mathsf{P}_{\mathcal{C}_{F;M}^{\mathcal{G},%
\limfunc{shift}}}^{\omega ,\mathbf{b}^{\ast }}$ and $\mathsf{Q}_{\mathcal{C}%
_{F;M}^{\mathcal{G},\limfunc{shift}}}^{\omega ,\mathbf{b}^{\ast }}$ jointly
in $F$ and $M$.

In term $II_{B}$ the quantities $\left\Vert x-m_{J^{\prime }}\right\Vert
_{L^{2}\left( \mathbf{1}_{J^{\prime }}\omega \right) }^{2}$ are no longer
additive except when the intervals $J^{\prime }$ are pairwise disjoint. As a
result we will use (\ref{Haar trick}) in the form,%
\begin{eqnarray}
\sum_{J^{\prime }\subset J}\left( \frac{\mathrm{P}_{1+\delta }^{\alpha
}\left( J^{\prime },\nu \right) }{\left\vert J^{\prime }\right\vert }\right)
^{2}\left\Vert x-m_{J^{\prime }}\right\Vert _{L^{2}\left( \mathbf{1}%
_{J^{\prime }}\omega \right) }^{2} &\lesssim &\frac{1}{\gamma ^{2\delta
^{\prime }}}\left( \frac{\mathrm{P}_{1+\delta ^{\prime }}^{\alpha }\left(
J,\nu \right) }{\left\vert J\right\vert }\right) ^{2}\sum_{J^{\prime \prime
}\subset J}\left\Vert \bigtriangleup _{J^{\prime \prime }}^{\omega
}x\right\Vert _{L^{2}\left( \omega \right) }^{2}  \label{Haar trick'} \\
&\lesssim &\left( \frac{\mathrm{P}_{1+\delta ^{\prime }}^{\alpha }\left(
J,\nu \right) }{\left\vert J\right\vert }\right) ^{2}\left\Vert
x-m_{J}\right\Vert _{L^{2}\left( \mathbf{1}_{J}\omega \right) }^{2}\ , 
\notag
\end{eqnarray}%
and exploit the decay in the Poisson integral $\mathrm{P}_{1+\delta ^{\prime
}}^{\alpha }$ along with weak goodness of the intervals $J$. As a
consequence we will be able to bound $II_{B}$ \emph{directly} by the strong
energy condition (\ref{strong b* energy}), without having to invoke the more
difficult functional energy condition. For the decay we compute that for $%
J\in \mathcal{M}_{\left( \mathbf{\rho },\varepsilon \right) -\limfunc{deep},%
\mathcal{G}}\left( F\right) $%
\begin{eqnarray*}
\frac{\mathrm{P}_{1+\delta ^{\prime }}^{\alpha }\left( J,\mathbf{1}%
_{F^{c}}\Phi \sigma \right) }{\left\vert J\right\vert } &\approx
&\int_{F^{c}}\frac{\left\vert J\right\vert ^{\delta ^{\prime }}}{\left\vert
y-c_{J}\right\vert ^{2+\delta ^{\prime }-\alpha }}\Phi \left( y\right)
d\sigma \left( y\right) \\
&\leq &\sum_{t=0}^{\infty }\int_{\pi _{\mathcal{F}}^{t+1}F\setminus \pi _{%
\mathcal{F}}^{t}F}\left( \frac{\left\vert J\right\vert }{\limfunc{dist}%
\left( c_{J},\left( \pi _{\mathcal{F}}^{t}F\right) ^{c}\right) }\right)
^{\delta ^{\prime }}\frac{1}{\left\vert y-c_{J}\right\vert ^{2-\alpha }}\Phi
\left( y\right) d\sigma \left( y\right) \\
&\lesssim &\sum_{t=0}^{\infty }\left( \frac{\left\vert J\right\vert }{%
\limfunc{dist}\left( c_{J},\left( \pi _{\mathcal{F}}^{t}F\right) ^{c}\right) 
}\right) ^{\delta ^{\prime }}\frac{\mathrm{P}^{\alpha }\left( J,\mathbf{1}%
_{\pi _{\mathcal{F}}^{t+1}F\setminus \pi _{\mathcal{F}}^{t}F}\Phi \sigma
\right) }{\left\vert J\right\vert },
\end{eqnarray*}%
and then use the weak goodness inequality%
\begin{equation*}
\limfunc{dist}\left( c_{J},\left( \pi _{\mathcal{F}}^{t}F\right) ^{c}\right)
\geq 2\ell \left( \pi _{\mathcal{F}}^{t}F\right) ^{1-\varepsilon }\ell
\left( J\right) ^{\varepsilon }\geq 2\cdot 2^{t\left( 1-\varepsilon \right)
}\ell \left( F\right) ^{1-\varepsilon }\ell \left( J\right) ^{\varepsilon
}\geq 2^{t\left( 1-\varepsilon \right) +1}\ell \left( J\right) ,
\end{equation*}%
to conclude that%
\begin{eqnarray}
\left( \frac{\mathrm{P}_{1+\delta ^{\prime }}^{\alpha }\left( J,\mathbf{1}%
_{F^{c}}\Phi \sigma \right) }{\left\vert J\right\vert }\right) ^{2}
&\lesssim &\left( \sum_{t=0}^{\infty }2^{-t\delta ^{\prime }\left(
1-\varepsilon \right) }\frac{\mathrm{P}^{\alpha }\left( J,\mathbf{1}_{\pi _{%
\mathcal{F}}^{t+1}F\setminus \pi _{\mathcal{F}}^{t}F}\Phi \sigma \right) }{%
\left\vert J\right\vert }\right) ^{2}  \label{decay in t} \\
&\lesssim &\sum_{t=0}^{\infty }2^{-t\delta ^{\prime }\left( 1-\varepsilon
\right) }\left( \frac{\mathrm{P}^{\alpha }\left( J,\mathbf{1}_{\pi _{%
\mathcal{F}}^{t+1}F\setminus \pi _{\mathcal{F}}^{t}F}\Phi \sigma \right) }{%
\left\vert J\right\vert }\right) ^{2}.  \notag
\end{eqnarray}%
Now we first apply Cauchy-Schwarz and (\ref{Haar trick'}) to obtain%
\begin{eqnarray*}
II_{B} &=&\sum_{F\in \mathcal{F}}\sum_{J\in \mathcal{M}_{\left( \mathbf{\rho 
},\varepsilon \right) -\limfunc{deep},\mathcal{G}}\left( F\right)
}\sum_{J^{\prime }\in \mathcal{C}_{F;J}^{\mathcal{G},\limfunc{shift}}}\frac{%
\mathrm{P}_{1+\delta }^{\alpha }\left( J^{\prime },\mathbf{1}_{F^{c}}\Phi
\sigma \right) }{\left\vert J^{\prime }\right\vert }\left\Vert
x-m_{J^{\prime }}\right\Vert _{L^{2}\left( \mathbf{1}_{J^{\prime }}\omega
\right) }\left\Vert \square _{J^{\prime }}^{\omega ,\mathbf{b}^{\ast
}}g\right\Vert _{L^{2}\left( \omega \right) }^{\bigstar } \\
&\leq &\left( \sum_{F\in \mathcal{F}}\sum_{J\in \mathcal{M}_{\left( \mathbf{%
\rho },\varepsilon \right) -\limfunc{deep},\mathcal{G}}\left( F\right)
}\sum_{J^{\prime }\in \mathcal{C}_{F;J}^{\mathcal{G},\limfunc{shift}}}\left( 
\frac{\mathrm{P}_{1+\delta }^{\alpha }\left( J^{\prime },\mathbf{1}%
_{F^{c}}\Phi \sigma \right) }{\left\vert J^{\prime }\right\vert }\right)
^{2}\left\Vert x-m_{J^{\prime }}\right\Vert _{L^{2}\left( \mathbf{1}%
_{J^{\prime }}\omega \right) }^{2}\right) ^{\frac{1}{2}}\left[
\sum_{F}\left\Vert g_{F}\right\Vert _{L^{2}\left( \omega \right) }^{\bigstar
2}\right] ^{\frac{1}{2}} \\
&\leq &\left( \sum_{F\in \mathcal{F}}\sum_{J\in \mathcal{M}_{\left( \mathbf{%
\rho },\varepsilon \right) -\limfunc{deep},\mathcal{G}}\left( F\right)
}\left( \frac{\mathrm{P}_{1+\delta ^{\prime }}^{\alpha }\left( J,\mathbf{1}%
_{F^{c}}\Phi \sigma \right) }{\left\vert J\right\vert }\right)
^{2}\left\Vert x-m_{J}\right\Vert _{L^{2}\left( \mathbf{1}_{J}\omega \right)
}^{2}\right) ^{\frac{1}{2}}\left\Vert g\right\Vert _{L^{2}\left( \omega
\right) } \\
&\equiv &\sqrt{II_{\limfunc{energy}}}\left\Vert g\right\Vert _{L^{2}\left(
\omega \right) },
\end{eqnarray*}%
and it remains to estimate $II_{\limfunc{energy}}$. From (\ref{decay in t})
and the strong energy condition (\ref{strong b* energy}), we have%
\begin{eqnarray*}
&&II_{\limfunc{energy}}=\sum_{F\in \mathcal{F}}\sum_{J\in \mathcal{M}%
_{\left( \mathbf{\rho },\varepsilon \right) -\limfunc{deep},\mathcal{G}%
}\left( F\right) }\left( \frac{\mathrm{P}_{1+\delta ^{\prime }}^{\alpha
}\left( J,\mathbf{1}_{F^{c}}\Phi \sigma \right) }{\left\vert J\right\vert }%
\right) ^{2}\left\Vert x-m_{J}\right\Vert _{L^{2}\left( \mathbf{1}_{J}\omega
\right) }^{2} \\
&\leq &\sum_{F\in \mathcal{F}}\sum_{J\in \mathcal{M}_{\left( \mathbf{\rho }%
,\varepsilon \right) -\limfunc{deep},\mathcal{G}}\left( F\right)
}\sum_{t=0}^{\infty }2^{-t\delta ^{\prime }\left( 1-\varepsilon \right)
}\left( \frac{\mathrm{P}^{\alpha }\left( J,\mathbf{1}_{\pi _{\mathcal{F}%
}^{t+1}F\setminus \pi _{\mathcal{F}}^{t}F}\Phi \sigma \right) }{\left\vert
J\right\vert }\right) ^{2}\left\Vert x-m_{J}\right\Vert _{L^{2}\left( 
\mathbf{1}_{J}\omega \right) }^{2} \\
&=&\sum_{t=0}^{\infty }2^{-t\delta ^{\prime }\left( 1-\varepsilon \right)
}\sum_{G\in \mathcal{F}}\sum_{F\in \mathfrak{C}_{\mathcal{F}}^{\left(
t+1\right) }\left( G\right) }\sum_{J\in \mathcal{M}_{\left( \mathbf{\rho }%
,\varepsilon \right) -\limfunc{deep},\mathcal{G}}\left( F\right) }\left( 
\frac{\mathrm{P}^{\alpha }\left( J,\mathbf{1}_{G\setminus \pi _{\mathcal{F}%
}^{t}F}\Phi \sigma \right) }{\left\vert J\right\vert }\right) ^{2}\left\Vert
x-m_{J}\right\Vert _{L^{2}\left( \mathbf{1}_{J}\omega \right) }^{2} \\
&\lesssim &\sum_{t=0}^{\infty }2^{-t\delta ^{\prime }\left( 1-\varepsilon
\right) }\sum_{G\in \mathcal{F}}\alpha _{\mathcal{F}}\left( G\right)
^{2}\sum_{F\in \mathfrak{C}_{\mathcal{F}}^{\left( t+1\right) }\left(
G\right) }\sum_{J\in \mathcal{M}_{\left( \mathbf{\rho },\varepsilon \right) -%
\limfunc{deep}}\left( F\right) }\left( \frac{\mathrm{P}^{\alpha }\left( J,%
\mathbf{1}_{G\setminus \pi _{\mathcal{F}}^{t}F}\sigma \right) }{\left\vert
J\right\vert }\right) ^{2}\left\Vert x-m_{J}\right\Vert _{L^{2}\left( 
\mathbf{1}_{J}\omega \right) }^{2} \\
&\lesssim &\sum_{t=0}^{\infty }2^{-t\delta ^{\prime }\left( 1-\varepsilon
\right) }\sum_{G\in \mathcal{F}}\alpha _{\mathcal{F}}\left( G\right)
^{2}\left( \mathcal{E}_{2}^{\alpha }\right) ^{2}\left\vert G\right\vert
_{\sigma }\lesssim \left( \mathcal{E}_{2}^{\alpha }\right) ^{2}\left\Vert
f\right\Vert _{L^{2}\left( \sigma \right) }^{2}.
\end{eqnarray*}

This completes the proof of the Intertwining Proposition \ref{strongly
adapted}.
\end{proof}

The task of controlling functional energy is taken up in Appendix B below.

\subsection{Paraproduct, neighbour and broken forms}

In this subsection we reduce boundedness of the local below form $\mathsf{B}%
_{\Subset _{\mathbf{r},\varepsilon }}^{A}\left( f,g\right) $ defined in (\ref%
{def local}) to boundedness of the associated stopping form%
\begin{equation}
\mathsf{B}_{\limfunc{stop}}^{A}\left( f,g\right) \equiv \sum_{\substack{ %
I\in \mathcal{C}_{A}^{\mathcal{D}}\text{ and }J\in \mathcal{C}_{A}^{\mathcal{%
G},\limfunc{shift}}  \\ J^{\maltese }\subsetneqq I\text{ and }\ell \left(
J\right) \leq 2^{-\mathbf{r}}\ell \left( I\right) }}\left( E_{I_{J}}^{\sigma
}\widehat{\square }_{I}^{\sigma ,\flat ,\mathbf{b}}f\right) \left\langle
T_{\sigma }^{\alpha }\left( \mathbf{1}_{A\setminus I_{J}}b_{A}\right)
,\square _{J}^{\omega ,\mathbf{b}^{\ast }}g\right\rangle _{\omega }\ ,
\label{def stop}
\end{equation}%
where the modified difference $\widehat{\square }_{I}^{\sigma ,\flat ,%
\mathbf{b}}$ must be carefully chosen (see (\ref{flat box}) and (\ref{factor
b_A}) in Appendix A) in order to control the corresponding paraproduct form
below. Indeed, below we will decompose 
\begin{equation*}
\mathsf{B}_{\Subset _{\mathbf{r},\varepsilon }}^{A}\left( f,g\right) =%
\mathsf{B}_{\limfunc{paraproduct}}^{A}\left( f,g\right) -\mathsf{B}_{%
\limfunc{stop}}^{A}\left( f,g\right) +\mathsf{B}_{\limfunc{neighbour}%
}^{A}\left( f,g\right) +\mathsf{B}_{\limfunc{broken}}^{A}\left( f,g\right) ,
\end{equation*}%
and then prove in (\ref{est para}), (\ref{est neigh}) and (\ref{broken
vanish}) the estimate%
\begin{eqnarray*}
\left\vert \mathsf{B}_{\Subset _{\mathbf{r},\varepsilon }}^{A}\left(
f,g\right) +\mathsf{B}_{\limfunc{stop}}^{A}\left( f,g\right) \right\vert
&\leq &\left\vert \mathsf{B}_{\limfunc{paraproduct}}^{A}\left( f,g\right)
\right\vert +\left\vert \mathsf{B}_{\limfunc{neighbour}}^{A}\left(
f,g\right) \right\vert +\left\vert \mathsf{B}_{\limfunc{broken}}^{A}\left(
f,g\right) \right\vert \\
&\lesssim &\mathfrak{T}_{T^{\alpha }}^{\mathbf{b}}\ \alpha _{\mathcal{A}%
}\left( A\right) \ \sqrt{\left\vert A\right\vert _{\sigma }}\ \left\Vert 
\mathsf{P}_{\mathcal{C}_{A}^{\mathcal{G},\limfunc{shift}}}^{\omega ,\mathbf{b%
}^{\ast }}g\right\Vert _{L^{2}\left( \omega \right) }^{\bigstar } \\
&&+\sqrt{\mathfrak{A}_{2}^{\alpha }}\left( \left\Vert \mathsf{P}_{\mathcal{C}%
_{A}}^{\sigma }f\right\Vert _{L^{2}(\sigma )}^{\bigstar }+\sqrt{%
\sum_{A^{\prime }\in \mathfrak{C}_{A}\left( A\right) }\left\vert A^{\prime
}\right\vert _{\sigma }\alpha _{\mathcal{A}}\left( A^{\prime }\right) ^{2}}%
\right) \left\Vert \mathsf{P}_{\mathcal{C}_{A}^{\mathcal{G},\limfunc{shift}%
}}^{\omega ,\mathbf{b}^{\ast }}g\right\Vert _{L^{2}(\omega )}^{\bigstar }\ ,
\end{eqnarray*}%
which can of course then be summed in $A\in \mathcal{A}$ to conclude that%
\begin{eqnarray*}
&&\sum_{A\in \mathcal{A}}\left\vert \mathsf{B}_{\Subset _{\mathbf{r}%
,\varepsilon }}^{A}\left( f,g\right) +\mathsf{B}_{\limfunc{stop}}^{A}\left(
f,g\right) \right\vert \\
&\lesssim &\left( \mathfrak{T}_{T^{\alpha }}^{\mathbf{b}}+\sqrt{\mathfrak{A}%
_{2}^{\alpha }}\right) \sqrt{\sum_{A\in \mathcal{A}}\left\{ \alpha _{%
\mathcal{A}}\left( A\right) ^{2}\left\vert A\right\vert _{\sigma
}+\left\Vert \mathsf{P}_{\mathcal{C}_{A}}^{\sigma }f\right\Vert
_{L^{2}(\sigma )}^{\bigstar 2}+\sum_{A^{\prime }\in \mathfrak{C}_{\mathcal{A}%
}\left( A\right) }\alpha _{\mathcal{A}}\left( A^{\prime }\right)
^{2}\left\vert A^{\prime }\right\vert _{\sigma }\right\} }\sqrt{\sum_{A\in 
\mathcal{A}}\left\Vert \mathsf{P}_{\mathcal{C}_{A}^{\mathcal{G},\limfunc{%
shift}}}^{\omega ,\mathbf{b}^{\ast }}g\right\Vert _{L^{2}\left( \omega
\right) }^{\bigstar 2}} \\
&\lesssim &\left( \mathfrak{T}_{T^{\alpha }}^{\mathbf{b}}+\sqrt{\mathfrak{A}%
_{2}^{\alpha }}\right) \left\Vert f\right\Vert _{L^{2}\left( \sigma \right)
}\left\Vert g\right\Vert _{L^{2}\left( \omega \right) }\ .
\end{eqnarray*}%
The stopping form is the subject of the section following this one.

Note from (\ref{factor b_A}) and (\ref{telescoping})\ in Appendix A, that
the modified dual martingale differences $\square _{I}^{\sigma ,\flat ,%
\mathbf{b}}$ and $\widehat{\square }_{I}^{\sigma ,\flat ,\mathbf{b}}$, 
\begin{equation*}
\square _{I}^{\sigma ,\flat ,\mathbf{b}}f\equiv \square _{I}^{\sigma ,%
\mathbf{b}}f-\sum_{I^{\prime }\in \mathfrak{C}_{\limfunc{broken}}\left(
I\right) }\mathbb{F}_{I^{\prime }}^{\sigma ,\mathbf{b}}f=b_{A}\sum_{I^{%
\prime }\in \mathfrak{C}\left( I\right) }\mathbf{1}_{I^{\prime
}}E_{I^{\prime }}^{\sigma }\left( \widehat{\square }_{I}^{\sigma ,\flat ,%
\mathbf{b}}f\right) =b_{A}\widehat{\square }_{I}^{\sigma ,\flat ,\mathbf{b}%
}f,
\end{equation*}%
satisfy the following telescoping property for all $K\in \left( \mathcal{C}%
_{A}\setminus \left\{ A\right\} \right) \cup \left( \bigcup_{A^{\prime }\in 
\mathfrak{C}_{\mathcal{A}}\left( A\right) }A^{\prime }\right) $ and $L\in 
\mathcal{C}_{A}$ with $K\subset L$:%
\begin{equation*}
\sum_{I:\ \pi K\subset I\subset L}E_{I_{K}}^{\sigma }\left( \widehat{\square 
}_{I}^{\sigma ,\flat ,\mathbf{b}}f\right) =\left\{ 
\begin{array}{ccc}
-E_{L}^{\sigma }\widehat{\mathbb{F}}_{L}^{\sigma ,\mathbf{b}}f & \text{ if }
& K\in \mathfrak{C}_{\mathcal{A}}\left( A\right) \\ 
E_{K}^{\sigma }\widehat{\mathbb{F}}_{K}^{\sigma ,\mathbf{b}}f-E_{L}^{\sigma }%
\widehat{\mathbb{F}}_{L}^{\sigma ,\mathbf{b}}f & \text{ if } & K\in \mathcal{%
C}_{A}%
\end{array}%
\right. .
\end{equation*}%
Fix $I\in \mathcal{C}_{A}$ for the moment. We will use%
\begin{eqnarray*}
\mathbf{1}_{I} &=&\mathbf{1}_{I_{J}}+\mathbf{1}_{\theta \left( I_{J}\right)
}\ , \\
\mathbf{1}_{I_{J}} &=&\mathbf{1}_{A}-\mathbf{1}_{A\setminus I_{J}}\ ,
\end{eqnarray*}%
where $\theta \left( I_{J}\right) \in \mathfrak{C}_{\mathcal{D}}\left(
I\right) \setminus \left\{ I_{J}\right\} $ is the $\mathcal{D}$-child of $I$
other than the child $I_{J}$ that contains $J$. We begin with the splitting%
\begin{eqnarray*}
&&\left\langle T_{\sigma }^{\alpha }\square _{I}^{\sigma ,\mathbf{b}%
}f,\square _{J}^{\omega ,\mathbf{b}^{\ast }}g\right\rangle _{\omega
}=\left\langle T_{\sigma }^{\alpha }\left( \mathbf{1}_{I_{J}}\square
_{I}^{\sigma ,\mathbf{b}}f\right) ,\square _{J}^{\omega ,\mathbf{b}^{\ast
}}g\right\rangle _{\omega }+\left\langle T_{\sigma }^{\alpha }\left( \mathbf{%
1}_{\theta \left( I_{J}\right) }\square _{I}^{\sigma ,\mathbf{b}}f\right)
,\square _{J}^{\omega ,\mathbf{b}^{\ast }}g\right\rangle _{\omega } \\
&=&\left\langle T_{\sigma }^{\alpha }\left( \mathbf{1}_{I_{J}}\square
_{I}^{\sigma ,\flat ,\mathbf{b}}f\right) ,\square _{J}^{\omega ,\mathbf{b}%
^{\ast }}g\right\rangle _{\omega }+\left\langle T_{\sigma }^{\alpha }\left( 
\mathbf{1}_{I_{J}}\sum_{I^{\prime }\in \mathfrak{C}_{\limfunc{broken}}\left(
I\right) }\mathbb{F}_{I^{\prime }}^{\sigma ,\mathbf{b}}f\right) ,\square
_{J}^{\omega ,\mathbf{b}^{\ast }}g\right\rangle _{\omega }+\left\langle
T_{\sigma }^{\alpha }\left( \mathbf{1}_{\theta \left( I_{J}\right) }\square
_{I}^{\sigma ,\mathbf{b}}f\right) ,\square _{J}^{\omega ,\mathbf{b}^{\ast
}}g\right\rangle _{\omega } \\
&\equiv &I+II+III\ .
\end{eqnarray*}%
From (\ref{factor b_A}) we have%
\begin{eqnarray*}
I &=&\left\langle T_{\sigma }^{\alpha }\left( \mathbf{1}_{I_{J}}\square
_{I}^{\sigma ,\flat ,\mathbf{b}}f\right) ,\square _{J}^{\omega ,\mathbf{b}%
^{\ast }}g\right\rangle _{\omega }=\left\langle T_{\sigma }^{\alpha }\left[
b_{A}\left( \mathbf{1}_{I_{J}}\widehat{\square }_{I}^{\sigma ,\flat ,\mathbf{%
b}}f\right) \right] ,\square _{J}^{\omega ,\mathbf{b}^{\ast }}g\right\rangle
_{\omega } \\
&=&E_{I_{J}}^{\sigma }\left( \widehat{\square }_{I}^{\sigma ,\flat ,\mathbf{b%
}}f\right) \left\langle T_{\sigma }^{\alpha }\left( \mathbf{1}%
_{I_{J}}b_{A}\right) ,\square _{J}^{\omega ,\mathbf{b}^{\ast
}}g\right\rangle _{\omega } \\
&=&E_{I_{J}}^{\sigma }\left( \widehat{\square }_{I}^{\sigma ,\flat ,\mathbf{b%
}}f\right) \left\langle T_{\sigma }^{\alpha }b_{A},\square _{J}^{\omega ,%
\mathbf{b}^{\ast }}g\right\rangle _{\omega }-E_{I_{J}}^{\sigma }\left( 
\widehat{\square }_{I}^{\sigma ,\flat ,\mathbf{b}}f\right) \left\langle
T_{\sigma }^{\alpha }\left( \mathbf{1}_{A\setminus I_{J}}b_{A}\right)
,\square _{J}^{\omega ,\mathbf{b}^{\ast }}g\right\rangle _{\omega }\ .
\end{eqnarray*}%
Since the function $\mathbb{F}_{I_{J}}^{\sigma ,\mathbf{b}}f$ is a constant
multiple of $b_{I_{J}}$ on $I_{J}$, we can define $\widehat{\mathbb{F}}%
_{I_{J}}^{\sigma ,\mathbf{b}}f\equiv \frac{1}{b_{I_{J}}}\mathbb{F}%
_{I_{J}}^{\sigma ,\mathbf{b}}f$ (or simply use the $PLBP$ we are assuming)
and then%
\begin{equation*}
II=\left\langle T_{\sigma }^{\alpha }\left( \mathbf{1}_{I_{J}}\sum_{I^{%
\prime }\in \mathfrak{C}_{\limfunc{broken}}\left( I\right) }\mathbb{F}%
_{I^{\prime }}^{\sigma ,\mathbf{b}}f\right) ,\square _{J}^{\omega ,\mathbf{b}%
^{\ast }}g\right\rangle _{\omega }=\mathbf{1}_{\mathfrak{C}_{\mathcal{A}%
}\left( A\right) }\left( I_{J}\right) \ E_{I_{J}}^{\sigma }\left( \widehat{%
\mathbb{F}}_{I_{J}}^{\sigma ,\mathbf{b}}f\right) \ \left\langle T_{\sigma
}^{\alpha }b_{I_{J}},\square _{J}^{\omega ,\mathbf{b}^{\ast }}g\right\rangle
_{\omega }\ ,
\end{equation*}%
where the presence of the indicator function $\mathbf{1}_{\mathfrak{C}_{%
\mathcal{A}}\left( A\right) }\left( I_{J}\right) $ simply means that term $%
II $ vanishes unless $I_{J}$ is an $\mathcal{A}$-child of $A$. We now write
these terms as%
\begin{eqnarray*}
\left\langle T_{\sigma }^{\alpha }\square _{I}^{\sigma ,\mathbf{b}}f,\square
_{J}^{\omega ,\mathbf{b}^{\ast }}g\right\rangle _{\omega }
&=&E_{I_{J}}^{\sigma }\left( \widehat{\square }_{I}^{\sigma ,\flat ,\mathbf{b%
}}f\right) \left\langle T_{\sigma }^{\alpha }b_{A},\square _{J}^{\omega ,%
\mathbf{b}^{\ast }}g\right\rangle _{\omega } \\
&&-E_{I_{J}}^{\sigma }\left( \widehat{\square }_{I}^{\sigma ,\flat ,\mathbf{b%
}}f\right) \left\langle T_{\sigma }^{\alpha }\left( \mathbf{1}_{A\setminus
I_{J}}b_{A}\right) ,\square _{J}^{\omega ,\mathbf{b}^{\ast }}g\right\rangle
_{\omega } \\
&&+\left\langle T_{\sigma }^{\alpha }\left( \mathbf{1}_{\theta \left(
I_{J}\right) }\square _{I}^{\sigma ,\mathbf{b}}f\right) ,\square
_{J}^{\omega ,\mathbf{b}^{\ast }}g\right\rangle _{\omega } \\
&&+\mathbf{1}_{\left\{ I_{J}\in \mathfrak{C}_{\mathcal{A}}\left( A\right)
\right\} }\ E_{I_{J}}^{\sigma }\left( \widehat{\mathbb{F}}_{I_{J}}^{\sigma ,%
\mathbf{b}}f\right) \ \left\langle T_{\sigma }^{\alpha }b_{I_{J}},\square
_{J}^{\omega ,\mathbf{b}^{\ast }}g\right\rangle _{\omega }\ ,
\end{eqnarray*}%
where the four lines are respectively a paraproduct, stopping, neighbour and
broken term.

The corresponding NTV splitting of $\mathsf{B}_{\Subset _{\mathbf{r}%
,\varepsilon }}^{A}\left( f,g\right) $ using (\ref{def local}) and (\ref{def
shorthand}) becomes%
\begin{eqnarray*}
\mathsf{B}_{\Subset _{\mathbf{r},\varepsilon }}^{A}\left( f,g\right)
&=&\left\langle T_{\sigma }^{\alpha }\left( \mathsf{P}_{\mathcal{C}%
_{A}}^{\sigma }f\right) ,\mathsf{P}_{\mathcal{C}_{A}^{\mathcal{G},\limfunc{%
shift}}}^{\omega }g\right\rangle _{\omega }^{\Subset _{\mathbf{r}%
,\varepsilon }}=\sum_{\substack{ I\in \mathcal{C}_{A}\text{ and }J\in 
\mathcal{C}_{A}^{\mathcal{G},\limfunc{shift}}  \\ J^{\maltese }\subsetneqq I%
\text{ and }\ell \left( J\right) \leq 2^{-\mathbf{r}}\ell \left( I\right) }}%
\left\langle T_{\sigma }^{\alpha }\left( \square _{I}^{\sigma ,\mathbf{b}%
}f\right) ,\square _{J}^{\omega ,\mathbf{b}^{\ast }}g\right\rangle _{\omega }
\\
&=&\mathsf{B}_{\limfunc{paraproduct}}^{A}\left( f,g\right) -\mathsf{B}_{%
\limfunc{stop}}^{A}\left( f,g\right) +\mathsf{B}_{\limfunc{neighbour}%
}^{A}\left( f,g\right) +\mathsf{B}_{\limfunc{broken}}^{A}\left( f,g\right) ,
\end{eqnarray*}%
where%
\begin{eqnarray*}
\mathsf{B}_{\limfunc{paraproduct}}^{A}\left( f,g\right) &\equiv &\sum 
_{\substack{ I\in \mathcal{C}_{A}\text{ and }J\in \mathcal{C}_{A}^{\mathcal{G%
},\limfunc{shift}}  \\ J^{\maltese }\subsetneqq I\text{ and }\ell \left(
J\right) \leq 2^{-\mathbf{r}}\ell \left( I\right) }}E_{I_{J}}^{\sigma
}\left( \widehat{\square }_{I}^{\sigma ,\flat ,\mathbf{b}}f\right)
\left\langle T_{\sigma }^{\alpha }b_{A},\square _{J}^{\omega ,\mathbf{b}%
^{\ast }}g\right\rangle _{\omega }\ , \\
\mathsf{B}_{\limfunc{stop}}^{A}\left( f,g\right) &\equiv &\sum_{\substack{ %
I\in \mathcal{C}_{A}\text{ and }J\in \mathcal{C}_{A}^{\mathcal{G},\limfunc{%
shift}}  \\ J^{\maltese }\subsetneqq I\text{ and }\ell \left( J\right) \leq
2^{-\mathbf{r}}\ell \left( I\right) }}E_{I_{J}}^{\sigma }\left( \widehat{%
\square }_{I}^{\sigma ,\flat ,\mathbf{b}}f\right) \left\langle T_{\sigma
}^{\alpha }\left( \mathbf{1}_{A\setminus I_{J}}b_{A}\right) ,\square
_{J}^{\omega ,\mathbf{b}^{\ast }}g\right\rangle _{\omega }\ , \\
\mathsf{B}_{\limfunc{neighbour}}^{A}\left( f,g\right) &\equiv &\sum 
_{\substack{ I\in \mathcal{C}_{A}\text{ and }J\in \mathcal{C}_{A}^{\mathcal{G%
},\limfunc{shift}}  \\ J^{\maltese }\subsetneqq I\text{ and }\ell \left(
J\right) \leq 2^{-\mathbf{r}}\ell \left( I\right) }}\left\langle T_{\sigma
}^{\alpha }\left( \mathbf{1}_{\theta \left( I_{J}\right) }\square
_{I}^{\sigma ,\mathbf{b}}f\right) ,\square _{J}^{\omega ,\mathbf{b}^{\ast
}}g\right\rangle _{\omega }\ ,
\end{eqnarray*}%
correspond to the three original NTV forms associated with $1$-testing, and
where 
\begin{equation}
\mathsf{B}_{\limfunc{broken}}^{A}\left( f,g\right) \equiv \sum_{\substack{ %
I\in \mathcal{C}_{A}\text{ and }J\in \mathcal{C}_{A}^{\mathcal{G},\limfunc{%
shift}}  \\ J^{\maltese }\subsetneqq I\text{ and }\ell \left( J\right) \leq
2^{-\mathbf{r}}\ell \left( I\right) }}\mathbf{1}_{\left\{ I_{J}\in \mathfrak{%
C}_{\mathcal{A}}\left( A\right) \right\} }\ E_{I_{J}}^{\sigma }\left( 
\widehat{\mathbb{F}}_{I_{J}}^{\sigma ,\mathbf{b}}f\right) \ \left\langle
T_{\sigma }^{\alpha }b_{I_{J}},\square _{J}^{\omega ,\mathbf{b}^{\ast
}}g\right\rangle _{\omega }=0\ ,  \label{broken vanish}
\end{equation}%
since $J^{\maltese }\subsetneqq I$ and $I_{J}\in \mathfrak{C}_{\mathcal{A}%
}\left( A\right) $ imply that $J^{\maltese }\notin \mathcal{C}_{A}^{\mathcal{%
G},\limfunc{shift}}$, contradicting $J\in \mathcal{C}_{A}^{\mathcal{G},%
\limfunc{shift}}$.

\begin{remark}
The inquisitive reader will note that the pairs $\left( I,J\right) $ arising
in the above sum with $J^{\maltese }\subsetneqq I$ replaced by $J^{\maltese
}=I$ are handled in the probabilistic estimate (\ref{HM bad}) for the bad
form $\Theta _{2}^{\limfunc{bad}\natural }$ defined in (\ref{Theta_2^bad
sharp}).
\end{remark}

\subsubsection{The paraproduct form}

The paraproduct form $\mathsf{B}_{\limfunc{paraproduct}}^{A}\left(
f,g\right) $ is easily controlled by the testing condition for $T^{\alpha }$
together with weak Riesz inequalities for dual martingale differences.
Indeed, recalling the telescoping identity (\ref{telescoping}), and that the
collection $\left\{ I\in \mathcal{C}_{A}\text{:\ }\ell \left( J\right) \leq
2^{-\mathbf{r}}\ell \left( I\right) \right\} $ is tree connected for all $%
J\in \mathcal{C}_{A}^{\mathcal{G},\limfunc{shift}}$, we have%
\begin{eqnarray*}
\mathsf{B}_{\limfunc{paraproduct}}^{A}\left( f,g\right) &=&\sum_{\substack{ %
I\in \mathcal{C}_{A}\text{ and }J\in \mathcal{C}_{A}^{\mathcal{G},\limfunc{%
shift}}  \\ J^{\maltese }\subsetneqq I\text{ and }\ell \left( J\right) \leq
2^{-\mathbf{r}}\ell \left( I\right) }}E_{I_{J}}^{\sigma }\left( \widehat{%
\square }_{I}^{\sigma ,\flat ,\mathbf{b}}f\right) \left\langle T_{\sigma
}^{\alpha }b_{A},\square _{J}^{\omega ,\mathbf{b}^{\ast }}g\right\rangle
_{\omega } \\
&=&\sum_{J\in \mathcal{C}_{A}^{\mathcal{G},\limfunc{shift}}}\left\langle
T_{\sigma }^{\alpha }b_{A},\square _{J}^{\omega ,\mathbf{b}^{\ast
}}g\right\rangle _{\omega }\left\{ \sum_{I\in \mathcal{C}_{A}\text{:\ }%
J^{\maltese }\subsetneqq I\text{ and }\ell \left( J\right) \leq 2^{-\mathbf{r%
}}\ell \left( I\right) }E_{I_{J}}^{\sigma }\left( \widehat{\square }%
_{I}^{\sigma ,\flat ,\mathbf{b}}f\right) \right\} \\
&=&\sum_{J\in \mathcal{C}_{A}^{\mathcal{G},\limfunc{shift}}}\left\langle
T_{\sigma }^{\alpha }b_{A},\square _{J}^{\omega ,\mathbf{b}^{\ast
}}g\right\rangle _{\omega }\left\{ \mathbf{1}_{\left\{ J:I^{\natural }\left(
J\right) _{J}\in \mathcal{C}_{A}\right\} }E_{I^{\natural }\left( J\right)
_{J}}^{\sigma }\widehat{\mathbb{F}}_{I^{\natural }\left( J\right)
_{J}}^{\sigma ,\mathbf{b}}f-E_{A}^{\sigma }\widehat{\mathbb{F}}_{A}^{\sigma ,%
\mathbf{b}}f\right\} \\
&=&\left\langle T_{\sigma }^{\alpha }b_{A},\sum_{J\in \mathcal{C}_{A}^{%
\mathcal{G},\limfunc{shift}}}\left\{ \mathbf{1}_{\left\{ J:I^{\natural
}\left( J\right) _{J}\in \mathcal{C}_{A}\right\} }E_{I^{\natural }\left(
J\right) _{J}}^{\sigma }\widehat{\mathbb{F}}_{I^{\natural }\left( J\right)
_{J}}^{\sigma ,\mathbf{b}}f-E_{A}^{\sigma }\widehat{\mathbb{F}}_{A}^{\sigma ,%
\mathbf{b}}f\right\} \square _{J}^{\omega ,\mathbf{b}^{\ast }}g\right\rangle
_{\omega }\ ,
\end{eqnarray*}%
where $I^{\natural }\left( J\right) $ denotes the smallest interval $I\in 
\mathcal{C}_{A}$ such that $J^{\maltese }\subsetneqq I$ and $\ell \left(
J\right) \leq 2^{-\mathbf{r}}\ell \left( I\right) $, and of course $%
I^{\natural }\left( J\right) _{J}$ denotes its child containing $J$. Note
that by construction of the modified difference operator $\square
_{I}^{\sigma ,\flat ,\mathbf{b}}$, the only time the average $\widehat{%
\mathbb{F}}_{I^{\natural }\left( J\right) _{J}}^{\sigma }f$ appears in the
above sum is when $I^{\natural }\left( J\right) _{J}\in \mathcal{C}_{A}$,
since the case $I^{\natural }\left( J\right) _{J}\in \mathcal{A}$ has been
removed to the broken term. This is reflected above with the inclusion of
the indicator $\mathbf{1}_{\left\{ J:I^{\natural }\left( J\right) _{J}\in 
\mathcal{C}_{A}\right\} }$. It follows that we have the bound $\left\vert 
\mathbf{1}_{\left\{ J:I^{\natural }\left( J\right) _{J}\in \mathcal{C}%
_{A}\right\} }E_{I^{\natural }\left( J\right) _{J}}^{\sigma }\widehat{%
\mathbb{F}}_{I^{\natural }\left( J\right) _{J}}^{\sigma ,\mathbf{b}%
}f\right\vert +\left\vert E_{A}^{\sigma }\widehat{\mathbb{F}}_{A}^{\sigma ,%
\mathbf{b}}f\right\vert \lesssim E_{A}^{\sigma }\left\vert f\right\vert \leq
\alpha _{\mathcal{A}}\left( A\right) $.

Thus from Cauchy-Schwarz, the upper weak Riesz inequalities Proposition \ref%
{half Riesz} for the pseudoprojections $\square _{J}^{\omega ,\mathbf{b}%
^{\ast }}g$ and the bound on the coefficients $\lambda _{J}\equiv \left( 
\mathbf{1}_{\left\{ J:I^{\natural }\left( J\right) _{J}\in \mathcal{C}%
_{A}\right\} }E_{I^{\natural }\left( J\right) _{J}}^{\sigma }\widehat{%
\mathbb{F}}_{I^{\natural }\left( J\right) _{J}}^{\sigma ,\mathbf{b}%
}f-E_{A}^{\sigma }\widehat{\mathbb{F}}_{A}^{\sigma ,\mathbf{b}}f\right) $
given by $\left\vert \lambda _{J}\right\vert \lesssim \alpha _{\mathcal{A}%
}\left( A\right) $, we have%
\begin{eqnarray}
&&  \label{est para} \\
\left\vert \mathsf{B}_{\limfunc{paraproduct}}^{A}\left( f,g\right)
\right\vert &=&\left\vert \left\langle T_{\sigma }^{\alpha }b_{A},\sum_{J\in 
\mathcal{C}_{A}^{\mathcal{G},\limfunc{shift}}}\left\{ \left( \mathbf{1}%
_{\left\{ J:I^{\natural }\left( J\right) _{J}\in \mathcal{C}_{A}\right\}
}E_{I^{\natural }\left( J\right) _{J}}^{\sigma }\widehat{\mathbb{F}}%
_{I^{\natural }\left( J\right) _{J}}^{\sigma ,\mathbf{b}}f-E_{A}^{\sigma }%
\widehat{\mathbb{F}}_{A}^{\sigma ,\mathbf{b}}f\right) \right\} \square
_{J}^{\omega ,\mathbf{b}^{\ast }}g\right\rangle _{\omega }\right\vert  \notag
\\
&\leq &\left\Vert \mathbf{1}_{A}T_{\sigma }^{\alpha }b_{A}\right\Vert
_{L^{2}\left( \omega \right) }\left\Vert \sum_{J\in \mathcal{C}_{A}^{%
\mathcal{G},\limfunc{shift}}}\lambda _{J}\square _{J}^{\omega ,\mathbf{b}%
^{\ast }}g\right\Vert _{L^{2}\left( \omega \right) }  \notag \\
&\lesssim &\alpha _{\mathcal{A}}\left( A\right) \ \left\Vert \mathbf{1}%
_{A}T_{\sigma }^{\alpha }b_{A}\right\Vert _{L^{2}\left( \omega \right) }\
\left\Vert \sum_{J\in \mathcal{C}_{A}^{\mathcal{G},\limfunc{shift}}}\square
_{J}^{\omega ,\mathbf{b}^{\ast }}g\right\Vert _{L^{2}\left( \omega \right)
}^{\bigstar }  \notag \\
&\leq &\mathfrak{T}_{T^{\alpha }}^{\mathbf{b}}\ \alpha _{\mathcal{A}}\left(
A\right) \ \sqrt{\left\vert A\right\vert _{\sigma }}\ \left\Vert \mathsf{P}_{%
\mathcal{C}_{A}^{\mathcal{G},\limfunc{shift}}}^{\omega ,\mathbf{b}^{\ast
}}g\right\Vert _{L^{2}\left( \omega \right) }^{\bigstar }.  \notag
\end{eqnarray}

\subsubsection{The neighbour form}

Next, the neighbour form $\mathsf{B}_{\limfunc{neighbour}}^{A}\left(
f,g\right) $ is easily controlled by the $\mathfrak{A}_{2}^{\alpha }$
condition using the pivotal estimate in Energy Lemma \ref{ener} and the fact
that the intervals $J\in \mathcal{C}_{A}^{\mathcal{G},\limfunc{shift}}$ are
good in $I$ and beyond when the pair $\left( I,J\right) $ occurs in the sum.
In particular, the information encoded in the stopping tree $\mathcal{A}$
plays no role here, apart from appearing in the corona projections on the
right hand side of (\ref{est neigh}) below. We have%
\begin{equation}
\mathsf{B}_{\limfunc{neighbour}}^{A}\left( f,g\right) =\sum_{\substack{ I\in 
\mathcal{C}_{A}\text{ and }J\in \mathcal{C}_{A}^{\mathcal{G},\limfunc{shift}%
}  \\ J^{\maltese }\subsetneqq I\text{ and }\ell \left( J\right) \leq 2^{-%
\mathbf{r}}\ell \left( I\right) }}\left\langle T_{\sigma }^{\alpha }\left( 
\mathbf{1}_{\theta \left( I_{J}\right) }\square _{I}^{\sigma ,\mathbf{b}%
}f\right) ,\square _{J}^{\omega ,\mathbf{b}^{\ast }}g\right\rangle _{\omega
},  \label{def neighbour}
\end{equation}%
where we keep in mind that the pairs $\left( I,J\right) \in \mathcal{D}%
\times \mathcal{G}$ that arise in the sum for $\mathsf{B}_{\limfunc{neighbour%
}}^{A}\left( f,g\right) $ satisfy the property that $J^{\maltese
}\subsetneqq I$, so that $J$ is good with respect to all intervals $K$ of
size at least that of $J^{\maltese }$, which includes $I$. Recall that $%
I_{J} $ is the child of $I$ that contains $J$, and that $\theta \left(
I_{J}\right) $ denotes its sibling in $I$, i.e. $\theta \left( I_{J}\right)
\in \mathfrak{C}_{\mathcal{D}}\left( I\right) \setminus \left\{
I_{J}\right\} $. Fix $\left( I,J\right) $ momentarily, and an integer $s\geq 
\mathbf{r}$. Using $\square _{I}^{\sigma ,\mathbf{b}}=\square _{I}^{\sigma
,\flat ,\mathbf{b}}+\square _{I,\limfunc{broken}}^{\sigma ,\flat ,\mathbf{b}%
} $ and the fact that $\square _{I}^{\sigma ,\flat ,\mathbf{b}}f$ is a
constant multiple of $b_{\theta \left( I_{J}\right) }$ on the interval $%
\theta \left( I_{J}\right) $, we have the estimates 
\begin{eqnarray*}
\left\vert \mathbf{1}_{\theta \left( I_{J}\right) }\square _{I}^{\sigma
,\flat ,\mathbf{b}}f\right\vert &=&\left\vert \left( E_{\theta \left(
I_{J}\right) }^{\sigma }\widehat{\square }_{I}^{\sigma ,\flat ,\mathbf{b}%
}f\right) b_{\theta \left( I_{J}\right) }\right\vert \leq C_{\mathbf{b}%
}\left\vert E_{\theta \left( I_{J}\right) }^{\sigma }\widehat{\square }%
_{I}^{\sigma ,\flat ,\mathbf{b}}f\right\vert , \\
\left\vert \mathbf{1}_{\theta \left( I_{J}\right) }\square _{I,\limfunc{%
broken}}^{\sigma ,\flat ,\mathbf{b}}f\right\vert &\leq &\mathbf{1}_{%
\mathfrak{C}_{A}\left( A\right) }\left( \theta \left( I_{J}\right) \right) \
E_{\theta \left( I_{J}\right) }^{\sigma }\left\vert f\right\vert ,
\end{eqnarray*}%
and hence%
\begin{equation}
\mathbf{1}_{\theta \left( I_{J}\right) }\left\vert \square _{I}^{\sigma ,%
\mathbf{b}}f\right\vert \leq C\mathbf{1}_{\theta \left( I_{J}\right) }\left(
\left\vert E_{\theta \left( I_{J}\right) }^{\sigma }\widehat{\square }%
_{I}^{\sigma ,\flat ,\mathbf{b}}f\right\vert +\mathbf{1}_{\mathfrak{C}%
_{A}\left( A\right) }\left( \theta \left( I_{J}\right) \right) \ E_{\theta
\left( I_{J}\right) }^{\sigma }\left\vert f\right\vert \right) ,
\label{box bound}
\end{equation}%
which will be used below after an application of the Energy Lemma. We can
write%
\begin{equation*}
\mathsf{B}_{\limfunc{neighbour}}^{A}\left( f,g\right) =\sum_{\substack{ I\in 
\mathcal{C}_{A}\text{ and }J\in \mathcal{G}_{\left( \kappa \left(
I_{J},J\right) ,\varepsilon \right) -\limfunc{good}}^{\mathcal{D}}\cap 
\mathcal{C}_{A}^{\mathcal{G},\limfunc{shift}}\text{ and }J^{\maltese
}\subsetneqq I  \\ d\left( J,\theta \left( I_{J}\right) \right) >2\ell
\left( J\right) ^{\varepsilon }\ell \left( \theta \left( I_{J}\right)
\right) ^{1-\varepsilon }\text{ and }\ell \left( J\right) \leq 2^{-\mathbf{r}%
}\ell \left( I\right) }}\left\langle T_{\sigma }^{\alpha }\left( \mathbf{1}%
_{\theta \left( I_{J}\right) }\square _{I}^{\sigma ,\mathbf{b}}f\right)
,\square _{J}^{\omega ,\mathbf{b}^{\ast }}g\right\rangle _{\omega }
\end{equation*}%
where we have included the conditions $J\in \mathcal{G}_{\left( \kappa
\left( I_{J},J\right) ,\varepsilon \right) -\limfunc{good}}^{\mathcal{D}}$
and $d\left( J,\theta \left( I_{J}\right) \right) >2\ell \left( J\right)
^{\varepsilon }\ell \left( \theta \left( I_{J}\right) \right)
^{1-\varepsilon }$ in the summation since they are already implied the
remaining four conditions, and will be used in estimates below.

We will also use the following fractional analogue of the Poisson inequality
in \cite{Vol}.

\begin{lemma}
\label{Poisson inequality}Suppose $0\leq \alpha <1$ and $J\subset I\subset K$
and that $d\left( J,\partial I\right) >2\ell \left( J\right) ^{\varepsilon
}\ell \left( I\right) ^{1-\varepsilon }$ for some $0<\varepsilon <\frac{1}{%
2-\alpha }$. Then for a positive Borel measure $\mu $ we have%
\begin{equation}
\mathrm{P}^{\alpha }(J,\mu \mathbf{1}_{K\setminus I})\lesssim \left( \frac{%
\ell \left( J\right) }{\ell \left( I\right) }\right) ^{1-\varepsilon \left(
2-\alpha \right) }\mathrm{P}^{\alpha }(I,\mu \mathbf{1}_{K\setminus I}).
\label{e.Jsimeq}
\end{equation}
\end{lemma}

\begin{proof}
We have%
\begin{equation*}
\mathrm{P}^{\alpha }\left( J,\mu \mathbf{1}_{K\setminus I}\right) \approx
\sum_{k=0}^{\infty }2^{-k}\frac{1}{\left\vert 2^{k}J\right\vert ^{1-\alpha }}%
\int_{\left( 2^{k}J\right) \cap \left( K\setminus I\right) }d\mu ,
\end{equation*}%
and $\left( 2^{k}J\right) \cap \left( K\setminus I\right) \neq \emptyset $
requires%
\begin{equation*}
d\left( J,K\setminus I\right) \leq c2^{k}\ell \left( J\right) ,
\end{equation*}%
for some dimensional constant $c>0$. Let $k_{0}$ be the smallest such $k$.
By our distance assumption we must then have%
\begin{equation*}
2\ell \left( J\right) ^{\varepsilon }\ell \left( I\right) ^{1-\varepsilon
}\leq d\left( J,\partial I\right) \leq c2^{k_{0}}\ell \left( J\right) ,
\end{equation*}%
or%
\begin{equation*}
2^{-k_{0}+1}\leq c\left( \frac{\ell \left( J\right) }{\ell \left( I\right) }%
\right) ^{1-\varepsilon }.
\end{equation*}%
Now let $k_{1}$ be defined by $2^{k_{1}}\equiv \frac{\ell \left( I\right) }{%
\ell \left( J\right) }$. Then assuming $k_{1}>k_{0}$ (the case $k_{1}\leq
k_{0}$ is similar) we have%
\begin{eqnarray*}
\mathrm{P}^{\alpha }\left( J,\mu \mathbf{1}_{K\setminus I}\right) &\approx
&\left\{ \sum_{k=k_{0}}^{k_{1}}+\sum_{k=k_{1}}^{\infty }\right\} 2^{-k}\frac{%
1}{\left\vert 2^{k}J\right\vert ^{1-\alpha }}\int_{\left( 2^{k}J\right) \cap
\left( K\setminus I\right) }d\mu \\
&\lesssim &2^{-k_{0}}\frac{\left\vert I\right\vert ^{1-\alpha }}{\left\vert
2^{k_{0}}J\right\vert ^{1-\alpha }}\left( \frac{1}{\left\vert I\right\vert
^{1-\alpha }}\int_{\left( 2^{k_{1}}J\right) \cap \left( K\setminus I\right)
}d\mu \right) +2^{-k_{1}}\mathrm{P}^{\alpha }\left( I,\mu \mathbf{1}%
_{K\setminus I}\right) \\
&\lesssim &\left( \frac{\ell \left( J\right) }{\ell \left( I\right) }\right)
^{\left( 1-\varepsilon \right) \left( 2-\alpha \right) }\left( \frac{\ell
\left( I\right) }{\ell \left( J\right) }\right) ^{1-\alpha }\mathrm{P}%
^{\alpha }\left( I,\mu \mathbf{1}_{K\setminus I}\right) +\frac{\ell \left(
J\right) }{\ell \left( I\right) }\mathrm{P}^{\alpha }\left( I,\mu \mathbf{1}%
_{K\setminus I}\right) ,
\end{eqnarray*}%
which is the inequality (\ref{e.Jsimeq}).
\end{proof}

Now fix $I_{0}=I_{J},I_{\theta }=\theta \left( I_{J}\right) \in \mathfrak{C}%
_{\mathcal{D}}\left( I\right) $ and assume that $J\Subset _{\mathbf{r}%
,\varepsilon }I_{0}$. Let $\frac{\ell \left( J\right) }{\ell \left(
I_{0}\right) }=2^{-s}$ in the pivotal estimate from Energy Lemma \ref{ener}
with $J\subset I_{0}\subset I$ to obtain 
\begin{align*}
& \left\vert \langle T_{\sigma }^{\alpha }\left( \mathbf{1}_{\theta \left(
I_{J}\right) }\square _{I}^{\sigma ,\mathbf{b}}f\right) ,\square
_{J}^{\omega ,\mathbf{b}^{\ast }}g\rangle _{\omega }\right\vert \\
& \lesssim \left\Vert \square _{J}^{\omega ,\mathbf{b}^{\ast }}g\right\Vert
_{L^{2}\left( \omega \right) }\sqrt{\left\vert J\right\vert _{\omega }}%
\mathrm{P}^{\alpha }\left( J,\mathbf{1}_{\theta \left( I_{J}\right)
}\left\vert \square _{I}^{\sigma ,\mathbf{b}}f\right\vert \sigma \right) \\
& \lesssim \left\Vert \square _{J}^{\omega ,\mathbf{b}^{\ast }}g\right\Vert
_{L^{2}\left( \omega \right) }\sqrt{\left\vert J\right\vert _{\omega }}\cdot
2^{-\left( 1-\varepsilon \left( 2-\alpha \right) \right) s}\mathrm{P}%
^{\alpha }\left( I_{0},\mathbf{1}_{\theta \left( I_{J}\right) }\left\vert
\square _{I}^{\sigma ,\mathbf{b}}f\right\vert \sigma \right) \\
& \lesssim \left\Vert \square _{J}^{\omega ,\mathbf{b}^{\ast }}g\right\Vert
_{L^{2}\left( \omega \right) }\sqrt{\left\vert J\right\vert _{\omega }}\cdot
2^{-\left( 1-\varepsilon \left( 2-\alpha \right) \right) s}\mathrm{P}%
^{\alpha }\left( I_{0},\mathbf{1}_{\theta \left( I_{J}\right) }\left(
\left\vert E_{\theta \left( I_{J}\right) }^{\sigma }\widehat{\square }%
_{I}^{\sigma ,\flat ,\mathbf{b}}f\right\vert +\mathbf{1}_{\mathfrak{C}%
_{A}\left( A\right) }\left( \theta \left( I_{J}\right) \right) \ E_{\theta
\left( I_{J}\right) }^{\sigma }\left\vert f\right\vert \right) \sigma
\right) .
\end{align*}%
Here we are using (\ref{e.Jsimeq}) in the third line, which applies since $%
J\subset I_{0}$, and we have used (\ref{box bound}) in the fourth line. It
will be convenient to use the shorthand notation%
\begin{equation*}
\mathbf{E}_{\theta \left( I_{J}\right) }^{\sigma }f\equiv \left\vert
E_{\theta \left( I_{J}\right) }^{\sigma }\widehat{\square }_{I}^{\sigma
,\flat ,\mathbf{b}}f\right\vert +\mathbf{1}_{\mathfrak{C}_{A}\left( A\right)
}\left( \theta \left( I_{J}\right) \right) \ E_{\theta \left( I_{J}\right)
}^{\sigma }\left\vert f\right\vert
\end{equation*}%
where the intervals $I$ and $I_{J}$ on the right hand side are determined
uniquely by the interval $\theta \left( I_{J}\right) $.

In the sum below, we keep the side lengths of the intervals $J$ fixed at $%
2^{-s}$ times that of $I_{0}$, and of course take $J\subset I_{0}$. We also
keep the underlying assumptions that $J\in \mathcal{C}_{A}^{\mathcal{G},%
\limfunc{shift}}$ and that $J\in \mathcal{G}_{\left( \kappa \left(
I_{J},J\right) ,\varepsilon \right) -\limfunc{good}}^{\mathcal{D}}$ in mind
without necessarily pointing to them in the notation. Matters will shortly
be reduced to estimating the following term: 
\begin{align*}
A(I,I_{0},I_{\theta },s)& \equiv \sum_{J\;:\;2^{s+1}\ell \left( J\right)
=\ell \left( I\right) :J\subset I_{0}}\left\vert \langle T_{\sigma }^{\alpha
}\left( \mathbf{1}_{I_{\theta }}\square _{I}^{\sigma ,\mathbf{b}}f\right)
,\square _{J}^{\omega ,\mathbf{b}^{\ast }}g\rangle _{\omega }\right\vert \\
& \leq 2^{-\left( 1-\varepsilon \left( 2-\alpha \right) \right) s}\left( 
\mathbf{E}_{\theta \left( I_{J}\right) }^{\sigma }f\right) \ \mathrm{P}%
^{\alpha }(I_{0},\mathbf{1}_{\theta \left( I_{J}\right) }\sigma
)\sum_{J\;:\;2^{s+1}\ell \left( J\right) =\ell \left( I\right) :\ J\subset
I_{0}}\left\Vert \square _{J}^{\omega ,\mathbf{b}^{\ast }}g\right\Vert
_{L^{2}\left( \omega \right) }\sqrt{\left\vert J\right\vert _{\omega }} \\
& \leq 2^{-\left( 1-\varepsilon \left( 2-\alpha \right) \right) s}\left( 
\mathbf{E}_{\theta \left( I_{J}\right) }^{\sigma }f\right) \ \mathrm{P}%
^{\alpha }(I_{0},\mathbf{1}_{\theta \left( I_{J}\right) }\sigma )\sqrt{%
\left\vert I_{0}\right\vert _{\omega }}\Lambda (I,I_{0},I_{\theta },s), \\
& \text{where }\Lambda (I,I_{0},I_{\theta },s)^{2}\equiv \sum_{J\in \mathcal{%
C}_{A}^{\mathcal{G},\limfunc{shift}}:\;2^{s+1}\ell \left( J\right) =\ell
\left( I\right) :\ J\subset I_{0}}\left\Vert \square _{J}^{\omega ,\mathbf{b}%
^{\ast }}g\right\Vert _{L^{2}\left( \omega \right) }^{2}\,.
\end{align*}%
The last line follows upon using the Cauchy-Schwarz inequality and the fact
that $J\in \mathcal{C}_{A}^{\mathcal{G},\limfunc{shift}}$. We also note that
since $2^{s+1}\ell \left( J\right) =\ell \left( I\right) $, 
\begin{eqnarray}
\sum_{I_{0}\in \mathfrak{C}_{\mathcal{D}}\left( I\right) }\Lambda
(I,I_{0},I_{\theta },s)^{2} &\equiv &\sum_{J\in \mathcal{C}_{A}^{\mathcal{G},%
\limfunc{shift}}:\;2^{s+1}\ell \left( J\right) =\ell \left( I\right) :\
J\subset I}\left\Vert \square _{J}^{\omega ,\mathbf{b}^{\ast }}g\right\Vert
_{L^{2}\left( \omega \right) }^{2}\ ;  \label{g} \\
\sum_{I\in \mathcal{C}_{A}}\sum_{I_{0}\in \mathfrak{C}_{\mathcal{D}}\left(
I\right) }\Lambda (I,I_{0},I_{\theta },s)^{2} &\leq &\left\Vert \mathsf{P}_{%
\mathcal{C}_{A}^{\mathcal{G},\limfunc{shift}}}^{\omega ,\mathbf{b}^{\ast
}}g\right\Vert _{L^{2}(\omega )}^{\bigstar 2}\ .  \notag
\end{eqnarray}

Using (\ref{box hat bound}) we obtain 
\begin{equation}
\left\vert E_{I_{\theta }}^{\sigma }\left( \widehat{\square }_{I}^{\sigma
,\flat ,\mathbf{b}}f\right) \right\vert \leq \sqrt{E_{I_{\theta }}^{\sigma
}\left\vert \widehat{\square }_{I}^{\sigma ,\flat ,\mathbf{b}}f\right\vert
^{2}}\lesssim \left\Vert \square _{I}^{\sigma ,\mathbf{b}}f\right\Vert
_{L^{2}\left( \sigma \right) }^{\bigstar }\ \left\vert I_{\theta
}\right\vert _{\sigma }^{-\frac{1}{2}},  \label{e.haarAvg}
\end{equation}%
and hence%
\begin{equation*}
\mathbf{E}_{\theta \left( I_{J}\right) }^{\sigma }f\equiv \left\vert
E_{\theta \left( I_{J}\right) }^{\sigma }\widehat{\square }_{I}^{\sigma
,\flat ,\mathbf{b}}f\right\vert +\mathbf{1}_{\mathfrak{C}_{A}\left( A\right)
}\left( \theta \left( I_{J}\right) \right) \ E_{\theta \left( I_{J}\right)
}^{\sigma }\left\vert f\right\vert \lesssim \left( \left\Vert \square
_{I}^{\sigma ,\mathbf{b}}f\right\Vert _{L^{2}\left( \sigma \right)
}^{\bigstar }\ +\mathbf{1}_{\mathfrak{C}_{A}\left( A\right) }\left( \theta
\left( I_{J}\right) \right) \ \left\vert I_{\theta }\right\vert _{\sigma }^{%
\frac{1}{2}}E_{\theta \left( I_{J}\right) }^{\sigma }\left\vert f\right\vert
\right) \left\vert I_{\theta }\right\vert _{\sigma }^{-\frac{1}{2}},
\end{equation*}%
and we can thus estimate $A(I,I_{0},I_{\theta },s)$ as follows: 
\begin{eqnarray*}
&&A(I,I_{0},I_{\theta },s) \\
&\lesssim &2^{-\left( 1-\varepsilon \left( 2-\alpha \right) \right) s}\left(
\left\Vert \square _{I}^{\sigma ,\mathbf{b}}f\right\Vert _{L^{2}\left(
\sigma \right) }^{\bigstar }+\mathbf{1}_{\mathfrak{C}_{A}\left( A\right)
}\left( I_{\theta }\right) \ \left\vert I_{\theta }\right\vert _{\sigma }^{%
\frac{1}{2}}E_{I_{\theta }}^{\sigma }\left\vert f\right\vert \right) \Lambda
(I,I_{0},I_{\theta },s)\cdot \left\vert I_{\theta }\right\vert _{\sigma }^{-%
\frac{1}{2}}\mathrm{P}^{\alpha }(I_{0},\mathbf{1}_{\theta \left(
I_{J}\right) }\sigma )\sqrt{\left\vert I_{0}\right\vert _{\omega }} \\
&\lesssim &\sqrt{\mathfrak{A}_{2}^{\alpha }}2^{-\left( 1-\varepsilon \left(
2-\alpha \right) \right) s}\left( \left\Vert \square _{I}^{\sigma ,\mathbf{b}%
}f\right\Vert _{L^{2}\left( \sigma \right) }^{\bigstar }+\mathbf{1}_{%
\mathfrak{C}_{A}\left( A\right) }\left( I_{\theta }\right) \ \left\vert
I_{\theta }\right\vert _{\sigma }^{\frac{1}{2}}E_{I_{\theta }}^{\sigma
}\left\vert f\right\vert \right) \Lambda (I,I_{0},I_{\theta },s)\,,
\end{eqnarray*}%
since $\mathrm{P}^{\alpha }(I_{0},\mathbf{1}_{\theta \left( I_{J}\right)
}\sigma )\lesssim \frac{\left\vert I_{\theta }\right\vert _{\sigma }}{%
\left\vert I_{\theta }\right\vert ^{1-\alpha }}$ shows that 
\begin{equation*}
\left\vert I_{\theta }\right\vert _{\sigma }^{-\frac{1}{2}}\mathrm{P}%
^{\alpha }(I_{0},\mathbf{1}_{\theta \left( I_{J}\right) }\sigma )\ \sqrt{%
\left\vert I_{0}\right\vert _{\omega }}\lesssim \frac{\sqrt{\left\vert
I_{\theta }\right\vert _{\sigma }}\sqrt{\left\vert I_{0}\right\vert _{\omega
}}}{\left\vert I\right\vert ^{1-\alpha }}\lesssim \sqrt{\mathfrak{A}%
_{2}^{\alpha }}.
\end{equation*}

An application of Cauchy-Schwarz to the sum over $I$ using (\ref{g}) then
shows that 
\begin{eqnarray*}
&&\sum_{I\in \mathcal{C}_{A}}\sum_{\substack{ I_{0},I_{\theta }\in \mathfrak{%
C}_{\mathcal{D}}\left( I\right)  \\ I_{0}\neq I_{\theta }}}%
A(I,I_{0},I_{\theta },s) \\
&\lesssim &\sqrt{\mathfrak{A}_{2}^{\alpha }}2^{-\left( 1-\varepsilon \left(
2-\alpha \right) \right) s}\sqrt{\sum_{I\in \mathcal{C}_{A}}\left\Vert
\square _{I}^{\sigma ,\mathbf{b}}f\right\Vert _{L^{2}\left( \sigma \right)
}^{\bigstar 2}+\sum_{I_{\theta }\in \mathfrak{C}_{A}\left( A\right)
}\left\vert I_{\theta }\right\vert _{\sigma }\left( E_{I_{\theta }}^{\sigma
}\left\vert f\right\vert \right) ^{2}}\sqrt{\sum_{I\in \mathcal{C}%
_{A}}\left( \sum_{\substack{ I_{0},I_{\theta }\in \mathfrak{C}_{\mathcal{D}%
}\left( I\right)  \\ I_{0}\neq I_{\theta }}}\Lambda (I,I_{0},I_{\theta
},s)\right) ^{2}} \\
&\lesssim &\sqrt{\mathfrak{A}_{2}^{\alpha }}2^{-\left( 1-\varepsilon \left(
2-\alpha \right) \right) s}\sqrt{\left\Vert \mathsf{P}_{\mathcal{C}%
_{A}}^{\sigma }f\right\Vert _{L^{2}(\sigma )}^{\bigstar 2}+\sum_{A^{\prime
}\in \mathfrak{C}_{A}\left( A\right) }\left\vert A^{\prime }\right\vert
_{\sigma }\left( E_{A^{\prime }}^{\sigma }\left\vert f\right\vert \right)
^{2}}\sqrt{\sum_{I\in \mathcal{C}_{A}}\left( \sum_{\substack{ I_{0}\in 
\mathfrak{C}_{\mathcal{D}}\left( I\right)  \\ I_{0}\neq I_{\theta }}}\Lambda
(I,I_{0},I_{\theta },s)\right) ^{2}} \\
&\lesssim &\sqrt{\mathfrak{A}_{2}^{\alpha }}2_{{}}^{-\left( 1-\varepsilon
\left( 2-\alpha \right) \right) s}\left( \lVert \mathsf{P}_{\mathcal{C}%
_{A}}^{\sigma }f\rVert _{L^{2}(\sigma )}^{\bigstar }+\sqrt{\sum_{A^{\prime
}\in \mathfrak{C}_{A}\left( A\right) }\left\vert A^{\prime }\right\vert
_{\sigma }\left( E_{A^{\prime }}^{\sigma }\left\vert f\right\vert \right)
^{2}}\right) \left\Vert \mathsf{P}_{\mathcal{C}_{A}^{\mathcal{G},\limfunc{%
shift}}}^{\omega ,\mathbf{b}^{\ast }}g\right\Vert _{L^{2}(\omega
)}^{\bigstar }\,.
\end{eqnarray*}%
This estimate is summable in $s\geq \mathbf{r}$ since $\varepsilon <\frac{1}{%
2-\alpha }$, and so the proof of 
\begin{eqnarray}
\left\vert \mathsf{B}_{\limfunc{neighbour}}^{A}\left( f,g\right) \right\vert
&=&\left\vert \sum_{\substack{ I\in \mathcal{C}_{A}\text{ and }J\in \mathcal{%
C}_{A}^{\mathcal{G},\limfunc{shift}}  \\ J^{\maltese }\subsetneqq I\text{
and }\ell \left( J\right) \leq 2^{-\mathbf{r}}\ell \left( I\right) }}%
\left\langle T_{\sigma }^{\alpha }\left( \mathbf{1}_{\theta \left(
I_{J}\right) }\square _{I}^{\sigma ,\mathbf{b}}f\right) ,\square
_{J}^{\omega ,\mathbf{b}^{\ast }}g\right\rangle _{\omega }\right\vert
\label{est neigh} \\
&\leq &\sum_{I\in \mathcal{C}_{A}}\sum_{\substack{ I_{0},I_{\theta }\in 
\mathfrak{C}_{\mathcal{D}}\left( I\right)  \\ I_{0}\neq I_{\theta }}}\sum_{s=%
\mathbf{r}}^{\infty }A(I,I_{0},I_{\theta },s)  \notag \\
&\lesssim &\sqrt{\mathfrak{A}_{2}^{\alpha }}\left( \left\Vert \mathsf{P}_{%
\mathcal{C}_{A}}^{\sigma }f\right\Vert _{L^{2}(\sigma )}^{\bigstar }+\sqrt{%
\sum_{A^{\prime }\in \mathfrak{C}_{A}\left( A\right) }\left\vert A^{\prime
}\right\vert _{\sigma }\alpha _{\mathcal{A}}\left( A^{\prime }\right) ^{2}}%
\right) \left\Vert \mathsf{P}_{\mathcal{C}_{A}^{\mathcal{G},\limfunc{shift}%
}}^{\omega ,\mathbf{b}^{\ast }}g\right\Vert _{L^{2}(\omega )}^{\bigstar } 
\notag
\end{eqnarray}%
is complete since $E_{A^{\prime }}^{\sigma }\left\vert f\right\vert \lesssim
\alpha _{\mathcal{A}}\left( A^{\prime }\right) $.

\section{Stopping form\label{Sec stop}}

In this section, we modify our adaptation in \cite{SaShUr7}, \cite{SaShUr9}
and \cite{SaShUr10}\footnote{%
And correct an error in \cite{SaShUr7} related to the restricted norms of
stopping forms for admissible collections.} of the argument of M. Lacey in 
\cite{Lac} to apply in the setting of a local $Tb$ theorem for an $\alpha $%
-fractional Calder\'{o}n-Zygmund operator $T^{\alpha }$ in $\mathbb{R}$
using the Monotonicity Lemma \ref{mono}, the energy condition, and the weak
goodness of Hyt\"{o}nen and Martikainen \cite{HyMa}. Following Lacey in \cite%
{Lac}, we construct $\mathcal{L}\,$-coronas from the `bottom up' with
stopping times involving the energies $\left\Vert \bigtriangleup
_{J}^{\omega ,\mathbf{b}^{\ast }}x\right\Vert _{L^{2}\left( \omega \right)
}^{2}$, but then overlay this with an additional top/down `indented' corona
construction in order to accommodate the weaker goodness of Hyt\"{o}nen and
Martikainen. We directly control the pairs $\left( I,J\right) $ in the
stopping form `essentially' according to the $\mathcal{L}\,$-coronas to
which $I$ and $J^{\maltese }$ are associated, by absorbing the case when
both $I$ and $J^{\maltese }$ belong to the same $\mathcal{L}\,$-corona, and
by using the Straddling and Subtraddling Lemmas to control the case when $I$
and $J^{\maltese }$ lie in different coronas, with a geometric gain coming
from the separation of the coronas in the `indented' construction overlaying
Lacey's bottom/up construction (we actually use the grandchild $J^{\flat }$
of $J^{\maltese }$ that contains $J$ to distinguish aborption cases from
straddling cases). We also use a Corona-straddling Lemma to control certain
extremal pairs $\left( I,J\right) $ that straddle two $\mathcal{A}$-coronas.
As in \cite{Lac}, an Orthogonality Lemma proves useful in all cases.
Finally, since we are using two independent dyadic grids, we must enlarge
the $\limfunc{skeleton}$ of an interval to include an infinite sequence of
points we call the $\limfunc{body}$ of the interval.

Apart from these changes, the remaining modifications are more obvious, such
as

\begin{itemize}
\item the use of the weak goodness of Hyt\"{o}nen and Martikainen \cite{HyMa}
for pairs $\left( I,J\right) $ arising in the stopping form, rather than
goodness for all intervals $I$ and $J$ that was available in \cite{Lac}, 
\cite{SaShUr7}, \cite{SaShUr9} and \cite{SaShUr10}. For the most part
definitions such as admissible collections are modified to require $%
J^{\maltese }\subsetneqq I$. In paricular, Lacey's size functional is
enlarged to include more intervals $K\in \mathcal{D}$ that are not good;

\item the pseudoprojections $\square _{I}^{\sigma ,\mathbf{b}},\square
_{J}^{\omega ,\mathbf{b}^{\ast }}$ and Carleson averaging operators $\nabla
_{I}^{\sigma },\nabla _{J}^{\omega }$ are used in place of the orthogonal
Haar projections, and the frame and weak Riesz inequalities compensate for
the lack of orthogonality.
\end{itemize}

Fix grids $\mathcal{D}$ and $\mathcal{G}$. In Section 6 we reduced matters
to proving (\ref{local est}), i.e.%
\begin{equation}
\left\vert \mathsf{B}_{\limfunc{stop}}^{A}\left( f,g\right) \right\vert
\lesssim \mathcal{NTV}_{\alpha }\left( \left\Vert \mathsf{P}_{\mathcal{C}%
_{A}^{\mathcal{D}}}^{\sigma ,\mathbf{b}}f\right\Vert _{L^{2}\left( \sigma
\right) }^{\bigstar }+\alpha _{\mathcal{A}}\left( A\right) \sqrt{\left\vert
A\right\vert _{\sigma }}\right) \left\Vert \mathsf{P}_{\mathcal{C}_{A}^{%
\mathcal{G},\limfunc{shift}}}^{\omega ,\mathbf{b}^{\ast }}g\right\Vert
_{L^{2}\left( \omega \right) }^{\bigstar }\ ,  \label{B stop form 3}
\end{equation}%
where we recall that $\mathcal{NTV}_{\alpha }$ is defined in (\ref{def NTV}%
), and the nonstandard `norms' are given in Notation \ref{nonstandard norm}
by,%
\begin{eqnarray*}
\left\Vert \mathsf{P}_{\mathcal{C}_{A}^{\mathcal{D}}}^{\sigma ,\mathbf{b}%
}f\right\Vert _{L^{2}\left( \sigma \right) }^{\bigstar 2} &\equiv
&\sum_{I\in \mathcal{C}_{A}^{\mathcal{D}}}\left( \left\Vert \square
_{I}^{\sigma ,\mathbf{b}}f\right\Vert _{L^{2}\left( \sigma \right)
}^{2}+\left\Vert \nabla _{I}^{\sigma }f\right\Vert _{L^{2}\left( \sigma
\right) }^{2}\right) , \\
\left\Vert \mathsf{P}_{\mathcal{C}_{A}^{\mathcal{G},\limfunc{shift}%
}}^{\omega ,\mathbf{b}^{\ast }}g\right\Vert _{L^{2}\left( \omega \right)
}^{\bigstar 2} &\equiv &\sum_{J\in \mathcal{C}_{A}^{\mathcal{G},\limfunc{%
shift}}}\left( \left\Vert \square _{J}^{\omega ,\mathbf{b}^{\ast
}}g\right\Vert _{L^{2}\left( \omega \right) }^{2}+\left\Vert \nabla
_{J}^{\omega }g\right\Vert _{L^{2}\left( \omega \right) }^{2}\right) ,
\end{eqnarray*}%
and that the stopping form is given in (\ref{def stop}) by%
\begin{equation*}
\mathsf{B}_{\limfunc{stop}}^{A}\left( f,g\right) \equiv \sum_{\substack{ %
I\in \mathcal{C}_{A}^{\mathcal{D}}\text{ and }J\in \mathcal{C}_{A}^{\mathcal{%
G},\limfunc{shift}}  \\ J^{\maltese }\subsetneqq I\text{ and }\ell \left(
J\right) \leq 2^{-\mathbf{r}}\ell \left( I\right) }}\left( E_{I_{J}}^{\sigma
}\widehat{\square }_{I}^{\sigma ,\flat ,\mathbf{b}}f\right) \left\langle
T_{\sigma }^{\alpha }\left( \mathbf{1}_{A\setminus I_{J}}b_{A}\right)
,\square _{J}^{\omega ,\mathbf{b}^{\ast }}g\right\rangle _{\omega }\ .
\end{equation*}

It is important to note that $J^{\maltese }\subsetneqq I$ implies $%
J^{\maltese }\subset I_{J}$, and it follows that we \emph{cannot} have $%
I_{J}\in \mathfrak{C}_{\mathcal{A}}\left( A\right) $, i.e. we cannot have
that the child of $I$ containing $J$ is a stopping interval in $\mathcal{A}$%
, since this would then contradict the assumption that $J\in \mathcal{C}%
_{A}^{\mathcal{G},\limfunc{shift}}$. Furthermore, the pair $\left(
I,J\right) =\left( J^{\maltese },J\right) $ does not arise in the sum simply
because of the requirement $J^{\maltese }\subsetneqq I$. For convenience in
notation, and without loss of generality, we now reindex the stopping form
with this in mind by replacing the pairs $\left( I,J\right) $ in the sum
above with new pairs $\left( I^{\prime },J^{\prime }\right) \equiv \left(
I_{J},J\right) $ (recall that the child of $I$ that contains $J$ is denoted $%
I_{J}$). The result is that%
\begin{equation*}
\mathsf{B}_{\limfunc{stop}}^{A}\left( f,g\right) =\sum_{\substack{ I^{\prime
}\in \mathcal{C}_{A}^{\mathcal{D},\limfunc{restrict}}\text{ and }J^{\prime
}\in \mathcal{C}_{A}^{\mathcal{G},\limfunc{shift}}  \\ J^{\maltese }\subset
I^{\prime }\text{ and }\ell \left( J^{\prime }\right) \leq 2^{1-\mathbf{r}%
}\ell \left( I^{\prime }\right) }}\left( E_{I^{\prime }}^{\sigma }\widehat{%
\square }_{\pi I^{\prime }}^{\sigma ,\flat ,\mathbf{b}}f\right) \left\langle
T_{\sigma }^{\alpha }\left( b_{A}\mathbf{1}_{A\setminus I^{\prime }}\right)
,\square _{J^{\prime }}^{\omega ,\mathbf{b}^{\ast }}g\right\rangle _{\omega
},
\end{equation*}%
where 
\begin{equation}
\mathcal{C}_{A}^{\mathcal{D},\limfunc{restrict}}\equiv \mathcal{C}%
_{A}\setminus \left\{ A\right\}  \label{def restrict}
\end{equation}%
is the $A$-corona with its top interval $A$ removed. Now we simply drop the
primes from the dummy variables $I^{\prime }$ and $J^{\prime }$ and relabel $%
1-\mathbf{r}$ as $-\mathbf{r}$ to obtain 
\begin{equation}
\mathsf{B}_{\limfunc{stop}}^{A}\left( f,g\right) =\sum_{\substack{ I\in 
\mathcal{C}_{A}^{\mathcal{D},\limfunc{restrict}}\text{ and }J\in \mathcal{C}%
_{A}^{\mathcal{G},\limfunc{shift}}  \\ J^{\maltese }\subset I\text{ and }%
\ell \left( J\right) \leq 2^{-\mathbf{r}}\ell \left( I\right) }}\left(
E_{I}^{\sigma }\widehat{\square }_{\pi I}^{\sigma ,\flat ,\mathbf{b}%
}f\right) \left\langle T_{\sigma }^{\alpha }\left( b_{A}\mathbf{1}%
_{A\setminus I}\right) ,\square _{J}^{\omega ,\mathbf{b}^{\ast
}}g\right\rangle _{\omega }.  \label{dummy}
\end{equation}

\begin{definition}
Suppose that $A\in \mathcal{A}$ and that $\mathcal{P}\subset \mathcal{C}%
_{A}^{\mathcal{D},\limfunc{restrict}}\times \mathcal{C}_{A}^{\mathcal{G},%
\limfunc{shift}}$. We say that the collection of pairs $\mathcal{P}$ is $A$%
\emph{-admissible} if
\end{definition}

\begin{itemize}
\item (good and $\left( \mathbf{r},\varepsilon \right) $-deeply embedded) $%
\ell \left( J\right) \leq 2^{-\mathbf{r}}\ell \left( I\right) $ and $%
J^{\maltese }\subset I\varsubsetneqq A$ for every $\left( I,J\right) \in 
\mathcal{P},$

\item (tree-connected in the first component) if $I_{1}\subset I_{2}$ and
both $\left( I_{1},J\right) \in \mathcal{P}$ and $\left( I_{2},J\right) \in 
\mathcal{P}$, then $\left( I,J\right) \in \mathcal{P}$ for every $I$ in the
geodesic $\left[ I_{1},I_{2}\right] =\left\{ I\in \mathcal{D}:I_{1}\subset
I\subset I_{2}\right\} $.
\end{itemize}

From now on we often write $\mathcal{C}_{A}$ and $\mathcal{C}_{A}^{\limfunc{%
restrict}}$ in place of $\mathcal{C}_{A}^{\mathcal{D}}$ and $\mathcal{C}%
_{A}^{\mathcal{D},\limfunc{restrict}}$ respectively, i.e. we drop the
superscript $\mathcal{D}$, when there is no possiblility of confusion. The
basic example of an admissible collection of pairs is obtained from the
pairs of intervals summed in the stopping form $\mathsf{B}_{\limfunc{stop}%
}^{A}\left( f,g\right) $ in (\ref{dummy}) which occurs in the inequality (%
\ref{B stop form 3}) above; 
\begin{equation}
\mathcal{P}^{A}\equiv \left\{ \left( I,J\right) :I\in \mathcal{C}_{A}^{%
\limfunc{restrict}}\text{, }J\in \mathcal{C}_{A}^{\mathcal{G},\limfunc{shift}%
}\text{, }J^{\maltese }\subset I\text{ and\ }\ell \left( J\right) \leq 2^{-%
\mathbf{r}}\ell \left( I\right) \right\} .  \label{initial P}
\end{equation}

\begin{definition}
\label{def stop P}Suppose that $A\in \mathcal{A}$ and that $\mathcal{P}$ is
an $A$\emph{-admissible} collection of pairs. Define the associated \emph{%
stopping} form $\mathsf{B}_{\limfunc{stop}}^{A,\mathcal{P}}$ by%
\begin{equation*}
\mathsf{B}_{\limfunc{stop}}^{A,\mathcal{P}}\left( f,g\right) \equiv
\sum_{\left( I,J\right) \in \mathcal{P}}\left( E_{I}^{\sigma }\widehat{%
\square }_{\pi I}^{\sigma ,\flat ,\mathbf{b}}f\right) \ \left\langle
T_{\sigma }^{\alpha }\left( b_{A}\mathbf{1}_{A\setminus I}\right) ,\square
_{J}^{\omega ,\mathbf{b}^{\ast }}g\right\rangle _{\omega }\ .
\end{equation*}%
where $\widehat{\square }_{\pi I}^{\sigma ,\flat ,\mathbf{b}}$ is the
modified dual martingale difference defined in (\ref{flat box}) and (\ref%
{flat box hat}).
\end{definition}

Recall the strong energy condition constant $\mathcal{E}_{2}^{\alpha }$
defined by%
\begin{equation*}
\left( \mathcal{E}_{2}^{\alpha }\right) ^{2}\equiv \sup_{I=\dot{\cup}I_{r}}%
\frac{1}{\left\vert I\right\vert _{\sigma }}\sum_{r=1}^{\infty }\left( \frac{%
\mathrm{P}^{\alpha }\left( I_{r},\mathbf{1}_{I}\sigma \right) }{\left\vert
I_{r}\right\vert }\right) ^{2}\left\Vert x-m_{I_{r}}^{\omega }\right\Vert
_{L^{2}\left( \mathbf{1}_{I_{r}}\omega \right) }^{2}\ .
\end{equation*}

\begin{proposition}
\label{stopping bound}Suppose that $A\in \mathcal{A}$ and that $\mathcal{P}$
is an $A$\emph{-admissible} collection of pairs. Then the stopping form $%
\mathsf{B}_{\limfunc{stop}}^{A,\mathcal{P}}$ satisfies the bound%
\begin{equation}
\left\vert \mathsf{B}_{\limfunc{stop}}^{A,\mathcal{P}}\left( f,g\right)
\right\vert \lesssim \left( \mathcal{E}_{2}^{\alpha }+\sqrt{\mathfrak{A}%
_{2}^{\alpha }}\right) \left\Vert \mathsf{P}_{\mathcal{C}_{A}}^{\sigma ,%
\mathbf{b}}f\right\Vert _{L^{2}\left( \sigma \right) }^{\bigstar }\left\Vert 
\mathsf{P}_{\mathcal{C}_{A}^{\mathcal{G},\limfunc{shift}}}^{\omega ,\mathbf{b%
}^{\ast }}g\right\Vert _{L^{2}\left( \omega \right) }^{\bigstar }\ .
\label{B stop form}
\end{equation}
\end{proposition}

The above proposition proves (\ref{B stop form 3}) - even without the term $%
\alpha _{\mathcal{A}}\left( A\right) \sqrt{\left\vert A\right\vert _{\sigma }%
}$ on the right - with the choice $\mathcal{P}=\mathcal{P}^{A}$. To prove
Proposition \ref{stopping bound}, we begin by letting 
\begin{eqnarray*}
\Pi _{1}\mathcal{P} &\equiv &\left\{ I\in \mathcal{C}_{A}^{\mathcal{D},%
\limfunc{restrict}}:\left( I,J\right) \in \mathcal{P}\text{ for some }J\in 
\mathcal{C}_{A}^{\mathcal{G},\func{shift}}\right\} , \\
\Pi _{2}\mathcal{P} &\equiv &\left\{ J\in \mathcal{C}_{A}^{\mathcal{G},\func{%
shift}}:\left( I,J\right) \in \mathcal{P}\text{ for some }I\in \mathcal{C}%
_{A}^{\mathcal{D},\limfunc{restrict}}\right\} ,
\end{eqnarray*}%
consist of the first and second components respectively of the pairs in $%
\mathcal{P}$, and writing 
\begin{eqnarray}
\mathsf{B}_{\limfunc{stop}}^{A,\mathcal{P}}\left( f,g\right) &=&\sum_{J\in
\Pi _{2}\mathcal{P}}\left\langle T_{\sigma }^{\alpha }\varphi _{J}^{\mathcal{%
P}},\square _{J}^{\omega ,\mathbf{b}^{\ast }}g\right\rangle _{\omega };
\label{def phi P} \\
\text{where }\varphi _{J}^{\mathcal{P}} &\equiv &\sum_{I\in \mathcal{C}_{A}^{%
\limfunc{restrict}}:\ \left( I,J\right) \in \mathcal{P}}b_{A}\left(
E_{I}^{\sigma }\widehat{\square }_{\pi I}^{\sigma ,\flat ,\mathbf{b}%
}f\right) \ \mathbf{1}_{A\setminus I}\ ,  \notag
\end{eqnarray}%
where $E_{I}^{\sigma }h\equiv \frac{1}{\left\vert I\right\vert _{\sigma }}%
\int_{I}hd\sigma $ denotes the $\sigma $-average of $h$ on\thinspace $I$,
and where we note that the function $\widehat{\square }_{\pi I}^{\sigma
,\flat ,\mathbf{b}}f$ is constant on $I$, so that $E_{I}^{\sigma }\widehat{%
\square }_{\pi I}^{\sigma ,\flat ,\mathbf{b}}f$ simply picks out the value
of $\widehat{\square }_{\pi I}^{\sigma ,\flat ,\mathbf{b}}f$ on $I$. By the
tree-connected property of $\mathcal{P}$, and the telescoping property (\ref%
{telescoping}), together with the bound $\alpha _{\mathcal{A}}\left(
A\right) $ on the averages of $f$ in the corona $\mathcal{C}_{A}$, we have%
\begin{equation}
\left\vert \varphi _{J}^{\mathcal{P}}\right\vert \lesssim \alpha _{\mathcal{A%
}}\left( A\right) \mathbf{1}_{A\setminus I_{\mathcal{P}}\left( J\right) },
\label{phi bound}
\end{equation}%
where $I_{\mathcal{P}}\left( J\right) \equiv \dbigcap \left\{ I:\left(
I,J\right) \in \mathcal{P}\right\} $ is the smallest interval $I$ for which $%
\left( I,J\right) \in \mathcal{P}$. It is important to note that $J$ is good
with respect to $I_{\mathcal{P}}\left( J\right) $ and beyond by the infusion
of weak goodness above.

Another important property of these functions is the sublinearity:%
\begin{equation}
\left\vert \varphi _{J}^{\mathcal{P}}\right\vert \leq \left\vert \varphi
_{J}^{\mathcal{P}_{1}}\right\vert +\left\vert \varphi _{J}^{\mathcal{P}%
_{2}}\right\vert ,\ \ \ \ \ \mathcal{P}=\mathcal{P}_{1}\dot{\cup}\mathcal{P}%
_{2}\ ,  \label{phi sublinear}
\end{equation}%
which is an immediate consequence of 
\begin{equation*}
\varphi _{J}^{\mathcal{P}_{1}\dot{\cup}\mathcal{P}_{2}}=\sum_{I\in \mathcal{C%
}_{A}^{\limfunc{restrict}}:\ \left( I,J\right) \in \mathcal{P}_{1}\dot{\cup}%
\mathcal{P}_{2}}\left\{ ...\right\} =\sum_{I\in \mathcal{C}_{A}^{\limfunc{%
restrict}}:\ \left( I,J\right) \in \mathcal{P}_{1}}\left\{ ...\right\}
+\sum_{I\in \mathcal{C}_{A}^{\limfunc{restrict}}:\ \left( I,J\right) \in 
\mathcal{P}_{2}}\left\{ ...\right\} =\varphi _{J}^{\mathcal{P}_{1}}+\varphi
_{J}^{\mathcal{P}_{2}}.
\end{equation*}%
Now apply the Monotonicity Lemma \ref{mono} to the inner product $%
\left\langle T_{\sigma }^{\alpha }\varphi _{J},\square _{J}^{\omega ,\mathbf{%
b}^{\ast }}g\right\rangle _{\omega }$ (which applies since $J$ is good in $%
I_{\mathcal{P}}\left( J\right) $) to obtain%
\begin{eqnarray*}
\left\vert \left\langle T_{\sigma }^{\alpha }\varphi _{J},\square
_{J}^{\omega ,\mathbf{b}^{\ast }}g\right\rangle _{\omega }\right\vert
&\lesssim &\frac{\mathrm{P}^{\alpha }\left( J,\left\vert \varphi _{J}^{%
\mathcal{P}}\right\vert \mathbf{1}_{A\setminus I_{\mathcal{P}}\left(
J\right) }\sigma \right) }{\left\vert J\right\vert }\left\Vert
\bigtriangleup _{J}^{\omega ,\mathbf{b}^{\ast }}x\right\Vert _{L^{2}\left(
\omega \right) }^{\spadesuit }\left\Vert \square _{J}^{\omega ,\mathbf{b}%
^{\ast }}g\right\Vert _{L^{2}\left( \omega \right) }^{\bigstar } \\
&&+\frac{\mathrm{P}_{1+\delta }^{\alpha }\left( J,\left\vert \varphi _{J}^{%
\mathcal{P}}\right\vert \mathbf{1}_{A\setminus I_{\mathcal{P}}\left(
J\right) }\sigma \right) }{\left\vert J\right\vert }\left\Vert
x-m_{J}^{\omega }\right\Vert _{L^{2}\left( \mathbf{1}_{J}\omega \right)
}\left\Vert \square _{J}^{\omega ,\mathbf{b}^{\ast }}g\right\Vert
_{L^{2}\left( \omega \right) }^{\bigstar }.
\end{eqnarray*}%
Thus we have%
\begin{eqnarray}
\left\vert \mathsf{B}_{\limfunc{stop}}^{A,\mathcal{P}}\left( f,g\right)
\right\vert &\leq &\sum_{J\in \Pi _{2}\mathcal{P}}\frac{\mathrm{P}^{\alpha
}\left( J,\left\vert \varphi _{J}^{\mathcal{P}}\right\vert \mathbf{1}%
_{A\setminus I_{\mathcal{P}}\left( J\right) }\sigma \right) }{\left\vert
J\right\vert }\left\Vert \bigtriangleup _{J}^{\omega ,\mathbf{b}^{\ast
}}x\right\Vert _{L^{2}\left( \omega \right) }^{\spadesuit }\left\Vert
\square _{J}^{\omega ,\mathbf{b}^{\ast }}g\right\Vert _{L^{2}\left( \omega
\right) }^{\bigstar }  \label{def mod B} \\
&&+\sum_{J\in \Pi _{2}\mathcal{P}}\frac{\mathrm{P}_{1+\delta }^{\alpha
}\left( J,\left\vert \varphi _{J}^{\mathcal{P}}\right\vert \mathbf{1}%
_{A\setminus I_{\mathcal{P}}\left( J\right) }\sigma \right) }{\left\vert
J\right\vert }\left\Vert x-m_{J}^{\omega }\right\Vert _{L^{2}\left( \mathbf{1%
}_{J}\omega \right) }\left\Vert \square _{J}^{\omega ,\mathbf{b}^{\ast
}}g\right\Vert _{L^{2}\left( \omega \right) }^{\bigstar }  \notag \\
&\equiv &\left\vert \mathsf{B}\right\vert _{\limfunc{stop},1,\bigtriangleup
^{\omega }}^{A,\mathcal{P}}\left( f,g\right) +\left\vert \mathsf{B}%
\right\vert _{\limfunc{stop},1+\delta ,\mathsf{P}^{\omega }}^{A,\mathcal{P}%
}\left( f,g\right) ,  \notag
\end{eqnarray}%
where we have dominated the stopping form by two sublinear stopping forms
that involve the Poisson integrals of order $1$ and $1+\delta $
respectively, and where the smaller Poisson integral $\mathrm{P}_{1+\delta
}^{\alpha }$ is multiplied by the larger quantity $\left\Vert
x-m_{J}^{\omega }\right\Vert _{L^{2}\left( \mathbf{1}_{J}\omega \right) }$.
It remains to show the following two inequalities where we abbreviate $%
\left\vert \mathsf{B}\right\vert _{\limfunc{stop},1,\bigtriangleup ^{\omega
}}^{A,\mathcal{P}}$ to $\left\vert \mathsf{B}\right\vert _{\limfunc{stop}%
,\bigtriangleup ^{\omega }}^{A,\mathcal{P}}$ and $\left\vert \mathsf{B}%
\right\vert _{\limfunc{stop},1+\delta ,\mathsf{P}^{\omega }}^{A,\mathcal{P}}$
to $\left\vert \mathsf{B}\right\vert _{\limfunc{stop},1+\delta }^{A,\mathcal{%
P}}$:%
\begin{equation}
\left\vert \mathsf{B}\right\vert _{\limfunc{stop},\bigtriangleup ^{\omega
}}^{A,\mathcal{P}}\left( f,g\right) \lesssim \left( \mathcal{E}_{2}^{\alpha
}+\sqrt{\mathfrak{A}_{2}^{\alpha }}\right) \left\Vert \mathsf{P}_{\pi \left(
\Pi _{1}\mathcal{P}\right) }^{\sigma ,\mathbf{b}}f\right\Vert _{L^{2}\left(
\sigma \right) }^{\bigstar }\left\Vert \mathsf{P}_{\Pi _{2}\mathcal{P}%
}^{\omega ,\mathbf{b}^{\ast }}g\right\Vert _{L^{2}\left( \omega \right)
}^{\bigstar },  \label{First inequality}
\end{equation}%
for $f\in L^{2}\left( \sigma \right) $ satisfying $E_{I}^{\sigma }\left\vert
f\right\vert \lesssim \alpha _{\mathcal{A}}\left( A\right) $ for all $I\in 
\mathcal{C}_{A}$, and where $\pi \left( \Pi _{1}\mathcal{P}\right) \equiv
\left\{ \pi _{\mathcal{D}}I:I\in \Pi _{1}\mathcal{P}\right\} $; and%
\begin{equation}
\left\vert \mathsf{B}\right\vert _{\limfunc{stop},1+\delta }^{A,\mathcal{P}%
}\left( f,g\right) \lesssim \left( \mathcal{E}_{2}^{\alpha }+\sqrt{%
A_{2}^{\alpha }}\right) \left\Vert \mathsf{P}_{\mathcal{C}_{A}^{\mathcal{D}%
}}^{\sigma ,\mathbf{b}}f\right\Vert _{L^{2}\left( \sigma \right) }\left\Vert 
\mathsf{P}_{\mathcal{C}_{A}^{\mathcal{G},\func{shift}}}^{\omega ,\mathbf{b}%
^{\ast }}g\right\Vert _{L^{2}\left( \omega \right) },
\label{Second inequality}
\end{equation}%
where we only need the case $\mathcal{P}=\mathcal{P}^{A}$ in this latter
inequality as there is no absorption involved in treating this second
sublinear form. We consider first the easier inequality (\ref{Second
inequality}) that does not require absorption.

\subsection{The bound for the second sublinear inequality}

Here we turn to proving (\ref{Second inequality}), i.e.%
\begin{eqnarray*}
\left\vert \mathsf{B}\right\vert _{\limfunc{stop},1+\delta }^{A,\mathcal{P}%
}\left( f,g\right) &=&\sum_{J\in \Pi _{2}\mathcal{P}}\frac{\mathrm{P}%
_{1+\delta }^{\alpha }\left( J,\left\vert \varphi _{J}\right\vert \mathbf{1}%
_{A\setminus I_{\mathcal{P}}\left( J\right) }\sigma \right) }{\left\vert
J\right\vert }\left\Vert x-m_{J}\right\Vert _{L^{2}\left( \mathbf{1}%
_{J}\omega \right) }\left\Vert \square _{J}^{\omega ,\mathbf{b}^{\ast
}}g\right\Vert _{L^{2}\left( \omega \right) }^{\bigstar } \\
&\lesssim &\left( \mathcal{E}_{2}^{\alpha }+\sqrt{A_{2}^{\alpha }}\right)
\left\Vert \mathsf{P}_{\mathcal{C}_{A}^{\mathcal{D}}}^{\sigma ,\mathbf{b}%
}f\right\Vert _{L^{2}\left( \sigma \right) }\left\Vert \mathsf{P}_{\mathcal{C%
}_{A}^{\mathcal{G},\func{shift}}}^{\omega ,\mathbf{b}^{\ast }}g\right\Vert
_{L^{2}\left( \omega \right) },
\end{eqnarray*}%
where since 
\begin{equation*}
\left\vert \varphi _{J}\right\vert =\left\vert \sum_{I\in \mathcal{C}%
_{A}^{\prime }:\ \left( I,J\right) \in \mathcal{P}}\left( E_{I}^{\sigma
}\square _{\pi I}^{\sigma ,\flat ,\mathbf{b}}f\right) \ \mathbf{1}%
_{A\setminus I}\right\vert \leq \sum_{I\in \mathcal{C}_{A}^{\prime }:\
\left( I,J\right) \in \mathcal{P}}\left\vert E_{I}^{\sigma }\square _{\pi
I}^{\sigma ,\flat ,\mathbf{b}}\left( f\right) \ \mathbf{1}_{A\setminus
I}\right\vert ,
\end{equation*}%
the sublinear form $\left\vert \mathsf{B}\right\vert _{\limfunc{stop}%
,1+\delta }^{A,\mathcal{P}}$ can be decomposed by pigeonholing the ratio of
side lengths of $J$ and $I$:%
\begin{eqnarray*}
\left\vert \mathsf{B}\right\vert _{\limfunc{stop},1+\delta }^{A,\mathcal{P}%
}\left( f,g\right) &=&\sum_{J\in \Pi _{2}\mathcal{P}}\frac{\mathrm{P}%
_{1+\delta }^{\alpha }\left( J,\left\vert \varphi _{J}\right\vert \mathbf{1}%
_{A\setminus I_{\mathcal{P}}\left( J\right) }\sigma \right) }{\left\vert
J\right\vert ^{\frac{1}{n}}}\left\Vert x-m_{J}\right\Vert _{L^{2}\left( 
\mathbf{1}_{J}\omega \right) }\left\Vert \square _{J}^{\omega ,\mathbf{b}%
^{\ast }}g\right\Vert _{L^{2}\left( \omega \right) }^{\bigstar } \\
&\leq &\sum_{\left( I,J\right) \in \mathcal{P}}\frac{\mathrm{P}_{1+\delta
}^{\alpha }\left( J,\left\vert E_{I}^{\sigma }\left( \square _{\pi
I}^{\sigma ,\flat ,\mathbf{b}}f\right) \right\vert \mathbf{1}_{A\setminus
I}\sigma \right) }{\left\vert J\right\vert ^{\frac{1}{n}}}\left\Vert
x-m_{J}\right\Vert _{L^{2}\left( \mathbf{1}_{J}\omega \right) }\left\Vert
\square _{J}^{\omega ,\mathbf{b}^{\ast }}g\right\Vert _{L^{2}\left( \omega
\right) }^{\bigstar } \\
&\equiv &\sum_{s=0}^{\infty }\left\vert \mathsf{B}\right\vert _{\limfunc{stop%
},1+\delta }^{A,\mathcal{P};s}\left( f,g\right) ; \\
\left\vert \mathsf{B}\right\vert _{\limfunc{stop},1+\delta }^{A,\mathcal{P}%
;s}\left( f,g\right) &\equiv &\sum_{\substack{ \left( I,J\right) \in 
\mathcal{P}  \\ \ell \left( J\right) =2^{-s}\ell \left( I\right) }}\frac{%
\mathrm{P}_{1+\delta }^{\alpha }\left( J,\left\vert \left( E_{I}^{\sigma
}\square _{\pi I}^{\sigma ,\flat ,\mathbf{b}}f\right) \right\vert \mathbf{1}%
_{A\setminus I}\sigma \right) }{\left\vert J\right\vert ^{\frac{1}{n}}}%
\left\Vert x-m_{J}\right\Vert _{L^{2}\left( \mathbf{1}_{J}\omega \right)
}\left\Vert \square _{J}^{\omega ,\mathbf{b}^{\ast }}g\right\Vert
_{L^{2}\left( \omega \right) }^{\bigstar }.
\end{eqnarray*}%
We will now adapt the argument for the stopping term starting on page 42 of 
\cite{LaSaUr2}, where the geometric gain from the assumed Energy Hypothesis
there will be replaced by a geometric gain from the smaller Poisson integral 
$\mathrm{P}_{1+\delta }^{\alpha }$ used here.

We exploit the additional decay in the Poisson integral $\mathrm{P}%
_{1+\delta }^{\alpha }$ as follows. Suppose that $J$ is good in $I$ with $%
\ell \left( J\right) =2^{-s}\ell \left( I\right) $. We then obtain from (\ref%
{Poisson decay}) above that%
\begin{equation}
\left( \frac{\mathrm{P}_{1+\delta }^{\alpha }\left( J,\mathbf{1}_{A\setminus
I}\sigma \right) }{\left\vert J\right\vert ^{\frac{1}{n}}}\right) \lesssim
2^{-s\delta \left( 1-\varepsilon \right) }\frac{\mathrm{P}^{\alpha }\left( J,%
\mathbf{1}_{A\setminus I}\sigma \right) }{\left\vert J\right\vert ^{\frac{1}{%
n}}}.  \label{Poisson decay'}
\end{equation}%
We next claim that for $s\geq 0$ an integer,%
\begin{eqnarray*}
\left\vert \mathsf{B}\right\vert _{\limfunc{stop},1+\delta }^{A,\mathcal{P}%
;s}\left( f,g\right) &=&\sum_{\substack{ \left( I,J\right) \in \mathcal{P} 
\\ \ell \left( J\right) =2^{-s}\ell \left( I\right) }}\frac{\mathrm{P}%
_{1+\delta }^{\alpha }\left( J,\left\vert \left( E_{I}^{\sigma }\square
_{\pi I}^{\sigma ,\flat ,\mathbf{b}}f\right) \right\vert \mathbf{1}%
_{A\setminus I}\sigma \right) }{\left\vert J\right\vert }\left\Vert
x-m_{J}\right\Vert _{L^{2}\left( \mathbf{1}_{J}\omega \right) }\left\Vert
\square _{J}^{\omega ,\mathbf{b}^{\ast }}g\right\Vert _{L^{2}\left( \omega
\right) }^{\bigstar } \\
&\lesssim &2^{-s\delta \left( 1-\varepsilon \right) }\ \left( \mathcal{E}%
_{2}^{\alpha }+\sqrt{A_{2}^{\alpha }}\right) \ \left\Vert f\right\Vert
_{L^{2}\left( \sigma \right) }\ \left\Vert g\right\Vert _{L^{2}\left( \omega
\right) }\,,
\end{eqnarray*}%
from which (\ref{Second inequality}) follows upon summing in $s\geq 0$. Now
using both%
\begin{eqnarray*}
\left\vert E_{I}^{\sigma }\square _{\pi I}^{\sigma ,\flat ,\mathbf{b}%
}f\right\vert &=&\frac{1}{\left\vert I\right\vert _{\sigma }}%
\int_{I}\left\vert \square _{\pi I}^{\sigma ,\flat ,\mathbf{b}}f\right\vert
d\sigma \leq \left\Vert \square _{\pi I}^{\sigma ,\flat ,\mathbf{b}%
}f\right\Vert _{L^{2}\left( \sigma \right) }\frac{1}{\sqrt{\left\vert
I\right\vert _{\sigma }}}, \\
\sum_{I\in \mathcal{D}}\left\Vert \square _{\pi I}^{\sigma ,\flat ,\mathbf{b}%
}f\right\Vert _{L^{2}\left( \sigma \right) }^{2} &\lesssim &\sum_{I\in 
\mathcal{D}}\left( \left\Vert \square _{\pi I}^{\sigma ,\mathbf{b}%
}f\right\Vert _{L^{2}\left( \sigma \right) }^{2}+\left\Vert \nabla _{\pi
I}^{\sigma }f\right\Vert _{L^{2}\left( \sigma \right) }^{2}\right) \approx
\left\Vert f\right\Vert _{L^{2}(\sigma )}^{2}\ ,
\end{eqnarray*}%
where the second line uses frame inequalities in Proposition \ref{dual frame}
and displays (\ref{low frame}) and (\ref{corr upper}) from Appendix A below,
we apply Cauchy-Schwarz in the $I$ variable above to see that 
\begin{eqnarray*}
&&\left[ \left\vert \mathsf{B}\right\vert _{\limfunc{stop},1+\delta }^{A,%
\mathcal{P};s}\left( f,g\right) \right] ^{2}\lesssim \left\Vert f\right\Vert
_{L^{2}(\sigma )}^{2} \\
&&\times \left[ \sum_{I\in \mathcal{C}_{A}^{\limfunc{restrict}}}\left( \frac{%
1}{\sqrt{\left\vert I\right\vert _{\sigma }}}\sum_{\substack{ J:\ \left(
I,J\right) \in \mathcal{P}  \\ \ell \left( J\right) =2^{-s}\ell \left(
I\right) }}\frac{\mathrm{P}_{1+\delta }^{\alpha }\left( J,\mathbf{1}%
_{A\setminus I}\sigma \right) }{\left\vert J\right\vert }\left\Vert
x-m_{J}\right\Vert _{L^{2}\left( \mathbf{1}_{J}\omega \right) }\left\Vert
\square _{J}^{\omega ,\mathbf{b}^{\ast }}g\right\Vert _{L^{2}\left( \omega
\right) }^{\bigstar }\right) ^{2}\right] ^{\frac{1}{2}}.
\end{eqnarray*}%
We can then estimate the sum inside the square brackets by%
\begin{equation*}
\sum_{I\in \mathcal{C}_{A}^{\prime }}\left\{ \sum_{\substack{ J:\ \left(
I,J\right) \in \mathcal{P}  \\ \ell \left( J\right) =2^{-s}\ell \left(
I\right) }}\left\Vert \square _{J}^{\omega ,\mathbf{b}^{\ast }}g\right\Vert
_{L^{2}\left( \omega \right) }^{\bigstar 2}\right\} \sum_{\substack{ J:\
\left( I,J\right) \in \mathcal{P}  \\ \ell \left( J\right) =2^{-s}\ell
\left( I\right) }}\frac{1}{\left\vert I\right\vert _{\sigma }}\left( \frac{%
\mathrm{P}_{1+\delta }^{\alpha }\left( J,\mathbf{1}_{A\setminus I}\sigma
\right) }{\left\vert J\right\vert }\right) ^{2}\left\Vert x-m_{J}\right\Vert
_{L^{2}\left( \mathbf{1}_{J}\omega \right) }^{2}\lesssim \left\Vert
g\right\Vert _{L^{2}\left( \omega \right) }^{2}A\left( s\right) ^{2},
\end{equation*}%
where%
\begin{equation*}
A\left( s\right) ^{2}\equiv \sup_{I\in \mathcal{C}_{A}^{\prime }}\sum 
_{\substack{ J:\ \left( I,J\right) \in \mathcal{P}  \\ \ell \left( J\right)
=2^{-s}\ell \left( I\right) }}\frac{1}{\left\vert I\right\vert _{\sigma }}%
\left( \frac{\mathrm{P}_{1+\delta }^{\alpha }\left( J,\mathbf{1}_{A\setminus
I}\sigma \right) }{\left\vert J\right\vert }\right) ^{2}\left\Vert
x-m_{J}\right\Vert _{L^{2}\left( \mathbf{1}_{J}\omega \right) }^{2}\,.
\end{equation*}%
Finally then we turn to the analysis of the supremum in last display. From
the Poisson decay (\ref{Poisson decay'}) we have 
\begin{eqnarray*}
A\left( s\right) ^{2} &\lesssim &\sup_{I\in \mathcal{C}_{A}^{\prime }}\frac{1%
}{\left\vert I\right\vert _{\sigma }}2^{-2s\delta \left( 1-\varepsilon
\right) }\sum_{\substack{ J:\ \left( I,J\right) \in \mathcal{P}  \\ \ell
\left( J\right) =2^{-s}\ell \left( I\right) }}\left( \frac{\mathrm{P}%
^{\alpha }\left( J,\mathbf{1}_{A\setminus I}\sigma \right) }{\left\vert
J\right\vert }\right) ^{2}\left\Vert x-m_{J}\right\Vert _{L^{2}\left( 
\mathbf{1}_{J}\omega \right) }^{2} \\
&\lesssim &2^{-2s\delta \left( 1-\varepsilon \right) }\left[ \left( \mathcal{%
E}_{2}^{\alpha }\right) ^{2}+A_{2}^{\alpha }\right] \,,
\end{eqnarray*}%
where the last inequality is the stopping energy inquality (\ref{def
stopping energy 3}) in the corona $\mathcal{C}_{A}$. Indeed, from Definition~%
\ref{def energy corona 3}, as $(I,J)\in \mathcal{P}$, we have that $I$ is 
\emph{not} a stopping interval in $\mathcal{A}$, and hence that (\ref{def
stop 3}) \emph{fails} to hold, delivering the estimate above since the terms 
$\left\Vert x-m_{J}\right\Vert _{L^{2}\left( \mathbf{1}_{J}\omega \right)
}^{2}$ are additive, as the $J^{\prime }s$ are pigeonholed by $\ell \left(
J\right) =2^{-s}\ell \left( I\right) $ and hence pairwise disjoint.

\subsection{The bound for the first sublinear inequality}

Recall the definition of the sublinear form $\left\vert \mathsf{B}%
\right\vert _{\limfunc{stop},\bigtriangleup ^{\omega }}^{A,\mathcal{P}%
}\left( f,g\right) $ in display (\ref{def mod B}), and the inequality (\ref%
{First inequality}) that we wish to prove.

\begin{definition}
\label{Norm hat}Denote by $\widehat{\mathfrak{N}}_{\limfunc{stop}%
,\bigtriangleup ^{\omega }}^{A,\mathcal{P}}$ the best constant in 
\begin{equation}
\left\vert \mathsf{B}\right\vert _{\limfunc{stop},\bigtriangleup ^{\omega
}}^{A,\mathcal{P}}\left( f,g\right) \leq \widehat{\mathfrak{N}}_{\limfunc{%
stop},\bigtriangleup ^{\omega }}^{A,\mathcal{P}}\left\Vert \mathsf{P}_{\pi
\left( \Pi _{1}\mathcal{P}\right) }^{\sigma ,\mathbf{b}}f\right\Vert
_{L^{2}\left( \sigma \right) }^{\bigstar }\left\Vert \mathsf{P}_{\Pi _{2}%
\mathcal{P}}^{\omega ,\mathbf{b}^{\ast }}g\right\Vert _{L^{2}\left( \omega
\right) }^{\bigstar },  \label{best hat}
\end{equation}%
where $f\in L^{2}\left( \sigma \right) $ satisfies $E_{I}^{\sigma
}\left\vert f\right\vert \leq \alpha _{\mathcal{A}}\left( A\right) $ for all 
$I\in \mathcal{C}_{A}$, and $g\in L^{2}\left( \omega \right) $, and $\pi
\left( \Pi _{1}\mathcal{P}\right) =\left\{ \pi I:I\in \Pi _{1}\mathcal{P}%
\right\} $. We refer to $\widehat{\mathfrak{N}}_{\limfunc{stop}%
,\bigtriangleup ^{\omega }}^{A,\mathcal{P}}$ as the \emph{restricted} norm
relative to the collection $\mathcal{P}$.
\end{definition}

Inequality (\ref{First inequality}) will follow once we have shown that $%
\widehat{\mathfrak{N}}_{\limfunc{stop},\bigtriangleup ^{\omega }}^{A,%
\mathcal{P}}\lesssim \mathcal{E}_{2}^{\alpha }+\sqrt{\mathfrak{A}%
_{2}^{\alpha }}$. To this end, the following general result on mutually
orthogonal admissible collections, given as (3.5) in \cite{Lac}, will prove
very useful. Given a set $\left\{ \mathcal{Q}_{m}\right\} _{m=0}^{\infty }$
of admissible collections for $A$, we say that the collections $\mathcal{Q}%
_{m}$ are \emph{mutually orthogonal}, if each collection $\mathcal{Q}_{m}$
satisfies%
\begin{equation*}
\mathcal{Q}_{m}\subset \dbigcup\limits_{j=0}^{\infty }\left\{ \mathcal{A}%
_{m,j}\times \mathcal{B}_{m,j}\right\} \ ,
\end{equation*}%
where the sets $\left\{ \mathcal{A}_{m,j}\right\} _{m,j}$ and $\left\{ 
\mathcal{B}_{m,j}\right\} _{m,j}$ are each pairwise disjoint in their
respective dyadic grids $\mathcal{D}$ and $\mathcal{G}$: 
\begin{equation*}
\sum_{m,j=0}^{\infty }\mathbf{1}_{\mathcal{A}_{m,j}}\leq \mathbf{1}_{%
\mathcal{D}}\text{ and }\sum_{m,j=0}^{\infty }\mathbf{1}_{\mathcal{B}%
_{m,j}}\leq \mathbf{1}_{\mathcal{G}}.
\end{equation*}

\begin{lemma}
\label{mut orth}Suppose that $\left\{ \mathcal{Q}_{m}\right\} _{m=0}^{\infty
}$ is a set of admissible collections for $A$ that are \emph{mutually
orthogonal}. Then $\mathcal{Q}\equiv \dbigcup\limits_{m=0}^{\infty }\mathcal{%
Q}_{m}$ is admissible, and the sublinear stopping form $\left\vert \mathsf{B}%
\right\vert _{\limfunc{stop},\bigtriangleup ^{\omega }}^{A,\mathcal{Q}%
}\left( f,g\right) $ has its localized norm $\widehat{\mathfrak{N}}_{%
\limfunc{stop},\bigtriangleup ^{\omega }}^{A,\mathcal{Q}}$ controlled by the 
\emph{supremum} of the localized norms $\widehat{\mathfrak{N}}_{\limfunc{stop%
},\bigtriangleup ^{\omega }}^{A,\mathcal{Q}_{m}}$: 
\begin{equation*}
\widehat{\mathfrak{N}}_{\limfunc{stop},\bigtriangleup ^{\omega }}^{A,%
\mathcal{Q}}\leq \sup_{m\geq 0}\widehat{\mathfrak{N}}_{\limfunc{stop}%
,\bigtriangleup ^{\omega }}^{A,\mathcal{Q}_{m}}.
\end{equation*}
\end{lemma}

\begin{proof}
If $J\in \Pi _{2}\mathcal{Q}_{m}$, then $\varphi _{J}^{\mathcal{Q}}=\varphi
_{J}^{\mathcal{Q}_{m}}$ and $I_{\mathcal{Q}}\left( J\right) =I_{\mathcal{Q}%
_{m}}\left( J\right) $, since the collection $\left\{ \mathcal{Q}%
_{m}\right\} _{m=0}^{\infty }$ is mutually orthogonal. Thus we have 
\begin{eqnarray*}
\left\vert \mathsf{B}\right\vert _{\limfunc{stop},\bigtriangleup ^{\omega
}}^{A,\mathcal{Q}}\left( f,g\right) &=&\sum_{J\in \Pi _{2}\mathcal{Q}}\frac{%
\mathrm{P}^{\alpha }\left( J,\left\vert \varphi _{J}^{\mathcal{Q}%
}\right\vert \mathbf{1}_{A\setminus I_{\mathcal{Q}}\left( J\right) }\sigma
\right) }{\left\vert J\right\vert }\left\Vert \bigtriangleup _{J}^{\omega ,%
\mathbf{b}^{\ast }}x\right\Vert _{L^{2}\left( \omega \right) }^{\spadesuit
}\left\Vert \square _{J}^{\omega ,\mathbf{b}^{\ast }}g\right\Vert
_{L^{2}\left( \omega \right) }^{\bigstar } \\
&=&\sum_{m\geq 0}\sum_{J\in \Pi _{2}\mathcal{Q}_{m}}\frac{\mathrm{P}^{\alpha
}\left( J,\left\vert \varphi _{J}^{\mathcal{Q}_{m}}\right\vert \mathbf{1}%
_{A\setminus I_{\mathcal{Q}_{m}}\left( J\right) }\sigma \right) }{\left\vert
J\right\vert }\left\Vert \bigtriangleup _{J}^{\omega ,\mathbf{b}^{\ast
}}x\right\Vert _{L^{2}\left( \omega \right) }^{\spadesuit }\left\Vert
\square _{J}^{\omega ,\mathbf{b}^{\ast }}g\right\Vert _{L^{2}\left( \omega
\right) }^{\bigstar }=\sum_{m\geq 0}\left\vert \mathsf{B}\right\vert _{%
\limfunc{stop},\bigtriangleup ^{\omega }}^{A,\mathcal{Q}_{m}}\left(
f,g\right) ,
\end{eqnarray*}%
and we can continue with the definition of $\widehat{\mathfrak{N}}_{\limfunc{%
stop},\bigtriangleup ^{\omega }}^{A,\mathcal{Q}_{m}}$ and Cauchy-Schwarz to
obtain%
\begin{eqnarray*}
\left\vert \mathsf{B}\right\vert _{\limfunc{stop},\bigtriangleup ^{\omega
}}^{A,\mathcal{Q}}\left( f,g\right) &\leq &\sum_{m\geq 0}\widehat{\mathfrak{N%
}}_{\limfunc{stop},\bigtriangleup ^{\omega }}^{A,\mathcal{Q}_{m}}\left\Vert 
\mathsf{P}_{\pi \left( \Pi _{1}\mathcal{Q}_{m}\right) }^{\sigma ,\mathbf{b}%
}f\right\Vert _{L^{2}\left( \sigma \right) }^{\bigstar }\left\Vert \mathsf{P}%
_{\Pi _{2}\mathcal{Q}_{m}}^{\omega ,\mathbf{b}^{\ast }}g\right\Vert
_{L^{2}\left( \omega \right) }^{\bigstar } \\
&\leq &\left( \sup_{m\geq 0}\widehat{\mathfrak{N}}_{\limfunc{stop}%
,\bigtriangleup ^{\omega }}^{A,\mathcal{Q}_{m}}\right) \sqrt{\sum_{m\geq
0}\left\Vert \mathsf{P}_{\pi \left( \Pi _{1}\mathcal{Q}_{m}\right) }^{\sigma
,\mathbf{b}}f\right\Vert _{L^{2}\left( \sigma \right) }^{\bigstar 2}}\sqrt{%
\sum_{m\geq 0}\left\Vert \mathsf{P}_{\Pi _{2}\mathcal{Q}_{m}}^{\omega ,%
\mathbf{b}^{\ast }}g\right\Vert _{L^{2}\left( \omega \right) }^{\bigstar 2}}
\\
&\leq &\left( \sup_{m\geq 0}\widehat{\mathfrak{N}}_{\limfunc{stop}%
,\bigtriangleup ^{\omega }}^{A,\mathcal{Q}_{m}}\right) \sqrt{\left\Vert 
\mathsf{P}_{\pi \left( \Pi _{1}\mathcal{Q}\right) }^{\sigma ,\mathbf{b}%
}f\right\Vert _{L^{2}\left( \sigma \right) }^{\bigstar 2}}\sqrt{\left\Vert 
\mathsf{P}_{\Pi _{2}\mathcal{Q}}^{\omega ,\mathbf{b}^{\ast }}g\right\Vert
_{L^{2}\left( \omega \right) }^{\bigstar 2}}.
\end{eqnarray*}
\end{proof}

Now we turn to proving inequality (\ref{First inequality}) for the sublinear
form $\left\vert \mathsf{B}\right\vert _{\limfunc{stop},\bigtriangleup
^{\omega }}^{A,\mathcal{P}}\left( f,g\right) $, i.e.%
\begin{eqnarray*}
\left\vert \mathsf{B}\right\vert _{\limfunc{stop},\bigtriangleup ^{\omega
}}^{A,\mathcal{P}}\left( f,g\right) &\equiv &\sum_{J\in \Pi _{2}\mathcal{P}}%
\frac{\mathrm{P}^{\alpha }\left( J,\left\vert \varphi _{J}\right\vert 
\mathbf{1}_{A\setminus I_{\mathcal{P}}\left( J\right) }\sigma \right) }{%
\left\vert J\right\vert }\left\Vert \bigtriangleup _{J}^{\omega ,\mathbf{b}%
^{\ast }}x\right\Vert _{L^{2}\left( \omega \right) }^{\spadesuit }\left\Vert
\square _{J}^{\omega ,\mathbf{b}^{\ast }}g\right\Vert _{L^{2}\left( \omega
\right) }^{\bigstar } \\
&\lesssim &\left( \mathcal{E}_{2}^{\alpha }+\sqrt{\mathfrak{A}_{2}^{\alpha }}%
\right) \left\Vert \mathsf{P}_{\pi \left( \Pi _{1}\mathcal{P}\right)
}^{\sigma ,\mathbf{b}}f\right\Vert _{L^{2}\left( \sigma \right) }^{\bigstar
}\left\Vert \mathsf{P}_{\Pi _{2}\mathcal{P}}^{\omega ,\mathbf{b}^{\ast
}}g\right\Vert _{L^{2}\left( \omega \right) }^{\bigstar }; \\
\text{where }\varphi _{J} &\equiv &\sum_{I\in \mathcal{C}_{A}^{\limfunc{%
restrict}}:\ \left( I,J\right) \in \mathcal{P}}\left( E_{I}^{\sigma }%
\widehat{\square }_{\pi I}^{\sigma ,\flat ,\mathbf{b}}f\right) b_{A}\ 
\mathbf{1}_{A\setminus I}\ \text{\ is supported in }A\setminus I_{\mathcal{P}%
}\left( J\right) ,
\end{eqnarray*}%
and $I_{\mathcal{P}}\left( J\right) $ denotes the smallest interval $I\in 
\mathcal{D}$ for which $\left( I,J\right) \in \mathcal{P}$. We recall the
stopping energy from (\ref{def stopping energy 3}),%
\begin{equation*}
\mathbf{X}_{\alpha }\left( \mathcal{C}_{A}\right) ^{2}\equiv \sup_{I\in 
\mathcal{C}_{A}}\frac{1}{\left\vert I\right\vert _{\sigma }}\sup_{I\supset 
\overset{\cdot }{\cup }J_{r}}\sum_{r=1}^{\infty }\left( \frac{\mathrm{P}%
^{\alpha }\left( J_{r},\mathbf{1}_{A}\sigma \right) }{\left\vert
J_{r}\right\vert }\right) ^{2}\left\Vert x-m_{J_{r}}\right\Vert
_{L^{2}\left( \mathbf{1}_{J_{r}}\omega \right) }^{2}\ ,
\end{equation*}%
where the intervals $J_{r}\in \mathcal{G}$ are pairwise disjoint in $I$.

What now follows is an adaptation to our sublinear form $\left\vert \mathsf{B%
}\right\vert _{\limfunc{stop},\square ^{\omega }}^{A,\mathcal{P}}$ of the
arguments of M. Lacey in \cite{Lac}, together with an additional `indented'
corona construction. We have the following Poisson inequality for intervals $%
B\subset A\subset I$:%
\begin{eqnarray}
\frac{\mathrm{P}^{\alpha }\left( A,\mathbf{1}_{I\setminus A}\sigma \right) }{%
\left\vert A\right\vert } &\approx &\int_{I\setminus A}\frac{1}{\left(
\left\vert y-c_{A}\right\vert \right) ^{2-\alpha }}d\sigma \left( y\right)
\label{BAI} \\
&\lesssim &\int_{I\setminus A}\frac{1}{\left( \left\vert y-c_{B}\right\vert
\right) ^{2-\alpha }}d\sigma \left( y\right) \approx \frac{\mathrm{P}%
^{\alpha }\left( B,\mathbf{1}_{I\setminus A}\sigma \right) }{\left\vert
B\right\vert }.  \notag
\end{eqnarray}

Fix $A\in \mathcal{A}$. Following \cite{Lac} we will use a `decoupled'
modification of the stopping energy $\mathbf{X}_{\alpha }\left( \mathcal{C}%
_{A}\right) $ to define a `size functional' of an $A$-admissible collection $%
\mathcal{P}$. So suppose that $\mathcal{P}$ is an $A$-admissible collection
of pairs of intervals, and recall that $\Pi _{1}\mathcal{P}$ and $\Pi _{2}%
\mathcal{P}$ denote the intervals in the first and second components of the
pairs in $\mathcal{P}$ respectively.

\begin{definition}
\label{rest K}For an $A$-admissible collection of pairs of intervals $%
\mathcal{P}$, and an interval $K\in \Pi _{1}\mathcal{P}$, define the
projection of $\mathcal{P}$ `relative to $K$' by 
\begin{equation*}
\Pi _{2}^{K}\mathcal{P}\equiv \left\{ J\in \Pi _{2}\mathcal{P}:\ J^{\maltese
}\subset K\right\} ,
\end{equation*}%
where we have suppressed dependence on $A$.
\end{definition}

\begin{definition}
\label{Pi below}We will use as the `size testing collection' of intervals
for $\mathcal{P}$ the collection 
\begin{equation*}
\Pi _{1}^{\limfunc{below}}\mathcal{P}\equiv \left\{ K\in \mathcal{D}%
:K\subset I\text{ for some }I\in \Pi _{1}\mathcal{P}\right\} ,
\end{equation*}%
which consists of all intervals contained in an interval from $\Pi _{1}%
\mathcal{P}$.
\end{definition}

Continuing to follow Lacey \cite{Lac}, we define a `size functional' of $%
\mathcal{P}$, actually two of them, as follows. Recall from Notation \ref%
{nonstandard norm} that for a pseudoprojection $\mathsf{Q}_{\mathcal{H}%
}^{\omega }$ on $x$ we have 
\begin{equation}
\left\Vert \mathsf{Q}_{\mathcal{H}}^{\omega ,\mathbf{b}^{\ast }}x\right\Vert
_{L^{2}\left( \omega \right) }^{\spadesuit 2}=\sum_{J\in \mathcal{H}%
}\left\Vert \bigtriangleup _{J}^{\omega ,\mathbf{b}^{\ast }}x\right\Vert
_{L^{2}\left( \omega \right) }^{\spadesuit 2}=\sum_{J\in \mathcal{H}}\left(
\left\Vert \bigtriangleup _{J}^{\omega ,\mathbf{b}^{\ast }}x\right\Vert
_{L^{2}\left( \omega \right) }^{2}+\inf_{z\in \mathbb{R}}\sum_{J^{\prime
}\in \mathfrak{C}_{\limfunc{broken}}\left( J\right) }\left\vert J^{\prime
}\right\vert _{\omega }\left( E_{J^{\prime }}^{\omega }\left\vert
x-z\right\vert \right) ^{2}\right) .  \label{h}
\end{equation}

\begin{definition}
\label{def ext size}If $\mathcal{P}$ is $A$-admissible, define an \emph{%
initial} size condition $\mathcal{S}_{\limfunc{init}\limfunc{size}}^{\alpha
,A}\left( \mathcal{P}\right) $ by%
\begin{equation}
\mathcal{S}_{\limfunc{init}\limfunc{size}}^{\alpha ,A}\left( \mathcal{P}%
\right) ^{2}\equiv \sup_{K\in \Pi _{1}^{\limfunc{below}}\mathcal{P}}\frac{1}{%
\left\vert K\right\vert _{\sigma }}\left( \frac{\mathrm{P}^{\alpha }\left( K,%
\mathbf{1}_{A\setminus K}\sigma \right) }{\left\vert K\right\vert }\right)
^{2}\left\Vert \mathsf{Q}_{\Pi _{2}^{K}\mathcal{P}}^{\omega ,\mathbf{b}%
^{\ast }}x\right\Vert _{L^{2}\left( \omega \right) }^{\spadesuit 2}.
\label{def P stop energy 3}
\end{equation}
\end{definition}

The following key fact is essential.

\textbf{Key Fact \#1:}%
\begin{equation}
K\subset A\text{ and }K\notin \mathcal{C}_{A}\Longrightarrow \Pi _{2}^{K}%
\mathcal{P}=\emptyset \ .  \label{later use}
\end{equation}%
To see this, suppose that $K\subset A$ and $K\notin \mathcal{C}_{A}$. Then $%
K\subset A^{\prime }$ for some $A^{\prime }\in \mathfrak{C}_{\mathcal{A}%
}\left( A\right) $, and so if there is $J^{\prime }\in \Pi _{2}^{K}\mathcal{P%
}$, then $\left( J^{\prime }\right) ^{\maltese }\subset K\subset A^{\prime }$%
, which implies that $\left( J^{\prime }\right) ^{\maltese }\notin \mathcal{C%
}_{A}^{\mathcal{G},\limfunc{shift}}$, which contradicts $\Pi _{2}^{K}%
\mathcal{P}\subset \mathcal{C}_{A}^{\mathcal{G},\limfunc{shift}}$. We now
observe from (\ref{later use}) that we may also write the initial size
functional as%
\begin{equation}
\mathcal{S}_{\limfunc{init}\limfunc{size}}^{\alpha ,A}\left( \mathcal{P}%
\right) ^{2}\equiv \sup_{K\in \Pi _{1}^{\limfunc{below}}\mathcal{P}\cap
C_{A}^{\limfunc{restrict}}}\frac{1}{\left\vert K\right\vert _{\sigma }}%
\left( \frac{\mathrm{P}^{\alpha }\left( K,\mathbf{1}_{A\setminus K}\sigma
\right) }{\left\vert K\right\vert }\right) ^{2}\left\Vert \mathsf{Q}_{\Pi
_{2}^{K}\mathcal{P}}^{\omega ,\mathbf{b}^{\ast }}x\right\Vert _{L^{2}\left(
\omega \right) }^{\spadesuit 2}.  \label{rewrite size}
\end{equation}

However, we will also need to control certain pairs $\left( I,J\right) \in 
\mathcal{P}$ using testing intervals $K$ which are strictly smaller than $%
J^{\maltese }$, namely those $K\in \mathcal{C}_{A}$ such that $K\subset
J^{\maltese }\subset \pi _{\mathcal{D}}^{\left( 2\right) }K$. For this, we
need a second key fact regarding the intervals $J^{\maltese }$, that will
also play a crucial role in controlling pairs in the indented corona below,
and which is that $J$ is always contained in one of the \emph{inner} two
grandchildren of $J^{\maltese }$.

\textbf{Key Fact \#2:}%
\begin{eqnarray}
&&\text{\textbf{either} }3J\subset J_{-/+}^{\maltese }\text{ \textbf{or} }%
3J\subset J_{+/-}^{\maltese }\ ,  \label{indentation} \\
&&\text{where }\mathfrak{C}_{J^{\maltese }}^{\left( 1\right) }=\left\{
J_{-}^{\maltese },J_{+}^{\maltese }\right\} \text{ and }\mathfrak{C}%
_{J^{\maltese }}^{\left( 2\right) }=\left\{ J_{-/-}^{\maltese
},J_{-/+}^{\maltese },J_{+/-}^{\maltese },J_{+/+}^{\maltese }\right\} , 
\notag \\
&&\text{and the children and grandchildren are listed left to right.}  \notag
\end{eqnarray}%
To see this, suppose without loss of generality that the child $%
J_{J}^{\maltese }$ of $J^{\maltese }$ that contains $J$ is the left child $%
J_{-}^{\maltese }$ (which exists because $J$ is good in $J^{\maltese }$).
Then observe that $J$ is by definition $\varepsilon -\limfunc{bad}$ in $%
J_{-}^{\maltese }=J_{J}^{\maltese }$, i.e. $\limfunc{dist}\left( J,\limfunc{%
body}J_{-}^{\maltese }\right) \leq 2\left\vert J\right\vert ^{\varepsilon
}\left\vert J_{-}^{\maltese }\right\vert ^{1-\varepsilon }$, and so cannot
lie in the leftmost grandchild $J_{-/-}^{\maltese }$. Indeed, if $J\subset
J_{-/-}^{\maltese }$, then%
\begin{equation*}
\limfunc{dist}\left( J,\limfunc{body}J^{\maltese }\right) =\limfunc{dist}%
\left( J,\limfunc{body}J_{-}^{\maltese }\right) \leq 2\left\vert
J\right\vert ^{\varepsilon }\left\vert J_{-}^{\maltese }\right\vert
^{1-\varepsilon }=2^{\varepsilon }\left\vert J\right\vert ^{\varepsilon
}\left\vert J^{\maltese }\right\vert ^{1-\varepsilon }<2\left\vert
J\right\vert ^{\varepsilon }\left\vert J^{\maltese }\right\vert
^{1-\varepsilon },
\end{equation*}%
contradicting the fact that $J$ is $\varepsilon -\limfunc{good}$ in $%
J^{\maltese }$. Thus we must have $J\subset J_{-/+}^{\maltese }$ (where the
body of $J^{\maltese }$ does not intersect the interior of $%
J_{-/+}^{\maltese }$, thus permitting $J$ to be $\varepsilon -\limfunc{good}$
in $J^{\maltese }$). Finally, the fact that $J$ is $\varepsilon -\limfunc{%
good}$ in $J^{\maltese }$ implies that $3J\subset J_{-/+}^{\maltese }$. See
Figure \ref{gra}.

\FRAME{ftbpFU}{6.6293in}{2.352in}{0pt}{\Qcb{The interval $J$ lies in one of
the two inner grandchildren of $J^{\maltese }$, namely $J_{-/+}^{\maltese }$
or $J_{+/-}^{\maltese }$. The $\limfunc{body}$ of $J^{\maltese }$ consists
of the infinitely many red dots.}}{\Qlb{gra}}{grandchild.wmf}{\special%
{language "Scientific Word";type "GRAPHIC";maintain-aspect-ratio
TRUE;display "USEDEF";valid_file "F";width 6.6293in;height 2.352in;depth
0pt;original-width 7.1805in;original-height 10.469in;cropleft "0";croptop
"0.6394";cropright "0.5800";cropbottom "0.5030";filename
'Grandchild.wmf';file-properties "XNPEU";}}

This second key fact is what underlies the construction of the indented
corona below, and motivates the next definition of augmented projection, in
which we allow intervals $K$ satisfying $J\subset K\subsetneqq J^{\maltese
}\subset \pi _{\mathcal{D}}^{\left( 2\right) }K$, as well as $K\in C_{A}$,
to be tested over in the augmented size condition below.

\begin{definition}
\label{augs}Suppose $\mathcal{P}$ is an $A$-admissible collection.

\begin{enumerate}
\item For $K\in \Pi _{1}\mathcal{P}$, define the \emph{augmented} projection
of $\mathcal{P}$ relative to $K$ by%
\begin{equation*}
\Pi _{2}^{K,\limfunc{aug}}\mathcal{P}\equiv \left\{ J\in \Pi _{2}\mathcal{P}%
:J\subset K\text{ and }J^{\maltese }\subset \pi _{\mathcal{D}}^{\left(
2\right) }K\right\} =\left\{ J\in \Pi _{2}\mathcal{P}:J^{\flat }\subset
K\right\} .
\end{equation*}

\item Define the corresponding \emph{augmented} size functional $\mathcal{S}%
_{\limfunc{aug}\limfunc{size}}^{\alpha ,A}\left( \mathcal{P}\right) $ by 
\begin{equation*}
\mathcal{S}_{\limfunc{aug}\limfunc{size}}^{\alpha ,A}\left( \mathcal{P}%
\right) ^{2}\equiv \sup_{K\in \Pi _{1}^{\limfunc{below}}\mathcal{P}\cap
C_{A}^{\limfunc{restrict}}}\frac{1}{\left\vert K\right\vert _{\sigma }}%
\left( \frac{\mathrm{P}^{\alpha }\left( K,\mathbf{1}_{A\setminus K}\sigma
\right) }{\left\vert K\right\vert }\right) ^{2}\left\Vert \mathsf{Q}_{\Pi
_{2}^{K,\limfunc{aug}}\mathcal{P}}^{\omega ,\mathbf{b}^{\ast }}x\right\Vert
_{L^{2}\left( \omega \right) }^{\spadesuit 2}\ .
\end{equation*}
\end{enumerate}
\end{definition}

We note that the augmented projection $\Pi _{2}^{K,\limfunc{aug}}\mathcal{P}$
includes intervals $J$ for which $J\subset K\subsetneqq J^{\maltese }\subset
\pi _{\mathcal{D}}^{\left( 2\right) }K$, and hence $J$ need not be $%
\varepsilon -\limfunc{good}$ inside $K$. For $M\in \mathcal{D}$, denote by $%
M_{J}$ and $M_{\searrow J}$ the child and grandchild respectively of $M$
that contains $J$, provided they exist. Then by the second key fact (\ref%
{indentation}), and using that the endpoints of both $J_{-/+}^{\maltese }$
and $J_{+/-}^{\maltese }$ lie in the $\limfunc{body}$ of $J^{\maltese }$, we
have two consequences,%
\begin{equation*}
K\in \left\{ J_{J}^{\maltese },J_{\searrow J}^{\maltese }\right\} \text{ and 
}3J\subset J_{\searrow J}^{\maltese }\subset 3J_{\searrow J}^{\maltese
}\subset J^{\maltese }\text{ for }J\in \Pi _{2}^{K,\limfunc{nar}}\mathcal{P},
\end{equation*}%
which will play an important role below.

The augmented size functional $\mathcal{S}_{\limfunc{aug}\limfunc{size}%
}^{\alpha ,A}\left( \mathcal{P}\right) $ is a `decoupled' form of the
stopping energy $\mathbf{X}_{\alpha }\left( \mathcal{C}_{A}\right) $
restricted to $\mathcal{P}$, in which the intervals $J$ appearing in $%
\mathbf{X}_{\alpha }\left( \mathcal{C}_{A}\right) $ no longer appear in the
Poisson integral in $\mathcal{S}_{\limfunc{aug}\limfunc{size}}^{\alpha
,A}\left( \mathcal{P}\right) $, and it plays a crucial role in Lacey's
argument in \cite{Lac}. We note two essential properties of this definition
of size functional:

\begin{enumerate}
\item \textbf{Monotonicity} of size: $\mathcal{S}_{\limfunc{aug}\limfunc{size%
}}^{\alpha ,A}\left( \mathcal{P}\right) \leq \mathcal{S}_{\limfunc{aug}%
\limfunc{size}}^{\alpha ,A}\left( \mathcal{Q}\right) $ if $\mathcal{P}%
\subset \mathcal{Q}$,

\item \textbf{Control} by energy and Muckenhoupt conditions: $\mathcal{S}_{%
\limfunc{aug}\limfunc{size}}^{\alpha ,A}\left( \mathcal{P}\right) \lesssim 
\mathcal{E}_{2}^{\alpha }+\sqrt{\mathfrak{A}_{2}^{\alpha }}$.
\end{enumerate}

The monotonicity property follows from $\Pi _{1}^{\limfunc{below}}\mathcal{P}%
\subset \Pi _{1}^{\limfunc{below}}\mathcal{Q}$ and $\Pi _{2}^{K}\mathcal{P}%
\subset \Pi _{2}^{K}\mathcal{Q}$. The control property is contained in the
next lemma, which uses the stopping energy control for the form $\mathsf{B}_{%
\limfunc{stop}}^{A}\left( f,g\right) $ associated with $A$.

\begin{lemma}
\label{energy control}If $\mathcal{P}^{A}$ is as in (\ref{initial P}) and $%
\mathcal{P}\subset \mathcal{P}^{A}$, then 
\begin{equation*}
\mathcal{S}_{\limfunc{aug}\limfunc{size}}^{\alpha ,A}\left( \mathcal{P}%
\right) \leq \mathbf{X}_{\alpha }\left( \mathcal{C}_{A}\right) \lesssim 
\mathcal{E}_{2}^{\alpha }+\sqrt{A_{2}^{\alpha }}+\sqrt{A_{2}^{\alpha ,%
\limfunc{punct}}}\ .
\end{equation*}
\end{lemma}

\begin{proof}
We have%
\begin{eqnarray*}
\mathcal{S}_{\limfunc{aug}\limfunc{size}}^{\alpha ,A}\left( \mathcal{P}%
\right) ^{2} &=&\sup_{K\in \Pi _{1}^{\limfunc{below}}\mathcal{P}\cap 
\mathcal{C}_{A}^{\limfunc{restrict}}}\frac{1}{\left\vert K\right\vert
_{\sigma }}\left( \frac{\mathrm{P}^{\alpha }\left( K,\mathbf{1}_{A\setminus
K}\sigma \right) }{\left\vert K\right\vert }\right) ^{2}\left\Vert \mathsf{Q}%
_{\Pi _{2}^{K}\mathcal{P}\cup \Pi _{2}^{K,\limfunc{aug}}\mathcal{P}}^{\omega
,\mathbf{b}^{\ast }}x\right\Vert _{L^{2}\left( \omega \right) }^{\spadesuit
2} \\
&\lesssim &\sup_{K\in \mathcal{C}_{A}^{\limfunc{restrict}}}\frac{1}{%
\left\vert K\right\vert _{\sigma }}\left( \frac{\mathrm{P}^{\alpha }\left( K,%
\mathbf{1}_{A}\sigma \right) }{\left\vert K\right\vert }\right)
^{2}\left\Vert x-m_{K}\right\Vert _{L^{2}\left( \mathbf{1}_{K}\omega \right)
}^{2}\leq \mathbf{X}_{\alpha }\left( \mathcal{C}_{A}\right) ^{2},
\end{eqnarray*}%
which is the first inequality in the statement of the lemma. The second
inequality follows from (\ref{def stopping bounds 3}).
\end{proof}

There is an important special circumstance, introduced by M. Lacey in \cite%
{Lac}, in which we can bound our forms by the size functional, namely when
the pairs all straddle a subpartition of $A$, and we present this in the
next subsection. In order to handle the fact that the intervals in $\Pi
_{1}^{\limfunc{below}}\mathcal{P}\cap C_{A}^{\limfunc{restrict}}$ need no
longer enjoy any goodness, we will need to formulate a Substraddling Lemma
to deal with this situation as well. See \textbf{Remark on lack of usual
goodness} after (\ref{N_L}), where it is explained how this applies to the
proof of (\ref{rem}). Then in the following subsection, we use the bottom/up
stopping time construction of M. Lacey, together with an additional
`indented' top/down corona construction, to reduce control of the sublinear
stopping form $\left\vert \mathsf{B}\right\vert _{\limfunc{stop}%
,\bigtriangleup ^{\omega }}^{A,\mathcal{P}}\left( f,g\right) $ in inequality
(\ref{First inequality}) to the three special cases addressed by the
Orthogonality Lemma, the Straddling Lemma and the Substraddling Lemma.

\subsection{Straddling, Substraddling and Corona-straddling Lemmas}

We begin with the Corona-straddling Lemma in which the straddling collection
is the set of $\mathcal{A}$-children of $A$, and applies to the `corona
straddling' subcollection of the initial admissible collection $\mathcal{P}%
^{A}$ - see (\ref{initial P}). Define the `corona straddling' collection $%
\mathcal{P}_{\func{cor}}^{A}$ by%
\begin{equation}
\mathcal{P}_{\func{cor}}^{A}\equiv \dbigcup\limits_{A^{\prime }\in \mathfrak{%
C}_{\mathcal{A}}\left( A\right) }\left\{ \left( I,J\right) \in \mathcal{P}%
^{A}:J\subset A^{\prime }\varsubsetneqq J^{\maltese }\subset \pi _{\mathcal{D%
}}^{\left( 2\right) }A^{\prime }\right\} .  \label{def cor}
\end{equation}%
Note that $\mathcal{P}_{\func{cor}}^{A}$ is an $A$-admissible collection
that consists of just those pairs $\left( I,J\right) $ for which $%
J^{\maltese }$ is either the $\mathcal{D}$-parent or the $\mathcal{D}$%
-grandparent of a stopping interval $A^{\prime }\in \mathfrak{C}_{\mathcal{A}%
}\left( A\right) $. The bound for the norm of the corresponding form is
controlled by the energy condition.

\begin{lemma}
\label{cor strad 1}We have the sublinear form bound%
\begin{equation*}
\mathfrak{N}_{\limfunc{stop},\bigtriangleup ^{\omega }}^{A,\mathcal{P}_{%
\func{cor}}^{A}}\leq C\mathcal{E}_{2}^{\alpha }.
\end{equation*}
\end{lemma}

\begin{proof}
The key point here is our assumption that $J\subset A^{\prime
}\varsubsetneqq J^{\maltese }\subset \pi _{\mathcal{D}}^{\left( 2\right)
}A^{\prime }$ for $\left( I,J\right) \in \mathcal{P}_{\func{cor}}^{A}$,
which implies that in fact $3J\subset A^{\prime }$ since $J\cap \limfunc{body%
}\left( \pi _{\mathcal{D}}^{\left( 2\right) }A^{\prime }\right) =\emptyset $
because $J$ is $\varepsilon -\limfunc{good}$ in $\pi _{\mathcal{D}}^{\left(
2\right) }A^{\prime }$. We start with%
\begin{eqnarray*}
\left\vert \mathsf{B}\right\vert _{\limfunc{stop},\bigtriangleup ^{\omega
}}^{A,\mathcal{P}_{\func{cor}}^{A}}\left( f,g\right) &=&\sum_{J\in \Pi _{2}%
\mathcal{P}_{\func{cor}}^{A}}\frac{\mathrm{P}^{\alpha }\left( J,\left\vert
\varphi _{J}^{\mathcal{P}_{\func{cor}}^{A}}\right\vert \mathbf{1}%
_{A\setminus I_{\mathcal{P}_{\func{cor}}^{A}}\left( J\right) }\sigma \right) 
}{\left\vert J\right\vert }\left\Vert \bigtriangleup _{J}^{\omega ,\mathbf{b}%
^{\ast }}x\right\Vert _{L^{2}\left( \omega \right) }^{\spadesuit }\left\Vert
\square _{J}^{\omega ,\mathbf{b}^{\ast }}g\right\Vert _{L^{2}\left( \omega
\right) }^{\bigstar } \\
&=&\sum_{A^{\prime }\in \mathfrak{C}_{\mathcal{A}}\left( A\right)
}\sum_{J\in \Pi _{2}\mathcal{P}_{\func{cor}}^{A}:\ 3J\subset A^{\prime }}%
\frac{\mathrm{P}^{\alpha }\left( J,\left\vert \varphi _{J}^{\mathcal{P}_{%
\func{cor}}^{A}}\right\vert \mathbf{1}_{A\setminus I_{\mathcal{P}_{\func{cor}%
}^{A}}\left( J\right) }\sigma \right) }{\left\vert J\right\vert }\left\Vert
\bigtriangleup _{J}^{\omega ,\mathbf{b}^{\ast }}x\right\Vert _{L^{2}\left(
\omega \right) }^{\spadesuit }\left\Vert \square _{J}^{\omega ,\mathbf{b}%
^{\ast }}g\right\Vert _{L^{2}\left( \omega \right) }^{\bigstar }; \\
\text{where }\varphi _{J}^{\mathcal{P}_{\func{cor}}^{A}} &\equiv &\sum_{I\in
\Pi _{1}\mathcal{P}_{\func{cor}}^{A}:\mathcal{\ }\left( I,J\right) \in 
\mathcal{P}_{\func{cor}}^{A}}b_{A}E_{I}^{\sigma }\left( \widehat{\square }%
_{\pi I}^{\sigma ,\flat ,\mathbf{b}}f\right) \ \mathbf{1}_{A\setminus I}\ .
\end{eqnarray*}%
If $J\in \Pi _{2}\mathcal{P}_{\func{cor}}^{A}$ and $J\subset A^{\prime }\in 
\mathfrak{C}_{\mathcal{A}}\left( A\right) $, then either (1) $A^{\prime
}=J_{-/+}^{\maltese }$ or $J_{+/-}^{\maltese }$ or (2) $A^{\prime
}=J_{-}^{\maltese }$ or $J_{+}^{\maltese }$, and we have 
\begin{equation*}
\frac{\mathrm{P}^{\alpha }\left( J,\mathbf{1}_{A\setminus I_{\mathcal{P}_{%
\func{cor}}^{A}}\left( J\right) }\sigma \right) }{\left\vert J\right\vert }%
\approx \left\{ 
\begin{array}{ccc}
\frac{\mathrm{P}^{\alpha }\left( A^{\prime },\mathbf{1}_{A\setminus I_{%
\mathcal{P}_{\func{cor}}^{A}}}\sigma \right) }{\left\vert A^{\prime
}\right\vert }\leq \frac{\mathrm{P}^{\alpha }\left( A^{\prime },\mathbf{1}%
_{A}\sigma \right) }{\left\vert A^{\prime }\right\vert } & \text{ if } & 
A^{\prime }=J_{-/+}^{\maltese }\text{ or }J_{+/-}^{\maltese } \\ 
\frac{\mathrm{P}^{\alpha }\left( A_{J}^{\prime },\mathbf{1}_{A\setminus I_{%
\mathcal{P}_{\func{cor}}^{A}}}\sigma \right) }{\left\vert A_{J}^{\prime
}\right\vert }\lesssim \frac{\mathrm{P}^{\alpha }\left( A^{\prime },\mathbf{1%
}_{A}\sigma \right) }{\left\vert A^{\prime }\right\vert } & \text{ if } & 
A^{\prime }=J_{-}^{\maltese }\text{ or }J_{+}^{\maltese }%
\end{array}%
\right. .
\end{equation*}%
Since $\left\vert \varphi _{J}^{\mathcal{P}_{\func{cor}}^{A}}\right\vert
\lesssim \alpha _{\mathcal{A}}\left( A\right) \mathbf{1}_{A}$ by (\ref{phi
bound}), we can then bound $\left\vert \mathsf{B}\right\vert _{\limfunc{stop}%
,\bigtriangleup ^{\omega }}^{A,\mathcal{P}_{\func{cor}}^{A}}\left(
f,g\right) $ by%
\begin{eqnarray*}
&&\alpha _{\mathcal{A}}\left( A\right) \sum_{A^{\prime }\in \mathfrak{C}_{%
\mathcal{A}}\left( A\right) }\left( \frac{\mathrm{P}^{\alpha }\left(
A^{\prime },\mathbf{1}_{A}\sigma \right) }{\left\vert A^{\prime }\right\vert 
}\right) \left\Vert \mathsf{Q}_{\Pi _{2}\mathcal{P}_{\func{cor}%
}^{A};A^{\prime }}^{\omega ,\mathbf{b}^{\ast }}x\right\Vert _{L^{2}\left(
\omega \right) }^{\spadesuit }\left\Vert \mathsf{P}_{\Pi _{2}\mathcal{P}_{%
\func{cor}}^{A};A^{\prime }}^{\omega ,\mathbf{b}^{\ast }}g\right\Vert
_{L^{2}\left( \omega \right) }^{\bigstar } \\
&\leq &\alpha _{\mathcal{A}}\left( A\right) \left( \sum_{A^{\prime }\in 
\mathfrak{C}_{\mathcal{A}}\left( A\right) }\left( \frac{\mathrm{P}^{\alpha
}\left( A^{\prime },\mathbf{1}_{A}\sigma \right) }{\left\vert A^{\prime
}\right\vert }\right) ^{2}\left\Vert x-m_{A^{\prime }}^{\sigma }\right\Vert
_{L^{2}\left( \mathbf{1}_{A^{\prime }}\sigma \right) }^{\spadesuit 2}\right)
^{\frac{1}{2}} \\
&&\ \ \ \ \ \ \ \ \ \ \ \ \ \ \ \ \ \ \ \ \ \ \ \ \ \ \ \ \ \ \ \ \ \ \ \ \
\ \ \ \ \ \ \ \ \times \left( \sum_{A^{\prime }\in \mathfrak{C}_{\mathcal{A}%
}\left( A\right) }\left\Vert \mathsf{P}_{\Pi _{2}\mathcal{P}_{\func{cor}%
}^{A};A^{\prime }}^{\omega ,\mathbf{b}^{\ast }}g\right\Vert _{L^{2}\left(
\omega \right) }^{\bigstar 2}\right) ^{\frac{1}{2}} \\
&\leq &\mathcal{E}_{2}^{\alpha }\alpha _{\mathcal{A}}\left( A\right) \sqrt{%
\left\vert A\right\vert _{\sigma }}\left\Vert \mathsf{P}_{\Pi _{2}\mathcal{P}%
_{\func{cor}}^{A}}^{\omega ,\mathbf{b}^{\ast }}g\right\Vert _{L^{2}\left(
\omega \right) }^{\bigstar }\leq \mathcal{E}_{2}^{\alpha }\alpha _{\mathcal{A%
}}\left( A\right) \sqrt{\left\vert A\right\vert _{\sigma }}\left\Vert 
\mathsf{P}_{\mathcal{C}_{A}^{\func{shift}}}^{\omega ,\mathbf{b}^{\ast
}}g\right\Vert _{L^{2}\left( \omega \right) }^{\bigstar },
\end{eqnarray*}%
where in the last line we have used the strong energy constant $\mathcal{E}%
_{2}^{\alpha }$ in (\ref{strong b* energy}).
\end{proof}

\begin{definition}
We say that an admissible collection of pairs $\mathcal{P}$ is \emph{reduced}
if it contains no pairs from $\mathcal{P}_{\func{cor}}^{A}$, i.e.%
\begin{equation*}
\mathcal{P}\cap \mathcal{P}_{\func{cor}}^{A}=\emptyset .
\end{equation*}
\end{definition}

\begin{definition}
\label{def aug}We define $J^{\flat }=J_{\searrow J}^{\maltese }$ to be the
inner grandchild of $J^{\maltese }$ that contains $J$.
\end{definition}

Recall that in terms of $J^{\flat }$ we rewrite%
\begin{equation*}
\Pi _{2}^{K,\limfunc{aug}}\mathcal{P}=\left\{ J\in \Pi _{2}\mathcal{P}%
:J\subset K\text{ and }J^{\maltese }\subset \pi _{\mathcal{D}}^{\left(
2\right) }K\right\} =\left\{ J\in \Pi _{2}\mathcal{P}:J\subset K\text{ and }%
J^{\flat }\subset K\right\} .
\end{equation*}

\begin{definition}
\label{flat straddles}Given a \emph{reduced admissible} collection of pairs $%
\mathcal{Q}$ for $A$, and a subpartition $\mathcal{S}\subset \Pi _{1}^{%
\limfunc{below}}\mathcal{Q}\cap \mathcal{C}_{A}^{\limfunc{restrict}}$ of
pairwise disjoint intervals in $A$, we say that $\mathcal{Q}$ $\flat $%
\textbf{straddles} $\mathcal{S}$ if for every pair $\left( I,J\right) \in 
\mathcal{Q}$ there is $S\in \mathcal{S}\cap \left[ J,I\right] $ with $%
J^{\flat }\subset S$. To avoid trivialities, we further assume that for
every $S\in \mathcal{S}$, there is at least one pair $\left( I,J\right) \in 
\mathcal{Q}$ with $J^{\flat }\subset S\subset I$. Here $\left[ J,I\right] $
denotes the geodesic in the dyadic tree $\mathcal{D}$ that connects $J^{%
\mathcal{D}}$ to $I$, where $J^{\mathcal{D}}$ is the minimal interval in $%
\mathcal{D}$ that contains $J$.
\end{definition}

\begin{definition}
\label{def Whit}For any dyadic interval $S\in \mathcal{D}$, define the
Whitney collection $\mathcal{W}\left( S\right) $ to consist of the maximal
subintervals $K$ of $S$ whose triples $3K$ are contained in $S$. Then set $%
\mathcal{W}^{\ast }\left( S\right) \equiv \mathcal{W}\left( S\right) \cup
\left\{ S\right\} $.
\end{definition}

The following geometric proposition will prove useful in proving the $\flat $%
Straddling Lemma \ref{straddle 3 ref} below.

\begin{proposition}
\label{flatness}Suppose $\mathcal{Q}$ is reduced admissible and $\flat $%
straddles a subpartition $\mathcal{S}$ of $A$. Fix $S\in \mathcal{S}$.
Define 
\begin{equation*}
\varphi _{J}^{\mathcal{Q}^{S}}\left[ h\right] \equiv \sum_{I\in \Pi _{1}%
\mathcal{Q}^{S}:\mathcal{\ }\left( I,J\right) \in \mathcal{Q}%
^{S}}b_{A}E_{I}^{\sigma }\left( \widehat{\square }_{\pi I}^{\sigma ,\flat ,%
\mathbf{b}}h\right) \ \mathbf{1}_{A\setminus I}\ ,
\end{equation*}%
assume that $h\in L^{2}\left( \sigma \right) $ is supported in the interval $%
A$, and that there is an interval $H\in \mathcal{C}_{A}$ with $H\supset S$
such that 
\begin{equation*}
E_{I}^{\sigma }\left\vert h\right\vert \leq CE_{H}^{\sigma }\left\vert
h\right\vert ,\ \ \ \ \ \text{for all }I\in \Pi _{1}^{\limfunc{below}}%
\mathcal{Q}\cap \mathcal{C}_{A}^{\limfunc{restrict}}\text{ with }I\supset S.
\end{equation*}%
Then%
\begin{eqnarray*}
&&\sum_{J\in \Pi _{2}\mathcal{Q}:\ J^{\flat }\subset S}\frac{\mathrm{P}%
^{\alpha }\left( J,\left\vert \varphi _{J}^{\mathcal{Q}}\left[ h\right]
\right\vert \mathbf{1}_{A\setminus I_{\mathcal{Q}}\left( J\right) }\sigma
\right) }{\left\vert J\right\vert }\left\Vert \bigtriangleup _{J}^{\omega ,%
\mathbf{b}^{\ast }}x\right\Vert _{L^{2}\left( \omega \right) }^{\spadesuit
}\left\Vert \square _{J}^{\omega ,\mathbf{b}^{\ast }}g\right\Vert
_{L^{2}\left( \omega \right) }^{\bigstar } \\
&\lesssim &\alpha _{\mathcal{H}}\left( H\right) \frac{\mathrm{P}^{\alpha
}\left( S,\mathbf{1}_{A\setminus S}\sigma \right) }{\left\vert S\right\vert }%
\left\Vert \mathsf{Q}_{\Pi _{2}^{S,\limfunc{aug}}\mathcal{Q}}^{\omega ,%
\mathbf{b}^{\ast }}x\right\Vert _{L^{2}\left( \omega \right) }^{\spadesuit
}\left\Vert \mathsf{P}_{\Pi _{2}^{S,\limfunc{aug}}\mathcal{Q}}^{\omega ,%
\mathbf{b}^{\ast }}g\right\Vert _{L^{2}\left( \omega \right) }^{\bigstar } \\
&&+\alpha _{\mathcal{H}}\left( H\right) \dsum\limits_{K\in \mathcal{W}\left(
S\right) }\frac{\mathrm{P}^{\alpha }\left( K,\mathbf{1}_{A\setminus K}\sigma
\right) }{\left\vert K\right\vert }\left\Vert \mathsf{Q}_{\Pi _{2}^{K,%
\limfunc{aug}}\mathcal{Q}}^{\omega ,\mathbf{b}^{\ast }}x\right\Vert
_{L^{2}\left( \omega \right) }^{\spadesuit }\left\Vert \mathsf{P}_{\Pi
_{2}^{K,\limfunc{aug}}\mathcal{Q}}^{\omega ,\mathbf{b}^{\ast }}g\right\Vert
_{L^{2}\left( \omega \right) }^{\bigstar }\ .
\end{eqnarray*}%
The sum over Whitney intervals $K\in \mathcal{W}\left( S\right) $ is only
required to bound the sum of those terms on the left for which $J^{\flat
}\subset S^{\prime \prime }$ for some $S^{\prime \prime }\in \mathfrak{C}_{%
\mathcal{D}}^{\left( 2\right) }\left( S\right) $.
\end{proposition}

\begin{proof}
Suppose first that $J^{\flat }=S\in \mathcal{C}_{A}^{\limfunc{restrict}}$.
Then $3S=3J^{\flat }\subset J^{\maltese }\subset I_{\mathcal{Q}}\left(
J\right) $ and using (\ref{phi bound}) with $\alpha _{\mathcal{H}}\left(
H\right) $ in place of $\alpha _{\mathcal{A}}\left( A\right) $, we have%
\begin{eqnarray*}
\frac{\mathrm{P}^{\alpha }\left( J,\left\vert \varphi _{J}^{\mathcal{Q}%
}\right\vert \mathbf{1}_{A\setminus I_{\mathcal{Q}}\left( J\right) }\sigma
\right) }{\left\vert J\right\vert } &\lesssim &\alpha _{\mathcal{H}}\left(
H\right) \frac{\mathrm{P}^{\alpha }\left( J,\mathbf{1}_{A\setminus
J^{\maltese }}\sigma \right) }{\left\vert J\right\vert } \\
&\lesssim &\alpha _{\mathcal{H}}\left( H\right) \frac{\mathrm{P}^{\alpha
}\left( S,\mathbf{1}_{A\setminus J^{\maltese }}\sigma \right) }{\left\vert
S\right\vert }\leq \alpha _{\mathcal{H}}\left( H\right) \frac{\mathrm{P}%
^{\alpha }\left( S,\mathbf{1}_{A\setminus S}\sigma \right) }{\left\vert
S\right\vert }.
\end{eqnarray*}%
Suppose next that $J^{\flat }=S^{\prime }\in \mathfrak{C}_{\mathcal{D}%
}\left( S\right) $. Then $3S^{\prime }=3J^{\flat }\subset J^{\maltese
}\subset I_{\mathcal{Q}}\left( J\right) $ and (\ref{phi bound}) give%
\begin{eqnarray*}
\frac{\mathrm{P}^{\alpha }\left( J,\left\vert \varphi _{J}^{\mathcal{Q}%
}\right\vert \mathbf{1}_{A\setminus I_{\mathcal{Q}}\left( J\right) }\sigma
\right) }{\left\vert J\right\vert } &\lesssim &\alpha _{\mathcal{H}}\left(
H\right) \frac{\mathrm{P}^{\alpha }\left( J,\mathbf{1}_{A\setminus
J^{\maltese }}\sigma \right) }{\left\vert J\right\vert } \\
&\lesssim &\alpha _{\mathcal{H}}\left( H\right) \frac{\mathrm{P}^{\alpha
}\left( S^{\prime },\mathbf{1}_{A\setminus J^{\maltese }}\sigma \right) }{%
\left\vert S^{\prime }\right\vert } \\
&\leq &\alpha _{\mathcal{H}}\left( H\right) \frac{\mathrm{P}^{\alpha }\left(
S^{\prime },\mathbf{1}_{A\setminus S}\sigma \right) }{\left\vert S^{\prime
}\right\vert }\approx \alpha _{\mathcal{H}}\left( H\right) \frac{\mathrm{P}%
^{\alpha }\left( S,\mathbf{1}_{A\setminus S}\sigma \right) }{\left\vert
S\right\vert }.
\end{eqnarray*}%
Thus in these two cases, by Cauchy-Schwarz, the left hand side of our
conclusion is bounded by a multiple of%
\begin{eqnarray*}
&&\alpha _{\mathcal{H}}\left( H\right) \frac{\mathrm{P}^{\alpha }\left( S,%
\mathbf{1}_{A\setminus S}\sigma \right) }{\left\vert S\right\vert }\left(
\sum_{J\in \Pi _{2}\mathcal{Q}:\ J^{\flat }\subset S}\left\Vert
\bigtriangleup _{J}^{\omega ,\mathbf{b}^{\ast }}x\right\Vert _{L^{2}\left(
\omega \right) }^{\spadesuit 2}\right) ^{\frac{1}{2}}\left( \sum_{J\in \Pi
_{2}\mathcal{Q}:\ J^{\flat }\subset S}\left\Vert \square _{J}^{\omega ,%
\mathbf{b}^{\ast }}g\right\Vert _{L^{2}\left( \omega \right) }^{\bigstar
2}\right) ^{\frac{1}{2}} \\
&=&\alpha _{\mathcal{H}}\left( H\right) \frac{\mathrm{P}^{\alpha }\left( S,%
\mathbf{1}_{A\setminus S}\sigma \right) }{\left\vert S\right\vert }%
\left\Vert \mathsf{Q}_{\Pi _{2}^{S,\limfunc{aug}}\mathcal{Q}}^{\omega ,%
\mathbf{b}^{\ast }}x\right\Vert _{L^{2}\left( \omega \right) }^{\spadesuit
}\left\Vert \mathsf{P}_{\Pi _{2}^{S,\limfunc{aug}}\mathcal{Q}}^{\omega ,%
\mathbf{b}^{\ast }}g\right\Vert _{L^{2}\left( \omega \right) }^{\bigstar }\ .
\end{eqnarray*}%
Finally, suppose that $J^{\flat }\subset S^{\prime \prime }$ for some $%
S^{\prime \prime }\in \mathfrak{C}_{\mathcal{D}}^{\left( 2\right) }\left(
S\right) $. Then $J^{\maltese }\subset S$, and Key Fact \#2 in (\ref%
{indentation}) shows that $3J^{\flat }\subset J^{\maltese }$, so that $%
3J^{\flat }\subset J^{\maltese }\subset S\subset I_{\mathcal{Q}}\left(
J\right) $. Thus we have $J^{\flat }\subset K=K\left[ J\right] $ for some $%
K\in \mathcal{W}\left( S\right) $ and so by (\ref{phi bound}) again,%
\begin{eqnarray*}
\frac{\mathrm{P}^{\alpha }\left( J,\left\vert \varphi _{J}^{\mathcal{Q}%
}\right\vert \mathbf{1}_{A\setminus I_{\mathcal{Q}}\left( J\right) }\sigma
\right) }{\left\vert J\right\vert } &\lesssim &\alpha _{\mathcal{H}}\left(
H\right) \frac{\mathrm{P}^{\alpha }\left( J,\mathbf{1}_{A\setminus S}\sigma
\right) }{\left\vert J\right\vert } \\
&\lesssim &\alpha _{\mathcal{H}}\left( H\right) \frac{\mathrm{P}^{\alpha
}\left( K,\mathbf{1}_{A\setminus S}\sigma \right) }{\left\vert K\right\vert }%
\leq \alpha _{\mathcal{H}}\left( H\right) \frac{\mathrm{P}^{\alpha }\left( K,%
\mathbf{1}_{A\setminus K}\sigma \right) }{\left\vert K\right\vert }.
\end{eqnarray*}%
Now we apply Cauchy-Schwarz again, but noting that $J^{\flat }\subset K$
this time, to obtain that the left hand side of our conclusion is bounded by
a multiple of%
\begin{eqnarray*}
&&\alpha _{\mathcal{H}}\left( H\right) \dsum\limits_{K\in \mathcal{W}\left(
S\right) }\frac{\mathrm{P}^{\alpha }\left( K,\mathbf{1}_{A\setminus K}\sigma
\right) }{\left\vert K\right\vert }\left( \sum_{J\in \Pi _{2}\mathcal{Q}:\
J^{\flat }\subset K}\left\Vert \bigtriangleup _{J}^{\omega ,\mathbf{b}^{\ast
}}x\right\Vert _{L^{2}\left( \omega \right) }^{\spadesuit 2}\right) ^{\frac{1%
}{2}}\left( \sum_{J\in \Pi _{2}\mathcal{Q}:\ J^{\flat }\subset K}\left\Vert
\square _{J}^{\omega ,\mathbf{b}^{\ast }}g\right\Vert _{L^{2}\left( \omega
\right) }^{\bigstar 2}\right) ^{\frac{1}{2}} \\
&=&\alpha _{\mathcal{H}}\left( H\right) \dsum\limits_{K\in \mathcal{W}\left(
S\right) }\frac{\mathrm{P}^{\alpha }\left( K,\mathbf{1}_{A\setminus K}\sigma
\right) }{\left\vert K\right\vert }\left\Vert \mathsf{Q}_{\Pi _{2}^{K,%
\limfunc{aug}}\mathcal{Q}}^{\omega ,\mathbf{b}^{\ast }}x\right\Vert
_{L^{2}\left( \omega \right) }^{\spadesuit }\left\Vert \mathsf{P}_{\Pi
_{2}^{K,\limfunc{aug}}\mathcal{Q}}^{\omega ,\mathbf{b}^{\ast }}g\right\Vert
_{L^{2}\left( \omega \right) }^{\bigstar }\ .
\end{eqnarray*}%
This completes the proof of Proposition \ref{flatness}.
\end{proof}

Recall the family of operators $\left\{ \square _{I}^{\sigma ,\pi ,\mathbf{b}%
}\right\} _{I\in \mathcal{C}_{A}^{\mathcal{A}}}$, where for $I\in \mathcal{C}%
_{A}^{\mathcal{A}}$, the dual martingale difference $\square _{I}^{\sigma
,\pi ,\mathbf{b}}$ is defined in (\ref{def pi box}) of Appendix A below, and
satisfies%
\begin{equation*}
\square _{I}^{\sigma ,\pi ,\mathbf{b}}f=\left[ \sum_{I^{\prime }\in 
\mathfrak{C}\left( I\right) }\mathbb{F}_{I^{\prime }}^{\sigma ,\pi ,\mathbf{b%
}}f\right] -\mathbb{F}_{I}^{\sigma ,\mathbf{b}}f=\sum_{I^{\prime }\in 
\mathfrak{C}\left( I\right) }\mathbb{F}_{I^{\prime }}^{\sigma ,b_{A}}f-%
\mathbb{F}_{I}^{\sigma ,b_{A}}f\ .
\end{equation*}%
Since $\square _{I}^{\sigma ,\pi ,\mathbf{b}}$ is the transpose of $%
\triangle _{I}^{\sigma ,\pi ,\mathbf{b}}$ for $I\in \mathcal{C}_{A}^{%
\mathcal{A}}$, the first line of Lemma \ref{b proj} (where the superscript $%
\pi $ is suppressed for convenience) shows that $\left\{ \square
_{I}^{\sigma ,\pi ,\mathbf{b}}\right\} _{I\in \mathcal{C}_{A}^{\mathcal{A}}}$
is a family of projections, and the second line of Lemma \ref{b proj} shows
it is an orthogonal family, i.e. 
\begin{equation*}
\square _{I}^{\sigma ,\pi ,\mathbf{b}}\square _{K}^{\sigma ,\pi ,\mathbf{b}%
}=\left\{ 
\begin{array}{ccc}
\square _{I}^{\sigma ,\pi ,\mathbf{b}} & \text{ if } & I=K \\ 
0 & \text{ if } & I\not=K%
\end{array}%
\right. ,\ \ \ \ \ I,K\in \mathcal{C}_{A}^{\mathcal{A}}.
\end{equation*}%
The orthogonal projections 
\begin{eqnarray*}
\mathsf{P}_{\pi \left( \Pi _{1}\mathcal{Q}\right) }^{\sigma ,\pi ,\mathbf{b}%
} &\equiv &\sum_{I\in \pi \left( \Pi _{1}\mathcal{Q}\right) }\square
_{I}^{\sigma ,\pi ,\mathbf{b}}=\sum_{I\in \Pi _{1}\mathcal{Q}}\square _{\pi
I}^{\sigma ,\pi ,\mathbf{b}}, \\
\text{where }\pi \left( \Pi _{1}\mathcal{Q}\right) &\equiv &\left\{ \pi _{%
\mathcal{D}}I:I\in \Pi _{1}\mathcal{Q}\right\} \text{ and }\Pi _{1}\mathcal{Q%
}\subset \mathcal{C}_{A}^{\mathcal{A},\limfunc{restrict}}\ ,
\end{eqnarray*}%
thus satisfy the equalities%
\begin{equation}
\square _{\pi I}^{\sigma ,\pi ,\mathbf{b}}f=\square _{\pi I}^{\sigma ,\pi ,%
\mathbf{b}}\mathsf{P}_{\pi \left( \Pi _{1}\mathcal{Q}\right) }^{\sigma ,\pi ,%
\mathbf{b}}f\text{ and }\widehat{\square }_{\pi I}^{\sigma ,\pi ,\mathbf{b}%
}f=\widehat{\square }_{\pi I}^{\sigma ,\pi ,\mathbf{b}}\mathsf{P}_{\pi
\left( \Pi _{1}\mathcal{Q}\right) }^{\sigma ,\pi ,\mathbf{b}}f,\ \ \ \ \ 
\text{for }I\in \Pi _{1}\mathcal{Q}\subset \mathcal{C}_{A}^{\mathcal{A},%
\limfunc{restrict}},  \label{sat}
\end{equation}%
which will permit us to apply certain projection tricks used for Haar
projections in the proof of $T1$ theorems.

However, in our sublinear stopping form $\left\vert \mathsf{B}\right\vert _{%
\limfunc{stop},\bigtriangleup ^{\omega }}^{A,\mathcal{Q}}$, the dual
martingale projections in use in the function%
\begin{equation}
\varphi _{J}^{\mathcal{Q}^{S}}\equiv \sum_{I\in \Pi _{1}\mathcal{Q}^{S}:%
\mathcal{\ }\left( I,J\right) \in \mathcal{Q}^{S}}b_{A}E_{I}^{\sigma }\left( 
\widehat{\square }_{\pi I}^{\sigma ,\flat ,\mathbf{b}}f\right) \ \mathbf{1}%
_{A\setminus I}\ ,  \label{in use}
\end{equation}%
given in (\ref{def phi P}) above, are the modified pseudoprojections $%
\left\{ \widehat{\square }_{\pi I}^{\sigma ,\flat ,\mathbf{b}}\right\}
_{I\in \Pi _{1}\mathcal{Q}}$, where $\square _{\pi I}^{\sigma ,\flat ,%
\mathbf{b}}$ differs from the orthogonal projection $\square _{\pi
I}^{\sigma ,\pi ,\mathbf{b}}$ for $I\in \Pi _{1}\mathcal{Q}$ by%
\begin{equation*}
\square _{\pi I}^{\sigma ,\flat ,\mathbf{b}}f-\square _{\pi I}^{\sigma ,\pi ,%
\mathbf{b}}f=\left\{ \left( \sum_{I^{\prime }\in \mathfrak{C}_{\limfunc{%
natural}}\left( \pi I\right) }\mathbb{F}_{I^{\prime }}^{\sigma
,b_{A}}f\right) -\mathbb{F}_{\pi I}^{\sigma ,b_{A}}f\right\} -\left\{ \left(
\sum_{I^{\prime }\in \mathfrak{C}\left( \pi I\right) }\mathbb{F}_{I^{\prime
}}^{\sigma ,b_{A}}f\right) -\mathbb{F}_{\pi I}^{\sigma ,b_{A}}f\right\}
=-\sum_{I^{\prime }\in \mathfrak{C}_{\limfunc{broken}}\left( \pi I\right) }%
\mathbb{F}_{I^{\prime }}^{\sigma ,b_{A}}f.
\end{equation*}%
But the "box support" $\func{Supp}_{\limfunc{box}}$ of this last expression $%
\sum_{I^{\prime }\in \mathfrak{C}_{\limfunc{broken}}\left( \pi I\right) }%
\mathbb{F}_{I^{\prime }}^{\sigma ,b_{A}}f$ consists of the broken children
of $\pi I$, $\mathfrak{C}_{\limfunc{broken}}\left( \pi I\right) $, and is
contained in the set $\dbigcup\limits_{I\in \mathcal{C}_{A}^{\limfunc{%
restrict}}}\dbigcup\limits_{I^{\prime }\in \mathfrak{C}_{\mathcal{A}}\left(
A\right) \cap \mathfrak{C}_{\mathcal{D}}\left( \pi I\right) }\left\{
I^{\prime }\right\} $, i.e. 
\begin{eqnarray*}
\func{Supp}_{\limfunc{box}}\left( \sum_{I^{\prime }\in \mathfrak{C}_{%
\limfunc{broken}}\left( \pi I\right) }\mathbb{F}_{I^{\prime }}^{\sigma
,b_{A}}f\right) &\subset &\left\{ I^{\prime }\in \mathfrak{C}_{\mathcal{A}%
}\left( A\right) :I^{\prime }\in \mathfrak{C}_{\limfunc{broken}}\left( \pi
I\right) \text{ for some }I\in \mathcal{C}_{A}^{\limfunc{restrict}}\right\}
\\
&=&\dbigcup\limits_{I\in \mathcal{C}_{A}^{\limfunc{restrict}%
}}\dbigcup\limits_{I^{\prime }\in \mathfrak{C}_{\mathcal{A}}\left( A\right)
\cap \mathfrak{C}_{\mathcal{D}}\left( \pi I\right) }\left\{ I^{\prime
}\right\} .
\end{eqnarray*}%
But $I\in \Pi _{1}\mathcal{Q}^{S}\subset \mathcal{C}_{A}^{\limfunc{restrict}%
} $ is a \emph{natural} child of $\pi I$, and so%
\begin{equation*}
I\cap \func{Supp}_{\limfunc{box}}\left( \sum_{I^{\prime }\in \mathfrak{C}_{%
\limfunc{broken}}\left( \pi I\right) }\mathbb{F}_{I^{\prime }}^{\sigma
,b_{A}}f\right) =\emptyset .
\end{equation*}%
It now follows that we have%
\begin{equation}
E_{I}^{\sigma }\left( \widehat{\square }_{\pi I}^{\sigma ,\flat ,\mathbf{b}%
}f\right) =E_{I}^{\sigma }\left( \widehat{\square }_{\pi I}^{\sigma ,\pi ,%
\mathbf{b}}f\right) ,\ \ \ \ \ \text{for }I\in \mathcal{C}_{A}^{\limfunc{%
restrict}}.  \label{fol}
\end{equation}

Returning to (\ref{in use}), we have from (\ref{sat}) and (\ref{fol}) the
identity 
\begin{eqnarray}
\varphi _{J}^{\mathcal{Q}^{S}} &\equiv &\sum_{I\in \Pi _{1}\mathcal{Q}^{S}:%
\mathcal{\ }\left( I,J\right) \in \mathcal{Q}^{S}}b_{A}E_{I}^{\sigma }\left( 
\widehat{\square }_{\pi I}^{\sigma ,\pi ,\mathbf{b}}f\right) \ \mathbf{1}%
_{A\setminus I}  \label{iden} \\
&=&\sum_{I\in \Pi _{1}\mathcal{Q}^{S}:\mathcal{\ }\left( I,J\right) \in 
\mathcal{Q}^{S}}b_{A}E_{I}^{\sigma }\left( \widehat{\square }_{\pi
I}^{\sigma ,\pi ,\mathbf{b}}\left( \mathsf{P}_{\pi \left( \Pi _{1}\mathcal{Q}%
\right) }^{\sigma ,\pi ,\mathbf{b}}f\right) \right) \ \mathbf{1}_{A\setminus
I}\ ,  \notag
\end{eqnarray}%
which will play a critical role in proving the following $\flat $Straddling
and Substraddling lemmas. The $\flat $Straddling Lemma is an adaptation of
Lemmas 3.19 and 3.16 in \cite{Lac}.

\begin{lemma}
\label{straddle 3 ref}Let $\mathcal{Q}$ be a reduced admissible collection
of pairs for $A$, and suppose that $\mathcal{S}\subset \Pi _{1}^{\limfunc{%
below}}\mathcal{Q}\cap \mathcal{C}_{A}^{\limfunc{restrict}}$ is a
subpartition of $A$ such that $\mathcal{Q}$ $\flat $straddles $\mathcal{S}$.
Then we have the restricted sublinear norm bound%
\begin{equation}
\widehat{\mathfrak{N}}_{\limfunc{stop},\bigtriangleup ^{\omega }}^{A,%
\mathcal{Q}}\leq C_{\mathbf{r}}\sup_{S\in \mathcal{S}}\mathcal{S}_{\limfunc{%
loc}\limfunc{size}}^{\alpha ,A;S}\left( \mathcal{Q}\right) \leq C_{\mathbf{r}%
}\mathcal{S}_{\limfunc{aug}\limfunc{size}}^{\alpha ,A}\left( \mathcal{Q}%
\right) ,  \label{sub loc bound}
\end{equation}%
where $\mathcal{S}_{\limfunc{loc}\limfunc{size}}^{\alpha ,A;S}$ is an $S$%
-localized size condition with an $S$-hole given by%
\begin{equation}
\mathcal{S}_{\limfunc{loc}\limfunc{size}}^{\alpha ,A;S}\left( \mathcal{Q}%
\right) ^{2}\equiv \sup_{K\in \mathcal{W}^{\ast }\left( S\right) \cap 
\mathcal{C}_{A}^{\limfunc{restrict}}}\frac{1}{\left\vert K\right\vert
_{\sigma }}\left( \frac{\mathrm{P}^{\alpha }\left( K,\mathbf{1}_{A\setminus
S}\sigma \right) }{\left\vert K\right\vert }\right) ^{2}\sum_{J\in \Pi
_{2}^{K,\limfunc{aug}}\mathcal{Q}}\left\Vert \bigtriangleup _{J}^{\omega ,%
\mathbf{b}^{\ast }}x\right\Vert _{L^{2}\left( \omega \right) }^{\spadesuit
2}.  \label{localized size ref}
\end{equation}
\end{lemma}

\begin{proof}
For $S\in \mathcal{S}$ let $\mathcal{Q}^{S}\equiv \left\{ \left( I,J\right)
\in \mathcal{Q}:J^{\flat }\subset S\subset I\right\} $. We begin by using
that the reduced collection $\mathcal{Q}$ straddles $\mathcal{S}$, together
with the sublinearity property (\ref{phi sublinear}) of $\varphi _{J}^{%
\mathcal{Q}}$, and with $\left\vert \mathsf{B}\right\vert _{\limfunc{stop}%
,\bigtriangleup ^{\omega }}^{A,\mathcal{Q}}\left( f,g\right) $ as in (\ref%
{def mod B}), to write%
\begin{eqnarray*}
\left\vert \mathsf{B}\right\vert _{\limfunc{stop},\bigtriangleup ^{\omega
}}^{A,\mathcal{Q}}\left( f,g\right) &=&\sum_{J\in \Pi _{2}\mathcal{Q}}\frac{%
\mathrm{P}^{\alpha }\left( J,\left\vert \varphi _{J}^{\mathcal{Q}%
}\right\vert \mathbf{1}_{A\setminus I_{\mathcal{Q}}\left( J\right) }\sigma
\right) }{\left\vert J\right\vert }\left\Vert \bigtriangleup _{J}^{\omega ,%
\mathbf{b}^{\ast }}x\right\Vert _{L^{2}\left( \omega \right) }^{\spadesuit
}\left\Vert \square _{J}^{\omega ,\mathbf{b}^{\ast }}g\right\Vert
_{L^{2}\left( \omega \right) }^{\bigstar } \\
&\leq &\sum_{S\in \mathcal{S}}\sum_{J\in \Pi _{2}^{S,\limfunc{aug}}\mathcal{Q%
}}\frac{\mathrm{P}^{\alpha }\left( J,\left\vert \varphi _{J}^{\mathcal{Q}%
^{S}}\right\vert \mathbf{1}_{A\setminus I_{\mathcal{Q}}\left( J\right)
}\sigma \right) }{\left\vert J\right\vert }\left\Vert \bigtriangleup
_{J}^{\omega ,\mathbf{b}^{\ast }}x\right\Vert _{L^{2}\left( \omega \right)
}^{\spadesuit }\left\Vert \square _{J}^{\omega ,\mathbf{b}^{\ast
}}g\right\Vert _{L^{2}\left( \omega \right) }^{\bigstar }; \\
\text{where }\varphi _{J}^{\mathcal{Q}^{S}} &\equiv &\sum_{I\in \Pi _{1}%
\mathcal{Q}^{S}:\mathcal{\ }\left( I,J\right) \in \mathcal{Q}%
^{S}}b_{A}E_{I}^{\sigma }\left( \widehat{\square }_{\pi I}^{\sigma ,\flat ,%
\mathbf{b}}f\right) \ \mathbf{1}_{A\setminus I}\ .
\end{eqnarray*}

At this point we invoke the identity (\ref{iden}),%
\begin{equation*}
\varphi _{J}^{\mathcal{Q}^{S}}=\sum_{I\in \Pi _{1}\mathcal{Q}^{S}:\mathcal{\ 
}\left( I,J\right) \in \mathcal{Q}^{S}}b_{A}E_{I}^{\sigma }\left( \widehat{%
\square }_{\pi I}^{\sigma ,\pi ,\mathbf{b}}\left( \mathsf{P}_{\pi \left( \Pi
_{1}\mathcal{Q}\right) }^{\sigma ,\pi ,\mathbf{b}}f\right) \right) \ \mathbf{%
1}_{A\setminus I}\ ,
\end{equation*}%
so that%
\begin{equation*}
\left\vert \mathsf{B}\right\vert _{\limfunc{stop},\bigtriangleup ^{\omega
}}^{A,\mathcal{Q}}\left( f,g\right) =\left\vert \mathsf{B}\right\vert _{%
\limfunc{stop},\bigtriangleup ^{\omega }}^{A,\mathcal{Q}}\left( h,g\right)
,\ \ \ \ \ \text{where }h\equiv \mathsf{P}_{\pi \left( \Pi _{1}\mathcal{Q}%
\right) }^{\sigma ,\pi ,\mathbf{b}}f\ .
\end{equation*}%
We will treat the sublinear form $\left\vert \mathsf{B}\right\vert _{%
\limfunc{stop},\bigtriangleup ^{\omega }}^{A,\mathcal{Q}}\left( h,g\right) $
with $h=\mathsf{P}_{\pi \left( \Pi _{1}\mathcal{Q}\right) }^{\sigma ,\pi ,%
\mathbf{b}}f$ using a small variation on the corresponding argument in Lacey 
\cite{Lac}\footnote{%
There is a gap in the treatment of the Straddling Lemma 11.10 on page 166 of 
\cite{SaShUr7}. The wrong restricted norm is used there, but can be fixed by
using the corresponding argument of Lacey in \cite{Lac}, equivalently
adapting the argument here. See Appendix C for a full discussion.}. Namely,
we will apply a Calder\'{o}n-Zygmund stopping time decomposition to the
function $h=\mathsf{P}_{\pi \left( \Pi _{1}\mathcal{Q}\right) }^{\sigma ,\pi
,\mathbf{b}}f$ on the interval $A$ \ with `obstacle' $\mathcal{S}\cup 
\mathfrak{C}_{A}$ $\left( A\right) $, to obtain stopping times $\mathcal{H}$ 
$\subset \mathcal{C}_{A}$ with the property that for all $H\in \mathcal{H}%
\setminus \left\{ A\right\} $ we have 
\begin{eqnarray*}
&&H\in \mathcal{C}_{A}\text{ is not strictly contained in any interval from }%
\mathcal{S}, \\
&&E_{H}^{\sigma }\left\vert h\right\vert >\Gamma E_{\pi _{\mathcal{H}%
}H}^{\sigma }\left\vert h\right\vert , \\
&&E_{H^{\prime }}^{\sigma }\left\vert h\right\vert \leq \Gamma E_{\pi _{%
\mathcal{H}}H}^{\sigma }\left\vert h\right\vert \text{ for all }H\subsetneqq
H^{\prime }\subset \pi _{\mathcal{H}}H\text{ with }H^{\prime }\in \mathcal{C}%
_{A}.
\end{eqnarray*}%
More precisely, define generation $0$ of $\mathcal{H}$ to consist of the
single interval $A$. Having defined generation $n$, let generation $n+1$
consist of the union over all intervals $M$ in generation $n$ of the maximal
intervals $M^{\prime }$ in $\mathcal{C}_{A}$ that are contained in $M$ with $%
E_{M^{\prime }}^{\sigma }\left\vert h\right\vert >\Gamma E_{M}^{\sigma
}\left\vert h\right\vert $, but are \emph{not} strictly contained in any
interval $S$ from $\mathcal{S}$ or contained in any interval $A^{\prime }$
from $\mathfrak{C}_{A}$ $\left( A\right) $ - thus the construction stops at
the obstacle $\mathcal{S}\cup \mathfrak{C}_{A}$ $\left( A\right) $. Then $%
\mathcal{H}$ is the union of all generations $n\geq 0$.

Denote by 
\begin{equation*}
\mathcal{C}_{H}^{\mathcal{H}}\equiv \left\{ H^{\prime }\in \mathcal{C}%
_{A}:H^{\prime }\subset H\text{ but }H^{\prime }\not\subset H^{\prime \prime
}\text{ for any }H^{\prime \prime }\in \mathfrak{C}_{\mathcal{H}}\left(
H\right) \right\}
\end{equation*}%
the usual $\mathcal{H}$-corona associated with the stopping interval $H$,
but restricted to $\mathcal{C}_{A}$, and let $\alpha _{\mathcal{H}}\left(
H\right) =E_{H}^{\sigma }\left\vert f\right\vert $ as is customary for a
Calder\'{o}n-Zygmund corona. Since these coronas $\mathcal{C}_{H}^{\mathcal{H%
}}$ are all contained in $\mathcal{C}_{A}$, we have the stopping energy from
the $\mathcal{A}$-corona $\mathcal{C}_{A}$ at our disposal, which as in \cite%
{Lac}, is crucial for the argument. Furthermore, we denote by%
\begin{equation}
\mathcal{Q}_{H}\equiv \left\{ \left( I,J\right) \in \mathcal{Q}:J\in 
\mathcal{C}_{H}^{\mathcal{H},\flat \func{shift}}\right\} ,\ \ \ \ \ \text{%
with }\mathcal{C}_{H}^{\mathcal{H},\flat \func{shift}}\equiv \left\{ J\in
\Pi _{2}\mathcal{Q}:J^{\flat }\in \mathcal{C}_{H}^{\mathcal{H}}\right\} ,
\label{def Q H}
\end{equation}%
the restriction of the pairs $\left( I,J\right) $ in $\mathcal{Q}$ to those
for which $J$ lies in the flat shifted $\mathcal{H}$-corona $\mathcal{C}%
_{H}^{\mathcal{H},\flat \func{shift}}$. Since the $\mathcal{H}$-stopping
intervals satisfy a $\sigma $-Carleson condition for $\Gamma $ chosen large
enough, we have the quasiorthogonal inequality 
\begin{equation}
\sum_{H\in \mathcal{H}}\alpha _{\mathcal{H}}\left( H\right) ^{2}\left\vert
H\right\vert _{\sigma }\lesssim \left\Vert h\right\Vert _{L^{2}\left( \sigma
\right) }^{2},  \label{qor}
\end{equation}%
which below we will see reduces matters to proving inequality (\ref{sub loc
bound}) for the family of reduced admissible collections $\left\{ \mathcal{Q}%
_{H}\right\} _{H\in \mathcal{H}}$ with constants independent of $H$:%
\begin{equation*}
\widehat{\mathfrak{N}}_{\limfunc{stop},\bigtriangleup ^{\omega }}^{A,%
\mathcal{Q}_{H}}\leq C_{\mathbf{r}}\sup_{S\in \mathcal{S}}\mathcal{S}_{%
\limfunc{loc}\limfunc{size}}^{\alpha ,A;S}\left( \mathcal{Q}_{H}\right) \leq
C_{\mathbf{r}}\mathcal{S}_{\limfunc{aug}\limfunc{size}}^{\alpha ,A}\left( 
\mathcal{Q}_{H}\right) ,\ \ \ \ \ H\in \mathcal{H}.
\end{equation*}

Given $S\in \mathcal{S}$, define $H_{S}\in \mathcal{H}$ to be the minimal
interval in $\mathcal{H}$ that contains $S$, and then define 
\begin{equation*}
\mathcal{H}_{\mathcal{S}}\equiv \left\{ H_{S}\in \mathcal{H}:S\in \mathcal{S}%
\right\} .
\end{equation*}%
Note that a given $H\in \mathcal{H}_{\mathcal{S}}$ may have many intervals $%
S\in \mathcal{S}$ such that $H=H_{S}$, and we denote the collection of these
intervals by $\mathcal{S}_{H}\equiv \left\{ S\in \mathcal{S}:H_{S}=H\
\right\} $. We will organize the straddling intervals $\mathcal{S}$ as%
\begin{equation*}
\mathcal{S}=\dbigcup\limits_{H\in \mathcal{H}_{\mathcal{S}%
}}\dbigcup\limits_{S\in \mathcal{S}_{H}}
\end{equation*}%
where each $S\in \mathcal{S}$ occurs exactly once in the union on the right
hand side, i.e. the collections $\left\{ \mathcal{S}_{H}\right\} _{H\in 
\mathcal{H}_{\mathcal{S}}}$ are pairwise disjoint.

We now momentarily fix $H\in \mathcal{H}_{\mathcal{S}}$, and consider the
reduced admissible collection $\mathcal{Q}_{H}$, so that its projection onto
the second component $\Pi _{2}\mathcal{Q}_{H}$ of $\mathcal{Q}_{H}$\ is 
\emph{contained} in the corona $\mathcal{C}_{H}^{\mathcal{H},\flat \func{%
shift}}$. Then the collection $\mathcal{Q}_{H}$ $\flat $straddles the set $%
\mathcal{S}_{H}=\left\{ S\in \mathcal{S}:H_{S}=H\ \right\} $. Moreover, $%
\mathcal{Q}_{H}=\dbigcup\limits_{S\in \mathcal{S}:\ S\subset H}\mathcal{Q}%
_{H}^{S}$ and $\Pi _{2}\mathcal{Q}_{H}^{S}=\Pi _{2}^{S,\limfunc{aug}}%
\mathcal{Q}_{H}$.

Recall that a Whitney interval $K$ was required in the right hand side of
the conclusion of Proposition \ref{flatness} only in the case that $J^{\flat
}\subset S^{\prime \prime }$ for some $S^{\prime \prime }\in \mathfrak{C}_{%
\mathcal{D}}^{\left( 2\right) }\left( S\right) $, which of course implies $%
3J^{\flat }\subset J^{\maltese }\subset S$. In this case we claim that $K\in 
\mathcal{C}_{A}$. Indeed, suppose in order to derive a contradiction, that $%
K\not\in \mathcal{C}_{A}$. Then $J^{\maltese }\not\subset K$, and hence $%
3J^{\maltese }\not\subset S$. Since $J^{\maltese }\subset S$, it follows
that $J^{\maltese }$ shares an endpoint with $S$ (since if not, then $%
3J^{\maltese }\subset S$, a contradiction). Now Key Fact \#2 in (\ref%
{indentation}) implies that the inner grandchild containing $J$, either $%
J_{-/+}^{\maltese }$ or $J_{+/-}^{\maltese }$, is contained in $K$ where $%
K\not\in \mathcal{C}_{A}$. This then implies that the pair $\left(
I,J\right) $ belongs to the corona straddling subcollection $\mathcal{P}_{%
\func{cor}}^{A}$, contradicting the assumption that $\mathcal{Q}$ is reduced.

Thus we have $S\in \Pi _{1}^{\limfunc{below}}\mathcal{Q}\cap \mathcal{C}%
_{A}^{\limfunc{restrict}}$ and $K\in \mathcal{W}\left( S\right) \cap 
\mathcal{C}_{A}^{\limfunc{restrict}}$ and we can use Proposition \ref%
{flatness} with $H=H_{S}$ to bound $\left\vert \mathsf{B}\right\vert _{%
\limfunc{stop},\bigtriangleup ^{\omega }}^{A,\mathcal{Q}}\left( f,g\right) $
by first summing over $H\in \mathcal{H}_{\mathcal{S}}$ and then over $S\in 
\mathcal{S}_{H}$. Indeed, $\mathcal{Q}_{H}$ $\flat $straddles $\mathcal{S}%
_{H}\equiv \left\{ S\in \mathcal{S}:H_{S}=H\ \right\} $, so that $\left\vert
\varphi _{J}^{\mathcal{Q}_{H}}\right\vert \lesssim \alpha _{\mathcal{H}%
}\left( H\right) \mathbf{1}_{A\setminus I_{\mathcal{Q}_{H}}\left( J\right) }$
by (\ref{phi bound}), and so the sum over $S\in \mathcal{S}_{H}$ of the
first term on the right side of the conclusion of Proposition \ref{flatness}
is bounded by%
\begin{eqnarray*}
&&\alpha _{\mathcal{H}}\left( H\right) \sum_{S\in \mathcal{S}_{H}}\sqrt{%
\left\vert S\right\vert _{\sigma }}\frac{1}{\sqrt{\left\vert S\right\vert
_{\sigma }}}\left( \frac{\mathrm{P}^{\alpha }\left( S,\mathbf{1}_{A\setminus
S}\sigma \right) }{\left\vert S\right\vert }\right) \left\Vert \mathsf{Q}%
_{\Pi _{2}^{S,\limfunc{aug}}\mathcal{Q}_{H}}^{\omega ,\mathbf{b}^{\ast
}}x\right\Vert _{L^{2}\left( \omega \right) }^{\spadesuit }\left\Vert 
\mathsf{P}_{\Pi _{2}^{S,\limfunc{aug}}\mathcal{Q}_{H}}^{\omega ,\mathbf{b}%
^{\ast }}g\right\Vert _{L^{2}\left( \omega \right) }^{\bigstar } \\
&\leq &\alpha _{\mathcal{H}}\left( H\right) \left\{ \sup_{S\in \mathcal{S}%
_{H}}\frac{1}{\sqrt{\left\vert S\right\vert _{\sigma }}}\left( \frac{\mathrm{%
P}^{\alpha }\left( S,\mathbf{1}_{A\setminus S}\sigma \right) }{\left\vert
S\right\vert }\right) \left\Vert \mathsf{Q}_{\Pi _{2}^{S,\limfunc{aug}}%
\mathcal{Q}_{H}}^{\omega ,\mathbf{b}^{\ast }}x\right\Vert _{L^{2}\left(
\omega \right) }^{\spadesuit }\right\} \sum_{S\in \mathcal{S}_{H}}\sqrt{%
\left\vert S\right\vert _{\sigma }}\left\Vert \mathsf{P}_{\Pi _{2}^{S,%
\limfunc{aug}}\mathcal{Q}_{H}}^{\omega ,\mathbf{b}^{\ast }}g\right\Vert
_{L^{2}\left( \omega \right) }^{\bigstar } \\
&\leq &\alpha _{\mathcal{H}}\left( H\right) \left\{ \sup_{S\in \mathcal{S}%
_{H}}\mathcal{S}_{\limfunc{loc}\limfunc{size}}^{\alpha ,A;S}\left( \mathcal{Q%
}_{H}\right) \right\} \sqrt{\left\vert H\right\vert _{\sigma }}\left\Vert 
\mathsf{P}_{\Pi _{2}\mathcal{Q}_{H}}^{\omega ,\mathbf{b}^{\ast
}}g\right\Vert _{L^{2}\left( \omega \right) }^{\bigstar }\ ,
\end{eqnarray*}%
where $\Pi _{2}^{K,\limfunc{aug}}\mathcal{Q}_{H}$ is as in Definition \ref%
{def aug}, and the corresponding sum over $S\in \mathcal{S}_{H}$ of the
second term is bounded by%
\begin{eqnarray*}
&&\alpha _{\mathcal{H}}\left( H\right) \sum_{S\in \mathcal{S}_{H}}\sum_{K\in 
\mathcal{W}\left( S\right) \cap \mathcal{C}_{A}^{\limfunc{restrict}}}\sqrt{%
\left\vert K\right\vert _{\sigma }}\frac{1}{\sqrt{\left\vert K\right\vert
_{\sigma }}}\left( \frac{\mathrm{P}^{\alpha }\left( K,\mathbf{1}_{A\setminus
S}\sigma \right) }{\left\vert K\right\vert }\right) \left\Vert \mathsf{Q}%
_{\Pi _{2}^{K,\limfunc{aug}}\mathcal{Q}_{H}^{S}}^{\omega ,\mathbf{b}^{\ast
}}x\right\Vert _{L^{2}\left( \omega \right) }^{\spadesuit }\left\Vert 
\mathsf{P}_{\Pi _{2}^{K,\limfunc{aug}}\mathcal{Q}_{H}^{S}}^{\omega ,\mathbf{b%
}^{\ast }}g\right\Vert _{L^{2}\left( \omega \right) }^{\bigstar } \\
&\lesssim &\alpha _{\mathcal{H}}\left( H\right) \sup_{S\in \mathcal{S}_{%
\mathcal{H}}}\mathcal{S}_{\limfunc{loc}\limfunc{size}}^{\alpha ,A;S}\left( 
\mathcal{Q}_{H}\right) \left( \sum_{S\in \mathcal{S}}\sum_{K\in \mathcal{W}%
\left( S\right) }\left\vert K\right\vert _{\sigma }\right) ^{\frac{1}{2}%
}\left\Vert \mathsf{P}_{\Pi _{2}\mathcal{Q}_{H}}^{\omega ,\mathbf{b}^{\ast
}}g\right\Vert _{L^{2}\left( \omega \right) }^{\bigstar } \\
&\leq &\left\{ \sup_{S\in \mathcal{S}_{\mathcal{H}}}\mathcal{S}_{\limfunc{loc%
}\limfunc{size}}^{\alpha ,A;S}\left( \mathcal{Q}_{H}\right) \right\} \alpha
_{\mathcal{H}}\left( H\right) \sqrt{\left\vert H\right\vert _{\sigma }}%
\left\Vert \mathsf{P}_{\Pi _{2}\mathcal{Q}_{\mathcal{H}}}^{\omega ,\mathbf{b}%
^{\ast }}g\right\Vert _{L^{2}\left( \omega \right) }^{\bigstar }.
\end{eqnarray*}

Using the definition of $\left\vert \mathsf{B}\right\vert _{\limfunc{stop}%
,\bigtriangleup ^{\omega }}^{A,\mathcal{Q}}\left( f,g\right) $ in (\ref{def
mod B}), we now sum the previous inequalities over the intervals $H\in 
\mathcal{H}_{\mathcal{S}}$ to obtain the following string of inequalities
(explained in detail after the display)%
\begin{eqnarray}
\left\vert \mathsf{B}\right\vert _{\limfunc{stop},\bigtriangleup ^{\omega
}}^{A,\mathcal{Q}}\left( f,g\right) &\leq &\left\{ \sup_{S\in \mathcal{S}}%
\mathcal{S}_{\limfunc{loc}\limfunc{size}}^{\alpha ,A;S}\left( \mathcal{Q}%
\right) \right\} \sum_{H\in \mathcal{H}_{\mathcal{S}}}\alpha _{\mathcal{H}%
}\left( H\right) \sqrt{\left\vert H\right\vert _{\sigma }}\left\Vert \mathsf{%
P}_{\Pi _{2}\mathcal{Q}_{H}}^{\omega ,\mathbf{b}^{\ast }}g\right\Vert
_{L^{2}\left( \omega \right) }^{\bigstar }  \label{unfix} \\
&\leq &\left\{ \sup_{S\in \mathcal{S}}\mathcal{S}_{\limfunc{loc}\limfunc{size%
}}^{\alpha ,A;S}\left( \mathcal{Q}\right) \right\} \sqrt{\sum_{H\in \mathcal{%
H}_{\mathcal{S}}}\alpha _{\mathcal{H}}\left( H\right) ^{2}\left\vert
H\right\vert _{\sigma }}\sqrt{\sum_{H\in \mathcal{H}_{\mathcal{S}%
}}\left\Vert \mathsf{P}_{\Pi _{2}\mathcal{Q}_{H}}^{\omega ,\mathbf{b}^{\ast
}}g\right\Vert _{L^{2}\left( \omega \right) }^{\bigstar 2}}  \notag \\
&\lesssim &\left\{ \sup_{S\in \mathcal{S}}\mathcal{S}_{\limfunc{loc}\limfunc{%
size}}^{\alpha ,A;S}\left( \mathcal{Q}\right) \right\} \left\Vert
h\right\Vert _{L^{2}\left( \sigma \right) }\sqrt{\sum_{H\in \mathcal{H}_{%
\mathcal{S}}}\left\Vert \mathsf{P}_{\Pi _{2}\mathcal{Q}_{H}}^{\omega ,%
\mathbf{b}^{\ast }}g\right\Vert _{L^{2}\left( \omega \right) }^{\bigstar 2}}
\notag \\
&\leq &\left\{ \sup_{S\in \mathcal{S}}\mathcal{S}_{\limfunc{loc}\limfunc{size%
}}^{\alpha ,A;S}\left( \mathcal{Q}\right) \right\} \left\Vert \mathsf{P}%
_{\pi \left( \Pi _{1}\mathcal{Q}\right) }^{\sigma ,\pi ,\mathbf{b}%
}f\right\Vert _{L^{2}\left( \sigma \right) }\left\Vert \mathsf{P}_{\Pi _{2}%
\mathcal{Q}}^{\omega ,\mathbf{b}^{\ast }}g\right\Vert _{L^{2}\left( \omega
\right) }^{\bigstar }  \notag \\
&\lesssim &\left\{ \sup_{S\in \mathcal{S}}\mathcal{S}_{\limfunc{loc}\limfunc{%
size}}^{\alpha ,A;S}\left( \mathcal{Q}\right) \right\} \left\Vert \mathsf{P}%
_{\pi \left( \Pi _{1}\mathcal{Q}\right) }^{\sigma ,\mathbf{b}}f\right\Vert
_{L^{2}\left( \sigma \right) }^{\bigstar }\left\Vert \mathsf{P}_{\Pi _{2}%
\mathcal{Q}}^{\omega ,\mathbf{b}^{\ast }}g\right\Vert _{L^{2}\left( \omega
\right) }^{\bigstar }\ ,  \notag
\end{eqnarray}%
where in the first line we have used $\mathcal{Q}=\dbigcup\limits_{H\in 
\mathcal{H}_{\mathcal{S}}}\mathcal{Q}_{H}$, which follows from the fact that
each $J^{\flat }$ is contained in a unique $S\in \mathcal{S}$; in the third
line we have used the quasiorthogonal inequality (\ref{qor}); in the fourth
line we have used that the sets $\Pi _{2}\mathcal{Q}_{H}\subset \mathcal{C}%
_{H}^{\mathcal{H},\flat \func{shift}}$ are pairwise disjoint in $H$ and have
union $\Pi _{2}\mathcal{Q}=$ $\overset{\cdot }{\dbigcup }_{H\in \mathcal{H}_{%
\mathcal{S}}}\Pi _{2}\mathcal{Q}_{H}$. In the final line, we have used first
the equality (\ref{box pi equals}), second the fact that the functions $%
\square _{I,\limfunc{broken}}^{\sigma ,\pi ,\mathbf{b}}f$ have pairwise
disjoint supports, third the upper weak Riesz inequality in Proposition \ref%
{half Riesz}, and fourth the estimate (\ref{F est}) - which relies on the
reverse H\"{o}lder property for children in Lemma \ref{prelim control of
corona} - to obtain%
\begin{eqnarray}
\left\Vert \mathsf{P}_{\pi \left( \Pi _{1}\mathcal{Q}\right) }^{\sigma ,\pi ,%
\mathbf{b}}f\right\Vert _{L^{2}\left( \sigma \right) }^{2} &=&\left\Vert
\sum_{I\in \pi \left( \Pi _{1}\mathcal{Q}\right) }\square _{I}^{\sigma ,%
\mathbf{b}}f-\sum_{I\in \pi \left( \Pi _{1}\mathcal{Q}\right) }\square _{I,%
\limfunc{broken}}^{\sigma ,\pi ,\mathbf{b}}f\right\Vert _{L^{2}\left( \sigma
\right) }^{2}  \label{needed for unfix} \\
&\lesssim &\left\Vert \sum_{I\in \pi \left( \Pi _{1}\mathcal{Q}\right)
}\square _{I}^{\sigma ,\mathbf{b}}f\right\Vert _{L^{2}\left( \sigma \right)
}^{2}+\left\Vert \sum_{I\in \pi \left( \Pi _{1}\mathcal{Q}\right) }\square
_{I,\limfunc{broken}}^{\sigma ,\pi ,\mathbf{b}}f\right\Vert _{L^{2}\left(
\sigma \right) }^{2}  \notag \\
&\lesssim &\left\Vert \mathsf{P}_{\pi \left( \Pi _{1}\mathcal{Q}\right)
}^{\sigma ,\mathbf{b}}f\right\Vert _{L^{2}\left( \sigma \right)
}^{2}+\sum_{I\in \pi \left( \Pi _{1}\mathcal{Q}\right) }\left\Vert \square
_{I,\limfunc{broken}}^{\sigma ,\pi ,\mathbf{b}}f\right\Vert _{L^{2}\left(
\sigma \right) }^{2}  \notag \\
&\lesssim &\sum_{I\in \pi \left( \Pi _{1}\mathcal{Q}\right) }\left\Vert
\square _{I}^{\sigma ,\mathbf{b}}f\right\Vert _{L^{2}\left( \sigma \right)
}^{2}+\sum_{I\in \pi \left( \Pi _{1}\mathcal{Q}\right) }\left\Vert
\bigtriangledown _{I}^{\sigma }f\right\Vert _{L^{2}\left( \sigma \right)
}^{2}\lesssim \left\Vert \mathsf{P}_{\pi \left( \Pi _{1}\mathcal{Q}\right)
}^{\sigma ,\mathbf{b}}f\right\Vert _{L^{2}\left( \sigma \right) }^{\bigstar
2}\ .  \notag
\end{eqnarray}

\medskip

We now use the fact that the supremum in the definition of $\mathcal{S}_{%
\limfunc{loc}\limfunc{size}}^{\alpha ,A;S}\left( \mathcal{Q}\right) $ is
taken over $K\in \mathcal{W}^{\ast }\left( S\right) \cap \mathcal{C}_{A}^{%
\limfunc{restrict}}$ to conclude that 
\begin{equation*}
\sup_{S\in \mathcal{S}}\mathcal{S}_{\limfunc{loc}\limfunc{size}}^{\alpha
,A;S}\left( \mathcal{Q}\right) \leq \mathcal{S}_{\limfunc{aug}\limfunc{size}%
}^{\alpha ,A}\left( \mathcal{Q}\right) ,
\end{equation*}%
and this completes the proof of Lemma \ref{straddle 3 ref}.
\end{proof}

In a similar fashion we can obtain the following Substraddling Lemma.

\begin{definition}
\label{def substraddles}Given a \emph{reduced admissible} collection of
pairs $\mathcal{Q}$ for $A$, and a $\mathcal{D}$-interval $L$ contained in $%
A $, we say that $\mathcal{Q}$ \textbf{substraddles} $L$ if for every pair $%
\left( I,J\right) \in \mathcal{Q}$ there is $K\in \mathcal{W}\left( L\right)
\cap \mathcal{C}_{A}^{\limfunc{restrict}}$ with $J\subset K\subset 3K\subset
I\subset L$.
\end{definition}

\begin{lemma}
\label{substraddle ref}Let $L$ be a $\mathcal{D}$-interval contained in $A$,
and suppose that $\mathcal{Q}$ is an admissible collection of pairs that
substraddles $L$. Then we have the sublinear form bound%
\begin{equation*}
\widehat{\mathfrak{N}}_{\limfunc{stop},\bigtriangleup ^{\omega }}^{A,%
\mathcal{Q}}\leq C\mathcal{S}_{\limfunc{aug}\limfunc{size}}^{\alpha
,A}\left( \mathcal{Q}\right) .
\end{equation*}
\end{lemma}

\begin{proof}
We will show that $\mathcal{Q}$ $\flat $straddles the subset $\mathcal{W}%
_{L} $ of Whitney intervals for $L$ given by%
\begin{equation*}
\mathcal{W}^{\mathcal{Q}}\left( L\right) \equiv \left\{ K\in \mathcal{W}%
\left( L\right) \cap \mathcal{C}_{A}^{\limfunc{restrict}}:J\subset K\subset
3K\subset I\subset L\text{ for some }\left( I,J\right) \in \mathcal{Q}%
\right\} .
\end{equation*}%
It is clear that $\mathcal{W}^{\mathcal{Q}}\left( L\right) \subset \Pi _{1}^{%
\limfunc{below}}\mathcal{Q}\cap \mathcal{C}_{A}^{\limfunc{restrict}}$ is a
subpartition of $A$. It remains to show that for every pair $\left(
I,J\right) \in \mathcal{Q}$ there is $K\in \mathcal{W}^{\mathcal{Q}}\left(
L\right) \cap \left[ J,I\right] $ such that $J^{\flat }\subset K$. But our
hypothesis implies that there is $K\in \mathcal{W}^{\mathcal{Q}}\left(
L\right) $ with $J\subset K\subset 3K\subset I\subset L$. We now consider
two cases.

\textbf{Case 1}: If $\pi _{\mathcal{D}}^{\left( 3\right) }K\subset L$, then
by Key Fact \#2 in (\ref{indentation}), i.e. $3J$ is contained in an \emph{%
inner} grandchild of $J^{\maltese }$. But $K$ is contained in an \emph{outer}
grandchild of $\pi _{\mathcal{D}}^{\left( 3\right) }K$ since $\pi _{\mathcal{%
D}}^{\left( 1\right) }K$ shares an endpoint with $L$, and so then does $\pi
_{\mathcal{D}}^{\left( 3\right) }K$). We thus have $J^{\maltese }\subset \pi
_{\mathcal{D}}^{\left( 2\right) }K$, which implies that $J^{\flat }\subset K$%
.

\textbf{Case 2}: If $\pi _{\mathcal{D}}^{\left( 3\right) }K\varsupsetneqq L$%
, then $K\subset 3K\subset I\subset L$ implies that $I=L=\pi _{\mathcal{D}%
}^{\left( 2\right) }K$. Thus we have $J^{\maltese }\subset I=\pi _{\mathcal{D%
}}^{\left( 2\right) }K$, which again gives $J^{\flat }\subset K$.

Now that we know $\mathcal{Q}$ $\flat $straddles the subset $\mathcal{W}^{%
\mathcal{Q}}\left( L\right) $, we can apply Lemma \ref{straddle 3 ref} to
obtain the required bound $\widehat{\mathfrak{N}}_{\limfunc{stop}%
,\bigtriangleup ^{\omega }}^{A,\mathcal{Q}}\leq C\mathcal{S}_{\limfunc{aug}%
\limfunc{size}}^{\alpha ,A}\left( \mathcal{Q}\right) $.
\end{proof}

\subsection{The bottom/up stopping time argument of M. Lacey}

Before introducing Lacey's stopping times, we note that the
Corona-straddling Lemma \ref{cor strad 1} allows us to remove the `corona
straddling' collection $\mathcal{P}_{\func{cor}}^{A}$ of pairs of intervals
in (\ref{def cor}) from the collection $\mathcal{P}^{A}$ in (\ref{initial P}%
) used to define the stopping form $\mathsf{B}_{\limfunc{stop}}^{A}\left(
f,g\right) $. The collection $\mathcal{P}^{A}\setminus \mathcal{P}_{\func{cor%
}}^{A}$ is of course also $A$-admissible.

\begin{conclusion}
\label{assume}We assume for the remainder of the proof that all admissible
collections $\mathcal{P}$ are reduced, i.e. 
\begin{equation}
\mathcal{P}^{A}\cap \mathcal{P}_{\func{cor}}^{A}=\emptyset ,\text{ as well
as }\mathcal{P}\cap \mathcal{P}_{\func{cor}}^{A}=\emptyset \text{ for all }A%
\text{-admissible }\mathcal{P}.  \label{empty assumption}
\end{equation}
\end{conclusion}

We remind the reader again that we will generally use $\left\vert
J\right\vert $ in the Poisson integrals and estimates, but will usually use $%
\ell \left( J\right) $ when defining collections of intervals. For an
interval $K\in \mathcal{D}$, we define%
\begin{equation*}
\mathcal{G}\left[ K\right] \equiv \left\{ J\in \mathcal{G}:J\subset K\right\}
\end{equation*}%
to consist of all intervals $J$ in the other grid $\mathcal{G}$ that are
contained in $K$. For an $A$-admissible collection $\mathcal{P}$ of pairs,
define two atomic measures $\omega _{\mathcal{P}}$ and $\omega _{\flat 
\mathcal{P}}$ in the upper half space $\mathbb{R}_{+}^{2}$ by%
\begin{equation}
\omega _{\mathcal{P}}\equiv \sum_{J\in \Pi _{2}\mathcal{P}}\left\Vert
\bigtriangleup _{J}^{\omega ,\mathbf{b}^{\ast }}x\right\Vert _{L^{2}\left(
\omega \right) }^{\spadesuit 2}\ \delta _{\left( c_{J^{\maltese }},\ell
\left( J^{\maltese }\right) \right) }\text{ and }\omega _{\flat \mathcal{P}%
}\equiv \sum_{J\in \Pi _{2}\mathcal{P}}\left\Vert \bigtriangleup
_{J}^{\omega ,\mathbf{b}^{\ast }}x\right\Vert _{L^{2}\left( \omega \right)
}^{\spadesuit 2}\ \delta _{\left( c_{J^{\flat }},\ell \left( J^{\flat
}\right) \right) },  \label{def atomic}
\end{equation}%
where $J^{\flat }=J_{\searrow J}^{\maltese }$ is the inner grandchild of $%
J^{\maltese }$ that contains $J$, i.e. $J^{\flat }=\left\{ 
\begin{array}{ccc}
J_{+/-}^{\maltese } & \text{ if } & J\subset J_{+/-}^{\maltese } \\ 
J_{-/+}^{\maltese } & \text{ if } & J\subset J_{-/+}^{\maltese }%
\end{array}%
\right. $. Note that each interval $J\in \Pi _{2}\mathcal{P}$ has its
`energy' $\left\Vert \bigtriangleup _{J}^{\omega ,\mathbf{b}^{\ast
}}x\right\Vert _{L^{2}\left( \omega \right) }^{\spadesuit 2}$ in the measure 
$\omega _{\flat \mathcal{P}}$ assigned to exactly one of the two points $%
\left( c_{J_{-/+}^{\maltese }},\frac{1}{4}\ell \left( J^{\maltese }\right)
\right) $ and $\left( c_{J_{+/-}^{\maltese }},\frac{1}{4}\ell \left(
J^{\maltese }\right) \right) $ in the upper half plane $\mathbb{R}_{+}^{2}$
since $J$ is either contained in $J_{-/+}^{\maltese }$\ or in $%
J_{+/-}^{\maltese }$ by Key Fact \#2 in (\ref{indentation}). Note also that
the atomic measure $\omega _{\flat \mathcal{P}}$ differs from the measure $%
\mu $ in (\ref{def mu n}) in Appendix B below - which is used there to
control the functional energy condition - in that here we bundle together
all the $J^{\prime }s$ having a common $J^{\flat }$. This is in order to
rewrite the \emph{augmented} size functional in terms of the measure $\omega
_{\flat \mathcal{P}}$. We can get away with this here, as opposed to in
Appendix B, due to the `smaller and decoupled' nature of the augmented size
functional to which we will relate $\omega _{\flat \mathcal{P}}$.

Define the tent $\mathbf{T}\left( L\right) $ over an interval $L$ to be the
convex hull of the interval $L\times \left\{ 0\right\} $ and the point $%
\left( c_{L},\ell \left( L\right) \right) \in \mathbb{R}_{+}^{2}$. Then for $%
J\in \Pi _{2}\mathcal{P}$ we have $J\in \Pi _{2}^{K,\limfunc{aug}}\mathcal{P}
$ \emph{iff} $\left\{ J\subset K\text{ and }J^{\maltese }\subset \pi _{%
\mathcal{D}}^{\left( 2\right) }K\right\} $ \emph{iff} $J^{\flat
}=J_{\searrow J}^{\maltese }\subset K$ \emph{iff} $\left( c_{J^{\flat
}},\ell \left( J^{\flat }\right) \right) \in \mathbf{T}\left( K\right) $. We
can now rewrite the augmented size functional of $\mathcal{P}$ in Definition %
\ref{augs}\ as%
\begin{equation}
\mathcal{S}_{\limfunc{aug}\limfunc{size}}^{\alpha ,A}\left( \mathcal{P}%
\right) ^{2}\equiv \sup_{K\in \Pi _{1}^{\limfunc{below}}\mathcal{P}\cap 
\mathcal{C}_{A}^{\limfunc{restrict}}}\frac{1}{\left\vert K\right\vert
_{\sigma }}\left( \frac{\mathrm{P}^{\alpha }\left( K,\mathbf{1}_{A\setminus
K}\sigma \right) }{\left\vert K\right\vert }\right) ^{2}\omega _{\flat 
\mathcal{P}}\left( \mathbf{T}\left( K\right) \right) .
\label{def P stop energy' 3}
\end{equation}%
It will be convenient to write%
\begin{equation*}
\Psi ^{\alpha }\left( K;\mathcal{P}\right) ^{2}\equiv \left( \frac{\mathrm{P}%
^{\alpha }\left( K,\mathbf{1}_{A\setminus K}\sigma \right) }{\left\vert
K\right\vert }\right) ^{2}\omega _{\flat \mathcal{P}}\left( \mathbf{T}\left(
K\right) \right) ,
\end{equation*}%
so that we have simply%
\begin{equation*}
\mathcal{S}_{\limfunc{aug}\limfunc{size}}^{\alpha ,A}\left( \mathcal{P}%
\right) ^{2}=\sup_{K\in \Pi _{1}^{\limfunc{below}}\mathcal{P}\cap \mathcal{C}%
_{A}^{\limfunc{restrict}}}\frac{\Psi ^{\alpha }\left( K;\mathcal{P}\right)
^{2}}{\left\vert K\right\vert _{\sigma }}.
\end{equation*}

\begin{remark}
The functional $\omega _{\flat \mathcal{P}}\left( \mathbf{T}\left( K\right)
\right) $ is increasing in $K$, while the functional $\frac{\mathrm{P}%
^{\alpha }\left( K,\mathbf{1}_{A\setminus K}\sigma \right) }{\left\vert
K\right\vert }$ is `almost decreasing' in $K$: if $K_{0}\subset K$ then%
\begin{eqnarray*}
\frac{\mathrm{P}^{\alpha }\left( K,\mathbf{1}_{A\setminus K}\sigma \right) }{%
\left\vert K\right\vert } &=&\int_{A\setminus K}\frac{d\sigma \left(
y\right) }{\left( \left\vert K\right\vert +\left\vert y-c_{K}\right\vert
\right) ^{2-\alpha }} \\
&\lesssim &\int_{A\setminus K}\frac{d\sigma \left( y\right) }{\left(
\left\vert K_{0}\right\vert +\left\vert y-c_{K_{0}}\right\vert \right)
^{2-\alpha }} \\
&\leq &C_{\alpha }\int_{A\setminus K_{0}}\frac{d\sigma \left( y\right) }{%
\left( \left\vert K_{0}\right\vert +\left\vert y-c_{K_{0}}\right\vert
\right) ^{2-\alpha }}=C_{\alpha }\frac{\mathrm{P}^{\alpha }\left( K_{0},%
\mathbf{1}_{A\setminus K_{0}}\sigma \right) }{\left\vert K_{0}\right\vert },
\end{eqnarray*}%
since $\left\vert K_{0}\right\vert +\left\vert y-c_{K_{0}}\right\vert \leq
\left\vert K\right\vert +\left\vert y-c_{K}\right\vert +\frac{1}{2}\limfunc{%
diam}\left( K\right) $ for $y\in A\setminus K$.
\end{remark}

Recall that if $\mathcal{P}$ is an admissible collection for a dyadic
interval $A$, the corresponding sublinear form in (\ref{def mod B}) and (\ref%
{First inequality}) is given by%
\begin{eqnarray*}
\left\vert \mathsf{B}\right\vert _{\limfunc{stop},\bigtriangleup ^{\omega
}}^{A,\mathcal{P}}\left( f,g\right) &\equiv &\sum_{J\in \Pi _{2}\mathcal{P}}%
\frac{\mathrm{P}^{\alpha }\left( J,\left\vert \varphi _{J}^{\mathcal{P}%
}\right\vert \mathbf{1}_{A\setminus I_{\mathcal{P}}\left( J\right) }\sigma
\right) }{\left\vert J\right\vert }\left\Vert \bigtriangleup _{J}^{\omega ,%
\mathbf{b}^{\ast }}x\right\Vert _{L^{2}\left( \omega \right) }^{\spadesuit
}\left\Vert \square _{J}^{\omega ,\mathbf{b}^{\ast }}g\right\Vert
_{L^{2}\left( \omega \right) }^{\bigstar }; \\
\text{where }\varphi _{J}^{\mathcal{P}} &\equiv &\sum_{I\in \mathcal{C}_{A}^{%
\limfunc{restrict}}:\ \left( I,J\right) \in \mathcal{P}}b_{A}E_{I}^{\sigma
}\left( \widehat{\square }_{\pi I}^{\sigma ,\flat ,\mathbf{b}}f\right) \ 
\mathbf{1}_{A\setminus I}\ .
\end{eqnarray*}%
In the notation for $\left\vert \mathsf{B}\right\vert _{\limfunc{stop}%
,\bigtriangleup ^{\omega }}^{A,\mathcal{P}}$, we are omitting dependence on
the parameter $\alpha $, and to avoid clutter, we will often do so from now
on when the dependence on $\alpha $ is inconsequential.

Recall further that the `size testing collection' of intervals $\Pi _{1}^{%
\limfunc{below}}\mathcal{P}$ for the initial size testing functional $%
\mathcal{S}_{\limfunc{init}\limfunc{size}}^{\alpha ,A}\left( \mathcal{P}%
\right) $ is the collection of all subintervals of intervals in $\Pi _{1}%
\mathcal{P}$, and moreover, by Key Fact \#1 in (\ref{later use}), that we
can restrict the collection to $\Pi _{1}^{\limfunc{below}}\mathcal{P}\cap 
\mathcal{C}_{A}^{\limfunc{restrict}}$. This latter set is used for the
augmented size functional.

\begin{description}
\item[Assumption] We may assume that the corona $\mathcal{C}_{A}$ is finite,
and that each $A$-admissible collection $\mathcal{P}$ is a finite
collection, and hence so are $\Pi _{1}\mathcal{P}$, $\Pi _{1}^{\limfunc{below%
}}\mathcal{P}\cap \mathcal{C}_{A}^{\limfunc{restrict}}$ and $\Pi _{2}%
\mathcal{P}$, provided all of the bounds we obtain are independent of the
cardinality of these latter collections.
\end{description}

Consider $0<\varepsilon <1$, where $\rho =1+\varepsilon $ will be chosen
later in (\ref{choose rho}). Begin by defining the collection $\mathcal{L}%
_{0}$ to consist of the \emph{minimal} dyadic intervals $K$ in $\Pi _{1}^{%
\limfunc{below}}\mathcal{P}\cap \mathcal{C}_{A}^{\limfunc{restrict}}$ such
that%
\begin{equation*}
\frac{\Psi ^{\alpha }\left( K;\mathcal{P}\right) ^{2}}{\left\vert
K\right\vert _{\sigma }}\geq \varepsilon \mathcal{S}_{\limfunc{aug}\limfunc{%
size}}^{\alpha ,A}\left( \mathcal{P}\right) ^{2}.
\end{equation*}%
where we recall that%
\begin{equation*}
\Psi ^{\alpha }\left( K;\mathcal{P}\right) ^{2}\equiv \left( \frac{\mathrm{P}%
^{\alpha }\left( K,\mathbf{1}_{A\setminus K}\sigma \right) }{\left\vert
K\right\vert }\right) ^{2}\omega _{\flat \mathcal{P}}\left( \mathbf{T}\left(
K\right) \right) .
\end{equation*}%
Note that such minimal intervals exist when $0<\varepsilon <1$ because $%
\mathcal{S}_{\limfunc{aug}\limfunc{size}}^{\alpha ,A}\left( \mathcal{P}%
\right) ^{2}$ is the supremum over $K\in \Pi _{1}^{\limfunc{below}}\mathcal{P%
}\cap \mathcal{C}_{A}^{\limfunc{restrict}}$ of $\frac{\Psi ^{\alpha }\left(
K;\mathcal{P}\right) ^{2}}{\left\vert K\right\vert _{\sigma }}$. A key
property of the minimality requirement is that%
\begin{equation}
\frac{\Psi ^{\alpha }\left( K^{\prime };\mathcal{P}\right) ^{2}}{\left\vert
K^{\prime }\right\vert _{\sigma }}<\varepsilon \mathcal{S}_{\limfunc{aug}%
\limfunc{size}}^{\alpha ,A}\left( \mathcal{P}\right) ^{2},
\label{key property 3}
\end{equation}%
whenever there is $K^{\prime }\in \Pi _{1}^{\limfunc{below}}\mathcal{P}\cap 
\mathcal{C}_{A}^{\limfunc{restrict}}$ with $K^{\prime }\varsubsetneqq K$ and 
$K\in \mathcal{L}_{0}$.

We now perform a stopping time argument `from the bottom up' with respect to
the atomic measure $\omega _{\mathcal{P}}$ in the upper half space. This
construction of a stopping time `from the bottom up', together with the
subsequent applications of the Orthogonality Lemma and the Straddling Lemma,
comprise the key innovations in Lacey's argument \cite{Lac}. However, in our
situation the intervals $I$ belonging to $\Pi _{1}^{\limfunc{below}}\mathcal{%
P}$ are no longer `good' in any sense, and we must include an additional
top/down stopping criterion in the next subsection to accommodate this lack
of `goodness'. The argument in \cite{Lac} will apply to these special
stopping intervals, called `indented' intervals, and the remaining intervals
form towers with a common endpoint, that are controlled using all three
straddling lemmas.

We refer to $\mathcal{L}_{0}$ as the initial or level $0$ generation of
stopping intervals. Set%
\begin{equation}
\rho =1+\varepsilon .  \label{def rho}
\end{equation}%
As in \cite{SaShUr7}, \cite{SaShUr9} and \cite{SaShUr10}, we follow Lacey 
\cite{Lac} by recursively defining a finite sequence of generations $\left\{ 
\mathcal{L}_{m}\right\} _{m\geq 0}$ by letting $\mathcal{L}_{m}$ consist of
the \emph{minimal} dyadic intervals $L$ in $\Pi _{1}^{\limfunc{below}}%
\mathcal{P}\cap \mathcal{C}_{A}^{\limfunc{restrict}}$ that contain an
interval from some previous level $\mathcal{L}_{\ell }$, $\ell <m$, such that%
\begin{equation}
\omega _{\flat \mathcal{P}}\left( \mathbf{T}\left( L\right) \right) \geq
\rho \omega _{\flat \mathcal{P}}\left( \dbigcup\limits_{L^{\prime }\in
\dbigcup\limits_{\ell =0}^{m-1}\mathcal{L}_{\ell }:\ L^{\prime }\subset L}%
\mathbf{T}\left( L^{\prime }\right) \right) .  \label{up stopping condition}
\end{equation}%
Since $\mathcal{P}$ is finite this recursion stops at some level $M$. We
then let $\mathcal{L}_{M+1}$ consist of all the maximal intervals in $\Pi
_{1}^{\limfunc{below}}\mathcal{P}\cap \mathcal{C}_{A}^{\limfunc{restrict}}$
that are not already in some $\mathcal{L}_{m}$ with $m\leq M$. Thus $%
\mathcal{L}_{M+1}$ will contain either none, some, or all of the maximal
intervals in $\Pi _{1}^{\limfunc{below}}\mathcal{P}$. We do not of course
have (\ref{up stopping condition}) for $A^{\prime }\in \mathcal{L}_{M+1}$ in
this case, but we do have that (\ref{up stopping condition}) fails for
subintervals $K$ of $A^{\prime }\in \mathcal{L}_{M+1}$ that are not
contained in any other $L\in \mathcal{L}_{m}$ with $m\leq M$, and this is
sufficient for the arguments below.

We now decompose the collection of pairs $\left( I,J\right) $ in $\mathcal{P}
$ into collections $\mathcal{P}^{\flat small}$ and $\mathcal{P}^{\flat big}$
according to the location of $I$ and $J^{\flat }$, but only after
introducing below the indented corona $\mathcal{H}$. The collection $%
\mathcal{P}^{\flat big}$ will then essentially consist of those pairs $%
\left( I,J\right) \in \mathcal{P}$ for which there are $L^{\prime },L\in 
\mathcal{H}$ with $L^{\prime }\varsubsetneqq L$ and such that $J^{\flat }\in 
\mathcal{C}_{L^{\prime }}^{\mathcal{H}}$ and $I\in \mathcal{C}_{L}^{\mathcal{%
H}}$. The collection $\mathcal{P}^{\flat small}$ will consist of the
remaining pairs $\left( I,J\right) \in \mathcal{P}$ for which there is $L\in 
\mathcal{H}$ such that $J^{\flat },I\in \mathcal{C}_{L}^{\mathcal{H}}$,
along with the pairs $\left( I,J\right) \in \mathcal{P}$ such that $I\subset
I_{0}$ for some $I_{0}\in \mathcal{L}_{0}$. This will cover all pairs $%
\left( I,J\right) $ in $\mathcal{P}\subset \mathcal{P}_{A}$, since for such
pairs, $I\in \mathcal{C}_{A}^{\limfunc{restrict}}$\ and $J\in \mathcal{C}%
_{A}^{\mathcal{G}\func{shift}}$, which in turn implies $I\in \mathcal{C}%
_{L}^{\mathcal{H}}$ and $J^{\flat }\in \mathcal{C}_{L^{\prime }}^{\mathcal{H}%
}$ for some $L,L^{\prime }\in \mathcal{H}$. But a considerable amount of
further analysis is required to prove (\ref{First inequality}).

First recall that $\mathcal{L}\equiv \dbigcup\limits_{m=0}^{M+1}\mathcal{L}%
_{m}$ is the tree of stopping $\omega _{\mathcal{P}}$-energy intervals
defined above. By the construction above, the maximal elements in $\mathcal{L%
}$ are the maximal intervals in $\Pi _{1}^{\limfunc{below}}\mathcal{P}\cap 
\mathcal{C}_{A}^{\limfunc{restrict}}$. For $L\in \mathcal{L}$, denote by $%
\mathcal{C}_{L}^{\mathcal{L}}$ the \emph{corona} associated with $L$ in the
tree $\mathcal{L}$,%
\begin{equation*}
\mathcal{C}_{L}^{\mathcal{L}}\equiv \left\{ K\in \mathcal{D}:K\subset L\text{
and there is no }L^{\prime }\in \mathcal{L}\text{ with }K\subset L^{\prime
}\subsetneqq L\right\} ,
\end{equation*}%
and define the \emph{shifted} $\mathcal{L}$-corona and the $\flat $\emph{%
shifted} $\mathcal{L}$-corona by%
\begin{eqnarray*}
\mathcal{C}_{L}^{\mathcal{L},\limfunc{shift}} &\equiv &\left\{ J\in \mathcal{%
G}:J^{\maltese }\in \mathcal{C}_{L}^{\mathcal{L}}\text{ }\right\} , \\
\mathcal{C}_{L}^{\mathcal{L},\flat \limfunc{shift}} &\equiv &\left\{ J\in 
\mathcal{G}:J^{\flat }\in \mathcal{C}_{L}^{\mathcal{L}}\text{ }\right\} .
\end{eqnarray*}%
It is the second flat shifted corona $\mathcal{C}_{L}^{\mathcal{L},\flat 
\limfunc{shift}}$ that will be used in our decompositions below, but we
retain the more natural shifted coronas $\mathcal{C}_{L}^{\mathcal{L},%
\limfunc{shift}}$ in order to make useful comparisons. Now the parameter $m$
in $\mathcal{L}_{m}$ refers to the level at which the stopping construction
was performed, but for\thinspace $L\in \mathcal{L}_{m}$, the corona children 
$L^{\prime }$ of $L$ are \emph{not} all necessarily in $\mathcal{L}_{m-1}$,
but may be in $\mathcal{L}_{m-t}$ for $t$ large.

At this point we introduce the notion of geometric depth $d$ in the tree $%
\mathcal{L}$ by defining%
\begin{eqnarray}
\mathcal{G}_{0} &\equiv &\left\{ L\in \mathcal{L}:L\text{ is maximal}%
\right\} ,  \label{geom depth} \\
\mathcal{G}_{1} &\equiv &\left\{ L\in \mathcal{L}:L\text{ is maximal wrt }%
L\subsetneqq L_{0}\text{ for some }L_{0}\in \mathcal{G}_{0}\right\} ,  \notag
\\
&&\vdots  \notag \\
\mathcal{G}_{d+1} &\equiv &\left\{ L\in \mathcal{L}:L\text{ is maximal wrt }%
L\subsetneqq L_{d}\text{ for some }L_{d}\in \mathcal{G}_{d}\right\} ,  \notag
\\
&&\vdots  \notag
\end{eqnarray}%
We refer to $\mathcal{G}_{d}$ as the $d^{th}$ generation of intervals in the
tree $\mathcal{L}$, and say that the intervals in $\mathcal{G}_{d}$ are at
depth $d$ in the tree $\mathcal{L}$ (the generations $\mathcal{G}_{d}$ here
are \emph{not} related to the grid $\mathcal{G}$), and we write $d_{\limfunc{%
geom}}\left( L\right) $ for the geometric depth of $L$. Thus the intervals
in $\mathcal{G}_{d}$ are the stopping intervals in $\mathcal{L}$ that are $d$
levels in the \emph{geometric} sense below the top level. While the
geometric depth $d_{\limfunc{geom}}$ is about to be superceded by the
`indented' depth $d_{\limfunc{indent}}$ defined in the next subsection, we
will return to the geometric depth in order to iterate Lacey's bottom/up
stopping criterion when proving the second line in (\ref{rest bounds}) in
Proposition \ref{bottom up 3} below.

\subsection{The indented corona construction}

Now we address the lack of goodness in $\Pi _{1}^{\limfunc{below}}\mathcal{P}%
\cap \mathcal{C}_{A}^{\limfunc{restrict}}$. For this we introduce an
additional top/down stopping time $\mathcal{H}$ over the collection $%
\mathcal{L}$. Given the initial generation 
\begin{equation*}
\mathcal{H}_{0}\equiv \mathcal{L}_{M+1}=\left\{ \text{maximal }L\in \mathcal{%
L}\right\} =\left\{ \text{maximal }I\in \Pi _{1}^{\limfunc{below}}\mathcal{P}%
\right\} ,
\end{equation*}%
define subsequent generations $\mathcal{H}_{k}$ as follows. For $k\geq 1$
and each $H\in \mathcal{H}_{k-1}$, let 
\begin{equation*}
\mathcal{H}_{k}\left( H\right) \equiv \left\{ \text{maximal }L\in \mathcal{L}%
:3L\subset H\right\}
\end{equation*}%
consist of the next $\mathcal{H}$-generation of $\mathcal{L}$-intervals
below $H$, and set $\mathcal{H}_{k}\equiv \dbigcup\limits_{H\in \mathcal{H}%
_{k-1}}\mathcal{H}_{k}\left( H\right) $. Finally set $\mathcal{H}\equiv
\dbigcup\limits_{k=0}^{\infty }\mathcal{H}_{k}$. We refer to the stopping
intervals $H\in \mathcal{H}$ as \emph{indented} stopping intervals since $%
3H\subset \pi _{\mathcal{H}}H$ for all $H\in \mathcal{H}$ at indented
generation one or more, i.e. each successive such $H$ is `indented' in its $%
\mathcal{H}$-parent. This property of indentation is precisely what is
required in order to generate geometric decay in indented generations at the
end of the proof. We refer to $k$ as the \emph{indented depth} of the
stopping interval $H\in \mathcal{H}_{k}$, written $k=d_{\limfunc{indent}%
}\left( H\right) $, which is a refinement of the geometric depth $d_{%
\limfunc{geom}}$ introduced above. We will often revert to writing the dummy
variable for intervals in $\mathcal{H}$ as $L$ instead of $H$. For $L\in 
\mathcal{H}$ define the $\mathcal{H}$-corona $\mathcal{C}_{L}^{\mathcal{H},%
\limfunc{shift}}$ and the $\mathcal{H}$-shifted corona $\mathcal{C}_{L}^{%
\mathcal{H},\limfunc{shift}}$ and $\mathcal{H}$-$\flat $shifted corona $%
\mathcal{C}_{L}^{\mathcal{H},\flat \limfunc{shift}}$ by%
\begin{eqnarray*}
\mathcal{C}_{L}^{\mathcal{H}} &\equiv &\left\{ I\in \mathcal{D}:I\subset L%
\text{ and }I\not\subset L^{\prime }\text{ for any }L^{\prime }\in \mathfrak{%
C}_{\mathcal{H}}\left( L\right) \right\} , \\
\mathcal{C}_{L}^{\mathcal{H},\limfunc{shift}} &\equiv &\left\{ J\in \mathcal{%
G}:J^{\maltese }\in \mathcal{C}_{L}^{\mathcal{H}}\right\} , \\
\mathcal{C}_{L}^{\mathcal{H},\flat \limfunc{shift}} &\equiv &\left\{ J\in 
\mathcal{G}:J^{\flat }\in \mathcal{C}_{L}^{\mathcal{H}}\right\} .
\end{eqnarray*}%
We will also need recourse to the coronas $\mathcal{C}_{L}^{\mathcal{H}}$
restricted to intervals in $\mathcal{L}$, i.e.%
\begin{equation*}
\mathcal{C}_{L}^{\mathcal{H}}\left( \mathcal{L}\right) \equiv \mathcal{C}%
_{L}^{\mathcal{H}}\cap \mathcal{L}=\left\{ T\in \mathcal{L}:T\subset L\text{
and }T\not\subset L^{\prime }\text{ for any }L^{\prime }\in \mathcal{H}\text{
with }L^{\prime }\subsetneqq L\right\} .
\end{equation*}

Then for $L\in \mathcal{H}$ and $t\geq 0$ define%
\begin{equation}
\mathcal{P}_{L,t}^{\mathcal{H}}\equiv \left\{ \left( I,J\right) \in \mathcal{%
P}:I\in \mathcal{C}_{L}^{\mathcal{H}}\text{ and }J\in \mathcal{C}_{L^{\prime
}}^{\mathcal{H},\limfunc{shift}}\text{ for some }L^{\prime }\in \mathcal{H}%
_{d_{\limfunc{indent}}\left( L\right) +t}\text{ with }L^{\prime }\subset
L\right\} .  \label{def PHLt}
\end{equation}%
In particular, $\left( I,J\right) \in \mathcal{P}_{L,t}^{\mathcal{H}}$
implies that $I$ is in the corona $\mathcal{C}_{L}^{\mathcal{H}}$, and that $%
J$ is in a shifted corona $\mathcal{C}_{L^{\prime }}^{\mathcal{H},\limfunc{%
shift}}$ that is $t$ levels of indented generation \emph{below} $\mathcal{C}%
_{L}^{\mathcal{H}}$ (when $t=0$ we have $L^{\prime }=L$). We emphasize the
distinction `indented generation' as this refers to the indented depth
rather than either the level of initial stopping construction of $\mathcal{L}
$, or the geometric depth. The point of introducing the tree $\mathcal{H}$
of indented stopping intervals, is that the inclusion $3L\subset \pi _{%
\mathcal{H}}L$ for all $L\in \mathcal{H}$ with $d_{\limfunc{indent}}\left(
L\right) \geq 1$ turns out to be an adequate substitute for the standard
`goodness' lost in the process of infusing the weak goodness of Hyt\"{o}nen
and Martikainen in\ Subsection \ref{Subsec HM} above.

Now within the $\mathcal{H}$-corona $\mathcal{C}_{L}^{\mathcal{H}}\left( 
\mathcal{L}\right) $, there are in general further intervals $T\in \mathcal{L%
}$ in addition to $L\in \mathcal{H}$ itself, but all of these further
intervals are contained in the two endpoint $\mathcal{L}$-towers%
\begin{eqnarray*}
\mathcal{T}_{\limfunc{left}}\left( L\right) &\equiv &\left\{ L^{\prime }\in 
\mathcal{L}:L^{\prime }\subsetneqq L\text{ and }\limfunc{left}\func{end}%
\left( L^{\prime }\right) =\limfunc{left}\func{end}\left( L\right) \right\}
\\
\mathcal{T}_{\limfunc{right}}\left( L\right) &\equiv &\left\{ L^{\prime }\in 
\mathcal{L}:L^{\prime }\subsetneqq L\text{ and }\limfunc{right}\func{end}%
\left( L^{\prime }\right) =\limfunc{right}\func{end}\left( L\right) \right\}
\end{eqnarray*}%
where $\limfunc{left}\func{end}\left( I\right) $ and $\limfunc{right}\func{%
end}\left( I\right) $ denote the left and right hand endpoints of $I$
respectively. Thus $\mathcal{C}_{L}^{\mathcal{H},\limfunc{restrict}}\left( 
\mathcal{L}\right) \equiv \mathcal{C}_{L}^{\mathcal{H}}\left( \mathcal{L}%
\right) \setminus \left\{ L\right\} $ consists of two `connected' $\mathcal{L%
}$-towers (possibly one or both empty), one in $\mathcal{T}_{\limfunc{left}%
}\left( L\right) $ and the other in $\mathcal{T}_{\limfunc{right}}\left(
L\right) $. Set $\mathcal{T}\left( L\right) \equiv \mathcal{T}_{\limfunc{left%
}}\left( L\right) \dot{\cup}\mathcal{T}_{\limfunc{right}}\left( L\right) 
\dot{\cup}\left\{ L\right\} $. See Figure \ref{ind}.

\FRAME{ftbpFU}{6.8393in}{3.4006in}{0pt}{\Qcb{Line segments (not to scale)
are the bottom/up stopping intervals in Lacey's tree $\mathcal{L}$. Red
segments are the \emph{indented} intervals and green segments are the
intervals in \emph{towers}. The top indented interval is boxed in purple,
the first generation of indented intervals are boxed in orange, and the
second generation in blue. Vertical lines indicate common endpoints.}}{\Qlb{%
ind}}{indented.wmf}{\special{language "Scientific Word";type
"GRAPHIC";maintain-aspect-ratio TRUE;display "USEDEF";valid_file "F";width
6.8393in;height 3.4006in;depth 0pt;original-width 7.1805in;original-height
10.469in;cropleft "0.0901";croptop "0.6533";cropright "1";cropbottom
"0.3462";filename 'Indented.wmf';file-properties "XNPEU";}}

\subsubsection{Decomposition of coronas}

Here we describe the decomposition of admissible collections of pairs
according to the regular shifted coronas $\mathcal{C}_{L}^{\mathcal{H},%
\limfunc{shift}}$ and $\mathcal{C}_{L}^{\mathcal{L},\limfunc{shift}}$.
Strictly speaking, these decompositions will not be used in the sequel, but
they do help to provide insight via comparison with the flat shifted
decompositions introduced in the next subsubsection, which will be used to
finish the proof. For $L\in \mathcal{H}$ and $t=0$ we further decompose $%
\mathcal{P}_{L,0}^{\mathcal{H}}$ in (\ref{def PHLt}) with $t=0$, i.e. 
\begin{equation*}
\mathcal{P}_{L,0}^{\mathcal{H}}=\left\{ \left( I,J\right) \in \mathcal{P}%
:I\in \mathcal{C}_{L}^{\mathcal{H}}\text{ and }J\in \mathcal{C}_{L}^{%
\mathcal{H},\limfunc{shift}}\right\} ,
\end{equation*}%
as%
\begin{eqnarray*}
\mathcal{P}_{L,0}^{\mathcal{H}} &=&\mathcal{P}_{L,0}^{\mathcal{H}-small}\dot{%
\cup}\mathcal{P}_{L,0}^{\mathcal{H}-big}; \\
\mathcal{P}_{L,0}^{\mathcal{H}-small} &\equiv &\left\{ \left( I,J\right) \in 
\mathcal{P}_{L,0}^{\mathcal{H}}:\text{there is no }L^{\prime }\in \mathcal{T}%
\left( L\right) \text{ with }J^{\maltese }\subset L^{\prime }\subset
I\right\} \\
&=&\left\{ \left( I,J\right) \in \mathcal{P}_{L,0}^{\mathcal{H}}:I\in 
\mathcal{C}_{L^{\prime }}^{\mathcal{L}}\setminus \left\{ L^{\prime }\right\} 
\text{ and }J\in \mathcal{C}_{L^{\prime }}^{\mathcal{L},\limfunc{shift}}%
\text{ for some }L^{\prime }\in \mathcal{T}\left( L\right) \right\} , \\
\mathcal{P}_{L,0}^{\mathcal{H}-big} &\equiv &\left\{ \left( I,J\right) \in 
\mathcal{P}_{L,0}^{\mathcal{H}}:\text{there is }L^{\prime }\in \mathcal{T}%
\left( L\right) \text{ with }J^{\maltese }\subset L^{\prime }\subset
I\right\} ,
\end{eqnarray*}%
with one exception: if $L\in \mathcal{H}_{0}=\mathcal{L}_{M+1}$ we set $%
\mathcal{P}_{L,0}^{\mathcal{H}-small}\equiv \mathcal{P}_{L,0}^{\mathcal{H}}$
and $\mathcal{P}_{L,0}^{\mathcal{H}-big}\equiv \emptyset $ since in this
case $L$ fails to satisfy (\ref{up stopping condition}) as pointed out
above. Finally, for $L\in \mathcal{H}$ we further decompose $\mathcal{P}%
_{L,0}^{\mathcal{H}-small}$ as%
\begin{eqnarray*}
\mathcal{P}_{L,0}^{\mathcal{H}-small} &=&\overset{\cdot }{\dbigcup }%
_{L^{\prime }\in \mathcal{T}\left( L\right) }\mathcal{P}_{L^{\prime },0}^{%
\mathcal{L}-small}; \\
\mathcal{P}_{L^{\prime },0}^{\mathcal{L}-small} &\equiv &\left\{ \left(
I,J\right) \in \mathcal{P}:I\in \mathcal{C}_{L^{\prime }}^{\mathcal{L}%
}\setminus \left\{ L^{\prime }\right\} \text{ and }J\in \mathcal{C}%
_{L^{\prime }}^{\mathcal{L},\limfunc{shift}}\right\} .
\end{eqnarray*}%
Then we set%
\begin{eqnarray}
\mathcal{P}^{big} &\equiv &\left\{ \dbigcup\limits_{L\in \mathcal{H}}%
\mathcal{P}_{L,0}^{\mathcal{H}-big}\right\} \dbigcup \left\{
\dbigcup\limits_{t\geq 1}\dbigcup\limits_{L\in \mathcal{H}}\mathcal{P}%
_{L,t}^{\mathcal{H}}\right\} ;  \label{def big small} \\
\mathcal{P}^{small} &\equiv &\dbigcup\limits_{L\in \mathcal{L}}\mathcal{P}%
_{L,0}^{\mathcal{L}-small}\text{ }.  \notag
\end{eqnarray}%
Note that every pair $\left( I,J\right) \in \mathcal{P}$ is included in
either $\mathcal{P}^{small}$ or $\mathcal{P}^{big}$ since every $I\in \Pi
_{1}^{\limfunc{below}}\mathcal{P}$ is contained in some $L\in \mathcal{H}%
_{0}=\mathcal{L}_{M+1}$.

\subsubsection{Flat shifted coronas}

More importantly, we now define the corresponding $\flat $shifted admissible
collections of pairs $\mathcal{P}_{L,t}^{\flat \mathcal{H}}$, etc., in which
we replace $\mathcal{C}_{L}^{\mathcal{H},\limfunc{shift}}$ and $\mathcal{C}%
_{L}^{\mathcal{L},\limfunc{shift}}$ with%
\begin{equation*}
\mathcal{C}_{L}^{\mathcal{H},\flat \limfunc{shift}}\equiv \left\{ J\in \Pi
_{2}\mathcal{P}:J^{\flat }\in \mathcal{C}_{L}^{\mathcal{H}}\right\} \text{
and }\mathcal{C}_{L}^{\mathcal{L},\flat \limfunc{shift}}\equiv \left\{ J\in
\Pi _{2}\mathcal{P}:J^{\flat }\in \mathcal{C}_{L}^{\mathcal{L}}\right\} .
\end{equation*}%
In these flat shifted $\mathcal{H}$ and $\mathcal{L}$ coronas, we have
effectively shift the intervals $J$ two levels `up' by requiring $J^{\flat
}\in \mathcal{C}_{L}^{\mathcal{L}}$ instead of $J^{\maltese }\in \mathcal{C}%
_{L}^{\mathcal{L}}$, etc., but because $\mathcal{P}$ is admissible, we
always have $J^{\maltese }\in \mathcal{C}_{A}^{\mathcal{A},\limfunc{restrict}%
}$. We define 
\begin{eqnarray*}
\mathcal{P}_{L,t}^{\flat \mathcal{H}} &\equiv &\left\{ \left( I,J\right) \in 
\mathcal{P}:I\in \mathcal{C}_{L}^{\mathcal{H}}\text{ and }J\in \mathcal{C}%
_{L^{\prime }}^{\mathcal{H},\flat \limfunc{shift}}\text{ for some }L^{\prime
}\in \mathcal{H}_{d_{\limfunc{indent}}\left( L\right) +t}\text{ with }%
L^{\prime }\subset L\right\} , \\
\mathcal{P}_{L,0}^{\flat \mathcal{H}} &=&\left\{ \left( I,J\right) \in 
\mathcal{P}:I\in \mathcal{C}_{L}^{\mathcal{H}}\text{ and }J\in \mathcal{C}%
_{L}^{\mathcal{H},\flat \limfunc{shift}}\right\} ,
\end{eqnarray*}%
and%
\begin{eqnarray*}
\mathcal{P}_{L,0}^{\flat \mathcal{H}} &=&\mathcal{P}_{L,0}^{\flat \mathcal{H}%
-small}\dot{\cup}\mathcal{P}_{L,0}^{\flat \mathcal{H}-big}; \\
\mathcal{P}_{L,0}^{\flat \mathcal{H}-small} &\equiv &\left\{ \left(
I,J\right) \in \mathcal{P}_{L,0}^{\flat \mathcal{H}}:\text{there is no }%
L^{\prime }\in \mathcal{T}\left( L\right) \text{ with }J^{\flat }\subset
L^{\prime }\subset I\right\} \\
&=&\left\{ \left( I,J\right) \in \mathcal{P}_{L,0}^{\flat \mathcal{H}}:I\in 
\mathcal{C}_{L^{\prime }}^{\mathcal{L}}\setminus \left\{ L^{\prime }\right\} 
\text{ and }J\in \mathcal{C}_{L^{\prime }}^{\mathcal{L},\flat \limfunc{shift}%
}\text{ for some }L^{\prime }\in \mathcal{T}\left( L\right) \right\} , \\
\mathcal{P}_{L,0}^{\flat \mathcal{H}-big} &\equiv &\left\{ \left( I,J\right)
\in \mathcal{P}_{L,0}^{\flat \mathcal{H}}:\text{there is }L^{\prime }\in 
\mathcal{T}\left( L\right) \text{ with }J^{\flat }\subset L^{\prime }\subset
I\right\} ,
\end{eqnarray*}%
with one exception: if $L\in \mathcal{H}_{0}=\mathcal{L}_{M+1}$ we set $%
\mathcal{P}_{L,0}^{\flat \mathcal{H}-small}\equiv \mathcal{P}_{L,0}^{\flat 
\mathcal{H}}$ and $\mathcal{P}_{L,0}^{\flat \mathcal{H}-big}\equiv \emptyset 
$ since in this case $L$ fails to satisfy (\ref{up stopping condition}) as
pointed out above. Finally, for $L\in \mathcal{H}$ we further decompose $%
\mathcal{P}_{L,0}^{\flat \mathcal{H}-small}$ as%
\begin{eqnarray*}
\mathcal{P}_{L,0}^{\flat \mathcal{H}-small} &=&\overset{\cdot }{\dbigcup }%
_{L^{\prime }\in \mathcal{T}\left( L\right) }\mathcal{P}_{L^{\prime
},0}^{\flat \mathcal{L}-small}; \\
\mathcal{P}_{L^{\prime },0}^{\flat \mathcal{L}-small} &\equiv &\left\{
\left( I,J\right) \in \mathcal{P}:I\in \mathcal{C}_{L^{\prime }}^{\mathcal{L}%
}\setminus \left\{ L^{\prime }\right\} \text{ and }J\in \mathcal{C}%
_{L^{\prime }}^{\mathcal{L},\flat \limfunc{shift}}\right\} .
\end{eqnarray*}%
Then we set%
\begin{eqnarray}
\mathcal{P}^{\flat big} &\equiv &\left\{ \dbigcup\limits_{L\in \mathcal{H}}%
\mathcal{P}_{L,0}^{\flat \mathcal{H}-big}\right\} \dbigcup \left\{
\dbigcup\limits_{t\geq 1}\dbigcup\limits_{L\in \mathcal{H}}\mathcal{P}%
_{L,t}^{\flat \mathcal{H}}\right\} ;  \label{def big small flat} \\
\mathcal{P}^{\flat small} &\equiv &\dbigcup\limits_{L\in \mathcal{L}}%
\mathcal{P}_{L,0}^{\flat \mathcal{L}-small}\text{ }.  \notag
\end{eqnarray}%
We observed above that every pair $\left( I,J\right) \in \mathcal{P}$ is
included in either $\mathcal{P}^{small}$ or $\mathcal{P}^{big}$, and it
follows that every pair $\left( I,J\right) \in \mathcal{P}$ is thus included
in either $\mathcal{P}^{\flat small}$ or $\mathcal{P}^{\flat big}$, simply
because the pairs $\left( I,J\right) $ have been shifted up by two dyadic
levels in the interval $J$. Thus the coronas $\mathcal{P}_{L,0}^{\flat 
\mathcal{L}-small}$ are now even \emph{smaller} than the regular coronas $%
\mathcal{P}_{L,0}^{\mathcal{L}-small}$, which permits the estimate (\ref%
{small claim' 3}) below to hold for the larger augmented size functional. On
the other hand, the coronas $\mathcal{P}_{L,0}^{\flat \mathcal{H}-big}$ and $%
\mathcal{P}_{L,t}^{\flat \mathcal{H}}$ are now bigger than before, requiring
the stronger straddling lemmas above in order to obtain the estimates (\ref%
{rest bounds}) below. More specifically, we will see that stopping forms
with pairs in $\mathcal{P}^{\flat big}$ will be estimated using the $\flat $%
Straddling and Substraddling Lemmas (Substraddling applies to part of $%
\mathcal{P}_{L,0}^{\flat \mathcal{H}-big}$ and $\flat $Straddling applies to
the remaining part of $\mathcal{P}_{L,0}^{\flat \mathcal{H}-big}$ and to $%
\mathcal{P}_{L,t}^{\flat \mathcal{H}}$), and it is here that the removal of
the corona-straddling collection $\mathcal{P}_{\func{cor}}^{A}$ is
essential, while forms with pairs in $\mathcal{P}^{\flat small}$ will be
absorbed.

\subsection{Size estimates}

Now we turn to proving the \emph{size estimates} we need for these
collections. Recall that the \emph{restricted} norm $\widehat{\mathfrak{N}}_{%
\limfunc{stop},\bigtriangleup ^{\omega }}^{A,\mathcal{P}}$ is the best
constant in the inequality%
\begin{equation*}
\left\vert \mathsf{B}\right\vert _{\limfunc{stop},\bigtriangleup ^{\omega
}}^{A,\mathcal{P}}\left( f,g\right) \leq \widehat{\mathfrak{N}}_{\limfunc{%
stop},\bigtriangleup ^{\omega }}^{A,\mathcal{P}}\left\Vert \mathsf{P}_{\Pi
_{1}\mathcal{P}}^{\sigma ,\mathbf{b}}f\right\Vert _{L^{2}\left( \sigma
\right) }^{\bigstar }\left\Vert \mathsf{P}_{\Pi _{2}\mathcal{P}}^{\omega ,%
\mathbf{b}^{\ast }}g\right\Vert _{L^{2}\left( \omega \right) }^{\bigstar },
\end{equation*}%
where $f\in L^{2}\left( \sigma \right) $ satisfies $E_{I}^{\sigma
}\left\vert f\right\vert \leq \alpha _{\mathcal{A}}\left( A\right) $ for all 
$I\in \mathcal{C}_{A}$, and $g\in L^{2}\left( \omega \right) $.

\begin{proposition}
\label{bottom up 3}Suppose $\rho $ in (\ref{def rho}) is greater than $1$,
and $\mathcal{P}$ is a \emph{reduced admissible} collection of pairs for a
dyadic interval $A$. Let $\mathcal{P}=\mathcal{P}^{\flat big}\dot{\cup}%
\mathcal{P}^{\flat small}$ be the decomposition satisfying\ (\ref{def big
small}) above, i.e.%
\begin{equation*}
\mathcal{P}=\left\{ \dbigcup\limits_{L\in \mathcal{H}}\mathcal{P}%
_{L,0}^{\flat \mathcal{H}-big}\right\} \dbigcup \left\{
\dbigcup\limits_{t\geq 1}\dbigcup\limits_{L\in \mathcal{H}}\mathcal{P}%
_{L,t}^{\flat \mathcal{H}}\right\} \ \cup \ \left( \dbigcup_{L\in \mathcal{L}%
}\mathcal{P}_{L,0}^{\flat \mathcal{L}-small}\right) \ .
\end{equation*}%
Then all of these collections $\mathcal{P}_{L,0}^{\flat \mathcal{L}-small}$, 
$\mathcal{P}_{L,0}^{\flat \mathcal{H}-big}$ and $\mathcal{P}_{L,t}^{\flat 
\mathcal{H}}$ are reduced admissible, and we have the estimate 
\begin{equation}
\mathcal{S}_{\limfunc{aug}\limfunc{size}}^{\alpha ,A}\left( \mathcal{P}%
_{L,0}^{\flat \mathcal{L}-small}\right) ^{2}\leq \left( \rho -1\right) 
\mathcal{S}_{\limfunc{aug}\limfunc{size}}^{\alpha ,A}\left( \mathcal{P}%
\right) ^{2},\ \ \ \ \ L\in \mathcal{L},  \label{small claim' 3}
\end{equation}%
and the localized norm bounds,%
\begin{eqnarray}
\widehat{\mathfrak{N}}_{\limfunc{stop},\bigtriangleup ^{\omega
}}^{A,\dbigcup\limits_{L\in \mathcal{H}}\mathcal{P}_{L,0}^{\flat \mathcal{H}%
-big}} &\leq &C\mathcal{S}_{\limfunc{aug}\limfunc{size}}^{\alpha ,A}\left( 
\mathcal{P}\right) ,  \label{rest bounds} \\
\widehat{\mathfrak{N}}_{\limfunc{stop},\bigtriangleup ^{\omega
}}^{A,\dbigcup\limits_{L\in \mathcal{H}}\mathcal{P}_{L,t}^{\flat \mathcal{H}%
}} &\leq &C\rho ^{-\frac{t}{2}}\mathcal{S}_{\limfunc{aug}\limfunc{size}%
}^{\alpha ,A}\left( \mathcal{P}\right) ,\ \ \ \ \ t\geq 1.  \notag
\end{eqnarray}
\end{proposition}

Using this proposition on size estimates, we can finish the proof of (\ref%
{First inequality}), and hence the proof of (\ref{B stop form 3}).

\begin{corollary}
The sublinear stopping form inequality (\ref{First inequality}) holds.
\end{corollary}

\begin{proof}
Recall that $\widehat{\mathfrak{N}}_{\limfunc{stop},\bigtriangleup ^{\omega
}}^{A,\mathcal{P}}$ is the best constant in the inequality (\ref{best hat}).
Since $\left\{ \mathcal{P}_{L,0}^{\flat \mathcal{L}-small}\right\} _{L\in 
\mathcal{L}}$ is a mutually orthogonal family of $A$-admissible pairs, the
Orthogonality Lemma \ref{mut orth} implies that%
\begin{equation*}
\widehat{\mathfrak{N}}_{\limfunc{stop},\bigtriangleup ^{\omega
}}^{A,\dbigcup\limits_{L\in \mathcal{L}}\mathcal{P}_{L,0}^{\flat \mathcal{L}%
-small}}\leq \sup_{L\in \mathcal{L}}\widehat{\mathfrak{N}}_{\limfunc{stop}%
,\bigtriangleup ^{\omega }}^{A,\mathcal{P}_{L,0}^{\flat \mathcal{L}-small}}.
\end{equation*}%
Using this, together with the decomposition of $\mathcal{P}$ and (\ref{rest
bounds}) above, we obtain%
\begin{eqnarray*}
\widehat{\mathfrak{N}}_{\limfunc{stop},\bigtriangleup ^{\omega }}^{A,%
\mathcal{P}} &\leq &\sup_{L\in \mathcal{H}}\widehat{\mathfrak{N}}_{\limfunc{%
stop},\bigtriangleup ^{\omega }}^{A,\dbigcup\limits_{L\in \mathcal{H}}%
\mathcal{P}_{L,0}^{\flat \mathcal{H}-big}}+\sum_{t=1}^{M+1}\sup_{L\in 
\mathcal{H}}\widehat{\mathfrak{N}}_{\limfunc{stop},\bigtriangleup ^{\omega
}}^{A,\dbigcup\limits_{L\in \mathcal{H}}\mathcal{P}_{L,t}^{\flat \mathcal{H}%
}}+\widehat{\mathfrak{N}}_{\limfunc{stop},\bigtriangleup ^{\omega
}}^{A,\dbigcup\limits_{L\in \mathcal{L}}\mathcal{P}_{L,0}^{\flat \mathcal{L}%
-small}} \\
&\lesssim &\mathcal{S}_{\limfunc{aug}\limfunc{size}}^{\alpha ,A}\left( 
\mathcal{P}\right) +\left( \sum_{t=1}^{M+1}\rho ^{-\frac{t}{2}}\right) 
\mathcal{S}_{\limfunc{aug}\limfunc{size}}^{\alpha ,A}\left( \mathcal{P}%
\right) +\sup_{L\in \mathcal{L}}\widehat{\mathfrak{N}}_{\limfunc{stop}%
,\bigtriangleup ^{\omega }}^{A,\mathcal{P}_{L,0}^{\flat \mathcal{L}-small}}\
.
\end{eqnarray*}%
Since the admissible collection $\mathcal{P}^{A}$ in (\ref{initial P}) that
arises in the stopping form is finite, we can define $\mathfrak{L}$ to be
the best constant in the inequality%
\begin{equation*}
\widehat{\mathfrak{N}}_{\limfunc{stop},\bigtriangleup ^{\omega }}^{A,%
\mathcal{P}}\leq \mathfrak{L}\mathcal{S}_{\limfunc{aug}\limfunc{size}%
}^{\alpha ,A}\left( \mathcal{P}\right) \text{ for all }A\text{-admissible
collections }\mathcal{P}.
\end{equation*}%
Now choose $\mathcal{P}$ so that 
\begin{equation*}
\frac{\widehat{\mathfrak{N}}_{\limfunc{stop},\bigtriangleup ^{\omega }}^{A,%
\mathcal{P}}}{\mathcal{S}_{\limfunc{aug}\limfunc{size}}^{\alpha ,A}\left( 
\mathcal{P}\right) }>\frac{1}{2}\mathfrak{L=}\frac{1}{2}\sup_{\mathcal{Q}%
\text{ is }A\text{-admissible}}\frac{\widehat{\mathfrak{N}}_{\limfunc{stop}%
,\bigtriangleup ^{\omega }}^{A,\mathcal{Q}}}{\mathcal{S}_{\limfunc{aug}%
\limfunc{size}}^{\alpha ,A}\left( \mathcal{Q}\right) }\ .
\end{equation*}%
Then using $\sum_{t=1}^{M+1}\rho ^{-\frac{t}{2}}\leq \frac{1}{\sqrt{\rho }-1}
$ we have%
\begin{eqnarray*}
\mathfrak{L} &<&2\frac{\widehat{\mathfrak{N}}_{\limfunc{stop},\bigtriangleup
^{\omega }}^{A,\mathcal{P}}}{\mathcal{S}_{\limfunc{aug}\limfunc{size}%
}^{\alpha ,A}\left( \mathcal{P}\right) }\leq \frac{C\frac{1}{\sqrt{\rho }-1}%
\mathcal{S}_{\limfunc{aug}\limfunc{size}}^{\alpha ,A}\left( \mathcal{P}%
\right) +C\sup_{L\in \mathcal{L}}\widehat{\mathfrak{N}}_{\limfunc{stop}%
,\bigtriangleup ^{\omega }}^{A,\mathcal{P}_{L,0}^{\flat \mathcal{L}-small}}}{%
\mathcal{S}_{\limfunc{aug}\limfunc{size}}^{\alpha ,A}\left( \mathcal{P}%
\right) } \\
&\leq &C\frac{1}{\sqrt{\rho }-1}+C\sup_{L\in \mathcal{L}}\mathfrak{L}\frac{%
\mathcal{S}_{\limfunc{aug}\limfunc{size}}^{\alpha ,A}\left( \mathcal{P}%
_{L,0}^{\flat \mathcal{L}-small}\right) }{\mathcal{S}_{\limfunc{aug}\limfunc{%
size}}^{\alpha ,A}\left( \mathcal{P}\right) }\leq C\frac{1}{\sqrt{\rho }-1}+C%
\mathfrak{L}\sqrt{\rho -1}\ ,
\end{eqnarray*}%
where we have used (\ref{small claim' 3}) in the last line. If we choose $%
\rho >1$ so that 
\begin{equation}
C\sqrt{\rho -1}<\frac{1}{2},  \label{choose rho}
\end{equation}%
then we obtain $\mathfrak{L}\leq 2C\frac{1}{\sqrt{\rho }-1}$. Together with
Lemma \ref{energy control}, this yields%
\begin{equation*}
\widehat{\mathfrak{N}}_{\limfunc{stop},\bigtriangleup ^{\omega }}^{A,%
\mathcal{P}}\leq \mathfrak{L}\mathcal{S}_{\limfunc{aug}\limfunc{size}%
}^{\alpha ,A}\left( \mathcal{P}\right) \leq 2C\frac{1}{\sqrt{\rho }-1}\left( 
\mathcal{E}_{2}^{\alpha }+\sqrt{A_{2}^{\alpha }}+\sqrt{A_{2}^{\alpha ,%
\limfunc{punct}}}\right)
\end{equation*}%
as desired, and completes the proof of inequality (\ref{First inequality}).
\end{proof}

Thus, in view of Conclusion \ref{assume}, it remains only to prove
Proposition \ref{bottom up 3} using the Orthogonality and Straddling and
Substraddling Lemmas above, and we now turn to this task.

\begin{proof}[Proof of Proposition \protect\ref{bottom up 3}]
We split the proof into three parts.

\textbf{Proof of (\ref{small claim' 3})}: To prove the inequality (\ref%
{small claim' 3}), suppose first that $L\notin \mathcal{L}_{M+1}$. In the
case that $L\in \mathcal{L}_{0}$ is an initial generation interval, then
from (\ref{key property 3}) and the fact that every $I\in \mathcal{P}%
_{L,0}^{\flat \mathcal{L}-small}$ satisfies $I\subsetneqq L$, we obtain that%
\begin{eqnarray*}
&&\mathcal{S}_{\limfunc{aug}\limfunc{size}}^{\alpha ,A}\left( \mathcal{P}%
_{L,0}^{\flat \mathcal{L}-small}\right) ^{2}=\sup_{K^{\prime }\in \Pi _{1}^{%
\limfunc{below}}\mathcal{P}_{L,0}^{\flat \mathcal{L}-small}\cap \mathcal{C}%
_{A}^{\limfunc{restrict}}}\frac{\Psi ^{\alpha }\left( K^{\prime };\mathcal{P}%
_{L,0}^{\flat \mathcal{L}-small}\right) ^{2}}{\left\vert K^{\prime
}\right\vert _{\sigma }} \\
&\leq &\sup_{K^{\prime }\in \Pi _{1}^{\limfunc{below}}\mathcal{P}\cap 
\mathcal{C}_{A}^{\limfunc{restrict}}:\ K^{\prime }\varsubsetneqq L}\frac{%
\Psi ^{\alpha }\left( K^{\prime };\mathcal{P}_{L,0}^{\flat \mathcal{L}%
-small}\right) ^{2}}{\left\vert K^{\prime }\right\vert _{\sigma }}\leq
\varepsilon \mathcal{S}_{\limfunc{aug}\limfunc{size}}^{\alpha ,A}\left( 
\mathcal{P}\right) ^{2}.
\end{eqnarray*}%
Now suppose that $L\not\in \mathcal{L}_{0}$ in addition to $L\notin \mathcal{%
L}_{M+1}$. Pick a pair $\left( I,J\right) \in \mathcal{P}_{L,0}^{\flat 
\mathcal{L}-small}$. Then $I$ is in the restricted corona $\mathcal{C}_{L}^{%
\mathcal{L},\limfunc{restrict}}$ and $J$ is in the $\flat $\emph{shifted}
corona $\mathcal{C}_{L}^{\mathcal{L},\flat \limfunc{shift}}$. Since $%
\mathcal{P}_{L,0}^{\flat \mathcal{L}-small}$ is a finite collection, the
definition of $\mathcal{S}_{\limfunc{aug}\limfunc{size}}^{\alpha ,A}\left( 
\mathcal{P}_{L,0}^{\flat \mathcal{L}-small}\right) $ shows that there is an
interval $K\in \Pi _{1}^{\limfunc{below}}\mathcal{P}_{L,0}^{\flat \mathcal{L}%
-small}\cap \mathcal{C}_{A}^{\limfunc{restrict}}$ so that%
\begin{equation*}
\mathcal{S}_{\limfunc{aug}\limfunc{size}}^{\alpha ,A}\left( \mathcal{P}%
_{L,0}^{\flat \mathcal{L}-small}\right) ^{2}=\frac{1}{\left\vert
K\right\vert _{\sigma }}\left( \frac{\mathrm{P}^{\alpha }\left( K,\mathbf{1}%
_{A\setminus K}\sigma \right) }{\left\vert K\right\vert }\right) ^{2}\omega
_{\flat \mathcal{P}}\left( \mathbf{T}\left( K\right) \right) .
\end{equation*}%
Note that $K\subsetneqq L$ by definition of $\mathcal{P}_{L,0}^{\flat 
\mathcal{L}-small}$. Now let $t$ be such that $L\in \mathcal{L}_{t}$, and
define 
\begin{equation*}
t^{\prime }=t^{\prime }\left( K\right) \equiv \max \left\{ s:\text{there is }%
L^{\prime }\in \mathcal{L}_{s}\text{ with }L^{\prime }\subset K\right\} ,
\end{equation*}%
and note that $t^{\prime }<t$. First, suppose that $t^{\prime }=0$ so that $%
K $ does not contain any $L^{\prime }\in \mathcal{L}$. Then it follows from
the construction at level $\ell =0$ that%
\begin{equation*}
\frac{1}{\left\vert K\right\vert _{\sigma }}\left( \frac{\mathrm{P}^{\alpha
}\left( K,\mathbf{1}_{A\setminus K}\sigma \right) }{\left\vert K\right\vert }%
\right) ^{2}\omega _{\flat \mathcal{P}}\left( \mathbf{T}\left( K\right)
\right) <\varepsilon \mathcal{S}_{\limfunc{aug}\limfunc{size}}^{\alpha
,A}\left( \mathcal{P}\right) ^{2},
\end{equation*}%
and hence from $\rho =1+\varepsilon $ we obtain 
\begin{equation*}
\mathcal{S}_{\limfunc{aug}\limfunc{size}}^{\alpha ,A}\left( \mathcal{P}%
_{L,0}^{\flat \mathcal{L}-small}\right) ^{2}<\varepsilon \mathcal{S}_{%
\limfunc{aug}\limfunc{size}}^{\alpha ,A}\left( \mathcal{P}\right)
^{2}=\left( \rho -1\right) \mathcal{S}_{\limfunc{aug}\limfunc{size}}^{\alpha
,A}\left( \mathcal{P}\right) ^{2}.
\end{equation*}%
Now suppose that $t^{\prime }\geq 1$. Then $K$ fails the stopping condition (%
\ref{up stopping condition}) with $m=t^{\prime }+1$, since otherwise it
would contain an interval $L^{\prime \prime }\in \mathcal{L}_{t^{\prime }+1}$
contradicting our definition of $t^{\prime }$, and so%
\begin{equation*}
\omega _{\flat \mathcal{P}}\left( \mathbf{T}\left( K\right) \right) <\rho
\omega _{\flat \mathcal{P}}\left( \mathbf{V}\left( K\right) \right) \text{
where }\mathbf{V}\left( K\right) \equiv \dbigcup\limits_{L^{\prime }\in
\dbigcup\limits_{\ell =0}^{t^{\prime }}\mathcal{L}_{\ell }:\ L^{\prime
}\subset K}\mathbf{T}\left( L^{\prime }\right) .
\end{equation*}

Now we use the crucial fact that the positive measure $\omega _{\flat 
\mathcal{P}}$ is \emph{additive} and finite to obtain from this that%
\begin{equation}
\omega _{\flat \mathcal{P}}\left( \mathbf{T}\left( K\right) \setminus 
\mathbf{V}\left( K\right) \right) =\omega _{\flat \mathcal{P}}\left( \mathbf{%
T}\left( K\right) \right) -\omega _{\flat \mathcal{P}}\left( \mathbf{V}%
\left( K\right) \right) \leq \left( \rho -1\right) \omega _{\flat \mathcal{P}%
}\left( \mathbf{V}\left( K\right) \right) .  \label{additive}
\end{equation}%
Now recall that%
\begin{equation*}
\mathcal{S}_{\limfunc{aug}\limfunc{size}}^{\alpha ,A}\left( \mathcal{Q}%
\right) ^{2}\equiv \sup_{K\in \Pi _{1}^{\limfunc{below}}\mathcal{Q}\cap
C_{A}^{\limfunc{restrict}}}\frac{1}{\left\vert K\right\vert _{\sigma }}%
\left( \frac{\mathrm{P}^{\alpha }\left( K,\mathbf{1}_{A\setminus K}\sigma
\right) }{\left\vert K\right\vert }\right) ^{2}\left\Vert \mathsf{Q}_{\Pi
_{2}^{K,\limfunc{aug}}\mathcal{Q}}^{\omega ,\mathbf{b}^{\ast }}x\right\Vert
_{L^{2}\left( \omega \right) }^{\spadesuit 2}.
\end{equation*}%
We claim it follows that for each $J\in \Pi _{2}^{K,\limfunc{aug}}\mathcal{P}%
_{L,0}^{\flat \mathcal{L}-small}$, the support $\left( c_{J^{\flat }},\ell
\left( J^{\flat }\right) \right) $ of the atom $\delta _{\left( c_{J^{\flat
}},\ell \left( J^{\flat }\right) \right) }$ is contained in the set $\mathbf{%
T}\left( K\right) $, but not in the set 
\begin{equation*}
\mathbf{V}\left( K\right) \equiv \dbigcup \left\{ \mathbf{T}\left( L^{\prime
}\right) :L^{\prime }\in \dbigcup\limits_{\ell =0}^{t^{\prime }}\mathcal{L}%
_{\ell }:\ L^{\prime }\subset K\right\} .
\end{equation*}%
Indeed, suppose in order to derive a contradiction, that $\left( c_{J^{\flat
}},\ell \left( J^{\flat }\right) \right) \in \mathbf{T}\left( L^{\prime
}\right) $ for some $L^{\prime }\in \mathcal{L}_{\ell }$ with $0\leq \ell
\leq t^{\prime }$. Recall that $L\in \mathcal{L}_{t}$ with $t^{\prime }<t$
so that $L^{\prime }\subsetneqq L$. Thus $\left( c_{J^{\flat }},\ell \left(
J^{\flat }\right) \right) \in \mathbf{T}\left( L^{\prime }\right) $ implies $%
J^{\flat }\subset L^{\prime }$, which contradicts the fact that 
\begin{equation*}
J\in \Pi _{2}^{K}\mathcal{P}_{L,0}^{\flat \mathcal{L}-small}\subset \Pi _{2}%
\mathcal{P}_{L,0}^{\flat \mathcal{L}-small}=\left\{ \left( I,J\right) \in 
\mathcal{P}:I\in \mathcal{C}_{L}^{\mathcal{L}}\setminus \left\{ L\right\} 
\text{ and }J\in \mathcal{C}_{L}^{\mathcal{L},\flat \limfunc{shift}}\right\}
\end{equation*}%
implies $J^{\flat }\in \mathcal{C}_{L}^{\mathcal{L}}$ - because $L^{\prime
}\notin \mathcal{C}_{L}^{\mathcal{L}}$.

Thus from the definition of $\omega _{\flat \mathcal{P}}$ in (\ref{def
atomic}), the `energy' $\left\Vert \mathsf{Q}_{\Pi _{2}^{K,\limfunc{aug}}%
\mathcal{P}_{L,0}^{\flat \mathcal{L}-small}}^{\omega ,\mathbf{b}^{\ast
}}x\right\Vert _{L^{2}\left( \omega \right) }^{\spadesuit 2}$ is at most the 
$\omega _{\flat \mathcal{P}}$-measure of $\mathbf{T}\left( K\right)
\setminus \mathbf{V}\left( K\right) $. Using 
\begin{equation*}
\omega _{\flat \mathcal{P}_{L,0}^{\flat \mathcal{L}-small}}\left( \mathbf{T}%
\left( K\right) \right) \leq \omega _{\flat \mathcal{P}}\left( \mathbf{T}%
\left( K\right) \setminus \mathbf{V}\left( K\right) \right) ,
\end{equation*}%
and (\ref{additive}), we then have%
\begin{eqnarray*}
&&\mathcal{S}_{\limfunc{aug}\limfunc{size}}^{\alpha ,A}\left( \mathcal{P}%
_{L,0}^{\flat \mathcal{L}-small}\right) ^{2} \\
&\leq &\sup_{K\in \Pi _{1}^{\limfunc{below}}\mathcal{P}_{L,0}^{\flat 
\mathcal{L}-small}\cap \mathcal{C}_{A}^{\limfunc{restrict}}}\frac{1}{%
\left\vert K\right\vert _{\sigma }}\left( \frac{\mathrm{P}^{\alpha }\left( K,%
\mathbf{1}_{A\setminus K}\sigma \right) }{\left\vert K\right\vert }\right)
^{2}\omega _{\flat \mathcal{P}}\left( \mathbf{T}\left( K\right) \setminus 
\mathbf{V}\left( K\right) \right) \\
&\leq &\left( \rho -1\right) \sup_{K\in \Pi _{1}^{\limfunc{below}}\mathcal{P}%
_{L,0}^{\flat \mathcal{L}-small}\cap \mathcal{C}_{A}^{\limfunc{restrict}}}%
\frac{1}{\left\vert K\right\vert _{\sigma }}\left( \frac{\mathrm{P}^{\alpha
}\left( K,\mathbf{1}_{A\setminus K}\sigma \right) }{\left\vert K\right\vert }%
\right) ^{2}\omega _{\flat \mathcal{P}}\left( \mathbf{V}\left( K\right)
\right) ,
\end{eqnarray*}%
and we can continue with 
\begin{eqnarray*}
&&\mathcal{S}_{\limfunc{aug}\limfunc{size}}^{\alpha ,A}\left( \mathcal{P}%
_{L,0}^{\flat \mathcal{L}-small}\right) ^{2} \\
&\leq &\left( \rho -1\right) \sup_{K\in \Pi _{1}^{\limfunc{below}}\mathcal{P}%
\cap \mathcal{C}_{A}^{\limfunc{restrict}}}\frac{1}{\left\vert K\right\vert
_{\sigma }}\left( \frac{\mathrm{P}^{\alpha }\left( K,\mathbf{1}_{A\setminus
K}\sigma \right) }{\left\vert K\right\vert }\right) ^{2}\omega _{\flat 
\mathcal{P}}\left( \mathbf{T}\left( K\right) \right) \\
&\leq &\left( \rho -1\right) \mathcal{S}_{\limfunc{aug}\limfunc{size}%
}^{\alpha ,A}\left( \mathcal{P}\right) ^{2}.
\end{eqnarray*}

In the remaining case where $L\in \mathcal{L}_{M+1}$ we can include $L$ as a
testing interval $K$ and the same reasoning applies. This completes the
proof of (\ref{small claim' 3}).

\bigskip

To prove the other inequality (\ref{rest bounds}) in Proposition \ref{bottom
up 3}, we will use the Straddling and Substraddling Lemmas to bound the norm
of certain `straddled' stopping forms by the augmented size functional $%
\mathcal{S}_{\limfunc{aug}\limfunc{size}}^{\alpha ,A}$, and the
Orthogonality Lemma to bound sums of `mutually orthogonal' stopping forms.
Recall that 
\begin{eqnarray*}
\mathcal{P}^{\flat big} &=&\left\{ \dbigcup\limits_{L\in \mathcal{H}}%
\mathcal{P}_{L,0}^{\flat \mathcal{H}-big}\right\} \dbigcup \left\{
\dbigcup\limits_{t\geq 1}\dbigcup\limits_{L\in \mathcal{H}}\mathcal{P}%
_{L,t}^{\flat \mathcal{H}}\right\} \equiv \mathcal{Q}_{0}^{\flat \mathcal{H}%
-big}\dbigcup \mathcal{Q}_{1}^{\flat \mathcal{H}-big}; \\
\mathcal{Q}_{0}^{\flat \mathcal{H}-big} &\equiv &\dbigcup\limits_{L\in 
\mathcal{L}}\mathcal{P}_{L,0}^{\flat \mathcal{H}-big}\ ,\ \ \ \ \ \mathcal{Q}%
_{1}^{\flat \mathcal{H}-big}\equiv \dbigcup\limits_{t\geq 1}\mathcal{P}%
_{t}^{\flat \mathcal{H}-big},\ \ \ \ \ \mathcal{P}_{t}^{\flat \mathcal{H}%
-big}\equiv \dbigcup\limits_{L\in \mathcal{H}}\mathcal{P}_{L,t}^{\flat 
\mathcal{H}}.
\end{eqnarray*}

\bigskip

\textbf{Proof of the second line in (\ref{rest bounds})}: We first turn to
the collection%
\begin{eqnarray*}
\mathcal{Q}_{1}^{\flat \mathcal{H}-big} &=&\dbigcup\limits_{t\geq
1}\dbigcup\limits_{L\in \mathcal{H}}\mathcal{P}_{L,t}^{\flat \mathcal{H}%
}=\dbigcup\limits_{t\geq 1}\mathcal{P}_{t}^{\flat \mathcal{H}-big}; \\
\mathcal{P}_{t}^{\flat \mathcal{H}-big} &\equiv &\dbigcup\limits_{L\in 
\mathcal{L}}\mathcal{P}_{L,t}^{\flat \mathcal{H}}\ ,\ \ \ \ \ t\geq 1,
\end{eqnarray*}%
where%
\begin{equation*}
\mathcal{P}_{L,t}^{\flat \mathcal{H}}=\left\{ \left( I,J\right) \in \mathcal{%
P}:I\in \mathcal{C}_{L}^{\mathcal{H}}\text{ and }J\in \mathcal{C}_{L^{\prime
}}^{\mathcal{H},\flat \limfunc{shift}}\text{ for some }L^{\prime }\in 
\mathcal{H}_{d_{\limfunc{indent}}\left( L\right) +t}\text{ with }L^{\prime
}\subset L\right\} .
\end{equation*}%
We now claim that the second line in (\ref{rest bounds}) holds, i.e.%
\begin{equation}
\widehat{\mathfrak{N}}_{\limfunc{stop},\bigtriangleup ^{\omega }}^{A,%
\mathcal{P}_{t}^{\flat \mathcal{H}-big}}\leq C\rho ^{-\frac{t}{2}}\mathcal{S}%
_{\limfunc{aug}\limfunc{size}}^{\alpha ,A}\left( \mathcal{P}\right) ,\ \ \ \
\ t\geq 1,  \label{S big t 3}
\end{equation}%
which recovers the key geometric gain obtained by Lacey in \cite{Lac},
except that here we are only gaining this decay relative to the indented
subtree $\mathcal{H}$ of the tree $\mathcal{L}$.

The case $t=1$ can be handled with relative ease since decay is not relevant
here. Indeed, $\mathcal{P}_{L,1}^{\flat \mathcal{H}}$ straddles the
collection $\mathfrak{C}_{\mathcal{H}}\left( L\right) $ of $\mathcal{H}$%
-children of $L$, and so the localized $\flat $Straddling Lemma \ref%
{straddle 3 ref} applies to give%
\begin{equation*}
\widehat{\mathfrak{N}}_{\limfunc{stop},\bigtriangleup ^{\omega }}^{A,%
\mathcal{P}_{L,1}^{\flat \mathcal{H}}}\leq C\mathcal{S}_{\limfunc{aug}%
\limfunc{size}}^{\alpha ,A}\left( \mathcal{P}_{L,1}^{\flat \mathcal{H}%
}\right) \leq C\mathcal{S}_{\limfunc{aug}\limfunc{size}}^{\alpha ,A}\left( 
\mathcal{P}\right) ,
\end{equation*}%
and then the Orthogonality Lemma \ref{mut orth} applies to give%
\begin{equation*}
\widehat{\mathfrak{N}}_{\limfunc{stop},\bigtriangleup ^{\omega }}^{A,%
\mathcal{P}_{1}^{\flat \mathcal{H}-big}}\leq \sup_{L\in \mathcal{H}}%
\mathfrak{N}_{\limfunc{stop},\bigtriangleup ^{\omega }}^{A,\mathcal{P}%
_{L,1}^{\flat \mathcal{H}}}\leq C\mathcal{S}_{\limfunc{aug}\limfunc{size}%
}^{\alpha ,A}\left( \mathcal{P}\right) ,
\end{equation*}%
since $\left\{ \mathcal{P}_{L,1}^{\flat \mathcal{H}}\right\} _{L\in \mathcal{%
L}}$ is mutually orthogonal as $\mathcal{P}_{L,1}^{\flat \mathcal{H}}\subset 
\mathcal{C}_{L}^{\mathcal{H}}\times \mathcal{C}_{L^{\prime }}^{\mathcal{H}%
,\flat \limfunc{shift}}$ with $L\in \mathcal{H}_{k}$ and $L^{\prime }\in 
\mathcal{H}_{k+1}$ for indented depth $k=k\left( L\right) $. The case $t=2$
is equally easy.

Now we consider the case $t\geq 2$, where it is essential to obtain
geometric decay in $t$. We remind the reader that all of our admissible
collections $\mathcal{P}_{L,t}^{\flat \mathcal{H}}$ are \emph{reduced} by
Conclusion \ref{assume}. We again apply Lemma \ref{straddle 3 ref} to $%
\mathcal{P}_{L,t}^{\flat \mathcal{H}}$ with $\mathcal{S}=\mathfrak{C}_{%
\mathcal{H}}\left( L\right) $, so that for any $\left( I,J\right) \in 
\mathcal{P}_{L,t}^{\flat \mathcal{H}}$, there is $H^{\prime }\in \mathfrak{C}%
_{\mathcal{H}}\left( L\right) $ with $J^{\flat }\subset H^{\prime
}\subsetneqq I\in \mathcal{C}_{L}^{\mathcal{H}}$. But this time we must use
the stronger localized bounds $\mathcal{S}_{\limfunc{loc}\limfunc{size}%
}^{\alpha ,A;S}$ with an $S$-hole, that give%
\begin{eqnarray}
\widehat{\mathfrak{N}}_{\limfunc{stop},\bigtriangleup ^{\omega }}^{A,%
\mathcal{P}_{L,t}^{\flat \mathcal{H}}} &\leq &C\sup_{H^{\prime }\in 
\mathfrak{C}_{\mathcal{H}}\left( L\right) }\mathcal{S}_{\limfunc{loc}%
\limfunc{size}}^{\alpha ,A;H^{\prime }}\left( \mathcal{P}_{L,t}^{\flat 
\mathcal{H}}\right) ,\ \ \ \ \ t\geq 0;  \label{t,n 3} \\
\mathcal{S}_{\limfunc{loc}\limfunc{size}}^{\alpha ,A;H^{\prime }}\left( 
\mathcal{P}_{L,t}^{\flat \mathcal{H}}\right) ^{2} &=&\sup_{K\in \mathcal{W}%
^{\ast }\left( H^{\prime }\right) \cap \mathcal{C}_{A}^{\limfunc{restrict}}}%
\frac{1}{\left\vert K\right\vert _{\sigma }}\left( \frac{\mathrm{P}^{\alpha
}\left( K,\mathbf{1}_{A\setminus H^{\prime }}\sigma \right) }{\left\vert
K\right\vert }\right) ^{2}\sum_{J\in \Pi _{2}^{K,\limfunc{aug}}\mathcal{P}%
_{L,t}^{\flat \mathcal{H}}}\left\Vert \bigtriangleup _{J}^{\omega ,\mathbf{b}%
^{\ast }}x\right\Vert _{L^{2}\left( \omega \right) }^{\spadesuit 2}\ . 
\notag
\end{eqnarray}

It remains to show that%
\begin{eqnarray}
&&\sum_{J\in \Pi _{2}^{K,\limfunc{aug}}\mathcal{P}_{L,t}^{\flat \mathcal{H}%
}}\left\Vert \bigtriangleup _{J}^{\omega ,\mathbf{b}^{\ast }}x\right\Vert
_{L^{2}\left( \omega \right) }^{\spadesuit 2}\leq \rho ^{-\left( t-2\right)
}\omega _{\flat \mathcal{P}}\left( \mathbf{T}\left( K\right) \right) ,
\label{rem} \\
&&\text{for}\ t\geq 2,\ K\in \mathcal{W}^{\ast }\left( H^{\prime }\right)
\cap \mathcal{C}_{A}^{\limfunc{restrict}},\ H^{\prime }\in \mathfrak{C}_{%
\mathcal{H}}\left( L\right) .  \notag
\end{eqnarray}%
so that we then have%
\begin{eqnarray*}
&&\frac{1}{\left\vert K\right\vert _{\sigma }}\left( \frac{\mathrm{P}%
^{\alpha }\left( K,\mathbf{1}_{A\setminus H^{\prime }}\sigma \right) }{%
\left\vert K\right\vert }\right) ^{2}\sum_{J\in \Pi _{2}^{K,\limfunc{aug}}%
\mathcal{P}_{L,t}^{\flat \mathcal{H}}}\left\Vert \bigtriangleup _{J}^{\omega
,\mathbf{b}^{\ast }}x\right\Vert _{L^{2}\left( \omega \right) }^{\spadesuit
2} \\
&\leq &\rho ^{-\left( t-2\right) }\frac{1}{\left\vert K\right\vert _{\sigma }%
}\left( \frac{\mathrm{P}^{\alpha }\left( K,\mathbf{1}_{A\setminus K}\sigma
\right) }{\left\vert K\right\vert }\right) ^{2}\omega _{\flat \mathcal{P}%
}\left( \mathbf{T}\left( K\right) \right) \leq \rho ^{-\left( t-2\right) }%
\mathcal{S}_{\limfunc{aug}\limfunc{size}}^{\alpha ,A}\left( \mathcal{P}%
\right) ^{2},
\end{eqnarray*}%
by (\ref{def P stop energy' 3}), and hence conclude the required bound for $%
\mathfrak{N}_{\limfunc{stop},\bigtriangleup ^{\omega }}^{A,\mathcal{P}%
_{L,t}^{\flat \mathcal{H}}}$, namely that%
\begin{eqnarray}
&&  \label{N_L} \\
\widehat{\mathfrak{N}}_{\limfunc{stop},\bigtriangleup ^{\omega }}^{A,%
\mathcal{P}_{L,t}^{\flat \mathcal{H}}} &\leq &C\sup_{H^{\prime }\in 
\mathfrak{C}_{\mathcal{H}}\left( L\right) }\sup_{K\in \mathcal{W}^{\ast
}\left( H^{\prime }\right) \cap \mathcal{C}_{A}^{\limfunc{restrict}}}\sqrt{%
\frac{1}{\left\vert K\right\vert _{\sigma }}\left( \frac{\mathrm{P}^{\alpha
}\left( K,\mathbf{1}_{A\setminus H^{\prime }}\sigma \right) }{\left\vert
K\right\vert }\right) ^{2}\sum_{J\in \Pi _{2}^{K,\limfunc{aug}}\mathcal{P}%
_{L,t}^{\flat \mathcal{H}}}\left\Vert \bigtriangleup _{J}^{\omega ,\mathbf{b}%
^{\ast }}x\right\Vert _{L^{2}\left( \omega \right) }^{\spadesuit 2}}  \notag
\\
&\leq &C\sqrt{\rho ^{-\left( t-2\right) }}\mathcal{S}_{\limfunc{aug}\limfunc{%
size}}^{\alpha ,A}\left( \mathcal{P}\right) =C^{\prime }\rho ^{-\frac{t}{2}}%
\mathcal{S}_{\limfunc{aug}\limfunc{size}}^{\alpha ,A}\left( \mathcal{P}%
\right) .  \notag
\end{eqnarray}

\medskip

\textbf{Remark on lack of usual goodness}: To prove (\ref{rem}), it is
essential that the intervals $H^{k+2}\in \mathcal{H}_{k+2}$ at the next
indented level down from $H^{k+1}\in \mathfrak{C}_{\mathcal{H}}\left(
L\right) $ are each contained in one of the Whitney intervals $K\in \mathcal{%
W}\left( H^{k+1}\right) \cap \mathcal{C}_{A}^{\limfunc{restrict}}$ for some $%
H^{k+1}\in \mathfrak{C}_{\mathcal{H}}\left( L\right) $. And this is the
reason we introduced the indented corona - namely so that $3H^{k+2}\subset
H^{k+1}$ for some $H^{k+1}\in \mathfrak{C}_{\mathcal{H}}\left( L\right) $,
and hence $H^{k+2}\subset K$ for some $K\in \mathcal{W}\left( H^{k+1}\right) 
$. In the argument of Lacey in \cite{Lac}, the corresponding intervals were
good in the usual sense, and so the above triple property was automatic.

\medskip

So we begin by fixing $K\in \mathcal{W}^{\ast }\left( H^{k+1}\right) \cap 
\mathcal{C}_{A}^{\limfunc{restrict}}$ with $H^{k+1}\in \mathfrak{C}_{%
\mathcal{H}}\left( L\right) $, and noting from the above that each $J\in \Pi
_{2}^{K,\limfunc{aug}}\mathcal{P}_{L,t}^{\flat \mathcal{H}}$ satisfies 
\begin{equation*}
J^{\flat }\subset H^{k+t}\subset H^{k+t-1}\subset ...\subset H^{k+2}\subset K
\end{equation*}%
for $H^{k+j}\in \mathcal{H}_{k+j}$ uniquely determined by $J^{\flat }$. Thus
for $t\geq 2$ we have 
\begin{eqnarray*}
\sum_{J\in \Pi _{2}^{K,\limfunc{aug}}\mathcal{P}_{L,t}^{\flat \mathcal{H}%
}}\left\Vert \bigtriangleup _{J}^{\omega ,\mathbf{b}^{\ast }}x\right\Vert
_{L^{2}\left( \omega \right) }^{\spadesuit 2} &=&\sum_{H^{k+t}\in \mathcal{H}%
_{k+t}:\ H^{k+t}\subset K}\sum_{J\in \Pi _{2}^{K,\limfunc{aug}}\mathcal{P}%
_{L,t}^{\flat \mathcal{H}}:\ J^{\flat }\subset H^{k+t}}\left\Vert
\bigtriangleup _{J}^{\omega ,\mathbf{b}^{\ast }}x\right\Vert _{L^{2}\left(
\omega \right) }^{\spadesuit 2} \\
&\leq &\sum_{H^{k+t}\in \mathcal{H}_{k+t}:\ H^{k+t}\subset K}\omega _{\flat 
\mathcal{P}}\left( \mathbf{T}\left( H^{k+t}\right) \right) .
\end{eqnarray*}%
In the case $t=2$ we are done since the final sum above is at most $\omega
_{\flat \mathcal{P}}\left( \mathbf{T}\left( K\right) \right) $.

Now suppose $t\geq 3$. In order to obtain geometric gain in $t$, we will
apply the stopping criterion (\ref{up stopping condition}) in the following
form,%
\begin{equation}
\sum_{L^{\prime }\in \mathfrak{C}_{\mathcal{L}}\left( L_{0}\right) }\omega
_{\flat \mathcal{P}}\left( \mathbf{T}\left( L^{\prime }\right) \right)
=\omega _{\flat \mathcal{P}}\left( \dbigcup\limits_{L^{\prime }\in \mathfrak{%
C}_{\mathcal{L}}\left( L_{0}\right) }\mathbf{T}\left( L^{\prime }\right)
\right) \leq \frac{1}{\rho }\omega _{\flat \mathcal{P}}\left( \mathbf{T}%
\left( L_{0}\right) \right) ,\ \ \ \ \ \text{for all }L_{0}\in \mathcal{L},
\label{foll form}
\end{equation}%
where we have used the fact that the \emph{maximal} intervals $L^{\prime }$
in the collection $\dbigcup\limits_{\ell =0}^{m-1}\left\{ L^{\prime }\in 
\mathcal{L}_{\ell }:\ L^{\prime }\subset L_{0}\right\} $ for $L_{0}\in 
\mathcal{L}_{m}$ (that appears in (\ref{up stopping condition})) are
precisely the $\mathcal{L}$-children of $L_{0}$ in the tree $\mathcal{L}$
(the intervals $L^{\prime }$ above are strictly contained in $L_{0}$ since $%
\rho >1$ in (\ref{up stopping condition})), so that%
\begin{equation*}
\dbigcup\limits_{L^{\prime }\in \Gamma }L^{\prime
}=\dbigcup\limits_{L^{\prime }\in \mathfrak{C}_{\mathcal{L}}\left(
L_{0}\right) }L^{\prime }\text{ where }\Gamma \equiv \dbigcup\limits_{\ell
=0}^{m-1}\left\{ L^{\prime }\in \mathcal{L}_{\ell }:\ L^{\prime }\subset
L_{0}\right\} .
\end{equation*}

In order to apply (\ref{foll form}), we collect the pairwise disjoint
intervals $H^{k+t}\in \mathcal{H}_{k+t}$ such that$\ H^{k+t}\subset
H^{k+2}\subset K$, into groups according to which interval $L^{k^{\prime
}+t-2}\in \mathcal{G}_{k^{\prime }+t-2}$ they are contained in, where $%
k^{\prime }=d_{\limfunc{geom}}\left( H^{k+2}\right) $ is the geometric depth
of $H^{k+2}$ in the tree $\mathcal{L}$ introduced in (\ref{geom depth}). It
follows that each interval $H^{k+t}\in \mathcal{H}_{k+t}$ is contained in a
unique interval $L^{d_{\limfunc{geom}}\left( H^{k+2}\right) +t-2}\in 
\mathcal{G}_{d_{\limfunc{geom}}\left( H^{k+2}\right) +t-2}$. Thus we\ obtain
from the previous inequality that%
\begin{eqnarray*}
\sum_{J\in \Pi _{2}^{K,\limfunc{aug}}\mathcal{P}_{L,t}^{\flat \mathcal{H}%
}}\left\Vert \bigtriangleup _{J}^{\omega ,\mathbf{b}^{\ast }}x\right\Vert
_{L^{2}\left( \omega \right) }^{\spadesuit 2} &\leq &\sum_{H^{k+t}\in 
\mathcal{H}_{k+t}:\ H^{k+t}\subset K}\omega _{\flat \mathcal{P}}\left( 
\mathbf{T}\left( H^{k+t}\right) \right) \\
&\leq &\sum_{\substack{ H^{k+2}\in \mathcal{H}_{k+2}  \\ H^{k+2}\subset K}}%
\sum_{\substack{ L^{k^{\prime }+t-2}\in \mathcal{G}_{k^{\prime }+t-2}:\
L^{k^{\prime }+t-2}\subset H^{k+2}  \\ \text{where }k^{\prime }=d_{\limfunc{%
geom}}\left( H^{k+2}\right) }}\omega _{\flat \mathcal{P}}\left( \mathbf{T}%
\left( L^{k^{\prime }+t-2}\right) \right) .
\end{eqnarray*}%
In the case $t=2$ we are done since the final sum above is dominated by 
\begin{equation*}
\sum_{H^{k+2}\in \mathcal{H}_{k+2}:\ H^{k+2}\subset K}\omega _{\flat 
\mathcal{P}}\left( \mathbf{T}\left( H^{k+2}\right) \right) \leq \omega
_{\flat \mathcal{P}}\left( \mathbf{T}\left( K\right) \right) .
\end{equation*}

For $t\geq 3$, we have 
\begin{eqnarray*}
&&\sum_{J\in \Pi _{2}^{K,\limfunc{aug}}\mathcal{P}_{L,t}^{\flat \mathcal{H}%
}}\left\Vert \bigtriangleup _{J}^{\omega ,\mathbf{b}^{\ast }}x\right\Vert
_{L^{2}\left( \omega \right) }^{\spadesuit 2}\leq \sum_{\substack{ %
H^{k+2}\in \mathcal{H}_{k+2}  \\ H^{k+2}\subset K}}\sum_{\substack{ %
L^{k^{\prime }+t-2}\in \mathcal{G}_{k^{\prime }+t-2}:\ L^{k^{\prime
}+t-2}\subset H^{k+2}  \\ \text{where }k^{\prime }=d_{\limfunc{geom}}\left(
H^{k+2}\right) }}\omega _{\flat \mathcal{P}}\left( \mathbf{T}\left(
L^{k^{\prime }+t-2}\right) \right) \\
&=&\sum_{\substack{ H^{k+2}\in \mathcal{H}_{k+2}  \\ H^{k+2}\subset K}}\sum 
_{\substack{ L^{k^{\prime }+t-3}\in \mathcal{G}_{k^{\prime }+t-3}:\
L^{k^{\prime }+t-3}\subset H^{k+2}  \\ \text{where }k^{\prime }=d_{\limfunc{%
geom}}\left( H^{k+2}\right) }}\left\{ \sum_{\substack{ L^{k^{\prime
}+t-2}\in \mathcal{G}_{k^{\prime }+t-2}:\ L^{k^{\prime }+t-2}\subset
L^{k^{\prime }+t-3}  \\ \text{where }k^{\prime }=d_{\limfunc{geom}}\left(
H^{k+2}\right) }}\omega _{\flat \mathcal{P}}\left( \mathbf{T}\left(
L^{k^{\prime }+t-2}\right) \right) \right\} \\
&\leq &\sum_{\substack{ H^{k+2}\in \mathcal{H}_{k+2}  \\ H^{k+2}\subset K}}%
\sum_{\substack{ L^{k^{\prime }+t-3}\in \mathcal{G}_{k^{\prime }+t-3}:\
L^{k^{\prime }+t-3}\subset H^{k+2}  \\ \text{where }k^{\prime }=d_{\limfunc{%
geom}}\left( H^{k+2}\right) }}\left\{ \frac{1}{\rho }\omega _{\flat \mathcal{%
P}}\left( \mathbf{T}\left( L^{k^{\prime }+t-3}\right) \right) \right\} ,
\end{eqnarray*}%
where in the last line we have used (\ref{foll form}) with $%
L_{0}=L^{k^{\prime }+t-3}$ on the sum in braces. We then continue (if
necessary) with%
\begin{eqnarray*}
\sum_{J\in \Pi _{2}^{K,\limfunc{aug}}\mathcal{P}_{L,t}^{\flat \mathcal{H}%
}}\left\Vert \bigtriangleup _{J}^{\omega ,\mathbf{b}^{\ast }}x\right\Vert
_{L^{2}\left( \omega \right) }^{\spadesuit 2} &\leq &\frac{1}{\rho }\sum 
_{\substack{ H^{k+2}\in \mathcal{H}_{k+2}  \\ H^{k+2}\subset K}}\sum 
_{\substack{ L^{k^{\prime }+t-3}\in \mathcal{G}_{k^{\prime }+t-3}:\
L^{k^{\prime }+t-3}\subset H^{k+2}  \\ \text{where }k^{\prime }=d_{\limfunc{%
geom}}\left( H^{k+2}\right) }}\omega _{\flat \mathcal{P}}\left( \mathbf{T}%
\left( L^{k^{\prime }+t-3}\right) \right) \\
&\leq &\frac{1}{\rho ^{2}}\sum_{\substack{ H^{k+2}\in \mathcal{H}_{k+2}  \\ %
H^{k+2}\subset K}}\sum_{\substack{ L^{k^{\prime }+t-4}\in \mathcal{G}%
_{k^{\prime }+t-4}:\ L^{k^{\prime }+t-4}\subset H^{k+2}  \\ \text{where }%
k^{\prime }=d_{\limfunc{geom}}\left( H^{k+2}\right) }}\omega _{\flat 
\mathcal{P}}\left( \mathbf{T}\left( L^{k^{\prime }+t-4}\right) \right) \\
&&\vdots \\
&\leq &\frac{1}{\rho ^{t-2}}\sum_{\substack{ H^{k+2}\in \mathcal{H}_{k+2} 
\\ H^{k+2}\subset K}}\sum_{\substack{ L^{k^{\prime }}\in \mathcal{G}%
_{k^{\prime }}:\ L^{k^{\prime }}\subset H^{k+2}  \\ \text{where }k^{\prime
}=d_{\limfunc{geom}}\left( H^{k+2}\right) }}\omega _{\flat \mathcal{P}%
}\left( \mathbf{T}\left( L^{k^{\prime }}\right) \right) .
\end{eqnarray*}%
Since $L^{k^{\prime }}\subset H^{k+2}$ implies $L^{k^{\prime }}=H^{k+2}$, we
now obtain%
\begin{equation*}
\sum_{J\in \Pi _{2}^{K,\limfunc{aug}}\mathcal{P}_{L,t}^{\flat \mathcal{H}%
}}\left\Vert \bigtriangleup _{J}^{\omega ,\mathbf{b}^{\ast }}x\right\Vert
_{L^{2}\left( \omega \right) }^{\spadesuit 2}\leq \frac{1}{\rho ^{t-2}}%
\sum_{H^{k+2}\in \mathcal{H}_{k+2}:\ H^{k+2}\subset K}\omega _{\flat 
\mathcal{P}}\left( \mathbf{T}\left( H^{k+2}\right) \right) \leq \frac{1}{%
\rho ^{t-2}}\omega _{\flat \mathcal{P}}\left( \mathbf{T}\left( K\right)
\right) ,
\end{equation*}%
which completes the proof of (\ref{rem}), and hence that of (\ref{N_L}).
Finally, an application of the Orthogonality Lemma \ref{mut orth}\ proves (%
\ref{S big t 3}).

\bigskip

\textbf{Proof of the first line in (\ref{rest bounds})}: At last we turn to
proving the first line in (\ref{rest bounds}). Recalling that $\mathcal{T}%
\left( L\right) =\mathcal{T}_{\limfunc{left}}\left( L\right) \dot{\cup}%
\mathcal{T}_{\limfunc{right}}\left( L\right) \dot{\cup}\left\{ L\right\} $,
we consider the collection 
\begin{eqnarray*}
&&\mathcal{Q}_{0}^{\flat \mathcal{H}-big}=\dbigcup\limits_{L\in \mathcal{H}}%
\mathcal{P}_{L,0}^{\flat \mathcal{H}-big}, \\
\text{where } &&\mathcal{P}_{L,0}^{\flat \mathcal{H}-big}=\left\{ \left(
I,J\right) \in \mathcal{P}_{L,0}^{\flat \mathcal{H}}:\text{there is }%
L^{\prime }\in \mathcal{T}\left( L\right) \text{ with }J^{\flat }\subset
L^{\prime }\subset I\right\} ,\ \ \ L\in \mathcal{H}, \\
\text{and } &&\mathcal{P}_{L,0}^{\flat \mathcal{H}}=\left\{ \left(
I,J\right) \in \mathcal{P}:I\in \mathcal{C}_{L}^{\mathcal{H}}\text{ and }%
J\in \mathcal{C}_{L}^{\mathcal{H},\flat \limfunc{shift}}\text{ for some }%
L\in \mathcal{H}\right\} ,\ \ \ L\in \mathcal{H},
\end{eqnarray*}%
and begin by claiming that%
\begin{equation}
\widehat{\mathfrak{N}}_{\limfunc{stop},\bigtriangleup ^{\omega }}^{A,%
\mathcal{P}_{L,0}^{\flat \mathcal{H}-big}}\leq C\mathcal{S}_{\limfunc{aug}%
\limfunc{size}}^{\alpha ,A}\left( \mathcal{P}_{L,0}^{\flat \mathcal{H}%
-big}\right) \leq C\mathcal{S}_{\limfunc{aug}\limfunc{size}}^{\alpha
,A}\left( \mathcal{P}\right) ,\ \ \ \ \ L\in \mathcal{H}.  \label{big t 3}
\end{equation}%
To see this, we fix $L\in \mathcal{H}$ and order the `left' tower of
intervals $\mathcal{T}_{\limfunc{left}}\left( L\right) =\left\{
L^{k}\right\} _{k=1}^{\infty }$ that lie in the restricted corona $\mathcal{C%
}_{L}^{\mathcal{H},\limfunc{restricted}}$ by decreasing side length, i.e. $%
\ell \left( L^{k+1}\right) \leq \ell \left( L^{k}\right) $ for all $k\geq 1$%
, and set $L^{0}=L$ (of course the tower may be finite, but for convenience
in notation, we won't reflect this in the notation). The `right' tower of
intervals $\mathcal{T}_{\limfunc{right}}\left( L\right) $ is handled
similarly and so not considered further here. Then $\mathcal{P}_{L,0}^{\flat 
\mathcal{H}-big}$ can be decomposed as follows, remembering that $J^{\flat
}\subset I\subset L$ for $\left( I,J\right) \in \mathcal{P}_{L,0}^{\flat 
\mathcal{H}-big}\subset \mathcal{P}_{L,0}^{\flat \mathcal{H}}$:%
\begin{eqnarray*}
\mathcal{P}_{L,0}^{\flat \mathcal{H}-big} &=&\overset{\cdot }{\bigcup }%
_{k=1}^{\infty }\left\{ \mathcal{R}_{L_{\limfunc{left}}^{k}}^{\flat \mathcal{%
L}}\dot{\cup}\mathcal{R}_{L_{\limfunc{right}}^{k}}^{\flat \mathcal{L}}\dot{%
\cup}\mathcal{R}_{L_{\limfunc{disjoint}}^{k}}^{\flat \mathcal{L}}\right\} \\
&=&\left( \overset{\cdot }{\bigcup }_{k=1}^{\infty }\mathcal{R}_{L_{\limfunc{%
left}}^{k}}^{\flat \mathcal{L}}\right) \dot{\cup}\left( \overset{\cdot }{%
\bigcup }_{k=1}^{\infty }\mathcal{R}_{L_{\limfunc{right}}^{k}}^{\flat 
\mathcal{L}}\right) \dot{\cup}\left( \overset{\cdot }{\bigcup }%
_{k=1}^{\infty }\mathcal{R}_{L_{\limfunc{disjoint}}^{k}}^{\flat \mathcal{L}%
}\right) \ ; \\
\mathcal{R}_{L_{\limfunc{right}}^{k}}^{\flat \mathcal{L}} &\equiv &\left\{
\left( I,J\right) \in \mathcal{P}_{L,0}^{\flat \mathcal{H}-big}:I\in 
\mathcal{C}_{L^{k-1}}^{\mathcal{L}}\text{ and }J^{\flat }\subset L_{\limfunc{%
right}}^{k}\right\} , \\
\mathcal{R}_{L_{\limfunc{left}}^{k}}^{\flat \mathcal{L}} &\equiv &\left\{
\left( I,J\right) \in \mathcal{P}_{L,0}^{\flat \mathcal{H}-big}:I\in 
\mathcal{C}_{L^{k-1}}^{\mathcal{L}}\text{ and }J^{\flat }\subset L_{\limfunc{%
left}}^{k}\right\} \\
\mathcal{R}_{L_{\limfunc{disjoint}}^{k}}^{\flat \mathcal{L}} &\equiv
&\left\{ \left( I,J\right) \in \mathcal{P}_{L,0}^{\flat \mathcal{H}%
-big}:I\in \mathcal{C}_{L^{k-1}}^{\mathcal{L}}\text{ and }J^{\flat }\in 
\mathcal{C}_{L^{k-1}}^{\mathcal{L}}\text{ and }J^{\flat }\cap
L^{k}=\emptyset \right\} \\
&=&\left\{ \left( I,J\right) \in \mathcal{P}_{L,0}^{\flat \mathcal{H}%
-big}:I=L^{k-1}\text{ and }J^{\flat }\in \mathcal{C}_{L^{k-1}}^{\mathcal{L}}%
\text{ and }J^{\flat }\cap L^{k}=\emptyset \right\} ,
\end{eqnarray*}%
and where in the last line we have used the fact that if $I,J^{\flat }\in 
\mathcal{C}_{L^{k-1}}^{\mathcal{L}}$and there is $L^{\prime }\in \mathcal{T}%
\left( L\right) $ with $J^{\flat }\subset L^{\prime }\subset I$, then we
must have $I=L^{k-1}$. All of the pairs $\left( I,J\right) \in \mathcal{P}%
_{L,0}^{\flat \mathcal{H}-big}$ are included in either $\mathcal{R}_{L_{%
\limfunc{right}}^{k}}^{\flat \mathcal{L}}$, $\mathcal{R}_{L_{\limfunc{left}%
}^{k}}^{\flat \mathcal{L}}$ or $\mathcal{R}_{L_{\limfunc{disjoint}%
}^{k}}^{\flat \mathcal{L}}$ for some $k$, since if $J^{\flat }\supset L^{k}$%
, then $J^{\flat }$ shares an endpoint with $L$, which contradicts the fact
that $3J^{\flat }\subset J^{\maltese }\subset I\subset L$.

We can easily deal with the `disjoint' collection $\mathcal{Q}^{\limfunc{%
disjoint}}\equiv \overset{\cdot }{\bigcup }_{k=1}^{\infty }\mathcal{R}_{L_{%
\limfunc{disjoint}}^{k}}^{\flat \mathcal{L}}$ by applying a $\emph{trivial}$
case of the $\flat $Straddling Lemma to $\mathcal{R}_{L_{\limfunc{disjoint}%
}^{k}}^{\flat \mathcal{L}}$ with a single straddling interval, followed by
an application of the Orthogonality Lemma to $\mathcal{Q}^{\limfunc{disjoint}%
}$. More precisely, every pair $\left( I,J\right) \in \mathcal{R}_{L_{%
\limfunc{disjoint}}^{k}}^{\flat \mathcal{L}}$ satisfies $J^{\flat }\subset
L^{k-1}=I$, so that the reduced admissible collection $\mathcal{R}_{L_{%
\limfunc{disjoint}}^{k}}^{\flat \mathcal{L}}$ $\flat $straddles the trivial
choice $\mathcal{S}=\left\{ L^{k-1}\right\} $, the singleton consisting of
just the interval $L^{k-1}$. Then the inequality%
\begin{equation*}
\widehat{\mathfrak{N}}_{\limfunc{stop},\bigtriangleup ^{\omega }}^{A,%
\mathcal{R}_{L_{\limfunc{disjoint}}^{k}}^{\flat \mathcal{L}}}\leq C\mathcal{S%
}_{\limfunc{aug}\limfunc{size}}^{\alpha ,A}\left( \mathcal{R}_{L_{\limfunc{%
disjoint}}^{k}}^{\flat \mathcal{L}}\right) ,
\end{equation*}%
follows from $\flat $Straddling Lemma \ref{straddle 3 ref}. The collection $%
\left\{ \mathcal{R}_{L_{\limfunc{disjoint}}^{k}}^{\flat \mathcal{L}}\right\}
_{k=1}^{\infty }$ is mutually orthogonal since%
\begin{eqnarray*}
\mathcal{R}_{L_{\limfunc{disjoint}}^{k}}^{\flat \mathcal{L}} &\subset &%
\mathcal{C}_{L^{k-1}}^{\mathcal{L}}\times \mathcal{C}_{L^{k-1}}^{\mathcal{L}%
,\flat \func{shift}}\ , \\
\dsum\limits_{k=1}^{\infty }\mathbf{1}_{\mathcal{C}_{L^{k-1}}^{\mathcal{L}}}
&\leq &\mathbf{1}\text{ and }\dsum\limits_{k=1}^{\infty }\mathbf{1}_{%
\mathcal{C}_{L^{k-1}}^{\mathcal{L},\flat \limfunc{shift}}}\leq \mathbf{1}.
\end{eqnarray*}%
Since $\overset{\cdot }{\bigcup }_{k=1}^{\infty }\mathcal{R}_{L_{\limfunc{%
disjoint}}^{k}}^{\flat \mathcal{L}}$ is reduced and admissible (each $J\in
\Pi _{2}\left( \overset{\cdot }{\bigcup }_{k=1}^{\infty }\mathcal{R}_{L_{%
\limfunc{disjoint}}^{k}}^{\flat \mathcal{L}}\right) $ is paired with a
single $I$, namely the top of the $\mathcal{L}$-corona to which $J^{\flat }$
belongs), the Orthogonality Lemma \ref{mut orth} applies to obtain the
estimate%
\begin{equation}
\widehat{\mathfrak{N}}_{\limfunc{stop},\bigtriangleup ^{\omega }}^{A,\overset%
{\cdot }{\bigcup }_{k=1}^{\infty }\mathcal{R}_{L_{\limfunc{disjoint}%
}^{k}}^{\flat \mathcal{L}}}=\widehat{\mathfrak{N}}_{\limfunc{stop}%
,\bigtriangleup ^{\omega }}^{A,\mathcal{Q}^{\limfunc{disjoint}}}\leq
\sup_{k\geq 1}\widehat{\mathfrak{N}}_{\limfunc{stop},\bigtriangleup ^{\omega
}}^{A,\mathcal{R}_{L_{\limfunc{disjoint}}^{k}}^{\flat \mathcal{L}}}\leq
C\sup_{k\geq 1}\mathcal{S}_{\limfunc{aug}\limfunc{size}}^{\alpha ,A}\left( 
\mathcal{R}_{L_{\limfunc{disjoint}}^{k}}^{\flat \mathcal{L}}\right) \leq C%
\mathcal{S}_{\limfunc{aug}\limfunc{size}}^{\alpha ,A}\left( \mathcal{P}%
_{L,0}^{\flat \mathcal{H}-big}\right) .  \label{disjoint bound}
\end{equation}

Now we turn to estimating the norm of the `right' collection $\mathcal{Q}^{%
\limfunc{right}}\equiv \bigcup_{k=1}^{\infty }\mathcal{R}_{L_{\limfunc{right}%
}^{k}}^{\flat \mathcal{L}}$. First we note that $L_{\limfunc{right}}^{k}\in 
\mathcal{C}_{A}^{\mathcal{A},\limfunc{restrict}}$ if $\left( I,J\right) \in 
\mathcal{R}_{L_{\limfunc{right}}^{k}}^{\flat \mathcal{L}}$ since $\mathcal{R}%
_{L_{\limfunc{right}}^{k}}^{\flat \mathcal{L}}$ is reduced, i.e. doesn't
contain any pairs $\left( I,J\right) $ with $J^{\flat }\subset A^{\prime }$
for some $A^{\prime }\in \mathfrak{C}_{\mathcal{A}}\left( A\right) $. Next
we note that $\mathcal{Q}^{\limfunc{right}}$ is admissible since if $J\in
\Pi _{2}\mathcal{Q}^{\limfunc{right}}$, then $J\in \Pi _{2}\mathcal{R}_{L_{%
\limfunc{right}}^{k}}^{\flat \mathcal{L}}$ for a unique index $k$, and of
course $\mathcal{R}_{L_{\limfunc{right}}^{k}}^{\flat \mathcal{L}}$ is
admissible, so that the intervals $I$ that are paired with $J$ are
tree-connected. Thus we can apply the Straddling Lemma \ref{straddle 3 ref}
to the reduced admissible collection $\mathcal{Q}^{\limfunc{right}}$, with
the `straddling' set $\mathcal{S}\equiv \left\{ L_{\limfunc{right}%
}^{k}\right\} _{k=1}^{\infty }\cap \mathcal{C}_{A}^{\mathcal{A},\limfunc{%
restrict}}$, to obtain the estimate%
\begin{equation}
\widehat{\mathfrak{N}}_{\limfunc{stop},\bigtriangleup ^{\omega
}}^{A,\bigcup_{k=1}^{\infty }\mathcal{R}_{L_{\limfunc{right}}^{k}}^{\flat 
\mathcal{L}}}=\widehat{\mathfrak{N}}_{\limfunc{stop},\bigtriangleup ^{\omega
}}^{A,\mathcal{Q}^{\limfunc{right}}}\leq C\mathcal{S}_{\limfunc{aug}\limfunc{%
size}}^{\alpha ,A}\left( \mathcal{Q}^{\limfunc{right}}\right) \leq C\mathcal{%
S}_{\limfunc{aug}\limfunc{size}}^{\alpha ,A}\left( \mathcal{P}_{L,0}^{\flat 
\mathcal{H}-big}\right) \ .  \label{right bound}
\end{equation}

As for the remaining `left' form $\left\vert \mathsf{B}\right\vert _{%
\limfunc{stop},\bigtriangleup ^{\omega }}^{A,\bigcup_{k=0}^{\infty }\mathcal{%
R}_{L_{\limfunc{left}}^{k}}^{\flat \mathcal{L}}}\left( f,g\right) $, if the
interval pair $\left( I,J\right) \in \mathcal{R}_{L_{\limfunc{left}%
}^{k}}^{\flat \mathcal{L}}$, then either $J^{\flat }\subset L_{\limfunc{left}%
}^{k}\subsetneqq J^{\maltese }$ or $J^{\maltese }\subset L_{\limfunc{left}%
}^{k}$. But $J^{\flat }\subset L_{\limfunc{left}}^{k}\subsetneqq J^{\maltese
}$ implies that either $J^{\flat }=L_{\limfunc{left}}^{k}\subsetneqq
J^{\maltese }\subset I\subset L$, which is impossible since $J^{\flat }$
cannot share an endpoint with $L$, or that $J^{\flat }=L_{-/+}^{k}=\left( L_{%
\limfunc{left}}^{k}\right) _{\limfunc{right}}\ $and $J^{\maltese }=L^{k}$.
So we conclude that if $\left( I,J\right) \in \mathcal{R}_{L_{\limfunc{left}%
}^{k}}^{\flat \mathcal{L}}$, then 
\begin{equation}
\text{either }J^{\maltese }\subset L_{\limfunc{left}}^{k}\text{ or "}%
J^{\maltese }=L^{k}\text{ and }J\subset L_{\limfunc{left}}^{k}\text{"}.
\label{either or}
\end{equation}

In either case in (\ref{either or}), there is a unique interval $K=K\left[ J%
\right] \in \mathcal{W}\left( L\right) $ that contains $J$. It follows that
there are now two remaining cases:

\textbf{Case 1}: $K\left[ J\right] \in \mathcal{C}_{A}^{\limfunc{restrict}}$,

\textbf{Case 2}: $K\left[ J\right] \subset A^{\prime }\subsetneqq I$ for
some $A^{\prime }\in \mathfrak{C}_{\mathcal{A}}\left( A\right) $.

However, by Key Fact \#2 in (\ref{indentation}), and the fact that the \emph{%
great} grandparent $\pi _{\mathcal{D}}^{\left( 3\right) }A^{\prime }$ of $%
A^{\prime }$ contains $A^{\prime }$ inside its leftmost grandchild $\left(
\pi _{\mathcal{D}}^{\left( 3\right) }A^{\prime }\right) _{-/-}$, the pairs $%
\left( I,J\right) $ in \textbf{Case 2} lie in the `corona straddling'
collection $P_{\func{cor}}^{A}$ that was removed from all $A$-admissible
collections in (\ref{empty assumption}) of Conclusion \ref{assume}\ above,
and thus there are no pairs in \textbf{Case 2} here. (We note in passing
that a given $A^{\prime }\in \mathfrak{C}_{A}\left( A\right) $ can occur as
one of the Whitney intervals $K$ in $\mathcal{W}\left( L\right) $ for at
most one $L\in \mathcal{H}$, the indented corona.) Thus we conclude that $K%
\left[ J\right] \in \mathcal{C}_{A}^{\limfunc{restrict}}$.

We now claim that $3K\left[ J\right] \subset I$ for all pairs $\left(
I,J\right) \in \bigcup_{k=1}^{\infty }\mathcal{R}_{L_{\limfunc{left}%
}^{k}}^{\flat \mathcal{L}}$. To see this, suppose that $\left( I,J\right)
\in \mathcal{R}_{L_{\limfunc{left}}^{k}}^{\flat \mathcal{L}}$ for some $%
k\geq 1$. Then by (\ref{either or}) we have both that $K\left[ J\right]
\subset L_{\limfunc{left}}^{k}$ and $L^{k}\subsetneqq I$. But then $K\left[ J%
\right] \subset L_{\limfunc{left}}^{k}$ implies that $3K\left[ J\right]
\subset L^{k}\subset I$ as claimed. See Figure \ref{tow}.

\FRAME{ftbpFU}{6.714in}{2.7422in}{0pt}{\Qcb{The case when $\left( I,J\right)
\in \mathcal{R}_{L_{\limfunc{right}}^{k}}^{\mathcal{L}}$ and $J^{\maltese }$
shares an endpoint with $L$ and $K\in \mathcal{W}\left( L\right) $ equals $%
J_{-/+}^{\maltese }$. In this case $3K\subset I$.}}{\Qlb{tow}}{tower.wmf}{%
\special{language "Scientific Word";type "GRAPHIC";maintain-aspect-ratio
TRUE;display "USEDEF";valid_file "F";width 6.714in;height 2.7422in;depth
0pt;original-width 7.4039in;original-height 9.804in;cropleft "0";croptop
"0.6397";cropright "0.9205";cropbottom "0.3600";filename
'Tower.wmf';file-properties "XNPEU";}}

Now the `left' collection $\mathcal{Q}^{\limfunc{left}}\equiv
\bigcup_{k=1}^{\infty }\mathcal{R}_{L_{\limfunc{left}}^{k}}^{\flat \mathcal{L%
}}$ is admissible, since if $J\in \Pi _{2}\mathcal{Q}^{\limfunc{left}}$ and $%
I_{j}\in \Pi _{1}\mathcal{Q}^{\limfunc{left}}$ with $\left( I_{j},J\right)
\in \mathcal{Q}^{\limfunc{left}}$ for $j=1,2$, then $I_{j}\in \mathcal{C}%
_{L^{k_{j}-1}}^{\mathcal{L}}$ for some $k_{j}$ and all of the intervals $%
I\in \left[ I_{1},I_{2}\right] \ $lie in one of the coronas $\mathcal{C}%
_{L^{k-1}}^{\mathcal{L}}$ for $k$ between $k_{1}$ and $k_{2}$. Thus $\left(
I,J\right) \in \mathcal{R}_{L_{\limfunc{left}}^{k}}^{\flat \mathcal{L}%
}\subset \mathcal{Q}^{\limfunc{left}}$, and we have proved the required
connectedness. From the containment $3K\left[ J\right] \subset I\subset L$
for all $\left( I,J\right) \in \bigcup_{k=1}^{\infty }\mathcal{R}_{L_{%
\limfunc{left}}^{k}}^{\flat \mathcal{L}}$, we now see that the reduced
admissible collection $\mathcal{Q}^{\limfunc{left}}$ \emph{substraddles} the
interval $L$. Hence the Substraddling Lemma \ref{substraddle ref} yields the
bound%
\begin{equation}
\widehat{\mathfrak{N}}_{\limfunc{stop},\bigtriangleup ^{\omega
}}^{A,\bigcup_{k=1}^{\infty }\mathcal{R}_{L_{\limfunc{left}}^{k}}^{\flat 
\mathcal{L}}}=\widehat{\mathfrak{N}}_{\limfunc{stop},\bigtriangleup ^{\omega
}}^{A,\mathcal{Q}^{\limfunc{left}}}\leq C\mathcal{S}_{\limfunc{aug}\limfunc{%
size}}^{\alpha ,A}\left( \mathcal{Q}^{\limfunc{left}}\right) \leq C\mathcal{S%
}_{\limfunc{aug}\limfunc{size}}^{\alpha ,A}\left( \mathcal{P}_{L,0}^{\flat 
\mathcal{H}-big}\right) .  \label{left bound}
\end{equation}%
Combining the bounds (\ref{disjoint bound}), (\ref{right bound}) and (\ref%
{left bound}), we obtain (\ref{big t 3}).

Finally, we observe that the collections $\mathcal{P}_{L,0}^{\flat \mathcal{H%
}-big}$ themselves are \emph{mutually orthogonal}, namely 
\begin{eqnarray*}
\mathcal{P}_{L,0}^{\flat \mathcal{H}-big} &\subset &\mathcal{C}_{L}^{%
\mathcal{H}}\times \mathcal{C}_{L}^{\mathcal{H},\flat \limfunc{shift}}\ ,\ \
\ \ \ L\in \mathcal{H}\ , \\
\dsum\limits_{L\in \mathcal{H}}\mathbf{1}_{\mathcal{C}_{L}^{\mathcal{H}}}
&\leq &\mathbf{1}\text{ and }\dsum\limits_{L\in \mathcal{H}}\mathbf{1}_{%
\mathcal{C}_{L}^{\mathcal{H},\flat \limfunc{shift}}}\leq \mathbf{1}.
\end{eqnarray*}%
Thus an application of the Orthogonality Lemma \ref{mut orth} shows that%
\begin{equation*}
\widehat{\mathfrak{N}}_{\limfunc{stop},\bigtriangleup ^{\omega }}^{A,%
\mathcal{Q}_{0}^{\flat \mathcal{H}-big}}\leq \sup_{L\in \mathcal{L}}\widehat{%
\mathfrak{N}}_{\limfunc{stop},\bigtriangleup ^{\omega }}^{A,\mathcal{P}%
_{L,0}^{\flat \mathcal{H}-big}}\leq C\mathcal{S}_{\limfunc{aug}\limfunc{size}%
}^{\alpha ,A}\left( \mathcal{P}\right) .
\end{equation*}%
Altogether, the proof of Proposition \ref{bottom up 3} is now complete.
\end{proof}

This finishes the proofs of the inequalities (\ref{First inequality}) and (%
\ref{B stop form 3}).

\section{Wrapup of the proof\label{Sub wrapup}}

At this point we have controlled, either directly or probabilistically, the
norms of all of the forms in our decompositions - namely the disjoint,
nearby, far below, paraproduct, neighbour, broken and stopping forms - in
terms of the Muckenhoupt, energy and \emph{functional energy} conditions,
along with an arbitrarily small multiple of the operator norm. Thus it only
remains to control the functional energy condition by the Muckenhoupt and
energy conditions, since then, using $\int \left( T_{\sigma }^{\alpha
}f\right) gd\omega =\Theta \left( f,g\right) +\Theta ^{\ast }\left(
f,g\right) $ with the further decompositions above, we will have shown that
for any fixed tangent line truncation of the operator $T_{\sigma }^{\alpha }$%
, as defined in Definition \ref{truncated op}, we have 
\begin{eqnarray*}
&&\left\vert \int \left( T_{\sigma }^{\alpha }f\right) gd\omega \right\vert =%
\boldsymbol{E}_{\Omega }^{\mathcal{D}}\boldsymbol{E}_{\Omega }^{\mathcal{G}%
}\left\vert \int \left( T_{\sigma }^{\alpha }f\right) gd\omega \right\vert
\leq \boldsymbol{E}_{\Omega }^{\mathcal{D}}\boldsymbol{E}_{\Omega }^{%
\mathcal{G}}\sum_{i=1}^{3}\left( \left\vert \Theta _{i}\left( f,g\right)
\right\vert +\left\vert \Theta _{i}^{\ast }\left( f,g\right) \right\vert
\right) \\
&&\ \ \ \ \ \leq \left( C_{\eta }\mathcal{NTV}_{\alpha }+\eta \mathfrak{N}%
_{T^{\alpha }}\right) \left\Vert f\right\Vert _{L^{2}\left( \sigma \right)
}\left\Vert g\right\Vert _{L^{2}\left( \omega \right) },\ \ \ f\in
L^{2}\left( \sigma \right) \text{ and }g\in L^{2}\left( \omega \right) ,
\end{eqnarray*}%
for an arbitarily small positive constant $\eta >0$, and a correspondingly
large finite constant $C_{\eta }$. Note that the testing constants $%
\mathfrak{T}_{T^{\alpha }}$ and $\mathfrak{T}_{T^{\alpha ,\ast }}$ in $%
\mathcal{NTV}_{\alpha }$ already include the supremum over all tangent line
truncations of $T^{\alpha }$, while the operator norm $\mathfrak{N}%
_{T^{\alpha }}$ on the left refers to a \emph{fixed} tangent line truncation
of $T^{\alpha }$. This gives%
\begin{equation*}
\mathfrak{N}_{T^{\alpha }}=\sup_{\left\Vert f\right\Vert _{L^{2}(\sigma
)}=1}\sup_{\left\Vert g\right\Vert _{L^{2}(\omega )}=1}\left\vert \int
\left( T_{\sigma }^{\alpha }f\right) gd\omega \right\vert \leq C_{\eta }%
\mathcal{NTV}_{\alpha }+\eta \mathfrak{N}_{T^{\alpha }},
\end{equation*}%
and since the truncated operators have finite operator norm $\mathfrak{N}%
_{T^{\alpha }}$, we can absorb the term $\eta \mathfrak{N}_{T^{\alpha }}$
into the left hand side for $\eta <1$ and obtain $\mathfrak{N}_{T^{\alpha
}}\leq C_{\eta }^{\prime }\mathcal{NTV}_{\alpha }$ for each tangent line
truncation of $T^{\alpha }$. Taking the supremum over all such truncations
of $T^{\alpha }$ finishes the proof of Theorem \ref{dim one}.

The task of controlling functional energy is taken up in Appendix B, after
first establishing weak frame and weak Riesz inequalities for martingale and
dual martingale differences in Appendix A (except for the lower weak Riesz
inequality for the martingale difference $\bigtriangleup _{Q}^{\mu ,\mathbf{b%
}}$).

\section{Appendix A: Martingale differences}

Most of the material in this appendix is known, see e.g. \cite{NTV3} and 
\cite{HyMa}. First, we recall the construction in \cite{SaShUr9} of a Haar
basis in $\mathbb{R}$ that is adapted to a measure $\mu $ (c.f. \cite{NTV2}
where this type of construction is made explicit). Given a dyadic interval $%
Q\in \mathcal{D}$, where $\mathcal{D}$ is a dyadic grid of intervals from $%
\mathcal{P}$, let $\bigtriangleup _{Q}^{\mu }$ denote orthogonal projection
onto the \textbf{one}-dimensional subspace $L_{Q}^{2}\left( \mu \right) $ of 
$L^{2}\left( \mu \right) $ that consists of linear combinations of the
indicators of\ the children $\mathfrak{C}\left( Q\right) $ of $Q$ that have $%
\mu $-mean zero over $Q$:%
\begin{equation*}
L_{Q}^{2}\left( \mu \right) \equiv \left\{ f=\dsum\limits_{Q^{\prime }\in 
\mathfrak{C}\left( Q\right) }a_{Q^{\prime }}\mathbf{1}_{Q^{\prime
}}:a_{Q^{\prime }}\in \mathbb{R},\int_{Q}fd\mu =0\right\} .
\end{equation*}%
Then we have the important telescoping property for dyadic intervals $%
Q_{1}\subset Q_{2}$:%
\begin{equation}
\mathbf{1}_{Q_{0}}\left( x\right) \left( \dsum\limits_{Q\in \left[
Q_{1},Q_{2}\right] }\bigtriangleup _{Q}^{\mu }f\left( x\right) \right) =%
\mathbf{1}_{Q_{0}}\left( x\right) \left( \mathbb{E}_{Q_{0}}^{\mu }f-\mathbb{E%
}_{Q_{2}}^{\mu }f\right) ,\ \ \ \ \ Q_{0}\in \mathfrak{C}\left( Q_{1}\right)
,\ f\in L^{2}\left( \mu \right) ,  \label{telescope}
\end{equation}

\begin{notation}
Here $\mathbb{E}_{Q}^{\mu }f\left( x\right) \equiv \mathbf{1}_{Q}\left(
x\right) \frac{1}{\left\vert Q\right\vert _{\mu }}\int_{Q}fd\mu $ denotes
the projection of $f$ onto $\limfunc{Span}\left\{ \mathbf{1}_{Q}\right\} $,
the one-dimensional subspace of multiples of the indicator of $Q$. We will
also denote the average value itself by $E_{Q}^{\mu }f\equiv \frac{1}{%
\left\vert Q\right\vert _{\mu }}\int_{Q}fd\mu $.
\end{notation}

It is convenient at times to use a fixed normalized basis $\left\{
h_{Q}^{\mu }\right\} $ of $L_{Q}^{2}\left( \mu \right) $. Then $\left\{
h_{Q}^{\mu }\right\} _{Q\in \mathcal{D}}$ is an orthonormal basis for $%
L^{2}\left( \mu \right) $, with the understanding that we add the constant
function $\mathbf{1}$ if $\mu $ is a finite measure. In particular we have%
\begin{equation*}
\left\Vert f\right\Vert _{L^{2}\left( \mu \right) }^{2}=\sum_{Q\in \mathcal{D%
}}\left\Vert \bigtriangleup _{Q}^{\mu }f\right\Vert _{L^{2}\left( \mu
\right) }^{2}=\sum_{Q\in \mathcal{D}}\left\vert \widehat{f}\left( Q\right)
\right\vert ^{2},\ \ \ \ \ \widehat{f}\left( Q\right) \equiv \left\langle
f,h_{Q}^{\mu }\right\rangle _{\mu }\ ,
\end{equation*}%
where the measure is suppressed in the notation $\widehat{f}$. Indeed, this
follows from (\ref{telescope}) and Lebesgue's differentiation theorem for
dyadic intervals. We also record the following useful estimate. If $%
I^{\prime }$ is either of the two $\mathcal{D}$-children of $I$, then 
\begin{equation}
\left\vert \mathbb{E}_{I^{\prime }}^{\mu }h_{I}^{\mu }\right\vert \leq \sqrt{%
\mathbb{E}_{I^{\prime }}^{\mu }\left( h_{I}^{\mu }\right) ^{2}}\leq \frac{1}{%
\sqrt{\left\vert I^{\prime }\right\vert _{\mu }}}.  \label{useful Haar}
\end{equation}

In the next subsection, we introduce martingale and dual martingale
differences for various $p$-weakly $\mu $-accretive families $\mathbf{b}%
=\left\{ b_{Q}\right\} _{Q\in \mathcal{P}}$, and establish convergence
properties for their expansions. Then in later subsections, we turn to frame
inequalities, and weak Riesz inequalities.

\subsection{Convergence for weakly and controlled accretive families}

Supposes that $\mathbf{b}=\left\{ b_{Q}\right\} _{Q\in \mathcal{P}}$ is a $p$%
-weakly $\mu $-accretive family on $\mathbb{R}$. Define the $\mathbf{b}$%
-expectation operator $\mathbb{E}_{Q}^{\mu ,\mathbf{b}}$ and the dual $%
\mathbf{b}$-expectation operator $\mathbb{F}_{Q}^{\mu ,\mathbf{b}}$ by 
\begin{eqnarray}
\mathbb{E}_{Q}^{\mu ,\mathbf{b}}f\left( x\right) &\equiv &\mathbf{1}%
_{Q}\left( x\right) \frac{1}{\int_{Q}b_{Q}d\mu }\int_{Q}fb_{Q}d\mu ,\ \ \ \
\ Q\in \mathcal{P}\ ,  \label{def expectation} \\
\mathbb{F}_{Q}^{\mu ,\mathbf{b}}f\left( x\right) &\equiv &\mathbf{1}%
_{Q}\left( x\right) b_{Q}\left( x\right) \frac{1}{\int_{Q}b_{Q}d\mu }%
\int_{Q}fd\mu ,\ \ \ \ \ Q\in \mathcal{P}\ .  \notag
\end{eqnarray}%
Occasionally we will use the modification of $\mathbb{F}_{Q}^{\mu ,\mathbf{b}%
}$ given by `dividing out' the factor $b_{Q}$:%
\begin{equation}
\widehat{\mathbb{F}}_{Q}^{\mu ,\mathbf{b}}f\left( x\right) \equiv \mathbf{1}%
_{Q}\left( x\right) \frac{1}{\int_{Q}b_{Q}d\mu }\int_{Q}fd\mu ,\ \ \ \ \
Q\in \mathcal{P}\ .  \label{F hat}
\end{equation}%
Then define the corresponding martingale and dual martingale differences by%
\begin{eqnarray}
\bigtriangleup _{Q}^{\mu ,\mathbf{b}}f\left( x\right) &\equiv &\left(
\sum_{Q^{\prime }\in \mathfrak{C}\left( Q\right) }\mathbb{E}_{Q^{\prime
}}^{\mu ,\mathbf{b}}f\left( x\right) \right) -\mathbb{E}_{Q}^{\mu ,\mathbf{b}%
}f\left( x\right) =\sum_{Q^{\prime }\in \mathfrak{C}\left( Q\right) }\mathbf{%
1}_{Q^{\prime }}\left( x\right) \left( \mathbb{E}_{Q^{\prime }}^{\mu ,%
\mathbf{b}}f\left( x\right) -\mathbb{E}_{Q}^{\mu ,\mathbf{b}}f\left(
x\right) \right) ,  \label{def diff} \\
\square _{Q}^{\mu ,\mathbf{b}}f\left( x\right) &\equiv &\left(
\sum_{Q^{\prime }\in \mathfrak{C}\left( Q\right) }\mathbb{F}_{Q^{\prime
}}^{\mu ,\mathbf{b}}f\left( x\right) \right) -\mathbb{F}_{Q}^{\mu ,\mathbf{b}%
}f\left( x\right) =\sum_{Q^{\prime }\in \mathfrak{C}\left( Q\right) }\mathbf{%
1}_{Q^{\prime }}\left( x\right) \left( \mathbb{F}_{Q^{\prime }}^{\mu ,%
\mathbf{b}}f\left( x\right) -\mathbb{F}_{Q}^{\mu ,\mathbf{b}}f\left(
x\right) \right) .  \notag
\end{eqnarray}

\begin{description}
\item[Exception] In the special case that $\mathfrak{C}\left( Q\right)
=\left\{ Q_{1},Q_{2}\right\} $ where $\left\vert Q_{1}\right\vert _{\mu }=0$
and $\left\vert Q\right\vert _{\mu }>0$, we set both $\bigtriangleup
_{Q}^{\mu ,\mathbf{b}}\equiv 0$ and $\square _{Q}^{\mu ,\mathbf{b}}\equiv 0$%
, and redefine the parent differences $\bigtriangleup _{\pi Q}^{\mu ,\mathbf{%
b}}$ and $\square _{\pi Q}^{\mu ,\mathbf{b}}$ as follows. Let $\mathfrak{C}%
\left( \pi Q\right) =\left\{ Q,\widetilde{Q}\right\} $ and set 
\begin{eqnarray*}
\bigtriangleup _{\pi Q}^{\mu ,\mathbf{b}}f\left( x\right) &\equiv &\left( 
\mathbb{E}_{Q}^{\mu ,b_{Q_{2}}}f\left( x\right) +\mathbb{E}_{\widetilde{Q}%
}^{\mu ,\mathbf{b}}f\left( x\right) \right) -\mathbb{E}_{\pi Q}^{\mu ,%
\mathbf{b}}f\left( x\right) , \\
\square _{\pi Q}^{\mu ,\mathbf{b}}f\left( x\right) &\equiv &\left( \mathbb{F}%
_{Q}^{\mu ,b_{Q_{2}}}f\left( x\right) +\mathbb{F}_{\widetilde{Q}}^{\mu ,%
\mathbf{b}}f\left( x\right) \right) -\mathbb{F}_{\pi Q}^{\mu ,\mathbf{b}%
}f\left( x\right) ,
\end{eqnarray*}%
where we have used the test function $b_{Q_{2}}$ in place of the expected $%
b_{Q}$ (because $\left\vert Q_{1}\right\vert _{\mu }=0$). With the analogous
modification when only one grandchild at level $k$ below $Q$ is charged by $%
\mu $, the telescoping property holds for these differences, and the reader
can easily verify all of the convergence statements and formulas below. For
the sake of convenience only, we will ignore these exceptions in the sequel,
and proceed under the assumption that all intervals are charged by $\mu $.
\end{description}

Note that in \cite{NTV3} and \cite{HyMa} this notation is reversed - they
use $\bigtriangleup _{Q}^{\mu ,\mathbf{b}}$ for our $\square _{Q}^{\mu ,%
\mathbf{b}}$. Finally, define the dual $\mathbf{b}$-expectation operator $%
\mathbb{F}_{m}^{\mu ,\mathbf{b}}$ on a function $f$ at level $m$ by%
\begin{equation*}
\mathbb{F}_{m}^{\mu ,\mathbf{b}}f\left( x\right) \equiv \sum_{Q\in \mathcal{D%
}_{m}}\mathbb{F}_{Q}^{\mu ,\mathbf{b}}f\left( x\right) =\sum_{Q\in \mathcal{D%
}_{m}}\mathbf{1}_{Q}\left( x\right) b_{Q}\left( x\right) \frac{1}{%
\int_{Q}b_{Q}d\mu }\int_{Q}fd\mu ,
\end{equation*}%
where $\mathcal{D}_{m}\equiv \left\{ Q\in \mathcal{D}:\ell \left( Q\right)
=2^{-m}\right\} $, and define the operators $\square _{m}^{\mu ,\mathbf{b}}$
by%
\begin{equation*}
\square _{m}^{\mu ,\mathbf{b}}\equiv \mathbb{F}_{m}^{\mu ,\mathbf{b}}-%
\mathbb{F}_{m-1}^{\mu ,\mathbf{b}}\ .
\end{equation*}

\begin{definition}
\label{controlled accretive}Let $\mu $ be a locally finite positive Borel
measure on $\mathbb{R}$, let $\mathbf{b}=\left\{ b_{Q}\right\} _{Q\in 
\mathcal{P}}$ be a $p$-weakly $\mu $-accretive family on $\mathbb{R}$ with
reverse H\"{o}lder control (\ref{rev Hol con}) on children, and let $%
\mathcal{A}$ be\ subset of a dyadic grid $\mathcal{D}$. We say that the
subfamily $\mathbf{b}=\left\{ b_{Q}\right\} _{Q\in \mathcal{D}}$ is a $p$%
\emph{-weakly }$\mu $\emph{-controlled accretive} family on $\mathcal{D}\ $%
(note we omit dependence $\mathcal{A}$ on in this notation) if%
\begin{equation*}
b_{Q}=\mathbf{1}_{Q}b_{A}\ ,\ \ \ \ \ Q\in \mathcal{C}_{A}\mathcal{\ },A\in 
\mathcal{A\ },
\end{equation*}%
and the set $\mathcal{A}$ satisfies a Carleson condition:%
\begin{equation*}
\sum_{Q\in \mathcal{A}:\ Q\subset K}\left\vert Q\right\vert _{\mu }\leq
C\left\vert K\right\vert _{\mu }\ ,\ \ \ \ \ \text{for all }K\in \mathcal{D},
\end{equation*}%
equivalently%
\begin{equation*}
\sum_{Q\in \mathcal{A}:\ Q\subset \Omega }\left\vert Q\right\vert _{\mu
}\leq C^{\prime }\left\vert \Omega \right\vert _{\mu }\ ,\ \ \ \ \ \text{for
all open sets }\Omega \subset \mathbb{R}.
\end{equation*}
\end{definition}

Denote the coronas associated to $\mathcal{A}$ by $\left\{ \mathcal{C}%
_{A}\right\} _{A\in \mathcal{A}}$, and now suppose that $\mathbf{b}=\left\{
b_{Q}\right\} _{Q\in \mathcal{P}}$ is $p$-weakly $\mu $-controlled
accretive, i.e. 
\begin{equation*}
0<1\leq \left\vert \frac{1}{\left\vert Q\right\vert _{\mu }}%
\int_{Q}b_{A}d\mu \right\vert \leq \left\Vert b_{A}\right\Vert _{L^{\infty
}\left( \mu \right) }\leq C_{\mathbf{b}}<\infty ,\ \ \ \ \ Q\in \mathcal{C}%
_{A}\mathcal{\ },A\in \mathcal{A}\ ,
\end{equation*}%
where in addition we have reverse H\"{o}lder control (\ref{rev Hol con}) on
children, and $\mathcal{A}$ satisfies a $\mu $-Carleson condition. Now
decompose the children $Q^{\prime }\in \mathfrak{C}\left( Q\right) $ into
the collection $\mathfrak{C}_{\limfunc{broken}}\left( Q\right) $ of \emph{%
broken} children $Q^{\prime }\in \mathcal{A}$ and the remaining collection $%
\mathfrak{C}_{\limfunc{natural}}\left( Q\right) $ of \emph{natural} children 
$Q^{\prime }\not\in \mathcal{A}$, i.e.%
\begin{equation*}
\mathfrak{C}_{\limfunc{broken}}\left( Q\right) \equiv \mathfrak{C}\left(
Q\right) \cap \mathcal{A}\text{ and }\mathfrak{C}_{\limfunc{natural}}\left(
Q\right) \equiv \mathfrak{C}\left( Q\right) \setminus \mathcal{A}.
\end{equation*}%
Let $\left\{ \square _{I}^{\mu ,\mathbf{b}}\right\} _{I\in \mathcal{D}}$ be
associated to the coronas $\left\{ \mathcal{C}_{A}\right\} _{A\in \mathcal{A}%
}$ and the $p$-weakly $\mu $-controlled accretive family $\mathbf{b}=\left\{
b_{A}\right\} _{A\in \mathcal{A}}$ as in Definition \ref{controlled
accretive} above. We will refer to the collection of dual martingale
differences $\left\{ \square _{I}^{\mu ,\mathbf{b}}\right\} _{I\in \mathcal{D%
}}$ as a `broken corona decomposition' in light of the fact that the testing
functions $b_{Q}$ `break' when passing from one corona to another. For $A\in 
\mathcal{A}$, define the corona `pseudoprojections' $\mathsf{P}_{\mathcal{C}%
_{A}}^{\mu ,\mathbf{b}}$ by%
\begin{equation*}
\mathsf{P}_{\mathcal{C}_{A}}^{\mu ,\mathbf{b}}f=\sum_{I\in \mathcal{C}%
_{A}}\square _{I}^{\mu ,\mathbf{b}}f\ ,
\end{equation*}%
We have the \emph{broken corona decomposition},%
\begin{equation*}
f=\sum_{I\in \mathcal{D}}\square _{I}^{\mu ,\mathbf{b}}f=\sum_{A\in \mathcal{%
A}}\mathsf{P}_{\mathcal{C}_{A}}^{\mu ,\mathbf{b}}f\ ,
\end{equation*}%
whose convergence properties we investigate in the next subsubsection.

\subsubsection{Convergence of controlled martingale differences}

As shown by Hyt\"{o}nen and Martikainen \cite{HyMa}, in the setting of a $2$%
-weakly $\mu $-controlled accretive\emph{\ }family, we have strong
convergence in $L^{2}\left( \mu \right) $ for the dual martingale
differences - and also for the martingale differences under the stronger
assumption of a $p$-weakly $\mu $-controlled accretive\emph{\ }family for
some $p>2$ (the proofs given there carry over to general measures). We only
use the case $p=\infty $.

\begin{lemma}
\label{conv prop}(\cite{NTV3}, \cite[Lemma 3.5]{HyMa}) Suppose $\mathbf{b}$
is an $\infty $-weakly $\mu $-controlled accretive family on a grid $%
\mathcal{D}$. Then we have the dual martingale and martingale identities%
\begin{eqnarray*}
f &=&\sum_{I\in \mathcal{D}_{N}}\mathbb{F}_{I}^{\mu ,\mathbf{b}}f+\sum_{I\in 
\mathcal{D}:\ \ell \left( I\right) \geq N+1}\square _{I}^{\mu ,\mathbf{b}}f\
,\ \ \ \ \ N\in \mathbb{Z}, \\
f &=&\sum_{I\in \mathcal{D}_{N}}\mathbb{E}_{I}^{\mu ,\mathbf{b}}f+\sum_{I\in 
\mathcal{D}:\ \ell \left( I\right) \geq N+1}\bigtriangleup _{I}^{\mu ,%
\mathbf{b}}f\ ,\ \ \ \ \ N\in \mathbb{Z},
\end{eqnarray*}%
in the sense of pointwise $\mu $-almost everywhere convergence, and also in
the sense of strong convergence in $L^{2}\left( \mu \right) $.
\end{lemma}

\subsection{Frame inequalities}

Define the positive sublinear operators%
\begin{eqnarray}
\bigtriangledown _{Q}^{\mu }f &\equiv &\sum_{Q^{\prime }\in \mathfrak{C}_{%
\limfunc{broken}}\left( Q\right) }\left( \frac{1}{\left\vert Q^{\prime
}\right\vert _{\mu }}\int_{Q^{\prime }}\left\vert f\right\vert d\mu \right) 
\mathbf{1}_{Q^{\prime }},  \label{Carleson avg op} \\
\widehat{\bigtriangledown }_{Q}^{\mu }f &\equiv &\sum_{Q^{\prime }\in 
\mathfrak{C}_{\limfunc{broken}}\left( Q\right) }\left( \frac{1}{\left\vert
Q^{\prime }\right\vert _{\mu }}\int_{Q^{\prime }}\left\vert f\right\vert
d\mu +\frac{1}{\left\vert Q\right\vert _{\mu }}\int_{Q}\left\vert
f\right\vert d\mu \right) \mathbf{1}_{Q^{\prime }},  \notag
\end{eqnarray}%
where we are suppressing here the dependence of both $\bigtriangledown
_{Q}^{\mu }$ and its larger version $\widehat{\bigtriangledown }_{Q}^{\mu }$
on the breaking intervals. Note also that $\widehat{\bigtriangledown }%
_{Q}^{\mu }=0$ if $Q$ has no broken children. We also set $Q_{\limfunc{broken%
}}\equiv \bigcup_{Q^{\prime }\in \mathfrak{C}_{\limfunc{broken}}\left(
Q\right) }Q^{\prime }$. We now show that the Carleson condition on broken
children gives the inequality%
\begin{equation}
\sum_{Q\in \mathcal{D}}\left\Vert \widehat{\bigtriangledown }_{Q}^{\mu
}f\right\Vert _{L^{2}\left( \mu \right) }^{2}\lesssim \left\Vert
f\right\Vert _{L^{2}\left( \mu \right) }^{2}\ .  \label{Car embed}
\end{equation}%
Indeed,%
\begin{eqnarray*}
\sum_{Q\in \mathcal{D}}\left\Vert \widehat{\bigtriangledown }_{Q}^{\mu
}f\right\Vert _{L^{2}\left( \mu \right) }^{2} &=&\sum_{Q\in \mathcal{D}}\int
\left\vert \sum_{Q^{\prime }\in \mathfrak{C}_{\limfunc{broken}}\left(
Q\right) }\mathbf{1}_{Q^{\prime }}\left( \frac{1}{\left\vert Q^{\prime
}\right\vert _{\mu }}\int_{Q^{\prime }}\left\vert f\right\vert d\mu +\frac{1%
}{\left\vert Q\right\vert _{\mu }}\int_{Q}\left\vert f\right\vert d\mu
\right) \right\vert ^{2}d\mu \\
&\lesssim &\sum_{Q\in \mathcal{D}}\sum_{Q^{\prime }\in \mathfrak{C}_{%
\limfunc{broken}}\left( Q\right) }\left\vert Q^{\prime }\right\vert _{\mu
}\left( \frac{1}{\left\vert Q^{\prime }\right\vert _{\mu }}\int_{Q^{\prime
}}\left\vert f\right\vert d\mu +\frac{1}{\left\vert Q\right\vert _{\mu }}%
\int_{Q}\left\vert f\right\vert d\mu \right) ^{2} \\
&\lesssim &\sum_{A\in \mathcal{A}}\left\vert A\right\vert _{\mu }\left( 
\frac{1}{\left\vert A\right\vert _{\mu }}\int_{A}\left\vert f\right\vert
d\mu \right) ^{2}+\sum_{A\in \mathcal{A}}\left\vert A\right\vert _{\mu
}\left( \frac{1}{\left\vert \pi _{\mathcal{D}}A\right\vert _{\mu }}\int_{\pi
_{\mathcal{D}}A}\left\vert f\right\vert d\mu \right) ^{2},
\end{eqnarray*}%
where $\mathcal{A}$ is a collection of stopping times that satisfy a
Carleson condition as in Definition \ref{controlled accretive}. The first
term in the last line above is at most $C\left\Vert f\right\Vert
_{L^{2}\left( \mu \right) }^{2}$ by the Carleson embedding theorem. We claim
that the second term is controlled by a variant of the Carleson embedding
theorem, 
\begin{equation}
\sum_{A\in \mathcal{A}}\left\vert A\right\vert _{\mu }\left( \frac{1}{%
\left\vert \pi _{\mathcal{D}}A\right\vert _{\mu }}\int_{\pi _{\mathcal{D}%
}A}\left\vert f\right\vert d\mu \right) ^{2}\lesssim \int \left\vert
f\right\vert ^{2}d\mu ,\ \ \ \ \ f\in L^{2}\left( \mu \right) .
\label{parent CET}
\end{equation}%
Indeed, with the measure $\nu \left( A\right) \equiv \left\vert A\right\vert
_{\mu }$ on $\mathcal{A}$, the sublinear map $T$ defined by $Tf\left(
A\right) \equiv \frac{1}{\left\vert \pi _{\mathcal{D}}A\right\vert _{\mu }}%
\int_{\pi _{\mathcal{D}}A}\left\vert f\right\vert d\mu $ takes $L^{\infty
}\left( \mu \right) $ to $L^{\infty }\left( \nu \right) $, and also $%
L^{1}\left( \mu \right) $ to $L^{1,\infty }\left( \nu \right) $, since if $%
\left\{ M\right\} $ are the maximal dyadic intervals $\pi _{\mathcal{D}}A$
such that%
\begin{equation*}
\frac{1}{\left\vert \pi _{\mathcal{D}}A\right\vert _{\mu }}\int_{\pi _{%
\mathcal{D}}A}\left\vert f\right\vert d\mu >\lambda >0,
\end{equation*}%
then%
\begin{eqnarray*}
\left\vert \left\{ A\in \mathcal{A}:Tf\left( A\right) >\lambda \right\}
\right\vert _{\nu } &=&\sum_{A:\ Tf\left( A\right) >\lambda }\left\vert
A\right\vert _{\mu }\leq \sum_{M}\sum_{A:\ A\subset M}\left\vert
A\right\vert _{\mu }\lesssim \sum_{M}\left\vert M\right\vert _{\mu } \\
&<&\sum_{M}\frac{1}{\lambda }\int_{M}\left\vert f\right\vert d\mu \leq \frac{%
1}{\lambda }\int \left\vert f\right\vert d\mu ,
\end{eqnarray*}%
since the maximal intervals $M$ are pairwise disjoint. Now interpolation
shows that $T$ takes $L^{2}\left( \mu \right) $ to $L^{2}\left( \nu \right) $%
, which is (\ref{parent CET}). These inequalities provide the reason for
referring to $\bigtriangledown _{Q}^{\mu }$ and $\widehat{\bigtriangledown }%
_{Q}^{\mu }$\ as \emph{Carleson averaging operators}.

From \cite{NTV3} we have that in the case $\mathbf{b}$ is an $\infty $%
-weakly $\mu $-controlled accretive family on a grid $\mathcal{D}$, and the
underlying measure $\mu $ is upper doubling, then the following frame
equivalences hold:%
\begin{eqnarray*}
\left\Vert f\right\Vert _{L^{2}\left( \mu \right) }^{2} &\approx &\sum_{Q\in 
\mathcal{D}}\left\{ \left\Vert \square _{Q}^{\mu ,\mathbf{b}}f\right\Vert
_{L^{2}\left( \mu \right) }^{2}+\left\Vert \bigtriangledown _{Q}^{\mu
}f\right\Vert _{L^{2}\left( \mu \right) }^{2}\right\} \\
&\approx &\sum_{Q\in \mathcal{D}}\left\{ \left\Vert \bigtriangleup _{Q}^{\mu
,\mathbf{b}}f\right\Vert _{L^{2}\left( \mu \right) }^{2}+\left\Vert
\bigtriangledown _{Q}^{\mu }f\right\Vert _{L^{2}\left( \mu \right)
}^{2}\right\} \ .
\end{eqnarray*}%
It appears however that the arguments in \cite{NTV3} hold for more general
measures\footnote{%
But in \cite{NTV3}, a proof of their inequality (3.3), which is our lower
frame inequality for $\bigtriangleup _{Q}^{\mu ,\mathbf{b}}$, seems not to
be explicitly given.}, and in any case, we will extend most of these frame
inequalities to certain of the weak Riesz inequalities below for arbitrary
positive Borel measures $\mu $. For this it will be convenient to refer
directly to the proofs of the frame inequalities, and since notation in \cite%
{NTV3} is very different from that used in this paper, we will instead
extend a similar argument of Hyt\"{o}nen and Martikainen \cite[Proposition
3.10]{HyMa} to general measures, adding a small additional argument for the
martingale differences in term $III_{A}$ in the proof of Proposition \ref%
{dual frame} below. This uses the following unweighted square function
estimate, which is essentially just the orthogonality of the standard Haar
projections $\bigtriangleup _{I}^{\mu }$ adapted to the measure $\mu $. For
this, we recall the general Haar projections $\mathsf{P}_{\mathcal{C}%
_{A}}^{\mu }h=\sum_{Q\in \mathcal{C}_{A}}\bigtriangleup _{Q}^{\mu }h$.

\begin{lemma}
\label{unweighted square}%
\begin{equation*}
\sum_{A\in \mathcal{A}}\sum_{Q\in \mathcal{C}_{A}}\sum_{Q^{\prime }\in 
\mathfrak{C}_{\limfunc{natural}}\left( Q\right) }\left\vert Q^{\prime
}\right\vert _{\mu }\left\vert E_{Q^{\prime }}^{\mu }h-E_{Q}^{\mu
}h\right\vert ^{2}\lesssim \left\Vert h\right\Vert _{L^{2}\left( \mu \right)
}^{2},\ \ \ \ \ h\in L^{2}\left( \mu \right) .
\end{equation*}
\end{lemma}

\begin{proof}
Recall that the Haar projection $\bigtriangleup _{I}^{\mu }$ is given by 
\begin{equation*}
\bigtriangleup _{I}^{\mu }h=\left( \sum_{I^{\prime }\in \mathfrak{C}\left(
I\right) }\mathbb{E}_{I^{\prime }}^{\mu }h\right) -\mathbb{E}_{I}^{\mu
}h=\sum_{I^{\prime }\in \mathfrak{C}\left( I\right) }\mathbf{1}_{I^{\prime
}}\left( \mathbb{E}_{I^{\prime }}^{\mu }h-\mathbb{E}_{I}^{\mu }h\right) ,
\end{equation*}%
and so%
\begin{equation*}
\sum_{Q^{\prime }\in \mathfrak{C}\left( Q\right) }\left\vert Q^{\prime
}\right\vert _{\mu }\left\vert \mathbb{E}_{Q^{\prime }}^{\mu }h-\mathbb{E}%
_{Q}^{\mu }h\right\vert ^{2}=\int \left\vert \sum_{Q^{\prime }\in \mathfrak{C%
}\left( Q\right) }\mathbf{1}_{Q^{\prime }}\left( \mathbb{E}_{Q^{\prime
}}^{\mu }h-\mathbb{E}_{Q}^{\mu }h\right) \right\vert ^{2}d\mu =\left\Vert
\bigtriangleup _{Q}^{\mu }h\right\Vert _{L^{2}\left( \mu \right) }^{2}.
\end{equation*}%
Thus we have%
\begin{eqnarray*}
\sum_{A\in \mathcal{A}}\sum_{Q\in \mathcal{C}_{A}}\sum_{Q^{\prime }\in 
\mathfrak{C}_{\limfunc{natural}}\left( Q\right) }\left\vert Q^{\prime
}\right\vert _{\mu }\left\vert \mathbb{E}_{Q^{\prime }}^{\mu }h-\mathbb{E}%
_{Q}^{\mu }h\right\vert ^{2} &\leq &\sum_{A\in \mathcal{A}}\sum_{Q\in 
\mathcal{C}_{A}}\sum_{Q^{\prime }\in \mathfrak{C}\left( Q\right) }\left\vert
Q^{\prime }\right\vert _{\mu }\left\vert \mathbb{E}_{Q^{\prime }}^{\mu }h-%
\mathbb{E}_{Q}^{\mu }h\right\vert ^{2} \\
&=&\sum_{A\in \mathcal{A}}\sum_{Q\in \mathcal{C}_{A}}\left\Vert
\bigtriangleup _{Q}^{\mu }h\right\Vert _{L^{2}\left( \mu \right)
}^{2}=\left\Vert h\right\Vert _{L^{2}\left( \mu \right) }^{2}\ .
\end{eqnarray*}
\end{proof}

\begin{proposition}
\label{dual frame}(see \cite[Remark 3.11]{HyMa} for the case of a doubling
measure with $p>2$) Suppose that $\mathbf{b}$ is an $\infty $-weakly $\mu $%
-controlled accretive family on a grid $\mathcal{D}$ with corona tops $%
\mathcal{A\subset D}$. Then we have the lower frame inequality%
\begin{equation*}
\sum_{I\in \mathcal{D}}\left\Vert \bigtriangleup _{Q}^{\mu ,\mathbf{b}%
}f\right\Vert _{L^{2}\left( \mu \right) }^{2}\lesssim \left\Vert
f\right\Vert _{L^{2}\left( \mu \right) }^{2}\ .
\end{equation*}
\end{proposition}

\begin{proof}
Given $A\in \mathcal{A}$, we begin\ with%
\begin{eqnarray*}
\sum_{Q\in \mathcal{C}_{A}}\left\Vert \bigtriangleup _{Q}^{\mu ,\mathbf{b}%
}f\right\Vert _{L^{2}\left( \mu \right) }^{2} &=&\sum_{Q\in \mathcal{C}%
_{A}}\sum_{Q^{\prime }\in \mathfrak{C}\left( Q\right) }\int_{Q^{\prime
}}\left\vert \mathbb{E}_{Q^{\prime }}^{\mu ,\mathbf{b}}f\left( x\right) -%
\mathbb{E}_{Q}^{\mu ,\mathbf{b}}f\left( x\right) \right\vert ^{2}d\mu \left(
x\right) \\
&=&\sum_{Q\in \mathcal{C}_{A}}\sum_{Q^{\prime }\in \mathfrak{C}_{\limfunc{%
natural}}\left( Q\right) }\int_{Q^{\prime }}\left\vert \frac{\int_{Q^{\prime
}}fb_{A}d\mu }{\int_{Q^{\prime }}b_{A}d\mu }-\frac{\int_{Q}fb_{A}d\mu }{%
\int_{Q}b_{A}d\mu }\right\vert ^{2}d\mu \left( x\right) \\
&&+\sum_{Q\in \mathcal{C}_{A}}\sum_{Q^{\prime }\in \mathfrak{C}_{\limfunc{%
broken}}\left( Q\right) }\int_{Q^{\prime }}\left\vert \frac{\int_{Q^{\prime
}}fb_{Q^{\prime }}d\mu }{\int_{Q^{\prime }}b_{Q^{\prime }}d\mu }-\frac{%
\int_{Q}fb_{A}d\mu }{\int_{Q}b_{A}d\mu }\right\vert ^{2}d\mu \left( x\right)
\\
&=&\sum_{Q\in \mathcal{C}_{A}}\sum_{Q^{\prime }\in \mathfrak{C}_{\limfunc{%
natural}}\left( Q\right) }\left\vert \frac{\int_{Q^{\prime }}fb_{A}d\mu }{%
\int_{Q^{\prime }}b_{A}d\mu }-\frac{\int_{Q}fb_{A}d\mu }{\int_{Q}b_{A}d\mu }%
\right\vert ^{2}\left\vert Q^{\prime }\right\vert _{\mu } \\
&&+\sum_{Q\in \mathcal{C}_{A}}\sum_{Q^{\prime }\in \mathfrak{C}_{\limfunc{%
broken}}\left( Q\right) }\left\vert \frac{\int_{Q^{\prime }}fb_{Q^{\prime
}}d\mu }{\int_{Q^{\prime }}b_{Q^{\prime }}d\mu }-\frac{\int_{Q}fb_{A}d\mu }{%
\int_{Q}b_{A}d\mu }\right\vert ^{2}\left\vert Q^{\prime }\right\vert _{\mu }
\\
&\equiv &I_{A}+II_{A}.
\end{eqnarray*}%
To estimate term $II_{A}$ we write%
\begin{eqnarray}
&&\left\vert \frac{\int_{Q^{\prime }}fb_{Q^{\prime }}d\mu }{\int_{Q^{\prime
}}b_{Q^{\prime }}d\mu }-\frac{\int_{Q}fb_{A}d\mu }{\int_{Q}b_{A}d\mu }%
\right\vert  \label{analogue} \\
&=&\frac{\left\vert \left( \frac{1}{\left\vert Q^{\prime }\right\vert _{\mu }%
}\int_{Q^{\prime }}fb_{Q^{\prime }}d\mu \right) \left( \frac{1}{\left\vert
Q\right\vert _{\mu }}\int_{Q}b_{A}d\mu \right) -\left( \frac{1}{\left\vert
Q\right\vert _{\mu }}\int_{Q}fb_{A}d\mu \right) \left( \frac{1}{\left\vert
Q^{\prime }\right\vert _{\mu }}\int_{Q^{\prime }}b_{Q^{\prime }}d\mu \right)
\right\vert }{\left\vert \frac{1}{\left\vert Q^{\prime }\right\vert _{\mu }}%
\int_{Q^{\prime }}b_{Q^{\prime }}d\mu \right\vert \left\vert \frac{1}{%
\left\vert Q\right\vert _{\mu }}\int_{Q}b_{A}d\mu \right\vert }  \notag \\
&\lesssim &\left\vert E_{Q^{\prime }}^{\mu }\left( fb_{Q^{\prime }}\right)
E_{Q}^{\mu }\left( b_{A}\right) -E_{Q}^{\mu }\left( fb_{A}\right)
E_{Q^{\prime }}^{\mu }\left( b_{Q^{\prime }}\right) \right\vert  \notag \\
&=&\left\vert E_{Q^{\prime }}^{\mu }\left( fb_{Q^{\prime }}\right)
E_{Q}^{\mu }\left( b_{A}\right) -E_{Q^{\prime }}^{\mu }\left( fb_{A}\right)
E_{Q^{\prime }}^{\mu }\left( b_{Q^{\prime }}\right) -\left[ E_{Q}^{\mu
}\left( fb_{A}\right) -E_{Q^{\prime }}^{\mu }\left( fb_{A}\right) \right]
E_{Q^{\prime }}^{\mu }\left( b_{Q^{\prime }}\right) \right\vert  \notag \\
&\lesssim &\left\vert E_{Q^{\prime }}^{\mu }\left( fb_{Q^{\prime }}\right)
\right\vert +\left\vert E_{Q^{\prime }}^{\mu }\left( fb_{A}\right)
\right\vert +\left\vert E_{Q}^{\mu }\left( fb_{A}\right) -E_{Q^{\prime
}}^{\mu }\left( fb_{A}\right) \right\vert ,  \notag
\end{eqnarray}%
and then 
\begin{eqnarray*}
II_{A} &\lesssim &\sum_{Q\in \mathcal{C}_{A}}\sum_{Q^{\prime }\in \mathfrak{C%
}_{\limfunc{broken}}\left( Q\right) }\left( \left\vert E_{Q^{\prime }}^{\mu
}\left( fb_{Q^{\prime }}\right) \right\vert ^{2}+\left\vert E_{Q^{\prime
}}^{\mu }\left( fb_{A}\right) \right\vert ^{2}+\left\vert E_{Q}^{\mu }\left(
fb_{A}\right) \right\vert ^{2}\right) \left\vert Q^{\prime }\right\vert
_{\mu } \\
&\lesssim &\sum_{Q\in \mathcal{C}_{A}}\sum_{Q^{\prime }\in \mathfrak{C}_{%
\limfunc{broken}}\left( Q\right) }\left( \left( E_{Q^{\prime }}^{\mu
}\left\vert f\right\vert \right) ^{2}+\left( E_{Q}^{\mu }\left\vert
f\right\vert \right) ^{2}\right) \left\vert Q^{\prime }\right\vert _{\mu }\ .
\end{eqnarray*}

Now we turn to term $I_{A}$ and use the analogue of (\ref{analogue}) for
natural children, along with the natural child bounds $\left\vert
\int_{Q^{\prime }}b_{A}d\mu \right\vert \gtrsim \left\vert Q^{\prime
}\right\vert _{\mu }$\ and $\left\vert \int_{Q}b_{A}d\mu \right\vert \gtrsim
\left\vert Q\right\vert _{\mu }$, to obtain%
\begin{eqnarray*}
I_{A} &=&\sum_{Q\in \mathcal{C}_{A}}\sum_{Q^{\prime }\in \mathfrak{C}_{%
\limfunc{natural}}\left( Q\right) }\left\vert \frac{\int_{Q^{\prime
}}fb_{A}d\mu }{\int_{Q^{\prime }}b_{A}d\mu }-\frac{\int_{Q}fb_{A}d\mu }{%
\int_{Q}b_{A}d\mu }\right\vert ^{2}\left\vert Q^{\prime }\right\vert _{\mu }
\\
&=&\sum_{Q\in \mathcal{C}_{A}}\sum_{Q^{\prime }\in \mathfrak{C}_{\limfunc{%
natural}}\left( Q\right) }\left\vert \frac{\left( \int_{Q^{\prime
}}fb_{A}d\mu \right) \left( \int_{Q}b_{A}d\mu \right) -\left(
\int_{Q^{\prime }}b_{A}d\mu \right) \left( \int_{Q}fb_{A}d\mu \right) }{%
\left( \int_{Q^{\prime }}b_{A}d\mu \right) \left( \int_{Q}b_{A}d\mu \right) }%
\right\vert ^{2}\left\vert Q^{\prime }\right\vert _{\mu } \\
&\lesssim &\sum_{Q\in \mathcal{C}_{A}}\sum_{Q^{\prime }\in \mathfrak{C}_{%
\limfunc{natural}}\left( Q\right) }\left\vert E_{Q}^{\mu }\left(
b_{A}\right) \right\vert ^{2}\left\vert E_{Q^{\prime }}^{\mu }\left(
fb_{A}\right) -E_{Q}^{\mu }\left( fb_{A}\right) \right\vert ^{2}\left\vert
Q^{\prime }\right\vert _{\mu } \\
&&+\sum_{Q\in \mathcal{C}_{A}}\sum_{Q^{\prime }\in \mathfrak{C}_{\limfunc{%
natural}}\left( Q\right) }\left\vert E_{Q}^{\mu }\left( b_{A}\right)
-E_{Q^{\prime }}^{\mu }b_{A}\right\vert ^{2}\left\vert E_{Q}^{\mu }\left(
fb_{A}\right) \right\vert ^{2}\left\vert Q^{\prime }\right\vert _{\mu } \\
&\equiv &III_{A}+IV_{A}.
\end{eqnarray*}

Now we have%
\begin{equation*}
III_{A}\lesssim \sum_{Q\in \mathcal{C}_{A}}\sum_{Q^{\prime }\in \mathfrak{C}%
_{\limfunc{natural}}\left( Q\right) }\left\vert E_{Q^{\prime }}^{\mu }\left(
fb_{A}\right) -E_{Q}^{\mu }\left( fb_{A}\right) \right\vert ^{2}\left\vert
Q^{\prime }\right\vert _{\mu }\ ,
\end{equation*}%
and for term $IV_{A}$, we introduce the quantities 
\begin{equation*}
\gamma _{Q}\equiv \sum_{Q^{\prime }\in \mathfrak{C}_{\limfunc{natural}%
}\left( Q\right) }\left\vert E_{Q}^{\mu }b_{A}-E_{Q^{\prime }}^{\mu
}b_{A}\right\vert ^{2}\left\vert Q^{\prime }\right\vert _{\mu },\ \ \ \ \ 
\text{for }Q\in \mathcal{C}_{A},A\in \mathcal{A\ },
\end{equation*}%
so that%
\begin{eqnarray*}
IV_{A} &=&\sum_{Q\in \mathcal{C}_{A}}\sum_{Q^{\prime }\in \mathfrak{C}_{%
\limfunc{natural}}\left( Q\right) }\left\vert E_{Q}^{\mu }\left(
b_{A}\right) -E_{Q^{\prime }}^{\mu }b_{A}\right\vert ^{2}\left\vert
Q^{\prime }\right\vert _{\mu }\left\vert E_{Q}^{\mu }\left( fb_{A}\right)
\right\vert ^{2} \\
&=&\sum_{Q\in \mathcal{C}_{A}}\gamma _{Q}\left\vert E_{Q}^{\mu }\left(
fb_{A}\right) \right\vert ^{2}\lesssim \sum_{Q\in \mathcal{C}_{A}}\gamma
_{Q}\left( E_{Q}^{\mu }\left\vert f\right\vert \right) ^{2}.
\end{eqnarray*}

Now note that the coefficients $\left\{ \gamma _{Q}\right\} _{Q\in \mathcal{D%
}}$ satisfy the Carleson condition%
\begin{equation}
\dsum\limits_{Q\subset B}\gamma _{Q}\lesssim \left\vert B\right\vert
_{\sigma }.  \label{Car cond}
\end{equation}%
Indeed, if $B\in \mathcal{C}_{A}$, then using $E_{Q}^{\mu }b_{A}=E_{Q}^{\mu
}\left( \mathbf{1}_{B}b_{A}\right) $ for $Q\subset B$ and the unweighted
square function estimate, and denoting by $\mathfrak{G}_{\mathcal{A}%
}^{t}\left( A\right) $ the $\mathcal{A}$-grandchildren at level $t$ below $A$%
, we have 
\begin{eqnarray*}
\dsum\limits_{Q\subset B}\gamma _{Q} &=&\sum_{A\in \mathcal{A}%
}\dsum\limits_{Q\in \mathcal{C}_{A}:\ Q\subset B}\sum_{Q^{\prime }\in 
\mathfrak{C}_{\limfunc{natural}}\left( Q\right) }\left\vert E_{Q}^{\mu
}b_{A}-E_{Q^{\prime }}^{\mu }b_{A}\right\vert ^{2}\left\vert Q^{\prime
}\right\vert _{\mu } \\
&&+\sum_{t=1}^{\infty }\sum_{H\in \mathfrak{G}_{\mathcal{A}}^{t}\left(
A\right) :\ H\subset B}\dsum\limits_{Q\in \mathcal{C}_{H}}\sum_{Q^{\prime
}\in \mathfrak{C}_{\limfunc{natural}}\left( Q\right) }\left\vert E_{Q}^{\mu
}b_{H}-E_{Q^{\prime }}^{\mu }b_{H}\right\vert ^{2}\left\vert Q^{\prime
}\right\vert _{\mu } \\
&\lesssim &\int_{B}\left\vert b_{A}\right\vert ^{2}d\mu +\sum_{t=1}^{\infty
}\sum_{H\in \mathfrak{G}_{\mathcal{A}}^{t}\left( A\right) :\ H\subset B}\int
\left\vert b_{H}\right\vert ^{2}d\mu \\
&\lesssim &\left\vert B\right\vert _{\mu }+\sum_{t=1}^{\infty }\sum_{H\in 
\mathfrak{G}_{\mathcal{A}}^{t}\left( A\right) :\ H\subset B}\left\vert
H\right\vert _{\mu }\lesssim \left\vert B\right\vert _{\mu }\ .
\end{eqnarray*}

Altogether then, using the above inequalities for $II_{A},III_{A}$ and $%
IV_{A}$, and then using the Carleson embedding theorem on the sum of the
terms $II_{A}$, and the Carleson embedding theorem again on the sum of the
terms $IV_{A}$ with the Carleson condition (\ref{Car cond}), we conclude that%
\begin{eqnarray*}
\sum_{A\in \mathcal{A}}\sum_{Q\in \mathcal{C}_{A}}\left\Vert \bigtriangleup
_{Q}^{\mu ,\mathbf{b}}f\right\Vert _{L^{2}\left( \mu \right) }^{2} &\lesssim
&\sum_{A\in \mathcal{A}}\left( II_{A}+III_{A}+IV_{A}\right) \\
&\lesssim &\sum_{A\in \mathcal{A}}\left\vert A\right\vert _{\mu }\left\{
\left( \frac{1}{\left\vert A\right\vert _{\mu }}\int_{A}\left\vert
f\right\vert d\mu \right) ^{2}+\left( \frac{1}{\left\vert \pi A\right\vert
_{\mu }}\int_{\pi A}\left\vert f\right\vert d\mu \right) ^{2}\right\} \\
&&+\sum_{A\in \mathcal{A}}III_{A}+\sum_{Q\in \mathcal{C}_{A}}\gamma
_{Q}\left( E_{Q}^{\mu }\left\vert f\right\vert \right) ^{2}\lesssim \int
\left\vert f\right\vert ^{2}d\mu +\sum_{A\in \mathcal{A}}III_{A}.
\end{eqnarray*}%
Thus it remains only to estimate the final square function expression, which
requires an additional argument due to the fact that the functions $fb_{A}$
differ from one corona to the next. For this we define%
\begin{equation}
b\equiv \sum_{A\in \mathcal{A}}b_{A}\mathbf{1}_{A\setminus
\bigcup_{A^{\prime }\in \mathfrak{C}_{\mathcal{A}}\left( A\right) }}
\label{def b}
\end{equation}%
and write%
\begin{eqnarray*}
&&\sum_{A\in \mathcal{A}}III_{A}\lesssim \sum_{A\in \mathcal{A}}\sum_{Q\in 
\mathcal{C}_{A}}\sum_{Q^{\prime }\in \mathfrak{C}_{\limfunc{natural}}\left(
Q\right) }\left\vert E_{Q^{\prime }}^{\mu }\left( fb_{A}\right) -E_{Q}^{\mu
}\left( fb_{A}\right) \right\vert ^{2}\left\vert Q^{\prime }\right\vert
_{\mu } \\
&\lesssim &\sum_{A\in \mathcal{A}}\sum_{Q\in \mathcal{C}_{A}}\sum_{Q^{\prime
}\in \mathfrak{C}_{\limfunc{natural}}\left( Q\right) }\left\vert
E_{Q^{\prime }}^{\mu }\left( fb_{A}-fb\right) -E_{Q}^{\mu }\left(
fb_{A}-fb\right) \right\vert ^{2}\left\vert Q^{\prime }\right\vert _{\mu } \\
&&+\sum_{A\in \mathcal{A}}\sum_{Q\in \mathcal{C}_{A}}\sum_{Q^{\prime }\in 
\mathfrak{C}_{\limfunc{natural}}\left( Q\right) }\left\vert E_{Q^{\prime
}}^{\mu }\left( fb\right) -E_{Q}^{\mu }\left( fb\right) \right\vert
^{2}\left\vert Q^{\prime }\right\vert _{\mu }\ .
\end{eqnarray*}%
Now the unweighted square function estimate applies to the second sum to give%
\begin{equation*}
\sum_{A\in \mathcal{A}}\sum_{Q\in \mathcal{C}_{A}}\sum_{Q^{\prime }\in 
\mathfrak{C}_{\limfunc{natural}}\left( Q\right) }\left\vert E_{Q^{\prime
}}^{\mu }\left( fb\right) -E_{Q}^{\mu }\left( fb\right) \right\vert
^{2}\left\vert Q^{\prime }\right\vert _{\mu }\lesssim \int \left\vert
fb\right\vert ^{2}d\mu \lesssim \int \left\vert f\right\vert ^{2}d\mu .
\end{equation*}%
To handle the first sum we write both%
\begin{eqnarray}
E_{Q^{\prime }}^{\mu }\left( fb_{A}-fb\right) &=&E_{Q^{\prime }}^{\mu
}\left( \sum_{A^{\prime }\in \mathfrak{C}_{\mathcal{A}}\left( A\right) }%
\mathbf{1}_{A^{\prime }}\left( fb_{A}-fb\right) \right) =E_{Q^{\prime
}}^{\mu }\left( \sum_{A^{\prime }\in \mathfrak{C}_{\mathcal{A}}\left(
A\right) }\mathbb{E}_{A^{\prime }}^{\mu }\left( fb_{A}-fb\right) \right) ,
\label{write both} \\
E_{Q}^{\mu }\left( fb_{A}-fb\right) &=&E_{Q}^{\mu }\left( \sum_{A^{\prime
}\in \mathfrak{C}_{\mathcal{A}}\left( A\right) }\mathbf{1}_{A^{\prime
}}\left( fb_{A}-fb\right) \right) =E_{Q}^{\mu }\left( \sum_{A^{\prime }\in 
\mathfrak{C}_{\mathcal{A}}\left( A\right) }\mathbb{E}_{A^{\prime }}^{\mu
}\left( fb_{A}-fb\right) \right) ,  \notag
\end{eqnarray}%
and use the unweighted square function estimate on each corona $\mathcal{C}%
_{A}$, applied to the function $fb_{A}-fb$, to obtain%
\begin{eqnarray*}
&&\sum_{Q\in \mathcal{C}_{A}}\sum_{Q^{\prime }\in \mathfrak{C}_{\limfunc{%
natural}}\left( Q\right) }\left\vert E_{Q^{\prime }}^{\mu }\left(
fb_{A}-fb\right) -E_{Q}^{\mu }\left( fb_{A}-fb\right) \right\vert
^{2}\left\vert Q^{\prime }\right\vert _{\mu } \\
&\lesssim &\int_{A}\left\vert \sum_{A^{\prime }\in \mathfrak{C}_{\mathcal{A}%
}\left( A\right) }\mathbb{E}_{A^{\prime }}^{\mu }\left( fb_{A}-fb\right)
\right\vert ^{2}d\mu =\sum_{A^{\prime }\in \mathfrak{C}_{\mathcal{A}}\left(
A\right) }\left\vert E_{A^{\prime }}^{\mu }\left( fb_{A}-fb\right)
\right\vert ^{2}\left\vert A^{\prime }\right\vert _{\mu }\lesssim
\sum_{A^{\prime }\in \mathfrak{C}_{\mathcal{A}}\left( A\right) }\left(
E_{A^{\prime }}^{\mu }\left\vert f\right\vert \right) ^{2}\left\vert
A^{\prime }\right\vert _{\mu }\ ,
\end{eqnarray*}%
We can now sum over $A\in \mathcal{A}$ to obtain%
\begin{eqnarray*}
&&\sum_{A\in \mathcal{A}}\sum_{Q\in \mathcal{C}_{A}}\sum_{Q^{\prime }\in 
\mathfrak{C}_{\limfunc{natural}}\left( Q\right) }\left\vert E_{Q^{\prime
}}^{\mu }\left( fb_{A}-fb\right) -E_{Q}^{\mu }\left( fb_{A}-fb\right)
\right\vert ^{2}\left\vert Q^{\prime }\right\vert _{\mu } \\
&\lesssim &\sum_{A\in \mathcal{A}}\sum_{A^{\prime }\in \mathfrak{C}_{%
\mathcal{A}}\left( A\right) }\left( E_{A^{\prime }}^{\mu }\left\vert
f\right\vert \right) ^{2}\left\vert A^{\prime }\right\vert _{\mu }\lesssim
\int \left\vert f\right\vert ^{2}d\mu
\end{eqnarray*}%
by the Carleson embedding theorem yet again.
\end{proof}

Essentially the same proof gives the lower frame inequality for the dual
martingale difference $\square _{Q}^{\mu ,\mathbf{b}}f=\left( \bigtriangleup
_{Q}^{\mu ,\mathbf{b}}\right) ^{\ast }f$: 
\begin{equation}
\sum_{Q\in \mathcal{D}}\left\Vert \square _{Q}^{\mu ,\mathbf{b}}f\right\Vert
_{L^{2}\left( \mu \right) }^{2}\lesssim \left\Vert f\right\Vert
_{L^{2}\left( \mu \right) }^{2}\ ,  \label{low frame}
\end{equation}%
but matters are made slightly simpler by the fact that we do not need the
function $b$ introduced in (\ref{def b}) above because in the definition of $%
\square _{Q}^{\mu ,\mathbf{b}}f$ the function $b_{Q}$ doesn't multiply $f$,
rather it sits outside the integral where it can be estimated crudely,
leaving the unweighted square function applied to $f$ alone. We leave this
for the reader, pointing out that the general idea of the proof can be found
in the proof given below in Proposition \ref{reverse half Riesz dual} for
the lower weak Riesz inequality for $\square _{Q}^{\mu ,\mathbf{b}}$.

The corresponding upper weak frame inequalities 
\begin{eqnarray}
\left\Vert f\right\Vert _{L^{2}\left( \mu \right) }^{2} &\lesssim
&\sum_{Q\in \mathcal{D}}\left\Vert \bigtriangleup _{Q}^{\mu ,\mathbf{b}%
}f\right\Vert _{L^{2}\left( \mu \right) }^{2}+\sum_{Q\in \mathcal{D}%
}\left\Vert \nabla _{Q}^{\mu }f\right\Vert _{L^{2}\left( \mu \right) }^{2}\ ,
\label{corr upper} \\
\left\Vert f\right\Vert _{L^{2}\left( \mu \right) }^{2} &\lesssim
&\sum_{Q\in \mathcal{D}}\left\Vert \square _{Q}^{\mu ,\mathbf{b}%
}f\right\Vert _{L^{2}\left( \mu \right) }^{2}+\sum_{Q\in \mathcal{D}%
}\left\Vert \nabla _{Q}^{\mu }f\right\Vert _{L^{2}\left( \mu \right) }^{2}\ ,
\notag
\end{eqnarray}%
are proved by duality using the `Calderon reproducing formulas' 
\begin{equation*}
\bigtriangleup _{Q}^{\mu ,\mathbf{b}}=\left( \bigtriangleup _{Q}^{\mu ,\pi ,%
\mathbf{b}}\right) ^{2}+\bigtriangleup _{Q,\limfunc{broken}}^{\mu ,\mathbf{b}%
}\text{ and }\square _{Q}^{\mu ,\mathbf{b}}=\left( \square _{Q}^{\mu ,\pi ,%
\mathbf{b}}\right) ^{2}+\square _{Q,\limfunc{broken}}^{\mu ,\mathbf{b}}
\end{equation*}%
introduced in the following subsection. We will actually prove stronger
inequalities below, namely upper weak \emph{Riesz} inequalities.

\subsection{Unbroken projections, broken differences and a Calder\'{o}n
reproducing formula}

Here we assume that $\int_{Q^{\prime }}b_{Q}d\mu \neq 0$, noting that in
applications we will have the stronger inequality $\left\vert
\int_{Q^{\prime }}b_{Q}d\omega \right\vert \gtrsim \mathbf{1}_{Q^{\prime
}}\left\Vert b_{Q}\right\Vert _{L^{\infty }\left( \omega \right) }>0$ due to
the assumed reverse H\"{o}lder control on children (\ref{rev Hol con}).
Define 
\begin{equation}
\mathbb{E}_{Q}^{\mu ,\pi ,\mathbf{b}}f\equiv \mathbf{1}_{Q}\frac{1}{%
\int_{Q}b_{\pi Q}d\mu }\int_{Q}b_{\pi Q}fd\mu \text{ and }\mathbb{F}%
_{Q}^{\mu ,\pi ,\mathbf{b}}f\equiv \mathbf{1}_{Q}\frac{b_{\pi Q}}{%
\int_{Q}b_{\pi Q}d\mu }\int_{Q}fd\mu  \label{def pi exp}
\end{equation}%
and%
\begin{eqnarray}
\bigtriangleup _{Q}^{\mu ,\pi ,\mathbf{b}}f &=&\left[ \sum_{Q^{\prime }\in 
\mathfrak{C}\left( Q\right) }\mathbb{E}_{Q^{\prime }}^{\mu ,\pi ,\mathbf{b}}f%
\right] -\mathbb{E}_{Q}^{\mu ,\mathbf{b}}f=\sum_{Q^{\prime }\in \mathfrak{C}%
\left( Q\right) }\mathbb{E}_{Q^{\prime }}^{\mu ,b_{Q}}f-\mathbb{E}_{Q}^{\mu
,b_{Q}}f,  \label{def pi box} \\
\square _{Q}^{\mu ,\pi ,\mathbf{b}}f &=&\left[ \sum_{Q^{\prime }\in 
\mathfrak{C}\left( Q\right) }\mathbb{F}_{Q^{\prime }}^{\mu ,\pi ,\mathbf{b}}f%
\right] -\mathbb{F}_{Q}^{\mu ,\mathbf{b}}f=\sum_{Q^{\prime }\in \mathfrak{C}%
\left( Q\right) }\mathbb{F}_{Q^{\prime }}^{\mu ,b_{Q}}f-\mathbb{F}_{Q}^{\mu
,b_{Q}}f,  \notag
\end{eqnarray}%
where on the far right of (\ref{def pi box}) we are using the notation $%
\mathbb{E}_{Q^{\prime }}^{\mu ,b}=\mathbf{1}_{Q^{\prime }}\frac{1}{%
\int_{Q^{\prime }}bd\mu }\int_{Q^{\prime }}bfd\mu $ when $b$ is simply a
function, rather than a family of functions $\mathbf{b}$, in order to
specify the testing function $b$ we use if it differs from the function $%
b_{Q^{\prime }}$ that is selected in the notation $\mathbb{E}_{Q^{\prime
}}^{\mu ,\mathbf{b}}$ when the family $\mathbf{b}$ appears in boldface as an
exponent. Similarly for other pseudoprojections.

For convenience of notation, we set $\bigtriangleup _{Q}^{\mu ,\pi ,\mathbf{b%
}}=\bigtriangleup _{Q}^{\mu ,b}$ and $\square _{Q}^{\mu ,\pi ,\mathbf{b}%
}=\square _{Q}^{\mu ,b}$ with $b=b_{Q}$. Note that%
\begin{eqnarray*}
\mathbb{E}_{Q}^{\mu ,b}\left( 1\right) &=&\mathbf{1}_{Q}\frac{1}{%
\int_{Q}bd\mu }\int_{Q}\left( 1\right) bd\mu =\mathbf{1}_{Q}\ , \\
\mathbb{F}_{Q}^{\mu ,b}\left( b\right) &=&\mathbf{1}_{Q}\frac{b}{%
\int_{Q}bd\mu }\int_{Q}\left( b\right) d\mu =\mathbf{1}_{Q}\ b, \\
\triangle _{Q}^{\mu ,b}\left( 1\right) &=&\left[ \sum_{Q^{\prime \prime }\in 
\mathfrak{C}\left( Q\right) }\mathbb{E}_{Q^{\prime \prime }}^{\mu ,b}1\right]
-\mathbb{E}_{Q}^{\mu ,b}1=\left[ \sum_{Q^{\prime \prime }\in \mathfrak{C}%
\left( Q\right) }\mathbf{1}_{Q^{\prime \prime }}\right] -\mathbf{1}_{Q}=0, \\
\square _{Q}^{\mu ,b}b &=&\left[ \sum_{Q^{\prime \prime }\in \mathfrak{C}%
\left( Q\right) }\mathbb{F}_{Q^{\prime \prime }}^{\mu ,b}b\right] -\mathbb{F}%
_{Q}^{\mu ,b}b=\left[ \sum_{Q^{\prime \prime }\in \mathfrak{C}\left(
Q\right) }\mathbf{1}_{Q^{\prime \prime }}b\right] -\mathbf{1}_{Q}b=0.
\end{eqnarray*}%
The next two lemmas are in \cite[see page 193]{NTV2}.

\begin{lemma}
\label{b orth}For dyadic cubes $R$ and $Q$ we have 
\begin{equation*}
\mathbb{E}_{R}^{\mu ,b}\bigtriangleup _{Q}^{\mu ,b}=\left\{ 
\begin{array}{ccc}
0 & \text{ if } & R\supset Q\text{ or }R\cap Q=\emptyset \\ 
\mathbf{1}_{R}E_{Q_{R}}^{\mu ,b}\bigtriangleup _{Q}^{\mu ,b}f & \text{ if }
& R\subsetneqq Q%
\end{array}%
\right. .
\end{equation*}
\end{lemma}

\begin{proof}
If $R\supset Q$, then since $\square _{Q}^{\mu ,b}b=0$, we have 
\begin{equation*}
\mathbb{E}_{R}^{\mu ,b}\bigtriangleup _{Q}^{\mu ,b}f=\mathbf{1}_{R}\frac{1}{%
\int_{R}bd\mu }\int_{R}\left( \bigtriangleup _{Q}^{\mu ,b}f\right) bd\mu =%
\mathbf{1}_{R}\frac{\left\langle \bigtriangleup _{Q}^{\mu
,b}f,b\right\rangle _{L^{2}\left( \mu \right) }}{\int_{R}bd\mu }=\mathbf{1}%
_{R}\frac{\left\langle f,\square _{Q}^{\mu ,b}b\right\rangle _{L^{2}\left(
\mu \right) }}{\int_{R}bd\mu }=0.
\end{equation*}%
On the other hand, if $R\subsetneqq Q$, then $R\subset Q^{\prime }$ for some 
$Q^{\prime }\in \mathfrak{C}\left( Q\right) $, and since $\bigtriangleup
_{Q}^{\mu ,b}f$ equals the constant $A=E_{Q_{R}}^{\mu ,b}\left(
\bigtriangleup _{Q}^{\mu ,b}f\right) $ on $Q^{\prime }$, and $E_{I}^{\mu
,b}1=1$ for all cubes $I$, we have%
\begin{equation*}
\mathbb{E}_{R}^{\mu ,b}\bigtriangleup _{Q}^{\mu ,b}f=\mathbb{E}_{R}^{\mu
,b}A=\mathbf{1}_{R}A=\mathbf{1}_{R}E_{Q_{R}}^{\mu ,b}\left( \bigtriangleup
_{Q}^{\mu ,b}f\right) .
\end{equation*}
\end{proof}

\begin{lemma}
\label{b proj}For dyadic cubes $R$ and $Q$ we have%
\begin{equation*}
\bigtriangleup _{R}^{\mu ,b}\bigtriangleup _{Q}^{\mu ,b}=\left\{ 
\begin{array}{ccc}
\bigtriangleup _{Q}^{\mu ,b} & \text{ if } & R=Q \\ 
0 & \text{ if } & R\not=Q%
\end{array}%
\right. .
\end{equation*}
\end{lemma}

\begin{proof}
By the top line in Lemma \ref{b orth}, it suffices to consider the case $%
R\cap Q\neq \emptyset $. First we suppose that $R\subsetneqq Q$. Then $%
R\subset Q^{\prime }$ for some $Q^{\prime }\in \mathfrak{C}\left( Q\right) $%
, and since $\bigtriangleup _{Q}^{\mu ,b}f=A$ is constant on $Q^{\prime }$,
and $\mathbb{E}_{I}^{\mu ,b}1=\mathbf{1}_{I}$ for any cube $I$, we obtain 
\begin{equation*}
\bigtriangleup _{R}^{\mu ,b}\bigtriangleup _{Q}^{\mu ,b}f=\left[
\sum_{R^{\prime }\in \mathfrak{C}\left( R\right) }\mathbb{E}_{R^{\prime
}}^{\mu ,b}\left\{ \bigtriangleup _{Q}^{\mu ,b}f\right\} \right] -\mathbb{E}%
_{R}^{\mu ,b}\left\{ \bigtriangleup _{Q}^{\mu ,b}f\right\} =\left[
\sum_{R^{\prime }\in \mathfrak{C}\left( R\right) }\mathbf{1}_{R^{\prime }}A%
\right] -\mathbf{1}_{R}A=0.
\end{equation*}

Next we suppose that $R=Q$ and obtain%
\begin{eqnarray*}
&&\bigtriangleup _{Q}^{\mu ,b}\bigtriangleup _{Q}^{\mu ,b}f \\
&=&\left[ \sum_{Q^{\prime }\in \mathfrak{C}\left( Q\right) }\mathbb{E}%
_{Q^{\prime }}^{\mu ,b}\left\{ \bigtriangleup _{Q}^{\mu ,b}f\right\} \right]
-\mathbb{E}_{Q}^{\mu ,b}\left\{ \bigtriangleup _{Q}^{\mu ,b}f\right\} \\
&=&\left[ \sum_{Q^{\prime }\in \mathfrak{C}\left( Q\right) }\mathbb{E}%
_{Q^{\prime }}^{\mu ,b}\left\{ \left[ \sum_{Q^{\prime \prime }\in \mathfrak{C%
}\left( Q\right) }\mathbb{E}_{Q^{\prime \prime }}^{\mu ,b}f\right] -\mathbb{E%
}_{Q}^{\mu ,b}f\right\} \right] -\mathbb{E}_{Q}^{\mu ,b}\left[ \left\{ \left[
\sum_{Q^{\prime \prime }\in \mathfrak{C}\left( Q\right) }\mathbb{E}%
_{Q^{\prime \prime }}^{\mu ,b}f\right] -\mathbb{E}_{Q}^{\mu ,b}f\right\} %
\right] \\
&=&\left\{ \left[ \sum_{Q^{\prime }\in \mathfrak{C}\left( Q\right) }\mathbb{E%
}_{Q^{\prime }}^{\mu ,b}f\right] -\mathbb{E}_{Q}^{\mu ,b}f\right\} -\mathbb{E%
}_{Q}^{\mu ,b}\left[ \left\{ \left[ \sum_{Q^{\prime \prime }\in \mathfrak{C}%
\left( Q\right) }\mathbb{E}_{Q^{\prime \prime }}^{\mu ,b}f\right] -\mathbb{E}%
_{Q}^{\mu ,b}f\right\} \right] \\
&=&\bigtriangleup _{Q}^{\mu ,b}f-\mathbb{E}_{Q}^{\mu ,b}\bigtriangleup
_{Q}^{\mu ,b}f=\bigtriangleup _{Q}^{\mu ,b}f,
\end{eqnarray*}%
where we have used Lemma \ref{b orth} with $R=Q$ for the final equality.

Finally we suppose that $R\supsetneqq Q$. Then $R_{Q}\supset Q$, and so by
the top line in Lemma \ref{b orth} we have 
\begin{equation*}
\bigtriangleup _{R}^{\mu ,b}\bigtriangleup _{Q}^{\mu ,b}f=\sum_{R^{\prime
}\in \mathfrak{C}\left( R\right) }\mathbb{E}_{R^{\prime }}^{\mu
,b}\bigtriangleup _{Q}^{\mu ,b}f-\mathbb{E}_{R}^{\mu ,b}\bigtriangleup
_{Q}^{\mu ,b}f=0-0=0.
\end{equation*}
\end{proof}

Now since we are assuming that $\int_{Q^{\prime }}b_{Q}\neq 0$, we can define%
\begin{eqnarray*}
\bigtriangleup _{Q,\limfunc{broken}}^{\mu ,\pi ,\mathbf{b}}f
&=&\bigtriangleup _{Q}^{\mu ,\mathbf{b}}f-\bigtriangleup _{Q}^{\mu ,\pi ,%
\mathbf{b}}f \\
&=&\left( \sum_{Q^{\prime }\in \mathfrak{C}\left( Q\right) }\mathbb{E}%
_{Q^{\prime }}^{\mu ,b_{Q^{\prime }}}f-\mathbb{E}_{Q}^{\mu ,b_{Q}}f\right)
-\left( \sum_{Q^{\prime }\in \mathfrak{C}\left( Q\right) }\mathbb{E}%
_{Q^{\prime }}^{\mu ,b_{Q}}f-\mathbb{E}_{Q}^{\mu ,b_{Q}}f\right) \\
&=&\sum_{Q^{\prime }\in \mathfrak{C}_{\limfunc{broken}}\left( Q\right) }%
\mathbb{E}_{Q^{\prime }}^{\mu ,b_{Q^{\prime }}}f-\mathbb{E}_{Q^{\prime
}}^{\mu ,b_{Q}}f\ ,
\end{eqnarray*}%
with a similar definition for $\square _{Q,\limfunc{broken}}^{\mu ,\mathbf{b}%
}f$. Altogether, with $\square _{Q}^{\mu ,\mathbf{b}}=\left( \bigtriangleup
_{Q}^{\mu ,\mathbf{b}}\right) ^{\ast }$, $\square _{Q}^{\mu ,\pi ,\mathbf{b}%
}=\left( \bigtriangleup _{Q}^{\mu ,\pi ,\mathbf{b}}\right) ^{\ast }$ and $%
\square _{Q,\limfunc{broken}}^{\mu ,\mathbf{b}}=\left( \bigtriangleup _{Q,%
\limfunc{broken}}^{\mu ,\mathbf{b}}\right) ^{\ast }$, we have%
\begin{equation}
\bigtriangleup _{Q}^{\mu ,\mathbf{b}}=\bigtriangleup _{Q}^{\mu ,\pi ,\mathbf{%
b}}+\bigtriangleup _{Q,\limfunc{broken}}^{\mu ,\pi ,\mathbf{b}}\text{ and }%
\square _{Q}^{\mu ,\mathbf{b}}=\square _{Q}^{\mu ,\pi ,\mathbf{b}}+\square
_{Q,\limfunc{broken}}^{\mu ,\pi ,\mathbf{b}}  \label{box pi equals}
\end{equation}%
where $\bigtriangleup _{Q}^{\mu ,\pi ,\mathbf{b}}$ and $\square _{Q}^{\mu
,\pi ,\mathbf{b}}$ are projections and 
\begin{eqnarray*}
\bigtriangleup _{Q,\limfunc{broken}}^{\mu ,\pi ,\mathbf{b}}f
&=&\sum_{Q^{\prime }\in \mathfrak{C}_{\limfunc{broken}}\left( Q\right)
}\left( \mathbb{E}_{Q^{\prime }}^{\mu ,b_{Q^{\prime }}}f-\mathbb{E}%
_{Q^{\prime }}^{\mu ,b_{Q}}f\right) , \\
\left\vert \bigtriangleup _{Q,\limfunc{broken}}^{\mu ,\pi ,\mathbf{b}%
}f\right\vert &\leq &\sum_{Q^{\prime }\in \mathfrak{C}_{\limfunc{broken}%
}\left( Q\right) }\left( \frac{\left\Vert \mathbf{1}_{Q^{\prime
}}b_{Q}\right\Vert _{L^{\infty }\left( \mu \right) }}{\left\vert
\int_{Q^{\prime }}b_{Q}d\mu \right\vert }+\left\Vert b_{Q^{\prime
}}\right\Vert _{L^{\infty }\left( \mu \right) }\right) \frac{1}{\left\vert
Q^{\prime }\right\vert _{\mu }}\int_{Q^{\prime }}\left\vert f\right\vert
d\mu \leq C\sum_{Q^{\prime }\in \mathfrak{C}_{\limfunc{broken}}\left(
Q\right) }\frac{1}{\left\vert Q^{\prime }\right\vert _{\mu }}\int_{Q^{\prime
}}\left\vert f\right\vert d\mu , \\
\square _{Q,\limfunc{broken}}^{\mu ,\pi ,\mathbf{b}}f &=&\sum_{Q^{\prime
}\in \mathfrak{C}_{\limfunc{broken}}\left( Q\right) }\left( \mathbb{F}%
_{Q^{\prime }}^{\mu ,b_{Q^{\prime }}}f-\mathbb{F}_{Q^{\prime }}^{\mu
,b_{Q}}f\right) , \\
\left\vert \square _{Q,\limfunc{broken}}^{\mu ,\pi ,\mathbf{b}}f\right\vert
&\leq &\sum_{Q^{\prime }\in \mathfrak{C}_{\limfunc{broken}}\left( Q\right)
}\left( \frac{\left\Vert \mathbf{1}_{Q^{\prime }}b_{Q}\right\Vert
_{L^{\infty }\left( \mu \right) }}{\left\vert \int_{Q^{\prime }}b_{Q}d\mu
\right\vert }+\left\Vert b_{Q^{\prime }}\right\Vert _{L^{\infty }\left( \mu
\right) }\right) \frac{1}{\left\vert Q^{\prime }\right\vert _{\mu }}%
\int_{Q^{\prime }}\left\vert f\right\vert d\mu \leq C\sum_{Q^{\prime }\in 
\mathfrak{C}_{\limfunc{broken}}\left( Q\right) }\frac{1}{\left\vert
Q^{\prime }\right\vert _{\mu }}\int_{Q^{\prime }}\left\vert f\right\vert
d\mu ,
\end{eqnarray*}%
where $C$ depends on both $C_{\mathbf{b}}$ and the constant in the reverse H%
\"{o}lder condition on children in (\ref{rev Hol con}).

Altogether then we have when $\int_{Q_{i}}b_{Q}d\sigma \neq 0$ for both
children $Q_{i}$ of $Q$, the `Calder\'{o}n reproducing formula', 
\begin{equation}
\bigtriangleup _{Q}^{\mu ,\mathbf{b}}f=\bigtriangleup _{Q}^{\mu ,\pi ,%
\mathbf{b}}\bigtriangleup _{Q}^{\mu ,\pi ,\mathbf{b}}f+\bigtriangleup _{Q,%
\limfunc{broken}}^{\mu ,\pi ,\mathbf{b}}f\text{ and }\square _{Q}^{\mu ,%
\mathbf{b}}f=\square _{Q}^{\mu ,\pi ,\mathbf{b}}\square _{Q}^{\mu ,\pi ,%
\mathbf{b}}f+\square _{Q,\limfunc{broken}}^{\mu ,\pi ,\mathbf{b}}f,
\label{square of delta}
\end{equation}%
and the pointwise estimates%
\begin{equation*}
\left\vert \bigtriangleup _{Q}^{\mu ,\mathbf{b}}f\right\vert ,\left\vert
\square _{Q}^{\mu ,\mathbf{b}}f\right\vert \leq C_{\delta ,\mathbf{b}%
}\sum_{Q^{\prime }\in \mathfrak{C}\left( Q\right) }\mathbf{1}_{Q^{\prime
}}\left( \frac{1}{\left\vert Q^{\prime }\right\vert _{\mu }}\int_{Q^{\prime
}}\left\vert f\right\vert d\mu +\frac{1}{\left\vert Q\right\vert _{\mu }}%
\int_{Q}\left\vert f\right\vert d\mu \right) \ ,
\end{equation*}%
and%
\begin{equation}
\left\vert \bigtriangleup _{Q,\limfunc{broken}}^{\mu ,\pi ,\mathbf{b}%
}f\right\vert \lesssim \left\vert \widehat{\bigtriangledown }_{Q}^{\mu
}f\right\vert \text{ and }\left\vert \square _{Q,\limfunc{broken}}^{\mu ,\pi
,\mathbf{b}}f\right\vert \lesssim \left\vert \widehat{\bigtriangledown }%
_{Q}^{\mu }f\right\vert ,  \label{F est}
\end{equation}%
which follow from the reverse H\"{o}lder property in Lemma \ref{prelim
control of corona} of the children of $Q$,%
\begin{equation}
\left\Vert \mathbf{1}_{Q_{i}}b_{Q}\right\Vert _{L^{\infty }\left( \sigma
\right) }<\frac{16C_{b_{Q}}}{\delta }\left\vert \frac{1}{\left\vert
Q_{i}\right\vert _{\sigma }}\int_{Q_{i}}b_{Q}d\sigma \right\vert \ ,\ \ \ \
\ Q_{i}\in \mathfrak{C}\left( Q\right) .  \label{rev Hold}
\end{equation}%
Note again that the formulas in (\ref{square of delta}) always hold because
our reverse H\"{o}lder assumption (\ref{rev Hol con}) in the triple corona
construction implies in particular that $\left\Vert \mathbf{1}%
_{Q_{i}}b_{Q}\right\Vert _{L^{\infty }\left( \sigma \right) }>0$.

\subsubsection{Another modified dual martingale difference}

Define another modified dual martingale difference by 
\begin{equation}
\square _{I}^{\sigma ,\flat ,\mathbf{b}}f\equiv \square _{I}^{\sigma ,%
\mathbf{b}}f-\sum_{I^{\prime }\in \mathfrak{C}_{\limfunc{broken}}\left(
I\right) }\mathbb{F}_{I^{\prime }}^{\sigma ,\mathbf{b}}f=\left(
\sum_{I^{\prime }\in \mathfrak{C}_{\limfunc{natural}}\left( I\right) }%
\mathbb{F}_{I^{\prime }}^{\sigma ,\mathbf{b}}f\right) -\mathbb{F}%
_{I}^{\sigma ,\mathbf{b}}f,  \label{flat box}
\end{equation}%
where we have removed the averages over broken children from $\square
_{I}^{\sigma ,\mathbf{b}}f$, but left the average over $I$ intact. On any
child $I^{\prime }$ of $I$, the function $\square _{I}^{\sigma ,\flat ,%
\mathbf{b}}f$ is thus a constant multiple of $b_{I}$, and so we have%
\begin{eqnarray}
\square _{I}^{\sigma ,\flat ,\mathbf{b}}f &=&b_{I}\sum_{I^{\prime }\in 
\mathfrak{C}\left( I\right) }\mathbf{1}_{I^{\prime }}E_{I^{\prime }}^{\sigma
}\left( \frac{1}{b_{I}}\square _{I}^{\sigma ,\flat ,\mathbf{b}}f\right)
=b_{I}\ \sum_{I^{\prime }\in \mathfrak{C}\left( I\right) }\mathbf{1}%
_{I^{\prime }}E_{I^{\prime }}^{\sigma }\left( \widehat{\square }_{I}^{\sigma
,\flat ,\mathbf{b}}f\right) ;  \label{flat box hat} \\
\widehat{\square }_{I}^{\sigma ,\flat ,\mathbf{b}}f &\equiv &\sum_{I^{\prime
}\in \mathfrak{C}\left( I\right) }\mathbf{1}_{I^{\prime }}\ E_{I^{\prime
}}^{\sigma }\left( \frac{1}{b_{I}}\square _{I}^{\sigma ,\flat ,\mathbf{b}%
}f\right) ,  \notag
\end{eqnarray}%
where we have denoted the constants in question by the expressions $%
E_{I^{\prime }}^{\sigma }\left( \widehat{\square }_{I}^{\sigma ,\flat ,%
\mathbf{b}}f\right) $, and then defined $\widehat{\square }_{I}^{\sigma
,\flat ,\mathbf{b}}f$ to be the corresponding operator. We record the
precise formula,%
\begin{equation*}
\widehat{\square }_{I}^{\sigma ,\flat ,\mathbf{b}}f=\sum_{I^{\prime }\in 
\mathfrak{C}_{\limfunc{natural}}\left( I\right) }\mathbf{1}_{I^{\prime }}\ %
\left[ \frac{1}{\int_{I^{\prime }}b_{I}d\mu }\int_{I^{\prime }}fd\mu -\frac{1%
}{\int_{I}b_{I}d\mu }\int_{I}fd\mu \right] -\sum_{I^{\prime }\in \mathfrak{C}%
_{\limfunc{broken}}\left( I\right) }\mathbf{1}_{I^{\prime }}\ \left[ \frac{1%
}{\int_{I}b_{I}d\mu }\int_{I}fd\mu \right] .
\end{equation*}%
Thus for $I\in \mathcal{C}_{A}$ we have 
\begin{equation}
\square _{I}^{\sigma ,\flat ,\mathbf{b}}f=b_{A}\sum_{I^{\prime }\in 
\mathfrak{C}\left( I\right) }\mathbf{1}_{I^{\prime }}E_{I^{\prime }}^{\sigma
}\left( \widehat{\square }_{I}^{\sigma ,\flat ,\mathbf{b}}f\right) =b_{A}%
\widehat{\square }_{I}^{\sigma ,\flat ,\mathbf{b}}f,  \label{factor b_A}
\end{equation}%
where the averages $E_{I^{\prime }}^{\sigma }\left( \widehat{\square }%
_{I}^{\sigma ,\flat ,\mathbf{b}}f\right) $ satisfy the following telescoping
property for all $K\in \left( \mathcal{C}_{A}\setminus \left\{ A\right\}
\right) \cup \left( \bigcup_{A^{\prime }\in \mathfrak{C}_{\mathcal{A}}\left(
A\right) }A^{\prime }\right) $ and $L\in \mathcal{C}_{A}$ with $K\subset L$:%
\begin{equation}
\sum_{I:\ \pi K\subset I\subset L}E_{I_{K}}^{\sigma }\left( \widehat{\square 
}_{I}^{\sigma ,\flat ,\mathbf{b}}f\right) =\left\{ 
\begin{array}{ccc}
-E_{L}^{\sigma }\widehat{\mathbb{F}}_{L}^{\sigma }f & \text{ if } & K\in 
\mathfrak{C}_{\mathcal{A}}\left( A\right) \\ 
E_{K}^{\sigma }\widehat{\mathbb{F}}_{K}^{\sigma }f-E_{L}^{\sigma }\widehat{%
\mathbb{F}}_{L}^{\sigma }f & \text{ if } & K\in \mathcal{C}_{A}%
\end{array}%
\right. ,  \label{telescoping}
\end{equation}%
where $\widehat{\mathbb{F}}_{K}^{\sigma }$ is defined in (\ref{F hat})
above. Indeed, recalling that $I_{K}$ denotes the child of $I$ that contains 
$K$, this is evident if we write%
\begin{equation*}
\mathbf{1}_{K}\ b_{A}E_{I_{K}}^{\sigma }\left( \widehat{\square }%
_{I}^{\sigma ,\flat ,\mathbf{b}}f\right) =\mathbf{1}_{K}\ \square
_{I}^{\sigma ,\flat ,\mathbf{b}}f=\mathbf{1}_{K}\ \left( \square
_{I}^{\sigma ,\mathbf{b}}f-\sum_{I^{\prime }\in \mathfrak{C}_{\limfunc{broken%
}}\left( I\right) }\mathbf{1}_{I^{\prime }}\mathbb{F}_{I^{\prime }}^{\sigma ,%
\mathbf{b}}f\right) ,
\end{equation*}%
and use the telescoping properties of $\left\{ \square _{I}^{\sigma ,\mathbf{%
b}}\right\} _{I\in \mathcal{D}}$, together with the fact that the set 
\begin{equation*}
\left\{ I^{\prime }:I^{\prime }\in \mathfrak{C}_{\limfunc{broken}}\left(
I\right) \text{ for some }I\in \mathcal{C}_{A}\right\} =\mathfrak{C}_{%
\mathcal{A}}\left( A\right)
\end{equation*}%
is pairwise disjoint and lies beneath the corona $\mathcal{C}_{A}$. Of
course we can similarly define $\bigtriangleup _{I}^{\sigma ,\flat ,\mathbf{b%
}}$.

Finally, in analogy with the broken differences $\bigtriangleup _{Q,\limfunc{%
broken}}^{\mu ,\pi ,\mathbf{b}}$ and $\square _{Q,\limfunc{broken}}^{\mu
,\pi ,\mathbf{b}}$ introduced above, we define%
\begin{equation}
\bigtriangleup _{I,\limfunc{broken}}^{\mu ,\flat ,\mathbf{b}}f\equiv
\sum_{I^{\prime }\in \mathfrak{C}_{\limfunc{broken}}\left( I\right) }\mathbb{%
E}_{I^{\prime }}^{\sigma ,\mathbf{b}}f\text{ and }\square _{I,\limfunc{broken%
}}^{\mu ,\flat ,\mathbf{b}}f\equiv \sum_{I^{\prime }\in \mathfrak{C}_{%
\limfunc{broken}}\left( I\right) }\mathbb{F}_{I^{\prime }}^{\sigma ,\mathbf{b%
}}f\ ,  \label{def flat broken}
\end{equation}%
so that%
\begin{equation}
\bigtriangleup _{I}^{\mu ,\mathbf{b}}=\bigtriangleup _{I}^{\mu ,\flat ,%
\mathbf{b}}+\bigtriangleup _{I,\limfunc{broken}}^{\mu ,\flat ,\mathbf{b}}%
\text{ and }\square _{I}^{\mu ,\mathbf{b}}=\square _{I}^{\mu ,\flat ,\mathbf{%
b}}+\square _{I,\limfunc{broken}}^{\mu ,\flat ,\mathbf{b}}\ .
\label{flat broken}
\end{equation}%
These modified differences and the identities (\ref{factor b_A}) and (\ref%
{telescoping}) play a useful role in the analysis of the nearby and
paraproduct forms.

\subsection{Weak Riesz inequalities\label{two sided Riesz}}

We begin with a strengthening of the upper frame inequality for dual
martingale pseudoprojections,%
\begin{equation*}
\left\Vert f\right\Vert _{L^{2}\left( \mu \right) }^{2}\lesssim \sum_{Q\in 
\mathcal{D}}\left\Vert \square _{Q}^{\mu ,\mathbf{b}}f\right\Vert
_{L^{2}\left( \mu \right) }^{2}+\sum_{Q\in \mathcal{D}}\left\Vert \widehat{%
\bigtriangledown }_{Q}^{\mu }f\right\Vert _{L^{2}\left( \mu \right) }^{2},
\end{equation*}%
and then proceed to consider the corresponding lower inequalities and
martingale versions as well. We refer to these strengthened inequalities as
weak Riesz inequalities for the following reason. A family of
pseudoprojections $\left\{ \Psi _{Q}^{\mu ,\mathbf{b}}f\right\} _{Q\in 
\mathcal{D}}$ is said to be a \emph{Riesz basis} for $L^{2}\left( \mu
\right) $ if for all subsets $\mathcal{B}\subset \mathcal{D}$ of the dyadic
grid we have%
\begin{equation*}
\sum_{Q\in \mathcal{B}}\left\Vert \Psi _{Q}^{\mu ,\mathbf{b}}f\right\Vert
_{L^{2}\left( \mu \right) }^{2}\lesssim \left\Vert \sum_{Q\in \mathcal{B}%
}\Psi _{Q}^{\mu ,\mathbf{b}}f\right\Vert _{L^{2}\left( \mu \right)
}^{2}\lesssim \sum_{Q\in \mathcal{B}}\left\Vert \Psi _{Q}^{\mu ,\mathbf{b}%
}f\right\Vert _{L^{2}\left( \mu \right) }^{2}\ ,\ \ \ \ \ f\in L^{2}\left(
\mu \right) .
\end{equation*}%
We refer to the left (respectively right) hand inequality above as the lower
(respectively upper) Riesz inequality. The families $\left\{ \square
_{Q}^{\mu ,\mathbf{b}}f\right\} _{Q\in \mathcal{D}}$ and $\left\{
\bigtriangleup _{Q}^{\mu ,\mathbf{b}}f\right\} _{Q\in \mathcal{D}}$ can 
\emph{fail} to be a Riesz basis for $L^{2}\left( \mu \right) $, but in the
case $p=\infty $, we show that each of these families is a Riesz basis in a
certain weak sense,\ involving Carleson averaging operators, that is made
precise below.

\begin{definition}
\label{Psi op}For any subset $\mathcal{B}$ of the grid $\mathcal{D}$, and
any sequence of real numbers $\mathbf{\lambda }=\left\{ \lambda _{I}\right\}
_{I\in \mathcal{B}}$, define the linear operator%
\begin{equation*}
\Psi _{\mathcal{B},\mathbf{\lambda }}^{\mu ,\mathbf{b}}f\equiv \sum_{I\in 
\mathcal{B}}\lambda _{I}\square _{I}^{\mu ,\mathbf{b}}f,
\end{equation*}%
which, with abuse of notation, we will refer to as a `pseudoprojection'. In
the event that all $\lambda _{I}=1$, we denote the operator by $\Psi _{%
\mathcal{B}}^{\mu ,\mathbf{b}}$, and note that, despite the fact it is a sum
of dual martingale averages, it is typically \textbf{not} a projection on $%
L^{2}\left( \mu \right) $.
\end{definition}

Note that the failure of $\square _{I}^{\mu ,\mathbf{b}}$ to be a projection
in general is what motivated the introduction of the projections $\square
_{I}^{\mu ,\pi ,\mathbf{b}}$ above. These projections have already played a
significant role in the proof of the Monotonicity Lemma earlier, and will
continue to play a role in other `duality' situations below.

\subsubsection{An upper weak Riesz inequality}

\begin{proposition}
\label{half Riesz}Suppose that $\mathbf{b}$ is an $\infty $-weakly $\mu $%
-controlled accretive family on a grid $\mathcal{D}$. Then we have the
following `$\square _{I}^{\mu ,\mathbf{b}}$-upper weak Riesz' inequality:%
\begin{eqnarray*}
&&\left\Vert \Psi _{\mathcal{B},\mathbf{\lambda }}^{\mu ,\mathbf{b}%
}f\right\Vert _{L^{2}\left( \mu \right) }^{2}\leq C\left\Vert \mathbf{%
\lambda }\right\Vert _{\infty }^{2}\left( \sum_{I\in \mathcal{B}}\left\Vert
\square _{I}^{\mu ,\mathbf{b}}f\right\Vert _{L^{2}\left( \mu \right)
}^{2}+\sum_{I\in \mathcal{B}}\left\Vert \widehat{\bigtriangledown }_{I}^{\mu
}f\right\Vert _{L^{2}\left( \mu \right) }^{2}\right) \\
&&\ \ \ \ \ \text{for all }f\in L^{2}\left( \mu \right) \text{ and all
subsets }\mathcal{B}\text{ of the grid }\mathcal{D}\text{ and all sequences }%
\mathbf{\lambda },
\end{eqnarray*}%
where $\left\Vert \mathbf{\lambda }\right\Vert _{\infty }\equiv \sup_{I\in 
\mathcal{B}}\left\vert \lambda _{I}\right\vert $ and the positive constant $%
C $ is independent of the subset $\mathcal{B}$. In particular, the
pseudoprojection $\Psi _{\mathcal{B},\mathbf{\lambda }}^{\mu ,\mathbf{b}}$
is a bounded linear operator on $L^{2}\left( \mu \right) $ if $\left\Vert 
\mathbf{\lambda }\right\Vert _{\infty }<\infty $:%
\begin{equation}
\left\Vert \Psi _{\mathcal{B},\mathbf{\lambda }}^{\mu ,\mathbf{b}%
}f\right\Vert _{L^{2}\left( \mu \right) }^{2}\leq C\left\Vert \mathbf{%
\lambda }\right\Vert _{\infty }^{2}\left\Vert f\right\Vert _{L^{2}\left( \mu
\right) }^{2}\ .  \label{Psi bound}
\end{equation}
\end{proposition}

\begin{proof}
We may suppose that the subset $\mathcal{B}$ is finite provided the
estimates we get are independent of the size of $\mathcal{B}$. Now let $%
g=\Psi _{\mathcal{B},\mathbf{\lambda }}^{\mu ,\mathbf{b}}f=\sum_{I\in 
\mathcal{B}}\lambda _{I}\square _{I}^{\mu ,\mathbf{b}}f$. Then from (\ref%
{square of delta}) we have%
\begin{eqnarray*}
\left\Vert g\right\Vert _{L^{2}\left( \mu \right) }^{2} &=&\int \left(
\sum_{I\in \mathcal{B}}\lambda _{I}\square _{I}^{\mu ,\mathbf{b}}f\right) g\
d\mu =\int \left( \sum_{I\in \mathcal{B}}\lambda _{I}\left[ \square
_{I}^{\mu ,\pi ,\mathbf{b}}\square _{I}^{\mu ,\pi ,\mathbf{b}}f+\square _{I,%
\limfunc{broken}}^{\mu ,\pi ,\mathbf{b}}f\right] \right) g\ d\mu \\
&=&\sum_{I\in \mathcal{B}}\lambda _{I}\int \left( \square _{I}^{\mu ,\pi ,%
\mathbf{b}}f\right) \left( \square _{I}^{\mu ,\pi ,\mathbf{b}}\right) ^{\ast
}g\ d\mu +\sum_{I\in \mathcal{B}}\lambda _{I}\int \left( \square _{I,%
\limfunc{broken}}^{\mu ,\pi ,\mathbf{b}}f\right) g\ d\mu \\
&\lesssim &\left\Vert \mathbf{\lambda }\right\Vert _{\infty }\left(
\sum_{I\in \mathcal{B}}\left\Vert \square _{I}^{\mu ,\pi ,\mathbf{b}%
}f\right\Vert _{L^{2}\left( \mu \right) }^{2}\right) ^{\frac{1}{2}}\left(
\sum_{I\in \mathcal{B}}\left\Vert \bigtriangleup _{I}^{\mu ,\pi ,\mathbf{b}%
}g\right\Vert _{L^{2}\left( \mu \right) }^{2}\right) ^{\frac{1}{2}} \\
&&+\left\Vert \mathbf{\lambda }\right\Vert _{\infty }\left( \sum_{I\in 
\mathcal{B}}\left\Vert \widehat{\bigtriangledown }_{I}^{\mu }f\right\Vert
_{L^{2}\left( \mu \right) }^{2}\right) ^{\frac{1}{2}}\left( \sum_{I\in 
\mathcal{B}}\left\Vert \widehat{\bigtriangledown }_{I}^{\mu }g\right\Vert
_{L^{2}\left( \mu \right) }^{2}\right) ^{\frac{1}{2}},
\end{eqnarray*}%
where we have used (\ref{F est}) in the last line. Now using $\bigtriangleup
_{I}^{\mu ,\pi ,\mathbf{b}}=\bigtriangleup _{I}^{\mu ,\mathbf{b}%
}-\bigtriangleup _{I,\limfunc{broken}}^{\mu ,\pi ,\mathbf{b}}$ and $%
\left\vert \bigtriangleup _{I,\limfunc{broken}}^{\mu ,\pi ,\mathbf{b}%
}g\right\vert \leq C_{\mathbf{b}}\widehat{\nabla }_{I}^{\mu }\left\vert
g\right\vert $, and $\square _{I}^{\mu ,\pi ,\mathbf{b}}=\square _{I}^{\mu ,%
\mathbf{b}}-\square _{I,\limfunc{broken}}^{\mu ,\pi ,\mathbf{b}}$, together
with the lower frame inequalities in Proposition \ref{dual frame}\ (note
that we do \emph{not} use lower \emph{Riesz} inequalities for $%
\bigtriangleup $!)%
\begin{eqnarray*}
\sum_{I\in \mathcal{B}}\left\Vert \bigtriangleup _{I}^{\mu ,\mathbf{b}%
}g\right\Vert _{L^{2}\left( \mu \right) }^{2} &\leq &\sum_{I\in \mathcal{D}%
}\left\Vert \bigtriangleup _{I}^{\mu ,\mathbf{b}}g\right\Vert _{L^{2}\left(
\mu \right) }^{2}\lesssim \left\Vert g\right\Vert _{L^{2}\left( \mu \right)
}^{2}\ , \\
\sum_{I\in \mathcal{B}}\left\Vert \widehat{\bigtriangledown }_{I}^{\mu
}\left\vert g\right\vert \right\Vert _{L^{2}\left( \mu \right) }^{2} &\leq
&\sum_{I\in \mathcal{D}}\left\Vert \widehat{\bigtriangledown }_{I}^{\mu
}\left\vert g\right\vert \right\Vert _{L^{2}\left( \mu \right) }^{2}\lesssim
\left\Vert g\right\Vert _{L^{2}\left( \mu \right) }^{2}\ ,
\end{eqnarray*}%
we obtain 
\begin{equation*}
\left\Vert g\right\Vert _{L^{2}\left( \mu \right) }^{2}\lesssim \left\Vert 
\mathbf{\lambda }\right\Vert _{\infty }\left( \sum_{I\in \mathcal{B}%
}\left\Vert \square _{I}^{\mu ,\mathbf{b}}f\right\Vert _{L^{2}\left( \mu
\right) }^{2}+\sum_{I\in \mathcal{B}}\left\Vert \widehat{\bigtriangledown }%
_{I}^{\mu }f\right\Vert _{L^{2}\left( \mu \right) }^{2}\right) ^{\frac{1}{2}%
}\left\Vert g\right\Vert _{L^{2}\left( \mu \right) }\ ,
\end{equation*}%
which gives the desired upper weak Riesz inequality upon dividing through by
the finite positive number $\left\Vert g\right\Vert _{L^{2}\left( \mu
\right) }$, and then squaring the resulting inequality.
\end{proof}

An analogous argument yields the next proposition. Recall that $\Psi _{%
\mathcal{B},\mathbf{\lambda }}^{\mu ,\mathbf{b}}f\equiv \sum_{I\in \mathcal{B%
}}\lambda _{I}\square _{I}^{\mu ,\mathbf{b}}f$ and so 
\begin{equation*}
\left( \Psi _{\mathcal{B},\mathbf{\lambda }}^{\mu ,\mathbf{b}}\right) ^{\ast
}f\equiv \sum_{I\in \mathcal{B}}\lambda _{I}\left( \square _{I}^{\mu ,%
\mathbf{b}}\right) ^{\ast }f=\sum_{I\in \mathcal{B}}\lambda
_{I}\bigtriangleup _{I}^{\mu ,\mathbf{b}}f.
\end{equation*}

\begin{proposition}
\label{half Riesz dual}Suppose that $\mathbf{b}$ is an $\infty $-weakly $\mu 
$-controlled accretive family. Then we have the `$\bigtriangleup _{I}^{\mu ,%
\mathbf{b}}$-upper weak Riesz' inequality:%
\begin{eqnarray*}
&&\left\Vert \left( \Psi _{\mathcal{B},\mathbf{\lambda }}^{\mu ,\mathbf{b}%
}\right) ^{\ast }f\right\Vert _{L^{2}\left( \mu \right) }^{2}\leq
C\left\Vert \mathbf{\lambda }\right\Vert _{\infty }^{2}\left( \sum_{I\in 
\mathcal{B}}\left\Vert \bigtriangleup _{I}^{\mu ,\mathbf{b}}f\right\Vert
_{L^{2}\left( \mu \right) }^{2}+\sum_{I\in \mathcal{B}}\left\Vert \widehat{%
\bigtriangledown }_{I}^{\mu }f\right\Vert _{L^{2}\left( \mu \right)
}^{2}\right) , \\
&&\ \ \ \ \ \text{for all }f\in L^{2}\left( \mu \right) \text{ and all
subsets }\mathcal{B}\text{ of the grid }\mathcal{D}\text{ and all sequences }%
\mathbf{\lambda },
\end{eqnarray*}%
and where the positive constant $C$ is independent of the subset $\mathcal{B}
$ and the sequence $\mathbf{\lambda }$. In particular, the pseudoprojection $%
\left( \Psi _{\mathcal{B},\mathbf{\lambda }}^{\mu ,\mathbf{b}}\right) ^{\ast
}$ is a bounded linear operator on $L^{2}\left( \mu \right) $:%
\begin{equation}
\left\Vert \left( \Psi _{\mathcal{B},\mathbf{\lambda }}^{\mu ,\mathbf{b}%
}\right) ^{\ast }f\right\Vert _{L^{2}\left( \mu \right) }^{2}\leq
C\left\Vert \mathbf{\lambda }\right\Vert _{\infty }^{2}\left\Vert
f\right\Vert _{L^{2}\left( \mu \right) }^{2}\ .  \label{Psi * bound}
\end{equation}
\end{proposition}

\begin{proof}
We may again suppose that the subset $\mathcal{B}$ is finite provided the
estimates we get are independent of the size of $\mathcal{B}$. Now let $%
g=\sum_{I\in \mathcal{B}}\lambda _{I}\bigtriangleup _{I}^{\mu ,\mathbf{b}}f$%
. Recall from (\ref{square of delta}) that $\bigtriangleup _{Q}^{\mu ,%
\mathbf{b}}=\bigtriangleup _{Q}^{\mu ,\pi ,\mathbf{b}}\bigtriangleup
_{Q}^{\mu ,\pi ,\mathbf{b}}+\bigtriangleup _{Q,\limfunc{broken}}^{\mu ,%
\mathbf{b}}$, and hence we have%
\begin{eqnarray*}
\left\Vert g\right\Vert _{L^{2}\left( \mu \right) }^{2} &=&\int \left(
\sum_{I\in \mathcal{B}}\lambda _{I}\bigtriangleup _{I}^{\mu ,\mathbf{b}%
}f\right) g\ d\mu =\int \left( \sum_{I\in \mathcal{B}}\lambda _{I}\left[
\bigtriangleup _{I}^{\mu ,\pi ,\mathbf{b}}\bigtriangleup _{I}^{\mu ,\pi ,%
\mathbf{b}}f+\bigtriangleup _{I,\limfunc{broken}}^{\mu ,\pi ,\mathbf{b}}f%
\right] \right) g\ d\mu \\
&=&\sum_{I\in \mathcal{B}}\lambda _{I}\int \left( \bigtriangleup _{I}^{\mu
,\pi ,\mathbf{b}}f\right) \left( \bigtriangleup _{I}^{\mu ,\pi ,\mathbf{b}%
}\right) ^{\ast }g\ d\mu +\sum_{I\in \mathcal{B}}\lambda _{I}\int \left(
\bigtriangleup _{I,\limfunc{broken}}^{\mu ,\pi ,\mathbf{b}}f\right) g\ d\mu
\\
&\lesssim &\left\Vert \mathbf{\lambda }\right\Vert _{\infty }\left(
\sum_{I\in \mathcal{B}}\left\Vert \bigtriangleup _{I}^{\mu ,\pi ,\mathbf{b}%
}f\right\Vert _{L^{2}\left( \mu \right) }^{2}\right) ^{\frac{1}{2}}\left(
\sum_{I\in \mathcal{B}}\left\Vert \square _{I}^{\mu ,\pi ,\mathbf{b}%
}g\right\Vert _{L^{2}\left( \mu \right) }^{2}\right) ^{\frac{1}{2}} \\
&&+\left\Vert \mathbf{\lambda }\right\Vert _{\infty }\left( \sum_{I\in 
\mathcal{B}}\left\Vert \widehat{\bigtriangledown }_{I}^{\mu }f\right\Vert
_{L^{2}\left( \mu \right) }^{2}\right) ^{\frac{1}{2}}\left( \sum_{I\in 
\mathcal{B}}\left\Vert \nabla _{I}^{\mu }g\right\Vert _{L^{2}\left( \mu
\right) }^{2}\right) ^{\frac{1}{2}}.
\end{eqnarray*}%
Now we continue as in the proof of Propostion \ref{half Riesz}, but using
instead the lower frame inequality (\ref{low frame}), to obtain%
\begin{equation*}
\left\Vert g\right\Vert _{L^{2}\left( \mu \right) }^{2}\lesssim \left\Vert 
\mathbf{\lambda }\right\Vert _{\infty }\left( \sum_{I\in \mathcal{B}%
}\left\Vert \bigtriangleup _{I}^{\mu ,\mathbf{b}}f\right\Vert _{L^{2}\left(
\mu \right) }^{2}+\sum_{I\in \mathcal{B}}\left\Vert \widehat{%
\bigtriangledown }_{I}^{\mu }f\right\Vert _{L^{2}\left( \mu \right)
}^{2}\right) \left\Vert g\right\Vert _{L^{2}\left( \mu \right) }\ ,
\end{equation*}%
which completes the proof of Propostion \ref{half Riesz dual} upon dividing
through by $\left\Vert g\right\Vert _{L^{2}\left( \mu \right) }$ and
squaring.
\end{proof}

\begin{remark}
The boundedness of the pseudoprojections $\Psi _{\mathcal{B}}^{\mu ,\mathbf{b%
}}$ and $\left( \Psi _{\mathcal{B}}^{\mu ,\mathbf{b}}\right) ^{\ast }$ on $%
L^{2}\left( \mu \right) $ (where the absence of the sequence $\mathbf{%
\lambda }$ in the subscript implies all $\lambda _{I}=1$) given by (\ref{Psi
bound}) and (\ref{Psi * bound}), can fail for a $2$-weakly $\mu $-controlled
accretive family on a grid $\mathcal{D}$. Indeed, if (\ref{Psi * bound})
holds, then for a $2$-weakly $\mu $-controlled accretive family,%
\begin{equation*}
\left\Vert \sum_{I\in \mathcal{B}_{+}}\bigtriangleup _{I}^{\mu ,\mathbf{b}%
}f-\sum_{I\in \mathcal{B}_{-}}\bigtriangleup _{I}^{\mu ,\mathbf{b}%
}f\right\Vert _{L^{2}\left( \mu \right) }^{2}\leq 2\left\Vert \sum_{I\in 
\mathcal{B}_{+}}\bigtriangleup _{I}^{\mu ,\mathbf{b}}f\right\Vert
_{L^{2}\left( \mu \right) }^{2}+2\left\Vert \sum_{I\in \mathcal{B}%
_{-}}\bigtriangleup _{I}^{\mu ,\mathbf{b}}f\right\Vert _{L^{2}\left( \mu
\right) }^{2}\leq C\left\Vert f\right\Vert _{L^{2}\left( \mu \right) }^{2}
\end{equation*}%
holds for all decompositions of $\mathcal{B}$ into a disjoint union $%
\mathcal{B=B}_{+}\overset{\cdot }{\cup }\mathcal{B}_{-}$. Then%
\begin{equation*}
\sum_{I\in \mathcal{B}}\left\Vert \bigtriangleup _{I}^{\mu ,\mathbf{b}%
}f\right\Vert _{L^{2}\left( \mu \right) }^{2}=\mathbb{E}_{\pm }\left\Vert
\sum_{I\in \mathcal{B}}\pm \bigtriangleup _{I}^{\mu ,\mathbf{b}}f\right\Vert
_{L^{2}\left( \mu \right) }^{2}\leq \mathbb{E}_{\pm }C\left\Vert
f\right\Vert _{L^{2}\left( \mu \right) }^{2}=C\left\Vert f\right\Vert
_{L^{2}\left( \mu \right) }^{2}\ ,
\end{equation*}%
which contradicts the example of Hyt\"{o}nen and Martikainen in \cite[%
Section 3.9]{HyMa} if $p=2$. Thus in order to prove a two weight local $Tb$
theorem for $p=2$, one \textbf{cannot} appeal in general to the boundedness
of pseudoprojections $\left( \Psi _{\mathcal{B},\mathbf{\lambda }}^{\mu ,%
\mathbf{b}}\right) ^{\ast }=\sum_{I\in \mathcal{B}}\lambda
_{I}\bigtriangleup _{I}^{\mu ,\mathbf{b}}$ on $L^{2}\left( \mu \right) $
even when $\lambda _{I}=1$ for all $I\in \mathcal{B}$.
\end{remark}

\subsubsection{A lower weak Riesz inequality}

The next proposition also assumes an $\infty $-weakly $\mu $-controlled
accretive family on a grid $\mathcal{D}$. For a subset $\mathcal{B}$ of a
dyadic grid $\mathcal{D}$, let 
\begin{equation}
\mathsf{P}_{\mathcal{B}}^{\mu }f=\sum_{Q\in \mathcal{B}}\bigtriangleup
_{Q}^{\mu }f  \label{Haar proj}
\end{equation}%
denote the orthogonal projection of $f$ onto the closed linear span of the
collection $\left\{ \bigtriangleup _{Q}^{\mu }\right\} _{Q\in \mathcal{B}}$
of Haar projections $\bigtriangleup _{Q}^{\mu }$ with $Q\in \mathcal{B}$. We
obtain an appropriate form of a weak lower Riesz inequality for the dual
martingale differences $\square _{Q}^{\mu ,\mathbf{b}}$ since for these
operators we no longer need the disruptive device of introducing the
function $b$ in (\ref{def b}) above.

\begin{proposition}
\label{reverse half Riesz dual}Suppose that $\mathbf{b}$ is an $\infty $%
-weakly $\mu $-controlled accretive family. Then we have a `weak lower
Riesz' inequality for dual martingale differences:%
\begin{eqnarray*}
&&\sum_{Q\in \mathcal{B}}\left\Vert \square _{Q}^{\mu ,\mathbf{b}%
}f\right\Vert _{L^{2}\left( \mu \right) }^{2}\leq C\left( \left\Vert \mathsf{%
P}_{\mathcal{B}}^{\mu }f\right\Vert _{L^{2}\left( \mu \right)
}^{2}+\sum_{Q\in \mathcal{B}}\left\Vert \widehat{\bigtriangledown }_{Q}^{\mu
}f\right\Vert _{L^{2}\left( \mu \right) }^{2}+\sum_{Q\in \mathcal{B}}\gamma
_{Q}\left\vert E_{Q}^{\mu }f\right\vert ^{2}\right) , \\
&&\text{for all }f\in L^{2}\left( \mu \right) \text{ and all subsets }%
\mathcal{B}\text{ of the grid }\mathcal{D}\text{, and where }\left\Vert 
\mathsf{P}_{\mathcal{B}}^{\mu }f\right\Vert _{L^{2}\left( \mu \right)
}^{2}=\sum_{Q\in \mathcal{B}}\left\Vert \bigtriangleup _{Q}^{\mu
}f\right\Vert _{L^{2}\left( \mu \right) }^{2}\ .
\end{eqnarray*}%
and where the positive constant $C$ depends only on the accretivity
constants, but is \emph{independent} of the subset $\mathcal{B}$ and the
testing family $\mathbf{b}$. Here the coefficients $\left\{ \gamma
_{Q}\right\} _{Q\in \mathcal{D}}$ form a Carleson sequence indexed by $%
\mathcal{D}$, i.e.%
\begin{equation*}
\sum_{Q\in \mathcal{D}:\ Q\subset P}\gamma _{Q}\leq C\left\vert P\right\vert
_{\mu }\ ,\ \ \ \ \ \text{for all }P\in \mathcal{D}.
\end{equation*}
\end{proposition}

The third term on the right hand side above is additive in $\mathcal{B}$
and, by the Carleson embedding theorem, satisfies%
\begin{equation*}
\sum_{Q\in \mathcal{D}}\gamma _{Q}\left\vert E_{Q}^{\mu }f\right\vert
^{2}\lesssim \left\Vert f\right\Vert _{L^{2}\left( \mu \right) }^{2}\ .
\end{equation*}

\begin{proof}
The main modifications in the proof of Proposition \ref{dual frame} that are
needed here, are that we use $\square _{Q}^{\mu ,\mathbf{b}}$ in place of $%
\bigtriangleup _{Q}^{\mu ,\mathbf{b}}$, which results in the testing
functions $b_{Q}$ appearing outside the integrals rather than inside the
integrals, and that we restrict the sums over $Q$ to $\mathcal{B}$, which
results in the presence of the term $\sum_{Q\in \mathcal{B}}\gamma
_{Q}\left\vert E_{Q}^{\mu }f\right\vert ^{2}$ on the right hand side above.
Since we are working with $\square _{Q}^{\mu ,\mathbf{b}}$ we will not need
the extra complications arising from the introduction of the function $b$ in
(\ref{def b}). With these modifications in mind, we now describe the
estimates we obtain for the terms analogous to $II_{A}$, $III_{A}$ and $%
IV_{A}$ in the proof of Proposition \ref{dual frame}. Given $A\in \mathcal{A}
$, we begin\ with%
\begin{eqnarray*}
\sum_{Q\in \mathcal{B}\cap \mathcal{C}_{A}}\left\Vert \square _{Q}^{\mu ,%
\mathbf{b}}f\right\Vert _{L^{2}\left( \mu \right) }^{2} &=&\sum_{Q\in 
\mathcal{B}\cap \mathcal{C}_{A}}\sum_{Q^{\prime }\in \mathfrak{C}\left(
Q\right) }\int_{Q^{\prime }}\left\vert \mathbb{F}_{Q^{\prime }}^{\mu ,%
\mathbf{b}}f\left( x\right) -\mathbb{F}_{Q}^{\mu ,\mathbf{b}}f\left(
x\right) \right\vert ^{2}d\mu \left( x\right) \\
&=&\sum_{Q\in \mathcal{B}\cap \mathcal{C}_{A}}\sum_{Q^{\prime }\in \mathfrak{%
C}_{\limfunc{natural}}\left( Q\right) }\left\vert \frac{b_{A}\int_{Q^{\prime
}}fd\mu }{\int_{Q^{\prime }}b_{A}d\mu }-\frac{b_{A}\int_{Q}fd\mu }{%
\int_{Q}b_{A}d\mu }\right\vert ^{2}\left\vert Q^{\prime }\right\vert _{\mu }
\\
&&+\sum_{Q\in \mathcal{B}\cap \mathcal{C}_{A}}\sum_{Q^{\prime }\in \mathfrak{%
C}_{\limfunc{broken}}\left( Q\right) }\left\vert \frac{b_{Q^{\prime
}}\int_{Q^{\prime }}fd\mu }{\int_{Q^{\prime }}b_{Q^{\prime }}d\mu }-\frac{%
b_{A}\int_{Q}fd\mu }{\int_{Q}b_{A}d\mu }\right\vert ^{2}\left\vert Q^{\prime
}\right\vert _{\mu } \\
&\equiv &I_{A}+II_{A}.
\end{eqnarray*}%
To estimate term $II_{A}$ we write%
\begin{eqnarray}
&&\left\vert \frac{b_{Q^{\prime }}\int_{Q^{\prime }}fd\mu }{\int_{Q^{\prime
}}b_{Q^{\prime }}d\mu }-\frac{b_{A}\int_{Q}fd\mu }{\int_{Q}b_{A}d\mu }%
\right\vert  \label{analogue'} \\
&=&\frac{\left\vert \left( \frac{1}{\left\vert Q^{\prime }\right\vert _{\mu }%
}b_{Q^{\prime }}\int_{Q^{\prime }}fd\mu \right) \left( \frac{1}{\left\vert
Q\right\vert _{\mu }}\int_{Q}b_{A}d\mu \right) -\left( \frac{1}{\left\vert
Q\right\vert _{\mu }}b_{A}\int_{Q}fd\mu \right) \left( \frac{1}{\left\vert
Q^{\prime }\right\vert _{\mu }}\int_{Q^{\prime }}b_{Q^{\prime }}d\mu \right)
\right\vert }{\left\vert \frac{1}{\left\vert Q^{\prime }\right\vert _{\mu }}%
\int_{Q^{\prime }}b_{Q^{\prime }}d\mu \right\vert \left\vert \frac{1}{%
\left\vert Q\right\vert _{\mu }}\int_{Q}b_{A}d\mu \right\vert }  \notag \\
&\lesssim &\left\vert b_{Q^{\prime }}\left( E_{Q^{\prime }}^{\mu }f\right)
\left( E_{Q}^{\mu }b_{A}\right) -b_{A}\left( E_{Q}^{\mu }f\right) \left(
E_{Q^{\prime }}^{\mu }b_{Q^{\prime }}\right) \right\vert  \notag \\
&=&\left\vert b_{Q^{\prime }}\left( E_{Q^{\prime }}^{\mu }f\right) \left(
E_{Q}^{\mu }b_{A}\right) -b_{A}\left( E_{Q^{\prime }}^{\mu }f\right) \left(
E_{Q^{\prime }}^{\mu }b_{Q^{\prime }}\right) -\left[ b_{A}\left( E_{Q}^{\mu
}f\right) -b_{A}\left( E_{Q^{\prime }}^{\mu }f\right) \right] \left(
E_{Q^{\prime }}^{\mu }b_{Q^{\prime }}\right) \right\vert  \notag \\
&\lesssim &\left\vert b_{Q^{\prime }}\left( E_{Q^{\prime }}^{\mu }f\right)
\right\vert +\left\vert b_{A}\left( E_{Q^{\prime }}^{\mu }f\right)
\right\vert +\left\vert b_{A}\left( E_{Q}^{\mu }f\right) -b_{A}\left(
E_{Q^{\prime }}^{\mu }f\right) \right\vert ,  \notag \\
&\lesssim &\left\vert E_{Q^{\prime }}^{\mu }f\right\vert +\left\vert
E_{Q^{\prime }}^{\mu }f\right\vert +\left\vert \left( E_{Q}^{\mu }f\right)
-\left( E_{Q^{\prime }}^{\mu }f\right) \right\vert ,  \notag
\end{eqnarray}%
and then 
\begin{equation*}
II_{A}\lesssim \sum_{Q\in \mathcal{B}\cap \mathcal{C}_{A}}\sum_{Q^{\prime
}\in \mathfrak{C}_{\limfunc{broken}}\left( Q\right) }\left( \left\vert
E_{Q^{\prime }}^{\mu }f\right\vert ^{2}+\left\vert E_{Q}^{\mu }f\right\vert
^{2}\right) \left\vert Q^{\prime }\right\vert _{\mu }\approx \sum_{Q\in 
\mathcal{B}\cap \mathcal{C}_{A}}\left\Vert \widehat{\bigtriangledown }%
_{Q}^{\mu }f\right\Vert _{L^{2}\left( \mu \right) }^{2}\ .
\end{equation*}%
Now we write%
\begin{eqnarray*}
I_{A} &=&\sum_{Q\in \mathcal{B}\cap \mathcal{C}_{A}}\sum_{Q^{\prime }\in 
\mathfrak{C}_{\limfunc{natural}}\left( Q\right) }\left\vert \frac{%
b_{A}\int_{Q^{\prime }}fd\mu }{\int_{Q^{\prime }}b_{A}d\mu }-\frac{%
b_{A}\int_{Q}fd\mu }{\int_{Q}b_{A}d\mu }\right\vert ^{2}\left\vert Q^{\prime
}\right\vert _{\mu } \\
&=&\sum_{Q\in \mathcal{B}\cap \mathcal{C}_{A}}\sum_{Q^{\prime }\in \mathfrak{%
C}_{\limfunc{natural}}\left( Q\right) }\left\vert \frac{\left(
b_{A}\int_{Q^{\prime }}fd\mu \right) \left( \int_{Q}b_{A}d\mu \right)
-\left( \int_{Q^{\prime }}b_{A}d\mu \right) \left( b_{A}\int_{Q}fd\mu
\right) }{\left( \int_{Q^{\prime }}b_{A}d\mu \right) \left(
\int_{Q}b_{A}d\mu \right) }\right\vert ^{2}\left\vert Q^{\prime }\right\vert
_{\mu } \\
&\lesssim &\sum_{Q\in \mathcal{B}\cap \mathcal{C}_{A}}\sum_{Q^{\prime }\in 
\mathfrak{C}_{\limfunc{natural}}\left( Q\right) }\left\vert E_{Q}^{\mu
}b_{A}\right\vert ^{2}\left\vert E_{Q^{\prime }}^{\mu }f-E_{Q}^{\mu
}f\right\vert ^{2}\left\vert Q^{\prime }\right\vert _{\mu } \\
&&+\sum_{Q\in \mathcal{B}\cap \mathcal{C}_{A}}\sum_{Q^{\prime }\in \mathfrak{%
C}_{\limfunc{natural}}\left( Q\right) }\left\vert E_{Q}^{\mu
}b_{A}-E_{Q^{\prime }}^{\mu }b_{A}\right\vert ^{2}\left\vert E_{Q}^{\mu
}f\right\vert ^{2}\left\vert Q^{\prime }\right\vert _{\mu } \\
&\equiv &III_{A}+IV_{A}.
\end{eqnarray*}%
Then we use the following stronger form of an inequality used in the proof
of the unweighted square function estimate in Lemma \ref{unweighted square},
namely 
\begin{equation}
\sum_{Q\in \mathcal{B}\cap \mathcal{C}_{A}}\sum_{Q^{\prime }\in \mathfrak{C}%
_{\limfunc{natural}}\left( Q\right) }\left\vert Q^{\prime }\right\vert _{\mu
}\left\vert E_{Q^{\prime }}^{\mu }h-E_{Q}^{\mu }h\right\vert ^{2}\lesssim
\sum_{Q\in \mathcal{B}\cap \mathcal{C}_{A}}\left\Vert \bigtriangleup
_{Q}^{\mu }h\right\Vert _{L^{2}\left( \mu \right) }^{2}\lesssim \left\Vert 
\mathsf{P}_{\mathcal{B}\cap \mathcal{C}_{A}}^{\mu }h\right\Vert
_{L^{2}\left( \mu \right) }^{2},  \label{stronger unweighted}
\end{equation}%
to dominate term $III_{A}$ by%
\begin{equation*}
III_{A}\lesssim \left\Vert b_{A}\mathsf{P}_{\mathcal{B}\cap \mathcal{C}%
_{A}}^{\mu }f\right\Vert _{L^{2}\left( \mu \right) }^{2}\lesssim \left\Vert 
\mathsf{P}_{\mathcal{B}\cap \mathcal{C}_{A}}^{\mu }f\right\Vert
_{L^{2}\left( \mu \right) }^{2}\ ,
\end{equation*}%
and we write term $IV_{A}$ as%
\begin{eqnarray*}
IV_{A} &=&\sum_{Q\in \mathcal{B}\cap \mathcal{C}_{A}}\gamma _{Q}\left\vert
E_{Q}^{\mu }f\right\vert ^{2}\ , \\
\text{where }\gamma _{Q} &\equiv &\sum_{Q^{\prime }\in \mathfrak{C}_{%
\limfunc{natural}}\left( Q\right) }\left\vert E_{Q}^{\mu }b_{A}-E_{Q^{\prime
}}^{\mu }b_{A}\right\vert ^{2}\left\vert Q^{\prime }\right\vert _{\mu },\ \
\ \ \ \text{for }Q\in \mathcal{C}_{A},A\in \mathcal{A\ }.
\end{eqnarray*}%
Now we sum over $A\in \mathcal{A}$ to obtain%
\begin{eqnarray*}
\sum_{Q\in \mathcal{B}}\left\Vert \square _{Q}^{\mu ,\mathbf{b}}f\right\Vert
_{L^{2}\left( \mu \right) }^{2} &=&\sum_{A\in \mathcal{A}}\sum_{Q\in 
\mathcal{B\cap C}_{A}}\left\Vert \square _{Q}^{\mu ,\mathbf{b}}f\right\Vert
_{L^{2}\left( \mu \right) }^{2}\lesssim \sum_{A\in \mathcal{A}}\left(
II_{A}+III_{A}+IV_{A}\right) \\
&\lesssim &\sum_{Q\in \mathcal{B}}\left\Vert \widehat{\bigtriangledown }%
_{Q}^{\mu ,\mathbf{b}}f\right\Vert _{L^{2}\left( \mu \right)
}^{2}+\left\Vert \mathsf{P}_{\mathcal{B}}^{\mu }f\right\Vert _{L^{2}\left(
\mu \right) }^{2}+\sum_{Q\in \mathcal{B}}\gamma _{Q}\left\vert E_{Q}^{\mu
}f\right\vert ^{2}\ ,
\end{eqnarray*}%
where the sequence $\left\{ \gamma _{Q}\right\} _{Q\in \mathcal{D}}$
satisfies the Carleson condition by (\ref{Car cond}) in the proof of
Proposition \ref{dual frame}.
\end{proof}

\begin{remark}
We are unable to obtain a corresponding lower weak Riesz inequality for the
martingale differences $\bigtriangleup _{Q}^{\mu ,\mathbf{b}}$ due to the
need for introducing the function $b$ in (\ref{def b}) as in the proof of
Proposition \ref{dual frame}, which does not interact well with $\mathcal{B}$
- see the argument surrounding (\ref{write both}) in the proof of
Proposition \ref{dual frame}. However, lower weak Riesz inequalities for the
martingale differences $\bigtriangleup _{Q}^{\mu ,\mathbf{b}}$ are not
needed in this paper - in fact, only upper weak Riesz inequalities are
needed for both $\square _{Q}^{\mu ,\mathbf{b}}$ and $\bigtriangleup
_{Q}^{\mu ,\mathbf{b}}$.
\end{remark}

\section{Appendix B:\ Control of functional energy\label{equiv}}

Now we arrive at one of the main propositions used in the proof of our
theorem. This result is proved \emph{independently} of the main theorem, and
only using the results on dual martingale differences established in the
previous appendix. The organization of the proof is almost identical to that
of the corresponding result in\ \cite[pages 128-151]{SaShUr7}, together with
the modifications in \cite[pages 348-360]{SaShUr9} to accommodate common
point masses, but we repeat the organization here with modifications
required for the use of two independent grids, and the appearance of weak
goodness entering through the intervals $J^{\maltese }$. Recall that the
functional energy constant $\mathfrak{F}_{\alpha }=\mathfrak{F}_{\alpha }^{%
\mathbf{b}^{\ast }}\left( \mathcal{D},\mathcal{G}\right) $ in (\ref%
{e.funcEnergy n}), $0\leq \alpha <n$, namely the best constant in the
inequality (see (\ref{def M_r-deep}) below for the definition of $\mathcal{W}%
\left( F\right) $), 
\begin{equation}
\sum_{F\in \mathcal{F}}\sum_{M\in \mathcal{W}\left( F\right) }\left( \frac{%
\mathrm{P}^{\alpha }\left( M,h\sigma \right) }{\left\vert M\right\vert }%
\right) ^{2}\left\Vert \mathsf{Q}_{\mathcal{C}_{F}^{\mathcal{G},\limfunc{%
shift}};M}^{\omega ,\mathbf{b}^{\ast }}x\right\Vert _{L^{2}\left( \omega
\right) }^{\spadesuit 2}\leq \mathfrak{F}_{\alpha }\lVert h\rVert
_{L^{2}\left( \sigma \right) }\,,  \label{fec}
\end{equation}%
depends on the grids $\mathcal{D}$ and $\mathcal{G}$, the goodness parameter 
$\varepsilon >0$ used in the definition of $J^{\maltese }$ through the
shifted corona $\mathcal{C}_{F}^{\mathcal{G},\limfunc{shift}}$, and on the
family of martingale differences $\left\{ \bigtriangleup _{J}^{\omega ,%
\mathbf{b}^{\ast }}\right\} _{J\in \mathcal{G}}$ associated with $x\in
L_{loc}^{2}\left( \omega \right) $, but not on the family of dual martingale
differences $\left\{ \square _{I}^{\sigma ,\mathbf{b}}\right\} _{I\in 
\mathcal{D}}$, since the function $h\in L^{2}\left( \sigma \right) $
appearing in the definition of functional energy is not decomposed as a sum
of pseudoprojections $\square _{I}^{\sigma ,\mathbf{b}}h$. Finally, we
emphasize that the pseudoprojection 
\begin{equation}
\mathsf{Q}_{\mathcal{C}_{F}^{\mathcal{G},\limfunc{shift}};M}^{\omega ,%
\mathbf{b}^{\ast }}\equiv \sum_{J\in \mathcal{C}_{F}^{\mathcal{G},\limfunc{%
shift}}:\ J\subset M}\bigtriangleup _{J}^{\omega ,\mathbf{b}^{\ast }}
\label{def pseudo rest}
\end{equation}%
here uses the shifted restricted corona in%
\begin{eqnarray}
\mathcal{C}_{F}^{\mathcal{G},\limfunc{shift}} &=&\left\{ J\in \mathcal{G}%
:J^{\maltese }\in \mathcal{C}_{F}^{\mathcal{D}}\right\} ,
\label{def shift cor rest} \\
\mathcal{C}_{F}^{\mathcal{G},\limfunc{shift}};K &\equiv &\left\{ J\in 
\mathcal{C}_{F}^{\mathcal{G},\limfunc{shift}}:J\subset K\right\} ,  \notag
\end{eqnarray}%
where $J^{\maltese }$ is defined using the $\limfunc{body}$ of an interval
as in Definition \ref{def sharp cross}, and where we have defined here the
`restriction' $\mathcal{C}_{F}^{\mathcal{G},\limfunc{shift}};K$ to the
interval $K$ of the corona $\mathcal{C}_{F}^{\mathcal{G},\limfunc{shift}}$
(c.f. $\Pi _{2}^{K}\mathcal{P}$\ in Definition \ref{rest K}, which uses the
stronger requirement $J^{\maltese }\subset K$). Moreover, recall from
Notation \ref{nonstandard norm} and the definition of $\nabla _{J}^{\omega }$
in (\ref{Carleson avg op}), that for any subset $\mathcal{H}$ of the grid $%
\mathcal{G}$,%
\begin{equation*}
\left\Vert \mathsf{Q}_{\mathcal{H}}^{\omega ,\mathbf{b}^{\ast }}x\right\Vert
_{L^{2}\left( \omega \right) }^{\spadesuit 2}\equiv \sum_{J\in \mathcal{H}%
}\left\Vert \bigtriangleup _{J}^{\omega ,\mathbf{b}^{\ast }}x\right\Vert
_{L^{2}\left( \omega \right) }^{\spadesuit 2}=\sum_{J\in \mathcal{H}}\left(
\left\Vert \bigtriangleup _{J}^{\omega ,\mathbf{b}^{\ast }}x\right\Vert
_{L^{2}\left( \omega \right) }^{2}+\inf_{z\in \mathbb{R}}\left\Vert \widehat{%
\nabla }_{J}^{\omega }\left( x-z\right) \right\Vert _{L^{2}\left( \omega
\right) }^{2}\right) ,
\end{equation*}%
so that we never need to consider the norm squared $\left\Vert \mathsf{Q}_{%
\mathcal{C}_{F}^{\mathcal{G},\limfunc{shift}};M}^{\omega ,\mathbf{b}^{\ast
}}x\right\Vert _{L^{2}\left( \omega \right) }^{2}$ of the pseudoprojection $%
\mathsf{Q}_{\mathcal{C}_{F}^{\mathcal{G},\limfunc{shift}};M}^{\omega ,%
\mathbf{b}^{\ast }}x$, something for which we have no lower Riesz
inequality. Note moreover that for $J\in \mathcal{G}$ and an arbitrary
interval $K$, we have by the frame inequality in Proposition \ref{dual frame}%
, 
\begin{eqnarray}
\sum_{J\in \mathcal{G}:\ J\subset K}\left\Vert \bigtriangleup _{J}^{\omega ,%
\mathbf{b}^{\ast }}x\right\Vert _{L^{2}\left( \omega \right) }^{2} &\lesssim
&\left\Vert x-m_{K}^{\omega }\right\Vert _{L^{2}\left( \mathbf{1}_{K}\omega
\right) }^{2},  \label{note more} \\
\sum_{J\in \mathcal{G}:\ J\subset K}\inf_{z\in \mathbb{R}}\left\Vert 
\widehat{\nabla }_{J}^{\omega }\left( x-z\right) \right\Vert _{L^{2}\left(
\omega \right) }^{2} &\leq &\sum_{J\in \mathcal{G}:\ J\subset K}\left\Vert 
\widehat{\nabla }_{J}^{\omega }\left\{ \left( x-p\right) \mathbf{1}%
_{K}\left( x\right) \right\} \right\Vert _{L^{2}\left( \omega \right)
}^{2}\lesssim \left\Vert \left( x-p\right) \right\Vert _{L^{2}\left( \mathbf{%
1}_{K}\omega \right) }^{2},\ \ \ p\in K,  \notag
\end{eqnarray}%
where the second line follows from (\ref{Car embed}).

\begin{description}
\item[Important note] If $J\in \mathcal{C}_{F}^{\mathcal{G},\limfunc{shift}}$%
, then in particular $J\Subset _{\mathbf{\rho },\varepsilon }F$ with $%
\mathbf{\rho }=\left[ \frac{3}{\varepsilon }\right] $ by Lemma \ref{good
scale}, and so $J\cap M\neq \emptyset $ for a \emph{unique} $M\in \mathcal{W}%
\left( F\right) $.
\end{description}

We will show that, uniformly in pairs of grids $\mathcal{D}$ and $\mathcal{G}
$, the functional energy constants $\mathfrak{F}_{\alpha }\left( \mathcal{D},%
\mathcal{G}\right) $ in (\ref{e.funcEnergy n}) are controlled by $\mathcal{A}%
_{2}^{\alpha }$, $A_{2}^{\alpha ,\limfunc{punct}}$ and the large energy
constant $\mathfrak{E}_{2}^{\alpha }$ - actually the proof shows that we
have control by the Whitney plugged energy constant as defined in (\ref{def
deep plug}) below. More precisely this is our control of functional energy
proposition.

\begin{proposition}
\label{func ener control}For all grids $\mathcal{D}$ and $\mathcal{G}$, and $%
\varepsilon >0$ sufficiently small, we have%
\begin{eqnarray*}
\mathfrak{F}_{\alpha }^{\mathbf{b}^{\ast }}\left( \mathcal{D},\mathcal{G}%
\right) &\lesssim &\mathfrak{E}_{2}^{\alpha }+\sqrt{\mathcal{A}_{2}^{\alpha }%
}+\sqrt{\mathcal{A}_{2}^{\alpha ,\ast }}+\sqrt{A_{2}^{\alpha ,\limfunc{punct}%
}}\ , \\
\mathfrak{F}_{\alpha }^{\mathbf{b},\ast }\left( \mathcal{G},\mathcal{D}%
\right) &\lesssim &\mathfrak{E}_{2}^{\alpha ,\ast }+\sqrt{\mathcal{A}%
_{2}^{\alpha }}+\sqrt{\mathcal{A}_{2}^{\alpha ,\ast }}+\sqrt{A_{2}^{\alpha
,\ast ,\limfunc{punct}}}\ ,
\end{eqnarray*}%
with implied constants independent of the grids $\mathcal{D}$ and $\mathcal{G%
}$.
\end{proposition}

In order to prove this proposition, we first turn to recalling these more
refined notions of energy constants.

\subsection{Various energy conditions\label{sub various energy cond}}

In this subsection we recall various refinements of the strong energy
conditions appearing in the main theorem above. Variants of this material
already appear in earlier papers, but we repeat it here both for convenience
and in order to introduce some arguments we will use repeatedly later on.
These refinements represent the `weakest' energy side conditions that
suffice for use in our proof, but despite this, we will usually use the
large energy constant $\mathfrak{E}_{2}^{\alpha }$ in estimates to avoid
having to pay too much attention to which of the energy conditions we need
to use - leaving the determination of the weakest conditions in such
situations to the interested reader. We begin with the notion of `deeply
embedded'. Recall that the goodness parameter $\mathbf{r}\in \mathbb{N}$ is
determined by $\varepsilon >0$ in (\ref{choice of r}), and that $%
0<\varepsilon <\frac{1}{2}<\frac{1}{2-\alpha }$.

For arbitrary intervals in $\,J,K\in \mathcal{P}$, we say that $J$ is $%
\left( \mathbf{\rho },\varepsilon \right) $-\emph{deeply embedded} in $K$,
which we write as $J\Subset _{\mathbf{\rho },\varepsilon }K$, when $J\subset
K$ and both 
\begin{eqnarray}
\ell \left( J\right) &\leq &2^{-\mathbf{\rho }}\ell \left( K\right) ,
\label{def deep embed} \\
d\left( J,\partial K\right) &\geq &2\ell \left( J\right) ^{\varepsilon }\ell
\left( K\right) ^{1-\varepsilon }.  \notag
\end{eqnarray}%
Note that we use the \emph{boundary} of $K$ for the definition of $J\Subset
_{\mathbf{\rho },\varepsilon }K$, rather than the \emph{skeleton} or \emph{%
body} of $K$, which would result in a more restrictive notion of $\left( 
\mathbf{\rho },\varepsilon \right) $-deeply embedded. We will use this
notion for the purpose of grouping $\varepsilon -\limfunc{good}$ intervals
into the following collections. Fix grids $\mathcal{D}$ and $\mathcal{G}$.
For $K\in \mathcal{D}$, define the collections,%
\begin{eqnarray}
\mathcal{M}_{\left( \mathbf{\rho },\varepsilon \right) -\limfunc{deep},%
\mathcal{G}}\left( K\right) &\equiv &\left\{ J\in \mathcal{G}:J\text{ is
maximal w.r.t }J\Subset _{\mathbf{\rho },\varepsilon }K\right\} ,
\label{def M_r-deep} \\
\mathcal{M}_{\left( \mathbf{\rho },\varepsilon \right) -\limfunc{deep},%
\mathcal{D}}\left( K\right) &\equiv &\left\{ M\in \mathcal{D}:M\text{ is
maximal w.r.t }M\Subset _{\mathbf{\rho },\varepsilon }K\right\} ,  \notag \\
\mathcal{W}\left( K\right) &\equiv &\left\{ M\in \mathcal{D}:M\text{ is
maximal w.r.t }3M\subset K\right\}  \notag
\end{eqnarray}%
where the first two consist of \emph{maximal} $\left( \mathbf{\rho }%
,\varepsilon \right) $-deeply embedded dyadic $\mathcal{G}$-subintervals $J$%
, respectively $\mathcal{D}$-subintervals $M$, of a $\mathcal{D}$-interval $%
K $, and the third consists of the maximal $\mathcal{D}$-subintervals $M$
whose triples are contained in $K$.

Let $\gamma >1$. Then the following bounded overlap property holds where $%
\mathcal{M}_{\left( \mathbf{\rho },\varepsilon \right) -\limfunc{deep}%
}\left( K\right) $ can be taken to be either $\mathcal{M}_{\left( \mathbf{%
\rho },\varepsilon \right) -\limfunc{deep},\mathcal{G}}\left( K\right) $ or $%
\mathcal{M}_{\left( \mathbf{\rho },\varepsilon \right) -\limfunc{deep},%
\mathcal{D}}\left( K\right) $ or $\mathcal{W}\left( K\right) $ throughout.

\begin{lemma}
Let $0<\varepsilon \leq 1<\gamma \leq 1+4\cdot 2^{\mathbf{\rho }\left(
1-\varepsilon \right) }$. Then%
\begin{equation}
\sum_{J\in \mathcal{M}_{\left( \mathbf{\rho },\varepsilon \right) -\limfunc{%
deep}}\left( K\right) }\mathbf{1}_{\gamma J}\leq \beta \mathbf{1}_{\left[
\dbigcup\limits_{J\in \mathcal{M}_{\left( \mathbf{\rho },\varepsilon \right)
-\limfunc{deep}}\left( K\right) }\gamma J\right] }  \label{bounded overlap}
\end{equation}%
holds for some positive constant $\beta $ depending only on $\gamma ,\mathbf{%
\rho }$ and $\varepsilon $. In addition $\gamma J\subset K$ for all $J\in 
\mathcal{M}_{\left( \mathbf{\rho },\varepsilon \right) -\limfunc{deep}%
}\left( K\right) $, and consequently%
\begin{equation}
\sum_{J\in \mathcal{M}_{\left( \mathbf{\rho },\varepsilon \right) -\limfunc{%
deep}}\left( K\right) }\mathbf{1}_{\gamma J}\leq \beta \mathbf{1}_{K}\ .
\label{bounded overlap in K}
\end{equation}%
A similar result holds for $\mathcal{W}\left( K\right) $.
\end{lemma}

\begin{proof}
We suppose $0<\varepsilon <1$ and leave the simpler case $\varepsilon =1$
for the reader. To prove (\ref{bounded overlap}), we first note that there
are at most $2^{\mathbf{\rho }+1}$ intervals $J$ contained in $K$ for which $%
\ell \left( J\right) >2^{-\mathbf{\rho }}\ell \left( K\right) $. On the
other hand, the maximal $\left( \mathbf{\rho },\varepsilon \right) $-deeply
embedded subintervals $J$ of $K$ also satisfy the comparability condition%
\begin{equation*}
2\ell \left( J\right) ^{\varepsilon }\ell \left( K\right) ^{1-\varepsilon
}\leq d\left( J,\partial K\right) \leq d\left( \pi J,\partial K\right) -\ell
\left( J\right) \leq 2\left( 2\ell \left( J\right) \right) ^{\varepsilon
}\ell \left( K\right) ^{1-\varepsilon }-\ell \left( J\right) \leq 4\ell
\left( J\right) ^{\varepsilon }\ell \left( K\right) ^{1-\varepsilon }-\ell
\left( J\right) .
\end{equation*}%
Now with $0<\varepsilon <1$ and $\gamma >1$ fixed, let $y\in K$. Then if $%
y\in \gamma J$, we have%
\begin{eqnarray*}
2\ell \left( J\right) ^{\varepsilon }\ell \left( K\right) ^{1-\varepsilon }
&\leq &d\left( J,\partial K\right) \leq \gamma \ell \left( J\right) +d\left(
\gamma J,\partial K\right) \\
&\leq &\gamma \ell \left( J\right) +d\left( y,\partial K\right) .
\end{eqnarray*}%
Now assume that $\frac{\ell \left( J\right) }{\ell \left( K\right) }\leq
\left( \frac{1}{\gamma }\right) ^{\frac{1}{1-\varepsilon }}$. Then we have $%
\gamma \ell \left( J\right) \leq \ell \left( J\right) ^{\varepsilon }\ell
\left( K\right) ^{1-\varepsilon }$ and so 
\begin{equation*}
\ell \left( J\right) ^{\varepsilon }\ell \left( K\right) ^{1-\varepsilon
}\leq d\left( y,\partial K\right) .
\end{equation*}%
But we also have 
\begin{equation*}
d\left( y,\partial K\right) \leq \ell \left( J\right) +d\left( J,\partial
K\right) \leq \ell \left( J\right) +4\ell \left( J\right) ^{\varepsilon
}\ell \left( K\right) ^{1-\varepsilon }-\ell \left( J\right) \leq 4\ell
\left( J\right) ^{\varepsilon }\ell \left( K\right) ^{1-\varepsilon },
\end{equation*}%
and so altogether, under the assumption that $\frac{\ell \left( J\right) }{%
\ell \left( K\right) }\leq \left( \frac{1}{\gamma }\right) ^{\frac{1}{%
1-\varepsilon }}$, we have%
\begin{eqnarray*}
\frac{1}{4}d\left( y,\partial K\right) &\leq &\ell \left( J\right)
^{\varepsilon }\ell \left( K\right) ^{1-\varepsilon }\leq d\left( y,\partial
K\right) , \\
\text{i.e. }\left( \frac{1}{4}\frac{d\left( y,\partial K\right) }{\ell
\left( K\right) ^{1-\varepsilon }}\right) ^{\frac{1}{\varepsilon }} &\leq
&\ell \left( J\right) \leq \left( \frac{d\left( y,\partial K\right) }{\ell
\left( K\right) ^{1-\varepsilon }}\right) ^{\frac{1}{\varepsilon }},
\end{eqnarray*}%
which shows that the number of $J^{\prime }s$ satisfying $y\in \gamma J$ and 
$\frac{\ell \left( J\right) }{\ell \left( K\right) }\leq \left( \frac{1}{%
\gamma }\right) ^{\frac{1}{1-\varepsilon }}$ is at most $C^{\prime }\frac{1}{%
\varepsilon }$. On the other hand, the number of $J^{\prime }s$ contained in 
$K$ satisfying $y\in \gamma J$ and $\frac{\ell \left( J\right) }{\ell \left(
K\right) }>\left( \frac{1}{\gamma }\right) ^{\frac{1}{1-\varepsilon }}$ is
at most $C^{\prime }\frac{1}{1-\varepsilon }\left( 1+\log _{2}\gamma \right) 
$. This proves (\ref{bounded overlap}) with 
\begin{equation*}
\beta =2^{\mathbf{\rho }+1}+C^{\prime }\frac{1}{\varepsilon }+C^{\prime }%
\frac{1}{1-\varepsilon }\left( 1+\log _{2}\gamma \right) .
\end{equation*}

In order to prove (\ref{bounded overlap in K}) it suffices, by (\ref{bounded
overlap}), to prove $\gamma J\subset K$ for all $J\in \mathcal{M}_{\left( 
\mathbf{\rho },\varepsilon \right) -\limfunc{deep}}\left( K\right) $. But $%
J\in \mathcal{M}_{\left( \mathbf{\rho },\varepsilon \right) -\limfunc{deep}%
}\left( K\right) $ implies%
\begin{equation*}
2\ell \left( J\right) ^{\varepsilon }\ell \left( K\right) ^{1-\varepsilon
}\leq d\left( J,\partial K\right) =d\left( c_{J},\partial K\right) +\frac{1}{%
2}\ell \left( J\right) .
\end{equation*}%
We wish to show $\gamma J\subset K$, which is implied by 
\begin{equation*}
\gamma \frac{1}{2}\ell \left( J\right) \leq d\left( c_{J},K^{c}\right)
=d\left( J,\partial K\right) +\frac{1}{2}\ell \left( J\right) .
\end{equation*}%
But we have%
\begin{equation*}
d\left( J,\partial K\right) +\frac{1}{2}\ell \left( J\right) \geq 2\ell
\left( J\right) ^{\varepsilon }\ell \left( K\right) ^{1-\varepsilon }+\frac{1%
}{2}\ell \left( J\right) ,
\end{equation*}%
and so it suffices to show that%
\begin{equation*}
2\ell \left( J\right) ^{\varepsilon }\ell \left( K\right) ^{1-\varepsilon }+%
\frac{1}{2}\ell \left( J\right) \geq \gamma \frac{1}{2}\ell \left( J\right) ,
\end{equation*}%
which is equivalent to%
\begin{equation*}
\gamma -1\leq 4\ell \left( J\right) ^{\varepsilon -1}\ell \left( K\right)
^{1-\varepsilon }.
\end{equation*}%
But the smallest that $\ell \left( J\right) ^{\varepsilon -1}\ell \left(
K\right) ^{1-\varepsilon }$ can get for $J\in \mathcal{M}_{\left( \mathbf{%
\rho },\varepsilon \right) -\limfunc{deep}}\left( K\right) $ is $2^{\mathbf{%
\rho }\left( 1-\varepsilon \right) }\geq 1$, and so $\gamma \leq 1+4\cdot 2^{%
\mathbf{\rho }\left( 1-\varepsilon \right) }$ implies $\gamma -1\leq 4\ell
\left( J\right) ^{\varepsilon -1}\ell \left( K\right) ^{1-\varepsilon }$,
which completes the proof.

The reader can easily verify the same argument works for the Whitney
collection $\mathcal{W}\left( K\right) $.
\end{proof}

Now we recall the notion of \emph{alternate} dyadic intervals from \cite%
{SaShUr7}, which we rename \emph{augmented} dyadic intervals here.

\begin{definition}
\label{def dyadic}Given a dyadic grid $\mathcal{D}$, the \emph{augmented
dyadic grid} $\mathcal{AD}$ consists\ of those intervals $I$ whose dyadic
children $I^{\prime }$ belong to the grid $\mathcal{D}$.
\end{definition}

Of course an augmented grid is not actually a grid because the nesting
property fails, but this terminology should cause no confusion. These
augmented grids will be needed in order to use the `prepare to puncture'
argument (introduced in \cite{SaShUr9}) at several places below.

Now we proceed to recall certain of the definitions of various energy
conditions from \cite{SaShUr5} and \cite{SaShUr7}. While these definitions
are not explicitly used in the proof of functional energy, some of the
arguments we give to control them will be appealed to later, and so we take
the time to develop these definitions in detail.

\subsubsection{Whitney energy conditions}

The following definition of Whitney energy condition uses the \emph{Whitney}
decomposition $\mathcal{M}_{\left( \mathbf{\rho },1\right) -\limfunc{deep},%
\mathcal{D}}\left( I_{r}\right) $ into $\mathcal{D}$-dyadic intervals in
which $\varepsilon =1$, as well as the `large' pseudoprojections%
\begin{equation}
\mathsf{Q}_{K}^{\omega ,\mathbf{b}^{\ast }}\equiv \sum_{J\in \mathcal{G}:\
J\subset K}\bigtriangleup _{J}^{\omega ,\mathbf{b}^{\ast }}.
\label{large pseudo}
\end{equation}

\begin{definition}
\label{energy condition}Suppose $\sigma $ and $\omega $ are locally finite
positive Borel measures on $\mathbb{R}$ and fix $\gamma >1$. Then the\
Whitney energy condition constant $\mathcal{E}_{2}^{\alpha ,\func{Whitney}}$
is given by%
\begin{equation*}
\left( \mathcal{E}_{2}^{\alpha ,\func{Whitney}}\right) ^{2}\equiv \sup_{%
\mathcal{D},\mathcal{G}}\sup_{I=\dot{\cup}I_{r}}\frac{1}{\left\vert
I\right\vert _{\sigma }}\sum_{r=1}^{\infty }\sum_{M\in \mathcal{W}\left(
I_{r}\right) }\left( \frac{\mathrm{P}^{\alpha }\left( M,\mathbf{1}%
_{I\setminus \gamma M}\sigma \right) }{\left\vert M\right\vert }\right)
^{2}\left\Vert \mathsf{Q}_{M}^{\omega ,\mathbf{b}^{\ast }}x\right\Vert
_{L^{2}\left( \omega \right) }^{\spadesuit 2},
\end{equation*}%
where $\sup_{\mathcal{D},\mathcal{G}}\sup_{I=\dot{\cup}I_{r}}$ is taken over

\begin{enumerate}
\item all dyadic grids $\mathcal{D}$ and $\mathcal{G}$,

\item all $\mathcal{D}$-dyadic intervals $I$,

\item and all partitions $\left\{ I_{r}\right\} _{r=1}^{N\text{ or }\infty }$
of the interval $I$ into $\mathcal{D}$-dyadic subintervals $I_{r}$.
\end{enumerate}
\end{definition}

If the parameter $\gamma >1$ above is chosen sufficiently close to $1$, then
the collection of intervals $\left\{ \gamma M\right\} _{M\in \mathcal{W}%
\left( I_{r}\right) }$ has bounded overlap $\beta $ by (\ref{bounded overlap
in K}), and the Whitney energy constant $\mathcal{E}_{2}^{\alpha ,\func{%
Whitney}}$ is controlled by the strong energy constant $\mathcal{E}%
_{2}^{\alpha }$ in (\ref{strong b* energy}),%
\begin{equation}
\mathcal{E}_{2}^{\alpha ,\func{Whitney}}\lesssim \mathcal{E}_{2}^{\alpha }.
\label{en con}
\end{equation}%
Indeed, to see this, fix a decomposition of an interval 
\begin{equation}
I=\overset{\cdot }{\dbigcup }_{1\leq r<\infty }\overset{\cdot }{\dbigcup }%
_{M\in \mathcal{W}\left( I_{r}\right) }M  \label{decomp int}
\end{equation}%
as in Definition \ref{energy condition}. Then consider the \emph{sub}%
decomposition 
\begin{equation*}
I\supset \overset{\cdot }{\dbigcup }_{1\leq r<\infty }\overset{\cdot }{%
\dbigcup }_{M\in \mathcal{W}\left( I_{r}\right) }M
\end{equation*}%
of the interval $I$ given by the collection of intervals,%
\begin{equation*}
\mathcal{I}\equiv \overset{\cdot }{\dbigcup }_{1\leq r<\infty }\mathcal{W}%
\left( I_{r}\right) .
\end{equation*}%
We then have%
\begin{equation*}
\left( \mathcal{E}_{2}^{\alpha }\right) ^{2}\geq \frac{1}{\left\vert
I\right\vert _{\sigma }}\sum_{r=1}^{\infty }\sum_{M\in \mathcal{W}\left(
I_{r}\right) }\left( \frac{\mathrm{P}^{\alpha }\left( M,\mathbf{1}_{I}\sigma
\right) }{\left\vert M\right\vert }\right) ^{2}\left\Vert x-m_{M}^{\omega
}\right\Vert _{L^{2}\left( \mathbf{1}_{M}\omega \right) }^{2}\ .
\end{equation*}%
Now $\mathrm{P}^{\alpha }\left( M,\mathbf{1}_{I}\sigma \right) \geq \mathrm{P%
}^{\alpha }\left( M,\mathbf{1}_{I\setminus \gamma M}\sigma \right) $ and
from (\ref{note more}),%
\begin{equation*}
\left\Vert x-m_{M}^{\omega }\right\Vert _{L^{2}\left( \mathbf{1}_{M}\omega
\right) }^{2}\gtrsim \left\Vert \mathsf{Q}_{M}^{\omega ,\mathbf{b}^{\ast
}}x\right\Vert _{L^{2}\left( \omega \right) }^{\spadesuit 2},
\end{equation*}
and combining these two inequalities, we obtain that%
\begin{equation*}
\left( \mathcal{E}_{2}^{\alpha }\right) ^{2}\geq c\frac{1}{\left\vert
I\right\vert _{\sigma }}\sum_{r=1}^{\infty }\sum_{M\in \mathcal{W}\left(
I_{r}\right) }\left( \frac{\mathrm{P}^{\alpha }\left( M,\mathbf{1}%
_{I\setminus \gamma M}\sigma \right) }{\left\vert M\right\vert }\right)
^{2}\left\Vert \mathsf{Q}_{M}^{\omega ,\mathbf{b}^{\ast }}x\right\Vert
_{L^{2}\left( \omega \right) }^{\spadesuit 2}\ .
\end{equation*}%
Thus we conclude that%
\begin{equation*}
\frac{1}{\left\vert I\right\vert _{\sigma }}\sum_{r=1}^{\infty }\sum_{M\in 
\mathcal{W}\left( I_{r}\right) }\left( \frac{\mathrm{P}^{\alpha }\left( M,%
\mathbf{1}_{I\setminus \gamma M}\sigma \right) }{\left\vert M\right\vert }%
\right) ^{2}\left\Vert \mathsf{Q}_{M}^{\omega ,\mathbf{b}^{\ast
}}x\right\Vert _{L^{2}\left( \omega \right) }^{\spadesuit 2}\leq \frac{C}{c}%
\beta \left( \mathcal{E}_{2}^{\alpha }\right) ^{2},
\end{equation*}%
and taking the supremum over all decompositions (\ref{decomp int}) as in
Definition \ref{energy condition}, we obtain (\ref{en con}).

There is a similar definition for the dual (backward) Whitney energy
conditions that simply interchanges $\sigma $ and $\omega $ everywhere.
These definitions of\ the Whitney energy conditions depend on the choice of $%
\gamma >1$.

\begin{description}
\item[Commentary on proofs] We now introduce a number of results concerning
partial plugging of the hole for Whitney energy conditions.
\end{description}

Note that we can `partially' plug the $\gamma $-hole in the Poisson integral 
$\mathrm{P}^{\alpha }\left( J,\mathbf{1}_{I\setminus \gamma J}\sigma \right) 
$ for $\mathcal{E}_{2}^{\alpha ,\func{Whitney}}$ using the offset $%
A_{2}^{\alpha }$ condition and the bounded overlap property (\ref{bounded
overlap in K}). Indeed, define 
\begin{eqnarray}
&&  \label{plug} \\
&&\left( \mathcal{E}_{2}^{\alpha ,\func{Whitney}\limfunc{partial}}\right)
^{2}\equiv \sup_{\mathcal{D},\mathcal{G}}\sup_{I=\dot{\cup}I_{r}}\frac{1}{%
\left\vert I\right\vert _{\sigma }}\sum_{r=1}^{\infty }\sum_{M\in \mathcal{W}%
\left( I_{r}\right) }\left( \frac{\mathrm{P}^{\alpha }\left( M,\mathbf{1}%
_{I\setminus M}\sigma \right) }{\left\vert M\right\vert }\right)
^{2}\left\Vert \mathsf{Q}_{M}^{\omega ,\mathbf{b}^{\ast }}x\right\Vert
_{L^{2}\left( \omega \right) }^{\spadesuit 2}\ .  \notag
\end{eqnarray}%
Recall from (\ref{bounded overlap in K}) that%
\begin{equation*}
\gamma M\subset I_{r}\text{ for all }M\in \mathcal{W}\left( I_{r}\right) 
\text{ provided }\gamma \leq 5.
\end{equation*}%
At this point we need the following analogues of the `energy $A_{2}^{\alpha
} $ conditions' from \cite{SaShUr9}, which we denote by $A_{2}^{\alpha ,%
\limfunc{energy}}$ and $A_{2}^{\alpha ,\ast ,\limfunc{energy}}$, and define
by%
\begin{eqnarray}
A_{2}^{\alpha ,\limfunc{energy}}\left( \sigma ,\omega \right) &\equiv
&\sup_{Q\in \mathcal{P}}\frac{\left\Vert \mathsf{Q}_{Q}^{\omega ,\mathbf{b}%
^{\ast }}\frac{x}{\ell \left( Q\right) }\right\Vert _{L^{2}\left( \omega
\right) }^{\spadesuit 2}}{\left\vert Q\right\vert ^{1-\alpha }}\frac{%
\left\vert Q\right\vert _{\sigma }}{\left\vert Q\right\vert ^{1-\alpha }},
\label{def energy A2} \\
A_{2}^{\alpha ,\ast ,\limfunc{energy}}\left( \sigma ,\omega \right) &\equiv
&\sup_{Q\in \mathcal{P}}\frac{\left\vert Q\right\vert _{\omega }}{\left\vert
Q\right\vert ^{1-\alpha }}\frac{\left\Vert \mathsf{Q}_{Q}^{\sigma ,\mathbf{b}%
}\frac{x}{\ell \left( Q\right) }\right\Vert _{L^{2}\left( \sigma \right)
}^{\spadesuit 2}}{\left\vert Q\right\vert ^{1-\alpha }}.  \notag
\end{eqnarray}%
Then if $\gamma \leq 5$, we have%
\begin{eqnarray}
&&\left( \mathcal{E}_{2}^{\alpha ,\func{Whitney}\limfunc{partial}}\right)
^{2}  \label{plug the hole deep} \\
&\lesssim &\sup_{\mathcal{D},\mathcal{G}}\sup_{I=\dot{\cup}I_{r}}\frac{1}{%
\left\vert I\right\vert _{\sigma }}\sum_{r=1}^{\infty }\sum_{M\in \mathcal{W}%
\left( I_{r}\right) }\left( \frac{\mathrm{P}^{\alpha }\left( M,\mathbf{1}%
_{I\setminus \gamma M}\sigma \right) }{\left\vert M\right\vert }\right)
^{2}\left\Vert \mathsf{Q}_{M}^{\omega ,\mathbf{b}^{\ast }}x\right\Vert
_{L^{2}\left( \omega \right) }^{\spadesuit 2}  \notag \\
&&+\sup_{\mathcal{D},\mathcal{G}}\sup_{I=\dot{\cup}I_{r}}\frac{1}{\left\vert
I\right\vert _{\sigma }}\sum_{r=1}^{\infty }\sum_{M\in \mathcal{W}\left(
I_{r}\right) }\left( \frac{\mathrm{P}^{\alpha }\left( M,\mathbf{1}_{\gamma
M\setminus M}\sigma \right) }{\left\vert M\right\vert }\right)
^{2}\left\Vert \mathsf{Q}_{M}^{\omega ,\mathbf{b}^{\ast }}x\right\Vert
_{L^{2}\left( \omega \right) }^{\spadesuit 2}  \notag \\
&\lesssim &\left( \mathcal{E}_{2}^{\alpha ,\func{Whitney}}\right) ^{2}+\sup_{%
\mathcal{D},\mathcal{G}}\sup_{I=\dot{\cup}I_{r}}\frac{1}{\left\vert
I\right\vert _{\sigma }}\sum_{r=1}^{\infty }\sum_{M\in \mathcal{W}\left(
I_{r}\right) }A_{2}^{\alpha ,\func{energy}}\left\vert \gamma M\right\vert
_{\sigma }\lesssim \left( \mathcal{E}_{2}^{\alpha ,\limfunc{deep}}\right)
^{2}+\beta A_{2}^{\alpha ,\func{energy}}\ ,  \notag
\end{eqnarray}%
by (\ref{bounded overlap in K}).

\subsubsection{Plugged energy conditions}

We continue to recall some results from \cite{SaShUr9} and \cite{SaShUr10}
that we will use repeatedly here. For example, we will use the punctured
Muckenhoupt conditions $A_{2}^{\alpha ,\limfunc{punct}}$ and $A_{2}^{\alpha
,\ast ,\limfunc{punct}}$ introduced earlier in (\ref{puncture}) to control
the \emph{plugged }energy conditions, where the hole in the argument of the
Poisson term $\mathrm{P}^{\alpha }\left( M,\mathbf{1}_{I\setminus M}\sigma
\right) $ in the partially plugged energy condition above, is replaced with
the `plugged' term $\mathrm{P}^{\alpha }\left( M,\mathbf{1}_{I}\sigma
\right) $, for example%
\begin{equation}
\left( \mathcal{E}_{2}^{\alpha ,\func{Whitney}\limfunc{plug}}\right)
^{2}\equiv \sup_{\mathcal{D},\mathcal{G}}\sup_{I=\dot{\cup}I_{r}}\frac{1}{%
\left\vert I\right\vert _{\sigma }}\sum_{r=1}^{\infty }\sum_{M\in \mathcal{W}%
\left( I_{r}\right) }\left( \frac{\mathrm{P}^{\alpha }\left( M,\mathbf{1}%
_{I}\sigma \right) }{\left\vert M\right\vert }\right) ^{2}\left\Vert \mathsf{%
Q}_{M}^{\omega ,\mathbf{b}^{\ast }}x\right\Vert _{L^{2}\left( \omega \right)
}^{\spadesuit 2}\ .  \label{def deep plug}
\end{equation}%
By an argument similar to that in (\ref{plug the hole deep}), we obtain%
\begin{equation}
\mathcal{E}_{2}^{\alpha ,\func{Whitney}\limfunc{plug}}\lesssim \mathcal{E}%
_{2}^{\alpha ,\func{Whitney}\limfunc{partial}}+A_{2}^{\alpha ,\func{energy}}.
\label{plug the hole deep'}
\end{equation}

We first show that the punctured Muckenhoupt conditions $A_{2}^{\alpha ,%
\limfunc{punct}}$ and $A_{2}^{\alpha ,\ast ,\limfunc{punct}}$ control
respectively the `energy $A_{2}^{\alpha }$ conditions' in (\ref{def energy
A2}). We will make reference to the proof of the next lemma (for the $T1$
theorem this is from \cite[Lemma 3.2 on page 328.]{SaShUr9}) several times
in the sequel. We repeat the proof from \cite[Lemma 3.2 on page 328.]%
{SaShUr9} but with modifications to accommodate the differences that arise
here in the setting of a local $Tb$ theorem. Recall that $\mathfrak{P}%
_{\left( \sigma ,\omega \right) }$ is defined in (\ref{def common point mass}%
) above.

\begin{lemma}
\label{energy A2}For any positive locally finite Borel measures $\sigma
,\omega $ we have%
\begin{eqnarray*}
A_{2}^{\alpha ,\limfunc{energy}}\left( \sigma ,\omega \right) &\lesssim
&A_{2}^{\alpha ,\limfunc{punct}}\left( \sigma ,\omega \right) , \\
A_{2}^{\alpha ,\ast ,\limfunc{energy}}\left( \sigma ,\omega \right)
&\lesssim &A_{2}^{\alpha ,\ast ,\limfunc{punct}}\left( \sigma ,\omega
\right) .
\end{eqnarray*}
\end{lemma}

\begin{proof}
Fix an interval $Q\in \mathcal{D}$. Recall the definition of $\omega \left(
Q,\mathfrak{P}_{\left( \sigma ,\omega \right) }\right) $ in (\ref{puncture}%
). If $\omega \left( Q,\mathfrak{P}_{\left( \sigma ,\omega \right) }\right)
\geq \frac{1}{2}\left\vert Q\right\vert _{\omega }$, then we trivially have%
\begin{eqnarray*}
\frac{\left\Vert \mathsf{Q}_{Q}^{\omega ,\mathbf{b}^{\ast }}\frac{x}{\ell
\left( Q\right) }\right\Vert _{L^{2}\left( \omega \right) }^{\spadesuit 2}}{%
\left\vert Q\right\vert ^{1-\alpha }}\frac{\left\vert Q\right\vert _{\sigma }%
}{\left\vert Q\right\vert ^{1-\alpha }} &\lesssim &\frac{\left\vert
Q\right\vert _{\omega }}{\left\vert Q\right\vert ^{1-\alpha }}\frac{%
\left\vert Q\right\vert _{\sigma }}{\left\vert Q\right\vert ^{1-\alpha }} \\
&\leq &2\frac{\omega \left( Q,\mathfrak{P}_{\left( \sigma ,\omega \right)
}\right) }{\left\vert Q\right\vert ^{1-\alpha }}\frac{\left\vert
Q\right\vert _{\sigma }}{\left\vert Q\right\vert ^{1-\alpha }}\leq
2A_{2}^{\alpha ,\limfunc{punct}}\left( \sigma ,\omega \right) .
\end{eqnarray*}%
On the other hand, if $\omega \left( Q,\mathfrak{P}_{\left( \sigma ,\omega
\right) }\right) <\frac{1}{2}\left\vert Q\right\vert _{\omega }$ then there
is a point $p\in Q\cap \mathfrak{P}_{\left( \sigma ,\omega \right) }$ such
that%
\begin{equation*}
\omega \left( \left\{ p\right\} \right) >\frac{1}{2}\left\vert Q\right\vert
_{\omega }\ ,
\end{equation*}%
and consequently, $p$ is the largest $\omega $-point mass in $Q$. Thus if we
define $\widetilde{\omega }=\omega -\omega \left( \left\{ p\right\} \right)
\delta _{p}$, then we have%
\begin{equation*}
\omega \left( Q,\mathfrak{P}_{\left( \sigma ,\omega \right) }\right)
=\left\vert Q\right\vert _{\widetilde{\omega }}\ .
\end{equation*}%
Now we observe from the construction of martingale differences that%
\begin{equation*}
\bigtriangleup _{J}^{\widetilde{\omega },\mathbf{b}^{\ast }}=\bigtriangleup
_{J}^{\omega ,\mathbf{b}^{\ast }},\ \ \ \ \ \text{for all }J\in \mathcal{D}%
\text{ with }p\notin J.
\end{equation*}%
So for each $s\geq 0$ there is a unique interval $J_{s}\in \mathcal{D}$ with 
$\ell \left( J_{s}\right) =2^{-s}\ell \left( Q\right) $ that contains the
point $p$. Now observe that, just as for the Haar projection, the
one-dimensional projection $\bigtriangleup _{J_{s}}^{\omega ,\mathbf{b}%
^{\ast }}$ is given by $\bigtriangleup _{J_{s}}^{\omega ,\mathbf{b}^{\ast
}}f=\left\langle h_{J_{s}}^{\omega ,\mathbf{b}^{\ast }},f\right\rangle
_{\omega }h_{J_{s}}^{\omega ,\mathbf{b}^{\ast }}$ for a unique up to $\pm $
unit vector $h_{J_{s}}^{\omega ,\mathbf{b}^{\ast }}$. For this interval we
then have%
\begin{eqnarray*}
\left\Vert \bigtriangleup _{J_{s}}^{\omega ,\mathbf{b}^{\ast }}x\right\Vert
_{L^{2}\left( \omega \right) }^{2} &=&\left\vert \left\langle
h_{J_{s}}^{\omega ,\mathbf{b}^{\ast }},x\right\rangle _{\omega }\right\vert
^{2}=\left\vert \left\langle h_{J_{s}}^{\omega ,\mathbf{b}^{\ast
}},x-p\right\rangle _{\omega }\right\vert ^{2} \\
&=&\left\vert \int_{J_{s}}h_{J_{s}}^{\omega ,\mathbf{b}^{\ast }}\left(
x\right) \left( x-p\right) d\omega \left( x\right) \right\vert
^{2}=\left\vert \int_{J_{s}}h_{J_{s}}^{\omega ,\mathbf{b}^{\ast }}\left(
x\right) \left( x-p\right) d\widetilde{\omega }\left( x\right) \right\vert
^{2} \\
&\leq &\left\Vert h_{J_{s}}^{\omega ,\mathbf{b}^{\ast }}\right\Vert
_{L^{2}\left( \widetilde{\omega }\right) }^{2}\left\Vert \mathbf{1}%
_{J_{s}}\left( x-p\right) \right\Vert _{L^{2}\left( \widetilde{\omega }%
\right) }^{2}\leq \left\Vert h_{J_{s}}^{\omega ,\mathbf{b}^{\ast
}}\right\Vert _{L^{2}\left( \omega \right) }^{2}\left\Vert \mathbf{1}%
_{J_{s}}\left( x-p\right) \right\Vert _{L^{2}\left( \widetilde{\omega }%
\right) }^{2} \\
&\leq &\ell \left( J_{s}\right) ^{2}\left\vert J_{s}\right\vert _{\widetilde{%
\omega }}\leq 2^{-2s}\ell \left( Q\right) ^{2}\left\vert Q\right\vert _{%
\widetilde{\omega }},
\end{eqnarray*}%
as well as%
\begin{equation*}
\inf_{z\in \mathbb{R}}\left\Vert \widehat{\nabla }_{J_{s}}^{\omega }\left(
x-z\right) \right\Vert _{L^{2}\left( \omega \right) }^{2}\lesssim \left\Vert
\left( x-p\right) \right\Vert _{L^{2}\left( \mathbf{1}_{J_{s}}\omega \right)
}^{2}=\left\Vert \left( x-p\right) \right\Vert _{L^{2}\left( \mathbf{1}%
_{J_{s}}\widetilde{\omega }\right) }^{2}\leq \ell \left( J_{s}\right)
^{2}\left\vert J_{s}\right\vert _{\widetilde{\omega }}\leq 2^{-2s}\ell
\left( Q\right) ^{2}\left\vert Q\right\vert _{\widetilde{\omega }}\ ,
\end{equation*}%
from (\ref{note more}). Thus we can estimate%
\begin{eqnarray}
&&  \label{omega tilda} \\
\left\Vert \mathsf{Q}_{Q}^{\omega ,\mathbf{b}^{\ast }}\frac{x}{\ell \left(
Q\right) }\right\Vert _{L^{2}\left( \omega \right) }^{\spadesuit 2} &\leq &%
\frac{1}{\ell \left( Q\right) ^{2}}\left( \sum_{J\in \mathcal{D}:\ J\subset
Q}\left\Vert \bigtriangleup _{J}^{\omega ,\mathbf{b}^{\ast }}x\right\Vert
_{L^{2}\left( \omega \right) }^{2}+\inf_{z\in \mathbb{R}}\left\Vert \widehat{%
\nabla }_{J_{s}}^{\omega }\left( x-z\right) \right\Vert _{L^{2}\left( \omega
\right) }^{2}\right)  \notag \\
&=&\frac{1}{\ell \left( Q\right) ^{2}}\left( \sum_{J\in \mathcal{D}:\
p\notin J\subset Q}\left\Vert \bigtriangleup _{J}^{\widetilde{\omega },%
\mathbf{b}^{\ast }}x\right\Vert _{L^{2}\left( \widetilde{\omega }\right)
}^{2}+\sum_{s=0}^{\infty }\left\Vert \bigtriangleup _{J_{s}}^{\omega ,%
\mathbf{b}^{\ast }}x\right\Vert _{L^{2}\left( \omega \right)
}^{2}+\inf_{z\in \mathbb{R}}\left\Vert \widehat{\nabla }_{J_{s}}^{\omega
}\left( x-z\right) \right\Vert _{L^{2}\left( \omega \right) }^{2}\right) 
\notag \\
&\lesssim &\frac{1}{\ell \left( Q\right) ^{2}}\left( \left\Vert \mathsf{Q}%
_{Q}^{\widetilde{\omega },\mathbf{b}^{\ast }}x\right\Vert _{L^{2}\left( 
\widetilde{\omega }\right) }^{\spadesuit 2}+\sum_{s=0}^{\infty }2^{-2s}\ell
\left( Q\right) ^{2}\left\vert Q\right\vert _{\widetilde{\omega }}\right) 
\notag \\
&\lesssim &\frac{1}{\ell \left( Q\right) ^{2}}\left( \ell \left( Q\right)
^{2}\left\vert Q\right\vert _{\widetilde{\omega }}+\sum_{s=0}^{\infty
}2^{-2s}\ell \left( Q\right) ^{2}\left\vert Q\right\vert _{\widetilde{\omega 
}}\right)  \notag \\
&\leq &3\left\vert Q\right\vert _{\widetilde{\omega }}=3\omega \left( Q,%
\mathfrak{P}_{\left( \sigma ,\omega \right) }\right) ,  \notag
\end{eqnarray}%
and so 
\begin{equation*}
\frac{\left\Vert \mathsf{Q}_{Q}^{\omega ,\mathbf{b}^{\ast }}\frac{x}{\ell
\left( Q\right) }\right\Vert _{L^{2}\left( \omega \right) }^{\spadesuit 2}}{%
\left\vert Q\right\vert ^{1-\alpha }}\frac{\left\vert Q\right\vert _{\sigma }%
}{\left\vert Q\right\vert ^{1-\alpha }}\lesssim \frac{3\omega \left( Q,%
\mathfrak{P}_{\left( \sigma ,\omega \right) }\right) }{\left\vert
Q\right\vert ^{1-\alpha }}\frac{\left\vert Q\right\vert _{\sigma }}{%
\left\vert Q\right\vert ^{1-\alpha }}\leq 3A_{2}^{\alpha ,\limfunc{punct}%
}\left( \sigma ,\omega \right) .
\end{equation*}%
Now take the supremum over $Q\in \mathcal{D}$ to obtain $A_{2}^{\alpha ,%
\limfunc{energy}}\left( \sigma ,\omega \right) \lesssim A_{2}^{\alpha ,%
\limfunc{punct}}\left( \sigma ,\omega \right) $. The dual inequality follows
upon interchanging the measures $\sigma $ and $\omega $.
\end{proof}

We isolate a simple but key fact that will be used repeatedly in what
follows:%
\begin{equation}
\sum_{Q\in \mathcal{D}:\ Q\subset P}\ell \left( Q\right) ^{2}\left\vert
Q\right\vert _{\mu }\lesssim \ell \left( P\right) ^{2}\left\vert
P\right\vert _{\mu }\ ,\ \ \ \ \ \text{for }P\in \mathcal{D}\text{ and }\mu 
\text{ a positive measure}.  \label{key fact}
\end{equation}%
Indeed, to see (\ref{key fact}), simply pigeonhole the length of $Q$
relative to that of $P$ and sum. The next corollary follows immediately from
Lemma \ref{energy A2}, (\ref{plug the hole deep}) and (\ref{plug the hole
deep'}).

\begin{corollary}
\label{all plugged}Provided $1<\gamma \leq 5$,%
\begin{equation*}
\mathcal{E}_{2}^{\alpha ,\func{Whitney}\limfunc{plug}}\lesssim \mathcal{E}%
_{2}^{\alpha ,\func{Whitney}\limfunc{partial}}+A_{2}^{\alpha ,\limfunc{punct}%
}\lesssim \mathcal{E}_{2}^{\alpha ,\func{Whitney}}+A_{2}^{\alpha ,\limfunc{%
punct}}\ ,
\end{equation*}%
and similarly for the dual plugged energy condition.
\end{corollary}

\subsubsection{Plugged $\mathcal{A}_{2}^{\protect\alpha ,\limfunc{energy}%
\limfunc{plug}}$ conditions}

Using Lemma \ref{energy A2} we can control the `plugged' energy $\mathcal{A}%
_{2}^{\alpha }$ conditions:%
\begin{eqnarray*}
\mathcal{A}_{2}^{\alpha ,\limfunc{energy}\limfunc{plug}}\left( \sigma
,\omega \right) &\equiv &\sup_{Q\in \mathcal{P}}\frac{\left\Vert \mathsf{Q}%
_{Q}^{\omega ,\mathbf{b}^{\ast }}\frac{x}{\ell \left( Q\right) }\right\Vert
_{L^{2}\left( \omega \right) }^{\spadesuit 2}}{\left\vert Q\right\vert
^{1-\alpha }}\mathcal{P}^{\alpha }\left( Q,\sigma \right) , \\
\mathcal{A}_{2}^{\alpha ,\ast ,\limfunc{energy}\limfunc{plug}}\left( \sigma
,\omega \right) &\equiv &\sup_{Q\in \mathcal{P}}\mathcal{P}^{\alpha }\left(
Q,\omega \right) \frac{\left\Vert \mathsf{Q}_{Q}^{\sigma ,\mathbf{b}}\frac{x%
}{\ell \left( Q\right) }\right\Vert _{L^{2}\left( \sigma \right)
}^{\spadesuit 2}}{\left\vert Q\right\vert ^{1-\alpha }}.
\end{eqnarray*}

\begin{lemma}
\label{energy A2 plugged}We have
\end{lemma}

\begin{eqnarray*}
\mathcal{A}_{2}^{\alpha ,\limfunc{energy}\limfunc{plug}}\left( \sigma
,\omega \right) &\mathcal{\lesssim }&\mathcal{A}_{2}^{\alpha }\left( \sigma
,\omega \right) +A_{2}^{\alpha ,\limfunc{energy}}\left( \sigma ,\omega
\right) , \\
\mathcal{A}_{2}^{\alpha ,\ast ,\limfunc{energy}\limfunc{plug}}\left( \sigma
,\omega \right) &\mathcal{\lesssim }&\mathcal{A}_{2}^{\alpha ,\ast }\left(
\sigma ,\omega \right) +A_{2}^{\alpha ,\ast ,\limfunc{energy}}\left( \sigma
,\omega \right) .
\end{eqnarray*}

\begin{proof}
We have%
\begin{eqnarray*}
\frac{\left\Vert \mathsf{Q}_{Q}^{\omega ,\mathbf{b}^{\ast }}\frac{x}{\ell
\left( Q\right) }\right\Vert _{L^{2}\left( \omega \right) }^{\spadesuit 2}}{%
\left\vert Q\right\vert ^{1-\alpha }}\mathcal{P}^{\alpha }\left( Q,\sigma
\right) &=&\frac{\left\Vert \mathsf{Q}_{Q}^{\omega ,\mathbf{b}^{\ast }}\frac{%
x}{\ell \left( Q\right) }\right\Vert _{L^{2}\left( \omega \right)
}^{\spadesuit 2}}{\left\vert Q\right\vert ^{1-\alpha }}\mathcal{P}^{\alpha
}\left( Q,\mathbf{1}_{Q^{c}}\sigma \right) +\frac{\left\Vert \mathsf{Q}%
_{Q}^{\omega ,\mathbf{b}^{\ast }}\frac{x}{\ell \left( Q\right) }\right\Vert
_{L^{2}\left( \omega \right) }^{\spadesuit 2}}{\left\vert Q\right\vert
^{1-\alpha }}\mathcal{P}^{\alpha }\left( Q,\mathbf{1}_{Q}\sigma \right) \\
&\lesssim &\frac{\left\vert Q\right\vert _{\omega }}{\left\vert Q\right\vert
^{1-\alpha }}\mathcal{P}^{\alpha }\left( Q,\mathbf{1}_{Q^{c}}\sigma \right) +%
\frac{\left\Vert \mathsf{Q}_{Q}^{\omega ,\mathbf{b}^{\ast }}\frac{x}{\ell
\left( Q\right) }\right\Vert _{L^{2}\left( \omega \right) }^{\spadesuit 2}}{%
\left\vert Q\right\vert ^{1-\alpha }}\frac{\left\vert Q\right\vert _{\sigma }%
}{\left\vert Q\right\vert ^{1-\alpha }} \\
&\lesssim &\mathcal{A}_{2}^{\alpha }\left( \sigma ,\omega \right)
+A_{2}^{\alpha ,\limfunc{energy}}\left( \sigma ,\omega \right) .
\end{eqnarray*}
\end{proof}

\subsection{The Poisson formulation}

Recall from Definitions \ref{def sharp cross} and \ref{shifted corona} that%
\begin{equation*}
\mathcal{C}_{F}^{\mathcal{G},\limfunc{shift}}=\left\{ J\in \mathcal{G}%
:J^{\maltese }\in \mathcal{C}_{F}\right\} ,
\end{equation*}%
where $F\in \mathcal{F}$ is a stopping interval in the dyadic grid $\mathcal{%
D}$. For convenience we repeat here the main result of this section,
Proposition \ref{func ener control}.

\begin{proposition}
\label{func ener control'}For all grids $\mathcal{D}$ and $\mathcal{G}$, and 
$\varepsilon >0$ sufficiently small, we have%
\begin{eqnarray*}
\mathfrak{F}_{\alpha }^{\mathbf{b}^{\ast }}\left( \mathcal{D},\mathcal{G}%
\right) &\lesssim &\mathfrak{E}_{2}^{\alpha }+\sqrt{\mathcal{A}_{2}^{\alpha }%
}+\sqrt{\mathcal{A}_{2}^{\alpha ,\ast }}+\sqrt{A_{2}^{\alpha ,\limfunc{punct}%
}}\ , \\
\mathfrak{F}_{\alpha }^{\mathbf{b},\ast }\left( \mathcal{G},\mathcal{D}%
\right) &\lesssim &\mathfrak{E}_{2}^{\alpha ,\ast }+\sqrt{\mathcal{A}%
_{2}^{\alpha }}+\sqrt{\mathcal{A}_{2}^{\alpha ,\ast }}+\sqrt{A_{2}^{\alpha
,\ast ,\limfunc{punct}}}\ ,
\end{eqnarray*}%
with implied constants independent of the grids $\mathcal{D}$ and $\mathcal{G%
}$.
\end{proposition}

To prove Proposition \ref{func ener control'}, we fix grids $\mathcal{D}$
and $\mathcal{G}$ and a subgrid $\mathcal{F}$ of $\mathcal{D}$\ as in (\ref%
{e.funcEnergy n}), and set 
\begin{equation}
\mu \equiv \sum_{F\in \mathcal{F}}\sum_{M\in \mathcal{W}\left( F\right)
}\left\Vert \mathsf{Q}_{F,M}^{\omega ,\mathbf{b}^{\ast }}x\right\Vert
_{L^{2}\left( \omega \right) }^{\spadesuit 2}\cdot \delta _{\left(
c_{M},\ell \left( M\right) \right) }\text{ and }d\overline{\mu }\left(
x,t\right) \equiv \frac{1}{t^{2}}d\mu \left( x,t\right) \ ,  \label{def mu n}
\end{equation}%
where $\mathcal{W}\left( F\right) $ consists of the maximal $\mathcal{D}$%
-subintervals of $F$ whose triples are contained in $F$, and where $\delta
_{\left( c_{M},\ell \left( M\right) \right) }$ denotes the Dirac unit mass
at the point $\left( c_{M},\ell \left( M\right) \right) $ in the upper
half-space $\mathbb{R}_{+}^{2}$. Here $M\in \mathcal{D}$ is a dyadic
interval with center $c_{M}$ and side length $\ell \left( M\right) $, and
for any interval $K\in \mathcal{P}$, the shorthand notation $\mathsf{P}%
_{F,K}^{\omega ,\mathbf{b}^{\ast }}$ (resp. $\mathsf{Q}_{F,K}^{\omega ,%
\mathbf{b}^{\ast }}$) is used for the localized pseudoprojection $\mathsf{P}%
_{\mathcal{C}_{F}^{\mathcal{G},\limfunc{shift}};K}^{\omega ,\mathbf{b}^{\ast
}}$ (resp. $\mathsf{Q}_{\mathcal{C}_{F}^{\mathcal{G},\limfunc{shift}%
};K}^{\omega ,\mathbf{b}^{\ast }}$) given in (\ref{def localization}):%
\begin{equation}
\mathsf{P}_{F,K}^{\omega ,\mathbf{b}^{\ast }}\equiv \mathsf{P}_{\mathcal{C}%
_{F}^{\mathcal{G},\limfunc{shift}};K}^{\omega ,\mathbf{b}^{\ast
}}=\sum_{J\subset K:\ J\in \mathcal{C}_{F}^{\mathcal{G},\limfunc{shift}%
}}\square _{J}^{\omega ,\mathbf{b}^{\ast }}\text{ }\left( \text{resp. }%
\mathsf{Q}_{F,K}^{\omega ,\mathbf{b}^{\ast }}\equiv \mathsf{Q}_{\mathcal{C}%
_{F}^{\mathcal{G},\limfunc{shift}};K}^{\omega ,\mathbf{b}^{\ast
}}=\sum_{J\subset K:\ J\in \mathcal{C}_{F}^{\mathcal{G},\limfunc{shift}%
}}\bigtriangleup _{J}^{\omega ,\mathbf{b}^{\ast }}\right) .  \label{def F,K}
\end{equation}%
We emphasize that all the subintervals $J$ that arise in the projection $%
\mathsf{Q}_{F,M}^{\omega ,\mathbf{b}^{\ast }}$ are good inside the intervals 
$F$ and beyond since $J^{\maltese }\subset F$. Here $J^{\maltese }$ is
defined in Definition \ref{def sharp cross} using the body of an interval.
Thus every $J\in \mathsf{Q}_{F}^{\omega ,\mathbf{b}^{\ast }}$ is contained
in a unique $M\in \mathcal{W}\left( F\right) $, so that $\mathsf{Q}%
_{F}^{\omega ,\mathbf{b}^{\ast }}=\overset{\cdot }{\dbigcup }_{M\in \mathcal{%
W}\left( F\right) }\mathsf{Q}_{F,M}^{\omega ,\mathbf{b}^{\ast }}$. We can
replace $x$ by $x-c$ inside the projection for any choice of $c$ we wish;
the projection is unchanged. More generally, $\delta _{q}$ denotes a Dirac
unit mass at a point $q$ in the upper half-space $\mathbb{R}_{+}^{2}$.

We will prove the two-weight inequality 
\begin{equation}
\left\Vert \mathbb{P}^{\alpha }\left( f\sigma \right) \right\Vert _{L^{2}(%
\mathbb{R}_{+}^{2},\overline{\mu })}\lesssim \left( \mathfrak{E}_{2}^{\alpha
}+\sqrt{\mathcal{A}_{2}^{\alpha }}+\sqrt{\mathcal{A}_{2}^{\alpha ,\ast }}+%
\sqrt{A_{2}^{\alpha ,\limfunc{punct}}}\right) \lVert f\rVert _{L^{2}\left(
\sigma \right) }\,,  \label{two weight Poisson n}
\end{equation}%
for all nonnegative $f$ in $L^{2}\left( \sigma \right) $, noting that $%
\mathcal{F}$ and $f$ are \emph{not} related here. Above, $\mathbb{P}^{\alpha
}(\cdot )$ denotes the $\alpha $-fractional Poisson extension to the upper
half-space $\mathbb{R}_{+}^{2}$,

\begin{equation*}
\mathbb{P}^{\alpha }\rho \left( x,t\right) \equiv \int_{\mathbb{R}}\frac{t}{%
\left( t^{2}+\left\vert x-y\right\vert ^{2}\right) ^{\frac{2-\alpha }{2}}}%
d\rho \left( y\right) ,
\end{equation*}%
so that in particular 
\begin{equation*}
\left\Vert \mathbb{P}^{\alpha }(f\sigma )\right\Vert _{L^{2}(\mathbb{R}%
_{+}^{2},\overline{\mu })}^{2}=\sum_{F\in \mathcal{F}}\sum_{M\in \mathcal{W}%
\left( F\right) }\mathbb{P}^{\alpha }\left( f\sigma \right) (c(M),\ell
\left( M\right) )^{2}\left\Vert \mathsf{Q}_{F,M}^{\omega ,\mathbf{b}^{\ast }}%
\frac{x}{\left\vert M\right\vert }\right\Vert _{L^{2}\left( \omega \right)
}^{\spadesuit 2}\,,
\end{equation*}%
and so (\ref{two weight Poisson n}) proves the first line in Proposition \ref%
{func ener control} upon inspecting (\ref{e.funcEnergy n}). Note also that
we can equivalently write $\left\Vert \mathbb{P}^{\alpha }\left( f\sigma
\right) \right\Vert _{L^{2}(\mathbb{R}_{+}^{2},\overline{\mu })}=\left\Vert 
\widetilde{\mathbb{P}}^{\alpha }\left( f\sigma \right) \right\Vert _{L^{2}(%
\mathbb{R}_{+}^{2},\mu )}$ where $\widetilde{\mathbb{P}}^{\alpha }\nu \left(
x,t\right) \equiv \frac{1}{t}\mathbb{P}^{\alpha }\nu \left( x,t\right) $ is
the renormalized Poisson operator. Here we have simply shifted the factor $%
\frac{1}{t^{2}}$ in $\overline{\mu }$ to $\left\vert \widetilde{\mathbb{P}}%
^{\alpha }\left( f\sigma \right) \right\vert ^{2}$ instead, and we will do
this shifting often throughout the proof when it is convenient to do so.

One version of the characterization of the two-weight inequality for
fractional and Poisson integrals in \cite{Saw3} was stated in terms of a
fixed dyadic grid $\mathcal{D}$ of intervals in $\mathbb{R}$ with sides
parallel to the coordinate axes. Using this theorem for the two-weight
Poisson inequality, but adapted to the $\alpha $-fractional Poisson integral 
$\mathbb{P}^{\alpha }$\footnote{%
The proof for $0\leq \alpha <1$ is essentially identical to that for $\alpha
=0$ given in \cite{Saw3}.}, we see that inequality (\ref{two weight Poisson
n}) requires checking these two inequalities for dyadic intervals $I\in 
\mathcal{D}$ and boxes $\widehat{I}=I\times \left[ 0,\ell \left( I\right)
\right) $ in the upper half-space $\mathbb{R}_{+}^{2}$: 
\begin{equation}
\int_{\mathbb{R}_{+}^{2}}\mathbb{P}^{\alpha }\left( \mathbf{1}_{I}\sigma
\right) \left( x,t\right) ^{2}d\overline{\mu }\left( x,t\right) \equiv
\left\Vert \mathbb{P}^{\alpha }\left( \mathbf{1}_{I}\sigma \right)
\right\Vert _{L^{2}(\overline{\mu })}^{2}\lesssim \left( \left( \mathfrak{E}%
_{2}^{\alpha }\right) ^{2}+\mathcal{A}_{2}^{\alpha }+\mathcal{A}_{2}^{\alpha
,\ast }+A_{2}^{\alpha ,\limfunc{punct}}\right) \sigma (I)\,,  \label{e.t1 n}
\end{equation}%
\begin{equation}
\int_{\mathbb{R}}[\mathbb{Q}^{\alpha }(t\mathbf{1}_{\widehat{I}}\overline{%
\mu })]^{2}d\sigma (x)\lesssim \left( \left( \mathfrak{E}_{2}^{\alpha
}\right) ^{2}+\mathcal{A}_{2}^{\alpha }+A_{2}^{\alpha ,\limfunc{punct}%
}\right) \int_{\widehat{I}}t^{2}d\overline{\mu }(x,t),  \label{e.t2 n}
\end{equation}%
for all \emph{dyadic} intervals $I\in \mathcal{D}$, and where the dual
Poisson operator $\mathbb{Q}^{\alpha }$ is given by 
\begin{equation*}
\mathbb{Q}^{\alpha }(t\mathbf{1}_{\widehat{I}}\overline{\mu })\left(
x\right) =\int_{\widehat{I}}\frac{t^{2}}{\left( t^{2}+\lvert x-y\rvert
^{2}\right) ^{\frac{2-\alpha }{2}}}d\overline{\mu }\left( y,t\right) \,.
\end{equation*}%
It is important to note that we can choose for $\mathcal{D}$ any fixed
dyadic grid, the compensating point being that the integrations on the left
sides of (\ref{e.t1 n}) and (\ref{e.t2 n}) are taken over the entire spaces $%
\mathbb{R}_{+}^{2}$ and $\mathbb{R}$ respectively\footnote{%
There is a gap in the proof of the Poisson inequality at the top of page 542
in \cite{Saw3}. However, this gap can be fixed as in \cite{SaWh} or \cite%
{LaSaUr1}.}.

\subsection{Poisson testing}

We now turn to proving the Poisson testing conditions (\ref{e.t1 n}) and (%
\ref{e.t2 n}). Similar testing conditions have been considered in \cite%
{SaShUr5}, \cite{SaShUr7}, \cite{SaShUr9} and \cite{SaShUr10}, and the
proofs there essentially carry over to the situation here, but careful
attention must now be paid to the changed definition of functional energy
and the weaker notion of goodness. We continue to circumvent the difficulty
of permitting common point masses here by using the energy Muckenhoupt
constants $A_{2}^{\alpha ,\limfunc{energy}}$ and $A_{2}^{\alpha ,\ast ,%
\limfunc{energy}}$, which require control by the punctured Muckenhoupt
constants $A_{2}^{\alpha ,\limfunc{punct}}$ and $A_{2}^{\alpha ,\ast ,%
\limfunc{punct}}$. The following elementary Poisson inequalities (see e.g. 
\cite{Vol}) will be used extensively.

\begin{lemma}
\label{Poisson inequalities}Suppose that $J,K,I$ are intervals in $\mathbb{R}
$, and that $\mu $ is a positive measure supported in $\mathbb{R}\setminus I$%
. If $J\subset K\subset \beta K\subset I$ for some $\beta >1$, then%
\begin{equation*}
\frac{\mathrm{P}^{\alpha }\left( J,\mu \right) }{\left\vert J\right\vert }%
\approx \frac{\mathrm{P}^{\alpha }\left( K,\mu \right) }{\left\vert
K\right\vert },
\end{equation*}%
while if $J\subset \beta K$, then%
\begin{equation*}
\frac{\mathrm{P}^{\alpha }\left( K,\mu \right) }{\left\vert K\right\vert }%
\lesssim \frac{\mathrm{P}^{\alpha }\left( J,\mu \right) }{\left\vert
J\right\vert }.
\end{equation*}
\end{lemma}

\begin{proof}
We have%
\begin{equation*}
\frac{\mathrm{P}^{\alpha }\left( J,\mu \right) }{\left\vert J\right\vert }=%
\frac{1}{\left\vert J\right\vert }\int \frac{\left\vert J\right\vert }{%
\left( \left\vert J\right\vert +\left\vert x-c_{J}\right\vert \right)
^{2-\alpha }}d\mu \left( x\right) ,
\end{equation*}%
where $J\subset K\subset \beta K\subset I$ implies that%
\begin{equation*}
\left\vert J\right\vert +\left\vert x-c_{J}\right\vert \approx \left\vert
K\right\vert +\left\vert x-c_{K}\right\vert ,\ \ \ \ \ x\in \mathbb{R}%
\setminus I,
\end{equation*}%
and where $J\subset \beta K$ implies that%
\begin{equation*}
\left\vert J\right\vert +\left\vert x-c_{J}\right\vert \lesssim \left\vert
J\right\vert +\left\vert c_{K}-c_{J}\right\vert +\left\vert
x-c_{K}\right\vert \lesssim \left\vert K\right\vert +\left\vert
x-c_{K}\right\vert ,\ \ \ \ \ x\in \mathbb{R}.
\end{equation*}
\end{proof}

Recall that in the case\ of the $T1$ theorem in \cite{SaShUr7}, where we
assumed \emph{traditional} goodness in a single family of grids $\mathcal{D}$%
, we had a \emph{strong} bounded overlap property associated with the
projections $\mathsf{P}_{F,J}^{\omega ,\mathbf{b}^{\ast }}$ defined there;
namely, that for each interval $I_{0}\in \mathcal{D}$, there were a bounded
number of intervals $F\in \mathcal{F}$ with the property that $F\supsetneqq
I_{0}\supset J$ for some $J\in \mathcal{M}_{\left( \mathbf{\rho }%
,\varepsilon \right) -\limfunc{deep}}\left( F\right) $ with $\mathsf{P}%
_{F,J}^{\omega ,\mathbf{b}^{\ast }}\neq 0$ (see the first part of Lemma 10.4
in \cite{SaShUr7}). However, we no longer have this strong bounded overlap
property when ordinary goodness is replaced with the \emph{weak} goodness of
Hyt\"{o}nen and Martikainen. Indeed, there may now be an \emph{unbounded}
number of intervals $F\in \mathcal{F}$ with $F\supsetneqq I_{0}\supset J$
and $\mathsf{P}_{F,J}^{\omega ,\mathbf{b}^{\ast }}\neq 0$, simply because
there can be $J^{\prime }\in \mathcal{G}$ with both $J^{\prime }\subset
I_{0} $ and $\left( J^{\prime }\right) ^{\maltese }$ \emph{arbitrarily}
large.

What will save us in obtaining the following lemma is that the Whitney
intervals $M$ in $\mathcal{W}\left( F\right) $ that happen to lie in some $%
I\in \mathcal{D}$ with $I\subset F$ have one of just two different forms: if 
$I$ shares an endpoint with $F$ then the intervals $M$ near that endpoint
are the same as those in $\mathcal{W}\left( I\right) $ - note that $F$ has
been replaced with $I$ here - while otherwise there are a bounded number of
Whitney intervals $M$ in $I$, and each such $M$ has side length comparable
to $\ell \left( I\right) $.

The next lemma will be used in bounding both of the local Poisson testing
conditions. Recall from Definition \ref{def dyadic}\ that $\mathcal{AD}$
consists of all augmented $\mathcal{D}$-dyadic intervals where $K$ is an
augmented dyadic interval if it is a union of $2$ $\mathcal{D}$-dyadic
intervals $K^{\prime }$ with $\ell \left( K^{\prime }\right) =\frac{1}{2}%
\ell \left( K\right) $.

\begin{lemma}
\label{refined lemma}Let $\mathcal{D}$ and $\mathcal{G}$ and $\mathcal{%
F\subset D}$ be grids and let $\left\{ \mathsf{Q}_{F,M}^{\omega ,\mathbf{b}%
^{\ast }}\right\} _{\substack{ F\in \mathcal{F}  \\ M\in \mathcal{W}\left(
F\right) }}$ be as in (\ref{def F,K}) above. For any augmented interval $%
I\in \mathcal{AD}$ define%
\begin{equation}
B\left( I\right) \equiv \sum_{F\in \mathcal{F}:\ F\supsetneqq I^{\prime }%
\text{ for some }I^{\prime }\in \mathfrak{C}\left( I\right) }\sum_{M\in 
\mathcal{W}\left( F\right) :\ M\subset I}\left( \frac{\mathrm{P}^{\alpha
}\left( M,\mathbf{1}_{I}\sigma \right) }{\left\vert M\right\vert }\right)
^{2}\left\Vert \mathsf{Q}_{F,M}^{\omega ,\mathbf{b}^{\ast }}x\right\Vert
_{L^{2}\left( \omega \right) }^{\spadesuit 2}\ .  \label{term B}
\end{equation}%
Then%
\begin{equation}
B\left( I\right) \lesssim \left( \left( \mathfrak{E}_{2}^{\alpha }\right)
^{2}+A_{2}^{\alpha ,\limfunc{energy}}\right) \left\vert I\right\vert
_{\sigma }\ .  \label{B bound}
\end{equation}
\end{lemma}

\begin{proof}
We first prove the bound (\ref{B bound}) for $B\left( I\right) $ ignoring
for the moment the possible case when $M=I$ in the sum defining $B\left(
I\right) $. So suppose that $I\in \mathcal{AD}$ is an augmented $\mathcal{D}$%
-dyadic interval. Define%
\begin{equation*}
\Lambda ^{\ast }\left( I\right) \equiv \left\{ M\subsetneqq I:M\in \mathcal{W%
}\left( F\right) \text{ for some }F\supsetneqq I^{\prime }\text{, }I^{\prime
}\in \mathfrak{C}\left( I\right) \text{ with }\mathsf{Q}_{F,M}^{\omega ,%
\mathbf{b}^{\ast }}x\neq 0\right\} ,
\end{equation*}%
and pigeonhole this collection as $\Lambda ^{\ast }\left( I\right)
=\dbigcup\limits_{I^{\prime }\in \mathfrak{C}\left( I\right) }\Lambda \left(
I^{\prime }\right) $, where for each $I^{\prime }\in \mathfrak{C}\left(
I\right) $ we define 
\begin{equation*}
\Lambda \left( I^{\prime }\right) \equiv \left\{ M\subset I^{\prime }:M\in 
\mathcal{W}\left( F\right) \text{ for some }F\supsetneqq I^{\prime }\text{
with }\mathsf{Q}_{F,M}^{\omega ,\mathbf{b}^{\ast }}x\neq 0\right\} .
\end{equation*}%
Consider first the case when $3I^{\prime }\subset F$, so that $d\left(
I^{\prime },\partial F\right) \geq \ell \left( I^{\prime }\right) $. Then if 
$M\in \mathcal{W}\left( F\right) $ for some $F\supsetneqq I^{\prime }$ we
have $\ell \left( M\right) =d\left( M,\partial F\right) $, and if in
addition $M\subset I^{\prime }$, then $M=I^{\prime }$. Consider the sum over
all $F\supsetneqq I^{\prime }=M$:%
\begin{eqnarray*}
B_{M}\left( I\right) &\equiv &\sum_{F\in \mathcal{F}:\ F\supsetneqq M\text{
for some }M\in \mathfrak{C}\left( I\right) \cap \mathcal{W}\left( F\right)
}\left( \frac{\mathrm{P}^{\alpha }\left( M,\mathbf{1}_{I}\sigma \right) }{%
\left\vert M\right\vert }\right) ^{2}\left\Vert \mathsf{Q}_{F,M}^{\omega ,%
\mathbf{b}^{\ast }}x\right\Vert _{L^{2}\left( \omega \right) }^{\spadesuit 2}
\\
&\leq &\left( \frac{\mathrm{P}^{\alpha }\left( M,\mathbf{1}_{I}\sigma
\right) }{\left\vert M\right\vert }\right) ^{2}\left\Vert \mathsf{Q}%
_{M}^{\omega ,\mathbf{b}^{\ast }}x\right\Vert _{L^{2}\left( \omega \right)
}^{\spadesuit 2}\lesssim \left( \frac{\mathrm{P}^{\alpha }\left( I,\mathbf{1}%
_{I}\sigma \right) }{\left\vert I\right\vert }\right) ^{2}\left\Vert \mathsf{%
Q}_{I}^{\omega ,\mathbf{b}^{\ast }}x\right\Vert _{L^{2}\left( \omega \right)
}^{\spadesuit 2}\lesssim A_{2}^{\alpha ,\limfunc{energy}}\left\vert
I\right\vert _{\sigma }\ ,
\end{eqnarray*}%
where we have used the definitions (\ref{def F,K}) and (\ref{large pseudo}).
Thus we have obtained the bound%
\begin{equation*}
\sum_{F\in \mathcal{F}:\ F\supsetneqq M\text{ for some }M\in \mathfrak{C}%
\left( I\right) \cap \mathcal{W}\left( F\right) }\left( \frac{\mathrm{P}%
^{\alpha }\left( M,\mathbf{1}_{I}\sigma \right) }{\left\vert M\right\vert }%
\right) ^{2}\left\Vert \mathsf{Q}_{F,M}^{\omega ,\mathbf{b}^{\ast
}}x\right\Vert _{L^{2}\left( \omega \right) }^{\spadesuit 2}\lesssim
A_{2}^{\alpha ,\limfunc{energy}}\left\vert I\right\vert _{\sigma }\ .
\end{equation*}

Now we turn to the case $3I^{\prime }\not\subset F$, i.e. when $\partial
I^{\prime }\cap \partial F$ consists of exactly one boundary point. In this
case, if both $M\subset I^{\prime }$ and $M\in \mathcal{W}\left( F\right) $
for some $F\supsetneqq I^{\prime }$, then we must have either $M\in \mathcal{%
W}\left( I^{\prime }\right) $ or $M\in \mathfrak{C}\left( I^{\prime }\right) 
$, since both $M$ and $I^{\prime }$ are then close to the same boundary
point in $\partial F$. Note that it is here that we use the Whitney
decompositions to full advantage. So again we can estimate%
\begin{eqnarray*}
&&\sum_{\substack{ F\in \mathcal{F}:\ F\supsetneqq I^{\prime }\text{ for
some }I^{\prime }\in \mathfrak{C}\left( I\right)  \\ 3I^{\prime }\not\subset
F}}\sum_{M\in \mathcal{W}\left( F\right) :\ M\subset I^{\prime }}\left( 
\frac{\mathrm{P}^{\alpha }\left( M,\mathbf{1}_{I}\sigma \right) }{\left\vert
M\right\vert }\right) ^{2}\left\Vert \mathsf{Q}_{F,M}^{\omega ,\mathbf{b}%
^{\ast }}x\right\Vert _{L^{2}\left( \omega \right) }^{\spadesuit 2} \\
&\leq &\sum_{M\in \left\{ \mathcal{W}\left( I^{\prime }\right) \cup 
\mathfrak{C}\left( I^{\prime }\right) \right\} \cap \mathcal{W}\left(
F\right) }\left( \frac{\mathrm{P}^{\alpha }\left( M,\mathbf{1}_{I}\sigma
\right) }{\left\vert M\right\vert }\right) ^{2}\left\Vert \mathsf{Q}%
_{M}^{\omega ,\mathbf{b}^{\ast }}x\right\Vert _{L^{2}\left( \omega \right)
}^{\spadesuit 2}\lesssim \left( \mathfrak{E}_{2}^{\alpha }\right)
^{2}\left\vert I\right\vert _{\sigma }\ .
\end{eqnarray*}

Finally, we consider the case $M=I$. In this case $I\in \mathcal{D}$ and so $%
F\supsetneqq I^{\prime }$ implies $F\supset I$ and we can estimate%
\begin{equation*}
\sum_{F\in \mathcal{F}:\ F\supset I}\left( \frac{\mathrm{P}^{\alpha }\left(
I,\mathbf{1}_{I}\sigma \right) }{\left\vert I\right\vert }\right)
^{2}\left\Vert \mathsf{Q}_{F,I}^{\omega ,\mathbf{b}^{\ast }}x\right\Vert
_{L^{2}\left( \omega \right) }^{\spadesuit 2}\leq \left( \frac{\mathrm{P}%
^{\alpha }\left( I,\mathbf{1}_{I}\sigma \right) }{\left\vert I\right\vert }%
\right) ^{2}\left\Vert \mathsf{Q}_{I}^{\omega ,\mathbf{b}^{\ast
}}x\right\Vert _{L^{2}\left( \omega \right) }^{\spadesuit 2}\lesssim
A_{2}^{\alpha ,\limfunc{energy}}\left\vert I\right\vert _{\sigma }\ .
\end{equation*}%
This completes the proof of Lemma \ref{refined lemma}.
\end{proof}

\subsection{The forward Poisson testing inequality}

Fix $I\in \mathcal{D}$. We split the integration on the left side of (\ref%
{e.t1 n}) into a local and global piece:%
\begin{equation*}
\int_{\mathbb{R}_{+}^{2}}\mathbb{P}^{\alpha }\left( \mathbf{1}_{I}\sigma
\right) ^{2}d\overline{\mu }=\int_{\widehat{I}}\mathbb{P}^{\alpha }\left( 
\mathbf{1}_{I}\sigma \right) ^{2}d\overline{\mu }+\int_{\mathbb{R}%
_{+}^{2}\setminus \widehat{I}}\mathbb{P}^{\alpha }\left( \mathbf{1}%
_{I}\sigma \right) ^{2}d\overline{\mu }\equiv \mathbf{Local}\left( I\right) +%
\mathbf{Global}\left( I\right) ,
\end{equation*}%
where more explicitly,%
\begin{eqnarray}
&&\mathbf{Local}\left( I\right) \equiv \int_{\widehat{I}}\left[ \mathbb{P}%
^{\alpha }\left( \mathbf{1}_{I}\sigma \right) \left( x,t\right) \right] ^{2}d%
\overline{\mu }\left( x,t\right) ;\ \ \ \ \ \overline{\mu }\equiv \frac{1}{%
t^{2}}\mu ,  \label{def local forward} \\
\text{i.e. }\overline{\mu } &\equiv &\ \sum_{F\in \mathcal{F}}\sum_{M\in 
\mathcal{W}\left( F\right) }\left\Vert \mathsf{Q}_{F,M}^{\omega ,\mathbf{b}%
^{\ast }}\frac{x}{\ell \left( M\right) }\right\Vert _{L^{2}\left( \omega
\right) }^{\spadesuit 2}\cdot \delta _{\left( c_{M},\ell \left( M\right)
\right) },  \notag
\end{eqnarray}%
where we recall $\mathsf{Q}_{F,M}^{\omega ,\mathbf{b}^{\ast }}$ is defined
in (\ref{def F,K}) above. Here is a brief schematic diagram of the
decompositions, with bounds in $\fbox{}$, used in this subsection:%
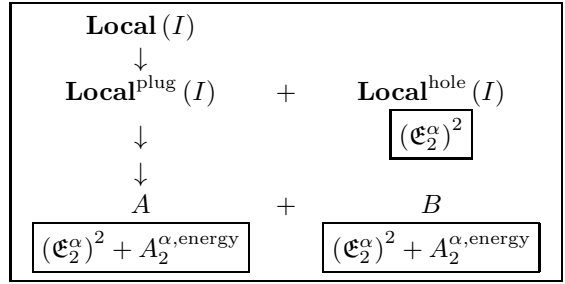
\begin{equation*}
\fbox{$%
\begin{array}{ccc}
\mathbf{Local}\left( I\right) &  &  \\ 
\downarrow &  &  \\ 
\mathbf{Local}^{\limfunc{plug}}\left( I\right) & + & \mathbf{Local}^{%
\limfunc{hole}}\left( I\right) \\ 
\downarrow &  & \fbox{$\left( \mathfrak{E}_{2}^{\alpha }\right) ^{2}$} \\ 
\downarrow &  &  \\ 
A & + & B \\ 
\fbox{$\left( \mathfrak{E}_{2}^{\alpha }\right) ^{2}+A_{2}^{\alpha ,\limfunc{%
energy}}$} &  & \fbox{$\left( \mathfrak{E}_{2}^{\alpha }\right)
^{2}+A_{2}^{\alpha ,\limfunc{energy}}$}%
\end{array}%
$}
\end{equation*}%
and%
\begin{equation*}
\fbox{$%
\begin{array}{ccccccc}
\mathbf{Global}\left( I\right) &  &  &  &  &  &  \\ 
\downarrow &  &  &  &  &  &  \\ 
A & + & B & + & C & + & D \\ 
\fbox{$A_{2}^{\alpha }$} &  & \fbox{$\left( \mathfrak{E}_{2}^{\alpha
}\right) ^{2}+A_{2}^{\alpha }+A_{2}^{\alpha ,\limfunc{energy}}$} &  & \fbox{$%
\mathcal{A}_{2}^{\alpha ,\ast }$} &  & \fbox{$\mathcal{A}_{2}^{\alpha ,\ast
}+A_{2}^{\alpha ,\limfunc{punct}}$}%
\end{array}%
$}.
\end{equation*}

As in our earlier papers \cite{SaShUr2}-\cite{SaShUr10} that used a single
family of random grids, we have the useful equivalence that%
\begin{equation}
\left( c\left( M\right) ,\ell \left( M\right) \right) \in \widehat{I}\text{ 
\textbf{if and only if} }M\subset I,  \label{tent consequence}
\end{equation}%
since $M$ and $I$ live in the common grid $\mathcal{D}$. We thus have

\begin{eqnarray*}
&&\mathbf{Local}\left( I\right) =\int_{\widehat{I}}\mathbb{P}^{\alpha
}\left( \mathbf{1}_{I}\sigma \right) \left( x,t\right) ^{2}d\overline{\mu }%
\left( x,t\right) \\
&=&\sum_{F\in \mathcal{F}}\sum_{M\in \mathcal{W}\left( F\right) :\ M\subset
I}\mathbb{P}^{\alpha }\left( \mathbf{1}_{I}\sigma \right) \left( c_{M},\ell
\left( M\right) \right) ^{2}\left\Vert \mathsf{Q}_{F,M}^{\omega ,\mathbf{b}%
^{\ast }}\frac{x}{\left\vert M\right\vert }\right\Vert _{L^{2}\left( \omega
\right) }^{\spadesuit 2} \\
&\approx &\sum_{F\in \mathcal{F}}\sum_{M\in \mathcal{W}\left( F\right) :\
M\subset I}\mathrm{P}^{\alpha }\left( M,\mathbf{1}_{I}\sigma \right)
^{2}\lVert \mathsf{Q}_{F,M}^{\omega ,\mathbf{b}^{\ast }}\frac{x}{\left\vert
M\right\vert }\rVert _{L^{2}\left( \omega \right) }^{\spadesuit 2} \\
&\approx &\mathbf{Local}^{\limfunc{plug}}\left( I\right) +\mathbf{Local}^{%
\func{hole}}\left( I\right) ,
\end{eqnarray*}%
where%
\begin{eqnarray*}
\mathbf{Local}^{\limfunc{plug}}\left( I\right) &\equiv &\sum_{F\in \mathcal{F%
}}\sum_{M\in \mathcal{W}\left( F\right) ):\ M\subset I}\left( \frac{\mathrm{P%
}^{\alpha }\left( M,\mathbf{1}_{F\cap I}\sigma \right) }{\left\vert
M\right\vert }\right) ^{2}\left\Vert \mathsf{Q}_{F,M}^{\omega ,\mathbf{b}%
^{\ast }}x\right\Vert _{L^{2}\left( \omega \right) }^{2}, \\
\mathbf{Local}^{\func{hole}}\left( I\right) &\equiv &\sum_{F\in \mathcal{F}%
}\sum_{M\in \mathcal{W}\left( F\right) :\ M\subset I}\left( \frac{\mathrm{P}%
^{\alpha }\left( M,\mathbf{1}_{I\setminus F}\sigma \right) }{\left\vert
M\right\vert }\right) ^{2}\left\Vert \mathsf{Q}_{F,M}^{\omega ,\mathbf{b}%
^{\ast }}x\right\Vert _{L^{2}\left( \omega \right) }^{\spadesuit 2}.
\end{eqnarray*}%
The `plugged' local sum $\mathbf{Local}^{\limfunc{plug}}\left( I\right) $
can be further decomposed into 
\begin{align*}
& \mathbf{Local}^{\limfunc{plug}}\left( I\right) =\left\{ \sum_{F\in 
\mathcal{F}:\ F\subset I}+\sum_{F\in \mathcal{F}:\ F\supsetneqq I}\right\}
\sum_{M\in \mathcal{W}\left( F\right) :\ M\subset I}\left( \frac{\mathrm{P}%
^{\alpha }\left( M,\mathbf{1}_{F\cap I}\sigma \right) }{\left\vert
M\right\vert }\right) ^{2}\left\Vert \mathsf{Q}_{F,M}^{\omega ,\mathbf{b}%
^{\ast }}x\right\Vert _{L^{2}\left( \omega \right) }^{\spadesuit 2} \\
& =A+B.
\end{align*}%
Then an application of the Whitney plugged energy condition gives 
\begin{eqnarray*}
A &=&\sum_{F\in \mathcal{F}:\ F\subset I}\sum_{M\in \mathcal{W}\left(
F\right) }\left( \frac{\mathrm{P}^{\alpha }\left( M,\mathbf{1}_{F\cap
I}\sigma \right) }{\left\vert M\right\vert }\right) ^{2}\left\Vert \mathsf{Q}%
_{F,M}^{\omega ,\mathbf{b}^{\ast }}x\right\Vert _{L^{2}\left( \omega \right)
}^{\spadesuit 2} \\
&\leq &\sum_{F\in \mathcal{F}:\ F\subset I}\left( \mathfrak{E}_{2}^{\alpha }+%
\sqrt{A_{2}^{\alpha ,\limfunc{energy}}}\right) ^{2}\left\vert F\right\vert
_{\sigma }\lesssim \left( \mathfrak{E}_{2}^{\alpha }+\sqrt{A_{2}^{\alpha ,%
\limfunc{energy}}}\right) ^{2}\left\vert I\right\vert _{\sigma }\,,
\end{eqnarray*}%
since $\left\Vert \mathsf{Q}_{F,M}^{\omega ,\mathbf{b}^{\ast }}x\right\Vert
_{L^{2}\left( \omega \right) }^{\spadesuit 2}\leq \left\Vert \mathsf{Q}%
_{M}^{\omega ,\mathbf{b}^{\ast }}x\right\Vert _{L^{2}\left( \omega \right)
}^{\spadesuit 2}$. We also used here that the stopping intervals $\mathcal{F}
$ satisfy a $\sigma $-Carleson measure estimate, 
\begin{equation*}
\sum_{F\in \mathcal{F}:\ F\subset F_{0}}\left\vert F\right\vert _{\sigma
}\lesssim \left\vert F_{0}\right\vert _{\sigma }.
\end{equation*}%
Lemma \ref{refined lemma} applies to the remaining term $B$ to obtain the
bound%
\begin{equation*}
B\lesssim \left( \left( \mathfrak{E}_{2}^{\alpha }\right) ^{2}+A_{2}^{\alpha
,\limfunc{energy}}\right) \left\vert I\right\vert _{\sigma }\ .
\end{equation*}

Next we show the inequality with `holes', where the support of $\sigma $ is
restricted to the complement of the interval $F$.

\begin{lemma}
\label{local hole}We have 
\begin{equation}
\mathbf{Local}^{\func{hole}}\left( I\right) \lesssim \left( \mathfrak{E}%
_{2}^{\alpha }\right) ^{2}\left\vert I\right\vert _{\sigma }\,.
\label{RTS n}
\end{equation}
\end{lemma}

\begin{proof}
Fix $I\in \mathcal{D}$ and define%
\begin{equation*}
\mathcal{F}_{I}\equiv \left\{ F\in \mathcal{F}:F\subset I\right\} \cup
\left\{ I\right\} ,
\end{equation*}%
and denote by $\pi F$, for this proof only, the parent of $F$ in the tree $%
\mathcal{F}_{I}$. Also denote by $d\left( F,F^{\prime }\right) \equiv d_{%
\mathcal{F}_{I}}\left( F,F^{\prime }\right) $ the distance from $F$ to $%
F^{\prime }$ in the tree $\mathcal{F}_{I}$, and denote by $d\left( F\right)
\equiv d_{\mathcal{F}_{I}}\left( F,I\right) $ the distance of $F$ from the
root $I$. Since $I\setminus F$ appears in the argument of the Poisson
integral, those $F\in \mathcal{F}\setminus \mathcal{F}_{I}$ do not
contribute to the sum and so we estimate%
\begin{equation*}
S\equiv \mathbf{Local}^{\func{hole}}\left( I\right) =\sum_{F\in \mathcal{F}%
_{I}}\sum_{M\in \mathcal{W}\left( F\right) :\ M\subset I}\left( \frac{%
\mathrm{P}^{\alpha }\left( M,\mathbf{1}_{I\setminus F}\sigma \right) }{%
\left\vert M\right\vert }\right) ^{2}\left\Vert \mathsf{Q}_{F,M}^{\omega ,%
\mathbf{b}^{\ast }}x\right\Vert _{L^{2}\left( \omega \right) }^{\spadesuit 2}
\end{equation*}%
by using $\sum_{F^{\prime }\in \mathcal{F}:\ F\subset F^{\prime }\subsetneqq
I}\frac{1}{d\left( F^{\prime }\right) ^{2}}\leq C$ to obtain\footnote{%
In \cite{SaShUr7} and \cite{SaShUr6} the first line of this display
incorrectly avoided the use of the Cauchy-Schwarz inequality. In the earlier
versions \cite{SaShUr5} and version \#2 of \cite{SaShUr6}, the argument was
correctly given by duality. The fix used here is taken from pages 94-95 of
version \#4 of \cite{SaShUr5}.} 
\begin{eqnarray*}
S &=&\sum_{F\in \mathcal{F}_{I}}\sum_{M\in \mathcal{W}\left( F\right) :\
M\subset I}\left( \sum_{F^{\prime }\in \mathcal{F}:\ F\subset F^{\prime
}\subsetneqq I}\frac{d\left( F^{\prime }\right) }{d\left( F^{\prime }\right) 
}\frac{\mathrm{P}^{\alpha }\left( M,\mathbf{1}_{\pi F^{\prime }\setminus
F^{\prime }}\sigma \right) }{\left\vert M\right\vert }\right) ^{2}\left\Vert 
\mathsf{Q}_{F,M}^{\omega ,\mathbf{b}^{\ast }}x\right\Vert _{L^{2}\left(
\omega \right) }^{\spadesuit 2} \\
&\leq &\sum_{F\in \mathcal{F}_{I}}\sum_{M\in \mathcal{W}\left( F\right) :\
M\subset I}\left( \sum_{F^{\prime }\in \mathcal{F}:\ F\subset F^{\prime
}\subsetneqq I}\frac{1}{d\left( F^{\prime }\right) ^{2}}\right) \\
&&\ \ \ \ \ \ \ \ \ \ \ \ \ \ \ \ \ \ \ \ \times \left( \sum_{F^{\prime }\in 
\mathcal{F}:\ F\subset F^{\prime }\subsetneqq I}d\left( F^{\prime }\right)
^{2}\left( \frac{\mathrm{P}^{\alpha }\left( M,\mathbf{1}_{\pi F^{\prime
}\setminus F^{\prime }}\sigma \right) }{\left\vert M\right\vert }\right)
^{2}\right) \left\Vert \mathsf{Q}_{F,M}^{\omega ,\mathbf{b}^{\ast
}}x\right\Vert _{L^{2}\left( \omega \right) }^{\spadesuit 2} \\
&\leq &C\sum_{F^{\prime }\in \mathcal{F}_{I}}d\left( F^{\prime }\right)
^{2}\sum_{F\in \mathcal{F}:\ F\subset F^{\prime }}\sum_{M\in \mathcal{W}%
\left( F\right) :\ M\subset I}\left( \frac{\mathrm{P}^{\alpha }\left( M,%
\mathbf{1}_{\pi F^{\prime }\setminus F^{\prime }}\sigma \right) }{\left\vert
M\right\vert }\right) ^{2}\left\Vert \mathsf{Q}_{F,M}^{\omega ,\mathbf{b}%
^{\ast }}x\right\Vert _{L^{2}\left( \omega \right) }^{\spadesuit 2} \\
&=&C\sum_{F^{\prime }\in \mathcal{F}_{I}}d\left( F^{\prime }\right)
^{2}\sum_{K\in \mathcal{W}\left( F^{\prime }\right) }\sum_{F\in \mathcal{F}%
:\ F\subset F^{\prime }}\sum_{M\in \mathcal{W}\left( F\right) :\ M\subset
I}\left( \frac{\mathrm{P}^{\alpha }\left( M,\mathbf{1}_{\pi F^{\prime
}\setminus F^{\prime }}\sigma \right) }{\left\vert M\right\vert }\right)
^{2}\left\Vert \mathsf{Q}_{F,M\cap K}^{\omega ,\mathbf{b}^{\ast
}}x\right\Vert _{L^{2}\left( \omega \right) }^{\spadesuit 2} \\
&\lesssim &\sum_{F^{\prime }\in \mathcal{F}_{I}}d\left( F^{\prime }\right)
^{2}\sum_{K\in \mathcal{W}\left( F^{\prime }\right) }\left( \frac{\mathrm{P}%
^{\alpha }\left( K,\mathbf{1}_{\pi F^{\prime }\setminus F^{\prime }}\sigma
\right) }{\left\vert K\right\vert }\right) ^{2}\sum_{F\in \mathcal{F}:\
F\subset F^{\prime }}\sum_{M\in \mathcal{W}\left( F\right) :\ M\subset
I}\left\Vert \mathsf{Q}_{F,M\cap K}^{\omega ,\mathbf{b}^{\ast }}x\right\Vert
_{L^{2}\left( \omega \right) }^{\spadesuit 2},
\end{eqnarray*}%
where in the fifth line we have used that each $J^{\prime }$ appearing in $%
\mathsf{Q}_{F,M}^{\omega ,\mathbf{b}^{\ast }}$ occurs in one of the $\mathsf{%
Q}_{F,M\cap K}^{\omega ,\mathbf{b}^{\ast }}$ since each $M$ is contained in
a unique $K$. We have also used there the Poisson inequalities in Lemma \ref%
{Poisson inequalities}.

We now use the lower frame inequality from Appendix A applied to the
function $\mathbf{1}_{K}\left( x-m_{K}^{\omega }\right) $ to obtain%
\begin{equation*}
\sum_{F\in \mathcal{F}:\ F\subset F^{\prime }}\sum_{M\in \mathcal{W}\left(
F\right) :\ M\subset I}\left\Vert \mathsf{Q}_{F,M\cap K}^{\omega ,\mathbf{b}%
^{\ast }}x\right\Vert _{L^{2}\left( \omega \right) }^{\spadesuit 2}\lesssim
\left\Vert \mathbf{1}_{K}\left( x-m_{K}^{\omega }\right) \right\Vert
_{L^{2}\left( \omega \right) }^{\spadesuit 2}\ .
\end{equation*}

Since the collection $\mathcal{F}_{I}$ satisfies a Carleson condition,
namely $\sum_{F\in \mathcal{F}_{I}}\left\vert F\cap I^{\prime }\right\vert
_{\sigma }\leq C\left\vert I^{\prime }\right\vert _{\sigma }$ for all
intervals $I^{\prime }$, we have geometric decay in generations:%
\begin{equation}
\sum_{F\in \mathcal{F}_{I}:\ d\left( F\right) =k}\left\vert F\right\vert
_{\sigma }\lesssim 2^{-\delta k}\left\vert I\right\vert _{\sigma }\ ,\ \ \ \
\ k\geq 0.  \label{geometric decay}
\end{equation}%
Indeed, with $m>2C$ we have for each $F^{\prime }\in \mathcal{F}_{I}$,%
\begin{equation}
\sum_{F\in \mathcal{F}_{I}:\ F\subset F^{\prime }\text{ and }d\left(
F,F^{\prime }\right) =m}\left\vert F\cap F^{\prime }\right\vert _{\sigma }<%
\frac{1}{2}\left\vert F^{\prime }\right\vert _{\sigma }\ ,  \label{half}
\end{equation}%
since otherwise%
\begin{equation*}
\sum_{F\in \mathcal{F}_{I}:\ F\subset F^{\prime }\text{ and }d\left(
F,F^{\prime }\right) \leq m}\left\vert F\cap F^{\prime }\right\vert _{\sigma
}\geq m\frac{1}{2}\left\vert F^{\prime }\right\vert _{\sigma }\ ,
\end{equation*}%
a contradiction. Now iterate (\ref{half}) to obtain (\ref{geometric decay}).

Thus we can write%
\begin{eqnarray*}
S &\lesssim &\sum_{F^{\prime }\in \mathcal{F}_{I}}d\left( F^{\prime }\right)
^{2}\sum_{K\in \mathcal{W}\left( F^{\prime }\right) }\left( \frac{\mathrm{P}%
^{\alpha }\left( K,\mathbf{1}_{\pi F^{\prime }\setminus F^{\prime }}\sigma
\right) }{\left\vert K\right\vert }\right) ^{2}\left\Vert \mathbf{1}%
_{K}\left( x-m_{K}^{\omega }\right) \right\Vert _{L^{2}\left( \omega \right)
}^{\spadesuit 2} \\
&=&\sum_{k=1}^{\infty }k^{2}\sum_{F^{\prime }\in \mathcal{F}_{I}:\ d\left(
F^{\prime }\right) =k}\sum_{K\in \mathcal{W}\left( F^{\prime }\right)
}\left( \frac{\mathrm{P}^{\alpha }\left( K,\mathbf{1}_{\pi F^{\prime
}\setminus F^{\prime }}\sigma \right) }{\left\vert K\right\vert }\right)
^{2}\left\Vert \mathbf{1}_{K}\left( x-m_{K}^{\omega }\right) \right\Vert
_{L^{2}\left( \omega \right) }^{\spadesuit 2}\equiv \sum_{k=1}^{\infty
}A_{k}\ ,
\end{eqnarray*}%
where $A_{k}$ is defined at the end of the above display. Hence using the
strong energy condition,%
\begin{eqnarray*}
A_{k} &=&k^{2}\sum_{F^{\prime }\in \mathcal{F}_{I}:\ d\left( F^{\prime
}\right) =k}\sum_{K\in \mathcal{W}\left( F^{\prime }\right) }\left( \frac{%
\mathrm{P}^{\alpha }\left( K,\mathbf{1}_{\pi F^{\prime }\setminus F^{\prime
}}\sigma \right) }{\left\vert K\right\vert }\right) ^{2}\left\Vert \mathbf{1}%
_{K}\left( x-m_{K}^{\omega }\right) \right\Vert _{L^{2}\left( \omega \right)
}^{\spadesuit 2} \\
&\lesssim &k^{2}\left( \mathfrak{E}_{2}^{\alpha }\right) ^{2}\sum_{F^{\prime
\prime }\in \mathcal{F}_{I}:\ d\left( F^{\prime \prime }\right)
=k-1}\left\vert F^{\prime \prime }\right\vert _{\sigma }\lesssim \left( 
\mathfrak{E}_{2}^{\alpha }\right) ^{2}k^{2}2^{-\delta k}\left\vert
I\right\vert _{\sigma }\ ,
\end{eqnarray*}%
where we have applied the strong energy condition for each $F^{\prime \prime
}\in \mathcal{F}_{I}$ with $d\left( F^{\prime \prime }\right) =k-1$ to obtain%
\begin{equation}
\sum_{F^{\prime }\in \mathcal{F}_{I}:\ \pi F^{\prime }=F^{\prime \prime
}}\sum_{K\in \mathcal{W}\left( F^{\prime }\right) }\left( \frac{\mathrm{P}%
^{\alpha }\left( K,\mathbf{1}_{F^{\prime \prime }\setminus F^{\prime
}}\sigma \right) }{\left\vert K\right\vert }\right) ^{2}\left\Vert \mathbf{1}%
_{K}\left( x-m_{K}^{\omega }\right) \right\Vert _{L^{2}\left( \omega \right)
}^{\spadesuit 2}\leq \left( \mathfrak{E}_{2}^{\alpha }\right) ^{2}\left\vert
F^{\prime \prime }\right\vert _{\sigma }\ .  \label{to obtain}
\end{equation}%
Finally then we obtain%
\begin{equation*}
S\lesssim \sum_{k=1}^{\infty }\left( \mathfrak{E}_{2}^{\alpha }\right)
^{2}k^{2}2^{-\delta k}\left\vert I\right\vert _{\sigma }\lesssim \left( 
\mathfrak{E}_{2}^{\alpha }\right) ^{2}\left\vert I\right\vert _{\sigma }\ ,
\end{equation*}%
which is (\ref{RTS n}).
\end{proof}

Altogether we have now proved the estimate $\mathbf{Local}\left( I\right)
\lesssim \left( \left( \mathfrak{E}_{2}^{\alpha }\right) ^{2}+A_{2}^{\alpha ,%
\limfunc{energy}}\right) \left\vert I\right\vert _{\sigma }$ when $I\in 
\mathcal{D}$, i.e. for every dyadic interval $I\in \mathcal{D}$,%
\begin{eqnarray}
&&  \label{local} \\
\mathbf{Local}\left( I\right) &\approx &\sum_{F\in \mathcal{F}}\sum_{M\in 
\mathcal{W}\left( F\right) :\ M\subset I}\left( \frac{\mathrm{P}^{\alpha
}\left( M,\mathbf{1}_{I}\sigma \right) }{\left\vert M\right\vert }\right)
^{2}\left\Vert \mathsf{Q}_{F,M}^{\omega ,\mathbf{b}^{\ast }}x\right\Vert
_{L^{2}\left( \omega \right) }^{\spadesuit 2}  \notag \\
&\lesssim &\left( \left( \mathfrak{E}_{2}^{\alpha }\right)
^{2}+A_{2}^{\alpha ,\limfunc{energy}}\right) \left\vert I\right\vert
_{\sigma },\ \ \ I\in \mathcal{D}.  \notag
\end{eqnarray}

\subsubsection{The augmented local estimate}

For future use in the `prepare to puncture' arguments below, we prove a
strengthening of the local estimate $\mathbf{Local}\left( I\right) $ to 
\emph{augmented} intervals $L\in \mathcal{AD}$.

\begin{lemma}
\label{shifted}With notation as above and $L\in \mathcal{AD}$ an augmented
interval, we have 
\begin{eqnarray}
&&  \label{shifted local} \\
\mathbf{Local}\left( L\right) &\equiv &\sum_{F\in \mathcal{F}}\sum_{M\in 
\mathcal{W}\left( F\right) :\ M\subset L}\left( \frac{\mathrm{P}^{\alpha
}\left( M,\mathbf{1}_{L}\sigma \right) }{\left\vert M\right\vert }\right)
^{2}\left\Vert \mathsf{Q}_{F,M}^{\omega ,\mathbf{b}^{\ast }}x\right\Vert
_{L^{2}\left( \omega \right) }^{\spadesuit 2}  \notag \\
&\lesssim &\left( \left( \mathfrak{E}_{2}^{\alpha }\right)
^{2}+A_{2}^{\alpha ,\limfunc{energy}}\right) \left\vert L\right\vert
_{\sigma },\ \ \ L\in \mathcal{AD}.  \notag
\end{eqnarray}
\end{lemma}

\begin{proof}
We prove (\ref{shifted local}) by repeating the above proof of (\ref{local})
and noting the points requiring change. First we decompose 
\begin{equation*}
\mathbf{Local}\left( L\right) \lesssim \mathbf{Local}^{\limfunc{plug}}\left(
L\right) +\mathbf{Local}^{\func{hole}}\left( L\right) +\mathbf{Local}^{%
\limfunc{offset}}\left( L\right)
\end{equation*}%
where $\mathbf{Local}^{\limfunc{plug}}\left( L\right) $ and $\mathbf{Local}^{%
\func{hole}}\left( L\right) $ are analogous to $\mathbf{Local}^{\limfunc{plug%
}}\left( I\right) $ and $\mathbf{Local}^{\func{hole}}\left( I\right) $
above, and where $\mathbf{Local}^{\limfunc{offset}}\left( L\right) $ is an
additional term arising because $L\setminus F$ need not be empty when $L\cap
F\neq \emptyset $ and $F$ is not contained in $L$:%
\begin{eqnarray*}
\mathbf{Local}^{\limfunc{plug}}\left( L\right) &\equiv &\sum_{F\in \mathcal{F%
}}\sum_{M\in \mathcal{W}\left( F\right) :\ M\subset L}\left( \frac{\mathrm{P}%
^{\alpha }\left( M,\mathbf{1}_{L\cap F}\sigma \right) }{\left\vert
M\right\vert }\right) ^{2}\left\Vert \mathsf{Q}_{F,M}^{\omega ,\mathbf{b}%
^{\ast }}x\right\Vert _{L^{2}\left( \omega \right) }^{\spadesuit 2}\ , \\
\mathbf{Local}^{\func{hole}}\left( L\right) &\equiv &\sum_{F\in \mathcal{F}%
:\ F\subset L}\sum_{M\in \mathcal{W}\left( F\right) :\ M\subset L}\left( 
\frac{\mathrm{P}^{\alpha }\left( M,\mathbf{1}_{L\setminus F}\sigma \right) }{%
\left\vert M\right\vert }\right) ^{2}\left\Vert \mathsf{Q}_{F,M}^{\omega ,%
\mathbf{b}^{\ast }}x\right\Vert _{L^{2}\left( \omega \right) }^{\spadesuit
2}\ , \\
\mathbf{Local}^{\limfunc{offset}}\left( L\right) &\equiv &\sum_{F\in 
\mathcal{F}:\ F\not\subset L}\sum_{M\in \mathcal{W}\left( F\right) :\
M\subset L}\left( \frac{\mathrm{P}^{\alpha }\left( M,\mathbf{1}_{L\setminus
F}\sigma \right) }{\left\vert M\right\vert }\right) ^{2}\left\Vert \mathsf{Q}%
_{F,M}^{\omega ,\mathbf{b}^{\ast }}x\right\Vert _{L^{2}\left( \omega \right)
}^{\spadesuit 2}\ .
\end{eqnarray*}%
We have%
\begin{align*}
& \mathbf{Local}^{\limfunc{plug}}\left( L\right) =\left\{ \sum_{F\in 
\mathcal{F}:\ F\subset \text{ some }L^{\prime }\in \mathfrak{C}\left(
L\right) }+\sum_{F\in \mathcal{F}:\ F\supsetneqq \text{ some }L^{\prime }\in 
\mathfrak{C}_{\mathcal{D}}\left( L\right) }\right\} \sum_{M\in \mathcal{W}%
\left( F\right) :\ M\subset L} \\
& \ \ \ \ \ \ \ \ \ \ \ \ \ \ \ \ \ \ \ \ \ \ \ \ \ \ \ \ \ \ \times \left( 
\frac{\mathrm{P}^{\alpha }\left( M,\mathbf{1}_{F\cap L}\sigma \right) }{%
\left\vert M\right\vert }\right) ^{2}\left\Vert \mathsf{Q}_{F,M}^{\omega ,%
\mathbf{b}^{\ast }}x\right\Vert _{L^{2}\left( \omega \right) }^{\spadesuit 2}
\\
& =A+B.
\end{align*}

Term $A$ satisfies%
\begin{equation*}
A\lesssim \left( \mathfrak{E}_{2}^{\alpha }+\sqrt{A_{2}^{\alpha ,\limfunc{%
energy}}}\right) ^{2}\left\vert L\right\vert _{\sigma }\ ,
\end{equation*}%
just as above using $\left\Vert \mathsf{Q}_{F,M}^{\omega }x\right\Vert
_{L^{2}\left( \omega \right) }^{2}\leq \left\Vert \mathsf{Q}_{M}^{\omega
}x\right\Vert _{L^{2}\left( \omega \right) }^{2}$, and the fact that the
stopping intervals $\mathcal{F}$ satisfy a $\sigma $-Carleson measure
estimate, 
\begin{equation*}
\sum_{F\in \mathcal{F}:\ F\subset L}\left\vert F\right\vert _{\sigma
}\lesssim \left\vert L\right\vert _{\sigma }.
\end{equation*}

Term $B$ is handled directly by Lemma \ref{refined lemma} with the augmented
interval $I=L$ to obtain%
\begin{equation*}
B\lesssim \left( \left( \mathfrak{E}_{2}^{\alpha }\right) ^{2}+A_{2}^{\alpha
,\limfunc{energy}}\right) \left\vert L\right\vert _{\sigma }\ .
\end{equation*}

To handle $\mathbf{Local}^{\func{hole}}\left( L\right) $, we define%
\begin{equation*}
\mathcal{F}_{L}\equiv \left\{ F\in \mathcal{F}:F\subset L\right\} \cup
\left\{ L\right\} ,
\end{equation*}%
and follow along the proof there with only trivial changes. The analogue of (%
\ref{to obtain}) is now%
\begin{equation*}
\sum_{F^{\prime }\in \mathcal{F}_{L}:\ \pi F^{\prime }=F^{\prime \prime
}}\sum_{K\in \mathcal{W}\left( F^{\prime }\right) }\left( \frac{\mathrm{P}%
^{\alpha }\left( K,\mathbf{1}_{F^{\prime \prime }\setminus F^{\prime
}}\sigma \right) }{\left\vert K\right\vert }\right) ^{2}\left\Vert \mathbf{1}%
_{K}\left( x-m_{K}^{\omega }\right) \right\Vert _{L^{2}\left( \omega \right)
}^{\spadesuit 2}\leq \left( \mathfrak{E}_{2}^{\alpha }\right) ^{2}\left\vert
F^{\prime \prime }\right\vert _{\sigma }\ ,
\end{equation*}%
the only change being that $\mathcal{F}_{L}$ now appears in place of $%
\mathcal{F}_{I}$, so that the energy condition still applies. We conclude
that 
\begin{equation*}
\mathbf{Local}^{\func{hole}}\left( L\right) \lesssim \left( \mathfrak{E}%
_{2}^{\alpha }\right) ^{2}\left\vert L\right\vert _{\sigma }\ .
\end{equation*}

Finally, the additional term $\mathbf{Local}^{\limfunc{offset}}\left(
L\right) $ is handled directly by Lemma \ref{refined lemma}, and this
completes the proof of the estimate (\ref{shifted local}) in Lemma \ref%
{shifted}.
\end{proof}

\subsubsection{The global estimate}

Now we turn to proving the following estimate for the global part of the
first testing condition \eqref{e.t1 n}:%
\begin{equation*}
\mathbf{Global}\left( I\right) =\int_{\mathbb{R}_{+}^{2}\setminus \widehat{I}%
}\mathbb{P}^{\alpha }\left( \mathbf{1}_{I}\sigma \right) ^{2}d\overline{\mu }%
\lesssim \left( \left( \mathfrak{E}_{2}^{\alpha }\right) ^{2}+\mathcal{A}%
_{2}^{\alpha ,\ast }+A_{2}^{\alpha ,\limfunc{punct}}\right) \left\vert
I\right\vert _{\sigma }.
\end{equation*}%
We begin by decomposing the integral above into four pieces. We have from (%
\ref{tent consequence}):%
\begin{eqnarray*}
&&\int_{\mathbb{R}_{+}^{2}\setminus \widehat{I}}\mathbb{P}^{\alpha }\left( 
\mathbf{1}_{I}\sigma \right) ^{2}d\overline{\mu }=\sum_{M:\ \left(
c_{M},\ell \left( M\right) \right) \in \mathbb{R}_{+}^{2}\setminus \widehat{I%
}}\mathbb{P}^{\alpha }\left( \mathbf{1}_{I}\sigma \right) \left( c_{M},\ell
\left( M\right) \right) ^{2}\sum_{\substack{ F\in \mathcal{F}:  \\ M\in 
\mathcal{W}\left( F\right) }}\left\Vert \mathsf{Q}_{F,M}^{\omega ,\mathbf{b}%
^{\ast }}\frac{x}{\left\vert M\right\vert }\right\Vert _{L^{2}\left( \omega
\right) }^{\spadesuit 2} \\
&=&\left\{ \sum_{\substack{ M\cap 3I=\emptyset  \\ \ell \left( M\right) \leq
\ell \left( I\right) }}+\sum_{M\subset 3I\setminus I}+\sum_{\substack{ M\cap
I=\emptyset  \\ \ell \left( M\right) >\ell \left( I\right) }}%
+\sum_{M\supsetneqq I}\right\} \mathbb{P}^{\alpha }\left( \mathbf{1}%
_{I}\sigma \right) \left( c_{M},\ell \left( M\right) \right) ^{2}\sum 
_{\substack{ F\in \mathcal{F}:  \\ M\in \mathcal{W}\left( F\right) }}%
\left\Vert \mathsf{Q}_{F,M}^{\omega ,\mathbf{b}^{\ast }}\frac{x}{\left\vert
M\right\vert }\right\Vert _{L^{2}\left( \omega \right) }^{\spadesuit 2} \\
&=&A+B+C+D.
\end{eqnarray*}

We further decompose term $A$ according to the length of $M$ and its
distance from $I$, and then use the pairwise disjointedness of the
projections $\mathsf{Q}_{F,M}^{\omega ,\mathbf{b}^{\ast }}$ in $F$ (see the
definition in (\ref{def F,K})) to obtain:%
\begin{eqnarray*}
A &\lesssim &\sum_{m=0}^{\infty }\sum_{k=1}^{\infty }\sum_{\substack{ %
M\subset 3^{k+1}I\setminus 3^{k}I  \\ \ell \left( M\right) =2^{-m}\ell
\left( I\right) }}\left( \frac{2^{-m}\left\vert I\right\vert }{d\left(
M,I\right) ^{2-\alpha }}\left\vert I\right\vert _{\sigma }\right)
^{2}\left\vert M\right\vert _{\omega } \\
&\lesssim &\sum_{m=0}^{\infty }2^{-2m}\sum_{k=1}^{\infty }\frac{\left\vert
I\right\vert ^{2}\left\vert I\right\vert _{\sigma }\left\vert
3^{k+1}I\setminus 3^{k}I\right\vert _{\omega }}{\left\vert 3^{k}I\right\vert
^{2\left( 2-\alpha \right) }}\left\vert I\right\vert _{\sigma } \\
&\lesssim &\sum_{m=0}^{\infty }2^{-2m}\sum_{k=1}^{\infty }3^{-2k}\left\{ 
\frac{\left\vert 3^{k+1}I\setminus 3^{k}I\right\vert _{\omega }\left\vert
3^{k}I\right\vert _{\sigma }}{\left\vert 3^{k}I\right\vert ^{2\left(
1-\alpha \right) }}\right\} \left\vert I\right\vert _{\sigma }\lesssim
A_{2}^{\alpha }\left\vert I\right\vert _{\sigma },
\end{eqnarray*}%
where the offset Muckenhoupt constant $A_{2}^{\alpha }$ applies because $%
3^{k+1}I$ has only three times the side length of $3^{k}I$.

\medskip

For term $B$ we first dispose of the nearby sum $B_{\limfunc{nearby}}$ that
consists of the sum over those $M$ which satisfy in addition $2^{-\mathbf{%
\rho }}\ell \left( I\right) \leq \ell \left( M\right) \leq \ell \left(
I\right) $. But it is a straightforward task to bound $B_{\limfunc{nearby}}$
by $CA_{2}^{\alpha ,\limfunc{energy}}\left\vert I\right\vert _{\sigma }$ as
there are at most $2^{\mathbf{\rho }+1}$ such intervals $M$. To bound $B_{%
\func{away}}\equiv B-B_{\limfunc{nearby}}$, we further decompose the sum
over $F\in \mathcal{F}$ according to whether or not $F\subset 3I\setminus I$:%
\begin{eqnarray*}
B_{\func{away}} &\approx &\sum_{M\subset 3I\setminus I\text{ and }\ell
\left( M\right) <2^{-\mathbf{\rho }}\ell \left( I\right) }\left( \frac{%
\mathrm{P}^{\alpha }\left( M,\mathbf{1}_{I}\sigma \right) }{\left\vert
M\right\vert }\right) ^{2}\sum_{\substack{ F\in \mathcal{F}:\ F\subset
3I\setminus I  \\ M\in \mathcal{W}\left( F\right) }}\left\Vert \mathsf{Q}%
_{F,M}^{\omega ,\mathbf{b}^{\ast }}x\right\Vert _{L^{2}\left( \omega \right)
}^{\spadesuit 2} \\
&&+\sum_{M\subset 3I\setminus I\text{ and }\ell \left( M\right) <2^{-\mathbf{%
\rho }}\ell \left( I\right) }\left( \frac{\mathrm{P}^{\alpha }\left( M,%
\mathbf{1}_{I}\sigma \right) }{\left\vert M\right\vert }\right) ^{2}\sum 
_{\substack{ F\in \mathcal{F}:\ F\not\subset 3I\setminus I  \\ M\in \mathcal{%
W}\left( F\right) }}\left\Vert \mathsf{Q}_{F,M}^{\omega ,\mathbf{b}^{\ast
}}x\right\Vert _{L^{2}\left( \omega \right) }^{\spadesuit 2} \\
&\equiv &B_{\func{away}}^{1}+B_{\func{away}}^{2}\ .
\end{eqnarray*}%
\ \ 

To estimate $B_{\func{away}}^{1}$, let 
\begin{equation}
\mathcal{J}^{\ast }\equiv \dbigcup\limits_{\substack{ F\in \mathcal{F}  \\ %
F\subset 3I\setminus I}}\dbigcup\limits_{\substack{ M\in \mathcal{W}\left(
F\right)  \\ M\subset 3I\setminus I\text{ and }\ell \left( M\right) <2^{-%
\mathbf{\rho }}\ell \left( I\right) }}\left\{ J\in \mathcal{C}_{F}^{\mathcal{%
G},\limfunc{shift}}:J\subset M\right\}  \label{def J*}
\end{equation}%
consist of all intervals $J\in \mathcal{G}$ for which the projection $%
\triangle _{J}^{\omega ,\mathbf{b}^{\ast }}$ occurs in one of the
projections $\mathsf{Q}_{F,M}^{\omega ,\mathbf{b}^{\ast }}$ in term $B_{%
\func{away}}^{1}$. In order to use $\mathcal{J}^{\ast }$ in the estimate for 
$B_{\func{away}}^{1}$ we need the following inequality. For any interval $%
M\in \mathcal{W}\left( F\right) $ we have%
\begin{eqnarray}
\left( \frac{\mathrm{P}^{\alpha }\left( M,\mathbf{1}_{I}\sigma \right) }{%
\left\vert M\right\vert }\right) ^{2}\left\Vert \mathsf{Q}_{F;M}^{\omega ,%
\mathbf{b}^{\ast }}x\right\Vert _{L^{2}\left( \omega \right) }^{\spadesuit
2} &=&\left( \frac{\mathrm{P}^{\alpha }\left( M,\mathbf{1}_{I}\sigma \right) 
}{\left\vert M\right\vert }\right) ^{2}\sum_{J\in \mathcal{C}_{F}^{\mathcal{G%
},\func{shift}}:\ J\subset M}\left\Vert \triangle _{J}^{\omega ,\mathbf{b}%
^{\ast }}x\right\Vert _{L^{2}\left( \omega \right) }^{\spadesuit 2}
\label{accomplished} \\
&\lesssim &\sum_{J\in \mathcal{C}_{F}^{\mathcal{G},\func{shift}}:\ J\subset
M}\left( \frac{\mathrm{P}^{\alpha }\left( J,\mathbf{1}_{I}\sigma \right) }{%
\left\vert J\right\vert }\right) ^{2}\left\Vert \triangle _{J}^{\omega ,%
\mathbf{b}^{\ast }}x\right\Vert _{L^{2}\left( \omega \right) }^{\spadesuit
2},  \notag
\end{eqnarray}%
since%
\begin{eqnarray*}
\frac{\mathrm{P}^{\alpha }\left( M,\mathbf{1}_{I}\sigma \right) }{\left\vert
M\right\vert } &=&\int_{I}\frac{1}{\left( \ell \left( M\right) +\left\vert
x-c_{M}\right\vert \right) ^{2-\alpha }}d\sigma \left( x\right) \\
&\lesssim &\int_{I}\frac{1}{\left( \ell \left( J\right) +\left\vert
x-c_{J}\right\vert \right) ^{2-\alpha }}d\sigma \left( x\right) =\frac{%
\mathrm{P}^{\alpha }\left( J,\mathbf{1}_{I}\sigma \right) }{\left\vert
J\right\vert }
\end{eqnarray*}%
for $J\subset M$ because%
\begin{equation*}
\ell \left( J\right) +\left\vert x-c_{J}\right\vert \lesssim \ell \left(
M\right) +\left\vert x-c_{M}\right\vert ,\ \ \ \ \ J\subset M\text{ and }%
x\in \mathbb{R}.
\end{equation*}%
We now use (\ref{accomplished}) to replace the sum over $M\in \mathcal{W}%
\left( F\right) $ in $B_{\func{away}}^{1}$, with a sum over $J\in \mathcal{J}%
^{\ast }$:%
\begin{eqnarray*}
B_{\func{away}}^{1} &=&\sum_{M\subset 3I\setminus I\text{ and }\ell \left(
M\right) <2^{-\mathbf{\rho }}\ell \left( I\right) }\left( \frac{\mathrm{P}%
^{\alpha }\left( M,\mathbf{1}_{I}\sigma \right) }{\left\vert M\right\vert }%
\right) ^{2}\sum_{\substack{ F\in \mathcal{F}:\ F\subset 3I\setminus I  \\ %
M\in \mathcal{W}\left( F\right) }}\left\Vert \mathsf{Q}_{F,M}^{\omega ,%
\mathbf{b}^{\ast }}x\right\Vert _{L^{2}\left( \omega \right) }^{\spadesuit 2}
\\
&\lesssim &\sum_{M\subset 3I\setminus I\text{ and }\ell \left( M\right) <2^{-%
\mathbf{\rho }}\ell \left( I\right) }\sum_{\substack{ F\in \mathcal{F}:\
F\subset 3I\setminus I  \\ M\in \mathcal{W}\left( F\right) }}\sum_{J\in 
\mathcal{C}_{F}^{\mathcal{G},\func{shift}}:\ J\subset M}\left( \frac{\mathrm{%
P}^{\alpha }\left( J,\mathbf{1}_{I}\sigma \right) }{\left\vert J\right\vert }%
\right) ^{2}\left\Vert \triangle _{J}^{\omega ,\mathbf{b}^{\ast
}}x\right\Vert _{L^{2}\left( \omega \right) }^{\spadesuit 2} \\
&\lesssim &\sum_{J\in \mathcal{J}^{\ast }}\left( \frac{\mathrm{P}^{\alpha
}\left( J,\mathbf{1}_{I}\sigma \right) }{\left\vert J\right\vert }\right)
^{2}\mathbf{\ }\left\Vert \bigtriangleup _{J}^{\omega ,\mathbf{b}^{\ast
}}x\right\Vert _{L^{2}\left( \omega \right) }^{\spadesuit 2}\ ,
\end{eqnarray*}%
where the final line follows since for each $J\in \mathcal{J}^{\ast }$ there
is a unique pair $\left( F,M\right) $ satisfying the conditions in the
second line.\bigskip

We will now exploit the smallness of $\varepsilon >0$ in the weak goodness
condition by decomposing the sum over $J\in \mathcal{J}^{\ast }$ according
to the length of $J$, and then using the fractional Poisson inequality (\ref%
{e.Jsimeq}) in Lemma \ref{Poisson inequality} on the neighbour $I^{\prime }$%
\ of $I$ containing $J$. Indeed, for $J\subset I^{\prime }\subset \mathbb{R}$
and $I\subset \mathbb{R}\setminus I^{\prime }$, we have 
\begin{equation}
\mathrm{P}^{\alpha }\left( J,\mathbf{1}_{I}\sigma \right) ^{2}\lesssim
\left( \frac{\ell \left( J\right) }{\ell \left( I\right) }\right)
^{2-2\left( 2-\alpha \right) \varepsilon }\mathrm{P}^{\alpha }\left( I,%
\mathbf{1}_{I}\sigma \right) ^{2},\ \ \ \ \ J\in \mathcal{J}^{\ast },
\label{Poisson inequalities 2}
\end{equation}%
where we have used that $\ell \left( I^{\prime }\right) =\ell \left(
I\right) $ and $\mathrm{P}^{\alpha }\left( I^{\prime },\mathbf{1}_{I}\sigma
\right) \approx \mathrm{P}^{\alpha }\left( I,\mathbf{1}_{I}\sigma \right) $,
and that the intervals $J\in \mathcal{J}^{\ast }$ are good in $I^{\prime }$
and beyond, and have side length at most $2^{-\mathbf{\rho }}\ell \left(
I\right) $, all because $J^{\maltese }\subset F\subset 3I\setminus I$ and we
have already dealt with the term $B_{\limfunc{nearby}}$. Moreover, we may
also assume here that the exponent $2-2\left( 2-\alpha \right) \varepsilon $
is positive, i.e.$\ \varepsilon <\frac{1}{2-\alpha }$, which is of course
implied by $0<\varepsilon <\frac{1}{2}$. We then obtain from (\ref{Poisson
inequalities 2}), the inequality $\left\Vert \bigtriangleup _{J}^{\omega ,%
\mathbf{b}^{\ast }}x\right\Vert _{L^{2}\left( \omega \right) }^{\spadesuit
2}\lesssim \left\vert J\right\vert ^{2}\left\vert J\right\vert _{\omega }$,
the pairwise disjointedness of the $M\in \mathcal{W}\left( F\right) $, the
uniqueness of $F$ with $J\in \mathcal{C}_{F}^{\mathcal{G},\limfunc{shift}}$,
and since $F\subset 3I\setminus I$ in the sum over $J\in \mathcal{J}^{\ast }$%
, that%
\begin{eqnarray*}
B_{\func{away}}^{1} &\lesssim &\sum_{J\in \mathcal{J}^{\ast }}\left( \frac{%
\mathrm{P}^{\alpha }\left( J,\mathbf{1}_{I}\sigma \right) }{\left\vert
J\right\vert }\right) ^{2}\mathbf{\ }\left\Vert \bigtriangleup _{J}^{\omega ,%
\mathbf{b}^{\ast }}x\right\Vert _{L^{2}\left( \omega \right) }^{\spadesuit
2}\lesssim \sum_{m=\mathbf{\rho }}^{\infty }\sum_{\substack{ J\in \mathcal{J}%
^{\ast }  \\ \ell \left( J\right) =2^{-m}\ell \left( I\right) }}\left(
2^{-m}\right) ^{2-2\left( 2-\alpha \right) \varepsilon }\mathrm{P}^{\alpha
}\left( I,\mathbf{1}_{I}\sigma \right) ^{2}\left\vert J\right\vert _{\omega }
\\
&\lesssim &\sum_{m=\mathbf{\rho }}^{\infty }\left( 2^{-m}\right) ^{2-2\left(
2-\alpha \right) \varepsilon }\left( \frac{\left\vert I\right\vert _{\sigma }%
}{\left\vert I\right\vert ^{1-\alpha }}\right) ^{2}\sum_{\substack{ J\subset
3I\setminus I  \\ \ell \left( J\right) =2^{-m}\ell \left( I\right) }}%
\left\vert J\right\vert _{\omega }\lesssim \sum_{m=\mathbf{\rho }}^{\infty
}\left( 2^{-m}\right) ^{2-2\left( 2-\alpha \right) \varepsilon }\frac{%
\left\vert I\right\vert _{\sigma }\left\vert 3I\setminus I\right\vert
_{\omega }}{\left\vert 3I\right\vert ^{2\left( 1-\alpha \right) }}\left\vert
I\right\vert _{\sigma }\lesssim A_{2}^{\alpha }\left\vert I\right\vert
_{\sigma }\ ,
\end{eqnarray*}%
since $2-2\left( 2-\alpha \right) \varepsilon >0$.

To complete the bound for term $B=B_{\limfunc{nearby}}+B_{\func{away}%
}^{1}+B_{\func{away}}^{2}$, it remains to estimate term $B_{\func{away}}^{2}$
in which we sum over $F\not\subset 3I\setminus I$. In this case $%
F\varsupsetneqq I^{\prime }$ for one of the two neighbours $I^{\prime }$ of $%
I$, and so we can apply Lemma \ref{refined lemma}, with $I$ there replaced
by the augmented intervals $I^{\prime }\cup I$, to obtain the estimate%
\begin{equation*}
B_{\func{away}}^{2}\lesssim \left( \left( \mathfrak{E}_{2}^{\alpha }\right)
^{2}+A_{2}^{\alpha ,\limfunc{energy}}\right) \left\vert I\right\vert
_{\sigma }\ .
\end{equation*}

\medskip

Next we turn to term $D$. The intervals $M$ occurring here are included in
the set of ancestors $A_{k}\equiv \pi _{\mathcal{D}}^{\left( k\right) }I$ of 
$I$, $1\leq k<\infty $.%
\begin{eqnarray*}
D &=&\sum_{k=1}^{\infty }\mathbb{P}^{\alpha }\left( \mathbf{1}_{I}\sigma
\right) \left( c\left( A_{k}\right) ,\left\vert A_{k}\right\vert \right)
^{2}\sum_{\substack{ F\in \mathcal{F}:  \\ A_{k}\in \mathcal{W}\left(
F\right) }}\left\Vert \mathsf{Q}_{F,A_{k}}^{\omega ,\mathbf{b}^{\ast }}\frac{%
x}{\lvert A_{k}\rvert }\right\Vert _{L^{2}\left( \omega \right)
}^{\spadesuit 2} \\
&=&\sum_{k=1}^{\infty }\mathbb{P}^{\alpha }\left( \mathbf{1}_{I}\sigma
\right) \left( c\left( A_{k}\right) ,\left\vert A_{k}\right\vert \right)
^{2}\sum_{\substack{ F\in \mathcal{F}:  \\ A_{k}\in \mathcal{W}\left(
F\right) }}\sum_{J^{\prime }\in \mathcal{C}_{F}^{\mathcal{G},\limfunc{shift}%
}:\ J^{\prime }\subset A_{k}\setminus I}\left\Vert \bigtriangleup
_{J^{\prime }}^{\omega ,\mathbf{b}^{\ast }}\frac{x}{\lvert A_{k}\rvert }%
\right\Vert _{L^{2}\left( \omega \right) }^{\spadesuit 2} \\
&&+\sum_{k=1}^{\infty }\mathbb{P}^{\alpha }\left( \mathbf{1}_{I}\sigma
\right) \left( c\left( A_{k}\right) ,\left\vert A_{k}\right\vert \right)
^{2}\sum_{\substack{ F\in \mathcal{F}:  \\ A_{k}\in \mathcal{W}\left(
F\right) }}\sum_{\substack{ J^{\prime }\in \mathcal{C}_{F}^{\mathcal{G},%
\limfunc{shift}}:\ J^{\prime }\subset A_{k}  \\ J^{\prime }\cap I\neq
\emptyset \text{ and }\ell \left( J^{\prime }\right) \leq \ell \left(
I\right) }}\left\Vert \bigtriangleup _{J^{\prime }}^{\omega ,\mathbf{b}%
^{\ast }}\frac{x}{\lvert A_{k}\rvert }\right\Vert _{L^{2}\left( \omega
\right) }^{\spadesuit 2} \\
&&+\sum_{k=1}^{\infty }\mathbb{P}^{\alpha }\left( \mathbf{1}_{I}\sigma
\right) \left( c\left( A_{k}\right) ,\left\vert A_{k}\right\vert \right)
^{2}\sum_{\substack{ F\in \mathcal{F}:  \\ A_{k}\in \mathcal{W}\left(
F\right) }}\sum_{\substack{ J^{\prime }\in \mathcal{C}_{F}^{\mathcal{G},%
\limfunc{shift}}:\ J^{\prime }\subset A_{k}  \\ J^{\prime }\cap I\neq
\emptyset \text{ and }\ell \left( J^{\prime }\right) >\ell \left( I\right) }}%
\left\Vert \bigtriangleup _{J^{\prime }}^{\omega ,\mathbf{b}^{\ast }}\frac{x%
}{\lvert A_{k}\rvert }\right\Vert _{L^{2}\left( \omega \right) }^{\spadesuit
2} \\
&\equiv &D_{\limfunc{disjoint}}+D_{\limfunc{descendent}}+D_{\limfunc{ancestor%
}}\ .
\end{eqnarray*}%
We thus have from the pairwise disjointedness of the projections $\mathsf{Q}%
_{F,A_{k}}^{\omega ,\mathbf{b}^{\ast }}$ in $F$ once again,%
\begin{eqnarray*}
D_{\limfunc{disjoint}} &=&\sum_{k=1}^{\infty }\mathbb{P}^{\alpha }\left( 
\mathbf{1}_{I}\sigma \right) \left( c\left( A_{k}\right) ,\left\vert
A_{k}\right\vert \right) ^{2}\sum_{\substack{ F\in \mathcal{F}:  \\ A_{k}\in 
\mathcal{W}\left( F\right) }}\sum_{J^{\prime }\in \mathcal{C}_{F}^{\mathcal{G%
},\limfunc{shift}}:\ J^{\prime }\subset A_{k}\setminus I}\left\Vert
\bigtriangleup _{J^{\prime }}^{\omega ,\mathbf{b}^{\ast }}\frac{x}{\lvert
A_{k}\rvert }\right\Vert _{L^{2}\left( \omega \right) }^{\spadesuit 2} \\
&\lesssim &\sum_{k=1}^{\infty }\left( \frac{\left\vert I\right\vert _{\sigma
}\left\vert A_{k}\right\vert }{\left\vert A_{k}\right\vert ^{2-\alpha }}%
\right) ^{2}\mathbf{\;}\left\vert A_{k}\setminus I\right\vert _{\omega
}=\left\{ \frac{\left\vert I\right\vert _{\sigma }}{\left\vert I\right\vert
^{1-\alpha }}\sum_{k=1}^{\infty }\frac{\left\vert I\right\vert ^{1-\alpha }}{%
\left\vert A_{k}\right\vert ^{2\left( 1-\alpha \right) }}\left\vert
A_{k}\setminus I\right\vert _{\omega }\right\} \left\vert I\right\vert
_{\sigma } \\
&\lesssim &\left\{ \frac{\left\vert I\right\vert _{\sigma }}{\left\vert
I\right\vert ^{1-\alpha }}\mathcal{P}^{\alpha }\left( I,\mathbf{1}%
_{I^{c}}\omega \right) \right\} \left\vert I\right\vert _{\sigma }\lesssim 
\mathcal{A}_{2}^{\alpha ,\ast }\left\vert I\right\vert _{\sigma },
\end{eqnarray*}%
since%
\begin{eqnarray*}
\sum_{k=1}^{\infty }\frac{\left\vert I\right\vert ^{1-\alpha }}{\left\vert
A_{k}\right\vert ^{2\left( 1-\alpha \right) }}\left\vert A_{k}\setminus
I\right\vert _{\omega } &=&\int \sum_{k=1}^{\infty }\frac{\left\vert
I\right\vert ^{1-\alpha }}{\left\vert A_{k}\right\vert ^{2\left( 1-\alpha
\right) }}\mathbf{1}_{A_{k}\setminus I}\left( x\right) d\omega \left(
x\right) \\
&=&\int \sum_{k=1}^{\infty }\frac{1}{2^{2\left( 1-\alpha \right) k}}\frac{%
\left\vert I\right\vert ^{1-\alpha }}{\left\vert I\right\vert ^{2\left(
1-\alpha \right) }}\mathbf{1}_{A_{k}\setminus I}\left( x\right) d\omega
\left( x\right) \\
&\lesssim &\int_{I^{c}}\left( \frac{\left\vert I\right\vert }{\left[
\left\vert I\right\vert +d\left( x,I\right) \right] ^{2}}\right) ^{1-\alpha
}d\omega \left( x\right) =\mathcal{P}^{\alpha }\left( I,\mathbf{1}%
_{I^{c}}\omega \right) ,
\end{eqnarray*}%
upon summing a geometric series with $2\left( 1-\alpha \right) >0$.

The next term $D_{\limfunc{descendent}}$ satisfies%
\begin{eqnarray*}
D_{\limfunc{descendent}} &\lesssim &\sum_{k=1}^{\infty }\left( \frac{%
\left\vert I\right\vert _{\sigma }\left\vert A_{k}\right\vert }{\left\vert
A_{k}\right\vert ^{2-\alpha }}\right) ^{2}\mathbf{\;}\left\Vert \mathsf{Q}%
_{3I}^{\omega ,\mathbf{b}^{\ast }}\frac{x}{2^{k}\lvert I\rvert }\right\Vert
_{L^{2}\left( \omega \right) }^{\spadesuit 2} \\
&=&\sum_{k=1}^{\infty }2^{-2k\left( 2-\alpha \right) }\left( \frac{%
\left\vert I\right\vert _{\sigma }}{\left\vert I\right\vert ^{1-\alpha }}%
\right) ^{2}\left\Vert \mathsf{Q}_{3I}^{\omega ,\mathbf{b}^{\ast }}\frac{x}{%
\lvert I\rvert }\right\Vert _{L^{2}\left( \omega \right) }^{\spadesuit 2} \\
&\lesssim &\left\{ \frac{\left\vert I\right\vert _{\sigma }\left\Vert 
\mathsf{Q}_{3I}^{\omega ,\mathbf{b}^{\ast }}\frac{x}{\lvert I\rvert }%
\right\Vert _{L^{2}\left( \omega \right) }^{\spadesuit 2}}{\left\vert
I\right\vert ^{2\left( 1-\alpha \right) }}\right\} \left\vert I\right\vert
_{\sigma }\lesssim A_{2}^{\alpha ,\limfunc{energy}}\left\vert I\right\vert
_{\sigma }\ .
\end{eqnarray*}

Lastly, for $D_{\limfunc{ancestor}}$ we note that there are at most two
intervals $K_{1}$ and $K_{2}$ in $\mathcal{G}$ having side length $\ell
\left( I\right) $ and such that $K_{i}\cap I\neq \emptyset $. Then each $%
J^{\prime }$ occurring in the sum in $D_{\limfunc{ancestor}}$ is of the form 
$J^{\prime }=A_{i}^{\ell }\equiv \pi _{\mathcal{G}}^{\left( \ell \right)
}K_{i}$ with $J^{\prime }\subset A_{k}$ for some $1\leq \ell \leq k$ and $%
i\in \left\{ 1,2\right\} $. Now we write%
\begin{eqnarray*}
D_{\limfunc{ancestor}} &=&\sum_{k=1}^{\infty }\mathbb{P}^{\alpha }\left( 
\mathbf{1}_{I}\sigma \right) \left( c\left( A_{k}\right) ,\left\vert
A_{k}\right\vert \right) ^{2}\sum_{\substack{ F\in \mathcal{F}:  \\ A_{k}\in 
\mathcal{W}\left( F\right) }}\sum_{\substack{ J^{\prime }\in \mathcal{C}%
_{F}^{\mathcal{G},\limfunc{shift}}:\ J^{\prime }\subset A_{k}  \\ J^{\prime
}\cap I\neq \emptyset \text{ and }\ell \left( J^{\prime }\right) >\ell
\left( I\right) }}\left\Vert \bigtriangleup _{J^{\prime }}^{\omega ,\mathbf{b%
}^{\ast }}\frac{x}{\lvert A_{k}\rvert }\right\Vert _{L^{2}\left( \omega
\right) }^{\spadesuit 2} \\
&\lesssim &\sum_{k=1}^{\infty }\left( \frac{\left\vert I\right\vert _{\sigma
}\left\vert A_{k}\right\vert }{\left\vert A_{k}\right\vert ^{2-\alpha }}%
\right) ^{2}\sum_{i=1}^{2}\sum_{\ell =1}^{k}\left\Vert \bigtriangleup
_{A_{i}^{\ell }}^{\omega ,\mathbf{b}^{\ast }}\frac{x}{\lvert A_{k}\rvert }%
\right\Vert _{L^{2}\left( \omega \right) }^{\spadesuit 2} \\
&\leq &2\sum_{k=1}^{\infty }\left( \frac{\left\vert I\right\vert _{\sigma
}\left\vert A_{k}\right\vert }{\left\vert A_{k}\right\vert ^{2-\alpha }}%
\right) ^{2}\left\Vert \mathsf{Q}_{A_{k}}^{\omega ,\mathbf{b}^{\ast }}\frac{x%
}{\lvert A_{k}\rvert }\right\Vert _{L^{2}\left( \omega \right) }^{\spadesuit
2}.
\end{eqnarray*}%
At this point we need a \emph{`prepare to puncture'} argument, as we will
want to derive geometric decay from $\left\Vert \mathsf{Q}_{J^{\prime
}}^{\omega ,\mathbf{b}^{\ast }}x\right\Vert _{L^{2}\left( \omega \right)
}^{\spadesuit 2}$ by dominating it by the `nonenergy' term $\left\vert
J^{\prime }\right\vert ^{2}\left\vert J^{\prime }\right\vert _{\omega }$, as
well as using the Muckenhoupt energy constant. For this we define $%
\widetilde{\omega }=\omega -\omega \left( \left\{ p\right\} \right) \delta
_{p}$ where $p$ is an atomic point in $I$ for which 
\begin{equation*}
\omega \left( \left\{ p\right\} \right) =\sup_{q\in \mathfrak{P}_{\left(
\sigma ,\omega \right) }:\ q\in I}\omega \left( \left\{ q\right\} \right) .
\end{equation*}%
(If $\omega $ has no atomic point in common with $\sigma $ in $I$ set $%
\widetilde{\omega }=\omega $.) Then we have $\left\vert I\right\vert _{%
\widetilde{\omega }}=\omega \left( I,\mathfrak{P}_{\left( \sigma ,\omega
\right) }\right) $ and%
\begin{equation*}
\frac{\left\vert I\right\vert _{\widetilde{\omega }}}{\left\vert
I\right\vert ^{1-\alpha }}\frac{\left\vert I\right\vert _{\sigma }}{%
\left\vert I\right\vert ^{1-\alpha }}=\frac{\omega \left( I,\mathfrak{P}%
_{\left( \sigma ,\omega \right) }\right) }{\left\vert I\right\vert
^{1-\alpha }}\frac{\left\vert I\right\vert _{\sigma }}{\left\vert
I\right\vert ^{1-\alpha }}\leq A_{2}^{\alpha ,\limfunc{punct}}.
\end{equation*}%
A key observation, already noted in the proof of Lemma \ref{energy A2}
above, is that%
\begin{equation}
\left\Vert \bigtriangleup _{K}^{\omega ,\mathbf{b}^{\ast }}x\right\Vert
_{L^{2}\left( \omega \right) }^{2}=\left\{ 
\begin{array}{ccc}
\left\Vert \bigtriangleup _{K}^{\omega ,\mathbf{b}^{\ast }}\left( x-p\right)
\right\Vert _{L^{2}\left( \omega \right) }^{2} & \text{ if } & p\in K \\ 
\left\Vert \bigtriangleup _{K}^{\omega ,\mathbf{b}^{\ast }}x\right\Vert
_{L^{2}\left( \widetilde{\omega }\right) }^{2} & \text{ if } & p\notin K%
\end{array}%
\right. \leq \ell \left( K\right) ^{2}\left\vert K\right\vert _{\widetilde{%
\omega }},\ \ \ \ \ \text{for all }K\in \mathcal{D}\ ,  \label{key obs}
\end{equation}%
and so, as in the proof of (\ref{omega tilda}) in Lemma \ref{energy A2},%
\begin{equation*}
\left\Vert \mathsf{Q}_{A_{k}}^{\omega ,\mathbf{b}^{\ast }}\frac{x}{%
\left\vert A_{k}\right\vert }\right\Vert _{L^{2}\left( \omega \right)
}^{\spadesuit 2}\lesssim \left\vert A_{k}\right\vert _{\widetilde{\omega }}\
.
\end{equation*}%
Then we continue with%
\begin{eqnarray*}
&&\sum_{k=1}^{\infty }\left( \frac{\left\vert I\right\vert _{\sigma
}\left\vert A_{k}\right\vert }{\left\vert A_{k}\right\vert ^{2-\alpha }}%
\right) ^{2}\left\Vert \mathsf{Q}_{A_{k}}^{\omega ,\mathbf{b}^{\ast }}\frac{x%
}{\lvert A_{k}\rvert }\right\Vert _{L^{2}\left( \omega \right) }^{\spadesuit
2} \\
&\lesssim &\sum_{k=1}^{\infty }\left( \frac{\left\vert I\right\vert _{\sigma
}\left\vert A_{k}\right\vert }{\left\vert A_{k}\right\vert ^{2-\alpha }}%
\right) ^{2}\left\vert A_{k}\right\vert _{\widetilde{\omega }} \\
&=&\sum_{k=1}^{\infty }\left( \frac{\left\vert I\right\vert _{\sigma }}{%
\left\vert A_{k}\right\vert ^{1-\alpha }}\right) ^{2}\left\vert
A_{k}\setminus I\right\vert _{\omega }+\sum_{k=1}^{\infty }\left( \frac{%
\left\vert I\right\vert _{\sigma }}{2^{k\left( 1-\alpha \right) }\left\vert
I\right\vert ^{1-\alpha }}\right) ^{2}\left\vert I\right\vert _{\widetilde{%
\omega }} \\
&\lesssim &\left( \mathcal{A}_{2}^{\alpha ,\ast }+A_{2}^{\alpha ,\limfunc{%
punct}}\right) \left\vert I\right\vert _{\sigma },
\end{eqnarray*}%
where the inequality $\sum_{k=1}^{\infty }\left( \frac{\left\vert
I\right\vert _{\sigma }}{\left\vert A_{k}\right\vert ^{1-\alpha }}\right)
^{2}\left\vert A_{k}\setminus I\right\vert _{\omega }\lesssim \mathcal{A}%
_{2}^{\alpha ,\ast }\left\vert I\right\vert _{\sigma }$ is already proved
above in the display estimating $D_{\limfunc{disjoint}}$.

\medskip

Finally, for term $C$ we will have to group the intervals $M$ into blocks $%
B_{i}$. We first split the sum according to whether or not $I$ intersects
the triple of $M$:%
\begin{eqnarray*}
C &\approx &\left\{ \sum_{\substack{ M:\ I\cap 3M=\emptyset  \\ \ell \left(
M\right) >\ell \left( I\right) }}+\sum_{\substack{ M:\ I\subset 3M\setminus
M  \\ \ell \left( M\right) >\ell \left( I\right) }}\right\} \left( \frac{%
\left\vert M\right\vert }{\left( \left\vert M\right\vert +d\left( M,I\right)
\right) ^{2-\alpha }}\left\vert I\right\vert _{\sigma }\right) ^{2}\sum 
_{\substack{ F\in \mathcal{F}:  \\ M\in \mathcal{W}\left( F\right) }}%
\left\Vert \mathsf{Q}_{F,M}^{\omega \mathbf{b}^{\ast }}\frac{x}{\left\vert
M\right\vert }\right\Vert _{L^{2}\left( \omega \right) }^{\spadesuit 2} \\
&=&C_{1}+C_{2}.
\end{eqnarray*}%
We first consider $C_{1}$. Let $\mathcal{M}$ consist of the maximal dyadic
intervals in the collection $\left\{ Q:3Q\cap I=\emptyset \right\} $, and
then let $\left\{ B_{i}\right\} _{i=1}^{\infty }$ be an enumeration of those 
$Q\in \mathcal{M}$ whose side length is at least $\ell \left( I\right) $.
Note in particular that $3B_{i}\cap I=\emptyset $. Now we further decompose
the sum in $C_{1}$ by grouping the intervals $M$ into the `Whitney'
intervals $B_{i}$, and then using the pairwise disjointedness of the
martingale supports of the pseudoprojections $\mathsf{Q}_{F,M}^{\omega ,%
\mathbf{b}^{\ast }}$ in $F$: 
\begin{eqnarray*}
C_{1} &\leq &\sum_{i=1}^{\infty }\sum_{M:\ M\subset B_{i}}\left( \frac{1}{%
\left( \left\vert M\right\vert +d\left( M,I\right) \right) ^{2-\alpha }}%
\left\vert I\right\vert _{\sigma }\right) ^{2}\sum_{\substack{ F\in \mathcal{%
F}:  \\ M\in \mathcal{W}\left( F\right) }}\left\Vert \mathsf{Q}%
_{F,M}^{\omega ,\mathbf{b}^{\ast }}x\right\Vert _{L^{2}\left( \omega \right)
}^{\spadesuit 2} \\
&\lesssim &\sum_{i=1}^{\infty }\left( \frac{1}{\left( \left\vert
B_{i}\right\vert +d\left( B_{i},I\right) \right) ^{2-\alpha }}\left\vert
I\right\vert _{\sigma }\right) ^{2}\sum_{M:\ M\subset B_{i}}\sum_{\substack{ %
F\in \mathcal{F}:  \\ M\in \mathcal{W}\left( F\right) }}\left\Vert \mathsf{Q}%
_{F,M}^{\omega ,\mathbf{b}^{\ast }}x\right\Vert _{L^{2}\left( \omega \right)
}^{\spadesuit 2} \\
&\lesssim &\sum_{i=1}^{\infty }\left( \frac{1}{\left( \left\vert
B_{i}\right\vert +d\left( B_{i},I\right) \right) ^{2-\alpha }}\left\vert
I\right\vert _{\sigma }\right) ^{2}\sum_{M:\ M\subset B_{i}}\left\vert
M\right\vert ^{2}\left\vert M\right\vert _{\omega } \\
&\lesssim &\sum_{i=1}^{\infty }\left( \frac{1}{\left( \left\vert
B_{i}\right\vert +d\left( B_{i},I\right) \right) ^{2-\alpha }}\left\vert
I\right\vert _{\sigma }\right) ^{2}\mathbf{\ }\left\vert B_{i}\right\vert
^{2}\left\vert B_{i}\right\vert _{\omega } \\
&\lesssim &\left\{ \sum_{i=1}^{\infty }\frac{\left\vert B_{i}\right\vert
_{\omega }\left\vert I\right\vert _{\sigma }}{\left\vert B_{i}\right\vert
^{2\left( 1-\alpha \right) }}\right\} \left\vert I\right\vert _{\sigma }\ ,
\end{eqnarray*}%
Now since $\left\vert B_{i}\right\vert \approx d\left( x,I\right) $ for $%
x\in B_{i}$, 
\begin{eqnarray*}
\sum_{i=1}^{\infty }\frac{\left\vert B_{i}\right\vert _{\omega }\left\vert
I\right\vert _{\sigma }}{\left\vert B_{i}\right\vert ^{2\left( 1-\alpha
\right) }} &=&\frac{\left\vert I\right\vert _{\sigma }}{\left\vert
I\right\vert ^{1-\alpha }}\sum_{i=1}^{\infty }\frac{\left\vert I\right\vert
^{1-\alpha }}{\left\vert B_{i}\right\vert ^{2\left( 1-\alpha \right) }}%
\left\vert B_{i}\right\vert _{\omega } \\
&\approx &\frac{\left\vert I\right\vert _{\sigma }}{\left\vert I\right\vert
^{1-\alpha }}\sum_{i=1}^{\infty }\int_{B_{i}}\frac{\left\vert I\right\vert
^{1-\alpha }}{d\left( x,I\right) ^{2\left( 1-\alpha \right) }}d\omega \left(
x\right) \\
&\approx &\frac{\left\vert I\right\vert _{\sigma }}{\left\vert I\right\vert
^{1-\alpha }}\sum_{i=1}^{\infty }\int_{B_{i}}\left( \frac{\left\vert
I\right\vert }{\left[ \left\vert I\right\vert +d\left( x,I\right) \right]
^{2}}\right) ^{1-\alpha }d\omega \left( x\right) \\
&\leq &\frac{\left\vert I\right\vert _{\sigma }}{\left\vert I\right\vert
^{1-\alpha }}\mathcal{P}^{\alpha }\left( I,\mathbf{1}_{I^{c}}\omega \right)
\leq \mathcal{A}_{2}^{\alpha ,\ast },
\end{eqnarray*}%
we obtain%
\begin{equation*}
C_{1}\lesssim \mathcal{A}_{2}^{\alpha ,\ast }\left\vert I\right\vert
_{\sigma }\ .
\end{equation*}

Next we turn to estimating term $C_{2}$ where the triple of $M$ contains $I$
but $M$ itself does not. Note that there are at most two such intervals $M$
of a given side length. So with this in mind, we sum over the intervals $M$
according to their lengths to obtain%
\begin{eqnarray*}
C_{2} &=&\sum_{m=1}^{\infty }\sum_{\substack{ M:\ I\subset 3M\setminus M  \\ %
\ell \left( M\right) =2^{m}\ell \left( I\right) }}\left( \frac{\left\vert
M\right\vert }{\left( \left\vert M\right\vert +\limfunc{dist}\left(
M,I\right) \right) ^{2-\alpha }}\left\vert I\right\vert _{\sigma }\right)
^{2}\sum_{\substack{ F\in \mathcal{F}:  \\ M\in \mathcal{W}\left( F\right) }}%
\left\Vert \mathsf{Q}_{F,M}^{\omega ,\mathbf{b}^{\ast }}\frac{x}{\left\vert
M\right\vert }\right\Vert _{L^{2}\left( \omega \right) }^{\spadesuit 2} \\
&\lesssim &\sum_{m=1}^{\infty }\left( \frac{\left\vert I\right\vert _{\sigma
}}{\left\vert 2^{m}I\right\vert ^{1-\alpha }}\right) ^{2}\mathbf{\ }%
\left\vert \left( 5\cdot 2^{m}I\right) \setminus I\right\vert _{\omega
}=\left\{ \frac{\left\vert I\right\vert _{\sigma }}{\left\vert I\right\vert
^{1-\alpha }}\sum_{m=1}^{\infty }\frac{\left\vert I\right\vert ^{1-\alpha
}\left\vert \left( 5\cdot 2^{m}I\right) \setminus I\right\vert _{\omega }}{%
\left\vert 2^{m}I\right\vert ^{2\left( 1-\alpha \right) }}\right\}
\left\vert I\right\vert _{\sigma } \\
&\lesssim &\left\{ \frac{\left\vert I\right\vert _{\sigma }}{\left\vert
I\right\vert ^{1-\alpha }}\mathcal{P}^{\alpha }\left( I,\mathbf{1}%
_{I^{c}}\omega \right) \right\} \left\vert I\right\vert _{\sigma }\leq 
\mathcal{A}_{2}^{\alpha ,\ast }\left\vert I\right\vert _{\sigma },
\end{eqnarray*}%
since in analogy with the corresponding estimate above,%
\begin{equation*}
\sum_{m=1}^{\infty }\frac{\left\vert I\right\vert ^{1-\alpha }\left\vert
\left( 5\cdot 2^{m}I\right) \setminus I\right\vert _{\omega }}{\left\vert
2^{m}I\right\vert ^{2\left( 1-\alpha \right) }}=\int \sum_{m=1}^{\infty }%
\frac{\left\vert I\right\vert ^{1-\alpha }}{\left\vert 2^{m}I\right\vert
^{2\left( 1-\alpha \right) }}\mathbf{1}_{\left( 5\cdot 2^{m}I\right)
\setminus I}\left( x\right) \ d\omega \left( x\right) \lesssim \mathcal{P}%
^{\alpha }\left( I,\mathbf{1}_{I^{c}}\omega \right) .
\end{equation*}

\subsection{The backward Poisson testing inequality\label{Subsec back test}}

The argument here follows the broad outline of the analogous argument in 
\cite{SaShUr7}, but using modifications from \cite{SaShUr9} that involve
`prepare to puncture arguments', using decompositions $\mathcal{W}\left(
F\right) $ in place of $\left( \mathbf{\rho },\varepsilon \right) $%
-decompositions, and using pseudoprojections $\mathsf{Q}_{F,M}^{\omega ,%
\mathbf{b}^{\ast }}x$ (see (\ref{def F,K}) for the definition). The final
change here is that there is no splitting into local and global parts as in 
\cite{SaShUr7} - instead, we follow the treatment in \cite{SaShUr6} in this
regard.

Fix $I\in \mathcal{D}$. It suffices to prove%
\begin{equation}
\mathbf{Back}\left( \widehat{I}\right) \equiv \int_{\mathbb{R}}\left[ 
\mathbb{Q}^{\alpha }\left( t\mathbf{1}_{\widehat{I}}\overline{\mu }\right)
\left( y\right) \right] ^{2}d\sigma (y)\lesssim \left\{ \mathcal{A}%
_{2}^{\alpha }+\left( \mathfrak{E}_{2}^{\alpha }+\sqrt{A_{2}^{\alpha ,%
\limfunc{energy}}}\right) \sqrt{A_{2}^{\alpha ,\limfunc{punct}}}\right\}
\int_{\widehat{I}}t^{2}d\overline{\mu }(x,t).  \label{e.t2 n'}
\end{equation}%
Note that for a `Poisson integral with holes' and a measure $\mu $ built
with Haar projections, Hyt\"{o}nen obtained in \cite{Hyt2} the simpler bound 
$A_{2}^{\alpha }$ for a term analogous to, but significantly smaller than, (%
\ref{e.t2 n'}). Here is a brief schematic diagram of the decompositions,
with bounds in $\fbox{}$, used in this subsection:%
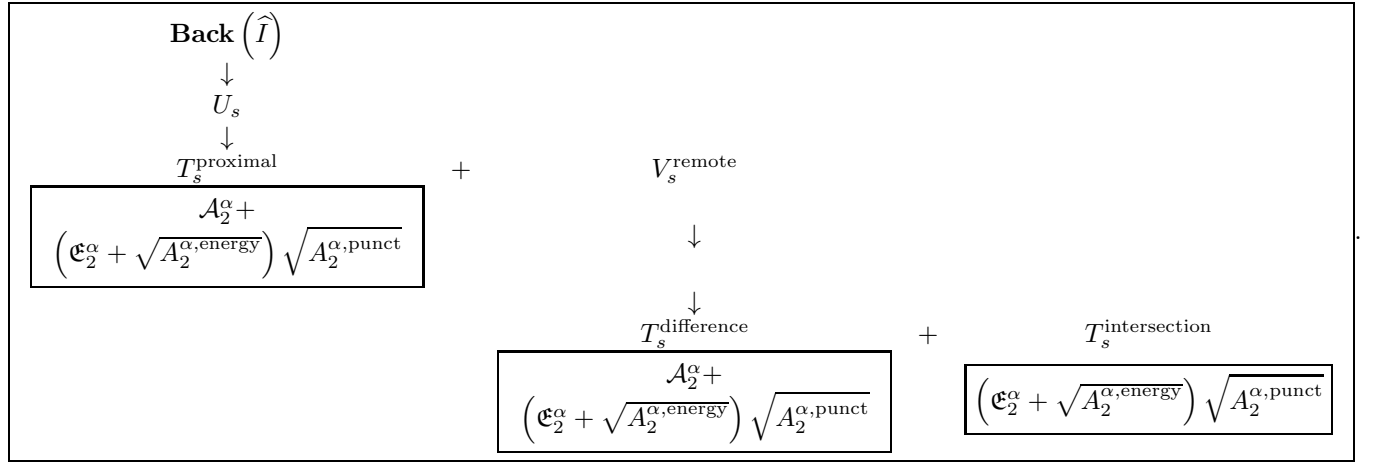
\begin{equation*}
\fbox{$%
\begin{array}{ccccc}
\mathbf{Back}\left( \widehat{I}\right) &  &  &  &  \\ 
\downarrow &  &  &  &  \\ 
U_{s} &  &  &  &  \\ 
\downarrow &  &  &  &  \\ 
T_{s}^{\limfunc{proximal}} & + & V_{s}^{\limfunc{remote}} &  &  \\ 
\fbox{$%
\begin{array}{c}
\mathcal{A}_{2}^{\alpha }+ \\ 
\left( \mathfrak{E}_{2}^{\alpha }+\sqrt{A_{2}^{\alpha ,\limfunc{energy}}}%
\right) \sqrt{A_{2}^{\alpha ,\limfunc{punct}}}%
\end{array}%
$} &  & \downarrow &  &  \\ 
&  & \downarrow &  &  \\ 
&  & T_{s}^{\limfunc{difference}} & + & T_{s}^{\limfunc{intersection}} \\ 
&  & \fbox{$%
\begin{array}{c}
\mathcal{A}_{2}^{\alpha }+ \\ 
\left( \mathfrak{E}_{2}^{\alpha }+\sqrt{A_{2}^{\alpha ,\limfunc{energy}}}%
\right) \sqrt{A_{2}^{\alpha ,\limfunc{punct}}}%
\end{array}%
$} &  & \fbox{$\left( \mathfrak{E}_{2}^{\alpha }+\sqrt{A_{2}^{\alpha ,%
\limfunc{energy}}}\right) \sqrt{A_{2}^{\alpha ,\limfunc{punct}}}$}%
\end{array}%
$}.
\end{equation*}%
Using (\ref{tent consequence}) we see that the integral on the right hand
side of (\ref{e.t2 n'}) is 
\begin{equation}
\int_{\widehat{I}}t^{2}d\overline{\mu }=\sum_{F\in \mathcal{F}}\sum_{M\in 
\mathcal{W}\left( F\right) :\ M\subset I}\lVert \mathsf{Q}_{F,M}^{\omega ,%
\mathbf{b}^{\ast }}x\rVert _{L^{2}\left( \omega \right) }^{\spadesuit 2}\,.
\label{mu I hat}
\end{equation}%
where $\mathsf{Q}_{F,M}^{\omega ,\mathbf{b}^{\ast }}$ was defined in (\ref%
{def F,K}).

We now compute using (\ref{tent consequence}) again that 
\begin{eqnarray}
\mathbb{Q}^{\alpha }\left( t\mathbf{1}_{\widehat{I}}\overline{\mu }\right)
\left( y\right) &=&\int_{\widehat{I}}\frac{t^{2}}{\left( t^{2}+\left\vert
x-y\right\vert ^{2}\right) ^{\frac{2-\alpha }{2}}}d\overline{\mu }\left(
x,t\right)  \label{PI hat} \\
&\approx &\sum_{F\in \mathcal{F}}\sum_{M\in \mathcal{W}\left( F\right) :\
M\subset I}\frac{\lVert \mathsf{Q}_{F,M}^{\omega ,\mathbf{b}^{\ast }}x\rVert
_{L^{2}\left( \omega \right) }^{\spadesuit 2}}{\left( \left\vert
M\right\vert +\left\vert y-c_{M}\right\vert \right) ^{2-\alpha }},  \notag
\end{eqnarray}%
and then expand the square and integrate to obtain that the term $\mathbf{%
Back}\left( \widehat{I}\right) $ is 
\begin{equation*}
\sum_{\substack{ F\in \mathcal{F}  \\ M\in \mathcal{W}\left( F\right)  \\ %
M\subset I}}\sum_{\substack{ F^{\prime }\in \mathcal{F}:  \\ M^{\prime }\in 
\mathcal{W}\left( F^{\prime }\right)  \\ M^{\prime }\subset I}}\int_{\mathbb{%
R}}\frac{\left\Vert \mathsf{Q}_{F,M}^{\omega ,\mathbf{b}^{\ast
}}x\right\Vert _{L^{2}\left( \omega \right) }^{\spadesuit 2}}{\left(
\left\vert M\right\vert +\left\vert y-c_{M}\right\vert \right) ^{2-\alpha }}%
\frac{\left\Vert \mathsf{Q}_{F^{\prime },M^{\prime }}^{\omega ,\mathbf{b}%
^{\ast }}x\right\Vert _{L^{2}\left( \omega \right) }^{\spadesuit 2}}{\left(
\left\vert M^{\prime }\right\vert +\left\vert y-c_{M^{\prime }}\right\vert
\right) ^{2-\alpha }}d\sigma \left( y\right) .
\end{equation*}

By symmetry we may assume that $\ell \left( M^{\prime }\right) \leq \ell
\left( M\right) $. We fix a nonnegative integer $s$, and consider those
intervals $M$ and $M^{\prime }$ with $\ell \left( M^{\prime }\right)
=2^{-s}\ell \left( M\right) $. For fixed $s$ we will control the expression 
\begin{eqnarray}
U_{s} &\equiv &\sum_{\substack{ F,F^{\prime }\in \mathcal{F}}}\sum 
_{\substack{ M\in \mathcal{W}\left( F\right) ,\ M^{\prime }\in \mathcal{W}%
\left( F^{\prime }\right)  \\ M,M^{\prime }\subset I,\ \ell \left( M^{\prime
}\right) =2^{-s}\ell \left( M\right) }}  \label{def Us} \\
&&\times \int_{\mathbb{R}}\frac{\left\Vert \mathsf{Q}_{F,M}^{\omega ,\mathbf{%
b}^{\ast }}x\right\Vert _{L^{2}\left( \omega \right) }^{\spadesuit 2}}{%
\left( \left\vert M\right\vert +\left\vert y-c_{M}\right\vert \right)
^{2-\alpha }}\frac{\left\Vert \mathsf{Q}_{F^{\prime },M^{\prime }}^{\omega ,%
\mathbf{b}^{\ast }}x\right\Vert _{L^{2}\left( \omega \right) }^{\spadesuit 2}%
}{\left( \left\vert M^{\prime }\right\vert +\left\vert y-c_{M^{\prime
}}\right\vert \right) ^{2-\alpha }}d\sigma \left( y\right) ,  \notag
\end{eqnarray}%
by proving that%
\begin{equation}
U_{s}\lesssim 2^{-\delta s}\left\{ \mathcal{A}_{2}^{\alpha }+\left( 
\mathfrak{E}_{2}^{\alpha }+\sqrt{A_{2}^{\alpha ,\limfunc{energy}}}\right) 
\sqrt{A_{2}^{\alpha ,\limfunc{punct}}}\right\} \int_{\widehat{I}}t^{2}d%
\overline{\mu },\ \ \ \ \ \text{where }\delta =\frac{1}{2}.  \label{Us bound}
\end{equation}%
With this accomplished, we can sum in $s\geq 0$ to control the term $\mathbf{%
Back}\left( \widehat{I}\right) $. We now decompose $U_{s}=T_{s}^{\limfunc{%
proximal}}+T_{s}^{\limfunc{difference}}+T_{s}^{\limfunc{intersection}}$ into
three pieces.

Our first decomposition is to write%
\begin{equation}
U_{s}=T_{s}^{\limfunc{proximal}}+V_{s}^{\limfunc{remote}}\ ,
\label{initial decomp}
\end{equation}%
where in the `proximal' term $T_{s}^{\limfunc{proximal}}$ we restrict the
summation over pairs of intervals $M,M^{\prime }$ to those satisfying $%
d\left( c_{M},c_{M^{\prime }}\right) <2^{s\delta }\ell \left( M\right) $;
while in the `remote' term $V_{s}^{\limfunc{remote}}$ we restrict the
summation over pairs of intervals $M,M^{\prime }$ to those satisfying the
opposite inequality $d\left( c_{M},c_{M^{\prime }}\right) \geq 2^{s\delta
}\ell \left( M\right) $. Then we further decompose 
\begin{equation*}
V_{s}^{\limfunc{remote}}=T_{s}^{\limfunc{difference}}+T_{s}^{\limfunc{%
intersection}},
\end{equation*}%
where in the `difference' term $T_{s}^{\limfunc{difference}}$ we restrict
integration in $y$ to the difference $\mathbb{R}\setminus B\left(
M,M^{\prime }\right) $ of $\mathbb{R}$ and 
\begin{equation}
B\left( M,M^{\prime }\right) \equiv B\left( c_{M},\frac{1}{2}d\left(
c_{M},c_{M^{\prime }}\right) \right) ,  \label{def BMM'}
\end{equation}%
the ball centered at $c_{M}$ with radius $\frac{1}{2}d\left(
c_{M},c_{M^{\prime }}\right) $; while in the `intersection' term $T_{s}^{%
\limfunc{intersection}}$ we restrict integration in $y$ to the intersection $%
\mathbb{R}\cap B\left( M,M^{\prime }\right) $ of $\mathbb{R}$ with the ball $%
B\left( M,M^{\prime }\right) $; i.e. 
\begin{eqnarray}
T_{s}^{\limfunc{intersection}} &\equiv &\sum_{\substack{ F,F^{\prime }\in 
\mathcal{F}}}\sum_{\substack{ M\in \mathcal{W}\left( F\right) ,\ M^{\prime
}\in \mathcal{W}\left( F^{\prime }\right)  \\ M,M^{\prime }\subset I,\ \ell
\left( M^{\prime }\right) =2^{-s}\ell \left( M\right)  \\ d\left(
c_{M},c_{M^{\prime }}\right) \geq 2^{s\left( 1+\delta \right) }\ell \left(
M^{\prime }\right) }}  \label{def Tints} \\
&&\times \int_{B\left( M,M^{\prime }\right) }\frac{\left\Vert \mathsf{Q}%
_{F,M}^{\omega ,\mathbf{b}^{\ast }}x\right\Vert _{L^{2}\left( \omega \right)
}^{\spadesuit 2}}{\left( \left\vert M\right\vert +\left\vert
y-c_{M}\right\vert \right) ^{2-\alpha }}\frac{\left\Vert \mathsf{Q}%
_{F^{\prime },M^{\prime }}^{\omega ,\mathbf{b}^{\ast }}x\right\Vert
_{L^{2}\left( \omega \right) }^{\spadesuit 2}}{\left( \left\vert M^{\prime
}\right\vert +\left\vert y-c_{M^{\prime }}\right\vert \right) ^{2-\alpha }}%
d\sigma \left( y\right) .  \notag
\end{eqnarray}%
Here is a schematic reminder of the these decompositions with the
distinguishing points of the definitions boxed:{}

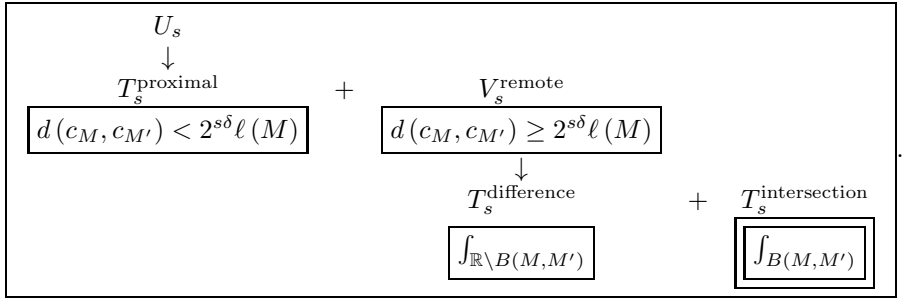
\begin{equation*}
\fbox{$%
\begin{array}{ccccc}
U_{s} &  &  &  &  \\ 
\downarrow &  &  &  &  \\ 
T_{s}^{\limfunc{proximal}} & + & V_{s}^{\limfunc{remote}} &  &  \\ 
\fbox{$d\left( c_{M},c_{M^{\prime }}\right) <2^{s\delta }\ell \left(
M\right) $} &  & \fbox{$d\left( c_{M},c_{M^{\prime }}\right) \geq 2^{s\delta
}\ell \left( M\right) $} &  &  \\ 
&  & \downarrow &  &  \\ 
&  & T_{s}^{\limfunc{difference}} & + & T_{s}^{\limfunc{intersection}} \\ 
&  & \fbox{$\int_{\mathbb{R}\setminus B\left( M,M^{\prime }\right) }$} &  & 
\fbox{$\fbox{$\int_{B\left( M,M^{\prime }\right) }$}$}%
\end{array}%
$}.
\end{equation*}

We will exploit the restriction of integration to $B\left( M,M^{\prime
}\right) $, together with the condition 
\begin{equation*}
d\left( c_{M},c_{M^{\prime }}\right) \geq 2^{s\left( 1+\delta \right) }\ell
\left( M^{\prime }\right) =2^{s\delta }\ell \left( M\right) ,
\end{equation*}%
which will then give an estimate for the term $T_{s}^{\limfunc{intersection}%
} $ using an argument dual to that used for the other terms $T_{s}^{\limfunc{%
proximal}}$ and $T_{s}^{\limfunc{difference}}$, to which we now turn.

\subsubsection{The proximal and difference terms}

We have%
\begin{align}
T_{s}^{\limfunc{proximal}}& \equiv \sum_{\substack{ F,F^{\prime }\in 
\mathcal{F}}}\sum_{\substack{ M\in \mathcal{W}\left( F\right) ,\ M^{\prime
}\in \mathcal{W}\left( F^{\prime }\right)  \\ M,M^{\prime }\subset I,\ \ell
\left( M^{\prime }\right) =2^{-s}\ell \left( M\right) \text{ and }d\left(
c_{M},c_{M^{\prime }}\right) <2^{s\delta }\ell \left( M\right) }}
\label{def Tproxs} \\
& \times \int_{\mathbb{R}}\frac{\left\Vert \mathsf{Q}_{F,M}^{\omega ,\mathbf{%
b}^{\ast }}x\right\Vert _{L^{2}\left( \omega \right) }^{\spadesuit 2}}{%
\left( \left\vert M\right\vert +\left\vert y-c_{M}\right\vert \right)
^{2-\alpha }}\frac{\left\Vert \mathsf{Q}_{F^{\prime },M^{\prime }}^{\omega ,%
\mathbf{b}^{\ast }}x\right\Vert _{L^{2}\left( \omega \right) }^{\spadesuit 2}%
}{\left( \left\vert M^{\prime }\right\vert +\left\vert y-c_{M^{\prime
}}\right\vert \right) ^{2-\alpha }}d\sigma \left( y\right)  \notag \\
& \leq \mathcal{M}_{s}^{\limfunc{proximal}}\sum_{F\in \mathcal{F}}\sum 
_{\substack{ M\in \mathcal{W}\left( F\right)  \\ M\subset I}}\lVert \mathsf{Q%
}_{F,M}^{\omega ,\mathbf{b}^{\ast }}z\rVert _{\omega }^{\spadesuit 2}=%
\mathcal{M}_{s}^{\limfunc{proximal}}\int_{\widehat{I}}t^{2}d\overline{\mu },
\notag
\end{align}%
where%
\begin{align*}
\mathcal{M}_{s}^{\limfunc{proximal}}& \equiv \sup_{F\in \mathcal{F}}\sup 
_{\substack{ M\in \mathcal{W}\left( F\right)  \\ M\subset I}}\mathcal{A}%
_{s}^{\limfunc{proximal}}\left( M\right) ; \\
\mathcal{A}_{s}^{\limfunc{proximal}}\left( M\right) & \equiv \sum_{F^{\prime
}\in \mathcal{F}}\sum_{\substack{ M^{\prime }\in \mathcal{W}\left( F^{\prime
}\right)  \\ M^{\prime }\subset I,\ \ell \left( M^{\prime }\right)
=2^{-s}\ell \left( M\right) \text{ and }d\left( c_{M},c_{M^{\prime }}\right)
<2^{s\delta }\ell \left( M\right) }}\int_{\mathbb{R}}S_{\left( M^{\prime
},M\right) }^{F^{\prime }}\left( y\right) d\sigma \left( y\right) ; \\
S_{\left( M^{\prime },M\right) }^{F^{\prime }}\left( x\right) & \equiv \frac{%
1}{\left( \left\vert M\right\vert +\left\vert y-c_{M}\right\vert \right)
^{2-\alpha }}\frac{\left\Vert \mathsf{Q}_{F^{\prime },M^{\prime }}^{\omega
}x\right\Vert _{L^{2}\left( \omega \right) }^{\spadesuit 2}}{\left(
\left\vert M^{\prime }\right\vert +\left\vert y-c_{M^{\prime }}\right\vert
\right) ^{2-\alpha }},
\end{align*}%
and similarly%
\begin{align}
T_{s}^{\limfunc{difference}}& \equiv \sum_{\substack{ F,F^{\prime }\in 
\mathcal{F}}}\sum_{\substack{ M\in \mathcal{W}\left( F\right) ,\ M^{\prime
}\in \mathcal{W}\left( F^{\prime }\right)  \\ M,M^{\prime }\subset I,\ \ell
\left( M^{\prime }\right) =2^{-s}\ell \left( M\right) \text{ and }d\left(
c_{M},c_{M^{\prime }}\right) \geq 2^{s\delta }\ell \left( M\right) }}
\label{def Tdiffs} \\
& \times \int_{\mathbb{R}\setminus B\left( M,M^{\prime }\right) }\frac{%
\left\Vert \mathsf{Q}_{F,M}^{\omega ,\mathbf{b}^{\ast }}x\right\Vert
_{L^{2}\left( \omega \right) }^{\spadesuit 2}}{\left( \left\vert
M\right\vert +\left\vert y-c_{M}\right\vert \right) ^{2-\alpha }}\frac{%
\left\Vert \mathsf{Q}_{F^{\prime },M^{\prime }}^{\omega ,\mathbf{b}^{\ast
}}x\right\Vert _{L^{2}\left( \omega \right) }^{\spadesuit 2}}{\left(
\left\vert M^{\prime }\right\vert +\left\vert y-c_{M^{\prime }}\right\vert
\right) ^{2-\alpha }}d\sigma \left( y\right)  \notag \\
& \leq \mathcal{M}_{s}^{\limfunc{difference}}\sum_{F\in \mathcal{F}}\sum 
_{\substack{ M\in \mathcal{W}\left( F\right)  \\ M\subset I}}\lVert \mathsf{Q%
}_{F,M}^{\omega ,\mathbf{b}^{\ast }}z\rVert _{\omega }^{\spadesuit 2}=%
\mathcal{M}_{s}^{\limfunc{difference}}\int_{\widehat{I}}t^{2}d\overline{\mu }%
;  \notag
\end{align}%
where%
\begin{eqnarray*}
\mathcal{M}_{s}^{\limfunc{difference}} &\equiv &\sup_{F\in \mathcal{F}}\sup 
_{\substack{ _{\substack{ M\in \mathcal{W}\left( F\right) }}  \\ M\subset I}}%
\mathcal{A}_{s}^{\func{difference}}\left( M\right) ; \\
\mathcal{A}_{s}^{\limfunc{difference}}\left( M\right) &\equiv
&\sum_{F^{\prime }\in \mathcal{F}}\sum_{\substack{ M^{\prime }\in \mathcal{W}%
\left( F^{\prime }\right)  \\ M^{\prime }\subset I,\ \ell \left( M^{\prime
}\right) =2^{-s}\ell \left( M\right) \text{ and }d\left( c_{M},c_{M^{\prime
}}\right) \geq 2^{s\delta }\ell \left( M\right) }}\int_{\mathbb{R}\setminus
B\left( M,M^{\prime }\right) }S_{\left( M^{\prime },M\right) }^{F^{\prime
}}\left( y\right) d\sigma \left( y\right) .
\end{eqnarray*}%
The restriction of integration in $\mathcal{A}_{s}^{\limfunc{difference}}$
to $\mathbb{R}\setminus B\left( M,M^{\prime }\right) $ will be used to
establish (\ref{vanishing close}) below.

\begin{notation}
\label{Sum *}Since the intervals $F,M,F^{\prime },M^{\prime }$ that arise in
all of the sums here satisfy 
\begin{equation*}
M\in \mathcal{W}\left( F\right) ,\ M^{\prime }\in \mathcal{W}\left(
F^{\prime }\right) \text{ and }\ell \left( M^{\prime }\right) =2^{-s}\ell
\left( M\right) \text{ and }M,M^{\prime }\subset I,
\end{equation*}%
we will often employ the notation $\overset{\ast }{\sum }$ to remind the
reader that, as applicable, these four conditions are in force even when
they are\ not explictly mentioned.
\end{notation}

Now fix $M$ as in $\mathcal{M}_{s}^{\limfunc{proximal}}$ respectively $%
\mathcal{M}_{s}^{\limfunc{difference}}$, and decompose the sum over $%
M^{\prime }$ in $\mathcal{A}_{s}^{\limfunc{proximal}}\left( M\right) $
respectively $\mathcal{A}_{s}^{\limfunc{difference}}\left( M\right) $ by%
\begin{eqnarray*}
&&\mathcal{A}_{s}^{\limfunc{proximal}}\left( M\right) =\sum_{F^{\prime }\in 
\mathcal{F}}\sum_{\substack{ M^{\prime }\in \mathcal{W}\left( F^{\prime
}\right)  \\ M^{\prime }\subset I,\ \ell \left( M^{\prime }\right)
=2^{-s}\ell \left( M\right) \text{ and }d\left( c_{M},c_{M^{\prime }}\right)
<2^{s\delta }\ell \left( M\right) }}\int_{\mathbb{R}}S_{\left( M^{\prime
},M\right) }^{F^{\prime }}\left( y\right) d\sigma \left( y\right) \\
&=&\sum_{F^{\prime }\in \mathcal{F}}\overset{\ast }{\sum_{\substack{ %
c_{M^{\prime }}\in 2M  \\ d\left( c_{M},c_{M^{\prime }}\right) <2^{s\delta
}\ell \left( M\right) }}}\int_{\mathbb{R}}S_{\left( M^{\prime },M\right)
}^{F^{\prime }}\left( y\right) d\sigma \left( y\right) +\sum_{F^{\prime }\in 
\mathcal{F}}\sum_{\ell =1}^{\infty }\overset{\ast }{\sum_{\substack{ %
c_{M^{\prime }}\in 2^{\ell +1}M\setminus 2^{\ell }M  \\ d\left(
c_{M},c_{M^{\prime }}\right) <2^{s\delta }\ell \left( M\right) }}}\int_{%
\mathbb{R}}S_{\left( M^{\prime },M\right) }^{F^{\prime }}\left( y\right)
d\sigma \left( y\right) \\
&\equiv &\sum_{\ell =0}^{\infty }\mathcal{A}_{s}^{\limfunc{proximal},\ell
}\left( M\right) ,
\end{eqnarray*}%
respectively%
\begin{eqnarray*}
&&\mathcal{A}_{s}^{\limfunc{difference}}\left( M\right) =\sum_{F^{\prime
}\in \mathcal{F}}\sum_{\substack{ M^{\prime }\in \mathcal{W}\left( F^{\prime
}\right)  \\ M^{\prime }\subset I,\ \ell \left( M^{\prime }\right)
=2^{-s}\ell \left( M\right) \text{ and }d\left( c_{M},c_{M^{\prime }}\right)
\geq 2^{s\delta }\ell \left( M\right) }}\int_{\mathbb{R}\setminus B\left(
M,M^{\prime }\right) }S_{\left( M^{\prime },M\right) }^{F^{\prime }}\left(
y\right) d\sigma \left( y\right) \\
&=&\sum_{F^{\prime }\in \mathcal{F}}\overset{\ast }{\sum_{\substack{ %
c_{M^{\prime }}\in 2M  \\ d\left( c_{M},c_{M^{\prime }}\right) \geq
2^{s\delta }\ell \left( M\right) }}}\int_{\mathbb{R}\setminus B\left(
M,M^{\prime }\right) }S_{\left( M^{\prime },M\right) }^{F^{\prime }}\left(
y\right) d\sigma \left( y\right) \\
&&+\sum_{\ell =1}^{\infty }\sum_{F^{\prime }\in \mathcal{F}}\overset{\ast }{%
\sum_{\substack{ c_{M^{\prime }}\in 2^{\ell +1}M\setminus 2^{\ell }M  \\ %
d\left( c_{M},c_{M^{\prime }}\right) \geq 2^{s\delta }\ell \left( M\right) }}%
}\int_{\mathbb{R}\setminus B\left( M,M^{\prime }\right) }S_{\left( M^{\prime
},M\right) }^{F^{\prime }}\left( y\right) d\sigma \left( y\right) \\
&\equiv &\sum_{\ell =0}^{\infty }\mathcal{A}_{s}^{\limfunc{difference},\ell
}\left( M\right) .
\end{eqnarray*}%
Let $m=2$ so that 
\begin{equation}
2^{-m}\leq \frac{1}{3}.  \label{smallest m}
\end{equation}%
Now decompose the integrals over $\mathbb{R}$ in $\mathcal{A}_{s}^{\limfunc{%
proximal},\ell }\left( M\right) $ by%
\begin{eqnarray*}
\mathcal{A}_{s}^{\limfunc{proximal},0}\left( M\right) &=&\sum_{F^{\prime
}\in \mathcal{F}}\overset{\ast }{\sum_{\substack{ c_{M^{\prime }}\in 2M  \\ %
d\left( c_{M},c_{M^{\prime }}\right) <2^{s\delta }\ell \left( M\right) }}}%
\int_{\mathbb{R}\setminus 4M}S_{\left( M^{\prime },M\right) }^{F^{\prime
}}\left( y\right) d\sigma \left( y\right) \\
&&+\sum_{F^{\prime }\in \mathcal{F}}\overset{\ast }{\sum_{\substack{ %
c_{M^{\prime }}\in 2M  \\ d\left( c_{M},c_{M^{\prime }}\right) <2^{s\delta
}\ell \left( M\right) }}}\int_{4M}S_{\left( M^{\prime },M\right)
}^{F^{\prime }}\left( y\right) d\sigma \left( y\right) \\
&\equiv &\mathcal{A}_{s,far}^{\limfunc{proximal},0}\left( M\right) +\mathcal{%
A}_{s,near}^{\limfunc{proximal},0}\left( M\right) ,
\end{eqnarray*}%
and%
\begin{eqnarray*}
\mathcal{A}_{s}^{\limfunc{proximal},\ell }\left( M\right) &=&\sum_{F^{\prime
}\in \mathcal{F}}\overset{\ast }{\sum_{\substack{ c_{M^{\prime }}\in 2^{\ell
+1}M\setminus 2^{\ell }M  \\ d\left( c_{M},c_{M^{\prime }}\right)
<2^{s\delta }\ell \left( M\right) }}}\int_{\mathbb{R}\setminus 2^{\ell
+2}M}S_{\left( M^{\prime },M\right) }^{F^{\prime }}\left( y\right) d\sigma
\left( y\right) \\
&&+\sum_{F^{\prime }\in \mathcal{F}}\overset{\ast }{\sum_{\substack{ %
c_{M^{\prime }}\in 2^{\ell +1}M\setminus 2^{\ell }M  \\ d\left(
c_{M},c_{M^{\prime }}\right) <2^{s\delta }\ell \left( M\right) }}}%
\int_{2^{\ell +2}M\setminus 2^{\ell -m}M}S_{\left( M^{\prime },M\right)
}^{F^{\prime }}\left( y\right) d\sigma \left( y\right) \\
&&+\sum_{F^{\prime }\in \mathcal{F}}\overset{\ast }{\sum_{\substack{ %
c_{M^{\prime }}\in 2^{\ell +1}M\setminus 2^{\ell }M  \\ d\left(
c_{M},c_{M^{\prime }}\right) <2^{s\delta }\ell \left( M\right) }}}%
\int_{2^{\ell -m}M}S_{\left( M^{\prime },M\right) }^{F^{\prime }}\left(
y\right) d\sigma \left( y\right) \\
&\equiv &\mathcal{A}_{s,far}^{\limfunc{proximal},\ell }\left( M\right) +%
\mathcal{A}_{s,near}^{\limfunc{proximal},\ell }\left( M\right) +\mathcal{A}%
_{s,close}^{\limfunc{proximal},\ell }\left( M\right) ,\ \ \ \ \ \ell \geq 1.
\end{eqnarray*}%
Similarly we decompose the integrals over the difference 
\begin{equation*}
B^{\ast }\equiv \mathbb{R}\setminus B\left( M,M^{\prime }\right)
\end{equation*}%
in $\mathcal{A}_{s}^{\limfunc{difference},\ell }\left( M\right) $ by%
\begin{eqnarray*}
\mathcal{A}_{s}^{\limfunc{difference},0}\left( M\right) &=&\sum_{F^{\prime
}\in \mathcal{F}}\overset{\ast }{\sum_{\substack{ c_{M^{\prime }}\in 2M  \\ %
d\left( c_{M},c_{M^{\prime }}\right) \geq 2^{s\delta }\ell \left( M\right) }}%
}\int_{B^{\ast }\setminus 4M}S_{\left( M^{\prime },M\right) }^{F^{\prime
}}\left( y\right) d\sigma \left( y\right) \\
&&+\sum_{F^{\prime }\in \mathcal{F}}\overset{\ast }{\sum_{\substack{ %
c_{M^{\prime }}\in 2M  \\ d\left( c_{M},c_{M^{\prime }}\right) \geq
2^{s\delta }\ell \left( M\right) }}}\int_{B^{\ast }\cap 4M}S_{\left(
M^{\prime },M\right) }^{F^{\prime }}\left( y\right) d\sigma \left( y\right)
\\
&\equiv &\mathcal{A}_{s,far}^{\limfunc{difference},0}\left( M\right) +%
\mathcal{A}_{s,near}^{\limfunc{difference},0}\left( M\right) ,
\end{eqnarray*}%
and%
\begin{eqnarray*}
\mathcal{A}_{s}^{\limfunc{difference},\ell }\left( M\right)
&=&\sum_{F^{\prime }\in \mathcal{F}}\overset{\ast }{\sum_{\substack{ %
c_{M^{\prime }}\in 2^{\ell +1}M\setminus 2^{\ell }M  \\ d\left(
c_{M},c_{M^{\prime }}\right) \geq 2^{s\delta }\ell \left( M\right) }}}%
\int_{B^{\ast }\setminus 2^{\ell +2}M}S_{\left( M^{\prime },M\right)
}^{F^{\prime }}\left( y\right) d\sigma \left( y\right) \\
&&+\sum_{F^{\prime }\in \mathcal{F}}\overset{\ast }{\sum_{\substack{ %
c_{M^{\prime }}\in 2^{\ell +1}M\setminus 2^{\ell }M  \\ d\left(
c_{M},c_{M^{\prime }}\right) \geq 2^{s\delta }\ell \left( M\right) }}}%
\int_{B^{\ast }\cap \left( 2^{\ell +2}M\setminus 2^{\ell -m}M\right)
}S_{\left( M^{\prime },M\right) }^{F^{\prime }}\left( y\right) d\sigma
\left( y\right) \\
&&+\sum_{F^{\prime }\in \mathcal{F}}\overset{\ast }{\sum_{\substack{ %
c_{M^{\prime }}\in 2^{\ell +1}M\setminus 2^{\ell }M  \\ d\left(
c_{M},c_{M^{\prime }}\right) \geq 2^{s\delta }\ell \left( M\right) }}}%
\int_{B^{\ast }\cap 2^{\ell -m}M}S_{\left( M^{\prime },M\right) }^{F^{\prime
}}\left( y\right) d\sigma \left( y\right) \\
&\equiv &\mathcal{A}_{s,far}^{\limfunc{difference},\ell }\left( M\right) +%
\mathcal{A}_{s,near}^{\limfunc{difference},\ell }\left( M\right) +\mathcal{A}%
_{s,close}^{\limfunc{difference},\ell }\left( M\right) ,\ \ \ \ \ \ell \geq
1.
\end{eqnarray*}

We now note the important point that the close terms $\mathcal{A}_{s,close}^{%
\limfunc{proximal},\ell }\left( M\right) $ and $\mathcal{A}_{s,close}^{%
\limfunc{difference},\ell }\left( M\right) $ both $\emph{vanish}$ for $\ell
>\delta s$ because of the decomposition (\ref{initial decomp}):%
\begin{equation}
\mathcal{A}_{s,close}^{\limfunc{proximal},\ell }\left( M\right) =\mathcal{A}%
_{s,close}^{\limfunc{difference},\ell }\left( M\right) =0,\ \ \ \ \ \ell
\geq 1+\delta s.  \label{vanishing close}
\end{equation}%
Indeed, if $c_{M^{\prime }}\in 2^{\ell +1}M\setminus 2^{\ell }M$, then we
have%
\begin{equation}
\frac{1}{2}2^{\ell }\ell \left( M\right) \leq d\left( c_{M},c_{M^{\prime
}}\right) .  \label{distJJ'}
\end{equation}%
Now the summands in $\mathcal{A}_{s,close}^{\limfunc{proximal},\ell }\left(
M\right) $ satisfy $d\left( c_{M},c_{M^{\prime }}\right) <2^{\delta s}\ell
\left( M\right) $, which by (\ref{distJJ'}) is impossible if $\ell \geq
1+\delta s$ - indeed, if $\ell \geq 1+\delta s$, we get the contradiction%
\begin{equation*}
2^{\delta s}\ell \left( M\right) =\frac{1}{2}2^{1+\delta s}\ell \left(
M\right) \leq \frac{1}{2}2^{\ell }\ell \left( M\right) \leq d\left(
c_{M},c_{M^{\prime }}\right) <2^{\delta s}\ell \left( M\right) .
\end{equation*}%
It now follows that $\mathcal{A}_{s,close}^{\limfunc{proximal},\ell }\left(
M\right) =0$. Thus we are left to consider the term $\mathcal{A}_{s,close}^{%
\limfunc{difference},\ell }\left( M\right) $, where the integration is taken
over the set $\mathbb{R}\setminus B\left( M,M^{\prime }\right) $. But we are
also restricted in $\mathcal{A}_{s,close}^{\limfunc{difference},\ell }\left(
M\right) $ to integrating over the interval $2^{\ell -m}M$, which is
contained in $B\left( M,M^{\prime }\right) $ by (\ref{distJJ'}). Indeed, the
smallest\ ball centered at $c_{M}$ that contains $2^{\ell -m}M$ has radius $%
\frac{1}{2}2^{\ell -m}\ell \left( M\right) $, which by (\ref{smallest m})
and (\ref{distJJ'}) is at most $\frac{1}{4}2^{\ell }\ell \left( M\right)
\leq \frac{1}{2}d\left( c_{M},c_{M^{\prime }}\right) $, the radius of $%
B\left( M,M^{\prime }\right) $. Thus the range of integration in the term $%
\mathcal{A}_{s,close}^{\limfunc{difference},\ell }\left( M\right) $ is the
empty set, and so $\mathcal{A}_{s,close}^{\limfunc{difference},\ell }\left(
M\right) =0$ as well as $\mathcal{A}_{s,close}^{\limfunc{proximal},\ell
}\left( M\right) =0$. This proves (\ref{vanishing close}).

From now on we treat $T_{s}^{\limfunc{proximal}}$ and $T_{s}^{\limfunc{%
difference}}$ in the same way since the terms $\mathcal{A}_{s,close}^{%
\limfunc{proximal},\ell }\left( M\right) $ and $\mathcal{A}_{s,close}^{%
\limfunc{difference},\ell }\left( M\right) $ both vanish for $\ell \geq
1+\delta s$. Thus we will suppress the superscripts $\limfunc{proximal}$ and 
$\limfunc{difference}$ in the $far$, $near$ and $close$ decomposition of $%
\mathcal{A}_{s}^{\limfunc{proximal},\ell }\left( M\right) $ and $\mathcal{A}%
_{s}^{\limfunc{difference},\ell }\left( M\right) $, and we will also
suppress the conditions $d\left( c_{M},c_{M^{\prime }}\right) <2^{s\delta
}\ell \left( M\right) $ and $d\left( c_{M},c_{M^{\prime }}\right) \geq
2^{s\delta }\ell \left( M\right) $ in the proximal and difference terms
since they no longer play a role. Using the pairwise disjointedness of the
shifted coronas $\mathcal{C}_{F}^{\mathcal{G},\limfunc{shift}}$, we have 
\begin{equation*}
\sum_{F^{\prime }\in \mathcal{F}}\left\Vert \mathsf{Q}_{F^{\prime
},A}^{\omega ,\mathbf{b}^{\ast }}x\right\Vert _{L^{2}\left( \omega \right)
}^{\spadesuit 2}\lesssim \left\vert A\right\vert ^{2}\left\vert A\right\vert
_{\omega }\ ,\ \ \ \ \ \text{for any interval }A.
\end{equation*}%
Note that if $c_{M^{\prime }}\in 2M$, then $M^{\prime }\subset 3M$. Then
with 
\begin{equation}
\mathcal{W}_{M}^{s}\equiv \dbigcup\limits_{F^{\prime }\in \mathcal{F}%
}\left\{ M^{\prime }\in \mathcal{W}\left( F^{\prime }\right) :M^{\prime
}\subset 3M\text{ and }\ell \left( M^{\prime }\right) =2^{-s}\ell \left(
M\right) \right\} ,  \label{def WMs}
\end{equation}%
we have%
\begin{eqnarray*}
\mathcal{A}_{s,far}^{0}\left( M\right) &\leq &\sum_{F^{\prime }\in \mathcal{F%
}}\overset{\ast }{\sum_{c_{M^{\prime }}\in 2M}}\int_{\mathbb{R}\setminus
4M}S_{\left( M^{\prime },M\right) }^{F^{\prime }}\left( y\right) d\sigma
\left( y\right) \\
&\lesssim &\sum_{A\in \mathcal{W}_{M}^{s}}\sum_{F^{\prime }\in \mathcal{F}:\
A\in \mathcal{W}\left( F^{\prime }\right) }\int_{\mathbb{R}\setminus 4M}%
\frac{\left\Vert \mathsf{Q}_{F^{\prime },M^{\prime }}^{\omega ,\mathbf{b}%
^{\ast }}x\right\Vert _{L^{2}\left( \omega \right) }^{\spadesuit 2}}{\left(
\left\vert M\right\vert +\left\vert y-c_{M}\right\vert \right) ^{2\left(
2-\alpha \right) }}d\sigma \left( y\right) \\
&\lesssim &\sum_{A\in \mathcal{W}_{M}^{s}}\int_{\mathbb{R}\setminus 4M}\frac{%
\left\vert A\right\vert ^{2}\left\vert A\right\vert _{\omega }}{\left(
\left\vert M\right\vert +\left\vert y-c_{M}\right\vert \right) ^{2\left(
2-\alpha \right) }}d\sigma \left( y\right) \\
&=&\left( \sum_{A\in \mathcal{W}_{M}^{s}}\left\vert A\right\vert
^{2}\left\vert A\right\vert _{\omega }\right) \int_{\mathbb{R}\setminus 4M}%
\frac{1}{\left( \left\vert M\right\vert +\left\vert y-c_{M}\right\vert
\right) ^{2\left( 2-\alpha \right) }}d\sigma \left( y\right) .
\end{eqnarray*}%
Now we use the standard pigeonholing of side length of $A$ to conclude that 
\begin{eqnarray}
\sum_{A\in \mathcal{W}_{M}^{s}}\left\vert A\right\vert ^{2}\left\vert
A\right\vert _{\omega } &=&\sum_{k=s}^{\infty }\sum_{A\in \mathcal{W}%
_{M}^{s}:\ \ell \left( A\right) =2^{-k}\ell \left( M\right) }\left\vert
A\right\vert ^{2}\left\vert A\right\vert _{\omega }\leq \sum_{k=s}^{\infty
}2^{-2k}\left\vert M\right\vert ^{2}\sum_{A\in \mathcal{W}_{M}^{s}:\ \ell
\left( A\right) =2^{-k}\ell \left( M\right) }\left\vert A\right\vert
_{\omega }  \label{stan pig} \\
&\leq &\sum_{k=s}^{\infty }2^{-2k}\left\vert M\right\vert ^{2}\left\vert
3M\right\vert _{\omega }\lesssim 2^{-2s}\left\vert M\right\vert
^{2}\left\vert 3M\right\vert _{\omega },  \notag
\end{eqnarray}%
so that combining the previous two displays we have%
\begin{eqnarray*}
\mathcal{A}_{s,far}^{0}\left( M\right) &\lesssim &2^{-2s}\left\vert
M\right\vert ^{2}\left\vert 3M\right\vert _{\omega }\int_{\mathbb{R}%
\setminus 4M}\frac{1}{\left( \left\vert M\right\vert +\left\vert
y-c_{M}\right\vert \right) ^{2\left( 2-\alpha \right) }}d\sigma \left(
y\right) \\
&\leq &2^{-2s}\left\vert 4M\right\vert _{\omega }\int_{\mathbb{R}\setminus
4M}\frac{1}{\left( \left\vert M\right\vert +\left\vert y-c_{M}\right\vert
\right) ^{2\left( 1-\alpha \right) }}d\sigma \left( y\right) \\
&\approx &2^{-2s}\frac{\left\vert 4M\right\vert _{\omega }}{\left\vert
4M\right\vert ^{1-\alpha }}\int_{\mathbb{R}\setminus 4M}\left( \frac{%
\left\vert M\right\vert }{\left( \left\vert M\right\vert +\left\vert
y-c_{M}\right\vert \right) ^{2}}\right) ^{1-\alpha }d\sigma \left( y\right)
\\
&\lesssim &2^{-2s}\frac{\left\vert 4M\right\vert _{\omega }}{\left\vert
4M\right\vert ^{1-\alpha }}\mathcal{P}^{\alpha }\left( 4M,\mathbf{1}_{%
\mathbb{R}\setminus 4M}\sigma \right) \lesssim 2^{-2s}\mathcal{A}%
_{2}^{\alpha }\ .
\end{eqnarray*}

To estimate the near term $\mathcal{A}_{s,near}^{0}\left( M\right) $, we
initially keep the energy $\left\Vert \mathsf{Q}_{F^{\prime },M^{\prime
}}^{\omega ,\mathbf{b}^{\ast }}z\right\Vert _{L^{2}\left( \omega \right)
}^{2}$ and write 
\begin{eqnarray}
\mathcal{A}_{s,near}^{0}\left( M\right) &\leq &\sum_{F^{\prime }\in \mathcal{%
F}}\overset{\ast }{\sum_{c_{M^{\prime }}\in 2M}}\int_{4M}S_{\left( M^{\prime
},M\right) }^{F^{\prime }}\left( y\right) d\sigma \left( y\right)
\label{A0snear} \\
&\approx &\sum_{F^{\prime }\in \mathcal{F}}\overset{\ast }{%
\sum_{c_{M^{\prime }}\in 2M}}\int_{4M}\frac{1}{\left\vert M\right\vert
^{2-\alpha }}\frac{\left\Vert \mathsf{Q}_{F^{\prime },M^{\prime }}^{\omega ,%
\mathbf{b}^{\ast }}x\right\Vert _{L^{2}\left( \omega \right) }^{\spadesuit 2}%
}{\left( \left\vert M^{\prime }\right\vert +\left\vert y-c_{M^{\prime
}}\right\vert \right) ^{2-\alpha }}d\sigma \left( y\right)  \notag \\
&=&\sum_{F^{\prime }\in \mathcal{F}}\frac{1}{\left\vert M\right\vert
^{2-\alpha }}\overset{\ast }{\sum_{c_{M^{\prime }}\in 2M}}\left\Vert \mathsf{%
Q}_{F^{\prime },M^{\prime }}^{\omega ,\mathbf{b}^{\ast }}x\right\Vert
_{L^{2}\left( \omega \right) }^{\spadesuit 2}\int_{4M}\frac{1}{\left(
\left\vert M^{\prime }\right\vert +\left\vert y-c_{M^{\prime }}\right\vert
\right) ^{2-\alpha }}d\sigma \left( y\right)  \notag \\
&=&\sum_{F^{\prime }\in \mathcal{F}}\frac{1}{\left\vert M\right\vert
^{2-\alpha }}\overset{\ast }{\sum_{c_{M^{\prime }}\in 2M}}\left\Vert \mathsf{%
Q}_{F^{\prime },M^{\prime }}^{\omega ,\mathbf{b}^{\ast }}x\right\Vert
_{L^{2}\left( \omega \right) }^{\spadesuit 2}\frac{\mathrm{P}^{\alpha
}\left( M^{\prime },\mathbf{1}_{4M}\sigma \right) }{\left\vert M^{\prime
}\right\vert }.  \notag
\end{eqnarray}%
In order to estimate the final sum above, we must invoke the `prepare to
puncture' argument above, as we will want to derive geometric decay from $%
\left\Vert \mathsf{Q}_{M^{\prime }}^{\omega ,\mathbf{b}^{\ast }}x\right\Vert
_{L^{2}\left( \omega \right) }^{\spadesuit 2}$ by dominating it by the
`nonenergy' term $\left\vert M^{\prime }\right\vert ^{2}\left\vert M^{\prime
}\right\vert _{\omega }$, as well as using the Muckenhoupt energy constant.
Choose an augmented interval $\widetilde{M}\in \mathcal{AD}$ satisfying $%
\dbigcup\limits_{c_{M^{\prime }}\in 2M}M^{\prime }\subset 4M\subset 
\widetilde{M}$ and $\ell \left( \widetilde{M}\right) \leq C\ell \left(
M\right) $. Define $\widetilde{\omega }=\omega -\omega \left( \left\{
p\right\} \right) \delta _{p}$ where $p$ is an atomic point in $\widetilde{M}
$ for which 
\begin{equation*}
\omega \left( \left\{ p\right\} \right) =\sup_{q\in \mathfrak{P}_{\left(
\sigma ,\omega \right) }:\ q\in \widetilde{M}}\omega \left( \left\{
q\right\} \right) .
\end{equation*}%
(If $\omega $ has no atomic point in common with $\sigma $ in $\widetilde{M}$%
, set $\widetilde{\omega }=\omega $.) Then we have $\left\vert \widetilde{M}%
\right\vert _{\widetilde{\omega }}=\omega \left( \widetilde{M},\mathfrak{P}%
_{\left( \sigma ,\omega \right) }\right) $ and%
\begin{equation*}
\frac{\left\vert \widetilde{M}\right\vert _{\widetilde{\omega }}}{\left\vert 
\widetilde{M}\right\vert ^{1-\alpha }}\frac{\left\vert \widetilde{M}%
\right\vert _{\sigma }}{\left\vert \widetilde{M}\right\vert ^{1-\alpha }}=%
\frac{\omega \left( \widetilde{M},\mathfrak{P}_{\left( \sigma ,\omega
\right) }\right) }{\left\vert \widetilde{M}\right\vert ^{1-\alpha }}\frac{%
\left\vert \widetilde{M}\right\vert _{\sigma }}{\left\vert \widetilde{M}%
\right\vert ^{1-\alpha }}\leq A_{2}^{\alpha ,\limfunc{punct}}.
\end{equation*}%
From (\ref{key obs}) and (\ref{key fact}) we also have%
\begin{equation*}
\sum_{F^{\prime }\in \mathcal{F}}\left\Vert \mathsf{Q}_{F^{\prime
},A}^{\omega ,\mathbf{b}^{\ast }}x\right\Vert _{L^{2}\left( \omega \right)
}^{\spadesuit 2}\lesssim \ell \left( A\right) ^{2}\left\vert A\right\vert _{%
\widetilde{\omega }}\ ,\ \ \ \ \ \text{for any interval }A.
\end{equation*}

Now by Cauchy-Schwarz and the augmented local estimate (\ref{shifted local})
in Lemma \ref{shifted} with $M=\widetilde{M}$ applied to the second line
below, and with $\mathcal{W}_{M}^{s}$ as in (\ref{def WMs}), and noting (\ref%
{stan pig}), the last sum in (\ref{A0snear}) is dominated by%
\begin{eqnarray}
&&\frac{1}{\left\vert M\right\vert ^{2-\alpha }}\left( \sum_{F^{\prime }\in 
\mathcal{F}}\overset{\ast }{\sum_{c\left( M^{\prime }\right) \in 2M}}%
\left\Vert \mathsf{Q}_{F^{\prime },M^{\prime }}^{\omega ,\mathbf{b}^{\ast
}}x\right\Vert _{L^{2}\left( \omega \right) }^{\spadesuit 2}\right) ^{\frac{1%
}{2}}  \label{dom by} \\
&&\ \ \ \ \ \ \ \ \ \ \times \left( \sum_{F^{\prime }\in \mathcal{F}}\overset%
{\ast }{\sum_{c_{M^{\prime }}\in 2M}}\left\Vert \mathsf{Q}_{F^{\prime
},M^{\prime }}^{\omega ,\mathbf{b}^{\ast }}x\right\Vert _{L^{2}\left( \omega
\right) }^{\spadesuit 2}\left( \frac{\mathrm{P}^{\alpha }\left( M^{\prime },%
\mathbf{1}_{4M}\sigma \right) }{\left\vert M^{\prime }\right\vert }\right)
^{2}\right) ^{\frac{1}{2}}  \notag \\
&\lesssim &\frac{1}{\left\vert M\right\vert ^{2-\alpha }}\left( \sum_{A\in 
\mathcal{W}_{M}^{s}}\left\vert A\right\vert ^{2}\left\vert A\right\vert _{%
\widetilde{\omega }}\right) ^{\frac{1}{2}}\sqrt{\left( \mathfrak{E}%
_{2}^{\alpha }\right) ^{2}+A_{2}^{\alpha ,\limfunc{energy}}}\sqrt{\left\vert 
\widetilde{M}\right\vert _{\sigma }}  \notag \\
&\lesssim &\frac{2^{-s}\left\vert M\right\vert }{\left\vert M\right\vert
^{2-\alpha }}\sqrt{\left\vert 4M\right\vert _{\widetilde{\omega }}}\sqrt{%
\left( \mathfrak{E}_{2}^{\alpha }\right) ^{2}+A_{2}^{\alpha ,\limfunc{energy}%
}}\sqrt{\left\vert \widetilde{M}\right\vert _{\sigma }}  \notag \\
&\lesssim &2^{-s}\sqrt{\left( \mathfrak{E}_{\alpha }\right)
^{2}+A_{2}^{\alpha ,\limfunc{energy}}}\sqrt{\frac{\left\vert \widetilde{M}%
\right\vert _{\widetilde{\omega }}}{\left\vert \widetilde{M}\right\vert
^{2-\alpha }}\frac{\left\vert \widetilde{M}\right\vert _{\sigma }}{%
\left\vert \widetilde{M}\right\vert ^{2-\alpha }}}\lesssim 2^{-s}\sqrt{%
\left( \mathfrak{E}_{2}^{\alpha }\right) ^{2}+A_{2}^{\alpha ,\limfunc{energy}%
}}\sqrt{A_{2}^{\alpha ,\limfunc{punct}}}\ .  \notag
\end{eqnarray}

Similarly, for $\ell \geq 1$, we can estimate the far term $\mathcal{A}%
_{s,far}^{\ell }\left( M\right) $ by the argument used for $\mathcal{A}%
_{s,far}^{0}\left( M\right) $ but applied to $2^{\ell }M$ in place of $M$.
For this need the following variant of $\mathcal{W}_{M}^{s}$ in (\ref{def
WMs}) given by%
\begin{equation}
\mathcal{W}_{M}^{s,\ell }\equiv \dbigcup\limits_{F^{\prime }\in \mathcal{F}%
}\left\{ M^{\prime }\in \mathcal{W}\left( F^{\prime }\right) :M^{\prime
}\subset 3\left( 2^{\ell }M\right) \text{ and }\ell \left( M^{\prime
}\right) =2^{-s-\ell }\ell \left( 2^{\ell }M\right) \right\} .
\label{def WMsell}
\end{equation}%
Then we have

\begin{eqnarray*}
\mathcal{A}_{s,far}^{\ell }\left( M\right) &\leq &\sum_{F^{\prime }\in 
\mathcal{F}}\overset{\ast }{\sum_{c_{M^{\prime }}\in \left( 2^{\ell
+1}M\right) \setminus \left( 2^{\ell }M\right) }}\int_{\mathbb{R}\setminus
2^{\ell +2}M}S_{\left( M^{\prime },M\right) }^{F^{\prime }}\left( y\right)
d\sigma \left( y\right) \\
&\lesssim &\sum_{A\in \mathcal{W}_{M}^{s,\ell }}\sum_{F^{\prime }\in 
\mathcal{F}:\ A\in \mathcal{W}\left( F^{\prime }\right) }\int_{\mathbb{R}%
\setminus 4\left( 2^{\ell }M\right) }\frac{\left\Vert \mathsf{Q}_{F^{\prime
},A}^{\omega ,\mathbf{b}^{\ast }}x\right\Vert _{L^{2}\left( \omega \right)
}^{\spadesuit 2}}{\left( \left\vert M\right\vert +\left\vert
y-c_{M}\right\vert \right) ^{2\left( 2-\alpha \right) }}d\sigma \left(
y\right) \\
&\lesssim &\sum_{A\in \mathcal{W}_{M}^{s,\ell }}\int_{\mathbb{R}\setminus
4\left( 2^{\ell }M\right) }\frac{\left\vert A\right\vert ^{2}\left\vert
A\right\vert _{\omega }}{\left( \left\vert M\right\vert +\left\vert
y-c_{M}\right\vert \right) ^{2\left( 2-\alpha \right) }}d\sigma \left(
y\right) \\
&=&\left( \sum_{A\in \mathcal{W}_{M}^{s,\ell }}\left\vert A\right\vert
^{2}\left\vert A\right\vert _{\omega }\right) \int_{\mathbb{R}\setminus
4\left( 2^{\ell }M\right) }\frac{1}{\left( \left\vert M\right\vert
+\left\vert y-c_{M}\right\vert \right) ^{2\left( 2-\alpha \right) }}d\sigma
\left( y\right) ,
\end{eqnarray*}%
where, just as for the sum over $A\in \mathcal{W}_{M}^{s,0}$, we have%
\begin{eqnarray}
&&\sum_{A\in \mathcal{W}_{M}^{s,\ell }}\left\vert A\right\vert
^{2}\left\vert A\right\vert _{\omega }  \label{stan pig'} \\
&=&\sum_{k=s}^{\infty }\sum_{A\in \mathcal{W}_{M}^{s,\ell }:\ \ell \left(
A\right) =2^{-k-\ell }\ell \left( 2^{\ell }M\right) }\left\vert A\right\vert
^{2}\left\vert A\right\vert _{\omega }\leq \sum_{k=s}^{\infty }2^{-2k-2\ell
}\left\vert 2^{\ell }M\right\vert ^{2}\sum_{A\in \mathcal{W}_{M}^{s,\ell }:\
\ell \left( A\right) =2^{-k-\ell }\ell \left( 2^{\ell }M\right) }\left\vert
A\right\vert _{\omega }  \notag \\
&\leq &\sum_{k=s}^{\infty }2^{-2k-2\ell }\left\vert 2^{\ell }M\right\vert
^{2}\left\vert 3\left( 2^{\ell }M\right) \right\vert _{\omega }\lesssim
2^{-2s-2\ell }\left\vert 2^{\ell }M\right\vert ^{2}\left\vert 3\left(
2^{\ell }M\right) \right\vert _{\omega }\ .  \notag
\end{eqnarray}%
Now using $\frac{\left\vert 2^{\ell }M\right\vert ^{2}}{\left( \left\vert
M\right\vert +\left\vert y-c_{M}\right\vert \right) ^{2\left( 2-\alpha
\right) }}\leq \frac{1}{\left( \left\vert 2^{\ell }M\right\vert +\left\vert
y-c_{2^{\ell }M}\right\vert \right) ^{2\left( 1-\alpha \right) }}$ for $%
y\notin 2^{\ell +2}M$, we can continue with%
\begin{eqnarray*}
&&\mathcal{A}_{s,far}^{\ell }\left( M\right) \lesssim 2^{-2s}2^{-2\ell
}\left\vert 2^{\ell +2}M\right\vert _{\omega }\int_{\mathbb{R}\setminus
2^{\ell +2}M}\frac{1}{\left( \left\vert 2^{\ell }M\right\vert +\left\vert
y-c_{2^{\ell }M}\right\vert \right) ^{2\left( 1-\alpha \right) }}d\sigma
\left( y\right) \\
&\approx &2^{-2s}2^{-2\ell }\frac{\left\vert 2^{\ell +2}M\right\vert
_{\omega }}{\left\vert 2^{\ell }M\right\vert ^{1-\alpha }}\int_{\mathbb{R}%
\setminus 2^{\ell +2}M}\left( \frac{\left\vert 2^{\ell }M\right\vert }{%
\left( \left\vert 2^{\ell }M\right\vert +\left\vert y-c_{2^{\ell
}M}\right\vert \right) ^{2}}\right) ^{1-\alpha }d\sigma \left( y\right) \\
&\lesssim &2^{-2s}2^{-2\ell }\left\{ \frac{\left\vert 2^{\ell
+2}M\right\vert _{\omega }}{\left\vert 2^{\ell }M\right\vert ^{1-\alpha }}%
\mathcal{P}^{\alpha }\left( 2^{\ell +2}M,1_{\mathbb{R}\setminus 2^{\ell
+2}M}\sigma \right) \right\} \lesssim 2^{-2s}2^{-2\ell }\mathcal{A}%
_{2}^{\alpha }\ .
\end{eqnarray*}

To estimate the near term $\mathcal{A}_{s,near}^{\ell }\left( M\right) $ we
must again invoke the \emph{`prepare to puncture'} argument. Choose an
augmented interval $\widetilde{M}\in \mathcal{AD}$ such that $%
\dbigcup\limits_{c_{M^{\prime }}\in 2^{\ell +1}M\setminus 2^{\ell
}M}M^{\prime }\subset 2^{\ell +2}M\subset \widetilde{M}$ and $\ell \left( 
\widetilde{M}\right) \leq C2^{\ell }\ell \left( M\right) $. Define $%
\widetilde{\omega }=\omega -\omega \left( \left\{ p\right\} \right) \delta
_{p}$ where $p$ is an atomic point in $\widetilde{M}$ for which 
\begin{equation*}
\omega \left( \left\{ p\right\} \right) =\sup_{q\in \mathfrak{P}_{\left(
\sigma ,\omega \right) }:\ q\in \widetilde{M}}\omega \left( \left\{
q\right\} \right) .
\end{equation*}%
(If $\omega $ has no atomic point in common with $\sigma $ in $\widetilde{M}$
set $\widetilde{\omega }=\omega $.) Then we have $\left\vert \widetilde{M}%
\right\vert _{\widetilde{\omega }}=\omega \left( \widetilde{M},\mathfrak{P}%
_{\left( \sigma ,\omega \right) }\right) $, and just as in the argument
above following (\ref{A0snear}), we have from (\ref{key obs}) and (\ref{key
fact}) that both%
\begin{equation*}
\frac{\left\vert \widetilde{M}\right\vert _{\widetilde{\omega }}}{\left\vert 
\widetilde{M}\right\vert ^{1-\alpha }}\frac{\left\vert \widetilde{M}%
\right\vert _{\sigma }}{\left\vert \widetilde{M}\right\vert ^{1-\alpha }}%
\leq A_{2}^{\alpha ,\limfunc{punct}}\text{ and }\sum_{F^{\prime }\in 
\mathcal{F}}\left\Vert \mathsf{Q}_{F^{\prime },M^{\prime }}^{\omega ,\mathbf{%
b}^{\ast }}x\right\Vert _{L^{2}\left( \omega \right) }^{\spadesuit
2}\lesssim \ell \left( M^{\prime }\right) ^{2}\left\vert M^{\prime
}\right\vert _{\widetilde{\omega }}\ .
\end{equation*}%
Thus using that $m=2$ in the definition of $A_{s,near}^{\ell }\left(
M\right) $, we see that

\begin{eqnarray*}
\mathcal{A}_{s,near}^{\ell }\left( M\right) &\leq &\sum_{F^{\prime }\in 
\mathcal{F}}\overset{\ast }{\sum_{c_{M^{\prime }}\in 2^{\ell +1}M\setminus
2^{\ell }M}}\int_{2^{\ell +2}M\setminus 2^{\ell -m}M}S_{\left( M^{\prime
},M\right) }^{F^{\prime }}\left( y\right) d\sigma \left( y\right) \\
&\approx &\sum_{F^{\prime }\in \mathcal{F}}\overset{\ast }{%
\sum_{c_{M^{\prime }}\in 2^{\ell +1}M\setminus 2^{\ell }M}}\int_{2^{\ell
+2}M\setminus 2^{\ell -m}M}\frac{1}{\left\vert 2^{\ell }M\right\vert
^{2-\alpha }}\frac{\left\Vert \mathsf{Q}_{F^{\prime },M^{\prime }}^{\omega ,%
\mathbf{b}^{\ast }}x\right\Vert _{L^{2}\left( \omega \right) }^{\spadesuit 2}%
}{\left( \left\vert M^{\prime }\right\vert +\left\vert y-c_{M^{\prime
}}\right\vert \right) ^{2-\alpha }}d\sigma \left( y\right) \\
&\lesssim &\frac{1}{\left\vert 2^{\ell }M\right\vert ^{2-\alpha }}%
\sum_{F^{\prime }\in \mathcal{F}}\overset{\ast }{\sum_{c_{M^{\prime }}\in
2^{\ell +1}M\setminus 2^{\ell }M}}\left\Vert \mathsf{Q}_{F^{\prime
},M^{\prime }}^{\omega ,\mathbf{b}^{\ast }}x\right\Vert _{L^{2}\left( \omega
\right) }^{\spadesuit 2} \\
&&\ \ \ \ \ \ \ \ \ \ \ \ \ \ \ \times \int_{2^{\ell +2}M}\frac{1}{\left(
\left\vert M^{\prime }\right\vert +\left\vert y-c_{M^{\prime }}\right\vert
\right) ^{2-\alpha }}d\sigma \left( y\right) ,
\end{eqnarray*}%
is dominated by%
\begin{eqnarray*}
&&\frac{1}{\left\vert 2^{\ell }M\right\vert ^{2-\alpha }}\sum_{F^{\prime
}\in \mathcal{F}}\overset{\ast }{\sum_{c_{M^{\prime }}\in 2^{\ell
+1}M\setminus 2^{\ell }M}}\left\Vert \mathsf{Q}_{F^{\prime },M^{\prime
}}^{\omega ,\mathbf{b}^{\ast }}x\right\Vert _{L^{2}\left( \omega \right)
}^{\spadesuit 2}\frac{\mathrm{P}^{\alpha }\left( M^{\prime },\mathbf{1}%
_{2^{\ell +2}M}\sigma \right) }{\left\vert M^{\prime }\right\vert } \\
&\leq &\frac{1}{\left\vert 2^{\ell }M\right\vert ^{2-\alpha }}\left(
\sum_{F^{\prime }\in \mathcal{F}}\overset{\ast }{\sum_{c_{M^{\prime }}\in
2^{\ell +1}M\setminus 2^{\ell }M}}\left\Vert \mathsf{Q}_{F^{\prime
},M^{\prime }}^{\omega ,\mathbf{b}^{\ast }}x\right\Vert _{L^{2}\left( \omega
\right) }^{\spadesuit 2}\right) ^{\frac{1}{2}} \\
&&\times \left( \sum_{F^{\prime }\in \mathcal{F}}\overset{\ast }{%
\sum_{c_{M^{\prime }}\in 2^{\ell +1}M\setminus 2^{\ell }M}}\left\Vert 
\mathsf{Q}_{F^{\prime },M^{\prime }}^{\omega ,\mathbf{b}^{\ast
}}x\right\Vert _{L^{2}\left( \omega \right) }^{\spadesuit 2}\left( \frac{%
\mathrm{P}^{\alpha }\left( M^{\prime },\mathbf{1}_{2^{\ell +2}M}\sigma
\right) }{\left\vert M^{\prime }\right\vert }\right) ^{2}\right) ^{\frac{1}{2%
}}.
\end{eqnarray*}

This can now be estimated as for the term $\mathcal{A}_{s,near}^{0}\left(
M\right) $, along with the augmented local estimate (\ref{shifted local}) in
Lemma \ref{shifted} with $M=\widetilde{M}$ applied to the final line above,
to get 
\begin{eqnarray*}
\mathcal{A}_{s,near}^{\ell }\left( M\right) &\lesssim &2^{-s}2^{-\ell }\frac{%
\left\vert 2^{\ell }M\right\vert }{\left\vert 2^{\ell }M\right\vert
^{2-\alpha }}\sqrt{\left\vert \widetilde{M}\right\vert _{\widetilde{\omega }}%
}\sqrt{\left( \mathfrak{E}_{2}^{\alpha }\right) ^{2}+A_{2}^{\alpha ,\limfunc{%
energy}}}\sqrt{\left\vert \widetilde{M}\right\vert _{\sigma }} \\
&\lesssim &2^{-s}2^{-\ell }\sqrt{\left( \mathfrak{E}_{2}^{\alpha }\right)
^{2}+A_{2}^{\alpha ,\limfunc{energy}}}\sqrt{\frac{\left\vert \widetilde{M}%
\right\vert _{\widetilde{\omega }}}{\left\vert \widetilde{M}\right\vert
^{1-\alpha }}\frac{\left\vert \widetilde{M}\right\vert _{\sigma }}{%
\left\vert \widetilde{M}\right\vert ^{1-\alpha }}} \\
&\lesssim &2^{-s}2^{-\ell }\sqrt{\left( \mathfrak{E}_{2}^{\alpha }\right)
^{2}+A_{2}^{\alpha ,\limfunc{energy}}}\sqrt{A_{2}^{\alpha ,\limfunc{punct}}}%
\ .
\end{eqnarray*}%
Each of the estimates for $\mathcal{A}_{s,far}^{\ell }\left( M\right) $ and $%
\mathcal{A}_{s,near}^{\ell }\left( M\right) $ is summable in both $s$ and $%
\ell $.

Now we turn to the terms $\mathcal{A}_{s,close}^{\ell }\left( M\right) $,
and recall from (\ref{vanishing close}) that $\mathcal{A}_{s,close}^{\ell
}\left( M\right) =0$ if $\ell \geq 1+\delta s$. So we now suppose that $\ell
\leq \delta s$. We have, with $m=2$ as in (\ref{smallest m}),%
\begin{eqnarray*}
&&\mathcal{A}_{s,close}^{\ell }\left( M\right) \leq \sum_{F^{\prime }\in 
\mathcal{F}}\overset{\ast }{\sum_{c_{M^{\prime }}\in 2^{\ell +1}M\setminus
2^{\ell }M}}\int_{2^{\ell -m}M}S_{\left( M^{\prime },M\right) }^{F^{\prime
}}\left( y\right) d\sigma \left( y\right)  \\
&\approx &\sum_{F^{\prime }\in \mathcal{F}}\overset{\ast }{%
\sum_{c_{M^{\prime }}\in 2^{\ell +1}M\setminus 2^{\ell }M}}\int_{2^{\ell
-m}M}\frac{1}{\left( \left\vert M\right\vert +\left\vert y-c_{M}\right\vert
\right) ^{2-\alpha }}\frac{\left\Vert \mathsf{Q}_{F^{\prime },M^{\prime
}}^{\omega ,\mathbf{b}^{\ast }}x\right\Vert _{L^{2}\left( \omega \right)
}^{\spadesuit 2}}{\left\vert 2^{\ell }M\right\vert ^{2-\alpha }}d\sigma
\left( y\right)  \\
&=&\left( \sum_{F^{\prime }\in \mathcal{F}}\overset{\ast }{%
\sum_{c_{M^{\prime }}\in 2^{\ell +1}M\setminus 2^{\ell }M}}\left\Vert 
\mathsf{Q}_{F^{\prime },M^{\prime }}^{\omega ,\mathbf{b}^{\ast
}}x\right\Vert _{L^{2}\left( \omega \right) }^{\spadesuit 2}\right) \frac{1}{%
\left\vert 2^{\ell }M\right\vert ^{2-\alpha }} \\
&&\ \ \ \ \ \ \ \ \ \ \ \ \ \ \ \ \ \ \ \ \times \int_{2^{\ell -m}M}\frac{1}{%
\left( \left\vert M\right\vert +\left\vert y-c_{M}\right\vert \right)
^{2-\alpha }}d\sigma \left( y\right) .
\end{eqnarray*}%
The argument used to prove (\ref{stan pig'}) gives the\ analogous inequality
with a hole $2^{\ell -1}M$, 
\begin{equation*}
\sum_{F^{\prime }\in \mathcal{F}}\overset{\ast }{\sum_{c_{M^{\prime }}\in
2^{\ell +1}M\setminus 2^{\ell }M}}\left\Vert \mathsf{Q}_{F^{\prime
},M^{\prime }}^{\omega ,\mathbf{b}^{\ast }}x\right\Vert _{L^{2}\left( \omega
\right) }^{\spadesuit 2}\lesssim 2^{-2s}\left\vert 2^{\ell }M\right\vert
^{2}\left\vert 2^{\ell +2}M\setminus 2^{\ell -1}M\right\vert _{\omega }\ .
\end{equation*}%
Thus we get%
\begin{eqnarray*}
&&\mathcal{A}_{s,close}^{\ell }\left( M\right)  \\
&\lesssim &2^{-2s}\left\vert 2^{\ell }M\right\vert ^{2}\left\vert 2^{\ell
+2}M\setminus 2^{\ell -1}M\right\vert _{\omega }\frac{1}{\left\vert 2^{\ell
}M\right\vert ^{2-\alpha }}\int_{2^{\ell -m}M}\frac{1}{\left( \left\vert
M\right\vert +\left\vert y-c_{M}\right\vert \right) ^{2-\alpha }}d\sigma
\left( y\right)  \\
&\lesssim &2^{-2s}\left\vert 2^{\ell }M\right\vert ^{2}\frac{\left\vert
2^{\ell +2}M\setminus 2^{\ell -1}M\right\vert _{\omega }}{\left\vert 2^{\ell
}M\right\vert ^{2-\alpha }}\frac{\left\vert 2^{\ell -m}M\right\vert _{\sigma
}}{\left\vert M\right\vert ^{2-\alpha }} \\
&\lesssim &2^{-2s}2^{\left( 2-\alpha \right) \ell }\frac{\left\vert 2^{\ell
+2}M\setminus 2^{\ell -1}M\right\vert _{\omega }}{\left\vert 2^{\ell
+2}M\right\vert ^{1-\alpha }}\frac{\left\vert 2^{\ell -m}M\right\vert
_{\sigma }}{\left\vert 2^{\ell -m}M\right\vert ^{1-\alpha }}\lesssim
2^{-2s}2^{\left( 2-\alpha \right) \ell }A_{2}^{\alpha },
\end{eqnarray*}%
provided that $m=2>1$. Note that we can use the offset Muckenhoupt constant $%
A_{2}^{\alpha }$ here since $2^{\ell +2}M\setminus 2^{\ell -1}M$ and $%
2^{\ell -m}M$ are disjoint. If $\ell \leq s$, then we have the relatively
crude estimate $\mathcal{A}_{s,close}^{\ell }\left( M\right) \lesssim
2^{-s}A_{2}^{\alpha }\ $without decay in $\ell $. But we are assuming $\ell
\leq \delta s$ here, and so we obtain a suitable estimate for $\mathcal{A}%
_{s,close}^{\ell }\left( M\right) $ provided we choose $0<\delta \leq \frac{1%
}{2-\alpha }$. Indeed, we then have%
\begin{equation*}
\sum_{l=1}^{\delta s}2^{-2s}2^{\left( 2-\alpha \right) \ell }A_{2}^{\alpha
}=2^{-2s}\left( \sum_{l=1}^{\delta s}2^{\left( 2-\alpha \right) \ell
}\right) A_{2}^{\alpha }\lesssim 2^{-2s}2^{\left( 2-\alpha \right) \delta
s}A_{2}^{\alpha }\leq 2^{-s}A_{2}^{\alpha }\ ,
\end{equation*}%
provided $\delta \leq \frac{1}{2-\alpha }$, and in particular we may take $%
\delta =\frac{1}{2}$. Altogether, the above estimates prove%
\begin{equation*}
T_{s}^{\limfunc{proximal}}+T_{s}^{\limfunc{difference}}\lesssim 2^{-s}\left( 
\mathcal{A}_{2}^{\alpha }+\sqrt{\left( \mathfrak{E}_{2}^{\alpha }\right)
^{2}+A_{2}^{\alpha ,\limfunc{energy}}}\sqrt{A_{2}^{\alpha ,\limfunc{punct}}}%
\right) ,
\end{equation*}%
which is summable in $s$.

\subsubsection{The intersection term}

Now we return to the term,%
\begin{eqnarray*}
T_{s}^{\limfunc{intersection}} &\equiv &\sum_{\substack{ F,F^{\prime }\in 
\mathcal{F}}}\sum_{\substack{ M\in \mathcal{W}\left( F\right) ,\ M^{\prime
}\in \mathcal{W}\left( F^{\prime }\right)  \\ M,M^{\prime }\subset I,\ \ell
\left( M^{\prime }\right) =2^{-s}\ell \left( M\right)  \\ d\left(
c_{M},c_{M^{\prime }}\right) \geq 2^{s\left( 1+\delta \right) }\ell \left(
M^{\prime }\right) }} \\
&&\times \int_{B\left( M,M^{\prime }\right) }\frac{\left\Vert \mathsf{Q}%
_{F,M}^{\omega ,\mathbf{b}^{\ast }}x\right\Vert _{L^{2}\left( \omega \right)
}^{\spadesuit 2}}{\left( \left\vert M\right\vert +\left\vert
y-c_{M}\right\vert \right) ^{2-\alpha }}\frac{\left\Vert \mathsf{Q}%
_{F^{\prime },M^{\prime }}^{\omega ,\mathbf{b}^{\ast }}x\right\Vert
_{L^{2}\left( \omega \right) }^{\spadesuit 2}}{\left( \left\vert M^{\prime
}\right\vert +\left\vert y-c_{M^{\prime }}\right\vert \right) ^{2-\alpha }}%
d\sigma \left( y\right) .
\end{eqnarray*}%
It will suffice to show that $T_{s}^{\limfunc{intersection}}$ satisfies the
estimate,%
\begin{eqnarray*}
T_{s}^{\limfunc{intersection}} &\lesssim &2^{-s\delta }\sqrt{\left( 
\mathfrak{E}_{2}^{\alpha }\right) ^{2}+A_{2}^{\alpha ,\limfunc{energy}}}%
\sqrt{A_{2}^{\alpha ,\limfunc{punct}}}\sum_{F^{\prime }\in \mathcal{F}%
^{\prime }}\sum_{\substack{ M^{\prime }\in \mathcal{M}_{\left( \mathbf{\rho }%
,\varepsilon \right) -\limfunc{deep}}\left( F^{\prime }\right)  \\ M^{\prime
}\subset I}}\lVert \mathsf{Q}_{F^{\prime },M^{\prime }}^{\omega ,\mathbf{b}%
^{\ast }}x\rVert _{L^{2}\left( \omega \right) }^{\spadesuit 2} \\
&=&2^{-s\delta }\sqrt{\left( \mathfrak{E}_{2}^{\alpha }\right)
^{2}+A_{2}^{\alpha ,\limfunc{energy}}}\sqrt{A_{2}^{\alpha ,\limfunc{punct}}}%
\int_{\widehat{I}}t^{2}\overline{\mu }\ .
\end{eqnarray*}%
Recalling $B\left( M,M^{\prime }\right) =B\left( c_{M},\frac{1}{2}d\left(
c_{M},c_{M^{\prime }}\right) \right) $, we can write (suppressing some
notation for clarity),%
\begin{eqnarray*}
&&T_{s}^{\limfunc{intersection}} \\
&=&\sum_{F,F^{\prime }}\sum_{\substack{ M,M^{\prime }}}\int_{B\left(
M,M^{\prime }\right) }\frac{\left\Vert \mathsf{Q}_{F,M}^{\omega ,\mathbf{b}%
^{\ast }}x\right\Vert _{L^{2}\left( \omega \right) }^{\spadesuit 2}}{\left(
\left\vert M\right\vert +\left\vert y-c_{M}\right\vert \right) ^{2-\alpha }}%
\frac{\left\Vert \mathsf{Q}_{F^{\prime },M^{\prime }}^{\omega ,\mathbf{b}%
^{\ast }}x\right\Vert _{L^{2}\left( \omega \right) }^{\spadesuit 2}}{\left(
\left\vert M^{\prime }\right\vert +\left\vert y-c_{M^{\prime }}\right\vert
\right) ^{2-\alpha }}d\sigma \left( y\right)  \\
&\approx &\sum_{F,F^{\prime }}\sum_{\substack{ M,M^{\prime }}}\left\Vert 
\mathsf{Q}_{F,M}^{\omega ,\mathbf{b}^{\ast }}x\right\Vert _{L^{2}\left(
\omega \right) }^{\spadesuit 2}\left\Vert \mathsf{Q}_{F^{\prime },M^{\prime
}}^{\omega ,\mathbf{b}^{\ast }}x\right\Vert _{L^{2}\left( \omega \right)
}^{\spadesuit 2}\frac{1}{\left\vert c_{M}-c_{M^{\prime }}\right\vert
^{2-\alpha }}\int_{B\left( M,M^{\prime }\right) }\frac{d\sigma \left(
y\right) }{\left( \left\vert M\right\vert +\left\vert y-c_{M}\right\vert
\right) ^{2-\alpha }} \\
&\leq &\sum_{F^{\prime }}\sum_{\substack{ M^{\prime }}}\left\Vert \mathsf{Q}%
_{F^{\prime },M^{\prime }}^{\omega ,\mathbf{b}^{\ast }}x\right\Vert
_{L^{2}\left( \omega \right) }^{\spadesuit 2}\sum_{F}\sum_{\substack{ M}}%
\frac{1}{\left\vert c_{M}-c_{M^{\prime }}\right\vert ^{2-\alpha }}\left\Vert 
\mathsf{Q}_{F,M}^{\omega ,\mathbf{b}^{\ast }}x\right\Vert _{L^{2}\left(
\omega \right) }^{\spadesuit 2}\int_{B\left( M,M^{\prime }\right) }\frac{%
d\sigma \left( y\right) }{\left( \left\vert M\right\vert +\left\vert
y-c_{M}\right\vert \right) ^{2-\alpha }} \\
&\equiv &\sum_{F^{\prime }}\sum_{\substack{ M^{\prime }}}\left\Vert \mathsf{Q%
}_{F^{\prime },M^{\prime }}^{\omega ,\mathbf{b}^{\ast }}x\right\Vert
_{L^{2}\left( \omega \right) }^{\spadesuit 2}S_{s}\left( M^{\prime }\right) ,
\end{eqnarray*}%
and since $\int_{B\left( M,M^{\prime }\right) }\frac{d\sigma \left( y\right) 
}{\left( \left\vert M\right\vert +\left\vert y-c_{M}\right\vert \right)
^{2-\alpha }}\approx \frac{\mathrm{P}^{\alpha }\left( M,\mathbf{1}_{B\left(
M,M^{\prime }\right) }\sigma \right) }{\left\vert M\right\vert }$, it
remains to show that for each fixed $M^{\prime }$,%
\begin{eqnarray*}
S_{s}\left( M^{\prime }\right)  &\approx &\sum_{F}\overset{\ast }{\sum
_{\substack{ M:\ d\left( c_{M},c_{M^{\prime }}\right) \geq 2^{s\left(
1+\delta \right) }\ell \left( M^{\prime }\right) }}}\frac{\left\Vert \mathsf{%
Q}_{F,M}^{\omega ,\mathbf{b}^{\ast }}x\right\Vert _{L^{2}\left( \omega
\right) }^{\spadesuit 2}}{\left\vert c_{M}-c_{M^{\prime }}\right\vert
^{2-\alpha }}\frac{\mathrm{P}^{\alpha }\left( M,\mathbf{1}_{B\left(
M,M^{\prime }\right) }\sigma \right) }{\left\vert M\right\vert } \\
&\lesssim &2^{-\delta s}\sqrt{\left( \mathfrak{E}_{2}^{\alpha }\right)
^{2}+A_{2}^{\alpha ,\limfunc{energy}}}\sqrt{A_{2}^{\alpha }}\ .
\end{eqnarray*}

We write%
\begin{eqnarray}
S_{s}\left( M^{\prime }\right)  &\approx &\sum_{k\geq s\left( 1+\delta
\right) }\frac{1}{\left( 2^{k}\left\vert M^{\prime }\right\vert \right)
^{2-\alpha }}\sum_{F}\overset{\ast }{\sum_{M:\ d\left( c_{M},c_{M^{\prime
}}\right) \approx 2^{k}\ell \left( M^{\prime }\right) }}\left\Vert \mathsf{Q}%
_{F,M}^{\omega ,\mathbf{b}^{\ast }}x\right\Vert _{L^{2}\left( \omega \right)
}^{\spadesuit 2}\frac{\mathrm{P}^{\alpha }\left( M,\mathbf{1}_{B\left(
M,M^{\prime }\right) }\sigma \right) }{\left\vert M\right\vert }
\label{def Sks} \\
&=&\sum_{k\geq s\left( 1+\delta \right) }\frac{1}{\left( 2^{k}\left\vert
M^{\prime }\right\vert \right) ^{2-\alpha }}S_{s}^{k}\left( M^{\prime
}\right) \ ;  \notag \\
S_{s}^{k}\left( M^{\prime }\right)  &\equiv &\sum_{F}\overset{\ast }{%
\sum_{M:\ d\left( c_{M},c_{M^{\prime }}\right) \approx 2^{k}\ell \left(
M^{\prime }\right) }}\left\Vert \mathsf{Q}_{F,M}^{\omega ,\mathbf{b}^{\ast
}}x\right\Vert _{L^{2}\left( \omega \right) }^{\spadesuit 2}\frac{\mathrm{P}%
^{\alpha }\left( M,\mathbf{1}_{B\left( M,M^{\prime }\right) }\sigma \right) 
}{\left\vert M\right\vert },  \notag
\end{eqnarray}%
where by $d\left( c_{M},c_{M^{\prime }}\right) \approx 2^{k}\ell \left(
M^{\prime }\right) $ we mean $2^{k}\ell \left( M^{\prime }\right) \leq
d\left( c_{M},c_{M^{\prime }}\right) \leq 2^{k+1}\ell \left( M^{\prime
}\right) $. Moreover, if $d\left( c_{M},c_{M^{\prime }}\right) \approx
2^{k}\ell \left( M^{\prime }\right) $, then from the fact that the radius of 
$B\left( M,M^{\prime }\right) $ is $\frac{1}{2}d\left( c_{M},c_{M^{\prime
}}\right) $, we obtain 
\begin{equation*}
B\left( M,M^{\prime }\right) \subset C_{0}2^{k}M^{\prime },
\end{equation*}%
where $C_{0}$ is a positive constant ($C_{0}=6$ works).

For fixed $k\geq s\left( 1+\delta \right) $, we invoke yet again the \emph{%
`prepare to puncture'} argument. Choose an augmented interval $\widetilde{%
M^{\prime }}\in \mathcal{AD}$ such that $C_{0}2^{k}M\subset \widetilde{%
M^{\prime }}$ and $\ell \left( \widetilde{M^{\prime }}\right) \leq
C2^{k}\ell \left( M^{\prime }\right) $. Define $\widetilde{\omega }=\omega
-\omega \left( \left\{ p\right\} \right) \delta _{p}$ where $p$ is an atomic
point in $\widetilde{M^{\prime }}$ for which 
\begin{equation*}
\omega \left( \left\{ p\right\} \right) =\sup_{q\in \mathfrak{P}_{\left(
\sigma ,\omega \right) }:\ q\in \widetilde{M^{\prime }}}\omega \left(
\left\{ q\right\} \right) .
\end{equation*}%
(If $\omega $ has no atomic point in common with $\sigma $ in $\widetilde{%
M^{\prime }}$, set $\widetilde{\omega }=\omega $.) Then we have $\left\vert 
\widetilde{M^{\prime }}\right\vert _{\widetilde{\omega }}=\omega \left( 
\widetilde{M^{\prime }},\mathfrak{P}_{\left( \sigma ,\omega \right) }\right) 
$ and so from (\ref{key obs}) and (\ref{key fact}),%
\begin{equation*}
\frac{\left\vert \widetilde{M^{\prime }}\right\vert _{\widetilde{\omega }}}{%
\left\vert \widetilde{M^{\prime }}\right\vert ^{1-\alpha }}\frac{\left\vert 
\widetilde{M^{\prime }}\right\vert _{\sigma }}{\left\vert \widetilde{%
M^{\prime }}\right\vert ^{1-\alpha }}\leq A_{2}^{\alpha ,\limfunc{punct}}%
\text{ and }\sum_{F\in \mathcal{F}}\left\Vert \mathsf{Q}_{F,A}^{\omega ,%
\mathbf{b}^{\ast }}x\right\Vert _{L^{2}\left( \omega \right) }^{\spadesuit
2}\lesssim \ell \left( A\right) ^{2}\left\vert A\right\vert _{\widetilde{%
\omega }}\ \text{for any interval }A.
\end{equation*}

Now we are ready to apply Cauchy-Schwarz and the augmented local estimate (%
\ref{shifted local}) in Lemma \ref{shifted} with $M=\widetilde{M^{\prime }}\ 
$to the second line below, and to apply the argument in (\ref{stan pig'}) to
the first line below, to get the following estimate for $S_{s}^{k}\left(
M^{\prime }\right) $ defined in (\ref{def Sks}) above:%
\begin{eqnarray*}
S_{s}^{k}\left( M^{\prime }\right)  &\leq &\left( \sum_{F}\sum_{M:\ d\left(
c_{M},c_{M^{\prime }}\right) \approx 2^{k}\ell \left( M^{\prime }\right)
}\left\Vert \mathsf{Q}_{F,M}^{\omega ,\mathbf{b}^{\ast }}x\right\Vert
_{L^{2}\left( \omega \right) }^{\spadesuit 2}\right) ^{\frac{1}{2}} \\
&&\times \left( \sum_{F}\sum_{M:\ d\left( c_{M},c_{M^{\prime }}\right)
\approx 2^{k}\ell \left( M^{\prime }\right) }\left\Vert \mathsf{Q}%
_{F,M}^{\omega ,\mathbf{b}^{\ast }}x\right\Vert _{L^{2}\left( \omega \right)
}^{\spadesuit 2}\left( \frac{\mathrm{P}^{\alpha }\left( M,\mathbf{1}%
_{B\left( M,M^{\prime }\right) }\sigma \right) }{\left\vert M\right\vert }%
\right) ^{2}\right) ^{\frac{1}{2}} \\
&\lesssim &\left( 2^{2s}\left\vert \widetilde{M^{\prime }}\right\vert
^{2}\left\vert \widetilde{M^{\prime }}\right\vert _{\widetilde{\omega }%
}\right) ^{\frac{1}{2}}\left( \left[ \left( \mathfrak{E}_{2}^{\alpha
}\right) ^{2}+A_{2}^{\alpha ,\limfunc{energy}}\right] \left\vert \widetilde{%
M^{\prime }}\right\vert _{\sigma }\right) ^{\frac{1}{2}} \\
&\lesssim &\sqrt{\left( \mathfrak{E}_{2}^{\alpha }\right) ^{2}+A_{2}^{\alpha
,\limfunc{energy}}}2^{s}\left\vert \widetilde{M^{\prime }}\right\vert \sqrt{%
\left\vert \widetilde{M^{\prime }}\right\vert _{\widetilde{\omega }}}\sqrt{%
\left\vert \widetilde{M^{\prime }}\right\vert _{\sigma }} \\
&\lesssim &\sqrt{\left( \mathfrak{E}_{2}^{\alpha }\right) ^{2}+A_{2}^{\alpha
,\limfunc{energy}}}\sqrt{A_{2}^{\alpha ,\limfunc{punct}}}2^{s}\left\vert 
\widetilde{M^{\prime }}\right\vert \left\vert \widetilde{M^{\prime }}%
\right\vert ^{1-\alpha } \\
&\approx &\sqrt{\left( \mathfrak{E}_{2}^{\alpha }\right) ^{2}+A_{2}^{\alpha ,%
\limfunc{energy}}}\sqrt{A_{2}^{\alpha ,\limfunc{punct}}}2^{s}2^{k\left(
1-\alpha \right) }\left\vert M^{\prime }\right\vert ^{2-\alpha },\ \ \ \ \ \ 
\text{since }\ell \left( \widetilde{M^{\prime }}\right) \approx 2^{k}\ell
\left( M^{\prime }\right) .
\end{eqnarray*}

Altogether then we have%
\begin{eqnarray*}
S_{s}\left( M^{\prime }\right)  &=&\sum_{k\geq \left( 1+\delta \right) s}%
\frac{1}{\left( 2^{k}\left\vert M^{\prime }\right\vert \right) ^{2-\alpha }}%
S_{s}^{k}\left( M^{\prime }\right)  \\
&\lesssim &\sqrt{\left( \mathfrak{E}_{2}^{\alpha }\right) ^{2}+A_{2}^{\alpha
,\limfunc{energy}}}\sqrt{A_{2}^{\alpha ,\limfunc{punct}}}\sum_{k\geq \left(
1+\delta \right) s}\frac{1}{\left( 2^{k}\left\vert M^{\prime }\right\vert
\right) ^{2-\alpha }}2^{s}2^{k\left( 1-\alpha \right) }\left\vert M^{\prime
}\right\vert ^{2-\alpha } \\
&=&\sqrt{\left( \mathfrak{E}_{2}^{\alpha }\right) ^{2}+A_{2}^{\alpha ,%
\limfunc{energy}}}\sqrt{A_{2}^{\alpha ,\limfunc{punct}}}\sum_{k\geq \left(
1+\delta \right) s}2^{s-k}\lesssim 2^{-\delta s}\sqrt{\left( \mathfrak{E}%
_{2}^{\alpha }\right) ^{2}+A_{2}^{\alpha ,\limfunc{energy}}}\sqrt{%
A_{2}^{\alpha ,\limfunc{punct}}},
\end{eqnarray*}%
which is summable in $s$. This completes the proof of (\ref{Us bound}), and
hence of the estimate for $\mathbf{Back}\left( \widehat{I}\right) $ in (\ref%
{e.t2 n'}).

The proof of Proposition \ref{func ener control} is now complete.

\section{Appendix C: Errata for the Revista paper}

The current paper adapts the arguments of our 2016 Revista paper \cite%
{SaShUr7} whenever possible (which in turn adapted arguments from many
earlier papers of various authors). To aid the reader in consulting \cite%
{SaShUr7}, we include here a list of up-to-date \emph{errata} for the
Revista paper \cite{SaShUr7}.

\bigskip

{\Large \#1:} Lemma 3.1 on page 90 and its proof should be replaced with the
following taken from \cite[\texttt{arXiv:1505.07816}v3]{SaShUr6}.

\textbf{Lemma 3.1} Given $\mathbf{r}\geq 3$, $\mathbf{\tau }\geq 1$ and $%
\frac{1}{\mathbf{r}}<\varepsilon <1-\frac{1}{\mathbf{r}}$, we have 
\begin{equation*}
\mathcal{D}_{\left( \mathbf{r}-1,\delta \right) -\limfunc{good}}\subset 
\mathcal{D}_{\left( \mathbf{r},\varepsilon \right) -\limfunc{good}}^{\mathbf{%
\tau }}\ ,
\end{equation*}%
provided%
\begin{equation*}
0<\delta \leq \frac{\mathbf{r}\varepsilon -1}{\mathbf{r}+\mathbf{\tau }}.
\end{equation*}

\begin{proof}
Suppose that $I\in \mathcal{D}_{\left( \mathbf{r}-1,\delta \right) -\limfunc{%
good}}$ where $\delta $ is as above. If $J$ is a child of $I$, then $J\in 
\mathcal{D}_{\left( \mathbf{r},\delta \right) -\limfunc{good}}$, and since $%
\delta <\varepsilon $, we also have $\mathcal{D}_{\left( \mathbf{r},\delta
\right) -\limfunc{good}}\subset \mathcal{D}_{\left( \mathbf{r},\varepsilon
\right) -\limfunc{good}}$. It remains to show that $\pi _{\mathcal{D}%
}^{\left( m\right) }I\in \mathcal{D}_{\left( \mathbf{r},\varepsilon \right) -%
\limfunc{good}}$ for $0\leq m\leq \mathbf{\tau }$. For this it suffices to
show that if $K\in \mathcal{D}$ satisfies\thinspace $\pi _{\mathcal{D}%
}^{\left( m\right) }I\subset K$ and $\ell \left( \pi _{\mathcal{D}}^{\left(
m\right) }I\right) \leq 2^{-\mathbf{r}}\ell \left( K\right) $ , then%
\begin{equation*}
\frac{1}{2}\left( \frac{\ell \left( \pi _{\mathcal{D}}^{\left( m\right)
}I\right) }{\ell \left( K\right) }\right) ^{\varepsilon }\ell \left(
K\right) \leq \limfunc{dist}\left( \pi _{\mathcal{D}}^{\left( m\right)
}I,K^{c}\right) .
\end{equation*}%
Now $\ell \left( I\right) =2^{-m}\ell \left( \pi _{\mathcal{D}}^{\left(
m\right) }I\right) \leq 2^{-\left( m+\mathbf{r}\right) }\ell \left( K\right)
\leq 2^{-\left( \mathbf{r}-1\right) }\ell \left( K\right) $ and $I\in 
\mathcal{D}_{\left( \mathbf{r}-1,\delta \right) -\limfunc{good}}$ imply that%
\begin{equation*}
\frac{1}{2}\left( \frac{\ell \left( I\right) }{\ell \left( K\right) }\right)
^{\delta }\ell \left( K\right) \leq \limfunc{dist}\left( I,K^{c}\right) ,
\end{equation*}%
and since the triangle inequality gives%
\begin{equation*}
\limfunc{dist}\left( I,K^{c}\right) \leq \limfunc{dist}\left( \pi _{\mathcal{%
D}}^{\left( m\right) }I,K^{c}\right) +2^{m}\ell \left( I\right) ,
\end{equation*}%
we see that it suffices to show%
\begin{equation*}
\frac{1}{2}\left( \frac{\ell \left( \pi _{\mathcal{D}}^{\left( m\right)
}I\right) }{\ell \left( K\right) }\right) ^{\varepsilon }\ell \left(
K\right) +2^{m}\ell \left( I\right) \leq \frac{1}{2}\left( \frac{\ell \left(
I\right) }{\ell \left( K\right) }\right) ^{\delta }\ell \left( K\right) ,\ \
\ \ \ 0\leq m\leq \mathbf{\tau }.
\end{equation*}%
This is equivalent to successively,%
\begin{eqnarray*}
\frac{1}{2}\left( \frac{2^{m}\ell \left( I\right) }{\ell \left( K\right) }%
\right) ^{\varepsilon }\ell \left( K\right) +2^{m}\ell \left( I\right) &\leq
&\frac{1}{2}\left( \frac{\ell \left( I\right) }{\ell \left( K\right) }%
\right) ^{\delta }\ell \left( K\right) ; \\
\left( \frac{2^{m}\ell \left( I\right) }{\ell \left( K\right) }\right)
^{\varepsilon }+2^{m+1}\frac{\ell \left( I\right) }{\ell \left( K\right) }
&\leq &\left( \frac{\ell \left( I\right) }{\ell \left( K\right) }\right)
^{\delta }; \\
2^{m\varepsilon }\left( \frac{\ell \left( I\right) }{\ell \left( K\right) }%
\right) ^{\varepsilon -\delta }+2^{m+1}\left( \frac{\ell \left( I\right) }{%
\ell \left( K\right) }\right) ^{1-\delta } &\leq &1,\ \ \ \ \ 0\leq m\leq 
\mathbf{\tau }.
\end{eqnarray*}%
Since $0<\delta <\varepsilon <1$ by our restriction on $\varepsilon $ and
our choice of $\delta $, and since $\frac{\ell \left( I\right) }{\ell \left(
K\right) }\leq 2^{-\left( m+\mathbf{r}\right) }$, it thus suffices to show
that%
\begin{eqnarray*}
2^{m\varepsilon }\left( 2^{-\left( m+\mathbf{r}\right) }\right)
^{\varepsilon -\delta }+2^{m+1}\left( 2^{-\left( m+\mathbf{r}\right)
}\right) ^{1-\delta } &\leq &1, \\
\text{i.e. }2^{m\varepsilon -\left( m+\mathbf{r}\right) \left( \varepsilon
-\delta \right) }+2^{m+1-\left( m+\mathbf{r}\right) \left( 1-\delta \right)
} &\leq &1,
\end{eqnarray*}%
for $0\leq m\leq \mathbf{\tau }$. In particular then it suffices to show both%
\begin{eqnarray*}
m\varepsilon -\left( m+\mathbf{r}\right) \left( \varepsilon -\delta \right)
&\leq &-1, \\
m+1-\left( m+\mathbf{r}\right) \left( 1-\delta \right) &\leq &-1,
\end{eqnarray*}%
equivalently both%
\begin{eqnarray*}
\left( \mathbf{r}+m\right) \delta &\leq &\mathbf{r}\varepsilon -1, \\
\left( \mathbf{r}+m\right) \delta &\leq &\mathbf{r}-2,
\end{eqnarray*}%
for $0\leq m\leq \mathbf{\tau }$. Finally then it suffices to show both%
\begin{equation*}
\delta \leq \frac{\mathbf{r}\varepsilon -1}{\mathbf{r}+\mathbf{\tau }}\text{
and }\delta \leq \frac{\mathbf{r}-2}{\mathbf{r}+\mathbf{\tau }}.
\end{equation*}%
Because of the restriction that $\frac{1}{\mathbf{r}}<\varepsilon <1-\frac{1%
}{\mathbf{r}}$, we see that $0<\mathbf{r}\varepsilon -1<\mathbf{r}-2$, and
it is now clear that the above display holds for our choice of $\delta $.
\end{proof}

\bigskip

{\Large \#2:} In the integral in line $3$ of page 93, the factor $%
s_{I}\left( y\right) $ should be raised to the power $n-\alpha $.

\bigskip

{\large \#3:} The definition of $\varphi _{F}$ in line $-2$ of page 124
should read%
\begin{equation*}
\varphi _{F}\equiv \sum_{k,\theta :\ \theta \left( I_{k}\right) \in \mathcal{%
F}}\mathbf{1}_{\theta \left( I_{k}\right) }\left( \mathbb{E}_{\theta \left(
I_{k}\right) }^{\sigma }f-\mathbb{E}_{I_{K}}^{\sigma }f\right) ,
\end{equation*}%
and the display in line $3$ of page 125 should read%
\begin{equation*}
\left\vert \sum_{F\in \mathcal{F}}\left\langle T_{\sigma }^{\alpha }\varphi
_{F},g_{F}\right\rangle \right\vert \lesssim \sqrt{A_{2}^{\alpha }}%
\left\Vert f\right\Vert _{L^{2}\left( \sigma \right) }\left\Vert
g\right\Vert _{L^{2}\left( \omega \right) }\ ,
\end{equation*}%
and finally in line $-3$ of page 125, the constant $\mathcal{NTV}_{\alpha }$
should be replaced by $\mathfrak{T}_{T^{\alpha }}$ in both of its
appearances.

\bigskip

{\Large \#4:} The\ display beginning with the term $S$ in line $4$ of page
134 should be replaced with 
\begin{eqnarray*}
S &=&\sum_{F\in \mathcal{F}_{I}}\sum_{J\in \mathcal{M}_{\mathbf{r}-\limfunc{%
deep}}(F)}\left( \sum_{F^{\prime }\in \mathcal{F}:\ F\subset F^{\prime
}\subsetneqq I}\frac{d\left( F^{\prime }\right) }{d\left( F^{\prime }\right) 
}\frac{\mathrm{P}^{\alpha }\left( J,\mathbf{1}_{\pi _{\mathcal{F}%
_{I}}F^{\prime }\setminus F^{\prime }}\sigma \right) }{\left\vert
J\right\vert ^{1/n}}\right) ^{2}\left\Vert \mathsf{P}_{F,J}^{\omega }\mathbf{%
x}\right\Vert _{L^{2}\left( \omega \right) }^{2} \\
&\leq &\sum_{F^{\prime }\in \mathcal{F}_{I}}d\left( F^{\prime }\right)
^{2}\sum_{F\in \mathcal{F}:\ F\subset F^{\prime }}\sum_{J\in \mathcal{M}_{%
\mathbf{r}-\limfunc{deep}}(F)}\left( \sum_{F^{\prime }\in \mathcal{F}:\
F\subset F^{\prime }\subsetneqq I}\frac{1}{d\left( F^{\prime }\right) ^{2}}%
\right) \left( \frac{\mathrm{P}^{\alpha }\left( J,\mathbf{1}_{\pi _{\mathcal{%
F}_{I}}F^{\prime }\setminus F^{\prime }}\sigma \right) }{\left\vert
J\right\vert ^{1/n}}\right) ^{2}\left\Vert \mathsf{P}_{F,J}^{\omega }\mathbf{%
x}\right\Vert _{L^{2}\left( \omega \right) }^{2} \\
&\leq &C\sum_{F^{\prime }\in \mathcal{F}_{I}}d\left( F^{\prime }\right)
^{2}\sum_{K\in \mathcal{M}_{\mathbf{r}-\limfunc{deep}}\left( F^{\prime
}\right) }\sum_{F\in \mathcal{F}:\ F\subset F^{\prime }}\sum_{J\in \mathcal{M%
}_{\mathbf{r}-\limfunc{deep}}(F):\ J\subset K}\left( \frac{\mathrm{P}%
^{\alpha }\left( J,\mathbf{1}_{\pi _{\mathcal{F}_{I}}F^{\prime }\setminus
F^{\prime }}\sigma \right) }{\left\vert J\right\vert ^{1/n}}\right)
^{2}\left\Vert \mathsf{P}_{F,J}^{\omega }\mathbf{x}\right\Vert _{L^{2}\left(
\omega \right) }^{2} \\
&\lesssim &\sum_{F^{\prime }\in \mathcal{F}_{I}}d\left( F^{\prime }\right)
^{2}\sum_{K\in \mathcal{M}_{\mathbf{r}-\limfunc{deep}}\left( F^{\prime
}\right) }\left( \frac{\mathrm{P}^{\alpha }\left( K,\mathbf{1}_{\pi _{%
\mathcal{F}_{I}}F^{\prime }\setminus F^{\prime }}\sigma \right) }{\left\vert
K\right\vert ^{1/n}}\right) ^{2}\sum_{F\in \mathcal{F}:\ F\subset F^{\prime
}}\sum_{J\in \mathcal{M}_{\mathbf{r}-\limfunc{deep}}(F):\ J\subset
K}\left\Vert \mathsf{P}_{F,J}^{\omega }\mathbf{x}\right\Vert _{L^{2}\left(
\omega \right) }^{2},
\end{eqnarray*}%
and then the display at the bottom of page 134 should be replaced with%
\begin{eqnarray*}
S &\lesssim &\sum_{F^{\prime }\in \mathcal{F}_{I}}d\left( F^{\prime }\right)
^{2}\sum_{K\in \mathcal{M}_{\mathbf{r}-\limfunc{deep}}\left( F^{\prime
}\right) }\left( \frac{\mathrm{P}^{\alpha }\left( K,\mathbf{1}_{\pi _{%
\mathcal{F}_{I}}F^{\prime }\setminus F^{\prime }}\sigma \right) }{\left\vert
K\right\vert ^{1/n}}\right) ^{2}\left\Vert \widehat{\mathsf{P}}_{F^{\prime
},K}^{\omega }\mathbf{x}\right\Vert _{L^{2}\left( \omega \right) }^{2} \\
&=&\sum_{k=0}^{\infty }k^{2}\sum_{F^{\prime }\in \mathcal{F}_{I}:\ d\left(
F^{\prime }\right) =k}\sum_{K\in \mathcal{M}_{\mathbf{r}-\limfunc{deep}%
}\left( F^{\prime }\right) }\left( \frac{\mathrm{P}^{\alpha }\left( K,%
\mathbf{1}_{\pi _{\mathcal{F}_{I}}F^{\prime }\setminus F^{\prime }}\sigma
\right) }{\left\vert K\right\vert ^{1/n}}\right) ^{2}\left\Vert \widehat{%
\mathsf{P}}_{F^{\prime },K}^{\omega ,\mathbf{b}^{\ast }}\mathbf{x}%
\right\Vert _{L^{2}\left( \omega \right) }^{2}\equiv \sum_{k=0}^{\infty
}A_{k}.
\end{eqnarray*}

\bigskip

{\Large \#5:} Beginning two lines above the display at the bottom of page
144, and ending after the display at the top of page 145, replace with the
following:%
\begin{equation*}
\sum_{F^{\prime }\in \mathcal{F}}\left\Vert \mathsf{Q}_{F^{\prime
},A}^{\omega }x\right\Vert _{L^{2}\left( \omega \right) }^{2}\lesssim 
\mathbf{\tau }\left\vert A\right\vert ^{\frac{2}{n}}\left\vert A\right\vert
_{\omega }\ ,\ \ \ \ \ \text{for any cube }A.
\end{equation*}%
Note that if $c_{J^{\prime }}\in 2J$ and $\ell \left( J^{\prime }\right)
<\ell \left( J\right) $, then $J^{\prime }\subset \frac{5}{2}J$. Then with 
\begin{equation*}
\mathcal{W}_{J}^{s}\equiv \dbigcup\limits_{F^{\prime }\in \mathcal{F}%
}\left\{ J^{\prime }\in \mathcal{M}_{\mathbf{r}-\limfunc{deep}}\left(
F^{\prime }\right) :J^{\prime }\subset \frac{5}{2}J\text{ and }\ell \left(
J^{\prime }\right) =2^{-s}\ell \left( J\right) \right\} ,
\end{equation*}%
we have%
\begin{eqnarray*}
\mathcal{A}_{s,far}^{0}\left( J\right) &\leq &\sum_{F^{\prime }\in \mathcal{F%
}}\overset{\ast }{\sum_{c_{J^{\prime }}\in 2J}}\int_{I\setminus \left(
3J\right) }S_{\left( J^{\prime },J\right) }^{F^{\prime }}\left( y\right)
d\sigma \left( y\right) \\
&\lesssim &\sum_{A\in \mathcal{W}_{J}^{s}}\sum_{F^{\prime }\in \mathcal{F}:\
A\in \mathcal{M}_{\mathbf{r}-\limfunc{deep}}\left( F^{\prime }\right)
}\int_{I\setminus \left( 3J\right) }\frac{\left\Vert \mathsf{Q}_{F^{\prime
},A}^{\omega }x\right\Vert _{L^{2}\left( \omega \right) }^{\spadesuit 2}}{%
\left( \left\vert J\right\vert ^{\frac{2}{n}}+\left\vert y-c_{J}\right\vert
\right) ^{2\left( 2-\alpha \right) }}d\sigma \left( y\right) \\
&\lesssim &\sum_{A\in \mathcal{W}_{J}^{s}}\int_{I\setminus \left( 3J\right) }%
\frac{\left\vert A\right\vert ^{\frac{2}{n}}\left\vert A\right\vert _{\omega
}}{\left( \left\vert J\right\vert ^{\frac{2}{n}}+\left\vert
y-c_{M}\right\vert \right) ^{2\left( 2-\alpha \right) }}d\sigma \left(
y\right) \\
&=&\left( \sum_{A\in \mathcal{W}_{J}^{s}}\left\vert A\right\vert ^{\frac{2}{n%
}}\left\vert A\right\vert _{\omega }\right) \int_{I\setminus \left(
3J\right) }\frac{1}{\left( \left\vert J\right\vert ^{\frac{2}{n}}+\left\vert
y-c_{J}\right\vert \right) ^{2\left( 2-\alpha \right) }}d\sigma \left(
y\right) .
\end{eqnarray*}%
Now we use the standard pigeonholing of side length of $A$ to conclude that 
\begin{eqnarray*}
\sum_{A\in \mathcal{W}_{J}^{s}}\left\vert A\right\vert ^{\frac{2}{n}%
}\left\vert A\right\vert _{\omega } &=&\sum_{k=s}^{\infty }\sum_{A\in 
\mathcal{W}_{J}^{s}:\ \ell \left( A\right) =2^{-k}\ell \left( J\right)
}\left\vert A\right\vert ^{\frac{2}{n}}\left\vert A\right\vert _{\omega
}\leq \sum_{k=s}^{\infty }2^{-2k}\left\vert J\right\vert ^{\frac{2}{n}%
}\sum_{A\in \mathcal{W}_{J}^{s}:\ \ell \left( A\right) =2^{-k}\ell \left(
J\right) }\left\vert A\right\vert _{\omega } \\
&\leq &\sum_{k=s}^{\infty }2^{-2k}\left\vert J\right\vert ^{2}\left\vert 
\frac{5}{2}J\right\vert _{\omega }\lesssim 2^{-2s}\left\vert J\right\vert
^{2}\left\vert \frac{5}{2}J\right\vert _{\omega },
\end{eqnarray*}%
so that combining the previous two displays we have%
\begin{eqnarray*}
\mathcal{A}_{s,far}^{0}\left( J\right) &\lesssim &2^{-2s}\left\vert
J\right\vert ^{2}\left\vert \frac{5}{2}J\right\vert _{\omega
}\int_{I\setminus \left( 3J\right) }\frac{1}{\left( \left\vert J\right\vert
^{\frac{2}{n}}+\left\vert y-c_{J}\right\vert \right) ^{2\left( 2-\alpha
\right) }}d\sigma \left( y\right) \\
&\leq &2^{-2s}\left\vert \frac{5}{2}J\right\vert _{\omega }\int_{I\setminus
\left( 3J\right) }\frac{1}{\left( \left\vert J\right\vert ^{\frac{2}{n}%
}+\left\vert y-c_{J}\right\vert \right) ^{2\left( 1-\alpha \right) }}d\sigma
\left( y\right) \\
&\approx &2^{-2s}\frac{\left\vert \frac{5}{2}J\right\vert _{\omega }}{%
\left\vert \frac{5}{2}J\right\vert ^{1-\alpha }}\int_{I\setminus \left(
3J\right) }\left( \frac{\left\vert J\right\vert }{\left( \left\vert
J\right\vert ^{\frac{2}{n}}+\left\vert y-c_{J}\right\vert \right) ^{2}}%
\right) ^{1-\alpha }d\sigma \left( y\right) \\
&\lesssim &2^{-2s}\frac{\left\vert \frac{5}{2}J\right\vert _{\omega }}{%
\left\vert \frac{5}{2}J\right\vert ^{1-\alpha }}\mathcal{P}^{\alpha }\left( 
\frac{5}{2}J,\mathbf{1}_{I\setminus \left( 3J\right) }\sigma \right)
\lesssim 2^{-2s}\mathcal{A}_{2}^{\alpha }\ .
\end{eqnarray*}

\bigskip

{\Large \#6:} Lines $3$ through $11$ on page 146 should be replaced with
this:

\bigskip

Similarly, for $\ell \geq 1$, we can estimate the far term $\mathcal{A}%
_{s,far}^{\ell }\left( J\right) $ by the argument used for $\mathcal{A}%
_{s,far}^{0}\left( J\right) $ but applied to $2^{\ell }J$ in place of $J$.
For this need the following variant of $\mathcal{W}_{J}^{s}$ given by%
\begin{equation*}
\mathcal{W}_{J}^{s,\ell }\equiv \dbigcup\limits_{F^{\prime }\in \mathcal{F}%
}\left\{ J^{\prime }\in \mathcal{W}\left( F^{\prime }\right) :J^{\prime
}\subset 3\left( 2^{\ell }J\right) \text{ and }\ell \left( J^{\prime
}\right) =2^{-s-\ell }\ell \left( 2^{\ell }J\right) \right\} .
\end{equation*}%
Then we have

\begin{eqnarray*}
\mathcal{A}_{s,far}^{\ell }\left( J\right) &\leq &\sum_{F^{\prime }\in 
\mathcal{F}}\overset{\ast }{\sum_{c_{J^{\prime }}\in \left( 2^{\ell
+1}J\right) \setminus \left( 2^{\ell }J\right) }}\int_{I\setminus 2^{\ell
+2}J}S_{\left( J^{\prime },J\right) }^{F^{\prime }}\left( y\right) d\sigma
\left( y\right) \\
&\lesssim &\sum_{A\in \mathcal{W}_{J}^{s,\ell }}\sum_{F^{\prime }\in 
\mathcal{F}:\ A\in \mathcal{M}_{\mathbf{r}-\limfunc{deep}}\left( F^{\prime
}\right) }\int_{I\setminus 4\left( 2^{\ell }J\right) }\frac{\left\Vert 
\mathsf{Q}_{F^{\prime },A}^{\omega ,\mathbf{b}^{\ast }}x\right\Vert
_{L^{2}\left( \omega \right) }^{\spadesuit 2}}{\left( \left\vert
J\right\vert ^{\frac{2}{n}}+\left\vert y-c_{J}\right\vert \right) ^{2\left(
2-\alpha \right) }}d\sigma \left( y\right) \\
&\lesssim &\sum_{A\in \mathcal{W}_{J}^{s,\ell }}\int_{I\setminus 4\left(
2^{\ell }J\right) }\frac{\left\vert A\right\vert ^{\frac{2}{n}}\left\vert
A\right\vert _{\omega }}{\left( \left\vert J\right\vert ^{\frac{2}{n}%
}+\left\vert y-c_{J}\right\vert \right) ^{2\left( 2-\alpha \right) }}d\sigma
\left( y\right) \\
&=&\left( \sum_{A\in \mathcal{W}_{J}^{s,\ell }}\left\vert A\right\vert ^{%
\frac{2}{n}}\left\vert A\right\vert _{\omega }\right) \int_{I\setminus
4\left( 2^{\ell }J\right) }\frac{1}{\left( \left\vert J\right\vert ^{\frac{2%
}{n}}+\left\vert y-c_{J}\right\vert \right) ^{2\left( 2-\alpha \right) }}%
d\sigma \left( y\right) ,
\end{eqnarray*}%
where, just as for the sum over $A\in \mathcal{W}_{J}^{s,0}$,%
\begin{eqnarray*}
\sum_{A\in \mathcal{W}_{J}^{s,\ell }}\left\vert A\right\vert ^{\frac{2}{n}%
}\left\vert A\right\vert _{\omega } &=&\sum_{k=s}^{\infty }\sum_{A\in 
\mathcal{W}_{J}^{s,\ell }:\ \ell \left( A\right) =2^{-k-\ell }\ell \left(
2^{\ell }J\right) }\left\vert A\right\vert ^{\frac{2}{n}}\left\vert
A\right\vert _{\omega }\leq \sum_{k=s}^{\infty }2^{-2k-2\ell }\left\vert
J\right\vert ^{\frac{2}{n}}\sum_{A\in \mathcal{W}_{M}^{s,\ell }:\ \ell
\left( A\right) =2^{-k-\ell }\ell \left( 2^{\ell }J\right) }\left\vert
A\right\vert _{\omega } \\
&\leq &\sum_{k=s}^{\infty }2^{-2k-2\ell }\left\vert 2^{\ell }J\right\vert ^{%
\frac{2}{n}}\left\vert 3\left( 2^{\ell }J\right) \right\vert _{\omega
}\lesssim 2^{-2s-2\ell }\left\vert 2^{\ell }J\right\vert ^{\frac{2}{n}%
}\left\vert 3\left( 2^{\ell }J\right) \right\vert _{\omega }\ .
\end{eqnarray*}%
Now we can continue with%
\begin{eqnarray*}
&&\mathcal{A}_{s,far}^{\ell }\left( J\right) \lesssim 2^{-2s}2^{-2\ell
}\left\vert 3\left( 2^{\ell }J\right) \right\vert _{\omega }\int_{I\setminus
4\left( 2^{\ell }J\right) }\frac{1}{\left( \left\vert 2^{\ell }J\right\vert
^{\frac{2}{n}}+\left\vert y-c_{2^{\ell }J}\right\vert \right) ^{2\left(
1-\alpha \right) }}d\sigma \left( y\right) \\
&\approx &2^{-2s}2^{-2\ell }\frac{\left\vert 3\left( 2^{\ell }J\right)
\right\vert _{\omega }}{\left\vert 3\left( 2^{\ell }J\right) \right\vert
^{1-\alpha }}\int_{I\setminus 4\left( 2^{\ell }J\right) }\left( \frac{%
\left\vert 2^{\ell }J\right\vert ^{\frac{2}{n}}}{\left( \left\vert 2^{\ell
}J\right\vert ^{\frac{2}{n}}+\left\vert y-c_{2^{\ell }J}\right\vert \right)
^{2}}\right) ^{1-\alpha }d\sigma \left( y\right) \\
&\lesssim &2^{-2s}2^{-2\ell }\left\{ \frac{\left\vert 3\left( 2^{\ell
}J\right) \right\vert _{\omega }}{\left\vert 3\left( 2^{\ell }J\right)
\right\vert ^{1-\alpha }}\mathcal{P}^{\alpha }\left( 3\left( 2^{\ell
}J\right) ,\mathbf{1}_{I\setminus 4\left( 2^{\ell }J\right) }\sigma \right)
\right\} \lesssim 2^{-2s}2^{-2\ell }\mathcal{A}_{2}^{\alpha }\ .
\end{eqnarray*}

\bigskip

{\Large \#7:} The three lines above Section 11 should be replaced with these
three lines in quotes\footnote{%
A detailed argument is given in Subsection 9.3 of \cite[v3]{SaShUr6}, and in
Subsection \ref{Subsec back test} of Appendix B above for the one
dimensional case adapted to $Tb$.}: "and then expanding the square and
calculating as in the proof of the local part given earlier to obtain the
bound $\mathcal{A}_{2}^{\alpha }+\left( \mathcal{E}_{\alpha }^{\func{plug}}+%
\sqrt{A_{2}^{\alpha ,\func{energy}}}\right) \sqrt{A_{2}^{\alpha ,\limfunc{%
punct}}}$. The details are similar and they are left to the reader."

\bigskip

{\Large \#8:} In line -3 of page 164 the final factor on the right hand side
should instead be $\omega _{\mathcal{P}_{L,0}^{\limfunc{small}}}\left( 
\mathbf{T}^{\mathbf{\tau }-\limfunc{deep}}\left( K\right) \right) $.

\bigskip

{\Large \#9:} There is a gap in the proof of the Orthogonality Lemma at the
top of page 170, where the restricted norm of the sublinear form used there
doesn't match the definition on page 161. The required change in definition
of the restricted norm, then forces an additional argument - due to Michael
Lacey in \cite{Lac} - that uses another Calder\'{o}n-Zygmund decomposition
in the proof of the Straddling Lemma on page 166. Here are the changes
required (see also Section 7 above, which adapts Lacey's additional argument
to the $Tb$ setting, in which the dual martingale differences are no longer
orthogonal projections).

First, the display on page 161 right before Proposition 11.8 should be
replaced with%
\begin{equation*}
\left\vert \mathsf{B}\right\vert _{\limfunc{stop},1,\bigtriangleup ^{\omega
}}^{A,\mathcal{P}}\left( f,g\right) \leq \mathfrak{N}_{\limfunc{stop}%
,1,\bigtriangleup ^{\omega }}^{A,\mathcal{P}}\left\Vert f\right\Vert
_{L^{2}\left( \sigma \right) }\left\Vert g\right\Vert _{L^{2}\left( \omega
\right) }\ ,
\end{equation*}%
in which the term $\alpha _{\mathcal{A}}\left( A\right) \sqrt{\left\vert
A\right\vert _{\sigma }}$ no longer appears on the right hand side.

Second, the sentence on page 169 right before Subsubsection 11.4.2 should be
replaced with the following material:

\bigskip

Now we sum these bounds in $s$ and $\ast $ and use $\sup_{S\in \mathcal{S}}%
\mathcal{S}_{\limfunc{size}}^{\alpha ,A;S}\left( \mathcal{Q}\right) \leq 
\mathcal{S}_{\limfunc{size}}^{\alpha ,A}\left( \mathcal{Q}\right) $ to obtain%
\begin{equation*}
\left\vert \mathsf{B}\right\vert _{\limfunc{stop},1,\bigtriangleup ^{\omega
}}^{A,\mathcal{Q}}\left( f,g\right) \leq \mathcal{S}_{\limfunc{size}%
}^{\alpha ,A}\left( \mathcal{Q}\right) \left\{ \alpha _{\mathcal{A}}\left(
A\right) \sqrt{\left\vert A\right\vert _{\sigma }}+\left\Vert f\right\Vert
_{L^{2}\left( \sigma \right) }\right\} \left\Vert g\right\Vert _{L^{2}\left(
\omega \right) }\ .
\end{equation*}%
However, this inequality has the unwanted term $\alpha _{\mathcal{A}}\left(
A\right) \sqrt{\left\vert A\right\vert _{\sigma }}$ included on the right
hand side, and we must apply an argument of M. Lacey \cite[see the proof of
Lemma 3.19]{Lac} to eliminate this term. We begin with%
\begin{eqnarray*}
\left\vert \mathsf{B}\right\vert _{\limfunc{stop},1,\bigtriangleup ^{\omega
}}^{A,\mathcal{Q}}\left( f,g\right) &=&\left\vert \mathsf{B}\right\vert _{%
\limfunc{stop},1,\bigtriangleup ^{\omega }}^{A,\mathcal{Q}}\left( h,g\right)
, \\
\text{where }h &\equiv &\mathsf{P}_{\pi \left( \Pi _{1}\mathcal{Q}\right)
}^{\sigma }f=\sum_{I\in \Pi _{1}\mathcal{Q}}\bigtriangleup _{\pi I}^{\sigma
}f,
\end{eqnarray*}%
which follows from the formula $\varphi _{J}^{\mathcal{Q}}\equiv \sum_{I\in 
\mathcal{C}_{A}^{\prime }:\ \left( I,J\right) \in \mathcal{Q}}\mathbb{E}%
_{I}^{\sigma }\left( \bigtriangleup _{\pi I}^{\sigma }f\right) \ \mathbf{1}%
_{A\setminus I}$ since $\bigtriangleup _{\pi I}^{\sigma }f=\bigtriangleup
_{\pi I}^{\sigma }h$ for $I\in \Pi _{1}\mathcal{Q}$, which holds in turn
because the Haar projections are orthogonal. Now define Calder\'{o}n-Zygmund
stopping times $\mathcal{H}$ for the function $h=\mathsf{P}_{\pi \left( \Pi
_{1}\mathcal{Q}\right) }^{\sigma }f\in L^{2}\left( \sigma \right) $, so that
the quasiorthogonal inequality%
\begin{equation*}
\sum_{H\in \mathcal{H}}\alpha _{\mathcal{H}}\left( H\right) ^{2}\left\vert
H\right\vert _{\sigma }\lesssim \left\Vert h\right\Vert _{L^{2}\left( \sigma
\right) }^{2}=\left\Vert \mathsf{P}_{\pi \left( \Pi _{1}\mathcal{Q}\right)
}^{\sigma }f\right\Vert _{L^{2}\left( \sigma \right) }^{2}
\end{equation*}%
holds with $\alpha _{\mathcal{H}}\left( H\right) =\mathbb{E}_{H}^{\sigma
}\left\vert h\right\vert $. We also define in the usual way coronas $%
\mathcal{C}_{H}^{\mathcal{H}}$ and $\mathcal{C}_{H}^{\mathcal{H},\mathbf{%
\tau }-\func{shift}}$.

Now we return to the previous inequalities we obtained for $\left\vert 
\mathsf{B}\right\vert _{\limfunc{stop},1,\bigtriangleup ^{\omega }}^{A,%
\mathcal{Q}_{s}}\left( f,g\right) $ and $\left\vert \mathsf{B}\right\vert _{%
\limfunc{stop},1,\bigtriangleup ^{\omega }}^{A,\mathcal{Q}_{\ast }}\left(
f,g\right) $, and replace the collection $\mathcal{Q}$ with the $A$%
-admissible collection $\mathcal{Q}_{H}\equiv \left\{ \left( I,J\right) \in 
\mathcal{Q}:J\in \mathcal{C}_{H}^{\mathcal{H},\mathbf{\tau }-\func{shift}%
}\right\} $, so that the arguments given there can be adapted to yield%
\begin{eqnarray*}
\left\vert \mathsf{B}\right\vert _{\limfunc{stop},1,\bigtriangleup ^{\omega
}}^{A,\left( \mathcal{Q}_{H}\right) _{s}}\left( f,g\right) &\leq &\sup_{S\in 
\mathcal{S}}\mathcal{S}_{\limfunc{size}}^{\alpha ,A;S}\left( \mathcal{Q}%
_{H}\right) \alpha _{\mathcal{H}}\left( H\right) \sqrt{\left\vert
H\right\vert _{\sigma }}\left\Vert \mathsf{P}_{\mathcal{C}_{H}^{\mathcal{H},%
\mathbf{\tau }-\func{shift}}}^{\omega }g\right\Vert _{L^{2}\left( \omega
\right) }\ ,\ \ \ \ \ \mathbf{\tau }\leq s\leq \mathbf{\rho }-1, \\
\left\vert \mathsf{B}\right\vert _{\limfunc{stop},1,\bigtriangleup ^{\omega
}}^{A,\left( \mathcal{Q}_{H}\right) _{\ast }}\left( f,g\right) &\leq
&\sup_{S\in \mathcal{S}}\mathcal{S}_{\limfunc{size}}^{\alpha ,A;S}\left( 
\mathcal{Q}_{H}\right) \alpha _{\mathcal{H}}\left( H\right) \sqrt{\left\vert
H\right\vert _{\sigma }}\left\Vert \mathsf{P}_{\mathcal{C}_{H}^{\mathcal{H},%
\mathbf{\tau }-\func{shift}}}^{\omega }g\right\Vert _{L^{2}\left( \omega
\right) }\ .
\end{eqnarray*}%
We then sum these improved bounds in $s$ and $\ast $ to obtain%
\begin{eqnarray*}
\left\vert \mathsf{B}\right\vert _{\limfunc{stop},1,\bigtriangleup ^{\omega
}}^{A,\mathcal{Q}}\left( f,g\right) &\leq &\sup_{S\in \mathcal{S}}\mathcal{S}%
_{\limfunc{size}}^{\alpha ,A;S}\left( \mathcal{Q}\right) \sum_{H\in \mathcal{%
H}}\alpha _{\mathcal{H}}\left( H\right) \sqrt{\left\vert H\right\vert
_{\sigma }}\left\Vert \mathsf{P}_{\mathcal{C}_{H}^{\mathcal{H},\mathbf{\tau }%
-\func{shift}}}^{\omega }g\right\Vert _{L^{2}\left( \omega \right) } \\
&\leq &\sup_{S\in \mathcal{S}}\mathcal{S}_{\limfunc{size}}^{\alpha
,A;S}\left( \mathcal{Q}\right) \sqrt{\sum_{H\in \mathcal{H}}\alpha _{%
\mathcal{H}}\left( H\right) ^{2}\left\vert H\right\vert _{\sigma }}\sqrt{%
\sum_{H\in \mathcal{H}}\left\Vert \mathsf{P}_{\mathcal{C}_{H}^{\mathcal{H},%
\mathbf{\tau }-\func{shift}}}^{\omega }g\right\Vert _{L^{2}\left( \omega
\right) }^{2}} \\
&\leq &\mathcal{S}_{\limfunc{size}}^{\alpha ,A}\left( \mathcal{Q}\right)
\left\Vert f\right\Vert _{L^{2}\left( \sigma \right) }\left\Vert
g\right\Vert _{L^{2}\left( \omega \right) }\ ,
\end{eqnarray*}%
where we have used the quasiorthogonal inequality and $\left\Vert \mathsf{P}%
_{\pi \left( \Pi _{1}\mathcal{Q}\right) }^{\sigma }f\right\Vert
_{L^{2}\left( \sigma \right) }^{2}\leq \left\Vert f\right\Vert _{L^{2}\left(
\sigma \right) }^{2}$ in the last line.

\section{Appendix D: Glossary}

\subsubsection{Section 1}

\begin{enumerate}
\item $C_{CZ}$; (\ref{sizeandsmoothness'})

\item $\mathfrak{N}_{T_{\sigma }^{\alpha }}$; (\ref{two weight'})

\item $T_{\sigma ,\delta ,R}^{\alpha }f\left( x\right) $; (\ref{def
truncation})

\item $p$\emph{-weakly }$\mu $\emph{-accretive} family $\mathbf{b}=\left\{
b_{Q}\right\} _{Q\in \mathcal{P}}$ of functions on $\mathbb{R}$; (\ref{local
accretive})

\item $\mathfrak{T}_{T^{\alpha }}^{\mathbf{b}},\mathfrak{T}_{T^{\alpha ,\ast
}}^{\mathbf{b}^{\ast }}$; (\ref{b testing cond})

\item $\mathrm{P}^{\alpha }\left( Q,\mu \right) ,\mathcal{P}^{\alpha }\left(
Q,\mu \right) $; (\ref{def Poisson})

\item $\mathcal{A}_{2}^{\alpha },\mathcal{A}_{2}^{\alpha ,\ast }$; (\ref{def
call A2})

\item $\mathfrak{P}_{\left( \sigma ,\omega \right) }$; (\ref{def common
point mass})

\item $A_{2}^{\alpha ,\limfunc{punct}},A_{2}^{\alpha ,\ast ,\limfunc{punct}}$%
; (\ref{def punct})

\item $\mathfrak{A}_{2}^{\alpha }$; (\ref{def A2})

\item $\mathcal{E}_{2}^{\alpha },\mathcal{E}_{2}^{\alpha ,\ast },\mathfrak{E}%
_{2}^{\alpha },\mathcal{NTV}_{\alpha }$;$\ $(\ref{strong b* energy}), (\ref%
{strong b energy}), (\ref{def frak energy}), (\ref{def NTV})
\end{enumerate}

\subsubsection{Section 2}

\begin{enumerate}
\item \emph{gradient elliptic} kernel $K^{\alpha }\left( x,y\right) $; (\ref%
{def grad elliptic})

\item $PLBP,$pointwise lower bound property; (\ref{plb})

\item reverse H\"{o}lder control of children (\ref{rev Hol con})

\item Calder\'{o}n-Zygmund stopping intervals; Definition \ref{CZ stopping
times}

\item $\mathbf{b}$-accretive/weak testing stopping intervals; Definition \ref%
{accretive stopping times gen}

\item $\sigma $-energy stopping intervals; Definition \ref{def energy corona
3}

\item $\mathbf{X}_{\alpha }\left( \mathcal{C}_{S}\right) $ energy stopping
times; (\ref{def stopping energy 3})

\item $\mathcal{D}_{\beta }$: (\ref{def dyadic grid})
\end{enumerate}

\subsubsection{Section 3}

\begin{enumerate}
\item $\limfunc{skel}K$,$\ \mathbb{S}_{x}^{\limfunc{dy}},\limfunc{body}K$
skeleton, spray, body of an interval; (\ref{skel spray body})

\item $\varepsilon -\limfunc{good}$ \emph{with respect to} an interval; (\ref%
{eps far})

\item $\mathcal{G}_{\left( k,\varepsilon \right) -\limfunc{good}}^{\mathcal{D%
}}$; Definition \ref{good two grids}

\item $k$-$\limfunc{bad}$ in a grid; (\ref{bad in grid})

\item $R^{\maltese },\kappa \left( R\right) $: \ref{def sharp cross}

\item $\Theta _{2}^{\limfunc{bad}\natural }\left( f,g\right) $; (\ref%
{Theta_2^bad sharp})

\item $\mathcal{G}_{k-\limfunc{bad}}^{\mathcal{D}}=\mathcal{G}_{\left(
k,\varepsilon \right) -\limfunc{bad}}^{\mathcal{D}}$; Definition \ref{def
Gbad}

\item $\Theta _{2}^{\limfunc{good}}\left( f,g\right) $; (\ref{def Theta 2
good})

\item $\mathfrak{3T}_{T^{\alpha }}^{\mathbf{b}}$; (\ref{triple b testing
cond})

\item $\mathfrak{FT}_{T^{\alpha }}^{\mathbf{b}}$; (\ref{full b testing})

\item $\mathcal{E}_{2}^{\alpha ,\limfunc{triple}}$; Definition \ref{def
triple energy}
\end{enumerate}

\subsubsection{Section 4}

\begin{enumerate}
\item $\Theta _{1}^{\limfunc{long}}\left( f,g\right) ,\Theta _{1}^{\limfunc{%
short}}\left( f,g\right) $; (\ref{decomp long short})
\end{enumerate}

\subsubsection{Section 5}

\begin{enumerate}
\item $\partial _{\eta }L$; (\ref{def halo})

\item $\mathcal{K}\left( I^{\prime },J^{\prime }\right) $; (\ref{def
K(I',J')})

\item $\left\{ E,F\right\} $; (\ref{def E,F})

\item $\mathrm{P}_{\delta }^{\alpha }\mathsf{Q}^{\omega }\left( J,\upsilon
\right) $; (\ref{def compact})

\item $\left\{ K_{\limfunc{out}},K_{\limfunc{in}}\right\} ^{\limfunc{orig}}$%
; (\ref{def orig})
\end{enumerate}

\subsubsection{Section 6}

\begin{enumerate}
\item $\mathcal{C}_{B}^{\mathcal{G},\limfunc{shift}}$; Definition \ref%
{shifted corona}

\item $\left\langle T_{\sigma }^{\alpha }\left( \mathsf{P}_{\mathcal{C}_{A}^{%
\mathcal{D}}}^{\sigma ,\mathbf{b}}f\right) ,\mathsf{P}_{\mathcal{C}_{B}^{%
\mathcal{G},\limfunc{shift}}}^{\omega ,\mathbf{b}^{\ast }}g\right\rangle
_{\omega }^{\Subset _{\mathbf{r},\varepsilon }}$; (\ref{def shorthand})

\item $\mathsf{B}_{\Subset _{\mathbf{r},\varepsilon }}^{A}\left( f,g\right) $%
; (\ref{def local})

\item \emph{Whitney} subintervals $\mathcal{W}\left( F\right) $; (\ref{def
Whitney})

\item $\mathfrak{F}_{\alpha }=\mathfrak{F}_{\alpha }^{\mathbf{b}^{\ast
}}\left( \mathcal{D},\mathcal{G}\right) $; (\ref{functional energy n})

\item $\mathsf{B}_{\limfunc{stop}}^{A}\left( f,g\right) $; (\ref{bounded
stopping form}), (\ref{def stop}), (\ref{dummy})

\item $\mathcal{C}_{F}^{\mathcal{G},\limfunc{shift}}$; Definition \ref{def
shift}

\item $\mathsf{Q}_{\mathcal{H}}^{\omega ,\mathbf{b}^{\ast }}$; (\ref{def
localization})

\item $\mathsf{B}_{\limfunc{broken}}^{A}\left( f,g\right) $; (\ref{broken
vanish})

\item $\mathsf{B}_{\limfunc{neighbour}}^{A}\left( f,g\right) $; (\ref{def
neighbour})

\item $\mathsf{B}_{\limfunc{stop}}^{A,\mathcal{P}}\left( f,g\right) $; (\ref%
{def stop P})
\end{enumerate}

\subsubsection{Section 7}

\begin{enumerate}
\item $\varphi _{J}^{\mathcal{P}}$; (\ref{def phi P})

\item $\widehat{\mathfrak{N}}_{\limfunc{stop},\bigtriangleup ^{\omega }}^{A,%
\mathcal{P}}$; (\ref{Norm hat})

\item $\Pi _{2}^{K}\mathcal{P}$; (\ref{rest K})

\item $\Pi _{1}^{\limfunc{below}}\mathcal{P}$; (\ref{Pi below})

\item $\mathcal{S}_{\limfunc{init}\limfunc{size}}^{\alpha ,A}\left( \mathcal{%
P}\right) ^{2}$; (\ref{def ext size})

\item $\Pi _{2}^{K,\limfunc{aug}}\mathcal{P}$; Definition \ref{augs}

\item $\mathcal{S}_{\limfunc{aug}\limfunc{size}}^{\alpha ,A}\left( \mathcal{P%
}\right) $; Definition \ref{augs}, (\ref{def P stop energy' 3})

\item $\mathcal{P}_{\func{cor}}^{A}$; (\ref{def cor})

\item $J^{\flat }=J_{\searrow J}^{\maltese }$; (\ref{def aug})

\item $\mathcal{Q}$ $\flat $\textbf{straddles} $\mathcal{S}$; Definition \ref%
{flat straddles}

\item $\mathcal{W}^{\ast }\left( S\right) $; (\ref{def Whit})

\item $\mathcal{S}_{\limfunc{loc}\limfunc{size}}^{\alpha ,A;S}\left( 
\mathcal{Q}\right) ^{2}$; (\ref{localized size ref})

\item $\mathcal{Q}$ \textbf{substraddles} $L$; (\ref{def substraddles})

\item $\omega _{\mathcal{P}}$ and $\omega _{\flat \mathcal{P}}$; (\ref{def
atomic})

\item $\mathcal{G}_{d}$; (\ref{geom depth})

\item $\mathcal{P}_{L,t}^{\mathcal{H}}$; (\ref{def PHLt})

\item $\mathcal{P}^{big}$, $\mathcal{P}^{small}$; (\ref{def big small})

\item $\mathcal{P}^{\flat big}$, $\mathcal{P}^{\flat small}$; (\ref{def big
small flat})
\end{enumerate}

\subsubsection{Section 9}

\begin{enumerate}
\item $\mathbb{E}_{Q}^{\mu ,\mathbf{b}}f\left( x\right) $; (\ref{def
expectation})

\item $\widehat{\mathbb{F}}_{Q}^{\mu ,\mathbf{b}}f\left( x\right) $; (\ref{F
hat})

\item $\bigtriangleup _{Q}^{\mu ,\mathbf{b}}f\left( x\right) $, $\square
_{Q}^{\mu ,\mathbf{b}}f\left( x\right) $; (\ref{def diff})

\item $\bigtriangledown _{Q}^{\mu }$, $\widehat{\bigtriangledown }_{Q}^{\mu
} $; (\ref{Carleson avg op})

\item $\mathbb{E}_{Q}^{\mu ,\pi ,\mathbf{b}}f$, $\mathbb{F}_{Q}^{\mu ,\pi ,%
\mathbf{b}}f$; (\ref{def pi exp})

\item $\bigtriangleup _{Q}^{\mu ,\pi ,\mathbf{b}}f$, $\square _{Q}^{\mu ,\pi
,\mathbf{b}}f$; (\ref{def pi box})

\item $\square _{I}^{\sigma ,\flat ,\mathbf{b}}f$; (\ref{flat box})

\item $\widehat{\square }_{I}^{\sigma ,\flat ,\mathbf{b}}f$; (\ref{flat box
hat})

\item $\bigtriangleup _{I,\limfunc{broken}}^{\mu ,\flat ,\mathbf{b}}f$, $%
\square _{I,\limfunc{broken}}^{\mu ,\flat ,\mathbf{b}}f$; (\ref{def flat
broken})

\item $\Psi _{\mathcal{B},\mathbf{\lambda }}^{\mu ,\mathbf{b}}f$: Definition %
\ref{Psi op}

\item $\mathsf{P}_{\mathcal{B}}^{\mu }f$; (\ref{Haar proj})
\end{enumerate}

\subsubsection{Section 10}

\begin{enumerate}
\item $\mathsf{Q}_{\mathcal{C}_{F}^{\mathcal{G},\limfunc{shift}};M}^{\omega ,%
\mathbf{b}^{\ast }}$; (\ref{def pseudo rest})

\item $\mathcal{C}_{F}^{\mathcal{G},\limfunc{shift}}$, $\mathcal{C}_{F}^{%
\mathcal{G},\limfunc{shift}};K$; (\ref{def shift cor rest})

\item $J\Subset _{\mathbf{\rho },\varepsilon }K$; (\ref{def deep embed})

\item $\mathcal{M}_{\left( \mathbf{\rho },\varepsilon \right) -\limfunc{deep}%
,\mathcal{G}}\left( K\right) $, $\mathcal{M}_{\left( \mathbf{\rho }%
,\varepsilon \right) -\limfunc{deep},\mathcal{D}}\left( K\right) $, $%
\mathcal{W}\left( K\right) $; (\ref{def M_r-deep})

\item augmented dyadic grid $\mathcal{AD}$; Definition \ref{def dyadic}

\item $\mathsf{Q}_{K}^{\omega ,\mathbf{b}^{\ast }}$; (\ref{large pseudo})

\item $\mathcal{E}_{2}^{\alpha ,\func{Whitney}\limfunc{partial}}$; (\ref%
{plug})

\item $A_{2}^{\alpha ,\limfunc{energy}}$, $A_{2}^{\alpha ,\ast ,\limfunc{%
energy}}$; (\ref{def energy A2})

\item $\mathcal{E}_{2}^{\alpha ,\func{Whitney}\limfunc{plug}}$; (\ref{def
deep plug})

\item $\mu $; (\ref{def mu n})

\item $\mathsf{P}_{F,K}^{\omega ,\mathbf{b}^{\ast }}\equiv \mathsf{P}_{%
\mathcal{C}_{F}^{\mathcal{G},\limfunc{shift}};K}^{\omega ,\mathbf{b}^{\ast
}} $; (\ref{def F,K})

\item $\mathbf{Local}\left( I\right) $, $\overline{\mu }$; (\ref{def local
forward})

\item $\mathcal{J}^{\ast }$; (\ref{def J*})

\item $\mathbf{Back}\left( \widehat{I}\right) $; (\ref{e.t2 n'})

\item $U_{s}$; (\ref{def Us})

\item $B\left( M,M^{\prime }\right) $; (\ref{def BMM'})

\item $T_{s}^{\limfunc{intersection}}$; (\ref{def Tints})

\item $T_{s}^{\limfunc{proximal}}$; (\ref{def Tproxs})

\item $T_{s}^{\limfunc{difference}}$; (\ref{def Tdiffs})

\item $\overset{\ast }{\sum }$; Notation \ref{Sum *}

\item $\mathcal{W}_{M}^{s}$; (\ref{def WMs})

\item $\mathcal{W}_{M}^{s,\ell }$; (\ref{def WMsell})
\end{enumerate}


\begin{thebibliography}{AuHoMuTaTh}
\bibitem[AAAHK]{AAAHK} \textsc{M. A. Alfonseca, P. Auscher, A. Axelsson, S.
Hofmann, and S. Kim,} \textit{Analyticity of layer potentials and }$L^{2}$%
\textit{\ solvability of boundary value problems for divergence form
ellliptic equations with complex }$L^{\infty }$\textit{\ coefficients}, 
\texttt{arXiv:0705.0836v1.}

\bibitem[AuHoMcTc]{AuHoLaMcTc} \textsc{Auscher, P., Hofmann, S., Lacey, M.,
McIntosh, A., and Tchamitchian, P.,} \textit{The Solution of the Kato Square
Root Problem for Second Order Elliptic Operators on }$\mathbb{R}^{n}$, Ann.
of Math. \textbf{156} (2002), 633--654.

\bibitem[AuHoMuTaTh]{AuHoMuTaTh} \textsc{P. Auscher, S. Hofmann, C. Muscalu,
T. Tao, and C. Thiele,} \textit{Carleson measures, trees, extrapolation, and 
}$T(b)$\textit{\ theorems}, Publ. Mat. \textbf{46} (2002), no. 2, 257--325.

\bibitem[Chr]{Chr} \textsc{Christ, M.,} \textit{A }$T(b)$\textit{\ theorem
with remarks on analytic capacity and the Cauchy integral}, Colloq. Math. 
\textbf{60/61} (1990), 601--628.

\bibitem[CoJoSe]{CoJoSe} \textsc{R. R. Coifman, P. W. Jones and S. Semmes,} 
\textit{Two elementary proofs of the }$L^{2}$\textit{\ boundedness of Cauchy
integrals on Lipschitz curves,} Journal of the A.M.S. \textbf{2} (1989), p.
553-564.

\bibitem[Dav1]{Dav1} \textsc{David, G.,} \textit{Unrectifiable 1-sets have
vanishing analytic capacity}, Rev. Mat. Iberoamericana \textbf{14} (2)
(1998), 369--479.

\bibitem[Dav2]{Dav2} \textsc{David, G.,} \textit{Analytic capacity, Calder%
\'{o}n-Zygmund operators, and rectifiability}, Publ. Mat. \textbf{43} (1)
(1999), 3--25.

\bibitem[DaJo]{DaJo} \textsc{David, Guy, Journ\'{e}, Jean-Lin,} \textit{A
boundedness criterion for generalized Calder\'{o}n-Zygmund operators,} Ann.
of Math. (2) \textbf{120} (1984), 371--397, MR763911 (85k:42041).

\bibitem[DaJoSe]{DaJoSe} \textsc{David,G.,Journ\'{e},J.-L.,andSemmes,S.,}Op%
\'{e}rateurs de Calder\'{o}n-Zygmund, fonctions para-accr\'{e}tives et
interpolation. Rev. Mat. Iberoamericana 1 (1985), 1--56.

\bibitem[Hof]{Hof} \textsc{Hofmann, S.,} \textit{A proof of the local }$Tb$%
\textit{\ theorem for standard Calder\'{o}n-Zygmund operators}, \texttt{%
arXiv:0705.0840v1}.

\bibitem[HoLaMc]{HoLaMc} \textsc{Hofmann, S., Lacey, M., and McIntosh, A.,}
The solution of the Kato problem for divergence form elliptic operators with
Gaussian heat kernel bounds. Ann. of Math. 156 (2002), 623--631.

\bibitem[HoMc]{HoMc} \textsc{Hofmann, S., and McIntosh, A.,} \textit{The
solution of the Kato problem in two dimensions}, In Proceedings of the
Conference on Harmonic Analysis and PDE (El Escorial, 2000), Publ. Mat.
Extra Vol. (2002), 143--160.

\bibitem[HuMuWh]{HuMuWh} \textsc{R. Hunt, B. Muckenhoupt} \textsc{and R. L.
Wheeden,} \textit{Weighted norm inequalities for the conjugate function and
the Hilbert transform}, Trans. Amer. Math. Soc. \textbf{176} (1973), 227-251.

\bibitem[Hyt]{Hyt} \textsc{Hyt\"{o}nen, Tuomas,} \textit{On Petermichl's
dyadic shift and the Hilbert transform}, C. R. Math. Acad. Sci. Paris 
\textbf{346} (2008), MR2464252.

\bibitem[Hyt2]{Hyt2} \textsc{Hyt\"{o}nen, Tuomas, }\textit{The two weight
inequality for the Hilbert transform with general measures, \texttt{%
arXiv:1312.0843v2}.}

\bibitem[HyMa]{HyMa} \textsc{Hyt\"{o}nen, Tuomas and H. Martikainen,} 
\textit{On general local }$Tb$\textit{\ theorems}, \texttt{arXiv:1011.0642v1.%
}

\bibitem[Lac]{Lac} \textsc{Lacey, Michael T.,}\textit{\ Two weight
inequality for the Hilbert transform: A real variable characterization, II},
Duke Math. J. Volume \textbf{163}, Number 15 (2014), 2821-2840.

\bibitem[LaMa]{LaMa} \textsc{M. T. Lacey and H. Martikainen,} \textit{Local }%
$Tb$\textit{\ theorem with }$L^{2}$\textit{\ testing conditions and general
measures: Calder\'{o}n--Zygmund operators}, \texttt{arXiv:1310.08531v1.}

\bibitem[LaSaUr1]{LaSaUr1} \textsc{Lacey, Michael T., Sawyer, Eric T.,
Uriarte-Tuero, Ignacio,} \textit{A characterization of two weight norm
inequalities for maximal singular integrals with one doubling measure,}
Analysis \& PDE, Vol. \textbf{5} (2012), No. 1, 1-60.

\bibitem[LaSaUr2]{LaSaUr2} \textsc{Lacey, Michael T., Sawyer, Eric T.,
Uriarte-Tuero, Ignacio,} \textit{A Two Weight Inequality for the Hilbert
transform assuming an energy hypothesis, } Journal of Functional Analysis,
Volume \textbf{263} (2012), Issue 2, 305-363.

\bibitem[LaSaShUr2]{LaSaShUr2} \textsc{Lacey, Michael T., Sawyer, Eric T.,
Shen, Chun-Yen, Uriarte-Tuero, Ignacio,} \textit{Two Weight Inequality for
the Hilbert Transform: A Real Variable Characterization, }\texttt{%
arXiv:1201.4319} (2012).

\bibitem[LaSaShUr3]{LaSaShUr3} \textsc{Lacey, Michael T., Sawyer, Eric T.,
Shen, Chun-Yen, Uriarte-Tuero, Ignacio,} \textit{Two weight inequality for
the Hilbert transform: A real variable characterization I}, Duke Math. J,
Volume \textbf{163}, Number 15 (2014), 2795-2820.

\bibitem[LaWi]{LaWi} \textsc{Lacey, Michael T., Wick, Brett D.,} \textit{Two
weight inequalities for Riesz transforms: uniformly full dimension weights}, 
\textit{\texttt{arXiv:1312.6163v3}}.

\bibitem[MaMeVe]{MaMeVe} \textsc{Mattila, P., Melnikov, M., and Verdera, J.,}
\textit{The Cauchy integral, analytic capacity, and uniform rectifiability},
Ann. of Math. (2) \textbf{144} (1) (1996), 127--136.

\bibitem[NTV1]{NTV1} \textsc{Nazarov, F., Treil, S., and Volberg, A.,} 
\textit{Weak type estimates and Cotlar inequalities for Calder\'{o}n-Zygmund
operators on nonhomogeneous spaces}, \textbf{1998} (9) IMRN (1997).

\bibitem[NTV2]{NTV2} \textsc{Nazarov, F., Treil, S. and Volberg, A.,} 
\textit{The }$Tb$\textit{-theorem on non-homogeneous spaces,} Acta Math. 
\textbf{190} (2003), no. 2, MR 1998349 (2005d:30053).

\bibitem[NTV3]{NTV3} \textsc{Nazarov, F., Treil, S., and Volberg, A.,} 
\textit{Accretive system }$Tb$\textit{-theorems on nonhomogeneous spaces},
Duke Math. J. \textbf{113} (2) (2002), 259--312.

\bibitem[NTV4]{NTV4} \textsc{F. Nazarov, S. Treil and A. Volberg,} \textit{%
Two weight estimate for the Hilbert transform and corona decomposition for
non-doubling measures}, preprint (2004) \texttt{arxiv:1003.1596}

\bibitem[Saw]{Saw3} \textsc{E. Sawyer,} \textit{A characterization of two
weight norm inequalities for fractional and Poisson integrals}, Trans.
A.M.S. \textbf{308} (1988), 533-545, MR\{930072 (89d:26009)\}.

\bibitem[SaShUr2]{SaShUr2} \textsc{Sawyer, Eric T., Shen, Chun-Yen,
Uriarte-Tuero, Ignacio,} A \textit{two weight theorem for }$\alpha $\textit{%
-fractional singular integrals with an energy side condition}, \texttt{%
arXiv:1302.5093v8.}

\bibitem[SaShUr3]{SaShUr3} \textsc{Sawyer, Eric T., Shen, Chun-Yen,
Uriarte-Tuero, Ignacio,} \textit{A geometric condition, necessity of energy,
and two weight boundedness of fractional Riesz transforms}, \texttt{%
arXiv:1310.4484v1.}

\bibitem[SaShUr4]{SaShUr4} \textsc{Sawyer, Eric T., Shen, Chun-Yen,
Uriarte-Tuero, Ignacio,} \textit{A note on failure of energy reversal for
classical fractional singular integrals}, IMRN, Volume \textbf{2015}, Issue
19, 9888-9920.

\bibitem[SaShUr5]{SaShUr5} \textsc{Sawyer, Eric T., Shen, Chun-Yen,
Uriarte-Tuero, Ignacio,} A \textit{two weight theorem for }$\alpha $\textit{%
-fractional singular integrals with an energy side condition and quasicube
testing}, \texttt{arXiv:1302.5093v10.}

\bibitem[SaShUr6]{SaShUr6} \textsc{Sawyer, Eric T., Shen, Chun-Yen,
Uriarte-Tuero, Ignacio,} A \textit{two weight theorem for }$\alpha $\textit{%
-fractional singular integrals with an energy side condition, quasicube
testing and common point masses}, \texttt{arXiv:1505.07816v2,v3.}

\bibitem[SaShUr7]{SaShUr7} \textsc{Sawyer, Eric T., Shen, Chun-Yen,
Uriarte-Tuero, Ignacio,} A \textit{two weight theorem for }$\alpha $\textit{%
-fractional singular integrals with an energy side condition}, Revista Mat.
Iberoam. \textbf{32} (2016), no. 1, 79-174.

\bibitem[SaShUr8]{SaShUr8} \textsc{Sawyer, Eric T., Shen, Chun-Yen,
Uriarte-Tuero, Ignacio,} The \textit{two weight }$T1$ \textit{theorem for
fractional Riesz transforms when one measure is supported on a curve}, 
\texttt{arXiv:1505.07822v4}.

\bibitem[SaShUr9]{SaShUr9} \textsc{Sawyer, Eric T., Shen, Chun-Yen,
Uriarte-Tuero, Ignacio,} A \textit{two weight fractional singular integral
theorem with side conditions, energy and }$k$\textit{-energy dispersed,}
Harmonic Analysis, Partial Differential Equations, Complex Analysis, Banach
Spaces, and Operator Theory (Volume 2) (Celebrating Cora Sadosky's life),
Springer 2017 (see also \texttt{arXiv:1603.04332v2}).

\bibitem[SaShUr10]{SaShUr10} \textsc{Sawyer, Eric T., Shen, Chun-Yen,
Uriarte-Tuero, Ignacio,} \textit{A good-}$\lambda $\textit{\ lemma, two
weight }$T1$\textit{\ theorems without weak boundedness, and a two weight
accretive global }$Tb$\textit{\ theorem,} Harmonic Analysis, Partial
Differential Equations and Applications (In Honor of Richard L. Wheeden),
Birkh\"{a}user 2017 (see also \texttt{arXiv:1609.08125v2}).

\bibitem[SaShUr11]{SaShUr11} \textsc{Sawyer, Eric T., Shen, Chun-Yen,
Uriarte-Tuero, Ignacio,} \textit{A counterexample in the theory of Calder%
\'{o}n-Zygmund operators, }\texttt{arXiv:16079.06071v3v1}.

\bibitem[SaWh]{SaWh} \textsc{E. Sawyer and R. L. Wheeden,} Weighted
inequalities for fractional integrals on Euclidean and homogeneous spaces, 
\textit{Amer. J. Math. }\textbf{114} (1992), 813-874.

\bibitem[Ste]{Ste} \textsc{E. M. Stein,} \textit{Harmonic Analysis:
real-variable methods, orthogonality, and oscillatory integrals},\textit{\ }%
Princeton University Press, Princeton, N. J., 1993.

\bibitem[Tol]{Tol} \textsc{Tolsa, X.,} \textit{Painlev\'{e}'s problem and
the semiadditivity of analytic capacity}, Acta Math. \textbf{190} (1)
(2003), 105--149.

\bibitem[Vol]{Vol} \textsc{A. Volberg,} \textit{Calder\'{o}n-Zygmund
capacities and operators on nonhomogeneous spaces,} CBMS Regional Conference
Series in Mathematics (2003), MR\{2019058 (2005c:42015)\}.
\end{thebibliography}
\end{document}